%% file: main.tex
\documentclass[openany,oneside]{book}

\newcommand{\ver}{28}
\newcommand{\ever}{1}
\usepackage{adjustbox}
\usepackage[dvipsnames]{xcolor}

\usepackage{amsmath,amsthm,amssymb,array,changepage,enumerate,etoolbox,float,multicol,pifont,tikz,tikz-cd,url,stmaryrd}
\usepackage{titlesec}

\definecolor{linkcolor}{RGB}{22,128,110}
\usepackage[colorlinks=true,
            linkcolor=linkcolor,
            citecolor=linkcolor,
            urlcolor=linkcolor,
            pdftitle={Multiple Zeta Values},
            pdfauthor={Henrik Bachmann}]{hyperref}
\usepackage[capitalise,nameinlink]{cleveref}
\usepackage[arrow, matrix, curve]{xy}
\usepackage{wrapfig}
\usepackage{geometry,graphicx}
\usepackage{caption}
\usepackage{tikzsymbols}
\usepackage{shuffle}
\usepackage[vcentermath]{youngtab}
\usepackage{ytableau}
\usepackage[utf8]{inputenc}

\usepackage[overlap, CJK]{ruby}
\usepackage{CJKulem}

\usepackage{mathtools}
\mathtoolsset{showonlyrefs}

\usetikzlibrary{decorations.pathreplacing,calc,arrows, fit, arrows.meta,angles,quotes}
\usetikzlibrary{arrows.meta,shadows.blur}
\usetikzlibrary{shapes,decorations}
\usetikzlibrary{shapes.arrows}

\definecolor{mygreen}{RGB}{43,58,55}
\definecolor{myblue}{RGB}{70, 5, 128}
\definecolor{fcolor}{RGB}{0,255,255} 
\definecolor{cec1d24}{RGB}{236,29,36}
\definecolor{cffffff}{RGB}{255,255,255}

\definecolor{realfillcol}{RGB}{240,249,255}   %
\definecolor{reallinecol}{RGB}{  0,114,178}

\definecolor{finitefillcol}{RGB}{241,252,238} %
\definecolor{finitelinecol}{RGB}{ 34,139, 34}

\definecolor{mesfillcol}{RGB}{255,244,220}    %
\definecolor{meslinecol}{RGB}{233,123,  0}

\definecolor{mflinecol}{RGB}{150, 65, 95}

\definecolor{qmzvfillcol}{RGB}{247,241,253}   %
\definecolor{qmzvlinecol}{RGB}{125, 60,152}

\definecolor{relcol}{RGB}{65, 150, 149}  

\definecolor{arrowfillcol}{RGB}{6, 85, 99}

\newcounter{relcnt}
\newenvironment{rel}[1][]{\refstepcounter{relcnt}\par\medskip\noindent    \textbf{\cmini{\therelcnt #1}}\, \rmfamily}{\medskip}
\newcommand\cmini[1]{\tikz[baseline=-.7ex]{\draw[thick,fill=white] (0,0) circle (0.26);
		\node[text=numcol] at (0,0) {{\bf \small #1}};}}
	
\definecolor{algrelcol}{RGB}{251, 255, 240}
\definecolor{linrelcol}{RGB}{251, 255, 240}

\definecolor{numcol}{RGB}{14, 88, 94}
\definecolor{numcol2}{RGB}{137,0,126}
\definecolor{bigcolor}{RGB}{168, 20, 57}
\definecolor{pentacol}{RGB}{98, 196, 147}
\definecolor{relscol}{RGB}{172, 187, 196}
\newcommand{\h}{\mathfrak H}

\usepackage{fancyhdr}
\usepackage{lastpage}
\renewcommand{\chaptermark}[1]{%
  \markboth{\S\thechapter\, #1}{}%
}
\renewcommand{\sectionmark}[1]{%
  \markright{#1}%
}
\newcommand{\zetasim}{\raisebox{0.1em}{\rotatebox[origin=c]{90}{{\color{darkgray}$\scriptstyle\zeta$}}}}

\newcommand{\zetasimm}{\raisebox{0.1em}{\rotatebox[origin=c]{270}{{\color{darkgray}$\scriptstyle\zeta$}}}}

\fancypagestyle{textpage}{
\renewcommand{\headrulewidth}{1pt}
	\lhead{} 
	\chead{\color{mygreen} \leftmark \,{\bf $:$}\,\rightmark} 
	\rhead{} 
	\lfoot{{\small \color{gray} Version \ver \,\,(\today)}} 
	\cfoot{\zetasim \thinspace \thepage \,\zetasimm} 
}

\fancypagestyle{plain}{
\renewcommand{\headrulewidth}{0pt}
	\lhead{} 
	\chead{} 
	\rhead{} 
	\lfoot{{\small \color{gray} Version \ver \,\,(\today)}} 
	\cfoot{\zetasim \thinspace \thepage \,\zetasimm} 
}

\fancypagestyle{titlepage}{
	\renewcommand{\headrulewidth}{0pt}
	\lhead{} 
	\chead{} 
	\rhead{} 
	\lfoot{{\small \color{gray} Version \ver \,\,(\today)}} 
	\cfoot{} 
}

\fancypagestyle{exercisepage}{
	\renewcommand{\headrulewidth}{0pt}
	\lhead{} 
	\chead{} 
	\rhead{} 
	\lfoot{{\small \color{gray} Exercise version \ever \,\,(\today)}} 
	\cfoot{\zetasim \thinspace \thepage \,\zetasimm} 
}

\newtheorem{theorem}{Theorem}[chapter]
\newtheorem{proposition}[theorem]{Proposition}
\newtheorem{definition}[theorem]{Definition}
\newtheorem{lemma}[theorem]{Lemma}
\newtheorem{corollary}[theorem]{Corollary}
\newtheorem{conjecture}[theorem]{Conjecture}
\newtheorem{problem}[theorem]{Open problem}

\newtheorem{exer}{Exercise}[chapter]
\theoremstyle{remark}
\newtheorem{ex}[theorem]{Example}
\newtheorem{exs}[theorem]{Examples}
\newtheorem{remark}[theorem]{Remark}
\newtheorem{constr}[theorem]{Construction}
\numberwithin{equation}{chapter}

\author{Henrik Bachmann}

\input{commands}

\begin{document}
\thispagestyle{titlepage}

\vspace*{-2cm} 
\begin{center}
{\color{mygreen}\huge\bfseries Multiple Zeta Values}\\[0.3cm]
{\Large Henrik Bachmann}\\[0.1cm]
\vspace{-0.1cm}
\PRLsep
\vspace{-0.2cm}
\end{center}

\makeatletter
\renewcommand\tableofcontents{%
    \@starttoc{toc}%
}
\makeatother

{\setlength{\parskip}{0pt}
\makeatletter
\patchcmd{\l@chapter}
  {\vskip 1.0em \@plus\p@}
  {\vskip .6em \@plus.6\p@}
  {}{}
\makeatother
\renewcommand{\contentsname}{}
\tableofcontents
}
\linespread{1.3}

\makeatletter
\renewcommand*\env@matrix[1][*\c@MaxMatrixCols c]{%
  \hskip -\arraycolsep
  \let\@ifnextchar\new@ifnextchar
  \renewcommand{\arraystretch}{0.77}%
  \array{#1}}
\makeatother

\pagestyle{plain}

\input{chap_Preface}

\input{chap_Introduction}

\clearpage
\pagestyle{textpage}

\input{chap_OverviewBasics}

\input{chap_AlgebraicSetup}

\input{chap_FiniteSymMZV}

\input{chap_FamiliesOfRelations}

\input{chap_FormalDoubleZetaSpace}

\input{chap_MES}

\input{chap_FormalSpaces}

\newpage 

\pagestyle{plain}
	
\newpage

\phantomsection
\addcontentsline{toc}{chapter}{References}
\renewcommand{\bibname}{References}

\end{document}

%% file: commands.tex
\renewcommand{\headrule}{\vspace{-0.2cm}\hbox to\headwidth{%
		\color{lightgray}\leaders\hrule height \headrulewidth\hfill}}

\renewcommand{\Re}{\operatorname{Re}}
\newcommand{\Der}{\operatorname{Der}}
\newcommand{\Aut}{\operatorname{Aut}}
\newcommand{\End}{\operatorname{End}}
\newcommand{\Hom}{\operatorname{Hom}}
\newcommand{\reg}{\operatorname{reg}}
\newcommand{\ad}{\operatorname{ad}}
\newcommand{\wt}{\operatorname{wt}}
\newcommand{\dep}{\operatorname{dep}}
\newcommand{\htt}{\operatorname{ht}}
\newcommand{\dropone}{\operatorname{Drop}_1}

\renewcommand{\emph}[1]{{\bf #1}}

\newcommand{\PRLsep}{\noindent\makebox[\linewidth]{\resizebox{\linewidth}{1pt}{$\bullet$}}\bigskip}

\newcommand\quotient[2]{
	\mathchoice
	{\text{\raise1ex\hbox{$#1$}\Big/\lower1ex\hbox{$#2$}}}
	{#1\,/\,#2}
	{#1\,/\,#2}
	{#1\,/\,#2}
}

\newcommand{\pmatr}[1]{\begin{pmatrix}#1\end{pmatrix}}
\newcommand{\sabcd}{\left( \begin{smallmatrix}
a & b \\
c & d
\end{smallmatrix} \right)}
\newcommand{\Q}{\mathbb{Q}}
\newcommand{\Z}{\mathbb{Z}}
\newcommand{\N}{\mathbb{N}}
\newcommand{\C}{\mathbb{C}}
\newcommand{\R}{\mathbb{R}}
\newcommand{\Ha}{\mathbb{H}}

\newcommand{\kk}{{\bf k}}
\newcommand{\kl}{{\bf l}}
\newcommand{\mz}{\mathcal{Z}}
\newcommand{\sh}{\shuffle}
\newcommand{\HH}{\mathfrak{H}}
\newcommand{\qsh}{\ast_\diamond}
\newcommand{\gdiamond}{{\hat\diamond}}
\newcommand{\gqsh}{\hat\ast}
\newcommand{\gshu}{\hat\shuffle}
\newcommand{\len}{\ell}
\newcommand{\ko}{\mathsf{k}}
\newcommand{\QA}{\ko\langle A \rangle}
\newcommand{\ds}{\operatorname{ds}}
\newcommand{\eds}{\mathsf{eds}}

\newcommand{\mes}{\mathcal{E}}
\newcommand{\rmes}{\mathcal{G}}
\newcommand{\aqmf}{\widetilde{\mf}}

\newcommand{\fdes}{\mathcal{D\!E}}
\DeclareRobustCommand{\Gde}{G^{\mathrm{DE}}\genfrac{(}{)}{0pt}{}}
\DeclareRobustCommand{\Pde}{P^{\mathrm{DE}}\genfrac{(}{)}{0pt}{}}

\newcommand{\fmes}{\rmes^{f}}
\newcommand{\Gf}{G^{f}}
\newcommand{\lwt}{\operatorname{lwt}}
\newcommand{\fil}[2]{\operatorname{Fil}^{#1}_{#2}}
\newcommand{\fmesz}{\fil{\lwt}{0}\fmes}
\newcommand{\fames}{\mes^{f}}
\newcommand{\fmz}{\mz^{f}}
\newcommand{\fzeta}{\zeta^{f}}
\newcommand{\fqmf}{\aqmf^{f}}
\newcommand{\fmf}{\mf^{f}}
\newcommand{\fcusp}{\mathsf{S}^{f}}
\newcommand{\fDelta}{\Delta^{f}}
\newcommand{\bist}{\ast_{\mathrm{bi}}}

\newcommand{\gG}{\mathfrak{G}}
\newcommand{\gbr}{\widetilde{\gb}}
\newcommand{\gL}{\mathfrak{L}}
\DeclareRobustCommand{\Gb}{G\genfrac{(}{)}{0pt}{}}

\newcommand{\GD}{\operatorname{D}}
\newcommand{\GW}{\operatorname{W}}
\newcommand{\Gdelta}{\delta}
\newcommand{\sltwo}{\mathfrak{sl}_2}

\newcommand{\g}{\operatorname{g}}
\newcommand{\gs}{\mz_q^\circ}
\newcommand{\gsk}{\mz_{q,\leq k}^\circ}
\newcommand{\gsev}{\mz_q^{\text{ev}}}
\newcommand{\gstwo}{\mz_q^{\geq 2}}
\newcommand{\gsadm}{\mz_q^{\text{adm}}}
\newcommand{\gsadmk}{\mz_{q,\leq k}^{\text{adm}}}
\newcommand{\gsadmj}{\mz_{q,\leq j}^{\text{adm}}}

\newcommand{\li}{\operatorname{Li}}

\newcommand{\ggen}{\mathfrak{g}}
\newcommand{\zm}{\operatorname{Z}}
\newcommand{\qdq}{q \frac{d}{dq}}

\newcommand{\fp}{\Z / p \Z}
\newcommand{\fpq}{\Z_{(p)}[q] / ([p]_q)}
\newcommand{\fmza}{{\mz_{\ma}^{f}}}
\newcommand{\fza}{\zeta_{\ma}^{f}}
\newcommand{\fzs}{\zeta_{\ms}^{f}}
\newcommand{\fZ}{\mathfrak{z}^{f}}
\newcommand{\ms}{\mathcal{S}}
\newcommand{\ma}{\mathcal{A}}
\newcommand{\mq}{\mathcal{Q}}

\newcommand{\zs}{\zeta_{\ms}}
\newcommand{\za}{\zeta_{\ma}}
\newcommand{\zq}{\zeta_{\mq}}
\newcommand{\mza}{\mz^{\ma}}

\newcommand{\mf}{\mathcal{M}}
\newcommand{\qmf}{\widetilde{\mf}^\Q}
\newcommand{\aG}{\mathbb{G}}
\newcommand{\ag}{\widehat{\g}}

\newcommand{\RR}{{\color{red}{R}}}
\newcommand{\UU}{{\color{blue}{U}}}

\newcommand{\SL}{\operatorname{SL}}
\newcommand{\SLZ}{\operatorname{SL}_2(\Z)}
\newcommand{\abcd}{\begin{pmatrix} a & b \\ c & d \end{pmatrix}}

\newcommand{\mabcd}{\frac{a\tau +b}{c\tau+d}}
\DeclareMathOperator{\MZB}{\mz[\pi i][\![q]\!]}

\newcommand{\hz}{\h^2}
\newcommand{\ON}{\overline{\N}}
\newcommand{\CC}{\mathcal{C}}
\newcommand{\HHH}{\widehat{\HH}}
\newcommand{\shq}{\shuffle_q}
\newcommand{\GL}{\operatorname{GL}}
\newcommand{\PGL}{\operatorname{PGL}}
\newcommand{\emptyword}{1}
\newcommand{\vsmall}{\rotatebox[origin=c]{-90}{$<$}}

\newcommand{\dz}{\mathcal{D}}
\newcommand{\gz}{\mathfrak{Z}} 
\newcommand{\gb}{\mathfrak{b}} 
\newcommand{\gp}{\mathfrak{P}} 
\newcommand{\gr}{\mathfrak{R}} 
\newcommand{\gt}{\mathfrak{T}} 
\newcommand{\gh}{\mathfrak{h}} 

\DeclareRobustCommand{\bi}{\genfrac{(}{)}{0pt}{}}

\DeclareRobustCommand{\mt}{\ggen \genfrac{(}{)}{0pt}{}}
\DeclareRobustCommand{\mb}{\g\genfrac{(}{)}{0pt}{}}
\DeclareRobustCommand{\ai}{\genfrac{[}{]}{0pt}{}}
\newcommand{\azbi}{A^{\text{bi}}_z}

\newcommand{\mtt}[3] {{\mathcal{T}\begin{psmallmatrix}
			#1\\
			#2 \\
			#3
\end{psmallmatrix}}}

\newcommand*{\braceme}[6][]{%
	\draw[
	shift={(#3:#2)},
	right to reversed-right to reversed,
	shorten >=-.75\pgflinewidth,
	#1
	] (0,0)
	arc[radius=#2, start angle=#3, end angle=#3+(#4-#3)/2] node[rotate=#3+(#4-#3)/2-90,above=2pt] (#5) {#6};
	\draw[
	shift={({#3+(#4-#3)/2}:#2)},
	left to reversed-left to reversed,
	shorten <=-.75\pgflinewidth,
	#1
	] (0,0)
	arc[radius=#2, start angle=#3+(#4-#3)/2, end angle=#4];
}

\definecolor{mycolor}{RGB}{194, 8, 88}
\newcommand{\todo}[1]{%
    \message{LaTeX Warning: Unfinished work on input line \the\inputlineno}%
    {\color{mycolor}\big\{\Coffeecup\,{\bf Todo:} #1\,\Coffeecup\big\}}%
}

\newcommand{\II}{\mathbb{I}}
\newcommand{\id}{\operatorname{id}}
\newcommand{\HHCK}{\mathcal{H}}

\newcommand{\fild}{\operatorname{Fil}^{\operatorname{D}}}
\newcommand{\grd}{\operatorname{gr}^{\operatorname{D}}}
\DeclareMathOperator{\odd}{\mathsf{O}}
\DeclareMathOperator{\ev}{\mathsf{E}}
\DeclareMathOperator{\Sx}{\mathsf{S}}

\newcommand{\piod}{\pi^{\text{od}}}

%% file: chap_Preface.tex
\chapter*{Preface}
\addcontentsline{toc}{chapter}{Preface}
\vspace{-2\baselineskip}
\begingroup
\linespread{1.0}\selectfont

I first became interested in multiple zeta values at the beginning of my master's studies at the University of Hamburg around fifteen years ago. At that time my supervisor Ulf K\"uhn suggested that I read Zudilin's paper \cite{Zu3} and the paper by Gangl, Kaneko and Zagier \cite{GKZ}. These two papers, and also their authors and my supervisor, who all have their own unique way of doing mathematics, had a strong influence on my academic path. Since then, I have been very lucky to meet all of the authors several times. In particular, the paper by Gangl, Kaneko and Zagier inspired my master's thesis on the Fourier expansion of multiple Eisenstein series and their connection to modular forms. This became the starting point for my interest in $q$-analogues of multiple zeta values, which I then studied further during my PhD.

After my PhD, I was fortunate to obtain a JSPS Postdoctoral Fellowship with Hidekazu Furusho at Nagoya University. This brought me to Japan, the country of multiple zeta values, which seems to have by far the largest community working on the subject. The many kind mathematicians in this community gave me many new insights into the field. I learned about many more aspects of multiple zeta values and always tried to connect them back to multiple Eisenstein series, $q$-analogues and modular forms. More recently, finite multiple zeta values also entered this picture, but their connections in these directions are only beginning to appear.

Since coming to Japan, I have enjoyed telling people abroad about the beautiful mathematics of multiple zeta values, the large community working on them here and, of course, the country itself. With these notes, I would like to continue this. The goal is to bring together as much as possible of what is known about these objects and, in particular, to emphasize the structures and connections which put all of them into a common framework. I hope that they will be useful for future generations.

These notes emerged from three courses given at Nagoya University: {\itshape Multiple zeta values and modular forms} (Spring 2020), {\itshape $q$-analogues and finite multiple zeta values} (Spring 2022) and {\itshape Multiple zeta values} (Spring 2025). The first course was given online during the COVID-19 pandemic, and videos of its fifteen lectures are available at \url{https://www.henrikbachmann.com/mzv2020.html}.

These notes are still under construction and will be updated in the future. Comments, corrections and suggestions are very welcome.

Finally, I would like to thank the following people for general discussions on the topic, feedback, comments and support, and for pointing out typos in these notes and/or their earlier versions: Takumi Anzawa, Vic Austen, Olivier Bouillot, Benjamin Brindle, David Broadhurst, Johannes Broedel, Francis Brown, Annika Burmester, Steven Charlton, Niclas Confurius, Cl\'ement Dupont, Kurusch Ebrahimi-Fard, Hidekazu Furusho, Herbert Gangl, Runxuan Gao, Ryotaro Harada, Minoru Hirose, Michael Hoffman, Kentaro Ihara, Jan-Willem van Ittersum, David Jarossay, Shin-ya Kadota, Masanobu Kaneko, Hayato Kanno, Hanamichi Kawamura, Naho Kawasaki, Adam Keilthy, Katsumi Kina, Nao Komiyama, Ulf K\"uhn, Yunyang Luo, Takumi Maesaka, Dominique Manchon, Nils Matthes, Hideki Murahara, Maki Nakasuji, Yasuo Ohno, Masataka Ono, Erik Panzer, Risan, Simon Rutard, Shingo Saito, Nobuo Sato, Oliver Schlotterer, Leila Schneps, Shin-ichiro Seki, Takeshi Shinohara, Yuta Suzuki, Yoshihiro Takeyama, Tatsushi Tanaka, Koji Tasaka, Hirofumi Tsumura, Can Turan, Pierre-Emmanuel Wulfman, Khalef Yaddaden, Shuji Yamamoto, Yoshinori Yamasaki, Jinbo Yu, Don Zagier, Federico Zerbini, Jianqiang Zhao, Jiayu Zhao and Wadim Zudilin.

\vfill

\noindent\hfill
\begin{tabular}{r}
Henrik Bachmann\\
Nagoya, August 2026
\end{tabular}
\endgroup

%% file: chap_Introduction.tex
\chapter*{Introduction}
\addcontentsline{toc}{chapter}{Introduction}

These notes study various aspects of \emph{multiple zeta values (MZVs)}, which are special values of generalizations of the Riemann zeta function defined for integers $k_1\geq 2, k_2,\dots,k_r \geq 1$ by
\begin{align}\label{eq:defmzvintro}
    \zeta(k_1,\dots,k_r) = \sum_{m_1>\dots>m_r>0} \frac{1}{m_1^{k_1}\cdots m_r^{k_r}} \in \R.
\end{align}
These real numbers appear in several areas of mathematics and theoretical physics. They have been studied since at least the time of Euler ($\sim$1740s), who discovered many of their algebraic properties. After seemingly being forgotten for over 200 years, multiple zeta values were rediscovered by many mathematicians and theoretical physicists beginning in the 1980s in several different contexts, such as modular forms, mixed Tate motives, quantum groups, moduli spaces of genus-zero curves with marked points \cite{Br1}, scattering amplitudes, resurgence theory, etc.

The goal of these notes is to provide an overview of multiple zeta values and some of their recently introduced variants. We will deal with objects which live in four different worlds visualized in the diagram below: {\color{reallinecol}Real numbers $\R$}, {\color{finitelinecol} the ring of ``poor man's adeles'' $\mathcal{A}$}, {\color{meslinecol}holomorphic functions in the complex upper half-plane $\mathcal{O}(\Ha)$}, and {\color{qmzvlinecol}$q$-series with rational coefficients $\Q\llbracket q \rrbracket$}.
\begin{center}
  \input{TikZ/4worlds}
\end{center}
We will begin in Section~\ref{subsec:riemann} by revisiting the classical Riemann zeta function and stating in Section~\ref{subsec:mzv} some of the folklore conjectures in the field of (multiple) zeta values.
Some simple examples will also show how to obtain relations among multiple zeta values using the finite double shuffle relations.
After this, we will introduce in Section~\ref{subsec:fmzv} \emph{finite multiple zeta values}, which are defined by reducing the defining sums \eqref{eq:defmzvintro} modulo primes~$p$ and then collecting all these values together. Finite multiple zeta values are not real numbers, but will reside in the ``poor man's adeles'' ring $\mathcal{A}$.
Nevertheless, we will see that they have a surprising connection to multiple zeta values via the Kaneko--Zagier conjecture \cite[Main Conjecture]{KZ}.

More precisely, we will introduce \emph{symmetric multiple zeta values}, which are constructed out of multiple zeta values, and which conjecturally satisfy the same relations as finite multiple zeta values.
All relations among these objects are conjectured to arise from an analogue of the double shuffle relations.
Next, we will discuss \emph{modular forms} and Eisenstein series in Section~\ref{subsec:mfandbk}.
Here modular forms are mysteriously connected to relations among multiple zeta values.
Some of this mystery will be explained in Section~\ref{subsec:mesoverview} by introducing a hybrid of multiple zeta values and Eisenstein series: \emph{multiple Eisenstein series}.

Finally, we will consider in Section~\ref{subsec:qmzv} another variant of MZVs: \emph{$q$-analogues of multiple zeta values ($q$MZVs)}.
These are $q$-series which degenerate to multiple zeta values as $q \rightarrow 1$. By viewing $q$ not just as a parameter but as a function in a complex variable $\tau$, via $q=e^{2\pi i \tau}$, we will see that $q$MZVs give another natural bridge between the theory of modular forms and multiple zeta values. 
For both $q$MZVs and multiple Eisenstein series, we will also discuss what the analogues of the double shuffle relations are.
After this overview, we will go into more detail regarding the algebraic structures of all these objects in Chapter~\ref{sec:algebraicsetup}.
For this, we will study quasi-shuffle algebras, which form the backbone of all the objects mentioned above.

After developing this algebraic setup, we return in Chapter~\ref{sec:fmzv} to finite and symmetric multiple zeta values and then discuss in Chapter~\ref{sec:families} further families of relations and variants, including Schur multiple zeta values. The formal double zeta space in Chapter~\ref{sec:mdandmzv} then makes the connection with modular forms precise in depth two, and its last section shows that, with a shift in depth, the same period polynomials give relations among finite multiple zeta values. We next study multiple Eisenstein series in Chapter~\ref{sec:mes}, including their regularizations, their derivatives and their $\mathfrak{sl}_2$-algebra structure. Formal multiple zeta values, the formal double Eisenstein space and formal multiple Eisenstein series in Sections~\ref{sec:fmz}, \ref{sec:fdes} and~\ref{sec:fmes} will then bring the double shuffle relations and the derivative structures together and connect three of the four worlds in the diagram above. Section~\ref{sec:cmes} of the final chapter constructs combinatorial multiple Eisenstein series. These give an
explicit $q$-series realization of the formal multiple Eisenstein series,
interpolate between rational solutions of the extended double shuffle
relations and multiple zeta values, and make it possible to compare the
different depth-two constructions precisely. Finally, formal finite multiple zeta values in Section~\ref{sec:formalfinitemzv} give a universal version of the finite and symmetric worlds and a formal version of the Kaneko--Zagier map. A reduced quotient of the formal double zeta space also gives a new connection between period polynomials and finite multiple zeta values.

There are various other points of view from which one can study multiple zeta values that are not treated in these notes.
For more details on multiple zeta values, we refer in particular to the books and lecture notes of Arakawa--Kaneko \cite{AK}, Kaneko \cite{K4}, Burgos--Fres\'an \cite{BF}, Waldschmidt \cite{W}, Zhao \cite{Zh2}, Zudilin \cite{Zu6}, and to Hoffman's collection of research papers on multiple zeta values \cite{H5}.
For example, a somewhat dual approach to the study of multiple zeta values is given by the work of Racinet. See \cite{Bu1} for a good overview and detailed explanation of how Racinet's setup can be transferred to $q$-analogues.

Let us finish the introduction with a small warning about conventions. Unfortunately, the multiple zeta value community does not agree on the order of an index. Throughout these notes, we use the convention of \cite{GKZ,IKZ}, i.e. the summation order in \eqref{eq:defmzvintro}. Thus for us $\zeta(2,1)=\zeta(3)$, whereas with the opposite convention the same identity is written as $\zeta(1,2)=\zeta(3)$. To make things even more confusing, some of the authors of these papers use both conventions in different works. I apologize to all readers using the other convention. When comparing formulas, one therefore always needs to check the order of the summation variables.

%% file: TikZ/4worlds.tex
\tikzset{%
  block/.style={draw, top color=white!90!#1, bottom color=white, shading angle=45, text width=18pt, text centered, rounded corners, minimum height=10pt},
}

\newcommand{\arr}[3]{
\fill[arrowfillcol] #1 circle[radius = 3pt];
  \draw[
    -{Triangle[length=4pt,width=6pt]},   
    line width = 2pt,
    draw = arrowfillcol
  ] 
    #1 -- #2 
    node[midway, sloped,
         fill = white,
         draw = arrowfillcol,
         rounded corners = 2pt,
         line width = 0.5pt,
         inner sep = 3pt,
         text = arrowfillcol]{#3};
}
\newcommand{\arrr}[3]{
  \draw[
    arrows={{Triangle[length=4pt,width=6pt]}-{Triangle[length=4pt,width=6pt]}},   
    line width = 2pt,
    draw = arrowfillcol
  ] 
    #1 -- #2 
    node[midway, sloped,
         fill = white,
         draw = arrowfillcol,
         rounded corners = 2pt,
         line width = 0.5pt,
         inner sep = 3pt,
         text = arrowfillcol]{#3};
}

\newcommand{\arl}[3]{
\fill[arrowfillcol] #1 circle[radius = 3pt];
\fill[arrowfillcol] #2 circle[radius = 5pt];
  \draw[
    -,   
    line width = 2pt,
    draw = arrowfillcol
  ] 
    #1 -- #2 
    node[midway, sloped,
         fill = white,
         draw = arrowfillcol,
         rounded corners = 2pt,
         line width = 0.5pt,
         inner sep = 3pt,
         text = arrowfillcol]{#3};
}

\begin{tikzpicture}[scale=0.27, >=Latex,every node/.style={scale=0.77}]

\def\linewidth{1pt}
\def\linewidthbig{2pt}
\def\arrowwidth{1pt}

\def\circrad{0.6}

\coordinate (MZV) at (31,23);

\coordinate (SMZV) at (30,34);

\coordinate (FMZV) at (10,34);

\coordinate (MES) at (30,8);

\coordinate (MF) at (26.5,15);

\coordinate (qMZV) at (10,15);

\coordinate (qg) at (10,8);

\draw[rounded corners,line width=\linewidthbig,color=reallinecol,fill=realfillcol,blur shadow] (21,21) rectangle (39,39) {};
\node at (37,37) [color=reallinecol]{{\bf \Huge $\R$}};

\draw[rounded corners,line width=\linewidthbig,color=finitelinecol,fill=finitefillcol,blur shadow] (1,21) rectangle (19,39) {};
\node at (3,37) [color=finitelinecol]{{\bf \Huge $\mathcal{A}$}};

\draw[rounded corners,line width=\linewidthbig,color=meslinecol,fill=mesfillcol,blur shadow] (21,1) rectangle (39,19) {};

\draw[
    arrows={-{Triangle[length=4pt,width=6pt]}},   
    line width = 1pt,
    dashed,
    draw = mflinecol
  ] (22,15) -- (8,22);

\draw[rounded corners,line width=\linewidthbig,color=qmzvlinecol,fill=qmzvfillcol,blur shadow] (1,1) rectangle (19,19) {};
\node at (4,3) [color=qmzvlinecol]{{\bf \Huge $\Q\llbracket q \rrbracket$}};

\node at (36,3) [color=meslinecol]{{\bf \Huge $\mathcal{O}(\Ha)$}};
\draw[rounded corners,line width=1pt,dashed,color=mflinecol] (22,13) rectangle (31,18) {};

\node at ($(MZV)+(0,1.5)$) [align=center]{{\bf Multiple Zeta Values (MZV)}};
\node at (MZV) [align=center]{{\bf \Large $\zeta(k_1,\dots,k_r)$}};

\node at ($(SMZV)+(0,1.5)$) [align=center]{{\bf Symmetric MZV}};
\node at (SMZV) [align=center]{{\bf \Large $\zs(k_1,\dots,k_r)$}};

\node at ($(FMZV)+(0,1.5)$) [align=center]{{\bf Finite MZV}};
\node at (FMZV) [align=center]{{\bf \Large $\za(k_1,\dots,k_r)$}};

\node at ($(MES)+(0,1.5)$) [align=center]{{\bf Multiple Eisenstein Series}};
\node at (MES) [align=center]{{\bf \Large $\aG_{k_1,\dots,k_r}(\tau)$}};
\node at (qg) [align=center]{{\bf \Large $\g(k_1,\dots,k_r)$}};

\node at (28,14.5) [align=center]{{\bf \Large $\aG_{2k}$}};
\node at (24,14.5) [align=center]{{\bf \Large $\Delta$}};

\node at ($(qMZV)+(0,1.7)$) [align=center]{{\bf $q$-analogues of MZV}};
\node at (qMZV) [align=center]{{\bf \Large $\zeta_q(k_1,\dots,k_r ; Q_1 , \dots , Q_r)$}};

\node at ($(MF)+(0,1.5)$) [align=center]{{\color{mflinecol} \bf Modular forms}};

\draw[
    arrows={-{Triangle[length=4pt,width=6pt]}},   
    line width = 1pt,
    dashed,
    draw = mflinecol
  ] (22.2,18) -- (22.2,26.5);

\node at (30,27) [align=center]{{\scriptsize \color{mflinecol} $168 \zeta(5,7)+150 \zeta(7,5)+28 \zeta(9,3) = \frac{5197}{691} \zeta(12)$}};    

\node at (8,24) [align=center]{{\scriptsize \color{mflinecol} $2\,\zeta_{\mathcal{A}}(1, 2, 1, 8) -18 \,\zeta_{\mathcal{A}}(1, 4, 1, 6)$}};    
\node at (8,23) [align=center]{{\scriptsize \color{mflinecol} $= 9\, \zeta_{\mathcal{A}}(1, 6, 1, 4) + 16 \,\zeta_{\mathcal{A}}(1, 8, 1, 2)$}};    

\arrr{($(MES)+(-6.5,0)$)}{($(qg)+(6.5,0)$)}{$q = e^{2\pi i \tau}$}

\arr{($(MES)+(3,3)$)}{(33,22)}{$\tau \rightarrow i \infty $}

\arr{($(qMZV)+(6,2)$)}{(26,23)}{$q \rightarrow 1$}

\arl{($(qMZV)+(4,3)$)}{(20,29)}{$q = e^{\frac{2\pi i}{n}}$}
\arr{(20,29)}{(25,33)}{$n \rightarrow \infty$}
\arr{(20,29)}{(15,33)}{$n$ prime}

\arrr{($(FMZV)+(4.5,0)$)}{($(SMZV)-(4.5,0)$)}{{\small Kaneko-Zagier Conj.}}

\node at (10,11) [align=center]{{\color{relcol} $\g_{2,1} = q^3 + 2q^4 + 6q^5 + 7q^6 + \dots$}};

\node at (10,5.5) [align=center]{{\color{relcol} $\g(6) = 6 \g(3,3) - 3 \g(4,2) + \frac{\g(4)}{4}  - \frac{\g(2)}{180}$}};
\node at (30,5.5) [align=center]{{\color{relcol} $\aG_{6} = 6 \aG_{3,3} - 3 \aG_{4,2}$}};
\node at (31,32) [align=center]{{\color{relcol} $2\zs(4,1)  + \zs(3,2) = 0$}};
\node at (9,32) [align=center]{{\color{relcol} $2\za(4,1)  + \za(3,2) = 0$}};
\node at (9.9,29) [align=center]{{\color{relcol} $\za(2,1) = (0, 1, 1, 3, 4, 5, \dots)$}};
\node at (10,27.5) [align=center]{{\color{relcol} $= 3 Z(3) \overset{?}{\neq} 0 $}};

\node at (30,30) [align=center]{{\color{relcol} $\zeta(5) \equiv \frac{2}{5}\zs(4,1)+\frac{2}{15} \zs(2,2,1)$}};

\end{tikzpicture}

%% file: chap_OverviewBasics.tex
\chapter{Overview \& Basics}\label{sec:overview}
In this chapter, we will give an overview of the values of the Riemann zeta function, multiple zeta values and some of their variants. These include finite and symmetric multiple zeta values, multiple Eisenstein series, $q$-analogues and Schur multiple zeta values. We will also mention some of the main conjectures concerning their structure and give a glimpse of their connections to modular forms. The general picture of these concepts will be discussed in detail in the following chapters.

\section{The values of the Riemann zeta function}\label{subsec:riemann}

The Riemann zeta function is defined for a complex variable $s\in \C$ with $\Re(s)>1$ by 
\begin{align}\label{def:riemannzeta}
\zeta(s) = \sum_{m>0}\frac{1}{m^s}\,.
\end{align}
This function appears in various fields of mathematics and theoretical physics and can be studied from various points of view.\\

\begin{minipage}[t]{0.4\textwidth}
    \centering\raisebox{\dimexpr \topskip-\height}{%
        \includegraphics[width=0.8\textwidth]{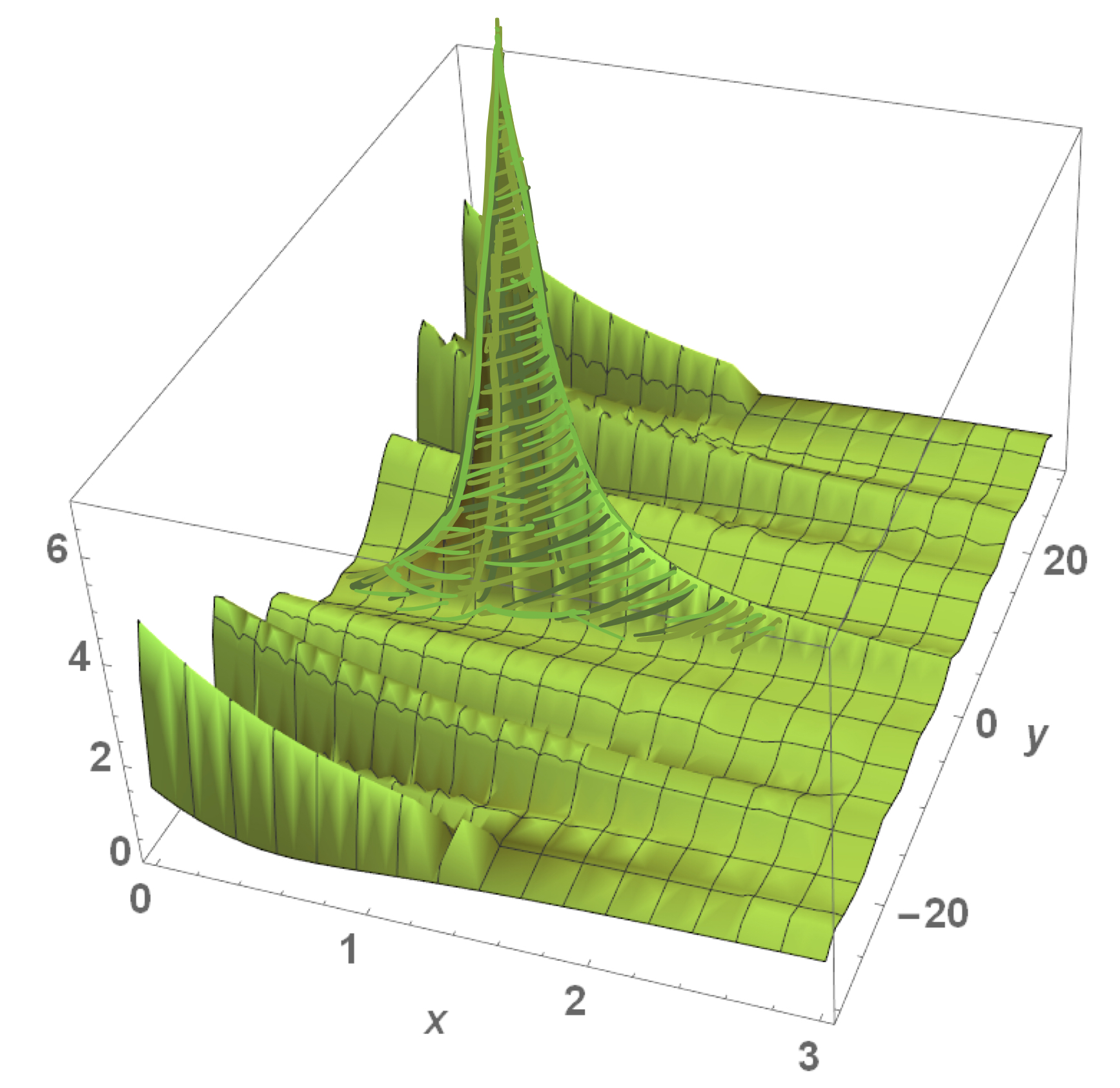}}
    
    The graph of $|\zeta(x+iy)|$ near \\the pole at $x+iy=1$.
\end{minipage}\hfill
\begin{minipage}[t]{0.6\textwidth}
    For example, it is well known that the Riemann zeta function can be analytically continued to the whole complex plane with a simple pole at $s=1$. 
    Although $\zeta(s)$ was already considered by L. Euler (1707 -- 1783), it was named after B. Riemann (1826 -- 1866), who proved its meromorphic continuation and functional equation and established a relation between its zeros and the distribution of prime numbers. 
    
    In particular, he gave his famous conjecture on the location of the zeros of the Riemann zeta function, stating that besides the trivial zeros at $s=-2,-4,-6,\dots$ all other zeros have real part  $\frac{1}{2}$.
\end{minipage}
\vspace{0.3cm}

The connection to prime numbers is given by the following product formula, which, to make things fair again, was named after Euler (Euler product formula) 
\begin{align*}
\zeta(s) = \prod_{p \text{ prime}}\frac{1}{1-p^{-s}}\,.\qquad (\Re(s)>1)
\end{align*}
In these notes, we will not study these analytic aspects but rather will be interested in the values of $\zeta(k)$, when $k\in \Z_{\geq 2}$ is a positive integer.  The first result is the famous formula by Euler for $\zeta(2)$, which states that 
\begin{align*}
\zeta(2) = \sum_{m>0}\frac{1}{m^2} = 1 + \frac{1}{4}+\frac{1}{9}+\frac{1}{16} + \dots = \frac{\pi^2}{6}\,.
\end{align*}

In general, Euler proved that $\zeta(2m)$ is always a rational multiple of $\pi^{2m}$, and he gave the following explicit formula in terms of Bernoulli numbers.
\begin{proposition}[Euler, 1735 {\cite{Eu}}]\label{prop:euler} For all $m\in \Z_{\geq 1}$ we have
    \begin{align*}
    \zeta(2m) = -\frac{B_{2m}}{2 (2m)!} (2\pi i)^{2m} \in \Q \pi^{2m}\,, 
    \end{align*}
    where $B_n$ denotes the $n$-th \emph{Seki-Bernoulli number\footnote{These numbers are often just called Bernoulli numbers. Due to the author's affiliation to Japan he prefers using a historically more fair/accurate naming.}} defined\footnote{A more natural definition takes the Seki-Bernoulli numbers \(B_n\) to be the constant terms of the Bernoulli polynomials \(B_n(x)\), uniquely characterised by the unit-interval identity \(\int_x^{x+1} B_n(t)\,dt = x^{\,n}\).}
 by the generating series
    \begin{align}\label{eq:bn}
    \sum_{n=0}^\infty B_n \frac{x^n}{n!} = \frac{x}{e^x-1}\,.
    \end{align}
\end{proposition}
\begin{proof}
    There are various ways to prove this fact and we will give the original approach due to Euler. First, consider the Weierstrass product of the sine function
    \begin{align}\label{eq:sinproduct}
    \frac{\sin(\pi x)}{\pi x} = \prod_{n\geq 1} \left(1 - \frac{x^2}{n^2}\right)\,.
    \end{align}
    For $x\in \C \backslash \Z$ we can take its logarithmic derivative to obtain the partial fraction expansion of the cotangent
    \begin{align*}
    \pi \cot(\pi x) = \frac{1}{x} + \sum_{n\geq 1}\left(\frac{1}{x+n} + \frac{1}{x-n}\right) \,.
    \end{align*}
    Expanding the right hand side in a geometric series $\frac{1}{x+n} = \frac{1}{n} \frac{1}{1+\frac{x}{n}}=\frac{1}{n}\sum_{m\geq 0} (-1)^m \frac{x^m}{n^m}$ gives 
    \begin{align*}
    \frac{1}{x} + \sum_{n\geq 1}\left(\frac{1}{x+n} + \frac{1}{x-n}\right) = \frac{1}{x} - \sum_{m=1}^\infty 2 \zeta(2m) x^{2m-1}\,.
    \end{align*}
    On the other hand the left hand side can be evaluated as 
    \begin{align*}
    \pi \cot(\pi x) = \pi i \frac{e^{\pi i x}+e^{-\pi i x}}{e^{\pi i x} - e^{-\pi i x}} =\pi i \left(1 + \frac{2}{e^{2\pi i x}- 1}\right) \overset{\eqref{eq:bn}}{=} \frac{1}{x} + \sum_{m=1}^\infty \frac{B_{2m} (2\pi i)^{2m}}{(2m)!} x^{2m-1}\,,
    \end{align*}
    where in the last equality we used $B_1=-\frac{1}{2}$.
\end{proof}
The first explicit values for $\zeta(2m)$ are given by the following
\begin{align*}
\zeta(2)=\frac{\pi^2}{6}\,,\quad \zeta(4)=\frac{\pi^4}{90}\,,\quad \zeta(6)=\frac{\pi^6}{945}\,,\quad \zeta(8)=\frac{\pi^8}{9450}\,,\quad \zeta(10)=\frac{\pi^{10}}{93555}\,,\quad \zeta(12)=\frac{691 \pi^{12}}{638512875}\,.
\end{align*}
Since $\pi$ is transcendental (Lindemann, 1882 \cite{Lin}), Proposition \ref{prop:euler} gives the only family of polynomial relations among even zeta values. On the other hand, one does not expect polynomial relations among odd zetas. This is part of the following folklore conjecture.
\begin{conjecture}\label{conj:singlezetaalgindep}
    The numbers $\pi^2, \zeta(3), \zeta(5), \zeta(7),\,\dots$ are algebraically independent over $\Q$.
\end{conjecture}
So far there is not much known regarding this conjecture. For the odd zeta values, the following theorem gives an overview of the known facts.
\begin{theorem}
    \begin{enumerate}[\textup{(}i\textup{)}] 
        \item $\zeta(3)$ is irrational.   \textup{(Ap\'ery, 1978 \cite{Ap})}
        \item For every sufficiently large odd integer $s$ we have
        \[
        \dim_{\Q}\langle 1,\zeta(3),\zeta(5),\ldots,\zeta(s)\rangle_{\Q}
        \geq 0.21\sqrt{\frac{s}{\log s}}\,.
        \]
        In particular, infinitely many of the values $\zeta(2m+1)$ are irrational. \textup{(Fischler, 2026 \cite{Fi})}
        \item At least one of the values $\zeta(5),\zeta(7),\zeta(9)$ and $\zeta(11)$ is irrational. \textup{(Zudilin, 2001 \cite{Zu1})}
        \item At least two of the values $\zeta(5),\zeta(7),\ldots,\zeta(35)$ are irrational. \textup{(Lai--Zhou, 2022 \cite{LZ})}
    \end{enumerate}
\end{theorem}

\begin{remark}
    The estimate in (ii) improves the earlier logarithmic bound of Ball and Rivoal \cite{BR}. In a different direction, Lai proved that for every sufficiently large even integer $s$, at least $1.284\sqrt{s/\log s}$ among $\zeta(3),\zeta(5),\ldots,\zeta(s-1)$ are irrational \cite{Lai}. Such results count irrational values and do not assert their linear independence. See \cite{FSZ} for an earlier counting result.

    It remains unknown whether any individual value $\zeta(2n+1)$ with $n\geq2$ is irrational, and no odd zeta value is known to be transcendental. See \cite{Zu4} for more details on (iii). Brown and Zudilin constructed infinitely many rational numbers $p/q$ with $0<\left|\zeta(5)-p/q\right|<q^{-0.86}$ in \cite{BrZ}.
\end{remark}

\section{Multiple zeta values}\label{subsec:mzv}
Due to Conjecture \ref{conj:singlezetaalgindep} we do not expect any linear relations among the values $\zeta(k_1)\zeta(k_2)$ if one of the $k_i$ is odd. But it turns out that certain parts of these products satisfy numerous relations among each other. Splitting the product $\zeta(k_1)\zeta(k_2)$ into the following three parts leads us to the definition of the double zeta values $\zeta(k_1,k_2)$
\begin{align}\label{eq:simplestuffle}
\begin{split}
\zeta(k_1)\zeta(k_2) = \sum_{m_1 > 0} \frac{1}{m_1^{k_1}}\sum_{m_2 > 0}\frac{1}{m_2^{k_2}} &= \left(\sum_{m_1 >m_2> 0} + \sum_{m_2 >m_1 > 0} + \sum_{m_1 = m_2 > 0}    \right) \frac{1}{m_1^{k_1} m_2^{k_2}} \qquad (k_1,k_2 \geq 2) \\
&=: \zeta(k_1,k_2) + \zeta(k_2,k_1) + \zeta(k_1+k_2)\,. 
\end{split}
\end{align}
For example, we have the following expressions for the products of Riemann zeta values 
\begin{align*}
\zeta(2)\zeta(5) = \zeta(2,5) + \zeta(5,2) + \zeta(7)\,,\qquad
\zeta(3)\zeta(4) = \zeta(3,4) + \zeta(4,3) + \zeta(7)\,.
\end{align*}
Even though we expect that there are no linear relations among $\zeta(2)\zeta(5)$ and     $\zeta(3)\zeta(4)$, their ``building blocks'' given by $\zeta(7), \zeta(2,5),\zeta(3,4),\zeta(4,3)$ and $ \zeta(5,2)$ satisfy various relations among each other. For example, we will see (Exercise \ref{ex1}) that
\begin{align}\label{eq:exrel}
\zeta(7) = 4 \zeta(3,4) + 3 \zeta(4,3) - 2 \zeta(5,2)\,.    %
\end{align}

Considering the product of more than just two zeta values and using the same idea as in \eqref{eq:simplestuffle} leads us for integers $k_1,\dots,k_r$ to sums of the form 
\begin{align}\label{eq:mzvsum}
\sum_{m_1>\dots>m_r>0} \frac{1}{m_1^{k_1}\cdots m_r^{k_r}}\,.
\end{align}
\begin{proposition}\label{prop:mzvconv}
    For integers $k_1 \geq 2$, $k_2,\dots,k_r\geq 1$ the sum \eqref{eq:mzvsum} converges. 
\end{proposition}
\begin{proof}
    It is enough to show the convergence for $k_1=2$ and $k_2=\dots=k_r=1$ for any $r$, since this gives an estimate for the other cases. Using the well-known inequality $\sum_{n=1}^m \frac{1}{n} \leq 1+ \log(m)$
    we obtain
    \begin{align*}
    \sum_{m_1>m_2>\cdots>m_r>0} \frac{1}{m_1^2 m_2\cdots m_r} = \sum_{m=1}^\infty \frac{1}{m^2}  \sum_{m>m_2>\dots>m_r>0} \frac{1}{m_2\cdots m_r} \leq \sum_{m=1}^\infty \frac{1}{m^2} (1+\log(m))^{r-1}
    \end{align*}
    and since $(1+\log(m))^{r-1} =  o(\sqrt{m})$ as $m\rightarrow \infty$ for any $r$, the above sum converges.
\end{proof}

Notice that the sum \eqref{eq:mzvsum} diverges in the case $k_1=1$ as for each $m_1 > r-1$, we obtain a lower bound of the multiple zeta value:
\begin{align}
\sum_{m_1 > \dots > m_r > 0} \frac{1}{ m_1 m_2^{k_2} \cdots m_r^{k_r}}
    \ge \frac{1}{(r-1)^{k_2} \cdots 1^{k_r}} \sum_{m_1 > r-1} \frac{1}{m_1} 
\end{align}
By the divergence of the harmonic series, the multiple zeta value diverges when $k_1 = 1$.

The multiple sum \eqref{eq:mzvsum} will give the definition of the multiple zeta values, which we will give after introducing the following notation.
\begin{definition} 
    \begin{enumerate}[\textup{(}i\textup{)}] 
        \item For $r\geq 0$ we call a tuple $\kk = (k_1,\dots,k_r) \in \Z_{\geq 1}^r$ of positive integers an \emph{index}. For $r=0$ we write $\kk = \emptyset$ and refer to it as the empty index.
        \item An index $\kk = (k_1,\dots,k_r)$ is called \emph{admissible} if $k_1\geq 2$ or  $\kk = \emptyset$.
        \item For an index $\kk = (k_1,\dots,k_r)$ we call $\wt(\kk)=k_1+\dots+k_r$ its \emph{weight} and $\dep(\kk) = r$ its \emph{depth}. We set $\wt(\emptyset)=\dep(\emptyset)=0$.
    \end{enumerate}
\end{definition}

\begin{definition} \label{def:mzv} For an admissible index $\kk = (k_1,\dots,k_r)$ define the \emph{multiple zeta value} $\zeta(\kk)$ by
    \begin{align*}
    \zeta(\kk) = \zeta(k_1,\dots,k_r) = \sum_{m_1>\dots>m_r>0} \frac{1}{m_1^{k_1}\cdots m_r^{k_r}} \in \R
    \end{align*}
    and $\zeta(\emptyset)=1$. In the case $r=1$ (resp. $r=2$) we refer to these as single (resp. double) zeta values.
\end{definition}
An early systematic study of these values was given by Hoffman in \cite{H1}, where they were called multiple harmonic series.

\begin{remark}
    \begin{enumerate}[\textup{(}i\textup{)}] 
        \item By Proposition \ref{prop:mzvconv}, for every admissible index $\kk$, $\zeta(\kk)$ is a real number. Although the notions of weight and depth may not be well defined for these real numbers (and in fact, we will already see in Proposition \ref{prop:z3z21} below that this is not the case for depth), we still refer to $\zeta(\kk) = \zeta(k_1,\dots,k_r)$ as having weight $\wt(\kk) = k_1 + \dots + k_r$ and depth $\dep(\kk) = r$.
        
        \item For $r = 1$, the multiple zeta values reduce to the values of the Riemann zeta function. One can also define multiple zeta functions $\zeta(s_1,\dots,s_r)$ for complex variables $s_1, \dots, s_r \in \C$ and study their analytic properties, similar to the classical case. See, for example, the thesis of Onozuka \cite{On} for a detailed and accessible survey, or \cite{Zh2}.
    \end{enumerate}
\end{remark}

As an analogue of Euler's result $\zeta(2) = \frac{\pi^2}{6}$ we have the following.
\begin{proposition}
For $n\geq 1$ we have 
\begin{align*}
    \zeta(\underbrace{2,\dots,2}_n) = \zeta(\{2\}^n) = \frac{\pi^{2n}}{(2n+1)!}\,.
\end{align*}
\end{proposition}
\begin{proof}Using the Taylor series and the product expansion of the sine function in \eqref{eq:sinproduct}, we can write the generating series of $\zeta(\{2\}^n)$ as
\begin{align}
1 + \sum_{n\geq 1} (-1)^n \zeta(\{2\}^n)  x^{2n} &= \prod_{n\geq 1} \left(1 - \frac{x^2}{n^2}\right) =\frac{\sin(\pi x)}{\pi x} = \sum_{n=0}^\infty (-1)^n \frac{(\pi x)^{2n}}{(2n+1)!}\,.
\end{align}
\end{proof}

As we have seen before in example \eqref{eq:exrel}, multiple zeta values satisfy various linear relations. The first one appears in weight $3$ and is originally due to Euler. We will see several ways to prove it, and the interested reader can find $32$ ways of doing so in \cite{BoBr}.
\begin{proposition}\label{prop:z3z21} We have $\zeta(3) = \zeta(2,1)$.
\end{proposition}
\begin{proof}
    One of the most elementary proofs known to the author is the following: Consider the sum
    \begin{align*}
    S = \sum_{m,n > 0} \frac{1}{m n(m+n)}= \sum_{m,n > 0} \frac{1}{n^2} \left(\frac{1}{m}  -  \frac{1}{m+n} \right) = \sum_{n=1}^\infty \frac{1}{n^2} \sum_{m=1}^n \frac{1}{m} = \zeta(3) + \zeta(2,1)\,.
    \end{align*}
    This sum can also be evaluated as follows
    \begin{align*}
    S= \sum_{m,n>0}\left(\frac{1}{n}+ \frac{1}{m}\right)\frac{1}{(m+n)^2} = \sum_{m,n > 0}\frac{1}{n(m+n)^2} + \sum_{m,n>0} \frac{1}{m(m+n)^2}  = 2\zeta(2,1)\,
    \end{align*}
    and therefore the relation $\zeta(3)=\zeta(2,1)$ follows.
\end{proof}
Another way to obtain relations among multiple zeta values is to evaluate the product $\zeta(k_1)\zeta(k_2)$ in two different ways. In \eqref{eq:simplestuffle} we saw that for  $k_1,k_2\geq 2$  
\begin{align}\label{eq:stuffle1}
\zeta(k_1)\zeta(k_2)=\zeta(k_1,k_2) + \zeta(k_2,k_1) + \zeta(k_1+k_2)\,,
\end{align}
which is often called the \emph{stuffle product} (also called \emph{harmonic product}). But we also have the following expression for the product, which is called the \emph{shuffle product}. 
\begin{proposition} \label{prop:shuffle1} For $k_1,k_2\geq 2$ we have
    \begin{align}\label{eq:shuffle1}
    \zeta(k_1) \zeta(k_2) = \sum_{j=2}^{k_1+k_2-1} \left( \binom{j-1}{k_1-1} + \binom{j-1}{k_2-1} \right) \zeta(j, k_1+k_2-j)\,,
    \end{align}
    where we use the usual convention $\binom{n}{k}=0$ for $n<k$.
\end{proposition}
\begin{proof}
    This is Exercise \ref{ex1} and can be done by using partial fraction decomposition. We will give another proof and an explanation of the word ``shuffle'' in Chapter \ref{sec:algebraicsetup} using iterated integrals.
\end{proof}
Comparing the right hand sides of \eqref{eq:stuffle1} and \eqref{eq:shuffle1} gives for $k_1,k_2\geq 2$ a linear relation among multiple zeta values, which is an example of a so-called \emph{(finite) double shuffle relation}. We will consider the stuffle/harmonic \& shuffle product and the resulting double shuffle relations in detail for arbitrary depth in Chapters \ref{sec:algebraicsetup} and \ref{sec:families}.

\begin{ex}\label{ex:z2z3}For $k_1=2$, $k_2=3$ equations  \eqref{eq:stuffle1} and \eqref{eq:shuffle1} give
    \begin{align*}
    \zeta(2)\zeta(3) &= \zeta(2,3)+\zeta(3,2)+\zeta(5)\,,\\
    \zeta(2)\zeta(3) &=\zeta(2,3) + 3 \zeta(3,2)+ 6 \zeta(4,1)\,,
    \end{align*}
    from which we deduce the linear relation $\zeta(5) = 2\zeta(3,2)+6\zeta(4,1)$.
\end{ex}

We will denote the $\Q$-vector space spanned by all multiple zeta values by
\begin{align*}
\mz = \langle \zeta(\kk) \mid \kk \text{ admissible} \rangle_\Q \,.
\end{align*}
For a fixed weight $k\geq 0$ we also define the space of weight $k$ multiple zeta values by
\begin{align*}
\mz_k = \langle \zeta(\kk) \mid \kk \text{ admissible}, \wt(\kk)=k \rangle_\Q \,.
\end{align*}
Notice that $\mz = \sum_{k\geq 0} \mz_k$. With the same idea as in \eqref{eq:simplestuffle}, where we showed that $\zeta(k_1)\zeta(k_2)$ is a linear combination of multiple zeta values of weight $k_1+k_2$, we will see in Chapter \ref{sec:algebraicsetup} that this is true for arbitrary products of multiple zeta values and we will show the following (see Corollary \ref{cor:zetashufflemap}).
\begin{proposition}\label{prop:mzissubalgebra}
    The space $\mz$ is a $\Q$-subalgebra of $\R$ and we have for $k_1,k_2\geq 0$ $$\mz_{k_1} \cdot \mz_{k_2} \subset \mz_{k_1+k_2}.$$
\end{proposition}

We saw in some examples that multiple zeta values of depth $r$ can be expressed in terms of multiple zeta values of lower depths, e.g. $\zeta(2,1)=\zeta(3)$. This reduction to lower depth always works if the parity of the weight and depths differ (here in the example the depth is $2$ and the weight is $3$). This so-called parity theorem will be explained in detail in Section \ref{sec:parity}, but we will mention it here since we will use a special case of it later in this overview.

\begin{theorem}[{Parity, \cite{Ts}}] If $r \not \equiv k \mod 2$ then any multiple zeta value of depth $r$ and weight $k$ can be written as products\footnote{Possibly including products with depth $0$, i.e. $\Q$-linear combinations.} of multiple zeta values of depth $<r$.
\end{theorem}

In the smallest depth, this was already known to Euler, who gave the following result.

\begin{proposition}\label{prop:paritydepth2}
For $k_1\geq 2, k_2\geq 1$ with $k=k_1+k_2$ odd we have
\begin{align*}
\zeta(k_1,k_2) =&  (-1)^{k_2} \sum_{\substack{j=2\\j \text{ even}}}^{k-3} \left( \binom{k-j-1}{k_1-1} + \binom{k-j-1}{k_2-1} + (-1)^{k_2}\delta_{j,k_1} \right)  \zeta(j)\zeta(k-j) \\
&+ \frac{1}{2}\left(  (-1)^{k_1} \binom{k_1+k_2}{k_2}  -1\right) \zeta(k)\,.
\end{align*}   
\end{proposition}

All the relations we saw so far: $\zeta(3)=\zeta(2,1)$,  the relation \eqref{eq:exrel}, and the finite double shuffle relations are relations among multiple zeta values of the same weight. Indeed it is expected that there exist no $\Q$-linear relations among multiple zeta values of different weights, which is part of the following conjecture.

\begin{conjecture} \label{conj:graded}The space $\mz$ is graded by weight, i.e. 
    \begin{align*}
    \mz = \bigoplus_{k\geq 0} \mz_k \,.
    \end{align*}
\end{conjecture}
This conjecture is very strong as it implies (Exercise \ref{ex2}) the transcendence of every multiple zeta value of non-zero weight. One of the main interests in the theory of multiple zeta values is to understand all of their $\Q$-linear relations. There are several families of relations which conjecturally give all linear relations among multiple zeta values in a fixed weight. We will describe some of them in Chapter~\ref{sec:families}. One of the most famous families of relations are the extended double shuffle relations, which we will describe in detail later, but which can be roughly explained as follows: In Chapter \ref{sec:algebraicsetup}, we will introduce \emph{stuffle regularized multiple zeta values} $\zeta^\ast(\kk ; T) \in \mz[T]$ and \emph{shuffle regularized multiple zeta values} $\zeta^\shuffle(\kk ; T) \in \mz[T]$, which are defined for arbitrary indices $\kk$. These assign to the non-admissible index $\kk=(1)$ the indeterminate $\zeta^\ast(1 ; T)=\zeta^\shuffle(1 ; T)=T$ and they satisfy for an admissible index $\kk$ that $\zeta^\ast(\kk ; T)=\zeta^\shuffle(\kk ; T)=\zeta(\kk)$. Regularized multiple zeta values satisfy for arbitrary indices analogues of the stuffle and shuffle product formulas. For example, in smallest depth we have as an analogue of \eqref{eq:stuffle1} and \eqref{eq:shuffle1} for all $k_1,k_2\geq 1$:
\begin{align}\label{eq:edsdep1}
\begin{split}
\zeta^\ast(k_1;T)\zeta^\ast(k_2;T) &= \zeta^\ast(k_1,k_2;T) + \zeta^\ast(k_2,k_1;T) + \zeta^\ast(k_1+k_2;T)\,,\\
\zeta^\shuffle(k_1;T) \zeta^\shuffle(k_2;T) &= \sum_{{\color{red}j=1}}^{k_1+k_2-1} \left( \binom{j-1}{k_1-1} + \binom{j-1}{k_2-1} \right) \zeta^\shuffle(j, k_1+k_2-j;T)\,.
\end{split}
\end{align}
For example, for $k_1 = 1$ and $k_2=2$ we get
\begin{align*}
    \zeta(2) T &= \zeta^\ast(1,2;T) + \zeta(2,1) + \zeta(3),\\
    \zeta(2) T &= \zeta^\shuffle(1,2;T) + 2\zeta(2,1)\,,
\end{align*}
from which we get $\zeta^\ast(1,2;T) = \zeta(2) T - \zeta(3)-\zeta(2,1)$ and $\zeta^\shuffle(1,2;T) = \zeta(2) T - 2 \zeta(2,1)$. Since we already proved that $\zeta(2,1) = \zeta(3)$ we see that $\zeta^\ast(1,2;T)=\zeta^\shuffle(1,2;T)$, but for general indices $\kk$ we have $\zeta^\ast(\kk ; T) \not = \zeta^\shuffle(\kk ; T)$. There exists an explicit $\R$-linear map $\rho: \R[T] \rightarrow \R[T]$, such that $\zeta^\shuffle(\kk ; T) = \rho(  \zeta^\ast(\kk ; T))$. In particular we will see that $\rho(1)=1$ and $\rho(T)=T$, which conversely shows that $\zeta(2,1)=\zeta(3)$, without using Proposition \ref{prop:z3z21}. Comparing $\zeta^\ast$ and $\zeta^\shuffle$ in general by this $\rho$ gives a large family of relations, which are called \emph{extended double shuffle relations}. Conjecturally these give all relations among multiple zeta values and we will give a precise version of the following conjecture later.

\begin{conjecture}[\cite{IKZ}]\label{conj:extdsh}All relations among multiple zeta values are a consequence of the extended double shuffle relations.
\end{conjecture}

 In particular, we have a conjecture for the dimension of the spaces $\mz_k$, which was first observed by Zagier based on extensive numerical calculations. To state the conjecture we first introduce the  integers $d_k$ given by the following generating series
\begin{align*}
\sum_{k\geq 0} d_k X^k = \frac{1}{1-X^2-X^3}\,,
\end{align*}
i.e. they are given by $d_0=1$, $d_1=0$, $d_2=1$ and the recursion $d_k = d_{k-2} + d_{k-3}$ for $k \geq 3$.

\begin{conjecture}[Zagier, 1994 {\cite[Section~9, p.~509]{Za2}}]\label{conj:zagier} We have $\dim_\Q \mz_k = d_k$ for all $k\geq 0$.
\end{conjecture}
This conjecture shows that multiple zeta values satisfy a lot of linear relations. For example, in weight $k=14$ there are $2^{12}=4096$ admissible indices (i.e. generators of $\mz_{14}$) and the conjectured dimension is $d_{14}=21$. In the following we give a table for the number of admissible indices, the conjectured number of linearly independent relations and the numbers $d_k$.
\begin{center}
\begin{tabular}{|c|c|c|c|c|c|c|c|c|c|c|c|c|c|c|c|}
    \hline
    weight $k$         & 0 & 1 & 2 & 3 & 4 & 5 & 6 & 7 & 8 & 9 & 10 & 11 & 12 & 13 & 14 \\ \hline
    \# of adm. ind.       & 1 & 0 & 1 & 2 & 4 & 8 & 16 & 32 & 64 & 128 & 256 & 512 & 1024 & 2048 & 4096 \\ \hline
    \# of relations  $\overset{?}{=}$      & 0& 0& 0& 1& 3& 6& 14& 29& 60& 123& 249& 503& 1012& 2032& 4075\\ \hline
    $\left( \dim_\Q\mathcal{Z}_k \overset{?}{=} \right) d_k$ & 1 & 0 & 1 & 1 & 1 & 2 & 2 & 3 & 4 & 5 & 7  & 9  & 12 & 16 & 21 \\ \hline
\end{tabular}
\end{center}

\vspace{0.1cm}

\begin{remark} One easy way to do numerical experiments for multiple zeta values is to use PARI/GP, which was actually the tool used by Zagier to come up with Conjecture \ref{conj:zagier}  (It is available for free at \url{https://pari.math.u-bordeaux.fr/}). 
    There the multiple zeta value $\zeta(k_1,\dots,k_r)$ is implemented as  $\operatorname{zetamult}([k_1,\dots,k_r])$ and the Riemann zeta function $\zeta(s)$ as $\operatorname{zeta}(s)$, which is of course the same as $\operatorname{zetamult}([s])$. Together with the function $\operatorname{lindep}([v_1,\dots,v_l])$ one can search for linear relations among the values $v_1,\dots,v_l$. For example to check if there is a relation between $\zeta(2,1)$ and $\zeta(3)$ one enters 
    \begin{verbatim}
    input: lindep([zetamult([2,1]),zeta(3)]) 
    output: [-1, 1]~
    \end{verbatim}
    \vspace{-\baselineskip}
    The output gives the coefficient of the relation $-1\cdot \zeta(2,1) + 1\cdot \zeta(3) =0$. It is unknown if Euler also used PARI/GP to come up with $\zeta(2)=\frac{\pi^2}{6}$, but he could have done so by using 
    \begin{verbatim}
    input: lindep([zeta(2),Pi^2]) 
    output: [-6, 1]~
    \end{verbatim}
    \vspace{-\baselineskip}
    Of course, these numerical calculations will not give a proof of any relations, since it is just a check up to a certain precision. But in any case an interested student should play around a little bit with these tools and maybe try to find nice relations and patterns, which he/she then could try to prove using the machinery developed in these notes. 
    
    If the coefficients in the output are extremely large, then this is an indication for the fact that there are no $\Q$-linear relations among the values in the input. For example to check if there is a $\Q$-linear relation among $\zeta(3),\pi^3$ and $1$ one enters
    \begin{verbatim}
    input: lindep([zeta(3),Pi^3,1])
    output: [-5229795329281686, 216810578846293, -435977217249266]~
    \end{verbatim}
\vspace{-\baselineskip}
    which indicates that there are (as expected) no relations among these values. The size of the numbers here depends on the current precision used, which can be changed by \verb+\p 50+ to, for example, set the precision to 50 significant digits.
\end{remark}

Conjecture \ref{conj:zagier} is out of reach at the moment and so far there is no weight $k$, for which we can actually prove that $\dim_{\Q} \mz_k > 1$, since for example it is not even known (but expected) that $\zeta(5)$ and $\zeta(2,3)$ are linearly independent and that $\dim_{\Q} \mz_5 =2$. Even though it seems to be impossible to give lower bounds for $\dim_{\Q} \mz_k$ so far, we know that the $d_k$ give upper bounds:

\begin{theorem}[Terasoma (2002) {\cite{Ter}}, Deligne--Goncharov (2005) {\cite{DG}}]\label{thm:terasoma} For all $k\geq 0$ we have $$\dim_{\Q} \mz_k \leq d_k.$$
\end{theorem}
There also is a conjecture on an explicit basis for $\mz$ due to Hoffman.
\begin{conjecture}[Hoffman \cite{H2}, 1997]\label{conj:hoffman} For $k\geq 0$ the multiple zeta values
    \[\{\zeta(k_1,\dots,k_r) \mid r\geq 0,\, k_1+\dots+k_r=k,\, k_1,\dots,k_r \in \{2,3\}\}\]
    form a basis of $\mz_k$. 
\end{conjecture}
Notice that this conjecture would imply  (Exercise \ref{ex2}) Zagier's dimension conjecture (Conjecture \ref{conj:zagier}) since $d_k$ counts exactly the number of indices of weight $k$ with only 2's and 3's.
Multiple zeta values with only 2's and 3's in their index will be called  \emph{Hoffman elements}. The linear independence of the Hoffman elements is unknown so far, but we know that these generate the whole space due to the following deep result of Brown.
\begin{theorem}[Brown \cite{Br2}, 2012]\label{thm:brown} For all $k\geq 0$ we have
    \[\mz_k = \langle \zeta(k_1,\dots,k_r) \mid r\geq 0,\, k_1+\dots+k_r=k,\, k_1,\dots,k_r \in \{2,3\}\rangle_\Q \,.\]    
\end{theorem}
Zagier's explicit evaluation of $\zeta(\{2\}^a,3,\{2\}^b)$ in \cite{Za4} is an important ingredient in the proof of Theorem \ref{thm:brown}.

The only known proofs of Theorem \ref{thm:terasoma} and \ref{thm:brown} use deep concepts from algebraic geometry, particularly the theory of mixed Tate motives.
\begin{remark}
    In his work \cite{Br2}, Brown shows that the motivic versions of Conjectures~\ref{conj:graded}, \ref{conj:zagier} and~\ref{conj:hoffman} hold for so-called motivic multiple zeta values $\mz^{\mathfrak{m}}$, which are conjecturally isomorphic as a $\Q$-algebra to $\mz$. Using the surjective period map $\operatorname{per}: \mz^{\mathfrak{m}} \twoheadrightarrow \mz$ Theorem \ref{thm:brown} is then just a consequence of his more general results. For details on this, we refer to the excellent book \cite{BF}.
\end{remark}

\newpage
\section{Finite and symmetric multiple zeta values}\label{subsec:fmzv}

The finite multiple zeta values are a variant of multiple zeta values, which are not given as elements of the real numbers. Instead, they are elements of the following ring $\mathcal{A}$ (see \cite[(2)]{KZ}) given by
\[  \mathcal{A} = \quotient{\prod_{p \text{ prime}} \fp \,}{\bigoplus_{p \text{ prime}} \fp }\,.\]
This ring was introduced by Zagier and he refers to it as the ``poor man's adeles'', due to its similarity to the usual finite adeles and its simpler definition. Elements in $\ma$ are infinite tuples $a=(a_p)_{p \text{ prime}} \in \ma$ with $a_p \in \fp$, such that for $a=(a_p)_{p},b=(b_p)_{p} \in \ma$ we have $a=b$ if and only if $a_p = b_p$ for $p \gg 0$. In other words, two such elements are the same if $a_p \neq  b_p$ for just finitely many primes $p$.
Due to this we have an embedding $ \iota : \Q \rightarrow \mathcal{A}$ (Exercise \ref{fex14}, see also \cite[\S1, Example~1]{KZ}), since for $\frac{a}{b} \in \Q$ we can get a solution $x_p$ of
\[  b\, x_p - a \equiv 0 \mod p \] 
for all but finitely many $p$ (those, where $p | b$). Choose $x_p$ arbitrary if it does not exist and define
\[   \iota \left(\frac{a}{b} \right) = ( x_2 , x_3 , x_5, x_7 , \dots ) \in \ma  \,.\]
For example, we have
\begin{align*}
    \iota \left(\frac{7}{15} \right) = (1, \ast_3 , \ast_5, 0, 10, 10, 5, 3, 2, 14,17, 35, 36, 32, 13,\dots)\,,
\end{align*}
where the $\ast_p$ can be any element in $\fp$. Notice that $\iota$ is injective and therefore $\mathcal{A}$ becomes a $\Q$-algebra. 
The finite multiple zeta value will be given by a collection of rational numbers depending on prime $p$. For this we first define the truncated version of multiple zeta values, called \emph{multiple harmonic sums}, for $\kk = (k_1,\dots,k_r) \in \Z^r$ and $m\geq 1$ by 
\begin{align}\label{eq:defhm}
 H_m(\kk) =    H_m(k_1,\dots,k_r) = \sum_{m\geq m_1>\dots>m_r>0} \frac{1}{m_1^{k_1}\cdots m_r^{k_r}} \in \Q\,.
\end{align}
The $H_m(1)$ are the classical harmonic numbers. Notice that for a prime $p$, the element $H_{p-1}(k_1,\dots,k_r) \mod p \in \fp$ is well-defined, since all $m_j$ in the denominator are smaller than $p$. Collecting these for all prime $p$ gives our main object for this section (see \cite[\S1, Example~7]{KZ}):

\begin{definition}
    For $\kk=(k_1,\ldots,k_r) \in \Z^r$ the {\bf finite multiple zeta value}  is defined by 	\begin{align*} 
\za(\kk) = \za(k_1,\dots,k_r) &= \left(  H_{p-1}(\kk) \mod p \right)_{p \text{ prime}}\!\!\!\! 
&= \left( \sum_{p > m_1>\cdots>m_r>0} \frac{1}{m_1^{k_1}\cdots m_r^{k_r}}  \mod p \right)_{p \text{ prime}} \!\!\!\!\in \mathcal{A}\,.
\end{align*}
Here we also set $\za(\emptyset)=1$.
\end{definition}
Notice that Kaneko and Zagier use the opposite convention for the order of an index in \cite{KZ}.

\begin{ex}
In the case $\kk=(2,1)$ the finite multiple zeta value is given by 
\begin{align}\begin{split}\label{eq:fmz21}
    \za(2,1) &= \left( \sum_{p > m_1>m_2>0} \frac{1}{m_1^2 m_2}  \mod p \right)_{p}\\
    &= \left( { \sum_{2 > m_1>m_2>0} \frac{1}{m_1^2 m_2} \!\!\!\!\mod 2, \sum_{3 > m_1>m_2>0} \frac{1}{m_1^2 m_2}  \!\!\!\!\mod 3, \dots } \right) \\
    &= (0, 1, 1, 3, 4, 5, 4, 15, 15, 27, 14, 2, 31, 15, 24, 49, 31, 38, 31, 54,\dots).
    \end{split}
\end{align}
For the index $\kk=(2,2)$ one gets
\begin{align*}
    \za(2,2) &= (0, 1, 3, 0, 0, 0, 0, 0, 0, 0, 0, 0, 0, 0, 0,\dots). 
\end{align*}
Here, it seems like at some points all entries become $0$. Since we ignore a finite number of places, this element would therefore be $0$ in $\mathcal{A}$. In fact, we will see that $\za(2,2)=0$ as $\za(k)=0$ (Proposition \ref{prop:fmzvdepth1}) and $\za(2,2) = \frac{1}{2}\left( \za(2)^2-\za(4)\right)$ by \eqref{eq:fmzvstuffle1}. In contrast, for \eqref{eq:fmz21} it seems that the entries do not become $0$. Actually, in addition to $p=2$ we only know two places where this is the case, namely for $p=16843$ and $p=2124679$, which are the only known ``Wolstenholme primes''. However, surprisingly, we have so far not been able to prove that $\za(2,1) \neq 0$ in $\mathcal{A}$. See Conjecture \ref{conj:zknonzero} and Open problem \ref{prob:nonzerofmzv} below.
\end{ex}

Notice that we also allowed negative $k_j$, since all the appearing sums are finite and there are no convergence issues. But as we will see in Proposition \ref{prop:negativeindexfmz} below, one can express the finite multiple zeta values with negative entries always in terms of those with positive entries. Therefore, we define the space of finite multiple zeta values as the following $\Q$-vector space
\begin{align*}
    \mza = \langle \za(k_1,\dots,k_r) \mid r\geq 0,\, k_1,\dots,k_r \geq 1 \rangle_\Q = \langle \za(\kk) \mid \kk \text{ index} \rangle_\Q \,.
\end{align*}
For a fixed weight $k\geq 0$ we also define the space of finite multiple zeta values of  weight $k$  as
\begin{align*}
\mza_k = \langle \za(\kk) \mid \kk \text{ index}, \wt(\kk)=k \rangle_\Q \,.
\end{align*}

Notice that $\Q \subset \mza$ since we also include the empty index $(r=0)$.

\begin{proposition}\label{prop:negativeindexfmz}For any $k_1,\dots,k_r \in \Z$ we have $\za(k_1,\dots,k_r) \in \mza$.
\end{proposition}
\begin{proof}
This is Exercise \ref{fex4}. The main idea is to use the Seki-Bernoulli formula in Lemma \ref{lem:fermatsekibernoulli}. See also \cite[Proposition~5]{KZ}.
\end{proof}

In the following, we will consider finite multiple zeta values in small depths. For this we will recall the following two statements, which will be used several times when doing calculations with finite multiple zeta values.
\begin{lemma}\label{lem:fermatsekibernoulli}
\begin{enumerate}[(i)]
\item {\normalfont (Fermat's little theorem)} For a prime $p$ and $m\in \Z$ we have
\begin{align*}
    m^p \equiv m \mod p\,.
\end{align*}
    \item {\normalfont(Seki-Bernoulli formula / Faulhaber's formula)} For $k, n\geq 1$ we have 
\begin{align*}
    \sum_{m=1}^{n-1} m^k = \frac{1}{k+1} \sum_{j=0}^k \binom{k+1}{j} B_j n^{k+1-j}\,.
\end{align*}
\end{enumerate}
\end{lemma}
\begin{proof}
There are various simple proofs for (i). For example, one can reduce the problem to $m\geq 0$ and prove it by induction on $m$. For the induction step one proves the `Freshman's dream'  $(a+b)^p \equiv a^p + b^p \mod p$ for $a,b\in \Z$, which can be shown by using that $\binom{p}{j}\equiv 0 \mod p$ for $0<j<p$.

For (ii) one can consider the exponential generating series of the left-hand side and gets
\begin{align*}
    \sum_{k=1}^\infty \left( \sum_{m=1}^{n-1} m^k  \right) \frac{X^{k-1}}{(k-1)!} &= \sum_{m=1}^{n-1} m e^{mX} =\frac{\partial}{\partial X} \frac{1-e^{nX}}{1-e^X} = \frac{\partial}{\partial X}  \frac{X}{e^X-1}  \frac{e^{nX}-1}{X} \\
    &= \frac{\partial}{\partial X} \left( \sum_{j=0}^\infty B_j \frac{X^j}{j!} \right) \left( \sum_{m=1}^\infty n^m \frac{X^{m-1}}{m!}\right) = \sum_{k=1}^\infty\left( k \sum_{j=0}^k  \frac{B_j}{j!} \frac{n^{k+1-j}}{(k+1-j)!}\right) X^{k-1} \,.
\end{align*}
Comparing the coefficient of $X^{k-1}$ then yields the result. 
\end{proof}

In particular, if $m$ is coprime to $p$ (e.g. $p>m>0$), then $m^{p-1} \equiv 1 \mod p$. Combining (i) and (ii) in Lemma \ref{lem:fermatsekibernoulli} yields for $n\geq 1$, $k\in \Z$ with $k\not = 0$ and primes $p > \max(|k|,n-1)$
\begin{align}\label{eq:hnkmodp}
   H_{n-1}(k) =  \sum_{m=1}^{n-1} \frac{1}{m^k} &\equiv  \sum_{m=1}^{n-1} m^{p-k-1} \equiv  -\frac{1}{k} \sum_{j=0}^{p-k-1} \binom{p-k}{j} B_j n^{p-k-j} \mod p.
\end{align}

The equation \eqref{eq:hnkmodp} can be used to evaluate finite multiple zeta values in depth one and two.

\begin{proposition}\label[proposition]{prop:fmzvdepth1} For $k\in \Z$ we have
\begin{align}
\za(k) = (H_{p-1}(k) \mod p)_p = \begin{cases} 0,&k\neq 0\\
-1,& k=0\end{cases}\,.
\end{align}
\end{proposition}
\begin{proof}
The $k=0$ case is clear and the other cases follow from \eqref{eq:hnkmodp} by using $p=n$. Alternatively, for primes $p>|k|+1$ one can choose a primitive root $a$ modulo $p$ to get for $k \not = 0$  (without using the Seki-Bernoulli formula)
\begin{align*}
    \sum_{m=1}^{p-1} \frac{1}{m^k} \equiv  \sum_{i=0}^{p-2} \frac{1}{a^{ik}} \equiv \frac{1 -  a^{-k (p-1)}}{1-a^{-k}}  \equiv 0 \mod p\,.
\end{align*}
For $k\geq1$, this is also contained in \cite[\S1, Example~7(1)]{KZ}. The remaining cases are contained in \cite[proof of Proposition~5]{KZ}.
\end{proof}
Proposition \ref{prop:fmzvdepth1} shows that finite multiple zeta values behave quite differently to multiple zeta values. But, as we will see later, there seems to be a conjectural deep connection between these two. In particular, there exists an element  $Z(k) \in \mza$ which serves as the correct \emph{analogue of $\zeta(k)$} instead of the naive choice $\za(k)$ (in a sense we will make clearer later). For $k\geq 2$ this element, introduced in \cite[\S1, Example~5]{KZ}, is defined by
\begin{align}\label{eq:defZ}
    Z(k) = \left( \frac{B_{p-k}}{k} \right)_p \in \ma.
\end{align}
Notice that this definition makes sense, since we can always ignore small primes $p$, i.e. those cases where the Seki-Bernoulli number $B_{p-k}$ is not defined. Notice that, since $B_{\text{odd}}$ vanish, we have $Z(2m)=0$ for all $m\geq 1$.

\begin{proposition}\label[proposition]{prop:fmzvdep2} For $k_1,k_2 \geq 1$ we have 
\begin{align*}
\za(k_1,k_2) =   (-1)^{k_1} \binom{k_1+k_2}{k_1} Z(k_1+k_2)\,.  
\end{align*}
\end{proposition}
\begin{proof}
For large enough primes $p$ we can use \eqref{eq:hnkmodp} to obtain
\begin{align*}
    \sum_{p>m_1 > m_2 > 0} \frac{1}{m_1^{k_1} m_2^{k_2}} &= \sum_{m_1=1}^{p-1}  \frac{1}{m_1^{k_1}}\sum_{m_2=1}^{m_1-1}\frac{1}{m_2^{k_2}} \equiv -\frac{1}{k_2} \sum_{m_1=1}^{p-1}  \sum_{j=0}^{p-k_2-1} \binom{p-k_2}{j} B_j m_1^{p-k_2-j-k_1} \\
    &\equiv \frac{1}{k_2} \binom{p-k_2}{p-k_1-k_2} B_{p-k_1-k_2} \equiv (-1)^{k_1} \binom{k_1+k_2}{k_1} \frac{B_{p-k_1-k_2}}{k_1+k_2} \mod p\,,
\end{align*}
where in the third equation we used $\sum_{m_1=1}^{p-1} m_1^{p-k_2-j-k_1} \equiv 0 \mod p$ except for $j=p-k_1-k_2$.
This is also \cite[(58)]{KZ}, after reversing the order of the index.
\end{proof}

Since the odd Seki-Bernoulli numbers vanish, Proposition \ref{prop:fmzvdep2} implies $\za(k_1,k_2)=0$ if $k_1+k_2$ is even. For odd weight we expect that these are non-zero, i.e. we have the following. 

\begin{conjecture}\label{conj:zknonzero}
For odd $k\geq 3$ we have $Z(k) \not = 0$ in $\ma$.
\end{conjecture}
This is the non-vanishing conjecture discussed in \cite[\S1, Example~5]{KZ}.

This conjecture is related to a classical conjecture of Siegel in algebraic number theory, which asserts that there are infinitely many regular primes, namely those primes $p$ that do not divide the class number of the $p$-th cyclotomic field. By a criterion of Kummer, a prime $p$ is regular if and only if it does not divide the numerator of the Bernoulli numbers $B_k$ for $k=2,4,\dots,p-3$. If this conjecture were true, we would also see that $Z(k) \not = 0$ in $\ma$. Conjecture \ref{conj:zknonzero} suggests that we expect certain elements in $\mza$ to be non-zero. By definition we also include $\za(\emptyset) = 1$, i.e. $\Q = \mza_0 \subset \mza$, and therefore have infinitely many non-zero elements. But the current situation is that we are actually not able to show that any of the $\za(\kk)$ for non-empty $\kk$ are non-vanishing. See also \cite[Introduction, Remark~1]{KZ}. In particular, we have the following:

\begin{problem}\label{prob:nonzerofmzv}
Show that there exists a non-empty index $\kk \in \Z^r_{\geq 1}$ such that  $\za(\kk) \neq 0$.
\end{problem}
See \cite{Se2} for a general discussion on the non-vanishing of finite multiple zeta values. \\

As for multiple zeta values, we can ask for relations among the finite multiple zeta values and the dimension of the space $\mza_k$. For example, as a direct consequence of the definition, we have  for $k_1,\dots,k_r \in \Z$ and any prime $p$
\begin{align}\begin{split}\label{eq:reversalfmzvcalc}
    \sum_{p > m_1 > \dots > m_r > 0} \frac{1}{m_1^{k_1} \cdots m_r^{k_r}} \equiv     &\sum_{p > m_1 > \dots > m_r > 0} \frac{ (-1)^{k_1+\dots+k_r}}{(p-m_1)^{k_1} \cdots (p-m_r)^{k_r}}\\
 \equiv&\sum_{0 < \tilde{m}_1 < \dots < \tilde{m}_r< p} \frac{ (-1)^{k_1+\dots+k_r}}{\tilde{m}_1^{k_1} \cdots \tilde{m}_r^{k_r}} \mod p\,,
 \end{split}
\end{align}

which implies the \emph{reversal formula} (see \cite[\S1, Example~7(3)]{KZ})
\begin{align}\label{eq:reversalformula}
    \za(k_1,\dots,k_r) = (-1)^{k_1+\dots+k_r} \za(k_r,\dots,k_1)\,.
\end{align}

In the case of multiple zeta values one source of relations (which conjecturally give all) is given by the double shuffle relations. We have an analogue of the \emph{stuffle product} \cite[Proposition~4]{KZ}, e.g. in smallest depth we have for any $k_1,k_2\geq 1$
\begin{align}\label{eq:fmzvstuffle1}
    \za(k_1) \za(k_2) = \za(k_1,k_2) + \za(k_2,k_1) + \za(k_1+k_2).
\end{align}
But in contrast to multiple zeta values we do not expect that the shuffle product formula \eqref{eq:shuffle1} holds. Instead, we have a family of linear relations, called the \emph{linear shuffle relations} \cite[Theorem~2]{KZ}. For example, as an analogue of Proposition \ref{prop:fmzvdep2}, we get for $k_1,k_2\geq1$ 
\begin{align}\label{eq:fmzvshuffle1}
 (-1)^{k_1} \za(k_1,k_2)=   \sum_{j=1}^{k_1+k_2-1} \left( \binom{j-1}{k_1-1} + \binom{j-1}{k_2-1} \right) \za(j, k_1+k_2-j)\,.
\end{align}
Due to Proposition \ref{prop:fmzvdepth1} and \ref{prop:fmzvdep2} the equations \eqref{eq:fmzvstuffle1} and \eqref{eq:fmzvshuffle1} do not produce really interesting relations, but for example we have (as a trivial consequence of Proposition \ref{prop:fmzvdep2}) 
\begin{align}\label{eq:za41z32}
2\za(4,1)  + \za(3,2) = 0\,.
\end{align}
But the higher-depth generalization, which we will present in Chapter \ref{sec:fmzv}, actually gives interesting and non-trivial relations. Moreover, we have the following conjecture.
\begin{conjecture}[Kaneko--Zagier {\cite[Conjecture~3]{KZ}}]All relations among finite multiple zeta values are a consequence of the stuffle product and the linear shuffle relations. 
\end{conjecture}

The really surprising aspect of the whole story is that there seems to be a connection to usual multiple zeta values, which we will describe now. For this we define a real counterpart of the finite multiple zeta values introduced by Kaneko and Zagier \cite[Introduction and \S8]{KZ}. This can also be seen as finite multiple zeta values at the infinite prime $p\rightarrow \infty$. 
The following idea is due to Kontsevich who communicated it to Zagier for the depth two case. In the definition of finite multiple zeta values we take a sum over $p > m_1 > \dots > m_r > 0$ and then consider it modulo $p$. Naively, therefore, one could say that we consider the sum over `$0 > -1 \geq m_1 > \dots > m_r \geq 1 > 0$', which we illustrate in the following picture. 
\begin{figure}[h!]
    \centering
\input{TikZ/tikz_infinitep}
\end{figure}

Inspired by this we consider the \emph{order $\succ$ on $\Z\backslash\{0\} \cup \{ \infty = -\infty\}$} defined by
\begin{align*}
    -1 \succ -2 \succ -3 \succ \dots \succ -\infty = \infty \succ \dots \succ 3 \succ 2 \succ 1\,.
\end{align*}
Then we define for an index $\kk =(k_1,\dots,k_r)$ and $m\geq 1$ as an analogue of the multiple harmonic sums the following rational numbers
\begin{align}\label{eq:defsm}
    S_m(\kk) = S_m(k_1,\dots,k_r) = \sum_{\substack{m_1 \succ \dots \succ m_r \\ m \geq |m_1|,\dots,|m_r| > 0}} \frac{1}{m_1^{k_1} \cdots m_r^{k_r}} \in \Q\,.
\end{align}

By direct calculation one then checks that we have the following.
\begin{proposition}\label{prop:smintermsofhm}  For all $m\geq 1$ and $k_1,\dots,k_r \geq 1$ we have 
\begin{align}\label{eq:smashm}
S_m(k_1,\dots,k_r)  = \sum_{j=0}^r (-1)^{k_1+\dots + k_j} H_m(k_j,k_{j-1},\dots,k_1) H_m(k_{j+1},\dots,k_r) \,.
\end{align}
Here we use $H_m(k_j,k_{j-1},\dots,k_1)=1$ (resp. $H_m(k_{j+1},\dots,k_r) =1$) for $j=0$  (resp. $j=r$).
\end{proposition}
\begin{proof}
This is Exercise \ref{fex5} or \ref{ex:hw22}.
\end{proof}

The limit of $H_m(\kk)$ for $m \rightarrow \infty$ exists only if $\kk$ is admissible. A priori it is therefore not clear that the limit of \eqref{eq:smashm} for $m\rightarrow \infty$ exists. Although this is the case (Proposition \ref{prop:symmzvstatements}), we will use \eqref{eq:smashm} as a motivation for the following definition, which is due to Kaneko and Zagier \cite[\S8, (87)]{KZ}.

\begin{definition}
For an index $\kk=(k_1,\ldots,k_r)$ and $\bullet \in \{\shuffle, \ast \}$ the {\bf $\bullet$-symmetric multiple zeta value}  is defined by 
\begin{align*}
      \zs^\bullet(\kk)=  \zs^\bullet(k_1,\dots,k_r) = \sum_{j=0}^r (-1)^{k_1+\dots+k_j} \zeta^\bullet(k_j,k_{j-1},\dots,k_1;T) \zeta^\bullet(k_{j+1},\dots,k_r;T)\,,
\end{align*}
where for $j=0$ (resp. $j=r$) we set $\zeta^\bullet(k_j,k_{j-1},\dots,k_1;T)=1$ (resp. $\zeta^\bullet(k_{j+1},\dots,k_r;T)=1$). 
\end{definition}

In depth one we have for $k\geq 1$
\begin{align*}
\zs^\bullet(k)=\zeta^\bullet(k;T) + (-1)^k \zeta^\bullet(k;T) = \left\{\begin{array}{ll}
2\zeta(k) & ,\,\hbox{$k$ is even} \\ 
0 & ,\,\hbox{$k$ is odd}
\end{array}
\right. \,.
\end{align*}

For the depth two case we have 
\begin{align*}
   \zs^\ast(1,1)= -\zeta(2),\qquad \zs^\shuffle(1,1)= 0,
\end{align*}
 and for $k_1,k_2 \geq 1$ with $k_1+k_2\geq 3$
\begin{align}\label{eq:smzvdepth2}
    \zs^\bullet(k_1,k_2) &= \zeta^\bullet(k_1,k_2;T) + (-1)^{k_1} \zeta^\bullet(k_1;T) \zeta^\bullet(k_2;T) + (-1)^{k_1+k_2} \zeta^\bullet(k_2,k_1;T)\\
    &= (1+(-1)^{k_1}) \zeta(k_1,k_2) + (-1)^{k_1}(1+(-1)^{k_2}) \zeta(k_2,k_1) + (-1)^{k_1} \zeta(k_1+k_2)\,.
\end{align}
Here we used that $\zeta^\ast(k_2,k_1;T)=\zeta^\shuffle(k_2,k_1;T)$ if $k_1+k_2\geq 3$. Since $(1+(-1)^{k_1})$ vanishes for $k_1=1$ (same for $k_2$) there are no non-admissible indices appearing and we therefore did not need any regularization in the last equation. In particular, we see that $\zs^\bullet(k_1,k_2)$ also does not depend on $T$. This is the case in general and the following proposition will be proved in Chapter \ref{sec:fmzv}.
\begin{proposition}\label{prop:symmzvstatements}
\begin{enumerate}[(i)]
    \item We have $ \zs^\bullet(\kk) \in \mz$ for $\bullet \in \{\shuffle, \ast \}$.\\
    (i.e. the $\bullet$-symmetric multiple zeta values are independent of $T$).
    \item For all indices $\kk$ we have 
\begin{align*}
    \lim_{m\rightarrow \infty }S_m(\kk)  = \zs^\ast(\kk)\,.
\end{align*}
    \item For all indices $\kk$ we have 
\begin{align*}
    \zs^\ast(\kk) \equiv  \zs^\shuffle(\kk) \mod \pi^2 \mz\,.
\end{align*}
\end{enumerate} 
\end{proposition}
Statements (i) and (iii) are \cite[Theorem~3]{KZ}, while (ii) is \cite[Proposition~9]{KZ}.

\begin{theorem}[{\cite[Theorem~6.1]{Yas}}]\label{thm:symspanmz} We have $\mz = \langle \zs^\ast(\kk) \mid \kk \text{ index} \rangle_\Q = \langle \zs^\shuffle(\kk) \mid \kk \text{ index} \rangle_\Q$.
\end{theorem}
This theorem is not obvious, since for example $\zs^\ast(k) = 0$ for odd $k$, i.e. it is not obvious why $\zeta(k)$ for odd $k$ can be written as a linear combination of $\zs^\ast(\kk)$. To express these, one needs to consider symmetric multiple zeta values of higher depth. For example, we have
\begin{align*}
\zeta(3) &= \frac{1}{3}\zs^\ast(2,1)\,,\\
    \zeta(5) &= \frac{2}{5}\zs^\ast(4,1)+\frac{2}{15} \zs^\ast(2,2,1)\,,\\
    \zeta(7) &= \frac{10}{7} \zs^\ast(6,1) - \frac{3}{7} \zs^\ast(2,2,1,2) + \frac{5}{21} \zs^\ast(2,2,2,1)\,,\\
    \zeta(9) &= \frac{2}{5} \zs^\ast(8,1) + \frac{2}{15} \zs^\ast(2,2,2,1,2)  + \frac{4}{15} \zs^\ast(2,2,2,3) + \frac{4}{15} \zs^\ast(2,2,2,2,1).
\end{align*}
One can obtain these expressions as follows: By using Proposition \ref{prop:paritydepth2} we have for odd $k\geq 3$
\begin{align*}
    \zeta(k-1,1) = \frac{1}{2}(k-1) \zeta(k) - \sum_{\substack{j=2\\j \text{ even}}}^{k-3} \zeta(j) \zeta(k-j)
\end{align*}
and together with \eqref{eq:smzvdepth2} we get $\zs^\ast(k-1,1) = 2 \zeta(k-1,1) + \zeta(k)$ from which we deduce
\begin{align}\label{eq:oddzetarecursion}
    \zeta(k) &= \frac{1}{k} \zs^\ast(k-1,1) + \frac{1}{k} \sum_{\substack{j=2\\j \text{ even}}}^{k-3}  \zs^\ast(j) \zeta(k-j).
\end{align}
As we will see later the $\zs^\ast$ satisfy the stuffle product formula and therefore \eqref{eq:oddzetarecursion} gives an inductive way of obtaining an explicit expression of $\zeta(k)$ for odd $k\geq 3$ in terms of $\zs^\ast$ by starting with $\zeta(3) =\frac{1}{3}\zs^\ast(2,1)$.

Theorem \ref{thm:symspanmz} shows that relations among multiple zeta values imply relations among $\ast$-symmetric multiple zeta values. For example, using the extended double shuffle relations one can show that 
\begin{align*}
2\zs^\ast(4,1)  + \zs^\ast(3,2) = \zeta(5) - 2 \zeta(2,3) + 4 \zeta(4,1) = 0\,.    
\end{align*}
Comparing this with \eqref{eq:za41z32} we see that the $\zs^\ast$ seem to satisfy the same relation as $\za$. But in general we know that this is not true, since for example $\zs^\ast(2) = 2\zeta(2) = \frac{\pi^2}{3}$ but $\za(2) = 0$. But, if we consider the $\bullet$-symmetric multiple zeta values modulo $\pi^2$, then we have that they are actually independent of $\bullet$ and they conjecturally satisfy the same relations as the finite multiple zeta values. This leads to the following definition. 
\begin{definition}
For $\kk=(k_1,\ldots,k_r)$ the {\bf symmetric multiple zeta value} $\zs(\kk) \in \mathcal{Z}/\pi^2\mathcal{Z}$ is defined by 
\begin{align*}
    \zs(\kk) \equiv  \zs^\ast(\kk) \equiv  \zs^\shuffle(\kk)  \mod \pi^2 \mz.
\end{align*}
\end{definition}
This is the definition in \cite[(88)]{KZ}.
Now we also have $\zs(k)=0$ for all $k\geq 1$ (as for $\za$), and the depth two case is given by the following. 

\begin{proposition} \label{prop:dep2smzv}For $k_1,k_2\geq 1$ we have 
\begin{align*}
    \zs(k_1,k_2) \equiv \begin{cases} 0,& k_1+k_2 \text{ even}\\
    (-1)^{k_1} \binom{k_1+k_2}{k_1}\zeta(k_1+k_2),& k_1+k_2 \text{ odd}
    \end{cases} \mod \pi^2 \mz\,.
\end{align*}
\end{proposition}
\begin{proof}
For even $k_1+k_2\geq 2$ we get by \eqref{eq:smzvdepth2} that
\begin{align*}
      \zs^\ast(k_1,k_2) = (1+(-1)^{k_1}) \zeta^\ast(k_1;T)\zeta^\ast(k_2;T) -\zeta(k_1+k_2)\,,
\end{align*}
which by Proposition \ref{prop:euler} is always a rational multiple of $\pi^{k_1+k_2}$. 
For the odd weight case we use the parity result in Proposition \ref{prop:paritydepth2}, which states that for $k_1\geq 2, k_2\geq 1$ with $k=k_1+k_2$ odd we have
\begin{align*}
\zeta(k_1,k_2) =&  (-1)^{k_2} \sum_{\substack{j=2\\j \text{ even}}}^{k-3} \left( \binom{k-j-1}{k_1-1} + \binom{k-j-1}{k_2-1} + (-1)^{k_2}\delta_{j,k_1} \right)  \zeta(j)\zeta(k-j) \\
&+ \frac{1}{2}\left(  (-1)^{k_1} \binom{k_1+k_2}{k_2}  -1\right) \zeta(k)\,.
\end{align*}
Plugging this into \eqref{eq:smzvdepth2} then yields the result, since the terms in the first sum vanish modulo $\pi^2 \mz$. 
This is also \cite[(89)]{KZ}, after reversing the order of the index.
\end{proof}
Notice that Proposition \ref{prop:dep2smzv} again has strong similarities with the analogous result for finite multiple zeta values (Proposition \ref{prop:fmzvdep2}), where we also had $\za(k_1,k_2)=0$ for even $k_1+k_2$. For odd $k_1+k_2$ we had the same formula except that $Z(k_1+k_2)$ takes the place of $\zeta(k_1+k_2)$, which underlines the comment that $Z(k)$ (defined in \eqref{eq:defZ}) can be seen as the finite analogue of $\zeta(k)$. In general we see that, at least in the example we have presented so far, the symmetric multiple zeta values satisfy the same relations as the finite multiple zeta values. This seems to be the case in general, which is part of the following surprising conjecture. 

\begin{conjecture}[Kaneko--Zagier conjecture {\cite[Main Conjecture]{KZ}}]\label{conj:kzconjecture} The following is an isomorphism of $\Q$-algebras
 \begin{align*}
\varphi_{KZ}: \mza &\longrightarrow \mathcal{Z}/\pi^2\mathcal{Z} \\
 \zeta_{\mathcal{A}}(\kk) & \longmapsto \zeta_{\mathcal{S}}(\kk) \,.
 \end{align*}
\end{conjecture}

In analogy, and according to Conjecture \ref{conj:kzconjecture}, we also have the following dimension conjecture.
\begin{conjecture}[Kaneko--Zagier {\cite[(6)]{KZ}}] For $k\geq 3$ we have $\dim_\Q \mza_k = d_{k-3}$.
\end{conjecture}

Conjecture \ref{conj:kzconjecture} in particular predicts that the finite and symmetric multiple zeta values satisfy the same $\Q$-linear relations. There are various families of linear relations that have been proven for finite as well as for symmetric multiple zeta values. One of these is given by the Hoffman duality. To state it, we first introduce the star versions of finite and symmetric multiple zeta values. 

For $\kk = (k_1,\dots,k_r) \in \Z^r$ and $m\geq 1$ define  
\begin{align}
 H^{\color{red}\star}_m(\kk) =    H^{\color{red}\star}_m(k_1,\dots,k_r) = \sum_{m \geq m_1 {\color{red}\geq} \dots  {\color{red}\geq} m_r>0} \frac{1}{m_1^{k_1}\cdots m_r^{k_r}} \in \Q\,,
\end{align}
and the \emph{finite multiple zeta star-values} 
\begin{align*} 
\za^\star(\kk) = \za^\star(k_1,\dots,k_r) &= \left(  H^\star_{p-1}(\kk) \mod p \right)_{p}\!\!\!\! 
&= \left( \sum_{p > m_1 \geq \cdots \geq m_r>0} \frac{1}{m_1^{k_1}\cdots m_r^{k_r}}  \mod p \right)_{p} \in \mathcal{A}\,.
\end{align*}
Here we also set $\za^\star(\emptyset)=1$. Notice that the star-values of both objects are just a linear combination of the usual values, since for example
\begin{align*}
     \za^\star(k_1,k_2,k_3) =  \za(k_1,k_2,k_3)+\za(k_1+k_2,k_3)+\za(k_1,k_2+k_3)+\za(k_1+k_2+k_3)\,.
\end{align*}
In general we have 
\begin{align*}
    \za^\star(k_1,\dots,k_r) = \sum_{\circ_j = `+` \text{ or } \circ_j=`,`}\za(k_1 \circ_1 \dots \circ_{r-1} k_r)\,.
\end{align*}
Using this description we define the {\bf star-version of the symmetric multiple zeta values} by 
\begin{align*}
    \zs^\star(k_1,\dots,k_r) = \sum_{\circ_j = `+` \text{ or } \circ_j=`,`}\zs(k_1 \circ_1 \dots \circ_{r-1} k_r) \in \mathcal{Z}/\pi^2\mathcal{Z}\,.
\end{align*}
By Conjecture \ref{conj:kzconjecture} the $ \za^\star$ and $\zs^\star$ should therefore also satisfy the same relations. 

Every non-empty index $\kk = (k_1, \dots, k_r)$ can be written as
\[(k_1, \dots, k_r) = (\overbrace{1+\dots+1}^{k_1}, \dots , \overbrace{1+\dots+1}^{k_r}) \,. \]
Define the {\bf Hoffman dual} $\kk^\vee$ by interchanging \,{ \bf $,$}\, and \,{ \bf +}\, in this representation.

For example, the Hoffman dual of $\kk=(3,2)$ is given by
\begin{align*}
\kk^\vee = (3,2)^\vee  = (1+1+1,1+1)^\vee = (1,1,1+1,1) = (1,1,2,1) \,.
\end{align*}
With this, we have the following family of linear relations which is true for finite as well as symmetric multiple zeta values. 

\begin{theorem}\label{thm:fmzvhoffmandual}
For all non-empty indices ${\bf k}$ and $\mathcal{F} \in \{ \mathcal{A},\mathcal{S}\}$ we have
\[ \zeta_\mathcal{F}^\star({\bf k}) = -\zeta_\mathcal{F}^\star({\bf k}^\vee) \,.\]
\end{theorem}
\begin{proof}
For $\mathcal{F}=\mathcal{A}$ this is \cite[Theorem 4.7]{H4} and for $\mathcal{F}=\mathcal{S}$ it can be found in \cite{J}
\end{proof}
We will see that Theorem \ref{thm:fmzvhoffmandual} can also be proven by using the results in \cite{BTT1}.

\section{(Quasi-)Modular forms}\label{subsec:mfandbk}
In this section, we want to give a glimpse of the connection of modular forms and multiple zeta values and state the Broadhurst-Kreimer conjecture, which is a refinement of Zagier's dimension conjecture (Conjecture \ref{conj:zagier}). For this we will give a naive argument why cusp forms give rise to relations among double zeta values, which we will make precise later in Chapter \ref{sec:mdandmzv}. We will not give a complete introduction to modular forms and just state the main structure theorems and definitions. For a complete introduction to the theory of modular forms we refer the reader to \cite{La} and \cite{Za3}.
\subsection{Basics of modular forms} \label{subsubsec:basicofmf}
Let $\Ha  =\big\{  x + i y \in \C \mid x,y \in \R \,, y > 0 \big\}$ denote the complex upper half plane. A holomorphic function $f\in \mathcal{O}(\Ha)$ is called a \emph{modular form of weight $k\in \Z$} (for $\SL_2(\Z)$) if it satisfies the following two conditions:
\begin{enumerate}[(i)]
    \item (Modular transformation): For all $\tau \in \Ha$ and all $\abcd \in \SL_2(\Z)$ we have
    \begin{align}\label{eq:modulartrans}
    f\left( \frac{a\tau +b}{c\tau +d}\right) = (c\tau +d)^k f(\tau).
    \end{align}
    \item (Holomorphy at $\infty$): $f(\tau)$ is bounded as $\tau \rightarrow i \infty$.
\end{enumerate}
Since $\SL_2(\Z)$ is generated by the two matrices $\begin{psmallmatrix}1&1\\0&1 \end{psmallmatrix}$ and $\begin{psmallmatrix}0&-1\\1&0\end{psmallmatrix}$ one can show that it suffices for (i) if $f$ satisfies
$f(\tau +1)=f(\tau)$ and $f(-\frac{1}{\tau}) = \tau^k f(\tau)$ for all $\tau\in \Ha$. Condition (ii) could also be replaced by saying that $f$ has a \emph{Fourier expansion} of the form 
\begin{equation}\label{eq:fourier}
f(\tau) = \sum_{n=0}^\infty a_n q^n \,\qquad (a_n \in \C)\,,\qquad q=q(\tau) = e^{2\pi i \tau}\,.
\end{equation}
The map $q\colon \tau \mapsto \exp(2\pi i \tau)$ is the holomorphic map that sends $\Ha$ to the punctured unit disc. The existence of \eqref{eq:fourier} states that $f$, as a function in $q$, can be analytically continued to $q=0$, and therefore gives rise to a holomorphic function in the whole open unit disc $\{q \in \C \mid |q|<1\}$.  This would therefore imply condition (ii) and explains the notion ``Holomorphy at $\infty$''.\\

By $\mf_k$ we denote the \emph{space of all modular forms of weight~$k$}, and we write $\mf=\sum_{k\geq 0}\mf_k$ for the space of all modular forms.  It is easy to see that $\mf$, equipped with the usual multiplication of holomorphic functions, forms a $\C$-algebra.  The coefficients $a_n\in\C$ in \eqref{eq:fourier} are called the \emph{Fourier coefficients} of~$f$, and a modular form with $a_0=0$ (that is, the sum in \eqref{eq:fourier} starts at $n=1$) is called a \emph{cusp form}.  We write
\[
\ms_k=\bigl\{f\in\mf_k \mid f=\sum_{n=1}^{\infty} a_n q^n\bigr\}
\]
for the \emph{space of all cusp forms of weight~$k$}. The first non-trivial examples of modular forms are the \emph{Eisenstein series}, which for even $k\geq 4$ are defined by
\begin{align}\label{eq:defclassicalgk}
\aG_k(\tau)=\frac12\sum_{\substack{m,n\in\Z\\(m,n)\ne(0,0)}}\frac{1}{(m\tau+n)^k}\,.
\end{align}
For all even $k\geq 4$ we have $\aG_k\in\mf_k$. Notice that \eqref{eq:defclassicalgk} vanishes for odd~$k$. In fact, for every odd~$k$ we have $\mf_k=0$. The Fourier expansion of Eisenstein series, which we will compute in Section~\ref{subsec:fourier} for a more general object, is
\begin{equation}\label{eq:gkfourier}
\aG_k(\tau)=\zeta(k)+\frac{(-2\pi i)^k}{(k-1)!}\sum_{n=1}^{\infty}\sigma_{k-1}(n)\,q^n,
\end{equation}
where $\sigma_{k-1}(n)=\sum_{d\mid n} d^{k-1}$ is the divisor sum. Formula \eqref{eq:gkfourier} makes sense for any $k\geq 2$, and therefore we define $\aG_k(\tau)$ for all $k\geq 2$ by \eqref{eq:gkfourier}. For $k=2$ and for odd~$k$, the series defined by \eqref{eq:gkfourier} is \emph{not} a modular form of weight~$k$. As we will see below, there are no non-trivial modular forms in these weights. Eisenstein series are the building blocks for all modular forms, and we summarize the main properties of modular forms in the following Theorem.
\begin{theorem} \label{thm:mf}
    \begin{enumerate}[\textup{(}i\textup{)}]            
        \item For even $k\geq 4$ we have $\mf_k = \C \cdot \aG_k \oplus \ms_k$. 
        \item For odd or negative $k$, as well as for $k=2$, we have $\mf_k = 0$, while $\mf_0 = \C$.
        \item For all $k_1,k_2 \geq 0$ we have $\mf_{k_1} \cdot \mf_{k_2} \subset \mf_{k_1+k_2}$.
        \item The Eisenstein series $\aG_4$ and $\aG_6$ are algebraically independent (over $\C$).
        \item We have 
        \[\mf = \bigoplus_{k=0}^\infty \mf_k = \C[\aG_4,\aG_6] \,,\]
        i.e. $\mf$ is a graded $\C$-algebra, which is isomorphic to the polynomial ring in two variables.
        \item For an even positive integer $k$ the dimension of $\mf_k$ is given by
        \begin{align} \label{eq:dimmk}
\dim_\C \mf_k = \begin{cases} \lfloor \frac{k}{12} \rfloor + 1&\,, \quad k \not\equiv 2 \mod 12 \\
 \lfloor \frac{k}{12} \rfloor &\,, \quad k \equiv 2 \mod 12
 \end{cases}\,.
 \end{align}
 The generating series for the dimension of modular forms in weight $k$ is given by
 \begin{align*}
     \mathsf{M}(X) &:= \sum_{k\geq 0} \dim_\C \mf_k X^k = \frac{1}{(1-X^4)(1-X^6)}\\
     &= 1+X^4+X^6+X^8+X^{10}+2 X^{12}+X^{14}+2 X^{16}+2 X^{18}+2 X^{20}+2 X^{22}+3 X^{24}+\dots .
 \end{align*}
        \item The space $\mf_k$ is generated by $\aG_k$ and products of two Eisenstein series, i.e. for even $k\geq 4$ 
        \begin{align*}
        \mf_k = \C \cdot  \aG_k  +  \langle \aG_{k_1} \aG_{k_2} \mid k_1,k_2 \geq 4 \text{ even}, k_1+k_2=k \rangle_\C \,.
        \end{align*}
    \end{enumerate}
\end{theorem}
\begin{proof}
    The statements (i), (ii) \& (iii) follow almost immediately from the definition. Statements (iv) and (v) need some complex analysis and can be found in any standard book of modular forms (e.g. \cite{Ba6}, \cite{La}, \cite{Za3}). The dimension formula in (vi) is a consequence of (iv) and (v). By these statements, we see that every modular form can be uniquely written as a polynomial in $\aG_4$ and $\aG_6$. In particular, the monomials $\aG_4^a \aG_6^b$ with $a,b\geq 0$ and $4a+6b=k$ form a basis of $\mf_k$. In order to count these one can consider the geometric series and obtains \begin{align*}
        \sum_{k\geq 0} \#\{\text{ pairs } a,b \text{ with } 4a+6b=k, a,b\geq 0 \} X^k = \sum_{a\geq 0} X^{4a}  \sum_{b\geq 0} X^{6b} = \frac{1}{(1-X^4)(1-X^6)}. 
    \end{align*}
    Statement (vii) follows from the work of Rankin and can be found in \cite{KoZ}.
\end{proof}

The first non-trivial cusp form is the \emph{discriminant function $\Delta$}  (a.k.a Ramanujan Delta function)
\begin{equation}\label{eq:deltadef}
\Delta(\tau) = q \prod_{n=1}^\infty \left( 1- q^n\right)^{24} = q - 24 q^2 +252 q^3 -1472 q^4  + \dots \,,
\end{equation}
which is a cusp form $\Delta \in \ms_{12}$ of weight $12$.

Since $\Delta$ has no zero in $\Ha$ and a zero of order one in $q=0$, one can show that every cusp form of weight $k$ can be written as a product of a modular form of weight $k-12$ and $\Delta$. Together with Theorem \ref{thm:mf} this gives the following well-known theorem.
\begin{theorem}\label{thm:mfdim}
    \begin{enumerate}[\textup{(}i\textup{)}] 
        \item For $k\geq 0$ the map $\mf_k \rightarrow \ms_{k+12}$ given by $f\mapsto \Delta \cdot f$ is an isomorphism of $\C$-vector spaces.
        \item The generating series for the dimension of cusp forms of weight $k$ is given by
        \begin{align*}
        \Sx(X) = \sum_{k\geq 0} \dim_\C \ms_k X^k &=X^{12} \sum_{k\geq 0} \dim_\C \mf_k X^k  =X^{12}\mathsf{M}(X) =  \frac{X^{12}}{(1-X^4)(1-X^6)}\,.
        \end{align*}
    \end{enumerate}
\end{theorem}

\subsection{The Broadhurst-Kreimer Conjecture}
We now want to state a refinement of Conjecture \ref{conj:zagier} given by Broadhurst and Kreimer, which indicates a connection of cusp forms and the dimension of the depth graded spaces of multiple zeta values.  The depth gives a filtration on the space $\mz$ and we write 
\[ \fild_r(\mz_k) = \langle \zeta(\kk) \mid \kk \text{ admissible}, \wt(\kk)=k\,, \dep(\kk) \leq r  \rangle_\Q \]
for its depth $r$ part and denote the associated graded part by $\grd_r(\mz_k)$. In other words, elements in $\grd_r(\mz_k)$ are multiple zeta values of weight $k$ and depth $r$ modulo multiple zeta values of lower depth. For example, the class of $\zeta(2,1)$ in $\grd_2(\mz_3)$ is zero, since $\zeta(2,1)=\zeta(3)$.

\begin{conjecture}[Broadhurst-Kreimer, 1997 \cite{BroK}] \label{conj:brokrei}
    The generating series of the dimensions of the weight- and depth-graded parts of multiple zeta values is given by
    \begin{align*}
    \sum_{k,r \geq 0} \dim_\Q \left(\grd_{r} \mz_k\right) X^k Y^r  &
    = 
    \frac{1 +\ev(X) Y}{ 1- \odd(X) Y + \Sx(X) Y^2 - \Sx(X) Y^4},
    \end{align*}
    where %
    \begin{align*}
    \ev(X) = \frac{X^2}{1-X^2},  \quad
    \odd(X) = \frac{X^3}{1-X^2} , \quad
    \Sx(X) = \sum_{k\geq 0} \dim_\C \ms_k X^k = \frac{X^{12}}{(1-X^4)(1-X^6)}.
    \end{align*}
\end{conjecture} 
The depth-graded motivic version of this picture is studied by Brown in \cite{Br3}.
This conjecture reduces to Zagier's dimension conjecture (Exercise \ref{ex2}) by setting $Y=1$. Observe that
\begin{align}\label{eq:bkeq}
\frac{1 +\ev(X) Y}{ 1- \odd(X) Y + \Sx(X) Y^2 - \Sx(X) Y^4} & = 1 + \left( \ev(X) + \odd(X) \right) Y + \left( \big(\ev(X) + \odd(X)\big) \odd(X) - \Sx(X) \right) Y^2  + \cdots \,,
\end{align}
which indicates that cusp forms give rise to relations in depth $2$. Before we make this more precise, we first want to give a naive reason why cusp forms give rise to linear relations among double zeta values.  Let $f=\sum_{n=1}^\infty a_n q^n \in \ms_k$ be a cusp form. By Theorem \ref{thm:mf} (vi) we know that $f$ can be written for some $\alpha, \beta_{a,b}  \in \C$ (usually we are interested in the case when $\alpha, \beta_{a,b}  \in \Q$) as
\begin{align}
f = \alpha \aG_k + \sum_{\substack{a,b \geq 4 \text{ even}\\a+b=k}} \beta_{a,b}  \aG_a \aG_b\,. \label{eq:fasprod}
\end{align}
Considering the constant term in the Fourier expansion of both sides then yields the relation
\begin{align*}
0 = \alpha \zeta(k)+ \sum_{\substack{a,b \geq 4 \text{ even}\\a+b=k}} \beta_{a,b}  \zeta(a) \zeta(b)\,.
\end{align*}
The products on the right-hand side can now be evaluated using either \eqref{eq:stuffle1} or \eqref{eq:shuffle1}  to obtain a linear relation among double zeta values. This approach is not really interesting, since the representation of a cusp form \eqref{eq:fasprod} is not unique and also the choice of expanding the product is arbitrary. Therefore, we cannot really relate a cusp form to a single relation among double zeta values. But there is a surprising 1:1 correspondence between certain relations and cusp forms, which we will explain now. \\

For even weight $k$, the Broadhurst-Kreimer conjecture predicts by \eqref{eq:bkeq} that $\dim_\Q \left(\grd_{2} \mz_k\right)$ is given by the coefficient of $X^k$ in $\odd(X) \odd(X) - \Sx(X)$. The coefficient of $\odd(X) \odd(X) $ counts the number of indices $(k_1,k_2)$ with $k_1,k_2\geq 3$ odd and $k_1+k_2=k$ for which we will write ``$(\text{odd},\text{odd})$'' in the following. If this were the only contribution to $\dim_\Q \left(\grd_{2} \mz_k\right)$ then a naive guess would be that the $\zeta(\text{odd},\text{odd})$ give a basis of $ \grd_{2} \mz_k$. Indeed, we will see in Chapter \ref{sec:mdandmzv}  that the $\zeta(\text{odd},\text{odd})$ span $ \grd_{2} \mz_k$. But the factor $\Sx(X)$ in the Broadhurst-Kreimer conjecture indicates that there are relations in weight $k$ between these values whenever there exist cusp forms of weight $k$.  The first relation between $\zeta(k_1,k_2)$, where both $k_1$ and $k_2$ are odd, appears in weight $k_1+k_2 = 12$ and is given by
\begin{align}\label{eq:firstoddrel}
-10394 \zeta(3,9)+ 47650\zeta(5,7)+ 41431\zeta(7,5) -720\zeta(9,3) -10394\zeta(11,1)=0\,.
\end{align}
As shown in \cite{GKZ} (see Chapter \ref{sec:mdandmzv}) we have for all even $k\geq 4$ the relation
\begin{align}\label{eq:oddsumform}
    \sum_{\substack{k_1\geq 3, k_2\geq 1 \text{ odd}\\k_1+k_2=k}} \zeta(k_1,k_2) = \frac{1}{4} \zeta(k)\,.
\end{align}
Combining the relations \eqref{eq:firstoddrel}  and \eqref{eq:oddsumform} gives 
    \begin{align*}
168 \zeta(5,7)+150 \zeta(7,5)+28 \zeta(9,3) = \frac{5197}{691} \zeta(12)\,.
\end{align*}
From this we get the following relation among $\zeta(\text{odd},\text{odd})$ in $ \grd_{2} \mz_{12}$
\begin{align*}
168 \zeta(5,7)+150 \zeta(7,5)+28 \zeta(9,3) \equiv 0 \mod \zeta(12)\,.
\end{align*}
In \cite{GKZ} Gangl-Kaneko-Zagier give an explicit construction of such relations for a given cusp form $f\in \ms_k$. For this they consider its period polynomial $p_f(X,Y)\in \C[X,Y]$ and show that one can obtain a relation among double zeta values explicitly from  the coefficients of this polynomial (see Chapter \ref{sec:mdandmzv}). For example the above relation can be obtained by taking for $f$ a certain multiple of the cusp form $\Delta$. A consequence of their results is the following.
\begin{theorem}[Gangl-Kaneko-Zagier, 2006 {\cite[Theorem~3]{GKZ}}]\label{thm:gkzperiodrel} For even $k\geq 4$ the number of (independent) $\Q$-linear relations among $\zeta(2a+1,2b-1)$ with $a,b\geq 1$ and $k =2(a+b)$ is at least $\dim_\C S_{k}$. 
\end{theorem}
Baumard and Schneps relate the Gangl--Kaneko--Zagier relations by duality to period-polynomial relations in the double shuffle Lie algebra \cite{BS}.
Conjecturally, the number of these relations is exactly $\dim_\C S_{k}$.
Another explanation for the occurrence of modular forms here was given by
Charlton and Keilthy in \cite{ChK}. They showed that, modulo products, the
period-polynomial relations among double zeta values are encoded by the
dihedral and differential relations for block-graded motivic multiple zeta
values in block degree $2$ (see \cite[Propositions~3.5 and~3.7]{ChK}).
Due to a recent result of Tasaka \cite{Tas}, also a somehow converse statement is known: Given a relation among $\zeta(\text{odd},\text{odd})$ (which follows from a certain set of relations) of weight $k$, we can construct explicitly a cusp form of weight $k$. His result uses double Eisenstein series which we will discuss later in this chapter and in more detail in Chapter \ref{sec:mes}.
Further connections between period polynomials and double and triple zeta values are studied by Ma and Tasaka in \cite{MaT}.

\subsection{Modular forms and finite multiple zeta values}\label{subsec:modularfmzv}
Similar to Theorem \ref{thm:gkzperiodrel}, Kaneko and Zagier observed (see \cite[\S3(g)]{KZ} and \cite[p.~184]{K2}) that there also seem to be relations among finite multiple zeta values corresponding to cusp forms. However, instead of having relations in depth $2$, they observed that (due to the shift in depth), one seems to obtain relations in depth $3$ and $4$. The relations in depth $4$ are among finite multiple zeta values of the form $\zeta_{\mathcal{A}}(1,2a,1,2b)$. In weight $12$ and $16$, the first such relations are given by
\begin{align*}
0 = 2\,\zeta_{\mathcal{A}}(1, 2, 1, 8) -18 \,\zeta_{\mathcal{A}}(1, 4, 1, 6) -9\, \zeta_{\mathcal{A}}(1, 6, 1, 4) - 16 \,\zeta_{\mathcal{A}}(1, 8, 1, 2)
\end{align*}
and 
\begin{align*}
0 =&\,\,\, 327 \,\zeta_{\mathcal{A}}(1, 2, 1, 12) + 111\, \zeta_{\mathcal{A}}(1, 4, 1, 10) + 53 \,\zeta_{\mathcal{A}}(1, 6, 1, 8) \\
&+ 59\, \zeta_{\mathcal{A}}(1, 8, 1, 6) + 120\, \zeta_{\mathcal{A}}(1, 10, 1, 4) + 366 \,\zeta_{\mathcal{A}}(1, 12, 1, 2).
\end{align*}
The relation in weight $16$ is also given in \cite[\S3(g)]{KZ}, where it is noted that it can be proved by the linear shuffle relations. Risan also conducted several numerical experiments in his master's thesis \cite{R}. Part of these observations, together with those made by Kaneko and Zagier in \cite[\S3(g), pp.~25--26]{KZ}, is now proved in joint work of the author and Risan.

As above, we write $\fild_r(\mza_k)$ for the subspace of $\mza_k$ spanned
by finite multiple zeta values of depth at most $r$.

\begin{theorem}[Bachmann--Risan, 2026
{\cite{BaR}}]\label{thm:fmzvmodularforms}
For every even $k\geq4$, we have
\begin{align}\label{eq:fmzvmodularspan}
\fild_4(\mza_k)
&=\sum_{\substack{a,b\geq1\\2(a+b+1)=k}}
\Q\za(1,2a,1,2b)
=\sum_{\substack{a,b\geq1\\2(a+b+1)=k}}
\Q\za(2a,1,2b,1)\,.
\end{align}
The dimension satisfies
\begin{align}\label{eq:fmzvmodularupperbound}
\dim_\Q\fild_4(\mza_k)
\leq \frac{k}{2}-2-\dim_\C\ms_k\,.
\end{align}
In particular, there are at least $\dim_\C\ms_k$ linearly independent
relations among the values in either sum in~\eqref{eq:fmzvmodularspan}.
Moreover, every even period polynomial gives an explicit, possibly trivial,
relation of this kind.
\end{theorem}

The remaining observations can be stated as follows.
\begin{conjecture}[Kaneko--Zagier, Risan]\label[conjecture]{conj:fmzvmodularforms}
For every even $k\geq4$, we have the following:
\begin{enumerate}[(i)]
\item Equality holds in~\eqref{eq:fmzvmodularupperbound}.
\item The two spaces in~\eqref{eq:fmzvmodularspan} are also equal to
\begin{align*}
\sum \Q\za(1,1,2a,2b)
&=\sum \Q\za(1,2a,2b,1)\\
&=\sum \Q\za(2a,1,1,2b)
=\sum \Q\za(2a,2b,1,1)\,,
\end{align*}
where the sums run over $a,b\geq1$ with $2(a+b+1)=k$.
\item We have
\begin{align*}
\fild_4(\mza_k)
&= \sum \Q \za(1, 2a, 2b+1) = \sum \Q \za(1, 2b+1, 2a) \\
&= \sum \Q \za(2a, 2b+1, 1) = \sum \Q \za(2b+1, 2a, 1)\,,
\end{align*}
where the sums run over $a,b\geq1$ with $2(a+b+1)=k$. The same statement
holds when we restrict to $a\geq2$ and $b\geq0$ instead.
\end{enumerate}
\end{conjecture}

If Conjecture~\ref{conj:fmzvmodularforms}~(i) holds, then the relations in
Theorem~\ref{thm:fmzvmodularforms} are exactly $\dim_\C\ms_k$ in number.
A direct construction of these finite relations will be given in \Cref{sec:fmzvmodularforms}, after we have introduced period polynomials. In \Cref{sec:formalfinitemzv} we will then study the map from the formal double zeta space and see that its formal version is an isomorphism, while the injectivity after the realization in actual finite multiple zeta values is still open.

\subsection{Quasimodular forms and $\sltwo$-algebras}\label{sec:qmfsl2algebras}
One can also define the Eisenstein series $\aG_2(\tau)$ by \eqref{eq:gkfourier}, i.e. 
\begin{align*}
\aG_2(\tau) = \zeta(2) + (2\pi i)^2 \sum_{n=1}^\infty \sigma_1(n) q^n.
\end{align*}
This is not a modular form anymore since one can check that for $\abcd \in \SLZ$ we have
\begin{equation*}
\aG_2\left( \mabcd \right) = (c\tau + d)^2 \aG_2(\tau) - \pi i c (c\tau +d)\,.
\end{equation*}
We will now introduce another normalization for Eisenstein series in order to consider modular forms with rational coefficients. Using Proposition \ref{prop:euler} we define for even $k\geq 2$ the \emph{normalized Eisenstein series} by 
\begin{equation}\label{eq:ek}
E_k(\tau) = \frac{1}{\zeta(k)} \aG_k(\tau) = 1 - \frac{2k}{B_k} \sum_{n=1}^\infty \sigma_{k-1}(n) q^n \,.
\end{equation}
In particular, we have 
\begin{align}\label{eq:e2e4e6}
\begin{split}
E_2(\tau) &= 1 - 24\sum_{n\geq 1} \sigma_1(n) q^n = 1 - 24 q  - 72 q^2 - 96 q^3 - 168 q^4 + \dots,\\
E_4(\tau) &= 1  + 240\sum_{n\geq 1} \sigma_3(n) q^n =1+240 q+2160 q^2+6720 q^3+17520 q^4+ \dots,\\
E_6(\tau) &= 1  -504\sum_{n\geq 1} \sigma_5(n) q^n = 1-504 q-16632 q^2-122976 q^3-532728 q^4 + \dots.
\end{split}
\end{align}
The space of modular forms with rational coefficients is then given by $\mf^\Q = \Q[E_4,E_6]$. Modular forms are holomorphic functions and therefore we can differentiate them with respect to $\tau$. It is convenient to consider the following notation for a modular form $f = \sum_{n=0}^\infty a_n q^n$:
\[ f^\prime := \qdq f = \frac{1}{2\pi i} \frac{d}{d\tau} f = \sum_{n=1}^\infty n a_n q^n \,.\]
Here the factor $2\pi i$ has been included in order to preserve the rationality properties of the Fourier coefficients. The derivative of a modular form is, in general, not a modular form anymore. By the definition of a modular form one can check that the derivative of a modular form $f \in \mf_k$ satisfies 
\begin{align}\label{eq:mfderivtransform}
    f^\prime\left( \mabcd \right) = (c\tau + d)^{k+2} f^\prime(\tau) + \frac{k}{2\pi i} c (c\tau +d)^{k+1} f(\tau) 
\end{align}
for all $\abcd \in \SLZ$. This, together with the transformation of the Eisenstein series $E_2$ leads to the definition of the \emph{Serre derivative}, given for $f \in \mf_k$ by
\[\partial_k f := f^\prime - \frac{k}{12} E_2 f \,.\]
\begin{proposition} \label{prop:derive2e4e6}
\begin{enumerate}[(i)]
    \item For a modular form $f \in \mf_k$ we have $\partial_k f \in \mf_{k+2}$.
    \item The ring $\qmf=\Q[E_2,E_4,E_6]$ is closed under differentiation and we have
\[ \qdq E_2 = \frac{E_2^2-E_4}{12}\,,\quad \qdq E_4=\frac{E_2 E_4 -E_6}{3}\,,\quad \qdq E_6 = \frac{E_2 E_6 - E_4^2}{2}\,. \]
\end{enumerate}
\end{proposition}
\begin{proof}
Statement (i) follows from a straightforward calculation together with the failure of modularity for $E_2$ and the transformation behavior of the derivative of a modular form \eqref{eq:mfderivtransform}.
For (ii) we can use (i) to get $\partial_4 E_4 = E_4^\prime - \frac{1}{3} E_2 E_4 \in M_6$ and $\partial_6 E_6 = E_6^\prime - \frac{1}{2} E_2 E_6 \in M_8$. Since both spaces are one-dimensional with basis $E_6$ and $E_4^2$ respectively we get the second and third equation after comparing the first Fourier coefficients. Using again the modularity formula of $E_2$ one can also show that $E_2^\prime - \frac{1}{12} E_2^2 \in M_4$. Therefore, this is also a multiple of $E_4$, which turns out to be $-\frac{1}{12}$ by comparing the Fourier coefficients.
\end{proof}

We will be interested in more derivations on $\qmf$ besides $\qdq$. For this, we first recall a few notations from algebra. For a field $K$ and $K$-algebra $A$ a \emph{derivation} is a $K$-linear map $d: A \rightarrow A$ satisfying the Leibniz rule $d(ab) = d(a) b + a d(b)$ for all $a,b \in A$. The set of all derivations on $A$ is denoted by $\operatorname{Der}(A)$.  A \emph{Lie algebra} is a $K$-vector space $V$ with a $K$-bilinear map $[\cdot,\cdot]: V \times V \rightarrow V$, the \emph{Lie bracket}, which satisfies $[x,x]=0$ for all $x\in V$ and the Jacobi identity $[x,[y,z]]+[z,[x,y]]+[y,[z,x]]=0$ for all $x,y,z \in V$. One can check by direct calculation that $\operatorname{Der}(A)$ is a Lie algebra with Lie bracket $[d_1,d_2] = d_1 \circ d_2 - d_2 \circ d_1$. Notice that by the Leibniz rule a derivation is already uniquely determined by its image on the algebra generators. In particular, for $\qmf=\Q[E_2,E_4,E_6]$ a derivation $d: \qmf \rightarrow \qmf$ is uniquely determined by $d(E_2), d(E_4)$ and $d(E_6)$, since we can then, for example, evaluate $$d(E_2^2 E_4) = d(E_2^2) E_4 + E_2^2 d(E_4) = 2 E_2 d(E_2) E_4 + E_2^2 d(E_4)\,.$$
Moreover, one can show that $E_2, E_4, E_6$ are algebraically independent over $\C$. Therefore, one can define a derivation $d$ by choosing arbitrary images for $d(E_k)$ for $k=2,4,6$. We define the derivations $W, \delta: \qmf \rightarrow \qmf$ by 
\begin{align*}
    W(E_2) &= 2 E_2, \quad W(E_4) = 4 E_4, \quad W(E_6) = 6 E_6, \\
    \delta(E_2) &= 12, \quad \delta(E_4) = \delta(E_6) =0\,.
\end{align*}
Notice that for $f \in \qmf_k$ we have $W(f) = k f$, which is called the \emph{weight operator}. The derivation $\delta$ can be seen as the derivative with respect to $\frac{E_2}{12}$. We will now see that the three derivations $\qdq,W, \delta$ make $\qmf$ into a so-called $\mathfrak{sl}_2$-algebra. Notice that the algebra of modular forms is given by $M=\ker \delta$.\\
 
The \emph{Lie algebra $\sltwo$} consists of the $2 \times 2$ complex matrices with trace zero. It is three-dimensional, and one possible basis is given by

\[
X = \begin{pmatrix} 0 & 1 \\ 0 & 0 \end{pmatrix}, \quad
H = \begin{pmatrix} 1 & 0 \\ 0 & -1 \end{pmatrix}, \quad
Y = \begin{pmatrix} 0 & 0 \\ -1 & 0 \end{pmatrix}.
\]

These fulfill the commutator relations
\begin{align}\label{eq:commutatorrel}
[H, X] = 2X, \quad [H, Y] = -2Y, \quad [Y, X] = H.
\end{align}
\begin{definition} A triple of operators $(X,H,Y)$ on an algebra $A$ is called an \emph{$\mathfrak{sl}_2$-triple} if they satisfy the commutator relations \eqref{eq:commutatorrel}. If furthermore the $X,H,Y$ are derivations on $A$ then $A$ is called an \emph{$\mathfrak{sl}_2$-algebra}.
\end{definition}

\begin{proposition}$(\qdq,W,\delta)$ is an $\mathfrak{sl}_2$-triple and thus
    $\qmf$ is an $\mathfrak{sl}_2$-algebra.
\end{proposition}
\begin{proof}
    By the above explanation, one only needs to check that $\eqref{eq:commutatorrel}$ is satisfied by acting on $E_2, E_4$ and $E_6$, which can be done by direct calculation using the explicit formulas for the action of $\qdq$ given in Proposition \ref{prop:derive2e4e6}.
\end{proof}

\section{Multiple Eisenstein series}\label{subsec:mesoverview}
In this section, we introduce multiple Eisenstein series and mention some conjectures and recent results concerning them. Multiple Eisenstein series were first introduced in \cite{GKZ} (in the depth $2$ case) and later studied further in \cite{Ba1} and \cite{Ba5}. Some recent developments can be found in \cite{BaBu}, \cite{BK2}, \cite{BI}, \cite{BIM}, \cite{BKM}, \cite{BT} and \cite{Bu1}. We will return to these developments in Chapters~\ref{sec:mes} and~\ref{sec:formalspaces}.

\begin{definition}
For $k_1, k_2,\dots,k_r \geq 2$ and $\tau \in \Ha$ the \emph{multiple Eisenstein series} are defined by\footnote{In the case $k_1=2$ the sum is not absolutely convergent and we need to use Eisenstein summation. We will make this precise in Chapter \ref{sec:mes}.}
\begin{align}
\mathbb{G}_{k_1,\dots,k_r}(\tau) := \sum_{\substack{\lambda_1 \succ \dots \succ \lambda_r \succ 0\\ \lambda_i \in \Z \tau + \Z}} \frac{1}{\lambda_1^{k_1} \dots \lambda_r^{k_r}}   \,,
\end{align}
where the order $\succ$ on the lattice $\Z \tau + \Z$ is defined as the lexicographical order given by
\begin{align*}
    m_1 \tau + n_1 \succ m_2 \tau + n_2 \quad :\Longleftrightarrow \quad  m_1 > m_2  \,\,\text{ or }\,\, m_1 = m_2 \wedge n_1 > n_2. 
\end{align*}
\end{definition}

Since $\mathbb{G}_{k_1,\dots,k_r}(\tau + 1) = \mathbb{G}_{k_1,\dots,k_r}(\tau)$ the multiple Eisenstein series possess a Fourier expansion, i.e. an expansion in $q=e^{2\pi i \tau}$, which was calculated in \cite{GKZ} for the $r=2$ case and for arbitrary depth by the author (\cite{Ba8}). We will explain the construction in detail in Chapter \ref{sec:mes} and just mention the final result here. For this we define for $k_1,\dots, k_r \geq 1$ the $q$-series
\begin{align}
\g(k_1,\dots,k_r)= \sum_{\substack{m_1 > \dots > m_r > 0\\ n_1, \dots , n_r > 0}} \frac{n_1^{k_1-1}}{(k_1-1)!} \dots \frac{n_r^{k_r-1}}{(k_r-1)!}  q^{m_1 n_1 + \dots + m_r n_r } \in \Q\llbracket q \rrbracket\,.
\end{align}
These $q$-series were studied in detail in \cite{Ba8}, \cite{BK1} and they are $q$-analogues of multiple zeta values, which we will discuss in general in the next section. For now we can just view them as the building blocks of the Fourier expansion of multiple Eisenstein series. In this context the $q$-series $\g$ always appear together with a power of $-2\pi i$ and therefore we set for $k_1,\dots,k_r \geq 1$
\begin{align*}
    \hat{\g}(k_1,\dots,k_r) := (-2\pi i)^{k_1+\dots + k_r} \g(k_1,\dots,k_r) \in \Q[\pi i]\llbracket  q \rrbracket\,.
\end{align*}
Notice that with that notation we can write the classical Eisenstein series as $$\mathbb{G}_k(\tau) = \zeta(k) + \hat{\g}(k)$$ and for multiple Eisenstein series we get the following generalization:
\begin{theorem}[$r=1,2$ \cite{GKZ}, $r\geq 1$ \cite{Ba8}]\label{thm:mesfourier}
 For $k_1,\dots,k_r \geq 2$ there exist explicit $\alpha^{k_1,\dots,k_r}_{l_1,\dots,l_r,j} \in \Z$, such that for $q=e^{2\pi i \tau}$ we have
\begin{align*}
    \mathbb{G}_{k_1,\dots,k_r}(\tau) = \zeta(k_1,\dots,k_r) +\!\!\!\!\!\!\!\!\!\!\sum_{\substack{0 < j < r\\l_1+\dots+l_r = k_1+\dots+k_r\\l_1\geq 2,l_2,\dots,l_r\geq 1}} \!\!\!\!\!\!\!\!\!\! \alpha^{k_1,\dots,k_r}_{l_1,\dots,l_r,j}\, \zeta(l_1,\dots,l_j)  \hat{\g}(l_{j+1},\dots,l_r) +  \hat{\g}(k_1,\dots,k_r)\,.
\end{align*}
In particular, $\mathbb{G}_{k_1,\dots,k_r}(\tau) = \zeta(k_1,\dots,k_r)+ \sum_{n> 0} a_{k_1,\dots,k_r}(n) q^n$ for some $a_{k_1,\dots,k_r}(n) \in \mz[\pi i]$.
\end{theorem}

In the case $r=2$ we get that for $k_1,k_2 \geq 2$ the Fourier expansion of the double Eisenstein series is given by (see \cite[Theorem 6]{GKZ} or Proposition \ref{prop:fourier23})
\begin{align*}
\aG_{k_1,k_2}(\tau) = \zeta(k_1,k_2)+ \!\!\!\!\!\!\!\!\!\sum_{\substack{l_1+l_2 = k_1+k_2 \\ l_1,l_2 \geq 2}}\!\!\! \left( (-1)^{k_2} \binom{l_1-1}{k_2-1} \!+\! (-1)^{l_1-k_1} \binom{l_1-1}{k_1-1} \! + \!\delta_{l_1,k_2}\right)\! \zeta(l_1) \ag(l_2) + \ag(k_1,k_2) \,.
\end{align*}

\begin{ex} Writing $P=-2\pi i$ we have 
\begin{align*}
\aG_6(\tau) &= \zeta(6) + \ag(6) \\
&= \zeta(6) + \frac{1}{120} P^6 q + \frac{11}{40} P^6 q^2 + \frac{61}{30} P^6 q^3 + \dots, \\[4pt]
    \aG_{4,2}(\tau) &= \zeta(4,2) + 2 \zeta(2)\, \ag(4) + 2 \zeta(3) \,\ag(3) + 4 \zeta(4) \,\ag(2) + \ag(4,2) \\
    &= \zeta(4,2) + \Big( \frac{1}{3} P^4 \zeta(2) + P^3 \zeta(3) + 4 P^2 \zeta(4) \Big) q +  \Big( 3 P^4 \zeta(2) + 5 P^3 \zeta(3) + 12 P^2 \zeta(4)\Big) q^2 + \dots,\\[4pt]
    \aG_{3,3}(\tau) &=  \zeta(3,3) + \zeta(3)\, \ag(3)  - 6 \zeta(4) \,\ag(2) +  \ag(3,3) \\
    &=\zeta(3,3)+ \Big( \frac{1}{2}P^3 \zeta(3) - 6 P^2 \zeta(4)  \Big) q +  \Big( \frac{5}{2} P^3 \zeta(3) - 18 P^2 \zeta(4)\Big) q^2 + \dots .
\end{align*}
Using the finite double shuffle relations (Exercise \ref{ex1}) one can show that $\zeta(6) + 3 \zeta(4,2) - 6 \zeta(3,3)=0$. Therefore, the constant term of $\aG_6 + 3 \aG_{4,2} - 6 \aG_{3,3}$ vanishes. But as it turns out, also the other coefficients vanish and we indeed have $\aG_6 + 3 \aG_{4,2} - 6 \aG_{3,3}=0$. This we can see for the first few coefficients as we have
\begin{align}
 \aG_6 + 3 \aG_{4,2} - 6 \aG_{3,3} &= \Big( \frac{P^6}{120} + P^4 \zeta(2) + 48 P^2 \zeta(4) \Big) q +   \Big( \frac{11}{40} P^6 + 9 P^4 \zeta(2) + 144 P^2 \zeta(4) \Big) q^2 + \dots .
\end{align}
And by Euler's result (Proposition \ref{prop:euler}) we have $\zeta(2m) = -\frac{B_{2m}}{2 (2m)!} P^{2m} $ from which we see that these coefficients all vanish. That this is indeed the case for all coefficients of $q^n$ will be a consequence of the relations we will prove among the $q$-series $\g$ later.
\end{ex}
Define the space of multiple Eisenstein series of weight $k\geq 0$ by
\begin{align*}
\mes_k = \langle \aG_{\kk} \mid \kk \in \Z_{\geq 2}^r, r\geq 0, \wt(\kk)=k \rangle_\Q \,
\end{align*}
and set $\mes = \sum_{k\geq 0} \mes_k$. As for multiple zeta values, we make the convention that $\aG_{\emptyset}=1$ is the only multiple Eisenstein series of weight $0$. We also set $\mes^\C = \C \otimes \mes$ and $\mes^\C_k = \C \otimes \mes_k$. Notice that the product of two multiple Eisenstein series can, by a combinatorial argument similar to the one used for multiple zeta values, be written again as a linear combination of multiple Eisenstein series. In particular, we have, for example,
\begin{align*}
    \aG_{k_1} \cdot \aG_{k_2} = \aG_{k_1,k_2} + \aG_{k_2,k_1} + \aG_{k_1+k_2}
\end{align*}
This structure leads to the following statement.

\begin{proposition}\label[proposition]{prop:mesissubalgebra}
The space $\mes$ is a $\Q$-algebra, and we have for $k_1,k_2\geq 0$,
$$ \mes_{k_1} \cdot \mes_{k_2} \subset \mes_{k_1+k_2}. $$
Moreover, the spaces of quasi-modular forms and modular forms are subalgebras of $\mes^\C$. The algebras $\Q[\aG_2,\aG_4,\aG_6]$ and $\Q[\aG_4,\aG_6]$ are subalgebras of $\mes$.
\end{proposition}

Assuming Conjecture \ref{conj:graded}, we expect the space of multiple zeta values to be graded by weight. Since multiple zeta values occur as the constant terms of multiple Eisenstein series, the same grading should extend to $\mes$. In contrast to the case of multiple zeta values, this was recently proved by the author and Kanno.
\begin{theorem}[{\cite[Main Theorem~B(i)]{BKa}}]\label{thm:mesgraded}
The space $\mes$ is graded by weight, i.e.,
\[
    \mes \;=\; \bigoplus_{k \ge 0} \mes_k.
\]
\end{theorem}
In \cite{BKa} it is also shown that the complex span of all multiple Eisenstein series is graded by weight and that it is isomorphic to $\C\otimes_\Q\mes$. We will explain the proof in \Cref{sec:messl2}.

The following statement was conjectured in \cite{BIM} and proved by the author and Kanno in \cite{BKa}.
\begin{theorem}[{\cite[Main Theorem~A]{BKa}}]\label{thm:messl2algebra}
\begin{enumerate}[(i)]
\item The maps $\GD,\GW,\Gdelta$ given on the generators of $\mes$ by
\begin{align*}
    \GD:\aG_{k_1,\dots,k_r} &\longmapsto (2\pi i) \frac{d}{d\tau} \aG_{k_1,\dots,k_r},\\
    \GW: \aG_{k_1,\dots,k_r} &\longmapsto (k_1+\dots+k_r) \aG_{k_1,\dots,k_r},
    \\
    \Gdelta: \aG_{k_1,\dots,k_r} &\longmapsto \begin{cases}
        -\frac{1}{2}\aG_{k_2,\dots,k_r},& k_1=2,\\
        0,& k_1>2,
    \end{cases}
\end{align*}
give well-defined $\Q$-linear maps $\mes \rightarrow \mes$.
\item The maps $\GD,\GW,\Gdelta$ are derivations on $\mes$.
\item $(\GD,\GW,\Gdelta)$ forms an $\mathfrak{sl}_2$-triple, i.e. we have the commutator relations
\begin{align}\label{eq:mescommutatorrel}
[\GW, \GD] = 2\GD, \quad [\GW, \Gdelta] = -2\Gdelta, \quad [\Gdelta, \GD] = \GW
\end{align}
and thus $\mes$ is an $\mathfrak{sl}_2$-algebra.
\end{enumerate}
\end{theorem}

The derivative part of this theorem was proved before in
\cite[Main Theorem~D(ii)]{BKM}. To state the result, extend the map $\aG$ linearly by
\begin{align*}
\aG(z_{k_1}\cdots z_{k_r})=\aG_{k_1,\dots,k_r}
\end{align*}
for $k_1,\dots,k_r\geq2$. The Drop1 operator from
Definition~\ref{def:dropone}, introduced in \cite{HMSW}, gives the map
\begin{align}\label{eq:mesderthetaoverview}
\theta(w)=\dropone\bigl(\ds(w,z_2)\bigr).
\end{align}
Theorem~\ref{thm:mesderformula} states that
\begin{align}\label{eq:mesderformulaoverview}
\GD\aG(w)=\aG\bigl(\theta(w)\bigr).
\end{align}
In particular, $\GD$ is a well-defined derivation on $\mes$ and $\GD\mes\subset\mes$. The well-definedness of $\GW$ is equivalent to \cref{thm:mesgraded}, and the well-definedness of $\Gdelta$ is proved in \cite{BKa} by a nice limit argument, which we will explain in \Cref{sec:messl2}. Also notice that an analogue of the map $\Gdelta$ (i.e., replacing $\aG$ by $\zeta$ in the definition) is not well-defined for multiple zeta values, since $\zeta(2)^2$ is a multiple of $\zeta(4)$.

\begin{theorem}[Theorem \ref{thm:turanthm}]
For $k \geq 2$ the operator $\GD = (2\pi i) \frac{d}{d\tau}$ acts on $\aG_k$ by
\begin{align*}
\GD \aG_{k} = (2k-1) \aG_{k+2} - \sum_{j=2}^k (k+j-1) \aG_{k+2-j,j}  - \aG_{2,k}\,. 
\end{align*}
\end{theorem}

There is also a completely explicit formula in depth two which does not involve the
Drop1 operator.
\begin{theorem}\label{thm:mesderdepthtwo}
For $k,l\geq2$, put $N=k+l+2$. Then
\begin{align}\label{eq:mesderdepthtwo}
\GD\aG_{k,l}={}&(2k-1)\aG_{k+2,l}+2k\aG_{k+1,l+1}
+(2l-1)\aG_{k,l+2}-2l\aG_{l+1,k+1}\notag\\
&+2l\left(\sum_{a=l+2}^{k+1}\aG_{a,N-a}
-\sum_{a=k+2}^{l+1}\aG_{a,N-a}\right)\notag\\
&-\sum_{\substack{a+b=k+2\\a,b\geq2}}(k+b-1)\aG_{a,b,l}
-\sum_{\substack{a+b=l+2\\a,b\geq2}}(l+b-1)\aG_{k,a,b}\notag\\
&-k\sum_{\substack{a+b=l+1\\a,b\geq2}}\aG_{k+1,a,b}
+2l\sum_{\substack{a+b=k+1\\a,b\geq2}}\aG_{l+1,a,b}\notag\\
&-2l\sum_{\substack{a+b=l+2\\a,b\geq2}}\aG_{a,k,b}
-\aG_{2,k,l}-\aG_{k,2,l}.
\end{align}
\end{theorem}
We will prove this formula in Section~\ref{sec:mesderiv} by using the ingredients of \cite{BKM}.

\begin{ex}
For example, Theorem~\ref{thm:mesderdepthtwo} gives
\begin{align*}
\GD\aG_{2,2}&=3\aG_{4,2}+3\aG_{2,4}-12\aG_{2,2,2},\\
\GD\aG_{3,2}&=-\aG_{3,4}+10\aG_{4,3}+5\aG_{5,2}
-10\aG_{2,3,2}-4\aG_{3,2,2},\\
\GD\aG_{2,3}&=5\aG_{2,5}+4\aG_{3,4}-9\aG_{4,3}
-16\aG_{2,2,3}-4\aG_{2,3,2}-8\aG_{3,2,2},\\
\GD\aG_{3,3}&=5\aG_{3,5}+5\aG_{5,3}-12\aG_{2,3,3}
-10\aG_{3,2,3}-10\aG_{3,3,2}+3\aG_{4,2,2}.
\end{align*}
The stuffle product
\begin{align*}
\aG_2\aG_{2,2}=3\aG_{2,2,2}+\aG_{2,4}+\aG_{4,2}
\end{align*}
shows that the first formula can also be written as
$\GD\aG_{2,2}=3\aG_2\aG_{2,2}-21\aG_{2,2,2}$. The middle two examples
also show that the order of the two indices matters.

The third identity gives the following check for the commutator relation $[\Gdelta,\GD] \aG_{2,3} = \GW \aG_{2,3}$ on these representatives:
    \begin{align*}
        \GD \aG_{3} &= 5\aG_{5}-4\aG_{3, 2}-6\aG_{2, 3},\\
        \big(\Gdelta \,\circ\, \GD\big)\aG_{2,3} &= -\frac{5}{2}\,\aG_{5}+2\,\aG_{3, 2}+8\,\aG_{2, 3},\\
        \big(\GD \,\circ\, \Gdelta\big) \aG_{2,3} &= -\frac{1}{2} \GD \aG_{3} =  -\frac{5}{2}\,\aG_{5}+2\,\aG_{3, 2}+3\,\aG_{2, 3},\\
        [\Gdelta,\GD] \aG_{2, 3} &=   \big(\Gdelta \,\circ\, \GD - \GD \,\circ\,\Gdelta \big) \aG_{2, 3} = 5 \aG_{2, 3} = \GW \aG_{2,3}.
    \end{align*}
Expressions of a multiple Eisenstein series in terms of the generators $\aG_{\kk}$ are not unique. Thus the general formula can give representatives different from those obtained by a weight-by-weight reduction.
\end{ex}
Notice that the algebra of quasimodular forms $\Q[\aG_2,\aG_4,\aG_6]$ carries the usual $\mathfrak{sl}_2$-structure described in Section~\ref{sec:qmfsl2algebras}, and the three operators in Theorem~\ref{thm:messl2algebra} restrict to the usual ones on this subalgebra. In particular, the restriction of $\Gdelta$ is $-\frac{1}{2}\frac{\partial}{\partial\aG_2}$. Surprisingly, this is true on all of $\mes$: by \cite[Main Theorem~A(iii)]{BKa}, the algebra $\mes$ is a polynomial ring in $\aG_2$ over the subalgebra spanned by $1$ and the $\aG_{k_1,\dots,k_r}$ with $k_1\geq3$, and $\Gdelta=-\frac{1}{2}\frac{\partial}{\partial\aG_2}$ (see \cref{thm:mespolynomialG2}).

Finally we want to mention a dimension conjecture for the spaces $\mes_k$. Recall that for multiple zeta values we have Zagier's dimension conjecture (Conjecture \ref{conj:zagier}), which states that
\begin{align}\label{eq:zagierconjagain}
       \sum_{k\geq 0} \dim_\Q \mz_k X^k  \overset{?}{=} \frac{1}{1-X^2-X^3}.
\end{align}
To propose an analogous conjecture for multiple Eisenstein series, we first give a different interpretation of the above conjecture. For this consider as before the following series
\begin{align}
    \mathsf{E}(X) = \frac{1}{1-X^2} = 1 + X^2 + X^4 + \dots, \quad \mathsf{O}(X) = \frac{X^3}{1-X^2} = X^3 + X^5 + \dots.
\end{align}
Then one can verify directly that \eqref{eq:zagierconjagain} can be rewritten as
\begin{align}\label{eq:zagierconjagainagain}
       \sum_{k\geq 0} \dim_\Q \mz_k X^k  \overset{?}{=} \frac{1}{1-X^2-X^3} = {\color{Blue} \mathsf{E}(X)} \cdot {\color{DarkOrchid} \frac{1}{1 -  \mathsf{O}(X)}}.
\end{align}
From the perspective of Hilbert--Poincar\'e series for graded algebras, this suggests 
\begin{align}
\mz \overset{?}{\cong} {\color{Blue}\Q[f_2]} \otimes {\color{DarkOrchid}\Q\langle f_3,f_5,\dots \rangle},
\end{align}
with the shuffle product on the right factor $\Q\langle f_3,f_5,\dots \rangle$ (see \cite{Bu1,BF} for details). The left factor $\Q[f_2]$ corresponds to all even single zeta values, which, by Proposition \ref{prop:euler}, satisfy $\zeta(2m) \in \Q[\zeta(2)]$.

In the case of multiple Eisenstein series, we know that not every single Eisenstein series of even weight is a rational multiple of a power of $\aG_2$, though all belong to the ring of quasimodular forms $\Q[\aG_2, \aG_4, \aG_6]$. Recall (Theorem \ref{thm:mfdim}) that the Hilbert--Poincar\'e series of modular forms, cusp forms and quasimodular forms are given by
\begin{align*}
\mathsf{M}(X) &= \frac{1}{(1-X^4)(1-X^6)},\\
\mathsf{S}(X) &= X^{12}\mathsf{M}(X) = \frac{X^{12}}{(1-X^4)(1-X^6)},\\
\widetilde{\mathsf{M}}(X) &= \mathsf{D}(X) \mathsf{M}(X) = \frac{1}{(1-X^2)(1-X^4)(1-X^6)}, 
\end{align*}
where $\mathsf{D}(X)= \frac{1}{1-X^2}$. One possible way to propose an analogue of Conjecture \eqref{eq:zagierconjagainagain} for multiple Eisenstein series would be therefore to replace $\mathsf{E}(X)$ by $\widetilde{\mathsf{M}}(X)$. Moreover, Theorem~\ref{thm:mesderformula} shows that the entire algebra $\mes$ is closed under the derivative $\GD$, which increases the weight by $2$. A natural guess is therefore to include the derivatives of generators in the right-hand factor of \eqref{eq:zagierconjagainagain}, replacing $\mathsf{O}(X)$ by $\mathsf{D}(X) \mathsf{O}(X)$. This seems almost correct to predict $\dim_\Q \mes_k$, except that there appear to be certain not yet understood relations among multiple Eisenstein series arising in pairs from cusp forms and their derivatives. 

\begin{conjecture}\label[conjecture]{conj:mesdim}
We have
\begin{align*}
    \sum_{k\geq 0} \dim_\Q \mes_k X^k &= \widetilde{\mathsf{M}}(X) \cdot  \frac{1}{1 - \mathsf{D}(X) \mathsf{O}(X) + 2 \mathsf{D}(X) \mathsf{S}(X)},\\
&= \mathsf{M}(X) \cdot  \frac{1}{1 - X^2 - \mathsf{O}(X) + 2 \mathsf{S}(X)},\\
    &= \frac{1}{1 - X^2 - X^3 - X^4 - X^5 + X^8 + X^9 + X^{10} + X^{11} + X^{12}}.
\end{align*}
\end{conjecture}

This conjecture is inspired by analogous conjectures by Okounkov \cite{Ok} and the author together with K\"uhn in \cite{BK2} for associated weight-graded spaces of certain $q$-analogues of multiple zeta values. Combined with the expectation that multiple Eisenstein series satisfy the same relations as these $q$-analogues modulo lower weight (\cite{Ba3}), we obtain Conjecture \ref{conj:mesdim}. This expectation is now a theorem of the author and Kanno \cite[Main Theorem~B(ii)]{BKa}, which identifies $\mes$ with the associated weight-graded algebra of the space of these $q$-analogues with all entries at least $2$ (see \cref{thm:mesgrqmzv}). Therefore Conjecture~\ref{conj:mesdim} is equivalent to the dimension conjecture of Okounkov and of the author and K\"uhn for this associated graded algebra. We now give a table of the expected dimensions, relations and generators of $\mes$. Notice that the number of indices of weight $k\geq 2$ with all entries $\geq 2$ and arbitrary depths $r \geq 1$ is given by the $(k-1)$-th Fibonacci number\footnote{Defined by $F_1=F_2=1$ and $F_{k}=F_{k-1}+F_{k-2}$ for $k\geq 3$.} $F_{k-1}$, since its generating series is given by
\begin{align*}
    \sum_{k\geq 2}\left(\text{\# of indices}\right) X^k = \sum_{r\geq 1} \Big( \sum_{k\geq 2} X^k \Big)^r = \sum_{r\geq 1}\Big( \frac{X^2}{1-X}\Big)^r = \frac{\frac{X^2}{1-X}}{1-\frac{X^2}{1-X}} = \frac{X^2}{1-X-X^2} = \sum_{k \geq 2} F_{k-1} X^k.
\end{align*}
We obtain the following table.
\begin{center}
\begin{tabular}{|c|c|c|c|c|c|c|c|c|c|c|c|c|c|c|c|}
    \hline
    weight $k$         & 0 & 1 & 2 & 3 & 4 & 5 & 6 & 7 & 8 & 9 & 10 & 11 & 12 & 13 & 14 \\ \hline
    \# of indices       & 1& 0& 1& 1& 2& 3& 5& 8& 13& 21& 34& 55& 89& 144 &233 \\ \hline
    \# of relations  $\overset{?}{=}$      & 0& 0& 0& 0& 0& 0& 1& 1& 4& 6& 13& 23& 42& 74& 129\\ \hline
    $ \dim_\Q\mes_k \overset{?}{=}$ & 1& 0& 1& 1& 2& 3& 4& 7& 9& 15& 21& 32& 47& 70& 104 \\ \hline
\end{tabular}
\end{center}

We will now give some examples of  relations among multiple Eisenstein series in low weight.

\begin{ex} The two relations among multiple Eisenstein series in weight $6$ and $7$ are
\begin{align*}
    \aG_{6} &= 6 \aG_{3,3} - 3 \aG_{4,2},\\
    \aG_{7} &= 4\aG_{3,4} + 3 \aG_{4,3} - 2 \aG_{5,2}.
\end{align*}
In weight $8$ we have the four relations
\begin{align*}
    \aG_{8} &= 12 \aG_{4,4}\\
     &= 10 \aG_{3,5} - 15 \aG_{4,4}+ 8 \aG_{5,3} - 5  \aG_{6,2} \\
     &= 28\aG_{5,3} + 30\aG_{4,2,2}- 20\aG_{3,3,2}\\
   &= -5 \aG_{2,6} + 48 \aG_{5,3} - 5 \aG_{6,2} + 15 \aG_{4,2,2} - 15 \aG_{2,4,2} + 30 \aG_{3,2,3} +30 \aG_{2,3,3}.
\end{align*}
These relations can be proven using the relations among the $q$-series $\g$ explained in Chapter \ref{sec:mes}. The first two relations, as well as the second relation in weight $8$ after using $\aG_8=12\aG_{4,4}$, also follow from the weighted sum formula in \cite[Main Theorem~A]{Ba12}. Both relations in weight $6$ and $7$ are also the first examples of a family of relations proved by the author, Kanno and Maesaka (see \cref{thm:mesRrelations}). Notice that applying $\delta$ to the last relation in weight $8$ gives the relation in weight $6$.
\end{ex}

The space $\mes$ was previously defined using multiple Eisenstein series $\aG_{k_1,\dots,k_r}$ where all indices satisfy $k_1,\dots,k_r \geq 2$. Since multiple zeta values are defined for indices satisfying $k_1\geq 2, k_2,\dots,k_r\geq 1$, a natural question arises: can the notion of multiple Eisenstein series be extended to include cases where some indices are equal to 1?

The answer is yes. There are indeed several ways to define regularized multiple Eisenstein series, which we will explain in Chapter \ref{sec:mes}. For these regularized multiple Eisenstein series, there also exist dimension conjectures and $\mathfrak{sl}_2$-algebra conjectures, analogous to the statements for the space $\mes$ above. We will explain the origin of these conjectures in Chapter \ref{sec:formalspaces}, when we discuss formal multiple Eisenstein series, for which the $\mathfrak{sl}_2$-algebra structure has been established in \cite{BIM}.

Let us finish this section by mentioning some variants of multiple Eisenstein series which we will not discuss in these notes. Double Eisenstein series of level $2$ were studied in \cite{KT}, and their level $N$ version in \cite{YZ}. A shuffle regularization for multiple Eisenstein series of arbitrary level was constructed in \cite{Ka}. Motivated by symmetric multiple zeta values, Hara, Sakugawa and Tasaka introduced symmetric multiple Eisenstein series in \cite{HaST}. They satisfy a linear shuffle relation by \cite[Proposition~3.9]{HaST}. After the usual normalization, in even weight $k\geq6$ the space spanned by the symmetric double Eisenstein series is given by the space of modular forms of weight $k$ together with the derivative of the Eisenstein series of weight $k-2$, while in odd weight $k\geq3$ its dimension is $\lfloor k/3\rfloor$. They also showed that every cusp form with rational Fourier coefficients is a rational linear combination of normalized symmetric triple Eisenstein series (see \cite[Theorems~1.1 and~1.2]{HaST}).

The Schur multiple zeta values in Section~\ref{subsec:schurmzvoverview} occur as constant terms of Schur multiple Eisenstein series, which are constructed in \cite[Section~4.3]{Yu}. The closely related Schur Eisenstein series and Schur MacMahon series are studied in \cite{BY2}. Connections between MacMahon's generalized sums-of-divisors and (odd) multiple Eisenstein series can also be found in \cite{Ba11}. Another, different family of multiple series of Eisenstein type is studied in \cite{BTs}.

\section{\texorpdfstring{$q$}{q}-analogues of multiple zeta values}\label{subsec:qmzv}
A $q$-analogue of a theorem, identity or expression is a generalization involving a new parameter $q$ that returns the original theorem, identity or expression in the limit as $q \rightarrow 1$\footnote{Here and in the following we mean by $q \rightarrow 1$ the limit of $q$ to $1$ on the real axis with $|q|<1$.}. The easiest example is the $q$-analogue of a natural number $m$ given by the $q$-integer
\begin{align} \label{eq:qinteger}
[m]_q = \frac{1-q^m}{1-q} = 1 + q + \dots + q^{m-1}\,, \qquad \lim\limits_{q\rightarrow 1} [m]_q = m\,.
\end{align}
There are various different models of $q$-analogues of multiple zeta values in the literature. We will consider a few of them and start with the most common model which was first independently studied by Bradley \cite{Bra} and Zhao \cite{Zh1}. For an admissible index $\kk = (k_1,\dots,k_r)$ these are defined by
\begin{align}
\zeta_q^{\rm{BZ}}(\kk)= \zeta_q^{\rm{BZ}}(k_1,\dots,k_r) = \sum_{m_1 > \cdots > m_r > 0} \frac{q^{(k_1-1)m_1} \cdots q^{(k_r-1)m_r}}{[m_1]_q^{k_1} \cdots [m_r]_q^{k_r} } \,.
\end{align}
By \eqref{eq:qinteger}, together with a justification that one can interchange summation and taking the limit (see proof of Proposition \ref{prop:gisqmzv}), it is easy to see that we have 

\[ \lim\limits_{q\rightarrow 1} \zeta_q^{\rm{BZ}}(\kk)= \zeta(\kk)\,. \] In Chapter \ref{sec:families} we will see that these $q$-series satisfy a lot of relations which are satisfied by multiple zeta values. In particular, this is the unique model of $q$-analogues (in the sense we will define later) which satisfies the duality relation, e.g.  $\zeta_q^{\rm{BZ}}(2,1)=\zeta_q^{\rm{BZ}}(3)$. But not all relations of multiple zeta values are also true for this model. In particular, the space of $\zeta_q^{\rm{BZ}}$ is not spanned by those with just entries with $2$ or $3$. Hirose, Maesaka and Watanabe recently proved the weaker statement that those with entries $\geq 2$ are enough \cite{HMW}. We will return to this in Section \ref{subsubsec:generalmopdifiedqana}. \\

We will mostly be interested in another model of $q$-analogues which is inspired by Eisenstein series and which was introduced by the author during his PhD studies (see \cite{Ba8}) and further studied in \cite{BK1}. These objects are not $q$-analogues in the above sense, but are often called modified $q$-analogues. By a {\bf modified $q$-analogue of weight $k$} we mean that we first need to multiply by $(1-q)^k$ before taking the limit $q\rightarrow 1$. One motivation to consider modified $q$-analogues is the following.

\begin{proposition} \label{prop:mfisqana} Let $f(q) = \sum_{n=0}^\infty a_n q^n \in \mf_k$ be a modular form of weight $k$. Then $f$ is, up to the factor $(2\pi i)^k$, a modified $q$-analogue of weight $k$ of its constant term $a_0$, i.e. we have
    \begin{align*}
    \lim_{q\rightarrow 1} (1-q)^k f(q) = (2\pi i)^k a_0\,.
    \end{align*}
\end{proposition}
\begin{proof} This is a consequence of Proposition \ref{prop:gisqmzv} below together with the fact that every modular form is a polynomial in $\aG_4$ and $\aG_6$. Another way to see this is by using the transformation property $f\left(-\frac{1}{\tau}\right) = \tau^k f(\tau)$. Taking the limit $q\rightarrow 1$ on the real axis corresponds to the limit $\tau \rightarrow 0$ on the positive imaginary axis, since $q=e^{2\pi i \tau}$. Together with $\lim_{\tau \rightarrow i\infty} f(\tau) = a_0$ we obtain
    \begin{align*}
            \lim_{q\rightarrow 1} (1-q)^k f(q) = \lim_{\tau\rightarrow 0} ((2\pi i \tau)^k +  O(\tau^{k+1})) f(\tau) = \lim_{\tau\rightarrow 0} (2\pi i)^k f\left(-\frac{1}{\tau}\right) = \lim_{\tau \rightarrow i\infty} (2\pi i)^k f(\tau) =  (2\pi i)^k a_0\,.
    \end{align*}
\end{proof}

\subsection{The $q$-series $\g$ as modified $q$-analogues of multiple zeta values}

For any $k\geq 2$ we defined in the previous section the Eisenstein series $\aG_k(\tau)$ by its Fourier expansion
\begin{equation*}
\aG_k(\tau) = \zeta(k) + \frac{(-2\pi i)^k}{(k-1)!} \sum_{n=1}^\infty \sigma_{k-1}(n) q^n \,.
\end{equation*}
Clearly we can obtain $\zeta(k)$ from these series in the case $q\rightarrow 0$, but as indicated in Proposition \ref{prop:mfisqana} this will also be possible by considering $q\rightarrow 1$ after some modification. From now on we will always consider $q$ as a formal variable or a fixed complex number with $|q|<1$. For $k\geq 0$ we set
\begin{align*}
\beta(k)=
\begin{cases}
-\dfrac{B_k}{2k!}\,,& k\text{ even},\\[0.2em]
0\,,& k\text{ odd}
\end{cases}
\,.
\end{align*}
In particular, $\beta(0)=-\frac12$ and $\beta(1)=0$. Define for even $k\geq 2$ the $q$-series $G_k = (-2\pi i)^{-k}\aG_k(\tau)$. By Euler's formula $\zeta(2m) = -\frac{B_{2m}}{2 (2m)!} (2\pi i)^{2m}$ (Proposition \ref{prop:euler}) we have
\begin{align}\label{eq:defgk}
G_k = \beta(k) + \frac{1}{(k-1)!} \sum_{n=1}^\infty \sigma_{k-1}(n) q^n   \in \Q[[q]]\,.
\end{align}
The right-hand side of \eqref{eq:defgk} makes sense for any $k\geq 1$ and we will use it to define $G_k$ for all $k\geq 1$. Its remaining part is given by the previously introduced $q$-series
\[ \g(k) = \frac{1}{(k-1)!} \sum_{n=1}^\infty \sigma_{k-1}(n) q^n \,,\]
i.e. $G_k = \beta(k) + \g(k)$.  The $\g(k)$ can be rewritten in the following way 
\begin{align}\label{eq:smallg}
\g(k) &= \frac{1}{(k-1)!} \sum_{n=1}^\infty \sigma_{k-1}(n) q^n = \sum_{m,d>0} \frac{d^{k-1}}{(k-1)!} q^{m d} = \sum_{m>0}     \frac{P_k(q^m)}{(1-q^m)^k} \,,
\end{align}
where for $k\geq 1$ the $P_k(X)\in \Q[X]$ are the \emph{Eulerian polynomials}\footnote{In the literature usually $(k-1)! X^{-1}P_k(X)$ is called Eulerian polynomial.}, defined by
\begin{align*}
\frac{P_k(X)}{(1-X)^k} = \sum_{d>0}\frac{d^{k-1}}{(k-1)!} X^{d}\,.
\end{align*} 
For $k=1,\dots,6$ these are given by 
\begin{align*}
P_1(X)&=P_2(X)=X,\,\quad P_3(X)=\frac{1}{2} X (X+1),\,\quad P_4(X)=\frac{1}{6} X (X^2+4X+1)\,,\\
P_5(X)&= \frac{1}{24} X (X+1) \left(X^2+10 X+1\right),\,\quad P_6(X)=\frac{1}{120} X \left(X^4+26 X^3+66 X^2+26 X+1\right)\,.
\end{align*}
\begin{lemma} \label{lem:eulpol} For all $k\geq 1$ we have $P_k(0)=0$ and $P_k(1)=1$.    
\end{lemma}
\begin{proof}
    This is Exercise \ref{ex3} (i).
\end{proof}

With the same argument as before we get in the arbitrary depth case the following. 
\begin{proposition} For any index $\kk=(k_1,\dots,k_r) \in \Z^r_{\geq 1}$ we have\footnote{These $q$-series are called brackets in \cite{BK1} and are denoted by $[k_1,\dots,k_r]$ there. }
    \begin{align*}
    \g(\kk)=\g(k_1,\dots,k_r) =  \sum_{m_1 > \cdots > m_r > 0} \frac{P_{k_1}(q^{m_1})}{(1-q^{m_1})^{k_1}} \cdots \frac{P_{k_r}(q^{m_r})}{(1-q^{m_r})^{k_r}} \in \Q[[q]]\,.
\end{align*}
\end{proposition}

Notice that this is a well-defined $q$-series for any index (even when $k_1=1$) since $P_k(X)\in X\Q[X]$ (Lemma \ref{lem:eulpol}). These series can be seen as modified $q$-analogues of multiple zeta values, where by modified we mean as before that we need to multiply by a power of $(1-q)$ before taking the limit $q\rightarrow 1$. 
\begin{proposition} \label{prop:gisqmzv} For any admissible index $\kk$ we have
    \begin{align*}
    \lim\limits_{q\rightarrow 1} (1-q)^{\wt(\kk)} \g(\kk) = \zeta(\kk)\,.
    \end{align*}
\end{proposition}
\begin{proof}
    First, we need to justify the interchange of summation and taking the limit. This is the argument of \cite[Proposition 6.4]{BK1}, and the required uniform convergence is proved in \cite[Lemma 6.6]{BK1}. Then this is an easy consequence of Lemma \ref{lem:eulpol} since for $k\geq 1$ we have
    \begin{align*}
    \lim\limits_{q\rightarrow 1} (1-q)^k\frac{P_k(q^m)}{(1-q^m)^k} =     \lim\limits_{q\rightarrow 1} \frac{P_k(q^m)}{[m]_q^k} = \frac{P_k(1)}{m^k} = \frac{1}{m^k}\,.
    \end{align*}
\end{proof}

We will denote the space spanned by all $\g(\kk)$ for any  (not necessarily admissible!) index $\kk$ by
\begin{align*}
\gs = \big\langle \,\g(\kk) \mid \kk \text{ index} \big\rangle_{\Q} \, \subset \, \Q[[q]]\,,
\end{align*}
where we also use the convention $\g(\emptyset)=1$.  In Chapter \ref{sec:algebraicsetup}  we will see (as a consequence of Lemma  \ref{lem:fmtoFM})  that this space is also closed under multiplication and we will prove the following. 

\begin{proposition}\label{prop:gisqalg} The space $\gs$ is a $\Q$-subalgebra of $\Q[[q]]$. \end{proposition}

As we did for multiple zeta values, we will show the lowest depth case of this proposition now. For $k_1,k_2\geq 1$ we have, as in \eqref{eq:simplestuffle}
\begin{align}\label{eq:gprodsimple}
\begin{split}
\g(k_1) \g(k_2)&= \sum_{m_1>0}     \frac{P_{k_1}(q^{m_1})}{(1-q^{m_1})^{k_1}} \sum_{m_2>0}     \frac{P_{k_2}(q^{m_2})}{(1-q^{m_2})^{k_2}} \\ &= \left(\sum_{m_1 >m_2> 0} + \sum_{m_2 >m_1 > 0} + \sum_{m_1 = m_2 > 0}    \right)  \frac{P_{k_1}(q^{m_1})}{(1-q^{m_1})^{k_1}} \frac{P_{k_2}(q^{m_2})}{(1-q^{m_2})^{k_2}} \\
&= \g(k_1,k_2)+\g(k_2,k_1) + \sum_{m>0}  \frac{P_{k_1}(q^{m})}{(1-q^{m})^{k_1}} \frac{P_{k_2}(q^{m})}{(1-q^{m})^{k_2}}\,.
\end{split}
\end{align}

That this is again an element in $\gs$ follows now from the following lemma, which can be proven by using generating series together with the definition of the Bernoulli numbers \eqref{eq:bn}.

\begin{lemma} \label{lem:rproduct}For $k\geq 1$ we set $R_k(X)=\frac{P_{k}(X)}{(1-X)^{k}} =\sum_{d>0}\frac{d^{k-1}}{(k-1)!} X^{d}$. Then for all $k_1,k_2 \geq 1$ we have
    \begin{align}\label{eq:rproduct}
    R_{k_1}(X) \cdot R_{k_2}(X) &=     R_{k_1+k_2}(X)+ \sum_{j=1}^{k_1+k_2-1}\left( \lambda^j_{k_1,k_2}+ \lambda^j_{k_2,k_1} \right) R_j(X)\,
    \end{align}
    where the rational numbers $\lambda^j_{k_1,k_2}$ are given by
    \begin{align*}
    \lambda^j_{k_1,k_2}  = (-1)^{k_2-1} \binom{k_1+k_2-1-j}{k_1-j}  \frac{B_{k_1+k_2-j}}{(k_1+k_2-j)!} \,,
    \end{align*}
    and where we use the convention $\binom{n}{k}=0$ for $k<0$.
\end{lemma}
\begin{proof}
    This is Exercise \ref{ex3} (ii).
\end{proof}
Lemma \ref{lem:rproduct} together with \eqref{eq:gprodsimple} gives the following analogue for the $q$-series $\g$ of the stuffle product formula $\zeta(k_1)\zeta(k_2)=\zeta(k_1,k_2) + \zeta(k_2,k_1) + \zeta(k_1+k_2)$ for multiple zeta values.

\begin{proposition} \label{prop:gstuffle1}For $k_1,k_2\geq 1$ we have
    \begin{align*}
    \g(k_1) \g(k_2) = &\g(k_1,k_2)+\g(k_2,k_1)+\g(k_1+k_2)+ \sum_{j=1}^{k_1+k_2-1}\left( \lambda^j_{k_1,k_2}+ \lambda^j_{k_2,k_1} \right) \g(j)\,.
    \end{align*}
\end{proposition}
\begin{proof}
    This follows immediately by plugging \eqref{eq:rproduct} into \eqref{eq:gprodsimple}.
\end{proof}
We see that the extra terms in the right-hand side are of lower weight and therefore these will vanish when multiplying both sides with $(1-q)^{k_1+k_2}$ and taking the limit $q\rightarrow 1$. In particular we see that this, together with Proposition \ref{prop:gisqmzv}, gives back the stuffle product formula for multiple zeta values. In contrast to $\mz$ (which is conjecturally graded by weight) the space $\gs$ is therefore  not graded by weight.
Section \ref{subsubsec:qshex} treats both of these products (for $\zeta$ and $\g$) simultaneously as examples for a quasi-shuffle product. 
The series $\g$ also satisfy an analogue of the shuffle product formula 
\begin{align*}
\zeta(k_1) \zeta(k_2) = \sum_{j=2}^{k_1+k_2-1} \left( \binom{j-1}{k_1-1} + \binom{j-1}{k_2-1} \right) \zeta(j, k_1+k_2-j)\,,
\end{align*}
which we saw in Proposition \ref{prop:shuffle1}. This formula involves not just the series $\g$ but also its derivative with respect to the differential operator $q \frac{d}{d q}$. 
\begin{proposition} \label{prop:gshuffle1} For $k_1,k_2\geq 1$ and $k=k_1+k_2>2$ we have
    \begin{align}\label{eq:gshuffle1}\begin{split}
    \g(k_1) \g(k_2) = &\sum_{j=1}^{k-1} \left( \binom{j-1}{k_1-1} + \binom{j-1}{k_2-1} \right) \g(j, k-j) \\
    &+ \binom{k-2}{k_1-1} \left(q\frac{d}{dq} \frac{\g(k-2)}{k-2} - \g(k-1) \right)\,.
    \end{split}
    \end{align}
    For $k_1=k_2=1$ we have
    \begin{align}\label{eq:gshuffle11}
    \g(1)^2=2\g(1,1)+\g(2)-\g(1)\,.
    \end{align}
\end{proposition}
We will give a proof of this formula by using generating series below. In general depth we will need, besides the derivative with respect to $q\frac{d}{dq}$, even more extra terms to get an analogue of the shuffle product. This will be discussed in Chapter \ref{sec:mes}. 

\begin{ex}\label{ex:gdshdep1}
    \begin{enumerate}[\textup{(}i\textup{)}] 
        \item 
Similar to the Example \ref{ex:z2z3}, where we showed
    \begin{align}        \label{eq:z2z3prod}
    \zeta(2)\zeta(3) &= \zeta(2,3)+\zeta(3,2)+\zeta(5)=\zeta(2,3) + 3 \zeta(3,2)+ 6 \zeta(4,1)\,,
    \end{align}
    we can use Propositions \ref{prop:gstuffle1} and \ref{prop:gshuffle1} to get the following formulas for the $q$-series $\g$
    \begin{align}\begin{split} \label{eq:g2g3prod}
    \g(2)\g(3) &= \g(2,3)+\g(3,2)+\g(5)-\frac{1}{12} \g(3)\,,\\
    \g(2)\g(3) &=\g(2,3) + 3 \g(3,2)+ 6 \g(4,1)-3 \g(4)+q\frac{d}{dq} \g(3)\,.\end{split}
    \end{align}
    From  this  we deduce the linear relation $\g(5) = 2\g(3,2)+6\g(4,1) + \frac{1}{12}\g(3) - 3 \g(4) +q\frac{d}{dq} \g(3)$. Multiplying equation \eqref{eq:g2g3prod} with $(1-q)^5$ and taking the limit $q \rightarrow 1$ gives the equation \eqref{eq:z2z3prod}. This will be explained in detail in Chapter \ref{sec:mes}.
    \item Since Propositions \ref{prop:gstuffle1} and \ref{prop:gshuffle1} are also valid for non admissible indices we get
    \begin{align*}
    \g(1) \g(2) &= \g(1,2) + \g(2,1) + \g(3) - \frac{1}{2} \g(2)\,,\\
    \g(1) \g(2) &= \g(1,2) + 2\g(2,1) + q \frac{d}{dq} \g(1) - \g(2)\,
    \end{align*}
    by using  $k_1=1, k_2=2$. This gives the following analogue of the relation $\zeta(3)=\zeta(2,1)$
    \begin{align*}
    \g(3) = \g(2,1) + q \frac{d}{dq} \g(1) - \frac{1}{2} \g(2)\,.
    \end{align*}
            \end{enumerate}
\end{ex}

Since the $\g(k)$ are essentially, up to a constant, the Eisenstein series of weight $k$, above formulas can be used to give purely combinatorial proofs of identities among modular forms. One simple example is the identity $\aG_4^2 = \frac{7}{6}\aG_8$, which is a consequence of Proposition %
\ref{prop:gstuffle1} and \ref{prop:gshuffle1} (Exercise \ref{eex3}). In Chapter \ref{sec:mdandmzv} we will elaborate on this combinatorial approach to modular forms. \\

\subsection{Generating series}
We now want to illustrate how the proofs of Proposition  \ref{prop:gstuffle1} and \ref{prop:gshuffle1} can be done by using generating series. This will be done in more generality in Chapter \ref{sec:mes}, but we want to satisfy the curious reader. 
The key point is that there are two different ways to write the generating series of $\g(\kk)$. Multiplying one of them leads to the stuffle product and the other one to the shuffle product. For $r\geq 1$ we will denote the generating series of $\g(k_1,\dots,k_r)$ by
\begin{align*}
\ggen(X_1,\dots,X_r) = \sum_{k_1,\dots,k_r\geq 1} \g(k_1,\dots,k_r) X_1^{k_1-1} \cdots X_r^{k_r-1}\,.
\end{align*}	

\begin{lemma} \label{lem:genexpr}We have 
    \begin{align}\label{eq:stufflegen}
\ggen(X_1,\dots,X_r) 
&= \sum_{m_1> \dots > m_r > 0}  \frac{e^{X_1} q^{m_1}}{1-e^{X_1}q^{m_1}} \cdots  \frac{e^{X_r} q^{m_r}}{1-e^{X_r}q^{m_r}} \\\label{eq:shufflegen}
&= \sum_{m_1> \dots > m_r > 0}  \frac{e^{m_1 X_r} q^{m_1}}{1-q^{m_1}} \frac{e^{m_2 (X_{r-1}-X_r)} q^{m_2}}{1-q^{m_2}} \cdots \frac{e^{m_r (X_{1}-X_2)} q^{m_r}}{1-q^{m_r}}\,.
    \end{align}
\end{lemma}
\begin{proof}
This is Exercise \ref{ex5} (i). The proof of \eqref{eq:stufflegen} follows directly from the definition. For \eqref{eq:shufflegen} a suitable change of summation variables is needed.
\end{proof}

Propositions \ref{prop:gstuffle1} and \ref{prop:gshuffle1} are a consequence of the following proposition by considering the coefficients of $X^{k_1-1} Y^{k_2-1}$. (Exercise \ref{ex5} (iii))
\begin{proposition} \label{prop:ggenproddep1}We have
\begin{align}
\ggen(X) \ggen(Y)  &= \ggen(X,Y)  + \ggen(Y,X) + \frac{1}{e^{X-Y}-1} \ggen(X) + \frac{1}{e^{Y-X}-1}\ggen(Y) \\
&= \ggen(X+Y,X) + \ggen(X+Y,Y) -\ggen(X+Y)  +  q \frac{d}{dq}\sum_{k\geq 1} \g(k) \frac{(X+Y)^k}{k} + \g(2)\,.
\end{align}
\end{proposition}
\begin{proof} Using  \eqref{eq:stufflegen} and \eqref{eq:shufflegen} in the smallest depth case together with the usual splitting of the summation of $m_1,m_2 >0$ into the cases $m_1 > m_2>0$, $m_2 > m_1 > 0$ and $m_1=m_2=m>0$ gives
\begin{align*}
\ggen(X) \ggen(Y) \overset{\eqref{eq:stufflegen}}{=} &\ggen(X,Y)  + \ggen(Y,X) + \sum_{m>0} \frac{e^X q^m}{1-e^X q^m} \frac{e^Y q^m}{1-e^Y q^m}  \,,\\
\ggen(X) \ggen(Y) \overset{\eqref{eq:shufflegen}}{=} &\ggen(X+Y,X) + \ggen(X+Y,Y) + \sum_{m>0} e^{m(X+Y)} \left(\frac{q^m}{1-q^m} \right)^2\,.
\end{align*}
It remains to evaluate the third term in both equations. For the first equation one can check by direct calculation that
\begin{align*}
\frac{e^X q^m}{1-e^X q^m} \frac{e^Y q^m}{1-e^Y q^m} &= \frac{1}{e^{X-Y}-1} \cdot \frac{e^X q^m}{1-e^X q^m} + \frac{1}{e^{Y-X}-1}\cdot   \frac{e^Y q^m}{1-e^Y q^m} \,,
\end{align*}
which then gives 
\begin{align*}
\ggen(X) \ggen(Y) = \ggen(X,Y)  + \ggen(Y,X) + \frac{1}{e^{X-Y}-1} \ggen(X) + \frac{1}{e^{Y-X}-1}\ggen(Y) \,.\\
\end{align*}
For the second equation one uses first $\left(\frac{q^m}{1-q^m} \right)^2 =\frac{q^m}{(1-q^m)^2} - \frac{q^m}{1-q^m}$, which gives
\begin{align*}
\sum_{m>0} e^{m(X+Y)} \left(\frac{q^m}{1-q^m} \right)^2 = \sum_{m>0} e^{m(X+Y)} \frac{q^m}{(1-q^m)^2}  - \ggen(X+Y)\,.
\end{align*}
The first sum on the right can then be evaluated as (Exercise \ref{ex5} (ii))
\begin{align}\label{eq:difgen}
\sum_{m>0} e^{m (X+Y)} \frac{q^m}{(1-q^m)^2} =\g(2) +  q \frac{d}{dq}\sum_{k\geq 1} \g(k) \frac{(X+Y)^k}{k} \,,
\end{align}
from which the claimed formula follows. 
\end{proof}

We end this section by giving another useful expression of the generating series of $\g(k)$. For this we introduce the $q$-Pochhammer symbol defined by 
\begin{align*}
    (a)_\infty = (a;q)_\infty = \prod_{k=0}^{\infty} (1-a q^k) = (1-a) (1-aq) (1-aq^2) \cdots .
\end{align*}
With this we get the following result, which was inspired by the work \cite{Q}.

\begin{proposition}
    We have
    \begin{align*}
        \frac{(e^X q)_\infty }{(q)_\infty} = \exp\left( -\sum_{k\geq 1} \frac{\g(k)}{k} X^k\right). 
    \end{align*}
\end{proposition}
\begin{proof}
By direct calculation we get 
\begin{align*}
    \log\left(  \frac{(e^X q)_\infty }{(q)_\infty}\right)  &= \sum_{m\geq 1} \log\left( \frac{1-e^X q^m }{1-q^m}\right) = \sum_{m\geq 1} \left( \log(1-e^X q^m) - \log(1-q^m) \right)  
    = - \sum_{m,n\geq 1} \frac{e^{nX}-1}{n}  q^{mn}\\
    &= - \sum_{n\geq 1} \frac{e^{nX}-1}{n} \frac{q^n}{1-q^n} = - \int_{0}^X \sum_{n\geq 1} \frac{e^{nZ} q^n}{1-q^n} dZ = - \int_{0}^X \mathfrak{g}(Z) dZ =  -\sum_{k\geq 1} \frac{\g(k)}{k} X^k.
\end{align*}
\end{proof}
As a direct consequence, we get (cf. \cite[equation~(1.2)]{Q})
\begin{align*}
        \frac{ (e^X q)_\infty (e^{-X} q)_\infty}{(q)^2_\infty } = \exp \left( - \sum_{k\geq 1} \frac{g(2k)}{k} X^{2k} \right) \in \Q[g(2),g(4),g(6)]\llbracket X\rrbracket
\end{align*}
and therefore the coefficients of $X^n$ of this expression are quasimodular forms (of mixed weight).

\subsection{General modified $q$-analogues of multiple zeta values}\label{subsubsec:generalmopdifiedqana}
As mentioned before, there are several models of $q$-analogues of multiple zeta values, and in \cite{Zh3} you can find a nice overview of some of them. Most of these have similar definitions as the $\g(\kk)$ with the difference that the Eulerian polynomials get replaced by other polynomials. Also the modified version (by which we mean that we factor out the factor $(1-q)^{\wt(\kk)}$) of the Bradley-Zhao model $\zeta_q^{\rm{BZ}}(\kk)$ is of this form, since instead of $P_k(X)$ the polynomials $X^{k-1}$ are used. In the following, we define a general type of $q$-analogue of multiple zeta values, which were introduced in \cite{BK2}.

\begin{definition}
    For $k_1,\dots,k_r \geq 1$ and polynomials $Q_1(X) \in X \Q[X]$ and $Q_2(X),\dots,Q_r(X) \in \Q[X]$ we define
    \begin{equation} \label{def:Z1}
    \zeta_q(k_1,\dots,k_r ; Q_1 , \dots , Q_r) = \sum_{m_1>\cdots> m_r >0} \frac{Q_1(q^{m_1}) \dots Q_r(q^{ m_r})}{(1-q^{m_1})^{k_1}\cdots (1-q^{m_r})^{k_r}} \,.
    \end{equation}
\end{definition}

As in Proposition \ref{prop:gisqmzv} these  series can be seen as (modified) $q$-analogues of $\zeta(k_1,\dots,k_r)$, since we have for $k_1 \geq 2$  \[ \lim\limits_{q\rightarrow 1} (1-q)^{k_1+\dots+k_r} \zeta_q(k_1,\ldots,k_r; Q_1 , \dots , Q_r) = Q_1(1) \cdots Q_r(1)\cdot \zeta(k_1,\ldots,k_r) \,. \]
We only consider the case where $\deg(Q_j) \leq k_j$ and consider the following $\Q$-vector space: 
\begin{align}\label{eq:defzq}
 \mz_q := \Big \langle   \zeta_q(k_1,\dots,k_r ; Q_1 , \dots , Q_r) \,\,\big|\,\, r\ge 0 ,\,k_1,\dots,k_r\geq1 ,\,\deg(Q_j) \leq k_j 
\Big\rangle_\Q\,,
\end{align}
where again $\zeta_q(\emptyset;\emptyset)=1$. It is again not hard to see that $\mz_q$ is a $\Q$-algebra, since it is again an example for a quasi-shuffle algebra (see Section \ref{subsubsec:qshex}), and we have for example
\[\zeta_q(k_1;Q_1)  \zeta_q(k_2;Q_2)  = \zeta_q(k_1, k_2;Q_1, Q_2)+ \zeta_q(k_2,k_1;Q_2,Q_1)+\zeta_q(k_1+k_2;Q_1 \cdot Q_2)\,.         \]

Since $\g(k_1,\dots,k_r) = \zeta_q(k_1,\dots,k_r; P_{k_1},\dots,P_{k_r})$ and $\deg(P_k) \leq k$ we have $\gs \subset \mz_q$. As we will see in Section \ref{subsubsec:qshex}, we can describe the analogue of the stuffle product for elements in $\gs$ and $\mz_q$ as examples of quasi-shuffle products. 
The reason to introduce the a priori bigger space $\mz_q$ is that we can describe the higher depth analogue of the shuffle product in this space. This we will do in Chapter \ref{sec:mes}. Even though we will not be able to describe the shuffle product analogue in the space $\gs$ explicitly, the two spaces are equal. Until recently, this was a conjecture of the author, discovered during his PhD thesis. It was proved by Hirose, Maesaka and Watanabe in \cite[Theorem~1.1]{HMW}.
\begin{theorem}[{Hirose--Maesaka--Watanabe \cite[Theorem~1.1]{HMW}}] We have $\gs = \mz_q$.
\end{theorem}
This theorem will be discussed in more detail in Chapter \ref{sec:mes}. In \cite{BK2} the authors also introduce the following subspaces of $\mz_q$:
\begin{align*}
 \mz_{q,1} &:= \Big \langle   \zeta_q(k_1,\dots,k_r ; Q_1 , \dots , Q_r) \,\,\big|\,\, r\ge 0 ,\,k_1,\dots,k_r\geq1 ,\,\deg(Q_j) \leq k_j-1
\Big\rangle_\Q\,,\\
 \mz^\circ_{q,1} &:= \Big \langle   \zeta_q(k_1,\dots,k_r ; Q_1 , \dots , Q_r) \,\,\big|\,\, r\ge 0 ,\,k_1,\dots,k_r\geq1 ,\,\deg(Q_j) \leq k_j-1, Q_j(X) \in X \Q[X]
\Big\rangle_\Q\,.
\end{align*}
Also here the $\circ$ can be removed. Until recently, this was another conjecture of the author based on numerical experiments.
\begin{theorem}[{Hirose--Maesaka--Watanabe \cite[Theorem~1.1]{HMW}}]We have $\mz^\circ_{q,1} = \mz_{q,1}$.
\end{theorem}
This theorem states that the modified Bradley-Zhao $q$-analogues, defined by  ($k_1\geq 2, k_2,\dots,k_r\geq 1$)
\begin{align*}
\overline{\zeta}^{\rm{BZ}}_q(k_1,\dots,k_r) &=  \zeta_q(k_1,\dots,k_r; X^{k_1-1}, \dots, X^{k_r-1})  \\
    &= \sum_{m_1>\dots > m_r >0}\frac{q^{(k_1-1)m_1} \cdots q^{(k_r-1)m_r}}{(1-q^{m_1})^{k_1} \cdots (1-q^{m_r})^{k_r}}
\end{align*}
can all be written in terms of those with entries $\geq 2$. This can be seen as a (weak) analogue of the Theorem of Brown (Theorem \ref{thm:brown}), stating that any multiple zeta value can be written in terms of those with entries $2$ or $3$.

\subsection{The BTT-philosophy}\label{subsec:bttphilo}
Finally, we explain the connection between $q$-analogues and finite and symmetric multiple zeta values which appears when the $q$-analogues are evaluated at roots of unity. We will refer to this viewpoint as the \emph{BTT philosophy}: finite and symmetric multiple zeta values arise as algebraic and analytic limits, respectively, when $q\rightarrow 1$.

In analogy to the multiple harmonic sum we define the \emph{multiple harmonic $q$-sum}, for $\kk = (k_1,\dots,k_r) \in \Z^r$ and $m\geq 1$ by 
\begin{align}
 H_m(\kk;q) =    H_m(k_1,\dots,k_r;q) = \sum_{m\geq m_1>\dots>m_r>0} \frac{q^{(k_1-1)m_1} \cdots q^{(k_r-1)m_r}}{[m_1]_q^{k_1} \cdots [m_r]_q^{k_r} } \in \Q[[q]]\,.
\end{align}

Notice that for an admissible index $\kk$ we have 
\begin{align*}
     \lim_{q \rightarrow 1}\lim_{m \rightarrow \infty}  H_m({\bf k};q) =  \lim_{q \rightarrow 1} \zeta_q({\bf k})  = \zeta({\bf k})\,.
\end{align*}

But some kind of magic happens, when one considers these limits in some sense at the same time. By this we mean that we make $q$ dependent on $m$ and then send $m\rightarrow \infty$.
For this we consider the values $H_{n-1}({\bf k};\zeta_n)  \in \Q(\zeta_n)$, where now $q=\zeta_n$ is a primitive $n$-th root of unity. For $n\rightarrow \infty$ we then get $\zeta_n \rightarrow 1$, i.e. in some sense we consider $q \rightarrow 1$ and $n \rightarrow \infty$ at the same time. Doing this for the explicit $n$-th root of unity $\zeta_n = e^{\frac{2\pi i}{n}}$ leads to the following result.

\begin{theorem}\label{thm:xilimit}
For any index set $\kk=(k_{1}, \ldots , k_{r})$ the limit $\lim\limits_{n \rightarrow \infty}    H_{n-1}(\kk ; e^{\frac{2\pi i}{n}})$ exists and we set
\begin{align*}
\xi(\kk) := \lim\limits_{n \rightarrow \infty}   H_{n-1}(\kk ; e^{\frac{2\pi i}{n}}) \in \C\,.
	\end{align*}
It is given by 
\begin{align*}
\xi(\kk) =\sum_{a=0}^{r}(-1)^{k_1+\cdots+k_a}
\zeta^*\big(k_{a}, k_{a-1}, \ldots , k_{1} ;\frac{\pi i}{2}\big)
\zeta^*\big(k_{a+1}, k_{a+2}, \ldots , k_{r};-\frac{\pi i}{2}\big)\,.
\end{align*}
In particular, we have 
\begin{align*}
     {\rm Re} \left(\xi(\kk) \right) \equiv \zeta_{\mathcal{S}} (\kk) \mod \pi^2 \mathcal{Z}\,.
\end{align*} 
\end{theorem}
Hirose's refined symmetric multiple zeta values give another interpretation of the values $\xi(\kk)$ in Theorem \ref{thm:xilimit} \cite[Remark~13]{Hi1}. Notice that Hirose uses the opposite convention for the order of an index.

\begin{theorem}\label{thm:bttfmzv}
For any primitive $p$-th root of unity $\zeta_p$, we have
	\[ ( H_{p-1}(\kk; \zeta_p) \mod \mathfrak{p})_p = \zeta_{\mathcal{A}} (\kk) \,,\]
where $\mathfrak{p}= (1-\zeta_p)$ is the prime ideal of $\Z[\zeta_p]$ generated by $1-\zeta_p$.
\end{theorem}
\begin{proof}
For $p$ prime we have $H_{p-1}(\kk;\zeta_p) \in \Z[\zeta_p]$. This follows from the fact that the $q$-integer $[m]_q$ at $q=\zeta_p$ is a cyclotomic unit, when $m$ is coprime with $p$, since in this case there exists a $t$ with $m \cdot t \equiv 1 \mod p$ and therefore
\[ \frac{1}{[m]_{\zeta_p}} = \frac{1-\zeta_p}{1-\zeta_p^m} = \frac{1-\zeta_p^{tm}}{1-\zeta_p^m} = 1+ \zeta_p^m + \dots + \zeta_p^{(t-1)m} \in \Z[\zeta_p] \,. \]
Moreover, we have $\Z[\zeta_p]/\mathfrak{p} \cong \fp$ (Exercise \ref{fex6}) and for $p>m>0$ we have $[m]_{\zeta_p}\equiv m \mod \mathfrak{p}$. Combining all these shows the desired result.
\end{proof}
The point of the following picture is that there are two different ways of letting
$q$ tend to $1$. Starting with a $q$-multiple zeta value, the usual analytic limit
$q\rightarrow1$ gives a multiple zeta value. If we first truncate at $p-1$ and then
work modulo $[p]_q$, we obtain the $\mq$-multiple zeta values which will be defined in
Section~\ref{sec:Qmzv}. On these objects there are two specializations. Algebraically,
we set $q=1$. Since $[p]_1=p$, this gives a finite multiple zeta value. Analytically,
we evaluate $q=e^{2\pi i/p}$ and let $p\rightarrow\infty$. Theorem
\ref{thm:xilimit} shows that this gives a symmetric multiple zeta value after taking
the real part modulo $\pi^2\mz$.

The green and orange outer arrows describe the corresponding operations without the
intermediate $\mq$-object: truncate and reduce modulo $p$ on the finite side, and
symmetrize on the symmetric side. Thus the arrows in the picture are not all maps of
the same kind. The main point is that the two red arrows start with the same
$\mq$-multiple zeta value. Therefore a $\Q$-linear relation among ordinary
$\mq$-multiple zeta values for which the analytic limits exist gives a finite
relation and a symmetric relation of the same shape. This gives an explanation
for the Kaneko--Zagier conjecture.
\begin{figure}[h!]
    \begin{center}
\includegraphics[width=0.80\textwidth]{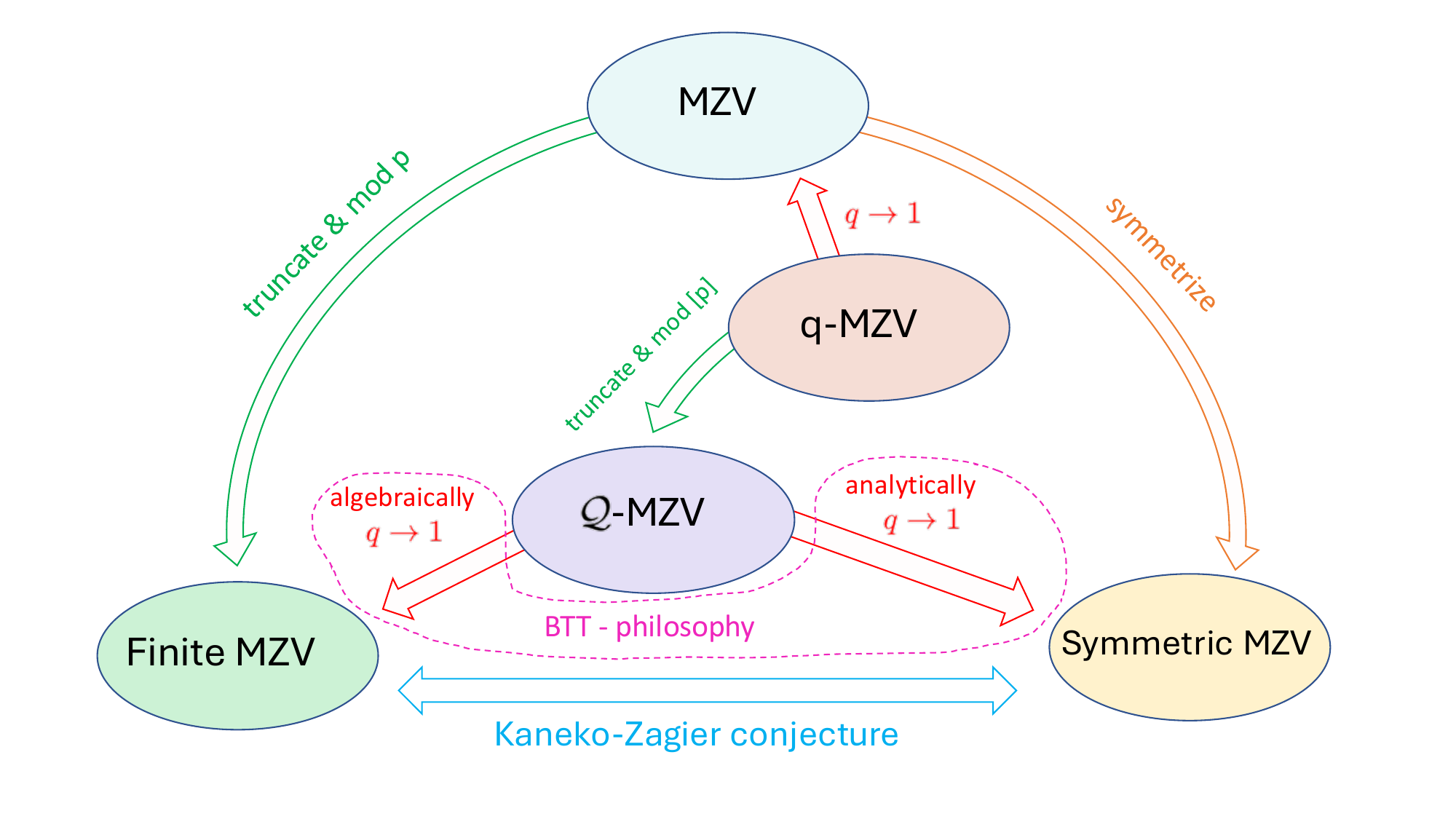}
    \end{center}
    \caption{Overview of the BTT philosophy}
    \label{fig:bttphilosophy}
\end{figure}

{\bf Goal:} Understand the values $H_{n-1}(\kk; \zeta_n)$ for primitive $n$-th roots of unity $\zeta_n$ (e.g. $\zeta_n=e^{\frac{2\pi i}{n}}$) and their relations to get simultaneous results on finite and symmetric multiple zeta values. \\

For this, we start by considering the depth one case again. Carlitz (1956) \cite{Car} introduced the {\bf degenerate Bernoulli} numbers $b_k(n) \in \Q[\frac{1}{n}]$ 	
\[\sum_{k=0}^\infty b_k(n) \frac{x^k}{k!} = \frac{x}{(1+\frac{x}{n})^n-1}\,.\]
These numbers can be seen as a degeneration of the Bernoulli numbers
\[ \lim\limits_{n \rightarrow \infty} b_k(n) = B_k\,,\qquad \sum_{k=0}^\infty B_k \frac{x^k}{k!} = \frac{x}{e^{x}-1}\,.\]

For even $k\geq 2$ we saw that $\zeta(k) = -\frac{B_k}{k!} \frac{(2\pi i)^k}{2}$ and as an analogue we have

\begin{proposition}\label{qhnindep1}
For all $n,k \geq 1$ we have
\[H_{n-1}(k; e^{\frac{2\pi i}{n}}) = -\frac{b_k\left( n \right)}{k!} \left(n (1 - e^{\frac{2\pi i}{n}})\right)^k \,.\]
\end{proposition}
\begin{proof}
This is Exercise \ref{fex7}.
\end{proof}

For an index set ${\bf k}=(k_1,\dots,k_r)$ we define its reverse by $\overline{\kk} = (k_r,\dots,k_1)$. 
\begin{theorem}\label{thm:btthoffmandual}
For all $n\geq 1$  and  all non-empty indices ${\bf k}$ we have for any primitive $n$-th root of unity $\zeta_n$
\[ H_{n-1}^\star({\bf k};\zeta_n) = (-1)^{\wt({\bf k})+1} H_{n-1}^\star(\overline{{\bf k}^\vee};\zeta_n) \,.\]
\end{theorem}
 
Theorem \ref{thm:fmzvhoffmandual} ($\zeta_\mathcal{A}^\star({\bf k}) = -\zeta_\mathcal{A}^\star({\bf k}^\vee)$)
 follows from this together with the reversal relation $\zeta^\star_{\mathcal{A}}(\kk) = (-1)^{\wt(\kk)}\zeta^\star_{\mathcal{A}}(\overline{\kk})$.

\section{Schur multiple zeta values}\label{subsec:schurmzvoverview}

We finish this overview with another generalization of multiple zeta values. So far, all versions of multiple zeta values in this chapter were indexed by tuples of positive integers. For Schur multiple zeta values, the tuple is replaced by a Young diagram whose boxes are filled by positive integers. These values were introduced by Nakasuji, Phuksuwan and Yamasaki in \cite{NPY}.

A Young diagram consists of finitely many left-aligned rows of boxes whose lengths weakly decrease from top to bottom. A skew Young diagram is obtained by removing a smaller Young diagram from its upper-left corner. We denote the set of boxes of such a diagram by $D$. For a filling $\mathbf m=(m_{i,j})_{(i,j)\in D}$ by positive integers we write $\mathbf m\in\operatorname{SSYT}(D)$ if the entries are weakly increasing along rows and strictly increasing down columns, i.e.
\begin{align*}
m_{i,j}\leq m_{i,j+1}\,,\qquad m_{i,j}<m_{i+1,j}
\end{align*}
whenever the two boxes in question belong to $D$.

A box is called a corner if neither the box immediately to its right nor the box immediately below it belongs to $D$.

\begin{definition}
For a filling $\mathbf h=(h_{i,j})_{(i,j)\in D}$ by positive integers which are at least $2$ in every corner, we define the corresponding \emph{Schur multiple zeta value} by
\begin{align*}
\zeta(\mathbf h)=\sum_{\mathbf m\in\operatorname{SSYT}(D)}
\prod_{(i,j)\in D}\frac{1}{m_{i,j}^{h_{i,j}}}\,.
\end{align*}
We call $\wt(\mathbf h)=\sum_{(i,j)\in D}h_{i,j}$ its \emph{weight}.
\end{definition}

The condition in the definition ensures that the series converges absolutely (see \cite[Lemma~2.1 and Section~4.1]{NPY}). For example, for $a\geq1$ and $b,c\geq2$ we have
\begin{align*}
\zeta\left(
{\ytableausetup{centertableaux,boxsize=1.4em}
\begin{ytableau}
a&b\\c
\end{ytableau}}
\right)
=\sum_{
\arraycolsep=1.4pt\def\arraystretch{0.8}
\begin{array}{ccc}
0<&m_a &\leq m_b \\
&\vsmall & \, \\
&m_c & \,
\end{array} }
\frac{1}{m_a^a\cdot m_b^b \cdot m_c^c}\,.
\end{align*}

For a non-empty admissible index $\kk=(k_1,\dots,k_r)$ we also define the \emph{multiple zeta-star value} by
\begin{align*}
\zeta^\star(\kk)=\zeta^\star(k_1,\dots,k_r)
=\sum_{m_1\geq\dots\geq m_r>0}\frac{1}{m_1^{k_1}\cdots m_r^{k_r}}\,.
\end{align*}
We also set $\zeta^\star(\emptyset)=1$. Multiple zeta values and multiple zeta-star values are the special cases of Schur multiple zeta values given by a column and a row, respectively. Because we use decreasing summation variables, the order of the entries is reversed. For $h_1,h_2\geq1$ and $h_3\geq2$ we have
\begin{align*}
\zeta\left(
{\ytableausetup{centertableaux,boxsize=1.4em}
\begin{ytableau}
h_1\\h_2\\h_3
\end{ytableau}}
\right)
&=\zeta(h_3,h_2,h_1),&
\zeta\left(
{\ytableausetup{centertableaux,boxsize=1.4em}
\begin{ytableau}
h_1&h_2&h_3
\end{ytableau}}
\right)
&=\zeta^\star(h_3,h_2,h_1).
\end{align*}

\begin{ex}\label{ex:schurmzvoverview}
\begin{enumerate}[\textup{(}i\textup{)}]
\item For $a,b\geq2$, separating the strict inequality from the equality gives
\begin{align*}
\zeta\left(
{\ytableausetup{centertableaux,boxsize=1.4em}
\begin{ytableau}
a&b
\end{ytableau}}
\right)
=\zeta^\star(b,a)=\zeta(b,a)+\zeta(a+b)
=\zeta(a)\zeta(b)-\zeta(a,b).
\end{align*}
Equivalently, the harmonic product formula \eqref{eq:simplestuffle} can be written in terms of Young diagrams as
\begin{align*}
\zeta\left(
{\ytableausetup{centertableaux,boxsize=1.4em}
\begin{ytableau}a\end{ytableau}}
\right)
\zeta\left(
{\ytableausetup{centertableaux,boxsize=1.4em}
\begin{ytableau}b\end{ytableau}}
\right)
=\zeta\left(
{\ytableausetup{centertableaux,boxsize=1.4em}
\begin{ytableau}a&b\end{ytableau}}
\right)
+\zeta\left(
{\ytableausetup{centertableaux,boxsize=1.4em}
\begin{ytableau}b\\a\end{ytableau}}
\right).
\end{align*}
In particular, we obtain
\begin{align*}
\zeta\left(
{\ytableausetup{centertableaux,boxsize=1.4em}
\begin{ytableau}2\\2\end{ytableau}}
\right)=\frac{3}{4}\zeta(4),\qquad
\zeta\left(
{\ytableausetup{centertableaux,boxsize=1.4em}
\begin{ytableau}2&2\end{ytableau}}
\right)=\frac{7}{4}\zeta(4).
\end{align*}

\item For $a,b,c\geq2$ we have the following relation (see also \cite[Example~2.5]{Ba4}, written in our convention):
\begin{align*}
\zeta\left(
{\ytableausetup{centertableaux,boxsize=1.4em}
\begin{ytableau}
a&b\\c
\end{ytableau}}
\right)
&=\zeta(c,a)\zeta(b)-\zeta(c,a,b)\\
&=\zeta(b,c,a)+\zeta(c,b,a)+\zeta(b+c,a)+\zeta(c,a+b).
\end{align*}
Indeed, the last line is obtained by splitting the summation into the four cases
\begin{align*}
m_b>m_c>m_a,\qquad m_c>m_b>m_a,\qquad m_b=m_c>m_a,\qquad m_c>m_b=m_a.
\end{align*}
The same four terms, together with $\zeta(c,a,b)$, appear when expanding $\zeta(c,a)\zeta(b)$, which proves the first equality. Setting $a=b=c=2$ and using the formula for $\zeta(\{2\}^n)$ gives
\begin{align*}
\zeta\left(
{\ytableausetup{centertableaux,boxsize=1.4em}
\begin{ytableau}
2&2\\2
\end{ytableau}}
\right)
=\zeta(2,2)\zeta(2)-\zeta(2,2,2)
=\frac{\pi^6}{840}=\frac{9}{8}\zeta(6).
\end{align*}

\item A checkerboard filling with alternating entries $1$ and $3$ gives the
following determinant of odd zeta values \cite[Introduction]{BY1}:
\begin{align*}
\zeta\left(
{\ytableausetup{centertableaux,boxsize=1.4em}
\begin{ytableau}
3&1&3\\
1&3\\
3
\end{ytableau}}
\right)
=\frac{1}{16}\det\begin{pmatrix}
\zeta(3)&\zeta(7)\\
\zeta(7)&\zeta(11)
\end{pmatrix}.
\end{align*}
The general checkerboard formulas will be discussed in
Section~\ref{sec:schurmzv}.
\end{enumerate}
\end{ex}

By splitting the summation domain in the same way, we see that every Schur multiple zeta value is a finite sum of ordinary multiple zeta values of the same weight: we order all tableau entries and combine the exponents belonging to equal entries.

For fillings which are constant along diagonals and whose entries are all at least $2$, Schur multiple zeta values also satisfy an analogue of the Jacobi--Trudi formula: they can be written as determinants whose entries are ordinary multiple zeta values. The skew version was proved in \cite[Theorem~4.3~(2)]{NPY}. An interpolated version for straight Young diagrams and its lattice path proof are given in \cite[Theorems~1.1 and 4.7]{Ba4}. We will return to Schur multiple zeta values in Section~\ref{sec:schurmzv}, where we will introduce the partition notation for skew Young diagrams and prove the Jacobi--Trudi formula.

\vspace{1cm}
\begin{center}
\bf{ \Large {\color{qmzvlinecol}\ding{94}}~ Exercises ~{\color{qmzvlinecol}\ding{94}}}\end{center}

The following is a collection of exercises intended to deepen the reader's understanding of this chapter.

\begin{exer}\label{ex1} 
    \begin{enumerate}[\textup{(}i\textup{)}] 
        \item 
        Prove Proposition \ref{prop:shuffle1}, i.e. show that for $k_1,k_2\geq 2$ we have
        \begin{align*}
        \zeta(k_1) \zeta(k_2) = \sum_{j=2}^{k_1+k_2-1} \left( \binom{j-1}{k_1-1} + \binom{j-1}{k_2-1} \right) \zeta(j, k_1+k_2-j).
        \end{align*}
        Hint: Use for $a,b \geq 1$ the partial fraction expansion
        \begin{align*}
        \frac{1}{x^a y^b} = \sum_{j=1}^{a+b-1} \left( \frac{\binom{j-1}{a-1}}{(x+y)^j y^{a+b-j}} + \frac{\binom{j-1}{b-1}}{(x+y)^j x^{a+b-j}} \right).
        \end{align*}
        (Optional: Give a proof of this formula)
        \item Use (i) together with $\zeta(k_1)\zeta(k_2)=\zeta(k_1,k_2) + \zeta(k_2,k_1) + \zeta(k_1+k_2)$ to prove the relations
        \begin{align*}
        \zeta(6) &= 6 \zeta(3,3) - 3 \zeta(4,2),\\
        \zeta(7) &= 4 \zeta(3,4) + 3 \zeta(4,3) - 2 \zeta(5,2).
        \end{align*}
    \end{enumerate}
\end{exer}

\vspace{0.5cm}

\begin{exer}\label{ex2}
    \begin{enumerate}[\textup{(}i\textup{)}] 
        \item Show that Conjecture \ref{conj:graded} together with Proposition \ref{prop:mzissubalgebra} would imply that all multiple zeta values (except for $\zeta(\emptyset)=1$) are transcendental.
        \item Show that Conjecture \ref{conj:hoffman} (Hoffman) would imply  Conjecture \ref{conj:zagier} (Zagier).
        \item Show that Conjecture \ref{conj:brokrei} (Broadhurst-Kreimer) would imply  Conjecture \ref{conj:zagier} (Zagier).
    \end{enumerate}
\end{exer}

\vspace{0.5cm}

\begin{exer}\label{fex14}
\begin{enumerate}[(i)]
    \item Show that there exists an injective $\Q$-algebra homomorphism 
    \begin{align*}
        \iota: \Q \rightarrow \mathcal{A},
    \end{align*}
    i.e., prove that $\mathcal{A}$ is a $\Q$-algebra.
    Furthermore, show that $\iota$ cannot be surjective.
    
    \item Show that $\sqrt{2} \not\in \mathcal{A}$, i.e., prove that there does not exist any $a \in \mathcal{A}$ such that $a^2 - 2 = 0$ in $\mathcal{A}$.
\end{enumerate}
\end{exer}

\vspace{0.5cm}

\begin{exer}\label{fex2} Find a formula for the number of indices and admissible indices of a given weight (and depth), i.e. for $r,k \geq 0$ find explicit expressions for the following four integers 
\begin{align*}
    I_k &=\sum_{r\geq 0} I_{k,r}\,,\qquad 
        I_{k,r} = \big| \{ \kk \in \Z^r_{\geq 1} \mid  \wt(\kk) = k\}\big|\,,\\
      I^0_k &=\sum_{r\geq 0} I^0_{k,r}\,,\qquad
        I^0_{k,r} = \big|\{ \kk \in  \Z^r_{\geq 1} \mid  \wt(\kk) = k, \kk \text{ is admissible}\}\big|\,.
\end{align*}
\end{exer}

\vspace{0.5cm}

\begin{exer}\label{fex4} Show that for any $k_1,\dots,k_r \in \Z$ we have $\za(k_1,\dots,k_r) \in \mza$, i.e. show that you can write $\za(k_1,\dots,k_r)$ as a linear combination of finite multiple zeta values with positive entries. 
\end{exer}

\vspace{0.5cm}

\begin{exer}\label{fex5}
Show that for all $m\geq 1$ and $k_1,\dots,k_r \geq 1$ we have 
\begin{align*}
    S_m(k_1,\dots,k_r)  = \sum_{j=0}^r (-1)^{k_1+\dots + k_j} H_m(k_j,k_{j-1},\dots,k_1) H_m(k_{j+1},\dots,k_r) \,,
\end{align*}
where $S_m$ is defined by \eqref{eq:defsm} and $H_m$ by \eqref{eq:defhm}.
\end{exer}

\vspace{0.5cm}

\begin{exer}\label{fex6} Let $\zeta_p$ be a primitive $p$-th root of unity and  consider the ideal $\mathfrak{p}= (1-\zeta_p)$ of $\Z[\zeta_p]$.
\begin{enumerate}[(i)]
    \item Show that $\mathfrak{p}$ is prime.
    \item Show that $\Z[\zeta_p]/\mathfrak{p} \cong \fp$.
\end{enumerate}
\end{exer}

\vspace{0.5cm}

\begin{exer}\label{fex7}
\begin{enumerate}[(i)]
\item Prove Proposition \ref{qhnindep1}, i.e. show that for all $n,k \geq 1$ we have
\[H_{n-1}(k; e^{\frac{2\pi i}{n}}) = -\frac{b_k\left( n \right)}{k!} \left(n (1 - e^{\frac{2\pi i}{n}})\right)^k \,.\]
\item Show that $b_k\left( n \right)n^k \in \Q[n]$.
\item Use (i) to give another proof of $\za(k)=0$ and $\zs(k)=0$ using Theorems \ref{thm:xilimit} and \ref{thm:bttfmzv}.
\end{enumerate}
\end{exer}

\vspace{0.5cm}

\begin{exer}\label{ex3} We defined for $k\geq 1$ the  Eulerian polynomials $P_k(X)$ and the power series $R_k(X)$ by  
    \begin{align*}
    R_{k}(X) =    \frac{P_k(X)}{(1-X)^k} = \sum_{d>0}\frac{d^{k-1}}{(k-1)!} X^{d}\,.
    \end{align*} 
    \begin{enumerate}[\textup{(}i\textup{)}] 
        \item Prove Lemma \ref{lem:eulpol}, i.e. show that we have $P_k(0)=0$ and $P_k(1)=1$ for all $k\geq 1$.
        \item Prove Lemma \ref{lem:rproduct}, i.e. show that for all $k_1,k_2 \geq 1$
        \begin{align*}
        R_{k_1}(X) \cdot R_{k_2}(X) &=     R_{k_1+k_2}(X)+ \sum_{j=1}^{k_1+k_2-1}\left( \lambda^j_{k_1,k_2}+ \lambda^j_{k_2,k_1} \right) R_j(X)\,,
        \end{align*}
        where the rational numbers $\lambda^j_{k_1,k_2}$ are given by
        \begin{align*}
        \lambda^j_{k_1,k_2}  = (-1)^{k_2-1} \binom{k_1+k_2-1-j}{k_1-j}  \frac{B_{k_1+k_2-j}}{(k_1+k_2-j)!} \,,
        \end{align*}
        and where we use the convention $\binom{n}{k}=0$ for $k<0$.
    \end{enumerate}
\end{exer}

\vspace{0.5cm}

\begin{exer} \label{eex3} Define for even $k\geq 4$ the normalized Eisenstein series $E_k$ by $E_k = \zeta(k)^{-1} \aG_k \in \mf_k$. Show that we have 
    \begin{align*}
    E_k = 1 - \frac{2 k!}{B_k} \g(k)
    \end{align*}
    and {\bf give two different proofs} of the identity 
    \begin{align}\label{eq:e4e8}
    E_4^2 = E_8
    \end{align}between the Eisenstein series
    \begin{align*}
    E_4&= 1+ 240 \sum_{n=1}^\infty \sigma_{3}(n) q^n \,\quad \text{and}     \quad  E_8= 1+  480\sum_{n=1}^\infty \sigma_{7}(n) q^n\,.
    \end{align*}
    \begin{enumerate}[\textup{(}i\textup{)}] 
        \item ``Modular'' proof of \eqref{eq:e4e8}: Use the theory of modular forms, i.e. use Theorem \ref{thm:mf} and \ref{thm:mfdim}.
        \item Combinatorial proof of \eqref{eq:e4e8}: You are just allowed to use Proposition \ref{prop:gstuffle1} and \ref{prop:gshuffle1}.
    \end{enumerate}
\end{exer}

\vspace{0.5cm}

\begin{exer} \label{ex5}
    \begin{enumerate}[\textup{(}i\textup{)}] 
        \item 
    Prove Lemma \ref{lem:genexpr}, i.e. show that 
        \begin{align*}
    \mathfrak{g}(X_1,\dots,X_r) = \sum_{k_1,\dots,k_r\geq 1} g(k_1,\dots,k_r) X_1^{k_1-1} \cdots X_r^{k_r-1}
    \end{align*}	
    can be written in the following two ways
    \begin{align*}
        \mathfrak{g}(X_1,\dots,X_r) 
        &= \sum_{m_1> \dots > m_r > 0}  \frac{e^{X_1} q^{m_1}}{1-e^{X_1}q^{m_1}} \cdots  \frac{e^{X_r} q^{m_r}}{1-e^{X_r}q^{m_r}} \\
        &= \sum_{m_1> \dots > m_r > 0}  \frac{e^{m_1 X_r} q^{m_1}}{1-q^{m_1}} \frac{e^{m_2 (X_{r-1}-X_r)} q^{m_2}}{1-q^{m_2}} \cdots \frac{e^{m_r (X_{1}-X_2)} q^{m_r}}{1-q^{m_r}}\,.
        \end{align*}
        \item Show equation \eqref{eq:difgen}, i.e. show that we have
        \begin{align*}
        \sum_{m>0} e^{m X} \frac{q^m}{(1-q^m)^2} =\g(2) +  q \frac{d}{dq}\sum_{k\geq 1} \g(k) \frac{X^k}{k} \,.
        \end{align*}
        \item Show that Propositions \ref{prop:gstuffle1} and \ref{prop:gshuffle1} are a consequence of Proposition \ref{prop:ggenproddep1}.
            \end{enumerate}
\end{exer}

\vspace{0.5cm}

%% file: TikZ/tikz_infinitep.tex
    \begin{tikzpicture}[
    scale=0.45,
       shorten > = 1pt,auto,
   node distance = 3cm,
      decoration = {snake,   %
                    pre length=3pt,post length=7pt,%
                    },
        every label/.append style={ font=\normalsize}]
        \foreach \cent in {0,10}
        {
        \draw[] (\cent + 3,0) arc (0:180:3);
        \draw[dotted] (\cent -3,0) arc (180:225:3);
        \draw[dotted] (\cent + 3,0) arc (0:-45:3);
        \draw[] (\cent -2.12,-2.12) arc (225:315:3);
        
        \node[text=gray] at (\cent-1.53,1.28){$m_1$};
        \node[text=gray] at (\cent-2,-.1){$m_2$};
        \node[text=gray] at (\cent + 1.53,1.28){$m_r$};
        \node[text=gray] at (\cent + 2,-.1){$m_{r-1}$};
        \node[anchor = north,text=gray, rotate = 70] at (\cent -2.377,0.77) {$<$}; 
        \node[anchor = north,text=gray, rotate = 120] at (\cent -2.165,-1.15) {$<$}; 
        \node[anchor = north,text=gray, rotate = -70] at (\cent + 2.377,0.77) {$<$}; 
        \node[anchor = north,text=gray, rotate = -120] at (\cent + 2.165,-1.15) {$<$}; 
        \draw[dotted, thick, gray] (\cent -1.414, -1.414) arc(225:315:2);
        }
        
        \draw[] (20 + 3,0) arc (0:180:3);
        \draw[dotted] (20 -3,0) arc (180:225:3);
        \draw[dotted] (20 + 3,0) arc (0:-45:3);
        \draw[] (20 -2.12,-2.12) arc (225:315:3);
        
        \node[text=gray] at (20-1.53,1.28){$m_1$};
        \node[text=gray] at (20-2,-.1){$m_2$};
        \node[text=gray] at (20 + 1.53,1.28){$m_r$};
        \node[text=gray] at (20 + 2,-.1){$m_{r-1}$};
        \node[anchor = north,text=gray, rotate = 70] at (20 -2.377,0.77) {$\prec$}; 
        \node[anchor = north,text=gray, rotate = 120] at (20 -2.165,-1.15) {$\prec$}; 
        \node[anchor = north,text=gray, rotate = -70] at (20 + 2.377,0.77) {$\prec$}; 
        \node[anchor = north,text=gray, rotate = -120] at (20 + 2.165,-1.15) {$\prec$}; 
        \draw[dotted, thick, gray] (20 -1.414, -1.414) arc(225:315:2);
        
        \node[circle, draw=black,    inner sep= 1pt, fill = black, label = {90:$0$}] at (0 ,3) {};
        \node[circle, draw=black,    inner sep= 1pt, fill = black, label = {126:$p-1$}] at (0 -1.76,2.427) {};        
        \node[circle, draw=black,    inner sep= 1pt, fill = black, label = {162:$p-2$}] at (0 -2.85,0.927) {};
        \node[circle, draw=black,    inner sep= 1pt, fill = black, label = {54:$1$}] at (0 + 1.76,2.427) {};        
        \node[circle, draw=black,    inner sep= 1pt, fill = black, label = {18:$2$}] at (0 + 2.85,0.927) {};
        \node[circle, draw=black,    inner sep= 1pt, fill = black, label = {288:$\frac{p-1}{2}$}] at (0 + 0.927, -2.85) {};
        \node[circle, draw=black,    inner sep= 1pt, fill = black, label = {262:$\frac{p+1}{2}$}] at (0 - 0.927, -2.85) {};
        
        \node[circle, draw=black,    inner sep= 1pt, fill = black, label = {90:$0$}] at (10 ,3) {};
        \node[circle, draw=black,    inner sep= 1pt, fill = black, label = {126:$-1$}] at (10 -1.76,2.427) {};        
        \node[circle, draw=black,    inner sep= 1pt, fill = black, label = {162:$-2$}] at (10 -2.85,0.927) {};
        \node[circle, draw=black,    inner sep= 1pt, fill = black, label = {54:$1$}] at (10 + 1.76,2.427) {};        
        \node[circle, draw=black,    inner sep= 1pt, fill = black, label = {18:$2$}] at (10 + 2.85,0.927) {};
        \node[circle, draw=black,    inner sep= 1pt, fill = black, label = {288:$\frac{p-1}{2}$}] at (10 + 0.927, -2.85) {};
        \node[circle, draw=black,    inner sep= 1pt, fill = black, label = {262:$-\frac{p-1}{2}$}] at (10 - 0.927, -2.85) {};
        
        \node[circle, draw=black,    inner sep= 1pt, fill = black, label = {90:$0$}] at (20 ,3) {};
        \node[circle, draw=black,    inner sep= 1pt, fill = black, label = {126:$-1$}] at (20 -1.76,2.427) {};        
        \node[circle, draw=black,    inner sep= 1pt, fill = black, label = {162:$-2$}] at (20 -2.85,0.927) {};
        \node[circle, draw=black,    inner sep= 1pt, fill = black, label = {54:$1$}] at (20 + 1.76,2.427) {};        
        \node[circle, draw=black,    inner sep= 1pt, fill = black, label = {18:$2$}] at (20 + 2.85,0.927) {};
        \node[circle, draw=black,    inner sep= 1pt, fill = black, label = {}] at (20 + 0.927, -2.85) {};
        \node[circle, draw=black,    inner sep= 1pt, fill = black, label = {}] at (20 - 0.927, -2.85) {};
        \node[circle, draw=black, inner sep=0pt, label = {270: $-\infty = \infty$}] at (20, -3){};
        \node[circle, inner sep=0pt, label = {270: mod $p$}] at (5,0) {};
        \path[draw=black, decorate,->] (3.75,0) -- (6.25,0);
        \path[draw=black, decorate,->] (13.75,0) -- (16.25,0);
        \node[circle, inner sep=0pt, label = {270: $p\to\infty$}] at (15,0) {};
    \end{tikzpicture}

%% file: chap_AlgebraicSetup.tex
\chapter{Algebra Setup}\label{sec:algebraicsetup}
In this chapter, we want to explain the algebraic structure of the spaces $\mz$ and $\mza$ (multiple and finite multiple zeta values), $\gs$ (the space of the $q$-analogues $\g(\kk)$) and $\mz_q$ ($q$-analogues of multiple zeta values). In particular, we will see why these spaces are all $\Q$-algebras. 
For this, we use Hoffman's algebraic setup from \cite{H2} and its generalization to quasi-shuffle algebras in \cite{HI}.

\section{Multiple polylogarithms, iterated integrals and duality}
In the following, we want to introduce the iterated integral expression for multiple zeta values. This will be used in the next subsection to give another explanation of the shuffle product formula in \eqref{eq:shuffle1}. We start by calculating one simple example by hand, before giving a general formula afterwards. Consider the following iterated integral
\begin{align}\begin{split}\label{eq:z2itint} 
\int_{0}^1 \frac{dt_1}{t_1} \int_{0}^{t_1}  \frac{dt_2}{1-t_2} &=     \int_{0}^1 \frac{dt_1}{t_1} \int_{0}^{t_1}  \sum_{n=0}^\infty t_2^n dt_2 = \int_{0}^1 \frac{dt_1}{t_1} \left[   \sum_{n=0}^\infty \frac{t_2^{n+1}}{n+1} \right]_0^{t_1}\\ &= \int_{0}^1   \sum_{n=0}^\infty \frac{t_1^{n}}{n+1}  dt_1 = \left[  \sum_{n=0}^\infty \frac{t_1^{n+1}}{(n+1)^2} \right]_0^{1} = \sum_{m>0} \frac{1}{m^2} = \zeta(2)\,.
\end{split}
\end{align}
With the same idea, one can also show that we have (Exercise \ref{ex6} (i))
\begin{align}\label{eq:z32itint}
\zeta({\color{mygreen}2},{\color{myblue}3}) =     {\color{mygreen}\int_{0}^1 \frac{dt_1}{t_1} \int_{0}^{t_1}  \frac{dt_2}{1-t_2}} {\color{myblue} \int_{0}^{t_2} \frac{dt_3}{t_3} \int_{0}^{t_3} \frac{dt_4}{t_4} \int_{0}^{t_4}  \frac{dt_5}{1-t_5}}\,.
\end{align}
In general, we will see that an index $\kk=(k_1,\dots,k_r)$ corresponds to an iterated integral of length $\wt(\kk)$, where each $k_j$ gives a block of $k_j-1$ integrals over $\frac{dt}{t}$ and one integral over $\frac{dt}{1-t}$.
To prove these iterated integrals in general, we will introduce multiple polylogarithms, which can be seen as a simultaneous generalization of the polylogarithm ($r=1$) and multiple zeta values  ($z=1$).
\begin{definition} For $|z|<1$ and  $\kk =(k_1,\dots,k_r)\in \Z_{\geq 1}^r$ we define the \emph{multiple polylogarithm} by
\begin{align*}
\li_\kk(z) = \li_{k_1,\dots,k_r}(z) = \sum_{m_1 > \cdots > m_r > 0} \frac{z^{m_1}}{m_1^{k_1} \cdots m_r^{k_r}}\,
\end{align*}
and set $\li_{\emptyset}(z)=1$.        
\end{definition}
For an arbitrary index $\kk$ the $\li_\kk(z)$ are holomorphic functions in the open unit disc, but clearly when $\kk$ is admissible $\li_\kk(z)$ is also defined for $z=1$ and we have  
\begin{align*}
\li_\kk(1)=\zeta(\kk)\,.
\end{align*}
Multiple polylogarithms also have an iterated integral expression, and for example using the same calculation as in \eqref{eq:z2itint} we see for example that
\begin{align*}
\li_2(z) =     \int_{0}^z \frac{dt_1}{t_1} \int_{0}^{t_1}  \frac{dt_2}{1-t_2} \,.
\end{align*}
It becomes clear that we will deal with iterated integrals of two different differential forms. To describe the iterated integrals and the shuffle product, we will therefore introduce the following algebraic setup. 
\subsection{The spaces $\HH$, $\HH^1$ and $\HH^0$ }
We denote by $\HH=\Q\langle x,y\rangle$ the polynomial ring in the two non-commutative variables $x$ and $y$. A monomial in $x$ and $y$ will also be called a \emph{word}, and $\HH$ is therefore the $\Q$-vector space spanned by all words in the \emph{letters} $x$ and $y$. 
Further, we define the subspace $\HH^1 = \Q + \HH y$, which is spanned by the \emph{empty word} $1$ and all words in $x$ and $y$ which end in $y$. For $k\geq 1$ we define \[ z_k =  x^{k-1}y \,.\]
With this we see that $\HH^1 = \Q\langle z_1, z_2, \dots \rangle$, i.e. we could say that $\HH^1$ is spanned by all words in the letters $z_k$. For an index $\kk =(k_1,\dots,k_r)\in \Z_{\geq 1}^r$ we define 
\[ z_{\kk} = z_{k_1} z_{k_2} \cdots z_{k_r} \in \HH^1\]
and set $z_\emptyset = 1$. Now define the space $\HH^0 = \Q + x\HH y$, which is the subspace of $\HH^1$ generated by all words which start in $x$ and end in $y$. In other words, $\HH^0$ is spanned by all $z_\kk$ with admissible indices $\kk$. Summarizing everything we  have
\begin{align*}
\HH^0 = \langle z_\kk \mid \kk \text{ admissible index} \rangle_\Q\quad \subset \quad \HH^1 = \langle z_\kk \mid \kk \text{ index} \rangle_\Q \quad \subset \quad  \HH=\Q\langle x,y\rangle\,.
\end{align*}

\subsection{Iterated integral expression for $\li$ and $\zeta$}
From now on we will restrict to real $|z|<1$ and consider integrals on the real axis. 
Since $\li_\kk(z)$ is defined for any $\kk$, we can view $\li$ as a $\Q$-linear map from $\HH^1$ to the space of real valued continuous functions on $(0,1)$, i.e.  $C((0,1);\R)$, defined on the generators by
\begin{align}\begin{split}\label{eq:limap}
\li\colon \HH^1 &\longrightarrow C((0,1);\R)\\
z_\kk &\longmapsto \li_{\kk}(z)\,.
\end{split}
\end{align}

By abuse of notation we write $\li\colon w \mapsto \li_w(z)$ for any $w \in \HH^1$, which is defined by linearly extending the definition on the generators $z_\kk$. For example for $w=xyxxy + 2xxxxy = z_2 z_3 + 2z_5 \in \HH^1$ we have $\li_w(z) = \li_{2,3}(z) + 2 \li_{5}(z)$. Now we want to describe the iterated integral expression for the multiple polylogarithm using this setup. 

\begin{lemma}\label{lem:lidif} Let $w \in \HH^1$ be a linear combination of words all starting with the letter $a \in \{x,y\}$, i.e. $w=au$ for some $u\in \HH^1$. Then we have 
\begin{align*}
\frac{d}{dz} \li_w(z) =  \frac{d}{dz} \li_{au}(z) = \begin{cases}
\frac{1}{z}\li_{u}(z) \,, &a = x\\
\frac{1}{1-z}\li_{u}(z) \,, &a = y
\end{cases}\,.
\end{align*}
\end{lemma}
\begin{proof} Since $\li$ is linear it suffices to prove the statement for a word $w$.	Assuming $w=z_\kk$ for $\kk=(k_1,\dots,k_r)$, we have
\begin{align*}
\frac{d}{dz} \li_w(z) &=  \frac{d}{dz} \li_\kk(z) = \frac{d}{dz}\sum_{m_1 > \cdots > m_r > 0} \frac{z^{m_1}}{m_1^{k_1}m_2^{k_2}  \cdots m_r^{k_r}}= \sum_{m_1 > \cdots > m_r > 0} \frac{z^{m_1-1}}{m_1^{k_1-1} m_2^{k_2} \cdots m_r^{k_r}}\,.
\end{align*}
Let $a=x$, which is equivalent to $k_1>1$. In this case we obtain 
\begin{align*}
\frac{d}{dz} \li_w(z) = \frac{d}{dz} \li_{xu}(z) = \frac{1}{z} \sum_{m_1 > \cdots > m_r > 0} \frac{z^{m_1}}{m_1^{k_1-1} m_2^{k_2} \cdots m_r^{k_r}} = \frac{1}{z}  \li_{u}(z)\,.
\end{align*}
If $a=y$, then we have $k_1=1$ and
\begin{align*}
\frac{d}{dz} \li_w(z) &= \frac{d}{dz} \li_{yu}(z) =  \sum_{m_1 > \cdots > m_r > 0} \frac{z^{m_1-1}}{ m_2^{k_2} \cdots m_r^{k_r}} = \sum_{m_2 > \cdots > m_r > 0} \frac{1}{ m_2^{k_2} \cdots m_r^{k_r}} \sum_{m_1=m_2+1}^\infty z^{m_1-1}\\
&= \frac{1}{1-z} \sum_{m_2 > \cdots > m_r > 0} \frac{z^{m_2}}{ m_2^{k_2} \cdots m_r^{k_r}}= \frac{1}{1-z} \li_{k_2,\dots,k_r}(z) = \frac{1}{1-z} \li_{u}(z)\,.
\end{align*}
\end{proof}

Motivated by the iterated integrals \eqref{eq:z2itint}, \eqref{eq:z32itint}, and the above Lemma, we define
\begin{align*}
\omega_x(t) = \frac{dt}{t}\,, \qquad \omega_y(t) = \frac{dt}{1-t}\,. 
\end{align*}
With these differential forms we  can write the multiple polylogarithms as the following iterated integral.

\begin{proposition} \label{prop:liitint}For any word $w=a_1\dots a_k \in \HH^1$, with $a_1,\dots,a_k \in \{x,y\}$ and  $0\leq z < 1$ we have
\begin{align*}
\li_w(z) = \int_0^z \omega_{a_1}(t_1) \int_0^{t_1} \omega_{a_2}(t_2) \dots \int_0^{t_{k-1}} \omega_{a_k}(t_k)\,. 
\end{align*}
\end{proposition}
\begin{proof}
This follows from Lemma \ref{lem:lidif} by induction on $k$. In the case $k=1$ we have $w= y = z_1$, i.e.  
\begin{align}
\li_w(z)  =     \li_{1}(z)= \sum_{m>0} \frac{z^m}{m} = \int_0^z \frac{dt}{1-t} = \int_0^z \omega_{y}(t)\,. 
\end{align}
The induction step is then exactly the statement of  Lemma \ref{lem:lidif} since $\li_w(0)=0$ for non-empty $w$.
\end{proof}
For a real $z$ we will also use the following simplified notation for iterated integrals for $a_1\dots a_k \in \HH^1$
\begin{align*}
\int\limits_{z>t_1>\dots>t_k>0} \omega_{a_1}(t_1) \cdots \omega_{a_k}(t_k) := \int_0^z \omega_{a_1}(t_1) \int_0^{t_1} \omega_{a_2}(t_2) \cdots \int_0^{t_{k-1}} \omega_{a_k}(t_k)\,.
\end{align*}
Since $\zeta(\kk)$ is just defined for admissible indices, we can, similar to \eqref{eq:limap}, define a $\Q$-linear map from $\HH^0$ to the space of multiple zeta values $\mz$, defined on the generators by
\begin{align}\begin{split}\label{eq:zetamap}
\zeta\colon \HH^0 &\longrightarrow \mz\\
z_\kk &\longmapsto \zeta(\kk)\,.
\end{split}
\end{align}
Also here we write $\zeta\colon w \mapsto \zeta(w)$ for any $w \in \HH^0$. 
Since $\li_w(1) = \zeta(w)$ for any $w \in \HH^0$ we also get an iterated integral expression for multiple zeta values as a consequence of Proposition \ref{prop:liitint}. 
\begin{corollary}\label{cor:mzvitint}For any word $w=a_1\dots a_k \in \HH^0$, with $a_1,\dots,a_k \in \{x,y\}$ we have
\begin{align*}
\zeta(w) = \int\limits_{1>t_1>\dots>t_k>0} \omega_{a_1}(t_1) \cdots \omega_{a_k}(t_k) \,. 
\end{align*}
\end{corollary}

\subsection{Duality relation}
We now give a direct consequence of the iterated integral expression. Making the change of variables $s_j = 1 - t_{k-j+1}$ in the iterated integral expression gives a linear relation among multiple zeta values, which is called the duality relation.  For example if $k=3$ we can make the change of variables $s_1 = 1- t_3$, $s_2 = 1-t_2$, $s_3=1-t_1$ in the following iterated integral
\begin{align}\begin{split} \label{eq:z321duality}
\zeta(3) =     \int_{0}^1 \frac{dt_1}{t_1} \int_{0}^{t_1}  \frac{dt_2}{t_2} \int_{0}^{t_2}  \frac{dt_3}{1-t_3}  
&=     \int_{1}^0 \frac{-ds_3}{1-s_3} \int_{1}^{s_3}  \frac{-ds_2}{1-s_2} \int_{1}^{s_2}  \frac{-ds_1}{s_1}  \\
&=      \int_{0}^1 \frac{ds_1}{s_1} \int_{0}^{s_1}  \frac{ds_2}{1-s_2} \int_{0}^{s_2}  \frac{ds_3}{1-s_3}   = \zeta(2,1)\,.
\end{split}
\end{align} 
from which we again get the relation $\zeta(3) = \zeta(2,1)$ in Proposition \ref{prop:z3z21}. This change of variables can be described nicely in terms of an anti-automorphism on the space  $\HH$.
For this we denote by $\tau$ the anti-automorphism of $\HH$ which interchanges $x$ and $y$. Here we view $\HH$, and all its subspaces, as $\Q$-algebras where the product is given by the usual non-commutative product in $\Q\langle x,y\rangle$. That $\tau$ is an anti-automorphism just means that $\tau(u w)= \tau(w)\tau(u)$ for $u,w \in \HH$ and $\tau(1)=1$. For example, if $w=z_3=xxy$, then 
\[ \tau(z_3)=\tau(xxy)=\tau(y)\tau(xx)=\tau(y)\tau(x)\tau(x)=xyy=z_2 z_1\,. \] 
Notice that $\tau(\HH^0) \subset \HH^0$, since any non-empty word $w\in \HH^0$ is of the form $w = xuy$ for some $u\in \HH$ and therefore $\tau(w)=\tau(xuy)=\tau(y)\tau(u)\tau(x)=x\tau(u)y \in \HH^0$. In addition, note that $\tau$ is an involution, that is, $\tau^2 = \operatorname{id}_\HH$ and $\tau(\HH^0) = \HH^0$.

\begin{proposition} [Duality relation] \label{prop:duality} For all $w\in \HH^0$ we have 
\begin{align*}    
\zeta(\tau(w)) = \zeta(w)\,.
\end{align*}
\end{proposition}
\begin{proof}
This is just a generalization of the variable change $s_j = 1 - t_{k-j+1}$ in the iterated integral expression in Corollary \ref{cor:mzvitint} similar to \eqref{eq:z321duality}. Interchanging $x$ and $y$ corresponds to $\omega_a(t_{k-j+1}) = -\omega_{\tau(a)}(s_j) $ for $a\in \{x,y\}$. The property of $\tau$ being an anti-automorphism corresponds to changing the order/directions of the integrals, which also gets rid of the minus signs. 
\end{proof}

A few explicit examples of the duality relations are given by the following Corollary, which both can be seen as a generalization of the formula $\zeta(3) = \zeta(2,1)$. Here we use the common notation  $\{ k_1,\dots,k_r \}^n=\underbrace{k_1,\dots,k_r,\dots,k_1,\dots,k_r}_{r n}$ for $n$ copies of the string $k_1,\dots,k_r$.
\begin{corollary} \begin{enumerate}[\textup{(}i\textup{)}] 
    \item For all $k\geq 3$ we have
    \begin{align*}
    \zeta(k) = \zeta(2,\underbrace{1,\dots,1}_{k-2}) = \zeta(2,\{1\}^{k-2})\,.
    \end{align*}
    \item For all $n\geq 1$ we have
    \begin{align*}
    \zeta(\{2,1\}^n) =     \zeta(\{3\}^n)\,. 
    \end{align*}
\end{enumerate}
\end{corollary}
\begin{proof} Both statements are immediate consequences of the duality relations, since $\tau(z_k) = \tau(x^{k-1}y)=xy^{k-1}=z_2 z_1\cdots z_1$ and $\tau( (z_2z_1)^n) = \tau(z_2z_1)^n = z_3^n$.
\end{proof}

\begin{remark} In Chapter \ref{sec:families} we will see another proof of the duality relation, which is not using the iterated integral expression. This new proof is based on so-called connected sums, which were just recently introduced by Seki and Yamamoto in \cite{SY}. There we will also see that the duality is true for the $q$-analogue model of Bradley-Zhao \eqref{eq:defqbz} and that we have $\zeta_q^{\rm{BZ}}(\tau(w)) = \zeta_q^{\rm{BZ}}(w)$, when considering $\zeta_q^{\rm{BZ}}$ as a map from $\HH^1$ to $\Q[[q]]$.
\end{remark}

\section{Shuffle, Stuffle, and Finite Double Shuffle Relations}

In this subsection, we will introduce the shuffle product $\sh$ and stuffle product $\ast$ on the spaces $\HH$, $\HH^1$ and $\HH^0$. We will then show that the space $\mz$ is a $\Q$-algebra (i.e. give a proof of Proposition \ref{prop:mzissubalgebra}) and see that the map $\zeta$ in \eqref{eq:zetamap} is an algebra homomorphism from $\HH^0$ to $\R$ with respect to both products $\sh$ and $\ast$. This will then lead to families of linear relations, which are called finite double shuffle relations. 

\subsection{The shuffle product}
The iterated integral expressions give another way to obtain the shuffle product formula 
\begin{align*}
\zeta(k_1) \zeta(k_2) = \sum_{j=2}^{k_1+k_2-1} \left( \binom{j-1}{k_1-1} + \binom{j-1}{k_2-1} \right) \zeta(j, k_1+k_2-j)
\end{align*}
in Proposition \ref{prop:shuffle1}, which was proved by using partial fraction decomposition. For example we have
\[\zeta(2)\zeta(3) =\zeta(2,3) + 3 \zeta(3,2)+ 6 \zeta(4,1)\,.\]

We will now describe how this relation can also be obtained from the iterated integral expression in Corollary \ref{cor:mzvitint}. Using the iterated integral expression of $\zeta(2)$ and $\zeta(3)$ we get 
\begin{align*}
\zeta(2)\zeta(3) &= \int\limits_{1>{\color{blue}t_1}> {\color{red}t_2} >0} \omega_{x}({\color{blue}t_1}) \omega_{y}({\color{red}t_2}) \int\limits_{1>{\color{blue}s_1}> {\color{blue}s_2}> {\color{red}s_3} >0} \omega_{x}({\color{blue}s_1}) \omega_{x}( {\color{blue}s_2}) \omega_{y}({\color{red}s_3})\,.
\end{align*}
By {\color{blue}blue} we indicate the variables which correspond to the differential form $\omega_x$ and by {\color{red}red} the ones corresponding to $\omega_y$. This makes it easier to translate the iterated integrals below back to multiple zeta values. 
Using Fubini's theorem the right-hand side is the iterated integral of  $\omega_{x}({\color{blue}t_1}) \omega_{y}({\color{red}t_2}) \omega_{x}({\color{blue}s_1}) \omega_{x}({\color{blue}s_2}) \omega_{y}({\color{red}s_3})$ over the domain where $1>{\color{blue}t_1}> {\color{red}t_2} >0$ and $1>{\color{blue}s_1}> {\color{blue}s_2}>{\color{red}s_3} >0$. This can be decomposed into the following iterated integrals, where we can neglect the non-trivial intersections $t_j=s_i$ since they have measure zero.

\begin{align*}
\zeta(2)\zeta(3) =&\,\Bigg(\, \int\limits_{1>{\color{blue}t_1}> {\color{red}t_2} > {\color{blue}s_1}>{\color{blue}s_2}>{\color{red}s_3}>0} +  \int\limits_{1>{\color{blue}t_1}> {\color{blue}s_1} > {\color{red}t_2}>{\color{blue}s_2}>{\color{red}s_3}>0} +  \int\limits_{1> {\color{blue}s_1}> {\color{blue}t_1} > {\color{red}t_2}>{\color{blue}s_2}>{\color{red}s_3}>0}+\int\limits_{1> {\color{blue}s_1}> {\color{blue}t_1} >{\color{blue}s_2}> {\color{red}t_2}>{\color{red}s_3}>0}\\
&+\int\limits_{1> {\color{blue}s_1}> {\color{blue}s_2} >{\color{blue}t_1}> {\color{red}t_2}>{\color{red}s_3}>0}+\int\limits_{1> {\color{blue}s_1}> {\color{blue}s_2} >{\color{blue}t_1}> {\color{red}s_3}>{\color{red}t_2}>0}+\int\limits_{1> {\color{blue}s_1}> {\color{blue}s_2} >{\color{red}s_3}>{\color{blue}t_1}> {\color{red}t_2}>0}+\int\limits_{1>{\color{blue}t_1}> {\color{blue}s_1} > {\color{blue}s_2}>{\color{red}t_2}>{\color{red}s_3}>0}\\
&+\int\limits_{1>{\color{blue}t_1}> {\color{blue}s_1} > {\color{blue}s_2}>{\color{red}s_3}>{\color{red}t_2}>0}+\int\limits_{1>{\color{blue}s_1}> {\color{blue}t_1} > {\color{blue}s_2}>{\color{red}s_3}>{\color{red}t_2}>0}\Bigg) \omega_{x}({\color{blue}t_1}) \omega_{y}({\color{red}t_2}) \omega_{x}({\color{blue}s_1}) \omega_{x}({\color{blue}s_2}) \omega_{y}({\color{red}s_3})\\
= \,&\, \zeta(2,3)+ \zeta(3,2)+\zeta(3,2)+\zeta(4,1) +\zeta(4,1)  +\zeta(4,1) + \zeta(3,2) +\zeta(4,1)+\zeta(4,1)+\zeta(4,1) \\
= \,&\,  \zeta(2,3) + 3 \zeta(3,2) + 6 \zeta(4,1)\,.
\end{align*}
The above example shows the origin of the name shuffle product, since one can interpret a multiple zeta value as a deck of blue and red cards, which correspond to the differential forms $\omega_x$ and $\omega_y$. Taking the product of two multiple zeta values then corresponds just to the shuffle of these two decks of cards.  We will now describe this product on the space $\HH$ and its subspaces. 

\begin{definition} We define the \emph{shuffle product}  $\sh$ on $\HH$ as the $\Q$-bilinear product, which satisfies $1 \shuffle w = w \shuffle 1 = w$ for any word $w\in \HH$ and
\begin{align*}
a_1 w_1 \shuffle a_2 w_2 = a_1 (w_1 \shuffle a_2 w_2) + a_2 (a_1 w_1 \shuffle w_2)
\end{align*}
for any letters $a_1, a_2 \in \{x,y\}$ and words $w_1,w_2 \in \HH$.
\end{definition}

By induction on the lengths of words, one can show that $\sh$ is a commutative and associative product and  $\HH_\sh = (\HH, \sh)$ is therefore a commutative $\Q$-algebra. One can also see that the subspaces $\HH^1$ and $\HH^0$ are both closed under $\sh$ and therefore we have subalgebras $\HH^0_\sh \subset \HH^1_\sh \subset \HH_\sh$.
You can check that this definition corresponds exactly to multiplying iterated integrals as above, i.e. we have 
\begin{align*}
z_2 \sh z_3 = xy \sh xxy = xyxxy + 3 xxyxy + 6 xxxyy = z_2 z_3 + 3 z_3 z_2 + 6 z_4 z_1 \,.
\end{align*}
That this is true in general will be proven now.

\begin{proposition} \label{prop:liisalghom} For any $w,u \in \HH^1$ we have
\begin{align*}
\li_w(z) \li_u(z) = \li_{w \sh u}(z)\,,
\end{align*}
i.e. the map $\li$ is an algebra homomorphism from $\HH^1_\sh$ to $C((0,1);\R)$.
\end{proposition}
\begin{proof} It is sufficient to prove the statement for words $w,u \in \HH^1$.    We will do this by induction on the sum of the lengths of  $w$ and $u$. If one of them equals the empty word $1$, the statement is clear. So let's assume that $w=a w'$ and $u=b u'$ for words $w',u' \in \HH^1$ and letters $a,b \in \{x,y\}$. Then we have 
\begin{align*}
\frac{d}{dz}\left( \li_w(z) \li_u(z)  \right) &=     \frac{d}{dz}\left( \li_{a w'}(z) \li_{b u'}(z)  \right) = \left( \frac{d}{dz} \li_{a w'}(z)\right) \li_{b u'}(z)  +  \li_{a w'}(z) \left(  \frac{d}{dz}\li_{b u'}(z)   \right)\,.
\end{align*}
Using now Lemma \ref{lem:lidif} we get $\frac{d}{dz} \li_{a w'}(z) = f_a(z)\li_{w'}(z) $ with $f_x(z)= \frac{1}{z}$ and $f_y(z)= \frac{1}{1-z}$. Using this together with the induction hypothesis we have
\begin{align*}
\frac{d}{dz}\left( \li_w(z) \li_u(z)  \right) &= f_a(z) \li_{w'}(z)\li_{b u'}(z) + f_b(z) \li_{aw'}(z)\li_{u'}(z) = f_a(z) \li_{w' \sh b u'}(z) + f_b(z)\li_{aw' \sh u'}(z)\,.
\end{align*}
Applying Lemma \ref{lem:lidif}  again gives
\[    
\frac{d}{dz}\left( \li_w(z) \li_u(z)  \right) =    \frac{d}{dz}\li_{a(w' \sh b u')}(z) +     \frac{d}{dz}\li_{b(aw' \sh u')}(z) =     \frac{d}{dz} \li_{w \sh u}(z)\,,
\]
i.e. $\li_w(z) \li_u(z) = \li_{w \sh u}(z)+c$ for some constant $c$. But since both sides vanish for $z=0$, we conclude $c=0$.
\end{proof}
For $w,u \in  \HH^0$ we can also set $z=1$ in the Proposition above and obtain the following.
\begin{corollary}\label{cor:zetashufflemap} For any $w,u \in \HH^0$ we have
\begin{align*}
\zeta(w) \zeta(u) = \zeta(w \sh u)\,.
\end{align*}
In particular, the space $\mz$ is a $\Q$-subalgebra of $\R$ and $\zeta$ is an algebra homomorphism from $\HH^0_\sh$ to $\mz$. 
\end{corollary}

\subsection{The stuffle product}
In Chapter \ref{sec:overview} we saw that for $k_1, k_2 \geq 2$ we have the stuffle product formula
\begin{align*}
\zeta(k_1)\zeta(k_2) = \left(\sum_{m_1 >m_2> 0} + \sum_{m_2 >m_1 > 0} + \sum_{m_1 = m_2 > 0}    \right) \frac{1}{m_1^{k_1} m_2^{k_2}} = \zeta(k_1,k_2) + \zeta(k_2,k_1) + \zeta(k_1+k_2)\,. 
\end{align*}
With the same argument, i.e. splitting up the summation, we get for $k_1, k_2\geq 2$, $k_3\geq 1$
\begin{align*}
\zeta(k_1)\zeta(k_2,k_3) = \zeta(k_1,k_2,k_3) +  \zeta(k_2,k_1,k_3) + \zeta(k_2,k_3,k_1)  +  \zeta(k_1+k_2,k_3)  +\zeta(k_2,k_1+k_3) \,.
\end{align*}
Similar to the shuffle product we will now define the stuffle product $\ast$ on the space $\HH^1$ and $\HH^0$ and then show that $\zeta$ is also an algebra homomorphism with respect to this product. Recall that $\HH^1 = \Q\langle z_1, z_2,\dots \rangle$, i.e. every element in $\HH^1$ can be viewed as a linear combination of words in the letters $z_j$ instead of the letters $x,y$. Here and in the following we will use the terminology `word' in these two different ways and it will be clear from context if we talk about words in the $z_j$ or $x,y$. 
In the next subsection, we will see that the stuffle and the shuffle product are both examples of quasi-shuffle products over different alphabets. 

\begin{definition} We define the \emph{stuffle product}  $\ast$ on $\HH^1$ as the $\Q$-bilinear product, which satisfies $1 \ast w = w \ast 1 = w$ for any word $w\in \HH^1$ and
\begin{align*}
z_{i} w_1 \ast z_{j} w_2 = z_i (w_1 \ast z_{j} w_2) + z_j (z_i w_1 \ast w_2) + z_{i+j} (w_1 \ast  w_2) 
\end{align*}
for any $i,j \geq 1$ and words $w_1,w_2 \in \HH^1$.
\end{definition}
Notice that this product also replicated the above product formula of multiple zeta values, since
\[ z_{k_1} \ast z_{k_2}   = z_{k_1} z_{k_2} + z_{k_2} z_{k_1} + z_{k_1+k_2}\,.\]
This product is also called the harmonic product and one can check (see \cite{H2}) that it is commutative and associative and therefore $\HH^1_\ast = (\HH^1, \ast)$ is a commutative $\Q$-algebra. By definition it is easy to check that $\HH^0$ is also closed under $\ast$ and we get a subalgebra $\HH^0_\ast \subset \HH^1_\ast$.

In the case of the shuffle product we used the polylogarithm to prove that the product of multiple zeta values satisfies the shuffle product formula by considering $z=1$. In the case of the stuffle product, we will consider the  \emph{truncated multiple zeta values}\footnote{We earlier called them multiple harmonic sums $H_M$ and the notation/naming is not agreed on in the literature. In this section, we will use the notion truncated multiple zeta values $\zeta_M$. Notice that $H_M = \zeta_{M+1}$ as we allowed $M=m_1$ in the definition of $H_M$.}, which are for an integer $M\geq 1$ and any index $\kk = (k_1,\dots, k_r) \in \Z_{\geq 1}^r$ defined by
\begin{align*}
\zeta_M(\kk) = \zeta_M(k_1,\dots,k_r) = \sum_{M>m_1>\dots>m_r>0} \frac{1}{m_1^{k_1}\cdots m_r^{k_r}} \in \Q \,.
\end{align*}
Clearly if $\kk$ is admissible we have $\lim_{M\rightarrow \infty} \zeta_M(\kk)  = \zeta(\kk)$. For a fixed $M$ we can view $\zeta_M$ as a $\Q$-linear map from $\HH^1$ to $\Q$, defined on the generators by $\zeta_M\colon z_\kk \mapsto \zeta_M(\kk)$.

\begin{proposition} \label{prop:zmstuffle}For any $w,u \in \HH^1$ and $M \geq 1$ we have
\begin{align*}
\zeta_M(w) \zeta_M(u) = \zeta_M(w \ast u)\,,
\end{align*}
i.e. the map $\zeta_M$ is an algebra homomorphism from $\HH^1_\ast$ to $\Q$.
\end{proposition}
\begin{proof}
This can be done by induction on the depths or $M$ (Exercise \ref{ex6} (ii)). We will prove this in more general form for quasi-shuffle algebras in the next section (Lemma \ref{lem:fmtoFM}).
\end{proof}

For $w,u \in  \HH^0$ we can also take the limit $M\rightarrow \infty$ in the Proposition above and obtain the following.
\begin{corollary}\label{cor:zetastufflemap} For any $w,u \in \HH^0$ we have
\begin{align*}
\zeta(w) \zeta(u) = \zeta(w \ast u)\,,
\end{align*}
i.e. $\zeta$ is an algebra homomorphism from $\HH^0_\ast$ to $\mz$. 
\end{corollary}

\subsection{Finite double shuffle relations}
Since the map $\zeta\colon \HH^0 \rightarrow \mz$ is an algebra homomorphism with respect to the shuffle product $\sh$ and the stuffle product $\ast$, we get a large family of linear relations among multiple zeta values.

\begin{proposition}[Finite double shuffle relations]\label{prop:fds} For $w,u \in \HH^0$ we have
\begin{align*}
\zeta(w \sh u - w\ast u) = 0\,.
\end{align*}
\end{proposition}
But it is also clear, that these do not give all linear relations among multiple zeta values. For example the relation $\zeta(3)=\zeta(2,1)$ is not a consequence of the above Proposition. Counting the finite double shuffle relations, we get the following table, which comes from the survey article \cite{Tan2}. In this article, you can also find the numbers of other families of relations, such as the duality relation.
\begin{center}
\begin{tabular}{|c|c|c|c|c|c|c|c|c|c|c|c|c|c|}
    \hline
    weight $k$         & 3 & 4 & 5 & 6 & 7 & 8 & 9 & 10 & 11 & 12  \\ \hline
    \# all conjectured relations       &1& 3& 6& 14& 29& 60& 123& 249& 503& 1012\\ \hline
    \# finite double shuffle relations  & 0 & 1 & 2 & 7 & 16 & 40 & 92 & 200  & 429  & 902  \\ \hline
\end{tabular}
\end{center}

We see that the first possible finite double shuffle relation appears in weight $4$ by choosing $w=u=z_2$, which gives
\begin{align*}
w \sh u - w\ast u = (2 z_2 z_2 + 4 z_3 z_1) - (2 z_2z_2 + z_4) = 4 z_3 z_1 - z_4
\end{align*}
i.e. $4\zeta(3,1)= \zeta(4)$. This relation is a special case of the following family of linear relations which is a consequence of finite double shuffle relations.

\begin{proposition}\label{prop:z31z44} For all $n\geq 1$ we have 
\begin{align*}
4^n \zeta(\{3,1\}^n) = \zeta(\{4\}^n)\,.
\end{align*}
\end{proposition}
\begin{proof}
This can be done by proving the following equations in $\HH^0$ (Exercise \ref{ex7})
\begin{align*}
\sum_{j=-n}^n (-1)^j z_2^{n-j} \sh z_2^{n+j} = 4^n (z_3 z_1)^n\,,\qquad 
\sum_{j=-n}^n (-1)^j z_2^{n-j} \ast z_2^{n+j} = z_4^n\,.
\end{align*}
(Notice: Here $z_k^n$ means $z_k z_k \dots z_k$, i.e. the usual non-commutative product in $\HH$ and \underline{not} the shuffle or stuffle product.)
\end{proof}
Together with the explicit formula (Exercise \ref{ex7})
\begin{align}\label{eq:222}
\zeta(\{2\}^n) = \frac{\pi^{2n}}{(2n+1)!} 
\end{align}
the proof of Proposition \ref{prop:z31z44} can also be used to show 
\begin{align*}
\zeta(\{4\}^n) = \frac{4^n 2 \pi^{4n}}{(4n+2)!}\,,\qquad \zeta(\{3,1\}^n) = \frac{2 \pi^{4n}}{(4n+2)!}\,,
\end{align*}
where the second equation here is known as the $3$-$1$ formula for multiple zeta values.

\newpage 

\section{Quasi-shuffle algebras}\label{subsec:quasishuffle}

In this section, we want to generalize what we did in the previous section. On $\HH$, we defined the shuffle product on words in letters $x$ and $y$, and on $\HH^1$, we defined the stuffle product on words in letters $z_j$ for $j\geq 1$. This idea will be generalized now by allowing an arbitrary set of letters $A$ and then defining a product on the space of words in these letters. We will mainly follow the definitions and theorems in \cite{HI} and \cite{IKZ}, but will also introduce some definitions and theorems which cannot be found in the literature.

\subsection{The algebra of letters and words} 
In the following, we assume that $\ko$ is a field containing $\Q$, and $A$ is a countable set which we refer to as the \emph{set of letters}. Let $\ko A$ be the  $\ko$-vector space generated by $A$ and let $\diamond$ be a $\ko$-bilinear, associative  and commutative  product on $\ko A$. We obtain a (non-unital) $\ko$-algebra $(\ko A, \diamond)$, which we refer to as the \emph{algebra of letters}. Notice that in the work \cite{HI} the authors just assume that $\diamond$ is associative and commutative and do not consider $\ko A$ as an algebra, but it will be useful for our purposes (e.g., Lemma \ref{lem:fmtoFM}). For such a product $\diamond$ on letters, we want to assign a product $\ast_\diamond$ on the space of words $\QA$, which generalizes the stuffle and shuffle product we have seen before. Here a monomial $w = a_1 \cdots a_l$ in $\QA$ will again be called a \emph{word} and the unit $(l=0)$, denoted by  $1$, is again called the \emph{empty word}. By $\len(w)=l$ we denote the \emph{length of the word} $w$.

\begin{definition} Let $\diamond$ be a product on $\ko A$ as above. Then we define the \emph{quasi-shuffle product}  $\qsh$ on $\QA$ as the $\ko$-bilinear product, which satisfies $1 \qsh w = w \qsh 1 = w$ for any word $w\in \QA$ and
\begin{align}\label{eq:qshdef}
a w \qsh b v = a (w \qsh b v) + b (a w \qsh v) + (a \diamond b) (w \qsh  v) 
\end{align}
for any letters $a,b \in A$ and words $w, v \in \QA$.
\end{definition}

\begin{theorem} \label{thm:quasishufflealgebra}The space $\QA$ equipped with the product $\qsh$ becomes a commutative $\ko$-algebra.
\end{theorem}
\begin{proof} This is Theorem 2.1 in \cite{HI}. It suffices to show that $\qsh$ is commutative and associative, which can be done straightforward by induction on the lengths of words. 
\end{proof}

We will call $(\QA,\qsh)$ a \emph{quasi-shuffle algebra} or \emph{algebra of words}. Notice that this generalizes the $\Q$-algebra $\HH_\sh$ by choosing $\ko = \Q$, $A=\{x,y\}$ and $a \diamond b=0$ for any  letters $\{x,y\}$ and it generalizes the $\Q$-algebra $\HH^1_\ast$, by choosing $\ko = \Q$, $A=\{z_1,z_2,\dots\}$ and $z_i \diamond z_j=z_{i+j}$ (See section \ref{subsubsec:qshex} for details). Our purpose to introduce quasi-shuffle products is to also describe the product structure of the $q$-series $\g(\kk)$ and the (modified) $q$-analogues $\zeta_q$.

\begin{lemma} \label{lem:fmtoFM}Let $R$ be a $\ko$-algebra and $f_m \colon (\ko A, \diamond) \rightarrow R$ be $\ko$-algebra homomorphisms for $m\geq 1$. Then for all $M\geq 1$ the $\ko$-linear map $F_M : \QA \rightarrow R$ defined on a word $w=a_1 \cdots a_r \in \QA$ by
\begin{align*}
F_M(w) = \sum_{M > m_1 > \dots > m_r> 0 } f_{m_1}(a_1) \cdots f_{m_r}(a_r)\,
\end{align*}
and $F_M(1)=1$ is a $\ko$-algebra homomorphism from $(\QA,\qsh)$ to $R$. 
\end{lemma}

\begin{proof} It suffices to show that for any $M\geq1$ and words $w,v \in \QA$ we have 
\begin{align*}
F_M(w) F_M(v) = F_M(w \qsh v)\,.
\end{align*}
We will prove this by induction on $M$. The case $M=1$ is trivial, since $F_1(w)=0$ for all non-empty $w$ and $1=F_1(1)F_1(1) = F_1(1\qsh 1)$. Notice that we have $F_M(aw) = \sum_{M>m>0} f_m(a) F_m(w)$ for a letter $a$ and a word $w$. For $w=aw'$, $v=b v'$ with letters $a,b \in A$ and words $w',v' \in \QA$ we therefore get
\begin{align*}
&F_M(w) F_M(v) = \sum_{M > m > 0 } f_{m}(a) F_m(w') \sum_{M > n > 0 } f_{n}(b) F_n(v')\\
&= \left( \sum_{M > m>n > 0 } +\sum_{M > n>m > 0 } + \sum_{M > m=n > 0 }   \right) f_{m}(a) F_m(w') f_{n}(b) F_n(v') \\
&= \sum_{M > m > 0 } f_{m}(a) F_m(w')  F_m(bv') + \sum_{M > n > 0 } f_{n}(b) F_n(aw') F_n(v')   + \sum_{M > m > 0 } f_{m}(a)  f_{m}(b)  F_m(w') F_m(v') \\
&= \sum_{M > m > 0 } f_{m}(a) F_m(w' \qsh bv') + \sum_{M > n > 0 } f_{n}(b) F_n(aw' \qsh v')   + \sum_{M > m > 0 } f_{m}(a \diamond b)  F_m(w' \qsh v') \\
&= F_M(a(w' \qsh bv'))+F_M(b(aw' \qsh v') ) + F_M((a \diamond b) (w' \qsh v')) = F_M(w \qsh v)\,.
\end{align*}
Here we used that $f_m$ is an algebra homomorphism together with the induction hypothesis in the fourth equation.
\end{proof}

\subsection{Examples of quasi-shuffle products \& (sub)algebras}\label{subsubsec:qshex}
In the following, we give a few explicit examples for quasi-shuffle products and algebras which appear in these notes.  We will also consider certain subalgebras, and the first statement we want to show now is that  subalgebras of the algebra of letters give subalgebras of the algebra of words.

\begin{proposition}\label{prop:lettersubgiveswordsub}
If $B \subset A$ is a subset of letters such that $(\ko B, \diamond)$ is a subalgebra of $(\ko A,\diamond)$, then $(\ko\langle B\rangle,\qsh)$ is a subalgebra of $(\ko\langle A\rangle,\qsh)$.
\end{proposition}
\begin{proof}
We need to show that $\ko\langle B\rangle$ is closed under $\qsh$. This follows by induction on the sum of the lengths of two words. Namely, if $a,b\in B$, then $a\diamond b\in\ko B$ and the right-hand side of \eqref{eq:qshdef} belongs to $\ko\langle B\rangle$ by bilinearity and the induction hypothesis.
\end{proof}

Most of our objects depend on some index $\kk = (k_1,\dots,k_r)$ and therefore most of our examples use the set of letters, the ``$z$-alphabet'', defined by 
\begin{align*}
A_z :=\{z_1,z_2,\dots\}\,. 
\end{align*}
Notice that we have an abuse of notation here, since $z_k=x^{k-1}y$ denoted elements in $\HH$ previously. But from the context it should always be clear if we talk about the formal elements $z_k$ in $A_z$ or the elements in $\HH$.  

\begin{enumerate} [\textup{(}i\textup{)}] 
\item  {\bf Shuffle product} For any field $\ko$ and any set of letters $A$, one can define the trivial product $a \diamond b=0$ for $a,b \in A$. The resulting quasi-shuffle product is then just the shuffle product $\sh = \qsh$. As a special case we considered $\ko = \Q$, $A=\{x,y\}$ before but will also deal with the shuffle product on $A_z$ later, which is sometimes also called the ``index-shuffle product''. 
\item {\bf Stuffle product} Another example we considered before is the stuffle product. Choosing $\ko = \Q$, $A_z =\{z_1,z_2,\dots\}$ and $z_i \diamond z_j=z_{i+j}$ we write $\ast = \qsh$. Here the $z_j$ are considered as variables themselves, but we see that $(\Q \langle A_z \rangle,\ast)$ is isomorphic to $\HH^1_\ast$ as a $\Q$-algebra, when sending $z_k$ to $x^{k-1}y$. Notice that the algebra of letters $(\Q A_z,\diamond)$ is isomorphic to $X\Q[X]$, by sending $z_k$ to $X^k$. Defining for $m\geq 1$ the algebra homomorphisms $f_m: X\Q[X] \rightarrow \Q$ by $p \mapsto p\left(m^{-1}\right)$ gives the truncated multiple zeta values $\zeta_M$ as the function $F_M$ in Lemma  \ref{lem:fmtoFM}. In particular, we obtain Proposition \ref{prop:zmstuffle} as a consequence. Considering the subsets $A_z^{\geq 2} \subset A_z$ and $A_z^{\text{ev}} \subset A_z$ defined by 
\begin{align*}
A_z^{\geq 2} :=A_z \backslash \{z_1\} = \{z_2,z_3,z_4,\dots\}\,,\qquad	A_z^{\text{ev}} :=\{z_2,z_4,z_6,\dots\}\,,
\end{align*}
we clearly have that $(\Q A_z^{\geq 2}, \diamond)$ and $(\Q A_z^{\text{ev}}, \diamond)$ are subalgebras of $(\Q A_z, \diamond)$. As a consequence of Proposition \ref{prop:lettersubgiveswordsub} we see that 
\begin{align*}
\mz^{\geq 2} &= \Big \langle   \zeta(k_1,\dots,k_r) \,\,\big|\,\, r\ge 0 ,\,k_1,\dots,k_r\geq 2 
\Big\rangle_\Q \,,\\
\mz^{\text{ev}} &= \Big \langle   \zeta(k_1,\dots,k_r) \,\,\big|\,\, r\ge 0 ,\,k_1,\dots,k_r\geq 2  \text{ even} 
\Big\rangle_\Q
\end{align*}
are subalgebras of $\mz$. Notice that by Brown's Theorem \ref{thm:brown} we actually have $\mz^{\geq 2}  = \mz$.

\item {\bf The $q$-series $\g(\kk)$:} Recall that we defined for an index $\kk = (k_1,\dots,k_r)$ the modified $q$-analogues 
\begin{align*}
\g(\kk)=\g(k_1,\dots,k_r) =  \sum_{m_1 > \cdots > m_r > 0} \frac{P_{k_1}(q^{m_1})}{(1-q^{m_1})^{k_1}} \cdots \frac{P_{k_r}(q^{m_r})}{(1-q^{m_r})^{k_r}} \in \Q[[q]]\,.
\end{align*}
Inspired by Lemma \ref{lem:rproduct} we define on $\Q A_z$ the product
\begin{align}\label{eq:defgdiamond}
z_{k_1 }\gdiamond z_{k_2}&=     z_{k_1+k_2}+ \sum_{j=1}^{k_1+k_2-1}\left( \lambda^j_{k_1,k_2}+ \lambda^j_{k_2,k_1} \right) z_j\,
\end{align}
where the rational numbers $\lambda^j_{k_1,k_2}$ are given by
\begin{align*}
\lambda^j_{k_1,k_2}  = (-1)^{k_2-1} \binom{k_1+k_2-1-j}{k_1-j}  \frac{B_{k_1+k_2-j}}{(k_1+k_2-j)!} \,.
\end{align*}
It is easy to see that this product is indeed associative and commutative. 
By Lemma  \ref{lem:rproduct} we see that for $m\geq 1$ the map $f_m : \Q A_z \rightarrow \Q[[q]]$ defined on the generators by 
\[f_m(z_k) = \frac{P_{k}(q^m)}{(1-q^m)^{k}}\]
is a $\Q$-algebra homomorphism from $(\Q A_z ,\gdiamond)$ to $\Q[[q]]$. We will denote the corresponding quasi-shuffle product by $\gqsh = \ast_{\gdiamond}$. 
 Using Lemma  \ref{lem:fmtoFM}, we see, after taking the limit $M\rightarrow \infty$, that the space $\gs$, spanned by all $\g(\kk)$, is a $\Q$-subalgebra of $\Q[[q]]$ and we can view $\g$ as an algebra homomorphism from $(\Q \langle A_z \rangle, \gqsh)$ to $\gs$. This proves Proposition \ref{prop:gisqalg}.
\begin{proposition} \label{prop:gevensubspace}The subspaces $\gsev \subset \gstwo \subset \gs$, defined by
    \begin{align*}
    \gstwo &= \Big \langle   \g(k_1,\dots,k_r) \,\,\big|\,\, r\ge 0 ,\,k_1,\dots,k_r\geq 2 
    \Big\rangle_\Q \,,\\
    \gsev &= \Big \langle   \g(k_1,\dots,k_r) \,\,\big|\,\, r\ge 0 ,\,k_1,\dots,k_r\geq 2  \text{ even} 
    \Big\rangle_\Q\,,
    \end{align*}        
    are $\Q$-subalgebras of $\gs$. 
\end{proposition}
\begin{proof}
    This will be a consequence of Proposition \ref{prop:lettersubgiveswordsub} after showing that  $(\Q A_z^{\geq 2}, \gdiamond)$ and $(\Q A_z^{\text{ev}}, \gdiamond)$ are subalgebras of $(\Q A_z, \gdiamond)$. Assume $z_{k_1}, z_{k_2} \in A_z^{\geq 2}$, i.e. $k_1,k_2 \geq 2$, then we see that all elements on the right-hand side in \eqref{eq:defgdiamond} are in $\Q A_z^{\geq 2}$, since $\lambda^1_{k_1,k_2}+\lambda^1_{k_2,k_1} = 0$ in these cases. This can be seen by writing 
    \begin{align}\label{eq:l1l1}
    \lambda^1_{k_1,k_2}+\lambda^1_{k_2,k_1} = \left((-1)^{k_1-1} + (-1)^{k_2-1}\right)\binom{k_1+k_2-2}{k_1-1} \frac{B_{k_1+k_2-1}}{(k_1+k_2-1)!}\,.
    \end{align}
    If $k_1$ and $k_2$ have different parity then this is clearly zero and if $k_1+k_2\geq 4$ is even, then $k_1+k_2-1\geq 3$ is odd and therefore $B_{k_1+k_2-1}$ vanishes. This shows that $(\Q A_z^{\geq 2}, \gdiamond)$ is a subalgebra of $(\Q A_z, \gdiamond)$.
    
    Now assume that $z_{k_1}, z_{k_2} \in A_z^{\text{ev}}$, i.e. $k_1,k_2$ are even. Since $B_k=0$ for odd $k> 1$ we have that $\lambda^j_{k_1,k_2} = \lambda^j_{k_2,k_1}=0$ in the cases that $j$ is odd and $1 \leq j < k_1+k_2-1$.  But since $k_1,k_2\geq 2$ we also have that $\lambda^{k_1+k_2-1}_{k_1,k_2} = \lambda^{k_1+k_2-1}_{k_2,k_1}=0$. This shows that  all elements on the right-hand side in \eqref{eq:defgdiamond} are in $\Q A_z^{\text{ev}}$ from which we get that $(\Q A_z^{\text{ev}}, \gdiamond)$ is a subalgebra of $(\Q A_z, \gdiamond)$.
\end{proof}

\item {\bf Bradley-Zhao q-MZV:} We define for an admissible index $\kk=(k_1,\dots,k_r)$
\begin{align}\label{eq:defqbz}
\zeta_q^{\rm{BZ}}(\kk)= \zeta_q^{\rm{BZ}}(k_1,\dots,k_r) = \sum_{m_1 > \cdots > m_r > 0} \frac{q^{(k_1-1)m_1} \cdots q^{(k_r-1)m_r}}{[m_1]_q^{k_1} \cdots [m_r]_q^{k_r} } \,.
\end{align}
Since we have for $m\geq 1$ and $k_1,k_2 \geq 2$
\begin{align*}
\frac{q^{(k_1-1)m}}{[m]_q^{k_1}} \frac{q^{(k_2-1)m}}{[m]_q^{k_2}} = \frac{q^{(k_1+k_2-1)m}}{[m]_q^{k_1+k_2}} + (1-q) \frac{q^{(k_1+k_2-2)m}}{[m]_q^{k_1+k_2-1}} 
\end{align*}
we choose $\ko = \Q(1-q)$ and define on $\ko A_z$ the product 
\begin{align*}
z_{k_1 }\diamond z_{k_2}&=     z_{k_1+k_2}+ (1-q) z_{k_1+k_2-1}.
\end{align*}
As before we get a quasi-shuffle algebra $(\ko \langle A_z \rangle,\qsh)$ and for an $M\geq 1$ an algebra homomorphism $F_M$ to $\Q[[q]]$ by sending $z_{k_1}\dots z_{k_r}$ to the truncated version $\zeta_{q,M}^{\rm{BZ}}(k_1,\dots,k_r) $ defined in the obvious way.  Notice that one could also consider the modified version $\overline{\zeta}_q^{\rm{BZ}}(\kk) := (1-q)^{-\wt(\kk)}\zeta_q^{\rm{BZ}}(\kk)$. With this one can choose again $\ko=\Q$ and use the product $z_{k_1 }\diamond z_{k_2}=     z_{k_1+k_2}+  z_{k_1+k_2-1}$. The algebraic structure of these modified versions is for example studied in  \cite{Tak1}.
\item {\bf Generalized modified q-MZV:} Take $A_{z,Q}= \{ z^Q_k \mid k\geq 1, Q\in \Q[X], \deg(Q) \leq k \}$ and define 
\begin{align*}
z^{Q_1}_{k_1} \diamond z^{Q_2}_{k_2}  = z^{Q_1 \cdot Q_2}_{k_1+ k_2}\,. 
\end{align*}
Then clearly $f_m: (\Q A_z ,\diamond) \rightarrow \Q[[q]]$ given by $f_m(z^{Q}_{k})= \frac{Q(q^m)}{(1-q^m)^k}$ is a $\Q$-algebra homomorphism. Again as a consequence of  Lemma  \ref{lem:rproduct} we see that  the space $\mz_q $
is a $\Q$-subalgebra of $\Q[[q]]$.

\end{enumerate}

\subsection{Words of repeating letters}
In the following  subsection, we want to state some standard facts on quasi-shuffle algebras, which were established in \cite{HI}, \cite{H3}, and \cite{IKZ} (without using the notion of quasi-shuffle algebras). Some of these will be given without proofs and we refer the reader to the above references for details. \\

Let $f = \sum_{n=0}^\infty c_n T^n \in \ko[[T]]$ and $\bullet \in \{\diamond, \qsh \}$, then we define for $a \in A$
\begin{align*}
f_\bullet(aX) = \sum_{n=0}^\infty c_n \underbrace{a \bullet \cdots \bullet a}_{n} X^n = \sum_{n=0}^\infty c_n a^{\bullet n} X^n  \in \QA[[X]] \,.
\end{align*}
In other words $f_\bullet(aX)$ means that we plug $aX$ into the power series $f$, and then use the product $\bullet$ to evaluate the products of $a$ in $(aX)^n$. Therefore it also makes sense to consider $f_\bullet(zX)$ for any $z \in \ko A[[X]]$ and then evaluate $(zX)^n$ as an element in $(\QA,\qsh)[[X]]$ or $(\ko A,\diamond)[[X]]$. This we will do in  Proposition \ref{prop:expqshlogdiamond} below. Notice that we have two different products on $\QA$ given by the quasi-shuffle product $\qsh$ and given by the usual non-commutative multiplication. When calculating with elements in $\QA[[X]]$, we will always do it with respect to the usual non-commutative multiplication in $\QA$.
\begin{proposition} \label{prop:expqshlogdiamond} For all $z\in \ko A[[X]]$ we have
\begin{align}\label{eq:explogcor}
\exp_{\qsh}( \log_\diamond(1+zX)) = \frac{1}{1-zX}\,.
\end{align}
\end{proposition}
\begin{proof} This is \cite[Corollary 5.1]{HI}, but it will also be a consequence of Proposition \ref{prop:psiffdiamond} below.
\end{proof}

This proposition can be restated in more explicit terms as follows:

\begin{proposition}
    Assume that we have an algebra homomorphism $\varphi: (\QA, \qsh) \rightarrow R$ in some $\ko$-algebra $R$. Applying $\varphi$ to \eqref{eq:explogcor} with $z=a \in A$ gives the following equation in $R[[X]]$

\begin{align}\label{eq:expphi}
\exp\left( \sum_{n=1}^{\infty} (-1)^{n-1} \varphi( a^{\diamond n}) \frac{X^n}{n} \right) = 1+ \sum_{n=1}^{\infty} \varphi(a^n) X^n.
\end{align}

In particular, we have 
\begin{align}\label{eq:expphipolynom}
\varphi(a^n) \in \ko\left[  \varphi( a^{\diamond j}) \mid 1 \leq j \leq n \right]\,,
\end{align}
i.e. for $a\in A$ the $\varphi(a^n)=\varphi(aa\cdots a)  $ is a polynomial in $\varphi( a^{\diamond j})$ with $1 \leq j \leq n$. 
\end{proposition}

\begin{corollary}\label{cor:zkkk}
\begin{enumerate}[\textup{(}i\textup{)}] 
    \item For all $k, M\geq 1$ we have 
    \begin{align*}
    \exp\left( \sum_{n=1}^{\infty} (-1)^{n-1} \zeta_M(nk) \frac{X^n}{n} \right) = 1+ \sum_{n=1}^{\infty} \zeta_M(\{k\}^n) X^n
    \end{align*}
    and therefore $\zeta_M(\{k\}^n) \in \Q[\zeta_M(j k) \mid 1\leq j \leq n]$ for all $n\geq 1$. In particular, for $k\geq 2$ these statements also hold by replacing $\zeta_M$ with $\zeta$.
    \item For all $k, n\geq 1$ we have 
    \begin{align*}
    \g(\{k\}^n)   \in \Q[\g(j) \mid 1\leq j \leq kn]\,.
    \end{align*}
    In addition if $k\geq 2$ is even, then $\g(\{k\}^n)   \in \Q[\g(j) \mid 2\leq j \leq kn\,, j \text{ even}]$
\end{enumerate}
\end{corollary}
\begin{proof}
Statement (i) follows directly from \eqref{eq:expphi} by using the algebra homomorphism $\zeta_M: (\Q\langle A_z\rangle,\ast) \rightarrow \Q$, $a=z_k$ and the  fact that $z_k^{\diamond n} = z_{k n}$ if $z_i \diamond z_j=z_{i+j}$. For (ii) we use the algebra homomorphism $\g: (\Q\langle A_z\rangle,\gqsh) \rightarrow \Q[[q]]$ together with \eqref{eq:expphipolynom}. The second part of (ii) follows since $\Q A_z^{\text{ev}}$ is closed under $\gdiamond$ as we saw in Proposition \ref{prop:gevensubspace} and therefore $z_k^{\gdiamond n} \in \Q A_z^{\text{ev}}$ if $k$ is even.
\end{proof}

\subsection{Application: Quasi-modular forms}
By Proposition \ref{prop:euler} we know that $\zeta(2m) \in \Q \pi^{2m}$ for $m\geq1$ and with Corollary \ref{cor:zkkk} we get
\begin{align*}
\zeta(2m,\dots,2m) \in \Q[\zeta(2)] = \Q[\pi^2]\,.
\end{align*} 
For example, as we have already seen before, we have for all $n\geq 1$
\begin{align*}
\zeta(\{2\}^n) = \frac{\pi^{2n}}{(2n+1)!} \,,\qquad 
\zeta(\{4\}^n) = \frac{4^n 2 \pi^{4n}}{(4n+2)!}\,.
\end{align*}

A similar statement can also be shown for the $q$-series $\g$, and the $q$-analogue of $\Q[\zeta(2)]$ is given by the \emph{ring of quasi-modular forms} (with rational coefficients), defined by
\begin{align*}
\qmf := \Q[\g(2), \g(4),\g(6)]\,.
\end{align*}

As an analogue of Euler's formula for $\zeta(2m)$ we get the following.
\begin{theorem} \label{thm:eulerforg}For all $m\geq 1$ we have $\g(2m) \in \qmf$.
\end{theorem}
\begin{proof}
For $m=1$ this follows from the definition of $\qmf$. Let $m\geq2$.
By Proposition \ref{prop:euler}, \eqref{eq:ek} and \eqref{eq:defgk} we have
\begin{align*}
G_{2m}=\beta(2m)E_{2m}.
\end{align*}
Since $E_{2m}$ is a modular form with rational Fourier coefficients,
Theorem~\ref{thm:mf}~(v) shows that
\begin{align*}
E_{2m}\in\mf^\Q=\Q[E_4,E_6].
\end{align*}
This also follows from \cite[Proposition~4.1~(i)]{BK1}.
Notice that
\begin{align*}
G_4=\beta(4)E_4=\beta(4)+\g(4),\qquad
G_6=\beta(6)E_6=\beta(6)+\g(6).
\end{align*}
Therefore,
\begin{align*}
G_{2m}\in\Q[G_4,G_6]
=\Q[\g(4),\g(6)]\subset\qmf.
\end{align*}
Since $\beta(2m)\in\Q$, the identity
$\g(2m)=G_{2m}-\beta(2m)$ gives the claim.
\end{proof}

\begin{corollary} For all $m\geq 1$ we have
\begin{align*}
\g(2m,\dots,2m) \in \qmf\,.
\end{align*}
\end{corollary}
\begin{proof}
This now follows directly from Theorem  \ref{thm:eulerforg} and Corollary \ref{cor:zkkk} (ii).
\end{proof}

\begin{proposition}\label{prop:qmprop}
\begin{enumerate}[\textup{(}i\textup{)}]  \item The space $\qmf$ is closed under $q\frac{d}{dq}$.
    \item We have 
    \begin{align*}
    \qmf = \Q[\g(2),\g'(2), \g''(2) ]\,.
    \end{align*}
    where $\g'$ denotes the derivative with respect to $q\frac{d}{dq}$.
    \item Let $\kk=(k_1,\dots,k_r)$ be an index with $k_1,\dots,k_r \geq 2$ even. Then we have 
    \begin{align*}
    \g^{\text{sym}}(\kk) := \sum_{\sigma \in S_r} \g(k_{\sigma(1)},\dots,k_{\sigma(r)}) \in \qmf\,,
    \end{align*}
    where $S_r$ denotes the set of all permutations of $\{1,\dots,r\}$.
    \item We have 
    \begin{align*}
    \qmf = \Big \langle  \g^{\text{sym}}(k_1,\dots,k_r)  \,\,\big|\,\, r\ge 0 ,\,k_1\geq k_2\geq\dots \geq k_r\geq 2 \text{ even}
    \Big\rangle_\Q\,,
    \end{align*} 
    where we set $    \g^{\text{sym}}(\emptyset)=1$.
\end{enumerate}
\end{proposition}
\begin{proof}
    This is Exercise \ref{ex8}.
    The first statement can be proven with  Proposition \ref{prop:gstuffle1} and \ref{prop:gshuffle1}  by giving explicit formulas for $q\frac{d}{dq} \g(2), q\frac{d}{dq} \g(4)$ and $q\frac{d}{dq} \g(6)$ as polynomials in $\g(2), \g(4)$ and $\g(6)$. From this one also deduces (ii). Statement (iii) and (iv) can be proven by induction on $r$ and the weight respectively.
\end{proof}

\subsection{Linear maps induced by power series}
In this section, we want to illustrate an important tool for quasi-shuffle algebras, which was first established in \cite{H3}
and later generalized in \cite[Section 3]{HI}. In addition, there is a recent work of Yamamoto \cite{Yam}, which generalizes this construction even more and gives a nice reinterpretation of some of the results we will mention here. 
One motivation for studying these maps is the following. For a given quasi-shuffle algebra $(\QA, \qsh)$, one can always construct an explicit isomorphism (of $\ko$-algebras) to $(\QA, \sh)$. In other words, all quasi-shuffle algebras over the same alphabet are isomorphic. 

We will first illustrate this basic idea on multiple zeta values, which are the images of an algebra homomorphism from $(\Q A_z, \ast)$ to $\R$. For now, we ignore convergence issues, since everything we are going to do could also be done for the truncated version $\zeta_M$. 
Recall that stuffle product formulas in small depths are given by 
\begin{align}\begin{split}\label{eq:stuffle23mzv}
\zeta(k_1)\zeta(k_2) &= \zeta(k_1,k_2) + \zeta(k_2,k_1) + \zeta(k_1+k_2)\,\\
\zeta(k_1)\zeta(k_2,k_3) &= \zeta(k_1,k_2,k_3) +  \zeta(k_2,k_1,k_3) + \zeta(k_2,k_3,k_1)  +  \zeta(k_1+k_2,k_3)  +\zeta(k_2,k_1+k_3) \,.
\end{split}
\end{align}
Now, one could ask for a new object $S(k_1,\dots,k_r) \in \R$, which does not satisfy the stuffle product formula, but the index shuffle product formula, that is, which is an image of an algebra homomorphism from $(\Q A_z, \sh)$ to $\R$. By this we mean we want to construct something out of the multiple zeta values which satisfies in low depths
\begin{align}\begin{split}\label{eq:shuffleex}
S(k_1)S(k_2) &= S(k_1,k_2) + S(k_2,k_1) \,,\\
S(k_1)S(k_2,k_3) &= S(k_1,k_2,k_3) +  S(k_2,k_1,k_3) + S(k_2,k_3,k_1) \,.
\end{split}
\end{align}
One can easily check that $S(k)=\zeta(k)$ and 
\begin{align}\begin{split}\label{eq:Sdep23}
S(k_1,k_2) &= \zeta(k_1,k_2) + \frac{1}{2} \zeta(k_1+k_2)\,,\\
S(k_1,k_2,k_3) &= \zeta(k_1,k_2,k_3) + \frac{1}{2} \zeta(k_1+k_2,k_3) + \frac{1}{2}\zeta(k_1,k_2+k_3) + \frac{1}{6}\zeta(k_1+k_2+k_3) \,,
\end{split}
\end{align}
satisfy the index-shuffle product formula \eqref{eq:shuffleex} as a consequence of the stuffle product formula \eqref{eq:stuffle23mzv}.

One might guess that this works in arbitrary depths and that the coefficients are given by $\frac{1}{n!}$, whenever adding together $n$ indices. Since $\frac{1}{n!}$ is the  coefficient of $T^n$ in $\exp(T)$ one might think of the above constructions as some exponential map  $\exp : (\Q A_z, \sh) \rightarrow (\Q A_z, \ast)$, which gives us the algebra homomorphism $S: (\Q A_z, \sh) \rightarrow \R$ by setting $S =  \zeta \circ \exp $.  
In this section, we want to make this precise by associating to some power series $f \in T \ko[[T]]$ a linear map $\Psi_f: \ko\langle A \rangle \rightarrow \ko\langle A \rangle$. In the case $f(T) = \exp(T)-1$ we will get the map in the example above.

Let $w=a_1a_2 \cdots a_n$ be a word of length $\len(w)=n$ with letters $a_1,\dots,a_n \in A$. Let $I=(i_1,\dots,i_m)$ be a composition of $n$, i.e. $i_1+\dots+i_m = n$ with $m\geq 1$, $i_1,\dots,i_m\geq 1$. For such an $I$ we define 
\begin{align*}
I[w] = (a_1 \diamond \cdots \diamond a_{i_1})(a_{i_1+1} \diamond \cdots \diamond a_{i_1+i_2}) \dots (a_{i_1+\cdots + i_{m-1}+1} \diamond \cdots \diamond a_n)\,.
\end{align*}

For example for $w=a_1 a_2 a_3$ a composition of $n=3$ is given by $I=(1,2)$, and we get $I[w]=a_1(a_2 \diamond a_3)$. By $\mathcal{C}(n)$ we denote the set of all compositions of $n$ and usually $n$ will be given by the length $n=\len(w)$ of some word $w$.
\begin{definition} For a formal power series $f = \sum_{i=1}^\infty c_i T^i \in T\ko[[T]]$ we define the $\ko$-linear map $\Psi_f: \QA \rightarrow \QA$ by $\Psi_f(1)=1$ and 
\begin{align*}
\Psi_f(w) = \sum_{I=(i_1,\dots,i_m) \in \mathcal{C}(\len(w))} c_{i_1}\cdots c_{i_m} I[w]
\end{align*}
for a nonempty word $w$. 
\end{definition}  
For example for words of length $2$ and $3$  we have $\mathcal{C}(2) =  \left\{ (1,1), (2)  \right\}$ and $\mathcal{C}(3) = \left\{ (1,1,1), (2,1), (1,2), (3)  \right\}$ which gives for $a_1,a_2,a_3 \in A$
\begin{align*}
\Psi_f(a_1 a_2 ) &= c_1^2 \,a_1 a_2 + c_2 \,a_1 \diamond a_2\,,\\
\Psi_f(a_1 a_2 a_3) &= c_1^3\, a_1 a_2 a_3 +c_1 c_2 \,\big( (a_1 \diamond a_2) a_3+  a_1 (a_2 \diamond a_3)  \big) + c_3\, a_1\diamond a_2 \diamond a_3\,. 
\end{align*}
Observe the similarity of this to \eqref{eq:Sdep23} when $A=A_z$, $z_i \diamond z_j = z_{i+j}$ and $c_i = \frac{1}{i!}$, i.e. $f(T) = \exp(T)-1$. In the case $f=T$, we obtain the identity map on $\QA$, i.e.  $\Psi_f(w) =w$ for all $w\in \QA$.  Now for $f,g \in T\ko[[T]]$ denote by $f \circ g \in T\ko[[T]]$ the usual composition of formal power series, i.e. the power series given by $(f \circ g )(T) = f(g(T))$.

\begin{theorem}  \label{thm:fcircg}For $f,g \in T\ko[[T]]$ we have 
\begin{align*}
\Psi_f \Psi_g = \Psi_{f\circ g}\,.
\end{align*}
\end{theorem}
\begin{proof}
This can be found in \cite[Theorem 3.1]{HI} and follows from a straightforward but tedious calculation, which we will omit here. 
\end{proof}

One interesting case of the linear maps $\Psi_f$ and $\Psi_g$ is given by $f(T)=\exp(T)-1$ and $g(T)= \log(1+T)$. In these cases we have $(f \circ g)(T)= (g \circ f)(T)=T$ and therefore $\Psi_f$ and $\Psi_g$ are inverse to each other by Theorem \ref{thm:fcircg}. As in \cite{HI} we write $\exp := \Psi_f$ and $\log := \Psi_g$, i.e. $\exp, \log \in \Hom(\QA,\QA)$. As illustrated in the introduction of this subsection we have the following result, which was first proven in \cite{H2}.

\begin{theorem} \label{thm:explog}The map 
\begin{align*}
\exp: (\QA,\sh) \longrightarrow (\QA, \qsh)
\end{align*}is a $\ko$-algebra isomorphism with inverse 
\begin{align*}
\log: (\QA,\qsh) \longrightarrow (\QA, \sh)\,.
\end{align*}     
\end{theorem}
\begin{proof}
This is \cite[Theorem 2.5]{H2}.
\end{proof}

Theorem \ref{thm:explog} will be used in Chapter \ref{sec:mes} in the construction of shuffle regularized multiple Eisenstein series, which were introduced in \cite{BT}. Another application will be the proof of Proposition \ref{prop:expqshlogdiamond}, which we will do now. For this, we will first give a more general statement on the linear maps $\Psi_f$. 

\begin{proposition} \label{prop:psiffdiamond} For $f\in T \ko[[T]]$ and $z\in \ko A[[X]]$ we have the following equality in $\QA[[X]]$
\begin{align*}
\Psi_f\left( \frac{1}{1- z X} \right) = \frac{1}{1-f_\diamond(zX)}\,.
\end{align*}
Here $\Psi_f$ acts on $\QA[[X]]$ componentwise. 
\end{proposition}

\begin{proof}
Let $f(T) = \sum_{i=1}^\infty c_i T^i$, then the left-hand side is given by 
\begin{align*}
\Psi_f\left( \frac{1}{1- z X} \right) &= \Psi_f \left(1+ zX + z^2 X^2 + z^3 X^3 + \dots\right) \\
&= 1 + \sum_{n \geq 1}  \sum_{I=(i_1,\dots,i_m) \in \mathcal{C}(n)} c_{i_1}\cdots c_{i_m} I[z^n]  X^n\\
&= 1 + \sum_{n \geq 1}  \sum_{i_1 = 1}^\infty c_{i_1} \underbrace{z \diamond \cdots \diamond z}_{i_1} X^{i_1} \sum_{I=(i_2,\dots,i_m) \in \mathcal{C}(n-i_1)} c_{i_2}\cdots c_{i_m} I[z^{n-i_1}]  X^{n-i_1}\,,
\end{align*}
where the last sum on the right needs to be interpreted as $1$ in the case $n=i_1$ and $0$ if $i_1>n$. Also notice that $I$ has been extended linearly to $\QA$ and it acts on $\QA[[X]]$ componentwise. 
Now recall that $f_\diamond(zX) = \sum_{i=1}^\infty c_{i} \underbrace{z \diamond \cdots \diamond z}_{i} X^{i}$. With this the above equation gives
\begin{align*}
\Psi_f\left( \frac{1}{1- z X} \right) &=  1+ f_\diamond(zX)  \Psi_f\left( \frac{1}{1- z X} \right)  \,,
\end{align*}
from which the statement follows. 
\end{proof}
We can use this Proposition together with Theorem \ref{thm:explog} to prove Proposition \ref{prop:expqshlogdiamond}. 

\begin{proof}[Proof of Proposition \ref{prop:expqshlogdiamond}]
First notice that for $f(T)=\exp(T)-1$ the left-hand side of Proposition \ref{prop:psiffdiamond} is given by 
\begin{align*}
\Psi_f\left( \frac{1}{1- z X} \right) = \exp\left( \frac{1}{1- z X} \right).
\end{align*}
On the other hand we have 
\begin{align*}
\frac{1}{1- z X}  &= 1+ zX + z^2 X^2 + z^3 X^3 + \dots\\
&= 1+ zX + \frac{1}{2!} \left(z \sh z\right) X^2 + \frac{1}{3!} \left(z\sh z \sh z \right) X^3 + \dots  = \exp_{\sh}(zX)\,.
\end{align*}
By Theorem \ref{thm:explog} $\exp: (\QA,\sh) \longrightarrow (\QA, \qsh)$ is an algebra homomorphism and we get 
\begin{align*}
\exp\left( \frac{1}{1- z X} \right) =  \exp\left( \exp_{\sh}(zX) \right) = \exp_{\qsh}(zX)\,,
\end{align*}
which gives for any $z\in \ko A[[X]]$ by Proposition \ref{prop:psiffdiamond}
\begin{align*}
\exp_{\qsh}(zX)= \frac{1}{1-\exp_\diamond(zX)}\,.
\end{align*}
Since this holds for any $z\in \ko A[[X]]$, the $zX$ can be replaced by any power series in $X\ko A[[X]]$, i.e. in particular we can choose $\log_{\diamond}(1+zX) \in X\ko A[[X]]$ for any $z\in \ko A[[X]]$, to get
\begin{align*}
\exp_{\qsh}(\log_{\diamond}(1+zX))= \frac{1}{1-\exp_\diamond(\log_{\diamond}(1+zX))} = \frac{1}{1-zX} \,,
\end{align*}
which is exactly the statement of Proposition \ref{prop:expqshlogdiamond}.
\end{proof}

\subsection{Subalgebras of words with restricted first and last letters}
Now we will present a general statement for quasi-shuffle algebras, which we will use to regularize multiple zeta values in the next section. Since multiple zeta values $\zeta(k_1,\dots,k_r)$ are just defined for indices with $k_1\geq 2$, the map $\zeta$ was    just defined on $\HH^0$. The subspace $\HH^0$ is spanned by words starting in $x$ and ending in $y$, or, when viewed as words in $A_z$, spanned by words not starting in the letter $z_1$. We want to extend this map to all indices, i.e. to the space $\HH^1$. For example to make sense of $\zeta$ for the element $z_1 z_2$ one first notices that
\begin{align*}
z_2 \ast z_1  = z_2 z_1 + z_1 z_2 + z_3\,,
\end{align*}
i.e. $z_1 z_2  =  {\color{blue}z_2} {\color{red}\ast z_1} {\color{blue}-  z_2 z_1 - z_3}$ is a polynomial in  {\color{red}$z_1$} (with respect to $\ast$) with {\color{blue} coefficients in $\HH^0$}. We will then view $z_1$ as a variable $T$ and define the stuffle regularized multiple zeta value as the polynomial $\zeta^\ast(1,2; T) = \zeta(2) T - \zeta(2,1) - \zeta(3)$, which will give us an algebra homomorphism $\zeta^\ast: \HH^1_\ast \rightarrow \mz[T]$. This we will do in the next section after proving a general statement for quasi-shuffle algebras in the following, which assures that a polynomial representation as above is possible in certain cases.

For subsets $S,E \subset A$ we define the following subspace of our quasi-shuffle algebra $Q=\QA$
\begin{align*}
Q_S^E &= \Q + \langle a_1 a_2 \dots a_n \mid a_1 \in S, a_2,\dots,a_{n-1} \in A, a_n\in E\,, n\geq 1\rangle_\Q \,,
\end{align*}
i.e. this is the subspace of $\QA$ of words starting with letters in $S$ and ending with letters in $E$. In particular we have $Q^A_A=Q$ and we omit writing $S$ or $E$ if they  equal $A$, i.e. $Q_S = Q^A_S$, $Q^E = Q_A^E$. 

\begin{proposition}
If $(\ko S,\diamond)$ and $(\ko E,\diamond)$ are subalgebras of $(\ko A,\diamond)$, then $(Q_S^E,\qsh)$ is a subalgebra of $(Q,\qsh)$.
\end{proposition}
\begin{proof} This is again a direct consequence of the definition of the quasi-shuffle product 
\begin{align*}
a w \qsh b v = a (w \qsh b v) + b (a w \qsh v) + (a \diamond b) (w \qsh  v) 
\end{align*}
since the first letters (resp. last letters) of the elements in this product just come from the first letters (resp. last letters) and their $\diamond$ products. 
\end{proof}

\begin{theorem}\label{thm:regmaps}
Assume  $(\ko S,\diamond)$, $(\ko E,\diamond)$ are subalgebras of $(\ko A,\diamond)$ and we have an $a \in A$, such that $A \diamond A\backslash \{a\} \subset \ko \left(A\backslash \{a\}\right)$.
\begin{enumerate}[\textup{(}i\textup{)}] 
\item If $a\in E$ we have $Q^E = Q^E_{A\backslash \{a\} }[a]$ and the map\footnote{In other words, the polynomial     $\sum_{j=0}^m w_j \qsh a^{\qsh j} $ with $w_j \in Q^E_{A\backslash \{a\} }$ gets sent to $\sum_{j=0}^m w_j  T^{j}$.}
\begin{align*}
\operatorname{pol}_a : Q^E &\longrightarrow Q^E_{A\backslash \{a\} }[T]\\
a&\longmapsto T
\end{align*}
is an isomorphism of $\ko$-algebras. 
\item  If $a\in S$ we have $Q_S = Q^{A \backslash \{a\} }_{S}[a]$ and the map
\begin{align*}
    \operatorname{pol}^{a} : Q_S &\longrightarrow Q^{A \backslash \{a\} }_{S}[T]\\
a &\longmapsto T
\end{align*}
is an isomorphism of $\ko$-algebras.
\end{enumerate}
\end{theorem}
\begin{proof}
For (i) we first want to show that $Q^E = Q^E_{A\backslash \{a\} }[a]$, so we need to show that any word $w\in Q^E$ is a polynomial in $a$ (with respect to the product $\qsh$) with coefficients given by linear combinations of words not starting in $a$. We write $w=a^m v$ for $m\geq 0$ and $v = b_1 \cdots b_l \in Q^E_{A\backslash \{a\} }$ and prove the statement by induction on $m$, where the $m=0$ case is clear since $w=v \in Q^E_{A\backslash \{a\} }$. By the definition of the quasi-shuffle product we obtain 
\begin{align*}
a \qsh a^{m-1} v  = m  \,a^m v  + a^{m-1} b_1 ( a \sh b_2\cdots b_l) +\sum_{j=0}^{m-2} a^j (a\diamond a) a^{m-2-j} v + a^{m-1} \sum_{i=1}^l b_1 \cdots (a \diamond b_i) \cdots b_l\,.
\end{align*}
The last three terms are linear combinations of words starting with $a^j$ with $0 \leq j<m$. Here we used for the last two sums the condition $a \diamond b_1 \in \ko  A\backslash \{a\}$. Also all words end with letters in $E$, since $a, b_l\in E$ and therefore $a\diamond b_l \in \ko E$. Therefore by the induction hypothesis we get that $w=a^m v$ is also an element in $Q^E_{A\backslash \{a\} }[a]$. To show that the given maps are isomorphisms we therefore just need to show that the representation as such a polynomial is unique, i.e. the given maps are injective. But this follows from the fact that there are no linear relations among the elements in $Q^E_{A\backslash \{a\} }$. Indeed, assuming that $0=\sum_{j=0}^{m} w_j \qsh  a^{\qsh j}$ with $w_j \in Q^E_{A\backslash \{a\} }$
we immediately get that $w_m =0$ since it is the only part which gives words starting with $a^m$. 
We omit the proof of (ii), since the argument for (ii) is exactly the same, except that we consider words ending in $a^m$ instead of starting in $a^m$. 
\end{proof}

Theorem \ref{thm:regmaps}  can be used to prove the following statements, for the spaces $\HH^0$,$\HH^1$,$\HH$  and the products $\sh$, $\ast$ and $\ast_\gdiamond$. The first two statements are classical results and the last one can be found in \cite[Theorem 2.14]{BK1}.
\begin{corollary}  \label{cor:polexprhh} We have 
\begin{enumerate}[\textup{(}i\textup{)}] 
\item $\HH_\shuffle^1 = \HH_\shuffle^0[y]$  and $\HH_\shuffle = \HH_\shuffle^1[x] = \HH_\shuffle^0[x,y]$.
\item  $\HH^1_\ast = \HH_\ast^0[z_1]$.
\item $\HH^1_{\gqsh} = \HH_{\gqsh}^0[z_1]$.
\end{enumerate}
\end{corollary}
\begin{proof}In all cases we have $\ko =\Q$. For (i) we choose $A=\{x,y\}$, $S=\{x\}$, $E=\{y\}$ and use the trivial product for $\diamond$. With this we have $\HH = Q$, $\HH^1 = Q^E$ and $\HH^0=Q^E_S$. The statement then follows from Theorem  \ref{thm:regmaps} (i) and (ii).  For (ii) and (iii) we choose $A=A_z$, $a=z_1$ and the usual stuffle product for (ii) and the quasi-shuffle product $\gqsh$ for (iii). The latter one was defined by using the following product $\gdiamond$ on $\ko A_z$
\begin{align*}
z_{k_1 }\gdiamond z_{k_2}&=     z_{k_1+k_2}+ \sum_{j=1}^{k_1+k_2-1}\left( \lambda^j_{k_1,k_2}+ \lambda^j_{k_2,k_1} \right) z_j\,
\end{align*}
where the rational numbers $\lambda^j_{k_1,k_2}$ are given by
\begin{align*}
\lambda^j_{k_1,k_2}  = (-1)^{k_2-1} \binom{k_1+k_2-1-j}{k_1-j}  \frac{B_{k_1+k_2-j}}{(k_1+k_2-j)!} \,.
\end{align*}
By the definition of the $\lambda^j_{k_1,k_2}$ one checks that $\lambda^1_{k_1,k_2} + \lambda^1_{k_2,k_1}=0$ whenever $k_1+k_2 \geq 3$ (By using the equation \eqref{eq:l1l1} in the proof of Proposition \ref{prop:gevensubspace}).  This shows that $z_1 \gdiamond A\backslash \{z_1\} \subset  \Q A\backslash \{z_1\} $ and we can apply Theorem  \ref{thm:regmaps}.
\end{proof}

\begin{ex}\label{ex:z112reg}
As an example of Corollary \ref{cor:polexprhh} (i), (ii) and (iii) we give the following expressions of $z_1 z_1 z_2$ as a polynomial in $z_1=y$ having coefficients in $\HH^0$ with respect to the products $\sh$, $\ast$ and $\gqsh$
\begin{align*}
z_1 z_1 z_2 &= \frac{1}{2} z_2 \sh z_1^{\sh 2} - 2 z_2 z_1 \sh z_1+ 3 z_2 z_1 z_1\,,\\
z_1 z_1 z_2 &= \frac{1}{2} z_2 \ast z_1^{\ast 2}-(z_2 z_1+z_3) \ast z_1  + \left(z_2 z_1 z_1 + z_3 z_1 + \frac{1}{2}z_4\right)\,,\\
z_1 z_1 z_2 &= \frac{1}{2} z_2 \,\gqsh\, z_1^{\gqsh 2}-(z_2 z_1+z_3 - z_2) \,\gqsh\, z_1  + \left(z_2 z_1 z_1 + z_3 z_1 + \frac{1}{2}z_4- z_3 - \frac{3}{2}z_2 z_1 + \frac{5}{12} z_2 \right)\,.
\end{align*}
Notice that the product $\sh$ here is with respect to the alphabet $A=\{x,y\}$ and not $A_z$! So for the first statement one should rewrite $z_1 z_1 z_2 = yyxy$, $z_2 z_1 = xyy$, etc. 
If you want to create more examples and play around with $\sh$, $\ast$ and $\gqsh$ you can use the following online tool: \url{https://www.henrikbachmann.com/shuffle.html}, where  these three products are implemented on the space $A_z$. There one could check the above examples by entering  (sh = $\sh$, st = $\ast$, gst = $\,\gqsh$ )
\begin{verbatim}
1/2*[1] sh [1] sh [2] -2*[1] sh [2,1]+ 3*[2,1,1]
1/2*[1] st [1] st [2] -[1] st ([2,1]+[3]) + [2,1,1] + [3,1] + 1/2*[4]
1/2*[1]gst[1]gst[2]-[1]gst([2,1]+[3]-[2])+[2,1,1]+[3,1]+1/2*[4]-[3]-3/2*[2,1]+5/12*[2]
\end{verbatim}
which all give the output $(1,1,2)$.
\end{ex}
\section{Regularizations}

As mentioned already at the beginning of the last section we now want to make sense of multiple zeta values for non-admissible indices by using the previous results on quasi-shuffle algebras. 

\subsection{Stuffle and shuffle regularized multiple zeta values}

As a consequence of  Theorem \ref{thm:regmaps}  (Corollary \ref{cor:polexprhh}) 
we have for $\bullet \in \{\sh , \ast\}$ isomorphisms of $\Q$-algebras
\begin{align*}
\reg^T_\bullet: \HH^1_\bullet \rightarrow \HH^0_{\bullet}[T]\,,
\end{align*}
which send an element $w= \sum_{j=0}^m w_j \bullet z_1^{\bullet j}$ with $w_j \in \HH^0$ to $\reg^T_\bullet(w)= \sum_{j=0}^m w_j  T^{j}$. This enables us to extend the algebra homomorphism $\zeta : \HH^0_\bullet \rightarrow \mz$ to an algebra homomorphism $\zeta^\bullet : \HH^1_\bullet \rightarrow \mz[T]$  by extending $\zeta$ to $\HH^0_\bullet [T]$ and setting $\zeta^\bullet =  \zeta \circ \reg^T_\bullet$, i.e. we have the following commutative diagram of $\Q$-algebra homomorphisms

\adjustbox{scale=1.1,center}{%
\begin{tikzcd}
\HH^1_\bullet \arrow[r,"\reg^T_\bullet"]\arrow{rd}[swap]{\zeta^\bullet} & \HH^0_{\bullet}[T] \arrow[d,"\zeta"]\\
& \mz[T]
\end{tikzcd}
    }

\begin{definition} \label{def:regmzv}Let $\kk = (k_1,\dots,k_r) \in \Z_{\geq 1}^r$  be any index. 
    \begin{enumerate}[\textup{(}i\textup{)}] 
        \item We define the \emph{shuffle regularized multiple zeta value} 
        \begin{align*}
        \zeta^\sh(\kk;T) =\zeta^\sh(k_1,\dots,k_r;T) := \zeta^\sh(z_\kk) \in \mz[T]
        \end{align*}
        In the case $T=0$ we just write $\zeta^\sh(\kk)= \zeta^\sh(\kk;0)$.\footnote{In the literature often ``shuffle/stuffle regularized multiple zeta value'' refers to the $T=0$ case.}
        \item We define the \emph{stuffle regularized multiple zeta value} 
\begin{align*}
\zeta^\ast(\kk;T) =\zeta^\ast(k_1,\dots,k_r;T) := \zeta^\ast(z_\kk) \in \mz[T]
\end{align*}
In the case $T=0$ we just write $\zeta^\ast(\kk)= \zeta^\ast(\kk;0)$. 
    \end{enumerate}
\end{definition}

From Example \ref{ex:z112reg} we obtain 
\begin{align*}
\zeta^\sh(1,1,2;T) &= \frac{1}{2}\zeta(2) T^2 - 2\zeta(2,1)T+ 3\zeta(2,1,1)\,,\\
\zeta^\ast(1,1,2;T) &= \frac{1}{2}\zeta(2) T^2 -(\zeta(2,1)+\zeta(3))T+ \zeta(2,1,1) + \zeta(3,1)+ \frac{1}{2}\zeta(4)\,.
\end{align*}
Although the coefficients of $T^2$ and $T$ are the same (because we know $\zeta(2,1)=\zeta(3)$), the constant terms differ. In general $\zeta^\sh(\kk;T)$ and $\zeta^\ast(\kk;T)$ are different if $\kk$ is non-admissible. However, in the next section we will see that there is an explicit relationship between these polynomials. Further, Theorem \ref{thm:eds} will show that for $w\in \HH^0, v\in \HH^1$ and   $\bullet \in \{\sh, \ast\}$ we have the extended double shuffle relations
\begin{align*}
\zeta^\bullet(w \sh v - w \ast v)=0\,.
\end{align*}
Since $\zeta^\sh$ and $\zeta^\ast$ differ on $\HH^1$ this is not obvious at all.  For example, we have  (Exercise \ref{ex9})
\begin{align*}
\zeta^\bullet(z_2 \sh z_1 z_1 - z_2 \ast z_1 z_1) =0\,,
\end{align*} 
which implies linear relations among multiple zeta values in weight four.

\subsection{Comparison of $\zeta^\sh$ and $\zeta^\ast$}

As we saw before, the two regularizations $\zeta^\sh(\kk;T)$ and $\zeta^\ast(\kk;T)$ differ as elements in $\R[T]$. In this section, we will present the exact relationship between these two regularizations as it was done in \cite{IKZ}.
For this first consider the following series 
\begin{align*}
A(u) &=\exp\left( \sum_{n=2}^\infty \frac{(-1)^n}{n} \zeta(n) u^n \right)\\
&= 1+ \frac{\zeta(2)}{2} u^2 - \frac{\zeta(3)}{3} u^3 + \left( \frac{\zeta(4)}{4} + \frac{\zeta(2)^2}{8}\right) u^4 - \left( \frac{\zeta(5)}{5} + \frac{\zeta(2) \zeta(3)}{6}\right) u^5 +\dots \,\\
&=: \sum_{k \geq 0} \gamma_k \,u^k\,.
\end{align*}
Here the $\gamma_k \in \Q[\zeta(j) \mid j\geq 2 ]$ are polynomials of single zeta values which, considered as multiple zeta values, have homogeneous weight  $k$. 
Using this we define the  $\R$-linear map $\rho: \R[T]\rightarrow \R[T]$ by 
\begin{align}\label{eq:rhodef}
\rho(e^{Tu}) := A(u) e^{Tu}\,.
\end{align}
Notice that this defines the linear map $\rho$ uniquely by comparing the coefficients of $u^m$ on both sides. Since $\rho$ is linear, we get 
\begin{align*}
\rho(e^{Tu}) &= \rho(1) + \rho(T) u + \frac{1}{2}\rho(T^2)  u^2+ \frac{1}{3!} \rho(T^3) u^3 + \dots \\
&:= \left( 1 + \gamma_1 u + \gamma_2 u^2 + \dots \right) \left( 1+ Tu + \frac{1}{2}(Tu)^2 + \frac{1}{3!} (Tu)^3 + \dots \right) = A(u) e^{Tu}\,,
\end{align*}
and therefore we obtain for $m\geq 0$ the explicit formula 
\begin{align*}
\rho\left(T^m\right) = m! \sum_{k=0}^{m} \gamma_k \frac{T^{m-k}}{(m-k)!}\,.
\end{align*}
Also notice that $\rho$ is bijective. For the first values of $m$ we get
\begin{align*}
\rho(1) &= 1\,,\\
\rho(T) &= T\,,\\
\rho(T^2) &= T^2 + \zeta(2)\,,\\
\rho(T^3) &= T^3 + 3 \zeta(2) T - 2 \zeta(3)\,,\\
\rho(T^4) &= T^4 + 6 \zeta(2) T^2 - 8 \zeta(3) T + 6 \zeta(4)+ 3\zeta(2)^2\,,\\
\rho(T^5) &= T^5 + 10 \zeta(2) T^3 - 20 \zeta(3) T^2 + \left(30 \zeta(4) + 15 \zeta(2)^2\right) T- 24 \zeta(5)- 20 \zeta(2) \zeta(3)\,.
\end{align*}

The stuffle regularized multiple zeta values are elements in $\R[T]$ and for example we saw that      
\begin{align*}
\zeta^\ast(1,1,2;T) &= \frac{1}{2}\zeta(2) T^2 -(\zeta(2,1)+\zeta(3))T+ \zeta(2,1,1) + \zeta(3,1)+ \frac{1}{2}\zeta(4)\,.
\end{align*}
Applying the linear map $\rho$ to this, we see that we just get an additional contribution of $\frac{1}{2}\zeta(2)^2$, i.e.
\begin{align*}
\rho\left(\zeta^\ast(1,1,2;T) \right) = \frac{1}{2}\zeta(2) T^2 -(\zeta(2,1)+\zeta(3))T+ \zeta(2,1,1) + \zeta(3,1)+ \frac{1}{2}\zeta(4)+ \frac{1}{2}\zeta(2)^2\,.
\end{align*}
Using the known relations $\zeta(3)=\zeta(2,1)$ and $\zeta(4)=\zeta(2,1,1)$ (duality), $\zeta(3,1)= \frac{1}{4}\zeta(4)$ (finite double shuffle) and $\zeta(2)^2 = \frac{5}{2} \zeta(4)$ (Euler), we get
\begin{align*}
\rho\left(\zeta^\ast(1,1,2;T) \right) &= \frac{1}{2}\zeta(2) T^2 -2\zeta(2,1)T+ 3\zeta(2,1,1) \\
&= \zeta^\sh(1,1,2;T)\,,
\end{align*}
i.e. the $\rho$ sends the stuffle regularized multiple zeta value to the shuffle regularized multiple zeta value. In general, the map $\rho$ has this property and we have the following. 
\begin{theorem} \label{thm:comparezastzsh}For all $\kk \in \Z_{\geq 1}^r$ we have
    \begin{align*}
    \zeta^\sh(\kk;T)  = \rho\left( \zeta^\ast(\kk;T) \right) \,.
    \end{align*}
    Or equivalently, when viewed as maps from $\HH^1$ to $\R[T]$, we have $\zeta^\sh = \rho \circ \zeta^\ast$.
\end{theorem}
\begin{proof} This is \cite[Theorem 1]{IKZ} and we just give a sketch of the proof here. A really detailed version of this proof can also be found in the book of Zhao \cite[Section 3.3.2]{Zh2}.
    The main idea is to compare the behavior of the truncated multiple zeta values $\zeta_M(\kk)$ (which satisfy the stuffle product formula for all $\kk$) and the multiple polylogarithm $\li_{\kk}(z)$ (which satisfies the shuffle product formula for all $\kk$) as $M\rightarrow \infty$ and $z \rightarrow 1$. We have the classical formula
    \begin{align*}
    \zeta_M(1) = 1 + \frac{1}{2} + \cdots + \frac{1}{M-1} = \log(M) + \gamma + \operatorname{O}\left(\frac{1}{M}\right),
    \end{align*}
    as $M \rightarrow \infty$, where $\gamma = 0.57721\dots$ denotes the Euler--Mascheroni constant.
    As a consequence of the stuffle product formula one can show by induction that for some $J$ we have (see \cite[Lemma 3.3.19]{Zh2})
    \begin{align}
    \zeta_M(\kk) = \zeta^\ast(\kk ; \log(M)+ \gamma) + \operatorname{O}( M^{-1} \log^J(M)) \qquad (\text{as } M \rightarrow \infty)\,.
    \end{align}
    Similarly by $\li_1(z) = \log\left(\frac{1}{1-z}\right)$ and using the shuffle product formula for $\li$ one can show, together with the fact that for admissible $\kk$ we have $\li_\kk(1) - \li_\kk(z) = \operatorname{O}(1-z)$, that  (see \cite[Lemma 3.3.20]{Zh2})
    \begin{align*}
    \li_\kk(z) = \zeta^\sh\left(\kk; \log\left(\frac{1}{1-z}\right) \right) + \operatorname{O}\left( (1-z) \log^J\left(\frac{1}{1-z} \right)\right)\qquad  (\text{as } z \rightarrow 1)\,.
    \end{align*}
    The connection between the multiple polylogarithm and the truncated multiple zeta values is that the latter basically give the Taylor coefficients of $\li_\kk$:
    \begin{align*}
    \li_\kk(z) = \li_{k_1,\dots,k_r}(z) &= \sum_{m_1 > \cdots > m_r > 0} \frac{z^{m_1}}{m_1^{k_1} \cdots m_r^{k_r}}= \sum_{m=1}^\infty \left( \sum_{m > m_2 > \cdots > m_r > 0} \frac{1}{m^{k_1} m_2^{k_2} \cdots m_r^{k_r}} \right) z^m  \\ 
    &= \sum_{m=1}^\infty \left( \zeta_{m+1}(\kk) - \zeta_m(\kk) \right) z^m  = (1-z) \sum_{m=1}^\infty \zeta_m(\kk)  z^{m-1} \,.
    \end{align*}
    The statement then follows from the following general fact, that for any polynomial $P(T)\in \R[T]$ and $Q(T)=\rho(P(T))$ one has  (see \cite[Lemma 1]{IKZ} or \cite[Lemma 3.3.17]{Zh2})
    \begin{align*}
    (1-z)\sum_{m=1}^\infty P\left( \log(m) + \gamma \right) z^{m-1} = Q\left( \log\left(\frac{1}{1-z}\right)  \right) + \operatorname{O}\left( (1-z) \log^J\left(\frac{1}{1-z} \right)\right)\qquad  (\text{as } z \rightarrow 1)
    \end{align*}
    for $J= \deg(P)-1$ and for $l\geq 0$
    \begin{align*}
    \sum_{m=1}^\infty \frac{\log^l(m)}{m} z^{m-1} = \operatorname{O}\left( \log^{l+1}\left(\frac{1}{1-z}\right)  \right) \qquad  (\text{as } z \rightarrow 1)\,.
    \end{align*}
    Choosing $P(T)= \zeta^\ast(\kk;T)$ and combining all the equations above gives $Q(T)=\zeta^\sh(\kk;T)$ and therefore $\zeta^\sh(\kk;T)=\rho(\zeta^\ast(\kk;T))$.
\end{proof}

\subsection{Stuffle regularization from the $q$-series $\g$}

The stuffle regularization can also be obtained as a limit of the $q$-series $\g$.
Recall that for any index $\kk=(k_1,\dots,k_r)$ we have
\begin{align*}
\g(\kk)&=\sum_{m_1>\cdots>m_r>0}
 \frac{P_{k_1}(q^{m_1})}{(1-q^{m_1})^{k_1}}\cdots
 \frac{P_{k_r}(q^{m_r})}{(1-q^{m_r})^{k_r}},
\end{align*}
and, if $\kk$ is admissible,
\begin{align}\label{eq:glimit}
 \lim_{q\rightarrow1}(1-q)^{\wt(\kk)}\g(\kk)=\zeta(\kk).
\end{align}
For $k\geq0$, set
\begin{align*}
\gsadm&=\big\langle\g(\kk)\mid \kk\text{ admissible}\big\rangle_\Q,\\
\gsadmk&=\big\langle\g(\kk)\mid \kk\text{ admissible},\ \wt(\kk)\leq k\big\rangle_\Q
\end{align*}
and define
\begin{align}\label{eq:defzkforg0}
\zm_k:\gsadmk&\longrightarrow\mz_k,\\
f&\longmapsto\lim_{q\rightarrow1}(1-q)^k f.
\end{align}
Equation \eqref{eq:glimit} shows that this map is surjective. We now extend it from
$\gsadm$ to
\begin{align*}
 \gsk=\big\langle\g(\kk)\mid \wt(\kk)\leq k\big\rangle_\Q.
\end{align*}
Recall that $\g:\HH^1_{\gqsh}\rightarrow\gs$ is an algebra homomorphism and
$\gsadm=\g(\HH^0)$. Corollary \ref{cor:polexprhh} therefore gives the following.

\begin{proposition}\label{prop:gg0}
\begin{enumerate}[\textup{(}i\textup{)}]
 \item We have $\gs=\gsadm[\g(1)]$.
 \item The series $\g(1)$ is algebraically independent over $\gsadm$.
\end{enumerate}
\end{proposition}
\begin{proof}
Statement (i) follows directly from Corollary \ref{cor:polexprhh}, and statement (ii) is
\cite[Theorem~2.14]{BK1}. The estimates used there are
\begin{align*}
 \g(1)\asymp\frac{-\log(1-q)}{1-q},\qquad
 \g(\kk)\asymp\frac{1}{(1-q)^{\wt(\kk)}}
\end{align*}
as $q\rightarrow1$ for admissible $\kk$. Here $f(q)\asymp g(q)$ means that the
absolute value of $f(q)/g(q)$ is bounded from above and below by positive constants
for $q$ sufficiently close to $1$. See also \cite[Lemma~2]{P}.
\end{proof}

The filtered form in the proof of \cite[Theorem~2.14]{BK1} shows that every
$f\in\gsk$ can be written uniquely as
\begin{align*}
 f=\sum_{j=0}^k f_j\g(1)^{k-j},\qquad f_j\in\gsadmj.
\end{align*}
We can therefore define
\begin{align*}
\zm_k^T:\gsk&\longrightarrow\mz[T],\\
\sum_{j=0}^k f_j\g(1)^{k-j}&\longmapsto
 \sum_{j=0}^k\zm_j(f_j)T^{k-j}.
\end{align*}

\begin{proposition}\label{prop:glimitregularization}
For every index $\kk$ we have
\begin{align*}
 \zm^T_{\wt(\kk)}\big(\g(\kk)\big)=\zeta^\ast(\kk;T).
\end{align*}
\end{proposition}
\begin{proof}
Let $w\in\HH^1$ be homogeneous of weight $k$. By Corollary
\ref{cor:polexprhh}, there are $u_j,v_j\in\HH^0$ such that
\begin{align*}
 w=\sum_{j=0}^k u_j\ast z_1^{\ast j}
  =\sum_{j=0}^k v_j\,\gqsh\,z_1^{\gqsh j},
\end{align*}
where $u_j$ is homogeneous of weight $k-j$ and $v_j$ has weight at most $k-j$.
The highest-weight part of
\begin{align*}
 z_a\gdiamond z_b=z_{a+b}+
 \sum_{l=1}^{a+b-1}\big(\lambda^l_{a,b}+\lambda^l_{b,a}\big)z_l
\end{align*}
is $z_{a+b}$. Thus the product induced by $\gqsh$ on the associated graded space is
the usual stuffle product. The uniqueness of the polynomial decomposition in
Corollary \ref{cor:polexprhh} now shows that $v_j-u_j$ has weight at most
$k-j-1$. Therefore,
\begin{align*}
 \zm_{k-j}\big(\g(v_j)\big)=\zm_{k-j}\big(\g(u_j)\big)=\zeta(u_j).
\end{align*}
Applying $\g$ to the second decomposition of $w$ and then using the definition of
$\zm_k^T$, we obtain
\begin{align*}
 \zm_k^T\big(\g(w)\big)
 &=\sum_{j=0}^k\zm_{k-j}\big(\g(v_j)\big)T^j
 =\sum_{j=0}^k\zeta(u_j)T^j
 =\zeta^\ast(w;T).
\end{align*}
This proves the statement. Compare also \cite[Theorem~2.14 and Proposition~6.4]{BK1}.
\end{proof}

\section{Extended double shuffle relations}
As an extension of the finite double shuffle relations (Proposition \ref{prop:fds}) we will now present a family of linear relations, which give conjecturally all linear relations among multiple zeta values. We define for $w,u \in \HH^1$ the element 
\begin{align*}
\ds(w,u) := w \sh u - w\ast u  \in \HH^1.
\end{align*} 
The statement of Proposition \ref{prop:fds} then was, that $\ds(w,u) \in \ker \zeta$ if $w,u \in \HH^0$. The extended version states, that one of the words $w$ and $u$ is allowed to be in $\HH^1$. In this case $\ds(w,u)$ is not necessarily in $\HH^0$ anymore, but after projecting to $\HH^0$ by the map $\reg^T_\bullet: \HH^1_\bullet \rightarrow \HH^0_{\bullet}[T]$ and then comparing the coefficients of $T$ (or setting $T=0$, for which we write $\reg_\bullet:=\reg^0_\bullet$) one still obtains a relation among multiple zeta values.  In other words $\ds(w,u)$ is in the kernel of the regularized multiple zeta value maps. 
\begin{theorem}[Extended double shuffle relations]\label{thm:eds} For $w\in \HH^1$, $u \in \HH^0$ and $\bullet \in \{\sh, \ast\}$  we have 
    \begin{align*}
    \zeta^\bullet(w \sh u - w\ast u; T) = 0\,,
    \end{align*}
    i.e. in particular $\reg_\bullet \left( \ds(w,u) \right) \in \ker \zeta$.
\end{theorem}
\begin{proof} By Theorem  \ref{thm:comparezastzsh} we have for all $w \in \HH^1$ 
    \begin{align*}
    \zeta^\sh(w;T)  = \rho\left( \zeta^\ast(w;T) \right) \,.
    \end{align*}
    Multiplying both sides with $\zeta^\sh(u)=\zeta^\ast(u)=\zeta(u) \in \R$, for $u\in \HH^0$, we get by the $\R$-linearity of $\rho$ 
    \begin{align*}
    \zeta^\sh(w \sh u;T)  &= \rho\left( \zeta^\ast(w \ast u;T) \right) \,\\
    &= \zeta^\sh(w \ast u;T)  \,.
    \end{align*}
    This gives $\zeta^\sh(w \sh u - w \ast u;T)=0$ and $\zeta^\ast(w \sh u - w \ast u;T)=0$ by applying the inverse of $\rho$.
\end{proof}

We can now give the precise version of Conjecture \ref{conj:extdsh}, which we announced in Chapter \ref{sec:overview}.

\begin{conjecture}\label[conjecture]{conj:extdshprecise} The kernel of $\zeta: \HH^0 \rightarrow \mz$ is given by
    \begin{align*}
    \ker \zeta &= \big\langle \reg_\sh \left( w \sh u - w \ast u \right)  \mid  w\in \HH^1,\, u \in \HH^0 \big\rangle_\Q \\
    &= \big\langle \reg_\ast \left( w \sh u - w \ast u \right)  \mid  w\in \HH^1,\, u \in \HH^0 \big\rangle_\Q\,,
    \end{align*}
    i.e. the extended double shuffle relations give all $\Q$-linear relations among multiple zeta values. 
\end{conjecture}

Since it is expected that the extended double shuffle relations give all relations among multiple zeta values, one could obtain upper bounds for the dimension of $\mz_k$, i.e., an alternative proof of Theorem \ref{thm:terasoma}, by counting these relations. This is still an open problem.

\begin{problem} Count the number of linearly independent extended double shuffle relations.
    Define for $k\geq 0$ and $\bullet \in \{\sh, \ast\}$ the spaces
    \begin{align*}
    \eds^\bullet_{k} = \big\langle \reg^\bullet \left( w \sh u - w \ast u \right)  \mid  w\in \HH^1,\, u \in \HH^0 , \wt(w)+\wt(u)=k \big\rangle_\Q\,.
    \end{align*}
    Show that for all $k\geq 2$ and any  $\bullet \in \{\sh, \ast\}$ we have
    \begin{align}\label{eq:dimeds}
    \dim_\Q \eds^\bullet_{k} = 2^{k-2} - d_k\,.
    \end{align}
    Here the $d_k$ are the conjectured dimensions (Conjecture \ref{conj:zagier}) of $\mz_k$, defined by     
    \begin{align*}
    \sum_{k\geq 0} d_k X^k = \frac{1}{1-X^2-X^3}\,,
    \end{align*}
    and $2^{k-2}$ is the number of admissible indices of weight $k$. 
\end{problem}
So far the equation \eqref{eq:dimeds} has been checked up to $k=21$ (T. Machide, T. Sonobe, 2020+) by extensive computer calculations. 

All the relations we obtained so far should be a consequence of the extended double shuffle relations. For the finite double shuffle relations this is obvious, but for the duality relation this is actually also an open problem. Recall that the duality relations (Proposition \ref{prop:duality}) stated that for all $v\in \HH^0$ we have 
\begin{align*}    
\zeta(\tau(v)) = \zeta(v)\,,
\end{align*}
where $\tau$ was the anti-automorphism defined on $\HH^0$ with $\tau(y)=x$ and $\tau(x)=y$.

\begin{problem}\label{prob:dualityeds} Show that for any $v \in \HH^0$ we have for  $\bullet \in \{\sh, \ast\}$
    \begin{align*}
    \tau(v) - v \in \big\langle \reg_\bullet \left( w \sh u - w \ast u \right)  \mid  w\in \HH^1,\, u \in \HH^0 \big\rangle_\Q\,,
    \end{align*}
    i.e. the duality relation is a consequence of the extended double shuffle relations.
\end{problem}

We now want to discuss a refinement of the extended double shuffle relations. For this, we first consider the following special case.

\begin{proposition}[Hoffman's relation (\cite{H1}, Theorem~5.1)] \label{prop:hoffmanrel} For an admissible index $\kk=(k_1,\dots,k_r)$ we have
    \begin{align*}
    \sum_{i=1}^r \zeta(k_1,\dots,k_{i-1}, k_i + 1, k_{i+1}, \dots, k_r) = \sum_{\substack{1 \leq i \leq r\\ k_i \geq 2}} \sum_{j=0}^{k_i-2} \zeta(k_1,\dots,k_{i-1}, k_i - j, j+1, k_{i+1},\dots,k_r)\,.
    \end{align*}
\end{proposition}
\begin{proof} 
    This is a special case of the extended double shuffle relation by choosing $w = z_1 = y$ and $u = z_\kk$. In particular, we have $\ds(z_1, z_\kk) \in \HH^0$ (Exercise \ref{ex10}) and therefore no regularization is necessary.  Another proof of this relation, using partial fraction decomposition, can be found in \cite[Theorem 1]{Zu3}.
\end{proof}

Hoffman's relation is a special case of the extended double shuffle relations, which are not a consequence of the finite double shuffle relations. Surprisingly, it seems that this is the only part of the extended double shuffle relations we need.

\begin{conjecture}[H. N. Minh, M. Petitot et al. (\cite{MJOP})] We have 
    \begin{align*}
    \ker \zeta &= \big\langle  w \sh u - w \ast u   \mid  w\in \HH^0 \cup \{z_1\},\, u \in \HH^0 \big\rangle_\Q \,,
    \end{align*}
    i.e. Hoffman's relation and the finite double shuffle relations give all linear relations among multiple zeta values. 
\end{conjecture}

But still, this set of relations seems to be too much, and the following gives an even better refinement.

\begin{conjecture}[M. Kaneko, M. Noro and K. Tsurumaki (\cite{KNT})] We have 
    \begin{align*}
    \ker \zeta &= \big\langle  w \sh u - w \ast u   \mid  w\in  \{z_1, z_2, z_3, z_2z_1\},\, u \in \HH^0 \big\rangle_\Q \,.
    \end{align*}
\end{conjecture}

The equivalence of these conjectures is not known and in particular we have the following open problem. 
\begin{problem} Show that 
    \begin{align*}
    \big\langle \reg_\sh \left( w \sh u - w \ast u \right)  \mid  w\in \HH^1,\, u \in \HH^0 \big\rangle_\Q  &= \big\langle \reg_\ast \left( w \sh u - w \ast u \right)  \mid  w\in \HH^1,\, u \in \HH^0 \big\rangle_\Q  \\
    &= \big\langle  w \sh u - w \ast u   \mid  w\in \HH^0 \cup \{z_1\},\, u \in \HH^0 \big\rangle_\Q \\
    &= \big\langle  w \sh u - w \ast u   \mid  w\in  \{z_1, z_2, z_3, z_2z_1\},\, u \in \HH^0 \big\rangle_\Q \,.
    \end{align*}
\end{problem}

Having an element $w\in \ker \zeta$ and $u \in \HH^0$, then clearly also $w \ast u, w \sh u \in \ker \zeta$. As far as the author knows, it is also still an open problem to show that the space of extended double shuffle relations is closed under the multiplication with elements in $\HH^0$. 

\begin{problem} Show that for $\bullet_1, \bullet_2 \in \{\sh, \ast\}$ we have 
    \begin{align*}
    \HH^0 \bullet_1 \big\langle \reg_{\bullet_2} \left( w \sh u - w \ast u \right)  \mid  w\in \HH^1,\, u \in \HH^0 \big\rangle_\Q   \subset \big\langle \reg_{\bullet_2} \left( w \sh u - w \ast u \right)  \mid  w\in \HH^1,\, u \in \HH^0 \big\rangle_\Q\,.
    \end{align*}
\end{problem}

\section{Hopf algebra structure}

In this section, we want to introduce the Hopf algebra structure of quasi-shuffle algebras. This can be useful to simplify some proofs on statements among multiple zeta values or can be used to understand multiple Eisenstein series and their regularizations. For this we first roughly recall what a Hopf algebra is.  Informally, a Hopf algebra is a vector space equipped simultaneously with algebra and coalgebra structures that interact compatibly, along with a map called the antipode acting as a generalized inverse. We now recall the precise definition:

\begin{definition}[Hopf algebra]
A \emph{Hopf algebra} over a field $\ko$ is a vector space $H$ equipped with:
\begin{enumerate}[(i)]
    \item an algebra structure $(H,m,u)$, where
    \[ m: H \otimes H \to H \quad \text{(multiplication)},\quad u: \ko \to H \quad \text{(unit)}, \]
    \item a coalgebra structure $(H,\Delta,\varepsilon)$, where
    \[ \Delta: H \to H \otimes H \quad \text{(coproduct)},\quad \varepsilon: H \to \ko \quad \text{(counit)}, \]
    \item and a linear map called the \emph{antipode}:
    \[ S: H \to H. \]
\end{enumerate}
These structures must satisfy the following axioms:
\begin{enumerate}[(i)]
    \item $(H,m,u)$ is an associative unital algebra.
    \item $(H,\Delta,\varepsilon)$ is a coassociative counital coalgebra, i.e.:
    \[(\Delta \otimes \mathrm{id}) \circ \Delta = (\mathrm{id} \otimes \Delta) \circ \Delta, \quad (\varepsilon \otimes \mathrm{id}) \circ \Delta = (\mathrm{id} \otimes \varepsilon) \circ \Delta = \mathrm{id}.\]
    \item Compatibility of structures: the coproduct $\Delta$ and counit $\varepsilon$ are algebra homomorphisms.
    \item Antipode axiom:
    \[m \circ (S \otimes \mathrm{id}) \circ \Delta = m \circ (\mathrm{id} \otimes S) \circ \Delta = u \circ \varepsilon.\]
\end{enumerate}
\end{definition}
The Antipode axiom can be interpreted nicely by introducing the \emph{convolution product}. Let $H$ be a Hopf algebra as above, and let $B$ be a $\ko$-algebra with multiplication $m_B$. Then, for $\ko$-algebra homomorphisms $f,g \in \Hom_{\ko\text{-alg}}(H,B)$, define their convolution product as
\begin{align*}
    f \star g = m_B \circ (f \otimes g) \circ \Delta
\end{align*}
This gives a $\ko$-linear map $f \star g: H \rightarrow B$ and one can check the following:

\begin{lemma} Let $H$ be a commutative Hopf algebra over $\ko$ and $B$ a commutative $\ko$-algebra. For $f,g \in \Hom_{\ko\text{-alg}}(H,B)$ we have $f \star g \in \Hom_{\ko\text{-alg}}(H,B)$.
\end{lemma}
\begin{proof}
Let $f, g \in \Hom_{\ko\text{-alg}}(H, B)$. Preservation of the unit follows since $\Delta(1_H) = 1_H \otimes 1_H$, so $(f \star g)(1_H) = m_B(f(1_H) \otimes g(1_H)) = 1_B$.
For multiplicativity, using Sweedler notation, 
\[
(f \star g)(h h') = \sum f(h_{(1)} h'_{(1)}) \cdot g(h_{(2)} h'_{(2)}) = \sum f(h_{(1)}) f(h'_{(1)}) \cdot g(h_{(2)}) g(h'_{(2)}),
\]
since $f$ and $g$ are multiplicative and $\Delta(h h') = \Delta(h) \Delta(h')$.

On the other hand,
\[
(f \star g)(h) \cdot (f \star g)(h') = \left( \sum f(h_{(1)}) \cdot g(h_{(2)}) \right) \cdot \left( \sum f(h'_{(1)}) \cdot g(h'_{(2)}) \right) = \sum f(h_{(1)}) \cdot g(h_{(2)}) \cdot f(h'_{(1)}) \cdot g(h'_{(2)}),
\]
where the last step expands the product. This matches the earlier expression because $B$ is commutative.    
\end{proof}

The quasi-shuffle algebra $\ko\langle A \rangle$ generated by words over an alphabet $A$ carries a natural Hopf algebra structure. For a word $w \in \ko\langle A \rangle$, let $\Delta:\ko\langle A \rangle \to \ko \langle A \rangle \otimes \ko \langle A \rangle$ be the $\ko$-linear map, called \emph{deconcatenation coproduct}, defined by
\begin{align*}
\Delta(w) = \sum_{uv=w} u \otimes v
\end{align*}
and $\varepsilon:\ko\langle A \rangle \to \ko$ be the $\ko$-linear map defined by
\begin{align*}
\varepsilon(w) = \begin{cases}1 &\text{if }w=1,\\0&\text{otherwise.}\end{cases}
\end{align*}
The canonical convolution product for $g, h:\ko\langle A \rangle \to B$ is given by
\begin{align}\label{eq:convolution}
(g \star h)(w) = \sum_{uv=w} g(u) h(v).
\end{align}

\begin{proposition}[\cite{HI}]
Equipped with $\Delta$ as a coproduct and $\varepsilon$ as a counit, $\ko\langle A \rangle$ is a coalgebra.
\end{proposition}
\begin{proof}
First, we will show that the diagram
\begin{center}
\begin{tikzcd}
&\ko\langle A \rangle \otimes \ko \langle A \rangle \otimes \ko \langle A \rangle&\\
\ko \langle A \rangle \otimes \ko \langle A \rangle \arrow{ru}{\Delta \otimes \id}&& \ko \langle A \rangle \otimes \ko \langle A \rangle \arrow[lu, swap, "\id \otimes \Delta"]\\
&\ko \langle A \rangle \arrow[lu, below, "\Delta"] \arrow[ru, swap, "\Delta"]
\end{tikzcd}
\end{center}
commutes.

Let $w \in \ko\langle A \rangle$ be a word. Then,
\begin{align*}
(\Delta \otimes \id)(\Delta(w)) &= (\Delta \otimes \id)\left(\sum_{uv=w} u \otimes v \right)
= \sum_{uv=w} \Delta(u) \otimes v
= \sum_{uv=w}\sum_{st=u} s \otimes t \otimes v
= \sum_{stv=w} s \otimes t \otimes v\\
&= \sum_{su=w} \sum_{tv=u} s\otimes t \otimes v
= \sum_{su=w} s \otimes \Delta(u)
= (\id \otimes \Delta)\left(\sum_{su=w}s \otimes u\right)
= (\id \otimes \Delta)(\Delta(w)).
\end{align*}

Next, we will show that for the bilinear maps
\begin{align*}
f:\ko\langle A \rangle \otimes \ko &\longrightarrow \ko \langle A \rangle\\
a \otimes 1 &\longmapsto a
\end{align*}
and
\begin{align*}
g:\ko \otimes \ko \langle A \rangle &\longrightarrow \ko \langle A \rangle\\
1 \otimes a &\longmapsto a,
\end{align*}
the diagram
\begin{center}
\begin{tikzcd}
\ko\langle A \rangle \otimes \ko \arrow[r, "f"] &\ko \langle A \rangle & \ko \otimes \ko \langle A \rangle \arrow[l, swap, "g"]\\
\ko\langle A \rangle \otimes \ko \langle A \rangle \arrow[u, "\id \otimes \varepsilon"] &\ko\langle A \rangle \arrow[l, "\Delta"]\arrow[r, swap, "\Delta"] \arrow[u, equals, "\id"] & \ko \langle A \rangle \otimes \ko \langle A \rangle \arrow[u, swap, "\varepsilon \otimes \id"]
\end{tikzcd}
\end{center}
commutes.

Let $w \in \ko \langle A \rangle$ be a word. We will first show that the left side commutes by
\begin{align*}
(f \circ (\id \otimes \varepsilon) \circ \Delta)(w) &= f\left((\id \otimes \varepsilon)\left(\sum_{uv=w}u \otimes v\right)\right)
= f\left(\sum_{uv=w}u \otimes \varepsilon(v)\right)
= f(w \otimes 1)
= w.
\end{align*}
The right-hand side of the diagram commutes with a similar argument.
\end{proof}

Define the following linear map
\begin{align*}
R:\ko\langle A \rangle &\longrightarrow \ko \langle A \rangle\\
a_1 \cdots a_n &\longmapsto a_n \cdots a_1\\
\emptyword &\longmapsto \emptyword.
\end{align*}
Clearly, $R$ is an involution. Moreover, we have the following proposition.

\begin{proposition}[{\cite[Proposition 4.2]{HI}}]\label{prop:R-automorphism}
$R$ is an automorphism of $(\ko \langle A \rangle, *_{\diamond})$.
\end{proposition}
\begin{proof}
Clearly, $R$ is a bijection. Let $S_{k,l}$ be the subgroup of the symmetric group $S_{k+l}$ consisting of all permutations $\sigma$ such that $\sigma(1) < \sigma(2) < \cdots < \sigma(k)$ and $\sigma(k+1) < \sigma(k+2) < \cdots < \sigma(k+l)$. Then, that $R$ is an endomorphism of $(\ko \langle A \rangle, \shuffle)$ follows directly from
\begin{align*}
R(a_1 a_2 \cdots a_k \shuffle a_{k+1} a_{k+2} \cdots a_{k+l}) &= \sum_{\sigma \in S_{k,l}} R(a_{\sigma(1)}a_{\sigma(2)} \cdots a_{\sigma(k+l)})=\sum_{\sigma \in S_{k,l}} a_{\sigma(k+l)} \cdots a_{\sigma(1)}\\
&= R(a_{k+1} \cdots a_{k+l}) \shuffle R(a_1 \cdots a_k).
\end{align*}
To obtain the statement for the quasi-shuffle product, we use Theorem \ref{thm:explog} and the fact that $R$ commutes with the maps $\exp$ and $\log$, which follows from Proposition \ref{prop:Rpsif} below. For $u,v\in\ko\langle A\rangle$ we then have
\begin{align*}
R(u\qsh v)&=R\bigl(\exp(\log(u)\shuffle\log(v))\bigr)
=\exp\bigl(R(\log(u))\shuffle R(\log(v))\bigr)\\
&=\exp\bigl(\log(R(u))\shuffle\log(R(v))\bigr)
=R(u)\qsh R(v)\,.
\end{align*}
\end{proof}

\begin{proposition}[{\cite[Proposition 4.3]{HI}}]\label{prop:Rpsif}
$R$ commutes with $\Psi_f$ for all $f \in T\ko[[T]]$.
\end{proposition}
\begin{proof}
For a composition $I = (i_1, \dots, i_l)$, we denote $\overline{I}=(i_l, \dots, i_1)$. For $f=c_1t + c_2t^2 + \cdots$,
\begin{align*}
R(\Psi_f(w)) &= \sum_{I=(i_1, \dots, i_m) \in \mathcal C(\ell(w))} c_{i_1} \cdots c_{i_m}R(I[w]) = \sum_{I=(i_1, \dots, i_m) \in \mathcal C(\ell(w))} c_{i_1} \cdots c_{i_m} \overline{I}[R(w)]\\
&= \sum_{\overline{I}=(i_m, \dots, i_1)} c_{i_m} \cdots c_{i_1} \overline{I}[R(w)] = \Psi_f(R(w)).
\end{align*}
\end{proof}

\begin{proposition}[{\cite[Theorem 4.2]{HI}}]\label{prop:antipode-relation}
The object $(\ko \langle A \rangle, *_\diamond, \Delta)$ is a Hopf algebra with $S=\Sigma TR$ as antipode. Here $R$ is the map defined above, $T$ is the $\ko$-linear map given by $T(w) = (-1)^{\len(w)} w$ and $\Sigma$ is the $\ko$-linear map given by $\Sigma(w) = \sum_{I \in \mathcal{C}(\len(w))} I[w]$.
\end{proposition}
\begin{proof}
We will first show that $(\ko \langle A \rangle, *_\diamond, \Delta)$ is a bialgebra by showing that $\Delta$ and $\varepsilon$ are algebra homomorphisms with respect to $*_\diamond$.

First,
\begin{align*}
\varepsilon(\emptyword *_\diamond \emptyword) = 1 = \varepsilon(\emptyword) \varepsilon(\emptyword)
\end{align*}
and for words $w, v \in \ko \langle A \rangle$ such that $\ell(w) + \ell(v) > 0$,
\begin{align*}
\varepsilon(w *_\diamond v) = 0 = \varepsilon(w) \varepsilon(v).
\end{align*}
Therefore, $\varepsilon$ is an algebra homomorphism.

Next, we will use mathematical induction on $\ell(w) + \ell(v)$. First,
\begin{align*}
\Delta(\emptyword *_\diamond \emptyword) &= \emptyword \otimes \emptyword= (\emptyword *_\diamond \emptyword) \otimes (\emptyword *_\diamond \emptyword)= (\emptyword \otimes \emptyword) *_\diamond (\emptyword \otimes \emptyword)= \Delta(\emptyword) *_\diamond \Delta(\emptyword).
\end{align*}

Suppose that for any words $w'', v'' \in \ko \langle A \rangle$ such that $\ell(w'') + \ell(v'') \leq k$,
\begin{align*}
\Delta(w'' *_\diamond v'') = \Delta(w'') *_\diamond \Delta(v'').
\end{align*}

Let $w'=aw, v'=bv \in \ko\langle A \rangle$ be words such that $\ell(w') + \ell(v')=k+1$. With Sweedler's notation,
\begin{align*}
\Delta(w) = \sum w_{(1)} \otimes w_{(2)}\text{ and }\Delta(v) = \sum v_{(1)} \otimes v_{(2)}.
\end{align*}
We get 
\begin{align*}
\Delta(w') *_\diamond \Delta(v') =& \Delta(aw) *_\diamond \Delta(bv)
= \left(\emptyword \otimes aw + \sum aw_{(1)} \otimes w_{(2)}\right) *_\diamond \left(\emptyword \otimes bv + \sum b v_{(1)} \otimes v_{(2)}\right)\\
=& \emptyword \otimes (aw *_\diamond bv) + \sum b v_{(1)} \otimes (aw *_\diamond v_{(2)}) + \sum aw_{(1)} \otimes (w_{(2)} *_\diamond bv) \\&+ \sum(aw_{(1)} *_\diamond bv_{(1)}) \otimes (w_{(2)} *_\diamond v_{(2)})\\
=& (\emptyword \otimes a(w *_\diamond bv)) + (\emptyword \otimes b(aw *_\diamond v)) + (\emptyword \otimes (a \diamond b) (w *_\diamond v))\\
&+\sum bv_{(1)} \otimes (aw *_\diamond v_{(2)}) + \sum aw_{(1)} \otimes (w_{(2)} *_\diamond bv) + \sum a(w_{(1)} *_\diamond bv_{(1)}) \otimes (w_{(2)} *_\diamond v_{(2)})\\
&+ \sum b(aw_{(1)} *_\diamond v_{(1)}) \otimes (w_{(2)} *_\diamond v_{(2)}) + \sum(a \diamond b)(w_{(1)} *_\diamond v_{(1)}) \otimes (w_{(2)} *_\diamond v_{(2)})\\
=& \emptyword \otimes a(w *_\diamond bv) + \emptyword \otimes b(aw *_\diamond v) + \emptyword \otimes (a \diamond b)(w *_\diamond v)\\
&+ \sum(a \otimes \emptyword) \left(\sum w_{(1)} \otimes (w_{(2)} *_\diamond bv) + \sum(w_{(1)} *_\diamond bv_{(1)}) \otimes (w_{(2)} *_\diamond v_{(2)})\right)\\
&+ \sum(b \otimes \emptyword)\left(\sum v_{(1)} \otimes (aw *_\diamond v_{(2)}) + \sum (a w_{(1)} *_\diamond v_{(1)}) \otimes (w_{(2)} *_\diamond v_{(2)})\right)\\
&+((a \diamond b) \otimes \emptyword)(\Delta(w) *_\diamond \Delta(v))\\
=& \emptyword \otimes a (w *_\diamond bv) + \emptyword \otimes b(aw *_\diamond v) + \emptyword \otimes (a \diamond b) (w*_\diamond v)\\
&+ (a \otimes \emptyword)\left(\left(\sum w_{(1)} \otimes w_{(2)}\right) *_\diamond \left(\emptyword \otimes bv + \sum bv_{(1)} \otimes v_{(2)}\right)\right)\\
&+ (b \otimes \emptyword)\left(\left(\sum v_{(1)} \otimes v_{(2)}\right) *_\diamond \left(\emptyword \otimes aw + \sum aw_{(1)} \otimes w_{(2)}\right)\right)\\
&+ ((a \diamond b) \otimes \emptyword) (\Delta(w)*_\diamond \Delta(v))\\
=& \emptyword \otimes a(w *_\diamond bv) + \emptyword \otimes b(aw *_\diamond v) + \emptyword \otimes (a \diamond b)(w *_\diamond v)\\
&+ (a \otimes \emptyword) (\Delta(w) *_\diamond \Delta(bv)) + (b \otimes \emptyword)(\Delta(v) *_\diamond \Delta(aw)) + ((a \diamond b) \otimes \emptyword) \Delta(w *_\diamond v)\\
=& \emptyword \otimes a (w *_\diamond bv) + (a \otimes \emptyword) \Delta(w *_\diamond bv) + \emptyword \otimes b(aw *_\diamond v) + (b \otimes \emptyword) \Delta(aw *_\diamond v)\\
&+ \emptyword \otimes (a \diamond b)(w *_\diamond v) + ((a \diamond b) \otimes \emptyword) \Delta(w *_\diamond v)\\
=& \Delta(a(w *_\diamond bv) + b(aw *_\diamond v) + (a \diamond b) (w *_\diamond v))\\
=& \Delta(aw*_\diamond bv).
\end{align*}
Therefore, $(\ko \langle A \rangle, *_\diamond, \Delta)$ is a bialgebra. Moreover, $\ko \langle A \rangle$ is filtered by word length and connected ($\ko \langle A \rangle^0=\ko \emptyword$). So, $(\ko \langle A \rangle, *_\diamond, \Delta)$ is a Hopf algebra with a unique antipode.

Finally, we will show that $S$ is the antipode of $(\ko \langle A \rangle, *_\diamond, \Delta)$. First, $S(\emptyword) = \emptyword$ and for $a_1, \dots, a_n \in A$,
\begin{align*}
S(a_1 \cdots a_n) = \Sigma TR(a_1 \cdots a_n) = \Sigma T(a_n \cdots a_1) = (-1)^n \Sigma(a_n \cdots a_1) = (-1)^n \sum_{I \in \mathcal C(n)} I[a_n \cdots a_1].
\end{align*}

Let $w \in \ko \langle A \rangle$ be a word. We will show the antipode relation
\begin{align}\label{eq:antipode-relation}
\sum_{uv=w} S(u) *_\diamond v= (\eta \circ \varepsilon)(w) = \sum_{uv=w} u *_\diamond S(v).
\end{align}

First, we have
\begin{align*}
\emptyword = S(\emptyword) *_\diamond \emptyword = (\eta \circ \varepsilon)(\emptyword) = \emptyword *_\diamond S(\emptyword).
\end{align*}

For word $w = a_1 \cdots a_n \in \ko \langle A \rangle$ such that $n \neq 0$,
\begin{align*}
\sum_{uv=w}S(u) *_\diamond v =& \sum_{i=0}^n S(a_1 \cdots a_i) *_\diamond a_{i+1} \cdots a_n\\
=& \sum_{i=0}^n (-1)^i \sum_{I \in \mathcal C(i)}(I[a_i \cdots a_1] *_\diamond a_{i+1} \cdots a_n)\\
=& \sum_{i=0}^n (-1)^i \sum_{k=1}^i \sum_{I \in \mathcal C(k-1)}((a_i \diamond \cdots \diamond a_k)I[a_{k-1} \cdots a_1]*_\diamond a_{i+1} \cdots a_n)\\
=& \sum_{i=0}^n (-1)^i \sum_{k=1}^i \sum_{I \in \mathcal C(k-1)}(a_i \diamond \cdots \diamond a_k)(I[a_{k-1} \cdots a_1] *_\diamond a_{i+1} \cdots a_n)\\
&+ \sum_{i=0}^{n-1}(-1)^i \sum_{k=1}^i \sum_{I \in \mathcal C(k-1)} a_{i+1}((a_i \diamond \cdots \diamond a_k)I[a_{k-1} \cdots a_1] *_\diamond a_{i+2} \cdots a_n)\\
&+ \sum_{i=0}^{n-1}(-1)^i \sum_{k=1}^i \sum_{I \in \mathcal C(k-1)}(a_{i+1} \diamond \cdots \diamond a_k) (I[a_{k-1} \cdots a_1]*_\diamond a_{i+2} \cdots a_n)\\
=& \sum_{i=0}^n(-1)^i \sum_{k=1}^i \sum_{I \in \mathcal C(k-1)}(a_i \diamond \cdots \diamond a_k) (I[a_{k-1} \cdots a_1] *_\diamond a_{i+1} \cdots a_n)\\
&+ \sum_{i=0}^{n-1} (-1)^i \sum_{k=1}^i \sum_{I \in \mathcal C(k-1)}(a_{i+1} \diamond \cdots \diamond a_k)(I[a_{k-1} \cdots a_1] *_\diamond a_{i+2} \cdots a_n)\\
&+ \sum_{i=0}^{n-1} (-1)^i \sum_{k=1}^i \sum_{I \in \mathcal C(k-1)} a_{i+1}((a_i \diamond \cdots \diamond a_k) I[a_{k-1} \cdots a_1] *_\diamond a_{i+2} \cdots a_n))\\
=& \sum_{i=0}^n (-1)^i \sum_{k=1}^i \sum_{I \in \mathcal C(k-1)} (a_i \diamond \cdots \diamond a_k) (I[a_{k-1} \cdots a_1] *_\diamond a_{i+1} \cdots a_n)\\
&- \sum_{i=1}^n (-1)^i \sum_{k=1}^{i-1} \sum_{I \in \mathcal C(k-1)} (a_i \diamond \cdots \diamond a_k) (I[a_{k-1} \cdots a_1] *_\diamond a_{i+1} \cdots a_n)\\
&+ \sum_{i=0}^{n-1} (-1)^i \sum_{k=1}^i \sum_{I \in \mathcal C(k-1)} a_{i+1}((a_i \diamond \cdots \diamond a_k) I[a_{k-1} \cdots a_1] *_\diamond a_{i+2} \cdots a_n)\\
=& \sum_{i=1}^n (-1)^i \sum_{I \in \mathcal C(i-1)} a_i(I[a_{i-1} \cdots a_1] *_\diamond a_{i+1} \cdots a_n)\\
&+ \sum_{i=1}^n (-1)^{i-1} \sum_{k=1}^{i-1} \sum_{I \in \mathcal C(k-1)} a_i ((a_{i-1} \diamond \cdots \diamond a_k) I[a_{k-1} \cdots a_1] *_\diamond a_{i+1} \cdots a_n)\\
&= 0\\
&= (\eta \circ \varepsilon)(w).
\end{align*}

Equation $(\eta \circ \varepsilon)(w)=\sum_{uv=w}u*_\diamond S(v)$ can be proven in a similar way.
\end{proof}

The convolution product gives a short proof of the stuffle product formula for the
symmetric sums introduced in Section~\ref{subsec:fmzv} and studied further in
Chapter~\ref{sec:fmzv}.

\begin{proposition}\label{prop:symmetricstuffle}
Extend both maps linearly to $\HH^1$. For $m\geq1$ and
$w,v\in\HH^1$, we have
\begin{align*}
 S_m(w\ast v)&=S_m(w)S_m(v),\\
 \zs^\ast(w\ast v)&=\zs^\ast(w)\zs^\ast(v).
\end{align*}
\end{proposition}
\begin{proof}
Define the signed reversal $\widetilde R:\HH^1\rightarrow\HH^1$ by
\begin{align*}
 \widetilde R(1)=1,\qquad
 \widetilde R(z_{k_1}\cdots z_{k_r})=(-1)^{k_1+\dots+k_r}z_{k_r}\cdots z_{k_1}.
\end{align*}
The stuffle product is homogeneous with respect to the weight. Therefore Proposition
\ref{prop:R-automorphism} shows that $\widetilde R$ is an automorphism of
$\HH^1_\ast$. For $m\geq1$, let $h_m:\HH^1_\ast\rightarrow\Q$ be the algebra
homomorphism given by $h_m(z_\kk)=H_m(\kk)$. Since $h_m$, $h_m\circ\widetilde R$,
$\zeta^\ast$ and $\zeta^\ast\circ\widetilde R$ are algebra homomorphisms, so are their
convolution products. Expanding them with \eqref{eq:convolution} gives
\begin{align*}
 \big((h_m\circ\widetilde R)\star h_m\big)(z_{k_1}\cdots z_{k_r})
 &=\sum_{j=0}^r(-1)^{k_1+\dots+k_j}
 H_m(k_j,\dots,k_1)H_m(k_{j+1},\dots,k_r)\\
 &=S_m(k_1,\dots,k_r)
\end{align*}
and, in the same way,
\begin{align*}
 \big((\zeta^\ast\circ\widetilde R)\star\zeta^\ast\big)
 (z_{k_1}\cdots z_{k_r})=\zs^\ast(k_1,\dots,k_r).
\end{align*}
The convolution lemma now proves both formulas. This is also the Hopf algebra proof
of \cite[Proposition~8]{KZ} (with the opposite convention for the order of an index).
\end{proof}

\section{Double shuffle relations for \texorpdfstring{$q$}{q}-analogues}\label{sec:qmzv}
In this section, we introduce the algebraic setup and the analogues of the double shuffle relations for $q$-analogues of multiple zeta values. This algebraic setup was, with a slightly different notation, introduced by Takeyama in \cite{Tak1} for an extended version of the Bradley-Zhao multiple zeta values. For a good overview how these relate to other models of $q$-analogues we refer to the overview article \cite{Bri}. We will use mainly the notations from \cite{Tak2} together with some self-made modification.\\

Recall that there is an easy description for the analogue of the harmonic product for the $q$-analogues $\zeta_q$. In order to give an analogue of the shuffle product, one needs to allow some kind of extended versions of $q$-analogues. Recall that for $k\geq 1$ we considered the term $\frac{q^{(k-1)m}}{[m]_q^k}$ in the definition of $\zeta_q$, in particular for $k=1$ the term $\frac{1}{[m]_q}$ appeared. We will extend this now by introducing a second term for the weight $1$ index, which will correspond to $\frac{q^m}{[m]_q}$. For this we define $\overline{\N} = \{ \bar{1} \} \cup \Z_{\geq 1}$ and then define for $k\in \overline{\N} $ and $m\geq 1$ the following $q$-series
\begin{align}\label{eq:deffm}
    f_k(m) = \begin{cases}  \frac{q^m}{[m]_q},&\text{ if } k = \bar{1}\\
    \frac{q^{(k-1)m}}{[m]_q^k},&\text{ if } k \in \Z_{\geq 1}
\end{cases}\,.
\end{align}
With this we extend the definition of $\zeta_q$ in the following way.

\begin{definition}
\begin{enumerate}[(i)]
    \item For $\kk = (k_1,\dots,k_r) \in \ON^r$ with $k_1 \neq 1$ define the (extended) $q$-multiple zeta value
    \begin{align*}
        \zeta_q(\kk) =  \zeta_q(k_1,\dots,k_r) = \sum_{m_1 > \dots > m_r > 0} \prod_{j=1}^r f_{k_j}(m_j) \in \Q[[q]]\,.
    \end{align*}
    \item For any $\kk = (k_1,\dots,k_r) \in \ON^r$ define the $q$-multiple polylogarithm by
    \begin{align*}
        \li^q_{\kk}(z) = \li^q_{k_1,\dots,k_r}(z) = \sum_{m_1 > \dots > m_r > 0} z^{m_1} \prod_{j=1}^r f_{k_j}(m_j) \in \Q[[q]][[z]]\,.
    \end{align*}
\end{enumerate}
\end{definition}

Notice that we view both objects as formal power series, but both can also be viewed as complex functions in $q,z \in \C$ with $|q|,|z| < 1$. In the case $k_1 \neq 1$ we can also consider the limit $\lim_{z\rightarrow 1} \li^q_{\kk}(z) =  \zeta_q(\kk)$ and in the case $k_1 \neq \bar{1},1$ we have $\lim_{q \rightarrow 1} \zeta_q(\kk) = \zeta(\kk)$. Moreover for any $\kk = (k_1,\dots,k_r) \in \ON^r$ we have $\lim_{q \rightarrow 1} \li^q_\kk = \li_\kk$.

The goal is to introduce an algebraic setup, which considers several different $q$-analogue models at the same time. We saw that the additional index $\bar{1}$ will correspond to a factor $\frac{q^m}{[m]_q}$ in the definition of $\zeta_q$. To generalize this further, we allow indices $k, \bar{k}$ and $\hat{k}$. For each of these we will introduce a letter $e_k, e_{\bar{k}}$ and $e_{\hat{k}}$. The rough idea is that we have the following correspondences, which indicate the relationship of letters to their corresponding factor in the definition of $\zeta_q$
\begin{align}\begin{split}\label{eq:ekmeaning}
    e_k &\longleftrightarrow  \frac{q^{(k-1)m}}{[m]^k_q}\\
    e_{\hat{k}} &\longleftrightarrow \frac{q^{k m}}{[m]^k_q}\\  
    e_{\bar{k}} &\longleftrightarrow\frac{q^{m}}{[m]^k_q}\,.
\end{split}
\end{align}
Notice that all these factors on the right-side evaluate to $\frac{1}{m^k}$ as $q\rightarrow 1$. Since we have
\begin{align}\label{eq:qmterms}
    \frac{q^{km}}{[m]_q^k} &= (-1)^{k-1} (1-q)^{k-1} \frac{q^m}{[m]_q} + \sum_{j=2}^k (-1)^{k-j} (1-q)^{k-j} \frac{q^{(j-1)m}}{[m]_q^j}\,,\qquad (k\geq 1)\\
    \frac{q^{m}}{[m]_q^k} &= \sum_{j=2}^k \binom{k-2}{j-2} (1-q)^{k-j} \frac{q^{(j-1)m}}{[m]_q^j}\,,\qquad (k\geq 2)
\end{align}
it is sufficient to just introduce $e_{\bar{1}}$ and $e_{k}$ for $k\geq 1$ in our setup, and then define $e_{\hat{k}}$ and $e_{\bar{k}}$ in terms of  $e_{\bar{1}}$ and $e_{k}$ by using the formulas \eqref{eq:qmterms}. In our algebraic setup the term $(1-q)$ will correspond to a formal variable $\hbar$ and instead of working with $\Q$-vector spaces we will work with $\CC$-modules, where $$\CC= \Q[\hbar, \hbar^{-1}]\,.$$

\begin{definition}\begin{enumerate}[(i)]
    \item We define $\HHH = \CC\langle a,b \rangle$ and its two $\CC$-submodules
\begin{align*}
    \HHH^0 = \CC + a \HHH b \quad \subset \quad \HHH^1 = \CC + \HHH b \quad \subset \quad \HHH\,.
\end{align*}
\item For $k\geq 1$ we write $e_k = a^{k-1}(a+\hbar)b$ and $e_{\bar{1}}=ab$. With this we have $b=\frac{1}{\hbar}(e_1-e_{\bar{1}})$ and 
\begin{align*}
    \HHH^1 &= \CC\langle e_{\bar{1}}, e_1, e_2, e_3, \dots \rangle,\\
    \HHH^0 &= \CC + \langle e_{k_1}\dots e_{k_r} \mid r\geq 1, k_1,\dots,k_r \in \ON, k_1 \neq 1\rangle_{\CC}\,.
\end{align*}
\item We set $e_{\widehat{1}} = e_{\bar{1}} = ab$ and for $k\geq 2$ define, inspired by \eqref{eq:ekmeaning} and \eqref{eq:qmterms},
\begin{align}\label{eq:defhatebare}
    e_{\widehat{k}} &:= (-1)^{k-1} \hbar^{k-1} e_{\bar{1}} + \sum_{j=2}^k (-1)^{k-j} \hbar^{k-j} e_{j} = a^k b\,,\\
    e_{\overline{k}} &:= \sum_{j=2}^k \binom{k-2}{j-2} \hbar^{k-j} e_j\,.
\end{align}
\end{enumerate}
\end{definition}

We define a $\CC$-module structure on $\Q[[q]]$ and $\Q[[q]][[z]]$ by defining the multiplication of $\hbar$ to be the multiplication with $(1-q)$. Since the words in $\HHH^0$ correspond exactly to those indices for which we defined $\zeta_q$, we can define the $\CC$-linear map 
\begin{align}\begin{split}\label{eq:zetaqmap}
\zeta_q \colon \HHH^0 &\longrightarrow \Q[[q]] \\
e_{\kk} := e_{k_1}\dots e_{k_r} &\longmapsto \zeta_q(\kk)\,,
\end{split}
\end{align}
where $\kk = (k_1,\dots,k_r) \in \ON^r$. Similarly we can define the $\CC$-linear map 
\begin{align}\begin{split}
\li^q \colon \HHH^1 &\longrightarrow \Q[[q]][[z]] \\
e_{\kk} &\longmapsto \li^q_\kk(z)\,.
\end{split}
\end{align}
Using \eqref{eq:defhatebare} we extend $\zeta_q(\kk)$ and $\li^q_\kk$ to indices $\kk = (k_1,\dots,k_r) \in \{k,\bar{k},\widehat{k} \mid k \in \Z_{\geq 1} \}$ by setting $\zeta_q(k_1,\dots,k_r) := \zeta_q(e_{k_1}\dots e_{k_r})$ and  $\li^q_{k_1,\dots,k_r}(z) := \li^q_{e_{k_1}\dots e_{k_r}}$. For example, using \eqref{eq:ekmeaning} and \eqref{eq:qmterms} we get
\begin{align*}
    \zeta_q(\widehat{9}, 3, \bar{4}) = \sum_{m_1 > m_2 > m_3 > 0} \frac{q^{9m_1}}{[m_1]_q^9} \frac{q^{2m_2}}{[m_2]_q^3} \frac{q^{m_3}}{[m_3]_q^4}\,.
\end{align*}

One can embed the classical $\HH$ into $\HHH$ in various different ways. For example, the embedding $\iota: x \mapsto a$, $\iota: y \mapsto (a+\hbar)b$ yields $\iota(z_k) = e_k$. Another variant would be $\widehat{\iota}: x \mapsto a$, $\widehat{\iota}: y \mapsto ab$ yielding $\widehat{\iota}(z_k) =   e_{\widehat{k}}$. To extend the $q$-stuffle product $\ast_q$, defined in Section \ref{subsubsec:qshex} on the space $\HH^1$, we first observe that the functions $f_k$ defined in \eqref{eq:deffm} satisfy for $k,k_1,k_2 \in \Z_{\geq 1}$ 
\begin{align*}
    f_{\bar{1}}^2 &= f_2 - (1-q) f_{\bar{1}},\qquad f_{\bar{1}} f_{k} = f_{k+1},\qquad f_{k_1} f_{k_2} = f_{k_1+k_2} + (1-q) f_{k_1+k_2-1}\,.
\end{align*}
Motivated by this we define $\ast_q = \qsh$ to be the quasi-shuffle product\footnote{Notice that the original definition of quasi-shuffle products is defined over fields and not rings. But all results we use here for quasi-shuffle products also work for algebras over rings.} on $\HHH^1 = \CC\langle \widehat{A}\rangle$ defined by the alphabet $\widehat{A}=\{e_{\bar{1}},e_1,e_2,\dots \}$ and the product $\diamond$ on $\CC \widehat{A}$ given by
\begin{align*}
    e_{\bar{1}} \diamond  e_{\bar{1}} &= e_2 - \hbar e_{\bar{1}},\qquad e_{\bar{1}}\diamond  e_{k} = e_{k+1},\qquad e_{k_1}\diamond  e_{k_2} = e_{k_1+k_2} + \hbar e_{k_1+k_2-1}\,.
\end{align*}
Notice that we have  $e_{\widehat{k_1}} \diamond e_{\widehat{k_2}} = e_{\widehat{k_1+k_2}}$ for any $k_1,k_2 \geq 1$.
This gives a commutative $\CC$-algebra $\HHH^1_{\ast_q}$. The subspace $\HHH^0$ is closed under $\ast_q$, which gives a subalgebra $\HHH^0_{\ast_q} \subset \HHH^1_{\ast_q}$, and as a consequence of Lemma \ref{lem:fmtoFM} we obtain:
\begin{proposition}The map $\zeta_q:\HHH^0_{\ast_q} \rightarrow \Q[[q]]$ is a $\CC$-algebra homomorphism.
\end{proposition}
Notice that $\ast_q$ is a natural extension of same-named quasi-shuffle product on $\HH^1$, since $\iota: \HH^1_{\ast_q} \rightarrow \HHH^1_{\ast_q}$ is an injective $\CC$-algebra homomorphism.\\

We now turn our attention to the analogue of the shuffle product for $q$-analogues. For this we define an action of $\HHH$ on the $\CC$-module $z \Q[[q,z]]$ for a $g \in z \Q[[q,z]]$ by $1 g = g$ and
\begin{align*}
    a g &= (1-q) \sum_{j=1}^\infty g_{|z=q^j z}\,\\
    b g &= \frac{z}{1-z} g\,.
\end{align*}
For any word $w=l_1\dots l_m \in \HHH$ with $l_j \in \{a,b\}$ we then define  recursively $w g= (l_1 \dots l_{m-1}) (l_m g)$ and extend this $\CC$-linearly to all of $\HHH$. Notice that  this action is just well defined on  $z \Q[[q,z]]$ since the action of $a$ would not be well-defined for power series having a constant term with respect to $z$. But we can extend the action to $\Q[[q,z]]$ for elements in $w \in \HHH^1$, since $b \Q[[q,z]] \subset z\Q[[q,z]]$.

\begin{ex} Consider the element $1 \in \Q[[q,z]]$ and $w=e_{\bar{1}}=ab$, then we have 
\begin{align*}
    e_{\bar{1}} 1 &= ab 1 = a \frac{z}{1-z} = (1-q) \sum_{j=1}^\infty \frac{q^j z}{1-q^j z} = (1-q) \sum_{j=1}^\infty \sum_{m=1}^\infty (q^j z)^m \\
    &= \sum_{m>0} \frac{(1-q) z^m q^m}{(1-q^m)}  = \sum_{m>0} \frac{z^m q^m}{[m]_q} = \li^q_{\bar{1}}(z) \,.
\end{align*}
If we want to calculate $e_{\hat{2}} 1 = aab 1$ we can use above calculations and get
\begin{align*}
    e_{\hat{2}} 1 &= a \li^q_{\bar{1}}(z)  = (1-q) \sum_{j=1}^\infty \sum_{m>0} \frac{(q^j z)^m q^m}{[m]_q} =  \sum_{m>0} \frac{q^{2m} z^m (1-q)}{(1-q^m)[m]_q} =  \sum_{m>0}\frac{z^m q^{2m}}{[m]_q^2} = \li^q_{\hat{2}}(z)\,.
\end{align*}
\end{ex}
In general, we obtain the following. 
\begin{proposition}\label{prop:qpolylogab} For all $w \in \HHH^1$ we have $\li^q_w = w 1$.
\end{proposition}
\begin{proof}
    This is Exercise \ref{fex11}.
\end{proof}

Now we observe (Exercise \ref{fex12}) that for $f,g \in z \Q[[q,z]]$ we have 
\begin{align}\begin{split}\label{eq:afag}
    (af) (ag) &= a ( (af) g + f(ag) + \hbar fg),\\
    (bf) g &= f (bg) = b (fg)\,.
    \end{split}
\end{align}

Motivated by Proposition \ref{prop:qpolylogab} and \eqref{eq:afag} we give the following definition.

\begin{definition} We define the $q$-shuffle product $\shq$ on $\HHH$ as the $\CC$-bilinear product, which satisfies $1 \shq w = w \shq 1 =w$ for any $w\in\HHH$  and 
\begin{align}
    aw \shq av &= a( (aw) \shq v + w \shq (av) + \hbar w \shq v)\,, \\
    (bw) \shq v &= w \shq (bv) = b (w \shq v)\,
\end{align}
for words $w,v \in \HHH$.
\end{definition}

\begin{ex} As an example we calculate 
\begin{align*}
    e_{\bar{1}} \shq e_{\bar{1}} = ab \shq ab = 2abab + \hbar abb = 2 e_{\bar{1}} e_{\bar{1}} + \hbar e_{\bar{1}} b = 2 e_{\bar{1}} e_{\bar{1}} + e_{\bar{1}} e_{1} - e_{\bar{1}} e_{\bar{1}}
\end{align*}
and
\begin{align*}
    e_{\hat{2}} \shq e_{\hat{2}} &= 4 e_{\hat{3}} e_{\hat{1}}  + 2 e_{\hat{2}} e_{\hat{2}} + \hbar (2 e_{\hat{3}} b + 4 e_{\hat{2}} e_{\hat{1}} ) + \hbar^2 e_{\hat{2}} b\\
    &= 4 e_{\hat{3}} e_{\hat{1}}  + 2 e_{\hat{2}} e_{\hat{2}} + 2 e_{\hat{3}}e_{1}-2 e_{\hat{3}}e_{\bar{1}} + \hbar (e_{\hat{2}} e_{1}  + 3e_{\hat{2}} e_{\bar{1}} ).
\end{align*}
Compare this to the classical shuffle products $z_1 \shuffle z_1 = 2  z_1 z_1$ and $z_2 \shuffle z_2 = 4 z_3 z_1 + 2 z_2 z_2$. 
\end{ex}

\begin{proposition} \begin{enumerate}[(i)]
    \item The space $\HHH_{\shq}=(\HHH, \shq)$ is a commutative $\CC$-algebra.
    \item The spaces $\HHH^1_{\shq}=(\HHH^1, \shq)$ and $\HHH^0_{\shq}=(\HHH^0, \shq)$ are subalgebras of $\HHH_{\shq}$.
\end{enumerate}
\end{proposition}
\begin{proof}
    The proofs are similar to the case of quasi-shuffle products. 
\end{proof}

\begin{remark} In general one can check that $\shq$ can be seen as the $q$-analogue of the shuffle product $\shuffle$, since the $q \rightarrow 1$ in our algebraic setup corresponds to setting $\hbar =0$. More precisely consider the following projection 
\begin{align*}
   \tilde{\HH}^1 := \Q[\hbar]\langle e_{\bar{1}},e_1,e_2,\dots \rangle \quad \underset{ \hbar = 0}{\longrightarrow}\quad \Q + \Q\langle a,ab\rangle ab \quad\underset{\substack{ab \mapsto y\\ a \mapsto x} }{\longrightarrow}\quad \HH^1 = \Q + \Q\langle x, y \rangle y = \Q\langle z_1,z_2,\dots \rangle\,.
\end{align*}
The space $\tilde{\HH}^1 \subset \HHH^1$ is closed under $\shq$ and the above map gives an algebra homomorphism from $\tilde{\HH}^1_{\shq}$ to $\HH^1_\shuffle$, which can be checked directly by the definition of $\shq$ after setting $\hbar = 0$ and observing that in this case $e_{\hat{k}} = e_{\bar{k}} = e_k = a^k b$.\end{remark}

\begin{proposition}\label{prop:qpolylogalghom} The maps $\li^q : \HHH^1_{\shq} \rightarrow \Q[[q,z]]$ and $\zeta_q: \HHH^0_{\shq} \rightarrow \Q[[q]]$ are $\CC$-algebra homomorphisms.
\end{proposition}
\begin{proof}
To show that $\li^q$ is a $\CC$-algebra homomorphism we need to check that for $w,v \in \HHH^1$ we have $\li^q(w) \li^q(v) = \li^q(w \shq v)$. This can be done as in Proposition \ref{prop:liisalghom}, i.e. by induction on $\len(w) + \len(v)$ together with Proposition \ref{prop:qpolylogab} and \eqref{eq:afag}. The second statement follows by taking the limit $z\rightarrow 1$.
\end{proof}

As an analogue of the finite double shuffle relations we obtain the following.

\begin{corollary}\label{cor:qfds}For $w,v \in \HHH^0$ we have 
\begin{align*}
    \zeta_q(w \shq v - w \ast_q v) = 0\,.
\end{align*}
\end{corollary}

\begin{exs} \label{ex:qdsh}\begin{enumerate}[(i)]
    \item For $w=v=e_{\bar{1}}$ we get 
    \begin{align*}
        e_{\bar{1}} \shq e_{\bar{1}} - e_{\bar{1}} \ast_q e_{\bar{1}} = 2 e_{\bar{1}} e_{\bar{1}} + \hbar \underbrace{e_{\bar{1}} b}_{abb} - \left( 2 e_{\bar{1}} e_{\bar{1}}  + \underbrace{e_2 - \hbar e_{\bar{1}}}_{aab} \right) = \hbar abb - aab = e_{\bar{1}} e_{1} - e_{\bar{1}} e_{\bar{1}} - e_{2} + \hbar e_{\bar{1}}\,,   \end{align*}
        which gives the relation 
        \begin{align*}
            0 = (1-q) \zeta_q(abb) - \zeta_q(aab) = \zeta_q(\bar{1}, 1) - \zeta_q(\bar{1},\bar{1}) - \zeta_q(2) + (1-q) \zeta_q(\bar{1})\,.
        \end{align*}
        \item Exercise \ref{fex13}: $w=e_{\bar{1}}$, $v=e_{\hat{2}}$.
    \item For $w=v=e_{\hat{2}}=aab$ we get 
    \begin{align*}
           e_{\hat{2}} \shq e_{\hat{2}} -    e_{\hat{2}} \ast_q e_{\hat{2}} =  4 e_{\hat{3}} e_{\hat{1}}  + 2 e_{\hat{2}} e_{\hat{2}} + 2 e_{\hat{3}}e_{1}-2 e_{\hat{3}}e_{\bar{1}} + \hbar (e_{\hat{2}} e_{1}  + 3e_{\hat{2}} e_{\bar{1}} ) - \left( 2 e_{\hat{2}}e_{\hat{2}} + e_{\hat{4}}\right)\,.
    \end{align*}
\end{enumerate}
\end{exs}

For some $w,v \in \HHH^0$ (but not all!) the limit $q\rightarrow 1$ of the individual terms in $\zeta_q(w \shq v - w \ast_q v)$ exists. Notice that the relations among $\zeta$ obtained from this are more than just the finite double shuffle relations in Proposition \ref{prop:fds}. In particular, we obtain the relation $\zeta(3) = \zeta(2,1)$. \\

It turns out that the relations in Corollary \ref{cor:qfds} do not give all $\CC$-linear relations among $\zeta_q$. It was first observed in \cite{Tak1}, that there is another family of linear relations among $\zeta_q$, called the resummation duality. The same family of linear relations was also discussed for other $q$-analogues, e.g. in \cite{Ba7} these are called partition relations and in \cite{Bri} the name $SZ$-duality is used. In fact, from an algebraic point of view, these dualities are similarly defined as for multiple zeta values. Example \ref{ex:qdsh} (i) shows that $\hbar abb - aab \in \ker \zeta_q$. In general, we define the anti-automorphism $\sigma: \HHH \rightarrow \HHH$ by $\sigma(1)=1$, $\sigma(a) = \hbar b$ and $\sigma(b) = \hbar^{-1} a$, e.g. $\sigma(aab) = \hbar abb$. With this we have the following:

\begin{proposition}\label{prop:sigmainvariance} For all $w \in \HHH^0$ we have 
\begin{align}\label{eq:sigmainvariance}
    \zeta_q(w) = \zeta_q(\sigma(w))\,.
\end{align}
\end{proposition}
\begin{proof}
This is a consequence of Theorem 4 in \cite{Tak1}.
\end{proof}

Even though the definitions of $\sigma$ and $\tau$ are quite similar, the relations obtained are essentially different.

Finally, we mention the following conjecture, which can be seen as the $q$-analogue version of Conjecture \ref{conj:extdsh}, which was first indirectly mentioned in \cite{Tak1} and then later explicitly stated in the thesis of the present author.

\begin{conjecture} \label{conj:qrelations1}All $\CC$-linear relations among $\zeta_q$ can be obtained by combining Corollary \ref{cor:qfds} and Proposition \ref{prop:sigmainvariance}.
\end{conjecture}

Another nice property of the involution $\sigma$ is that it translates the $q$-harmonic product $\ast_q$ into the $q$-shuffle product $\shq$.

\begin{proposition}For all $w,v \in \HHH^0$ we have
      \begin{align}
          w \shq v = \sigma( \sigma(w) \ast_q \sigma(v) )\,.
      \end{align}
\end{proposition}
\begin{proof}
A proof of a more general statement can be found in \cite[Corollary 3.22]{Sa}.
\end{proof}

With this Conjecture \ref{conj:qrelations1} has the following reinterpretation.

\begin{conjecture} \label{conj:qrelations2}All algebraic/linear relations among $\zeta_q$ are a consequence of the $q$-harmonic product $\ast_q$ and the $\sigma$-invariance \eqref{eq:sigmainvariance}. 
\end{conjecture}

\vspace{1cm}
\begin{center}
\bf{ \Large {\color{qmzvlinecol}\ding{94}}~ Exercises ~{\color{qmzvlinecol}\ding{94}}}\end{center}

The following is a collection of exercises intended to deepen the reader's understanding of this chapter.

\begin{exer} \label{ex6}
\begin{enumerate}[\textup{(}i\textup{)}] 
    \item Calculate \eqref{eq:z32itint} by hand, i.e. show 
        \begin{align*}
    \zeta(2,3) = 	\int_{0}^1 \frac{dt_1}{t_1} \int_{0}^{t_1}  \frac{dt_2}{1-t_2}  \int_{0}^{t_2} \frac{dt_3}{t_3} \int_{0}^{t_3} \frac{dt_4}{t_4} \int_{0}^{t_4}  \frac{dt_5}{1-t_5}\,.
    \end{align*}
    without using Corollary \ref{cor:mzvitint}.
    \item Prove Proposition \ref{prop:zmstuffle}, i.e. show that for any $w,u \in \HH^1$ and $M \geq 1$ we have
    \[ \zeta_M(w) \zeta_M(u) = \zeta_M(w \ast u)\,.\]
\end{enumerate}		
\end{exer}

\vspace{0.5cm}
\begin{exer}\label{ex:hw21}
\begin{enumerate}[(i)]
\item Calculate $z_{a}z_{b} \ast z_{c} z_{d}$ for $a,b,c,d\geq 1$.
    \item Show that $\HH^1_\ast := (\HH^1,\ast)$ is a commutative $\Q$-algebra.
\item Show that the $\Q$-linear map $d: \HH^1 \longrightarrow \HH^1$ defined on words by\footnote{Notice that for $r=1$ we view $z_{k_2}\cdots z_{k_r}$ as the empty word, i.e. we define $d(z_1)=1$.}
\begin{align*}
    d: z_{k_1}\cdots z_{k_r} &\longmapsto \begin{cases}
        z_{k_2}\cdots z_{k_r} &, \text{ if } k_1=1 ,\\
        0&, \text{ else}
    \end{cases}
\end{align*}
is a derivation on $\HH^1_\ast$.
\item Define $\HH^0 := \ker(d)$. Show that for every $w \in \HH^1$ there exist unique $w_j \in \HH^0$ such that 
\begin{align*}
w = \sum_{j=0}^l w_j \ast z_1^{\ast j}.\end{align*}In other words, show that $\HH^1_\ast=\HH^0_\ast[z_1]$, where $\HH^0_\ast := (\HH^0,\ast)$.
\item For $w=z_1 z_1 z_3 z_2$ determine the $w_j \in \HH^0$ in (iv).
\end{enumerate}
\end{exer}

\vspace{0.5cm}

\begin{exer} \label{ex:hw22}
\begin{enumerate}[(i)]
    \item Prove Proposition \ref{prop:smintermsofhm}, i.e., 
show that for all $m\geq 1$ and $k_1,\dots,k_r \geq 1$ we have 
\begin{align*}
    S_m(k_1,\dots,k_r)  = \sum_{j=0}^r (-1)^{k_1+\dots + k_j} H_m(k_j,k_{j-1},\dots,k_1) H_m(k_{j+1},\dots,k_r).
\end{align*}
\item Show that for a fixed $m\geq 1$ the $\Q$-linear map defined on a word by
\begin{align*}
  h_m: z_{k_1}\cdots z_{k_r} \longmapsto  H_m(k_1,\dots,k_r) 
\end{align*}
is a $\Q$-algebra homomorphism $h_m: \HH^1_\ast \rightarrow \Q$.
\item Use the results from Exercise \ref{ex:hw21} to show that for any $k_1,\dots,k_r \geq 1$ we have 
\begin{align*}
   \sum_{j=0}^r (-1)^{k_1+\dots + k_j} z_{k_j} z_{k_{j-1}} \cdots z_{k_1} \ast z_{k_{j+1}} \cdots z_{k_r} \in \HH^0.
\end{align*}
Conclude, by using (ii), that $\lim_{m\rightarrow \infty} S_m(k_1,\dots,k_r) \in \mz$ for all $k_1,\dots,k_r \geq 1$.
\end{enumerate} 
\end{exer}
\vspace{0.5cm}

\begin{exer}\label{fex8}
Show that the quasi-shuffle product $\qsh$ defined in \eqref{eq:qshdef} is associative. 
\end{exer}

\vspace{0.5cm}

\begin{exer} \label{ex7}
    \begin{enumerate}[\textup{(}i\textup{)}] 
        \item Prove Proposition \ref{prop:z31z44}, i.e. show that for $n\geq 1$ we have 
        \begin{align*}
        4^n \zeta(\{3,1\}^n) = \zeta(\{4\}^n)\,.
        \end{align*}
        \item Show \eqref{eq:222}, i.e. show that we have for $n\geq 1$
            \begin{align*}
            \zeta(\{2\}^n) = \frac{\pi^{2n}}{(2n+1)!} \,.
            \end{align*}
            (Hint: Calculate the coefficient of $x^{2n}$ in \eqref{eq:sinproduct}.)
    \end{enumerate}
\end{exer}				

\vspace{0.5cm}

\begin{exer} \label{ex8} Prove Proposition \ref{prop:qmprop}, i.e. show the following:
\begin{enumerate}[\textup{(}i\textup{)}]  \item The space $\qmf$ is closed under $q\frac{d}{dq}$.
\item We have 
\begin{align*}
\qmf = \Q[\g(2),\g'(2), \g''(2) ]\,.
\end{align*}
where $\g'$ denotes the derivative with respect to $q\frac{d}{dq}$.
\item Let $\kk=(k_1,\dots,k_r)$ be an index with $k_1,\dots,k_r \geq 2$ even. Then we have 
\begin{align*}
\g^{\text{sym}}(\kk) := \sum_{\sigma \in S_r} \g(k_{\sigma(1)},\dots,k_{\sigma(r)}) \in \qmf\,,
\end{align*}
where $S_r$ denotes the set of all permutations of $\{1,\dots,r\}$.
\item We have 
\begin{align*}
\qmf = \Big \langle  \g^{\text{sym}}(k_1,\dots,k_r)  \,\,\big|\,\, r\ge 0 ,\,k_1\geq k_2\geq\dots \geq k_r\geq 2 \text{ even}
\Big\rangle_\Q\,,
\end{align*} 
where we set $	\g^{\text{sym}}(\emptyset)=1$.
\end{enumerate}
\end{exer}			

\vspace{0.5cm}

\begin{exer} \label{ex9} Let $\bullet \in \{\sh, \ast\}$.
    \begin{enumerate}[\textup{(}i\textup{)}] 
        \item Calculate $\zeta^\bullet(1,2,1)$.
        \item Show that 
        \begin{align*}
        \zeta^\bullet(z_2 \sh z_1 z_1 - z_2 \ast z_1 z_1) =0\,,
        \end{align*} 
        by using the finite double shuffle relations and/or the duality relation and/or Euler's formula.
    \end{enumerate}
    
\end{exer}

\vspace{0.5cm}

\begin{exer}\label{ex10}
\begin{enumerate}[\textup{(}i\textup{)}] 
\item Show that for any admissible index $\kk$ we have $\ds(z_1, z_\kk) = z_1 \sh z_\kk - z_1 \ast z_\kk \in \HH^0$.
\item Prove Hoffman's relation (Proposition \ref{prop:hoffmanrel}) by using the extended double shuffle relations.
\end{enumerate}
\end{exer}

\vspace{0.5cm}

\begin{exer}\label{fex9}\begin{enumerate}[(i)]
    \item Show that for any $k\geq 1$ and $m\geq 1$ we have
	\begin{align*}
 1+ \sum_{n=1}^{\infty} H_m(\overbrace{k,\dots,k}^n) X^n = 	\exp\left( \sum_{n=1}^{\infty} (-1)^{n-1} H_m(nk) \frac{X^n}{n} \right)\,.
	\end{align*}
	\item Prove that for $k\geq 1$ we have $\za(k,\dots,k) =0$ and $\zeta(2k,\dots,2k) \in \Q[\pi^2]$.
\end{enumerate} 
\end{exer}

\vspace{0.5cm}

\begin{exer}\label{fex10} Show that $\HH^1_\ast = \HH_\ast^0[z_1]$.
\end{exer}

\vspace{0.5cm}

\begin{exer}\label{fex11} Prove Proposition \ref{prop:qpolylogab}, i.e. show that for any $w\in \HHH^1$ we have $\li^q_w = w 1$.
\end{exer}

\vspace{0.5cm}

\begin{exer}\label{fex12} Show that for $f,g \in z \Q[[q,z]]$ we have 
\begin{align*}
    (af) (ag) &= a ( (af) g + f(ag) + \hbar fg),\\
    (bf) g &= f (bg) = b (fg)\,.
\end{align*}
\end{exer}

\vspace{0.5cm}

\begin{exer}\label{fex13}
\begin{enumerate}[(i)]
    \item Determine $\zeta_q(e_{\bar{1}} \shq e_{\hat{2}} - e_{\bar{1}} \ast_q  e_{\hat{2}} )$ and consider the limit $q\rightarrow 1$ if possible. 
    \item Calculate $\sigma( \sigma(e_{\hat{2}}) \ast_q \sigma(e_{\hat{2}}))$ and compare it with $e_{\hat{2}} \shq e_{\hat{2}}$. 
\end{enumerate}
\end{exer}

%% file: chap_FiniteSymMZV.tex
\chapter{Finite \& Symmetric MZVs} \label{sec:fmzv}
In this chapter, we want to come back to the finite and symmetric multiple zeta values, which we introduced in Section \ref{subsec:fmzv}. We will first prove the linear shuffle relations for finite multiple zeta values and then have a closer look at symmetric multiple zeta values. The relations among finite multiple zeta values coming from modular forms will be studied later in \Cref{sec:fmzvmodularforms}, since we need the period polynomials of Chapter~\ref{sec:mdandmzv} for them. At the end, we want to make the BTT philosophy of Section \ref{subsec:bttphilo} precise. For this, we introduce the $\mathcal Q$-multiple zeta values of Takeyama and Tasaka \cite{TT}.

\section{Linear shuffle relations for finite MZVs}
Using the notation introduced in Section \ref{subsec:fmzv}, we can view the finite multiple zeta values as a $\Q$-linear map $\za: \HH^1 \rightarrow \mza$. By Lemma \ref{lem:fmtoFM} we see that $\za$ is an algebra homomorphism from $\HH^1_\ast$ to $\mza$, i.e. for any $w,v \in \HH^1$ we have
\begin{align*}
    \za(w) \za(v) = \za(w \ast v)\,.
\end{align*}
So we see that the product of two finite multiple zeta values can be expressed in the same way as multiple zeta values. For these we saw that we also have the shuffle product formula $\zeta(w) \zeta(v) = \zeta(w \shuffle v)$ for $w,v \in \HH^0$. We expect\footnote{Since we cannot prove that $\za(w) \neq 0$ for any non-empty word $w \in \HH^1$, we cannot say that the shuffle product formula is not satisfied.} that this formula is not true for finite multiple zeta values. Instead we have the following \emph{linear shuffle relations} for finite multiple zeta values. For this recall that for $w=z_{k_1}\dots z_{k_r} \in \HH^1$ we denote its reverse by $\overline{w} = z_{k_r}\dots z_{k_1}$. We can view the reverse as a $\Q$-linear map $\HH^1 \rightarrow \HH^1$ by extending it linearly.

\begin{theorem}[Kaneko--Zagier {\cite[Theorem~2]{KZ}}] \label{thm:linearshufflefmzv}For all $w,v\in \HH^1$ we have 
\begin{align*}
    \za(w \sh v) = (-1)^{\wt(w)} \za(\overline{w} v)\,,
\end{align*}
where $\overline{w} v$ denotes the concatenation of $\overline{w}$ and $v$.
\end{theorem}
\begin{proof}
For $n\geq 0$ we write $c_n(\sum_{n\geq 0} a_n z^n) = a_n$ to denote the $n$-th coefficient of a power series in $z$. Then we have for an index $\kk = (k_1,\dots,k_r) \in \Z_{\geq 1}^r$ and a prime $p$
\begin{align*}
    H_{p-1}(\kk) = \sum_{p > m_1 > \dots > m_r>0} \frac{1}{m_1^{k_1} \dots m_r^{k_r}} = \sum_{p > m > 0} c_m\left( \underbrace{\sum_{m_1 > \dots > m_r>0} \frac{z^{m_1}}{m_1^{k_1} \dots m_r^{k_r}}}_{\li_\kk(z)} \right)\,.
\end{align*}
In general, we have for $w \in \HH^1$
\begin{align*}
    H_{p-1}(w) = \sum_{p > m > 0}c_m(\li_w(z))\,.
\end{align*}
Using Proposition \ref{prop:liisalghom} ($\li_{w \shuffle v}(z) = \li_{w}(z)\li_{v}(z)$) we get for $w=z_{k_1}\dots z_{k_r}$ and $v=z_{l_1}\dots z_{l_s}$
\begin{align*}
    H_{p-1}(w \shuffle v) &=  \sum_{p > m > 0}c_m(\overbrace{\li_{w \shuffle v}(z)}^{\li_{w}(z)\li_{v}(z)})= \sum_{p > m > 0} \sum_{\substack{p>i,j>0\\i+j=m}} c_i(\li_{w}(z)) c_j(\li_{v}(z)) \\
    &= \sum_{\substack{p > i,j > 0 \\ p > i+j >0}} \sum_{i > m_2 > \dots > m_r >0 } \frac{1}{i^{k_1} m_2^{k_2}\cdots m_r^{k_r}} \sum_{j > n_2 > \dots > n_s >0 } \frac{1}{j^{l_1} n_2^{l_2}\cdots n_s^{l_s}}\,.
\end{align*}
If we consider this sum modulo $p$, reversing the order of the first sum, and doing the change of variables $i' = p-i$, $m'_t = p - m_t$ ($t=2,\dots,r$), we obtain
\begin{align*}
      H_{p-1}(w \shuffle v) &\equiv   \sum_{\substack{p > i,j > 0 \\ p > i+j >0}} \sum_{i > m_2 > \dots > m_r >0 } \frac{(-1)^{k_1+\dots+k_r}}{(p-i)^{k_1} (p-m_2)^{k_2}\dots (p-m_r)^{k_r}} \\
      &\qquad{}\cdot \sum_{j > n_2 > \dots > n_s >0 } \frac{1}{j^{l_1} n_2^{l_2}\dots n_s^{l_s}} \quad \mod p\\
      &\equiv \sum_{p > m'_{r} > \dots > m'_2 > i' > j > n_2 > \dots > n_s > 0} \frac{(-1)^{k_1+\dots+k_r}}{{m'}_r^{k_r} \cdots {m'}_2^{k_2} {i'}^{k_1} j^{l_1} n_2^{l_2} \cdots n_s^{l_s} } \quad \mod p\\
      &\equiv (-1)^{\wt(w)} H_{p-1}(\overline{w}v)\,\quad \mod p,
\end{align*}
which implies the statement in the theorem. 
Notice that \cite{KZ} uses the opposite convention for the order of an index. Reversing all indices in its Theorem~2 gives exactly the formula above.
\end{proof}

\begin{corollary}
For all $w\in \HH^1$ we have $\za(w) = (-1)^{\wt(w)} \za(\overline{w})$.
\end{corollary}
\begin{proof}
This is the reversal formula \eqref{eq:reversalformula}, but now can also be seen as a special case of the above Theorem, by using $v=1$.
\end{proof}
\begin{conjecture}[Kaneko--Zagier {\cite[Conjectures~3 and~6]{KZ}}]\label{conj:fmzvallrelations}
All algebraic and linear relations over $\Q$ among finite multiple zeta values can be deduced from $(w,v \in \HH^1)$
\begin{align}\tag{$*$}\label{eq:conjast}
     \za(w * v) &= \za(w) \za(v), \\
     \tag{$\shuffle$}\label{eq:conjshuffle}
    \za(w \shuffle v) &= (-1)^{\wt(w)}\za(\overline{w}v)
\end{align}
In other words: If we define $K$ to be the ideal in $\HH^1_*$ generated by $w\shuffle v - (-1)^{\wt(w)} \overline{w}v$ for all $w,v \in \HH^1$, then we expect that $\ker(\za) = K$.
\end{conjecture}

The universal quotient by this ideal will be studied in
Section~\ref{sec:formalfinitemzv}. There we will also construct a formal
version of the Kaneko--Zagier map.

In particular, all relations we have proved so far for finite multiple zeta values should be a consequence of \eqref{eq:conjast} and \eqref{eq:conjshuffle} (Exercise \ref{ffex14}: Show that $\za(k)=0$ for all $k\geq 1$). Since $\za(k)=0$, equation \eqref{eq:conjast} reduces in the special case $v=z_k$ to $\za(w * z_k)=0$. It seems that this equation is enough to reduce all linear relations together with \eqref{eq:conjshuffle}:
\begin{conjecture}[Kaneko--Zagier {\cite[Conjecture~4]{KZ}}]
All $\Q$-linear relations among finite multiple zeta values can be deduced from \eqref{eq:conjshuffle} and $\za(w * z_k)=0$ for $w\in \HH^1$ and $k\geq 1$.
\end{conjecture}

\section{Symmetric MZVs}
Recall from \cite[\S8, (87)]{KZ} that for
an index $\kk=(k_1,\ldots,k_r)$ and $\bullet \in \{\shuffle, \ast \}$ the {\bf $\bullet$-symmetric multiple zeta value}  is defined by
\begin{align*}
      \zs^\bullet(\kk)=  \zs^\bullet(k_1,\dots,k_r) = \sum_{j=0}^r (-1)^{k_1+\dots+k_j} \zeta^\bullet(k_j,k_{j-1},\dots,k_1;T) \zeta^\bullet(k_{j+1},\dots,k_r;T)\,.
\end{align*}
In this section, we want to sketch the proofs of Proposition \ref{prop:symmzvstatements}. Its parts (i) and (iii) are Theorem~3 in \cite{KZ}, and part (ii) is Proposition~9. In other words, we want to show the following:
\begin{enumerate}[(i)]
    \item We have $ \zs^\bullet(\kk) \in \mz$ for $\bullet \in \{\shuffle, \ast \}$.\\
    (i.e. the $\bullet$-symmetric multiple zeta values are independent of $T$).
    \item For all indices $\kk$ we have 
\begin{align*}
    \lim_{m\rightarrow \infty }S_m(\kk)  = \zs^\ast(\kk)\,.
\end{align*}
    \item For all indices $\kk$ we have 
\begin{align*}
    \zs^\ast(\kk) \equiv  \zs^\shuffle(\kk) \mod \pi^2 \mz\,.
\end{align*}
\end{enumerate} 
For this we will make use of the following Lemma.

\begin{lemma} \label{lem:Glem}For an admissible index    $\kk$ and $\bullet \in \{\shuffle, \ast \}$  we define the following generating series
\begin{align}\label{eq:defG}
    G^\bullet(\kk;X,T) := \sum_{j=0}^\infty \zeta^\bullet(\{1\}^j,\kk ; T) X^j\,
\end{align}
and write in the special case $T=0$ just $ G^\bullet(\kk;X):=  G^\bullet(\kk;X,0)$.
    \begin{enumerate}[(i)]
        \item  For $\bullet \in \{\shuffle, \ast \}$ we have 
        \begin{align*}
             G^\bullet(\kk;X,T) = e^{XT} G^\bullet(\kk;X)\,.
        \end{align*}
        \item We have 
        \begin{align*}
          G^\ast(\kk;X) = \Gamma(1+X)^{-1}e^{-\gamma X}   G^\shuffle(\kk;X)\,.
        \end{align*}
    \end{enumerate}
\end{lemma}
\begin{proof}
For (i), write $z_\kk=z_{k_1}\cdots z_{k_r}$. Since $\kk$ is admissible,
$z_\kk\in\HH^0$, and in $\HH^1[[X]]$ we have
\begin{align*}
 \frac{1}{1-z_1X}z_\kk=\sum_{j=0}^\infty z_1^jz_\kk X^j,
\end{align*}
where the powers on the right are taken with respect to concatenation. Proposition~10
in \cite{IKZ} states, for $\bullet\in\{\shuffle,\ast\}$, that
\begin{align*}
 \reg^T_\bullet\left(\frac{1}{1-z_1X}z_\kk\right)
 =e^{TX}\reg_\bullet\left(\frac{1}{1-z_1X}z_\kk\right).
\end{align*}
Applying $\zeta$ coefficientwise gives
$G^\bullet(\kk;X,T)=e^{TX}G^\bullet(\kk;X)$, which proves (i).

For (ii) recall that we defined the linear map $\rho: \R[T] \rightarrow \R[T]$ in \eqref{eq:rhodef} by $\rho(e^{Tu}) = A(u) e^{Tu}$, where 
\begin{align*}
    A(u) = \exp\left( \sum_{n=2}^\infty \frac{(-1)^n}{n} \zeta(n) u^n\right) = \Gamma(1+u) e^{\gamma u}\,.
\end{align*}
Here the second equation follows by considering the logarithmic derivative for the Gamma function and $\gamma = \lim_{n\rightarrow \infty}\left(\sum_{m=1}^n \frac{1}{m} - \log(n)\right)$ denotes Euler's constant. By Theorem \ref{thm:comparezastzsh} we have $\zeta^\shuffle = \rho \circ \zeta^\ast$. Applying $\rho$ to (i) with $\bullet = \ast$ and setting $T=0$ yields (ii), since $\rho$ is $\R$-linear and $\zeta^\ast(\{1\}^j, \kk; 0)\in \R$. 
\end{proof}

For the proof of  Proposition \ref{prop:symmzvstatements} (i) and (iii) we will consider certain partial sums of the terms appearing in the definition of $\zs$. First notice that the statements in (i) and (iii) are trivial for indices $\kk = (k_1,\dots,k_r)$ with $k_1,\dots,k_r\geq2$, since in this case no regularization is necessary and we have $\zs^\shuffle(\kk)=\zs^\ast(\kk)\in \R$. Now suppose that $\kk =(k_1,\dots,k_r)= (\overline{\kl}, \{1\}^h, {\bf m})$ for $h\geq 1$ and some admissible indices $\kl$ and ${\bf m}$. Then we claim that the statements in Proposition \ref{prop:symmzvstatements} (i) and (iii) are already true for the following partial sum of $\zeta^\bullet(\kk)$
\begin{align}\label{eq:partialsum}
 \sum_{i=0}^h (-1)^{i+\wt(\kl)} \zeta^\bullet(\{1\}^i,\kl;T) \zeta^\bullet(\{1\}^{h-i},{\bf m};T)\,.
\end{align}
Notice that \eqref{eq:partialsum} is the coefficient of $X^h$ in 
\begin{align*}
    P^\bullet(\kl,{\bf m}; X,T) := (-1)^{\wt(\kl)}G^\bullet(\kl;-X,T)G^\bullet({\bf m};X,T)\,.
\end{align*}
Now we are ready to give the proofs of  Proposition \ref{prop:symmzvstatements}.
\begin{proof}[Proof of Proposition \ref{prop:symmzvstatements} (i)]
By Lemma \ref{lem:Glem} (i) we have 
\begin{align*}
       P^\bullet(\kl,{\bf m}; X,T)  = (-1)^{\wt(\kl)} e^{-XT}G^\bullet(\kl;-X) e^{XT}G^\bullet({\bf m};X) = (-1)^{\wt(\kl)} G^\bullet(\kl;-X) G^\bullet({\bf m};X) 
\end{align*}
and therefore \eqref{eq:partialsum} is independent of $T$ and so is $\zeta^\bullet_S(\kk)$.
\end{proof}
\begin{proof}[Proof of Proposition \ref{prop:symmzvstatements} (iii)] Since $P$ is independent of $T$ we will just write $ P^\bullet(\kl,{\bf m}; X)$ in the following. By Lemma \ref{lem:Glem} (i),(ii) we obtain 
\begin{align*}
 P^\ast(\kl,{\bf m}; X) = \frac{1}{\Gamma(1+X)\Gamma(1-X)} P^\shuffle(\kl,{\bf m}; X).
\end{align*}
Using the well-known identity for the sine
\begin{align*}
    \frac{1}{\Gamma(1+X)\Gamma(1-X)} &= \frac{\sin(\pi X)}{\pi X} = \sum_{n=0}^\infty (-1)^n \zeta(\{2\}^n) X^{2n}= \sum_{n=0}^\infty (-1)^n \frac{\pi^{2n} X^{2n}}{(2n+1)!} = 1 - \frac{\pi^2}{6}X^2+\dots \,
\end{align*}
we get $ P^\ast(\kl,{\bf m}; X) -  P^\shuffle(\kl,{\bf m}; X) \in \pi^2 \mz[[X]]$\,. This shows that \eqref{eq:partialsum} for $\bullet = \ast$ and $\bullet = \shuffle$ are the same modulo $\pi^2 \mz$ and therefore the same also holds for $\zs^\bullet(\kk)$.
\end{proof}

As a direct consequence of the proof, we get:

\begin{corollary} If there are no consecutive $1$'s in $\kk$ we have $\zs^\ast(\kk) = \zs^\shuffle(\kk)$.
\end{corollary}
This is also \cite[Corollary to Theorem~3]{KZ}.

\begin{proof}[Proof of Proposition \ref{prop:symmzvstatements} (ii)]
Put $T_m=\log(m+1)+\gamma$. The asymptotic formula used in the proof of
Theorem~\ref{thm:comparezastzsh} gives, for every index $\kl$,
\begin{align*}
 H_m(\kl)=\zeta^\ast(\kl;T_m)
 +\operatorname{O}\big(m^{-1}\log^J(m+1)\big)
\end{align*}
for some $J$ depending on $\kl$. Here the shift by one comes from the convention that
$H_m$ has the weak upper bound $m$, whereas $\zeta_M$ in the proof of Theorem
\ref{thm:comparezastzsh} has the strict upper bound $M$. Inserting this asymptotic
into \eqref{eq:smashm}, and using that the regularized values are polynomials in $T_m$,
we obtain
\begin{align*}
 S_m(\kk)
 &=\sum_{j=0}^r(-1)^{k_1+\dots+k_j}
 \zeta^\ast(k_j,\dots,k_1;T_m)
 \zeta^\ast(k_{j+1},\dots,k_r;T_m)+o(1)\\
 &=\zs^\ast(\kk)+o(1).
\end{align*}
The last equality follows from the already proven independence of $T$. Taking the limit
proves the claim. This argument is also given in \cite[Proposition~9]{KZ}.
\end{proof}
\vspace{0.2cm}

\begin{theorem}[Kaneko--Zagier {\cite[Theorem~5]{KZ}}]\label{thm:linearshufflesmzv}
 For all words $w,v\in \HH^1$ we have 
\begin{align*}
    \zs^\shuffle(w \sh v) = (-1)^{\wt(w)} \zs^\shuffle(\overline{w} v)\,.
\end{align*}
\end{theorem}
\begin{proof}
On the space of indices of depth $r$, define the linear maps
\begin{align*}
 \mathcal S_i(z_{k_1}\cdots z_{k_r})
 =(-1)^{k_1+\dots+k_i}
 (z_{k_i}\cdots z_{k_1})\sh(z_{k_{i+1}}\cdots z_{k_r})
 \qquad(0\leq i\leq r)
\end{align*}
and put $\mathcal S=\sum_{i=0}^r\mathcal S_i$. Since $\zeta^\shuffle$ is a
shuffle algebra homomorphism, the definition of the symmetric value can be written as
\begin{align*}
 \zs^\shuffle(u)=\zeta^\shuffle(\mathcal S(u);T).
\end{align*}
The combinatorial identity
\begin{align}\label{eq:symmetricshuffleoperator}
 \mathcal S\circ\mathcal S_i=\mathcal S\qquad(0\leq i\leq r)
\end{align}
holds already in the formal space of indices (see \cite[Proposition~10]{KZ}, or \cite[Proposition~9.7]{K3}). One proof of \eqref{eq:symmetricshuffleoperator}
uses the generating series of all depth-$r$ indices. There $\mathcal S_i$ is represented by the signed
sum of the permutations which first increase and then decrease. The alternating sum
defining $\mathcal S$ is fixed by each of these elements. This is the classical
descent-algebra identity used in the cited proofs.

Now let $i=\dep(w)$ and apply \eqref{eq:symmetricshuffleoperator} to
$\overline wv$. We have
\begin{align*}
 \mathcal S_i(\overline wv)=(-1)^{\wt(w)}(w\sh v),
\end{align*}
and hence
\begin{align*}
 \zs^\shuffle(\overline wv)
 &=\zeta^\shuffle\big(\mathcal S(\mathcal S_i(\overline wv));T\big)
 =(-1)^{\wt(w)}\zs^\shuffle(w\sh v).
\end{align*}
This gives the statement. Notice that the references use the opposite convention
for the order of an index. Reversing all indices gives exactly the formula above.
\end{proof}

\section{BTT philosophy and \texorpdfstring{$\mathcal{Q}$}{Q}-multiple zeta values} \label{sec:Qmzv}

We will now come back to the multiple harmonic sums at roots of unity. For $\kk = (k_1,\dots,k_r) \in \Z^r$ and $m\geq 1$ these were defined by
\begin{align}
 H_m(\kk;q) =    H_m(k_1,\dots,k_r;q) = \sum_{m\geq m_1>\dots>m_r>0} \frac{q^{(k_1-1)m_1} \cdots q^{(k_r-1)m_r}}{[m_1]_q^{k_1} \cdots [m_r]_q^{k_r} } \in \Q[[q]]\,.
\end{align}
As explained in Section \ref{subsec:bttphilo}, we are interested in the values $H_{n-1}({\bf k};\zeta_n)  \in \Q(\zeta_n)$, since these can be related to symmetric and finite multiple zeta values. The connection to finite multiple zeta values was given by Theorem \ref{thm:bttfmzv}, which stated that for any primitive $p$-th root of unity $\zeta_p$, we have
	\[ ( H_{p-1}(\kk; \zeta_p) \mod \mathfrak{p})_p = \zeta_{\mathcal{A}} (\kk) \,,\]
where $\mathfrak{p}= (1-\zeta_p)$ is the prime ideal of $\Z[\zeta_p]$ generated by $1-\zeta_p$. In the next section we will prove the connection to symmetric multiple zeta values and then define $\mathcal{Q}$-multiple zeta values. 

\subsection{The connection to symmetric multiple zeta values}

 In this section, we want to sketch the proof of Theorem \ref{thm:xilimit}, i.e. we want to show that for any index set $\kk=(k_{1}, \ldots , k_{r})$ the limit $\xi(\kk) := \lim\limits_{n \rightarrow \infty}    H_{n-1}(\kk ; e^{\frac{2\pi i}{n}})$ exists and we have
\begin{align*}
\xi(\kk) = \lim\limits_{n \rightarrow \infty}    H_{n-1}(\kk ; e^{\frac{2\pi i}{n}})=\sum_{a=0}^{r}(-1)^{k_1+\cdots+k_a}
\zeta^*\big(k_{a}, k_{a-1}, \ldots , k_{1} ;\frac{\pi i}{2}\big)
\zeta^*\big(k_{a+1}, k_{a+2}, \ldots , k_{r};-\frac{\pi i}{2}\big)\,.
\end{align*}

In particular, this shows
\begin{align*}
     {\rm Re} \left(\xi(\kk) \right) \equiv \zeta_{\mathcal{S}} (\kk) \mod \pi^2 \mathcal{Z}\,.
\end{align*} 
The proof of Theorem \ref{thm:xilimit} was first given in \cite{BTT1}, which we will sketch in the following. Hirose's refined symmetric multiple zeta values give another interpretation of these limits \cite[Remark~13]{Hi1}, after reversing the order of the index.

\begin{lemma}\label{lem:tindepend}
For any index $\kk=(k_1,\dots,k_r)$ the polynomial 
\begin{align*}
R_\kk(X;T) = \sum_{j=0}^r (-1)^{k_1+\dots+k_j} \zeta^\ast(k_j,k_{j-1},\dots,k_1;T+X) \zeta^\ast(k_{j+1},\dots,k_r;T-X)
\end{align*}
does not depend on $T$, i.e.  $R_\kk(X;T)=R_\kk(X;0)$.
\end{lemma}
\begin{proof}
This is a generalization of the $\ast$-part of Proposition \ref{prop:symmzvstatements} (i) and can be proven in a similar way. Exercise \ref{fex15}.
\end{proof}

We rewrite the value $H_{n-1}(\kk; e^{2\pi i/n})$. 
Let $n$ be a positive integer. 
When $q=e^{2\pi i/n}$ we see that 
\begin{align*}
\frac{1}{[m]_q} = \frac{1-q}{1-q^{m}}=e^{-\frac{\pi i}{n}(m-1)}\frac{\sin{\frac{\pi}{n}}}{\sin{\frac{m \pi}{n}}} \qquad 
(n>m > 0).
\end{align*}
Therefore it holds that 
\begin{align*}
H_{n-1}(\kk;e^{\frac{2\pi i}{n}}) =
\left(e^{\frac{\pi i}{n}}\frac{n}{\pi}\sin{\frac{\pi}{n}}\right)^{\mathrm{wt}(\mathbf{k})}
\sum_{n>m_{1}>\cdots >m_{r}>0}
\prod_{j=1}^{r}\frac{e^{\frac{\pi i}{n}(k_{j}-2)m_{j}}}
{\left(\frac{n}{\pi}\sin{\frac{m_{j}\pi}{n}}\right)^{k_{j}}}
\end{align*}
for any non-empty index $\kk=(k_{1}, \ldots , k_{r})$. 
Decompose the set $\{(m_{1}, \ldots , m_{r}) \in \Z^{r} \, | \, n>m_{1}>\cdots >m_{r}>0\}$ into the disjoint union
\begin{align*}
\bigsqcup_{a=0}^{r}
\{(m_{1}, \ldots , m_{r}) \in \Z^{r} \, | \, n>m_{1}>\cdots >m_{a}>\frac{n}{2} \ge m_{a+1}>\cdots >m_{r}>0 \}
\end{align*}
and change the summation variables $m_{j}$ to $n_{j}=n-m_{a+1-j} \, (1\le j \le a)$ and 
$l_{j}=m_{a+j} \, (1 \le j \le r-a)$. 
Then we find that 
\begin{align*}
& 
H_{n-1}(\kk;e^{\frac{2\pi i}{n}}) =
\left(e^{\frac{\pi i}{n}}\frac{n}{\pi}\sin{\frac{\pi}{n}}\right)^{\mathrm{wt}(\mathbf{k})}\\ 
&\times\sum_{a=0}^{r}(-1)^{\sum_{j=1}^{a}k_{j}} 
\sum_{n/2>n_{1}>\cdots >n_{a}>0}
\prod_{j=1}^{a}\frac{e^{-\frac{\pi i}{n}(k_{a+1-j}-2)n_{j}}}
{\left(\frac{n}{\pi}\sin{\frac{n_{j}\pi}{n}}\right)^{k_{a+1-j}}}
\sum_{n/2\ge l_{1}>\cdots >l_{r-a}>0}
\prod_{j=1}^{r-a}\frac{e^{\frac{\pi i}{n}(k_{a+j}-2)l_{j}}}
{\left(\frac{n}{\pi}\sin{\frac{l_{j}\pi}{n}}\right)^{k_{a+j}}}. 
\end{align*}

Motivated by the above expression we introduce the following numbers. 
For an index $\kk=(k_{1}, \ldots , k_{r})$ and a positive integer $n$, we define 
\begin{align*}
A_{n}^{-}(\kk)&=\sum_{n/2>m_{1}>\cdots >m_{r}>0}
\prod_{j=1}^{r}\frac{e^{-\frac{\pi i}{n}(k_{j}-2)m_{j}}}
{\left(\frac{n}{\pi}\sin{\frac{m_{j}\pi}{n}}\right)^{k_{j}}}, \\ 
A_{n}^{+}(\kk)&=\sum_{n/2 \ge m_{1}>\cdots >m_{r}>0}
\prod_{j=1}^{r}\frac{e^{\frac{\pi i}{n}(k_{j}-2)m_{j}}}
{\left(\frac{n}{\pi}\sin{\frac{m_{j}\pi}{n}}\right)^{k_{j}}}. 
\end{align*}
Then we see that 
\begin{align*}
& 
H_{n-1}(\kk;e^{\frac{2\pi i}{n}}) =
\left(e^{\frac{\pi i}{n}}\frac{n}{\pi}\sin{\frac{\pi}{n}}\right)^{\mathrm{wt}(\mathbf{k})} 
\sum_{a=0}^{r}(-1)^{\sum_{j=1}^{a}k_{j}} 
A_{n}^{-}(k_{a}, k_{a-1}, \ldots , k_{1}) 
A_{n}^{+}(k_{a+1}, k_{a+2}, \ldots , k_{r}).  
\end{align*}

Notice that we have 
\renewcommand{\arraystretch}{1.5}
\begin{align*}
A_{n}^{-}(k_{1}, \ldots , k_{r})=\left\{ 
\begin{array}{ll}
\overline{A_{n}^{+}(k_{1}, \ldots , k_{r})} &  (\hbox{$n$: odd}), \\
\overline{A_{n}^{+}(k_{1}, \ldots , k_{r})}+(-\frac{\pi i}{n})^{k_{1}}\, 
\overline{A_{n}^{+}(k_{2}, \ldots , k_{r})} &  
(\hbox{$n$: even}),
\end{array}
\right. 
\end{align*}
where the bar on the right-hand side denotes complex conjugation. We want to understand the behavior of $A_n^+$ for $n\rightarrow \infty$. For this we will first consider the cases where the index $\kk$ is admissible. 

\begin{lemma}[\cite{BTT1} Lemma 2.7]\label{lem:T-asym-admissible}
Let $\kk$ be an admissible index. 
Then it holds that 
\begin{align*}
A_{n}^{+}(\kk)=\zeta(\kk)+O\left(\frac{(\log{n})^{J_{1}(\kk)}}{n}\right) \qquad (n \to +\infty),   
\end{align*} 
where $J_{1}(\kk)$ is a positive integer which depends on $\kk$. 
\end{lemma}

\begin{proof}
Set $\kk=(k_1,\ldots,k_r)$ and define for $k\geq 1$ the function
\begin{align*}
g_k(x)=e^{(k-2)ix} \left(\frac{x}{\sin x}\right)^k \,. 
\end{align*}
Then it holds that 
$\left|A_{n}^{+}(\kk)-\zeta(\kk)\right| \le I_{1}+I_{2}$, where 
\begin{align*}
I_{1}&=\sum_{n/2\ge m_{1}>\cdots >m_{r}>0}
\prod_{j=1}^{r}\frac{1}{m_{j}^{k_{j}}} 
\left| 
\prod_{j=1}^{r}g_{k_j}\left(\frac{m_{j}\pi}{n}\right)-1
\right|, \\ 
I_{2}&=\sum_{m>n/2}\frac{1}{m^{k_{1}}} 
\left(
\sum_{m>m_{2}>\cdots >m_{r}>0}\prod_{j=2}^{r}\frac{1}{m_{j}^{k_{j}}}  
\right). 
\end{align*} 
Since $g_k(x)=1+(k-2)ix + o(x) \ (x\rightarrow +0)$, 
there exists a positive constant $C$ depending on $k$ such that 
$|g_k(m\pi/n)-1| \le Cm/n$
for all integers $m$ and $n$ satisfying $n/2\ge m>0$.
Using the identity 
\begin{align*}
\left(\prod_{j=1}^{r}x_{j}\right)-1=\sum_{a=1}^{r}\, 
(\prod_{j=1}^{a-1}x_{j}) \, (x_{a}-1)
\end{align*} 
and the inequality $0<(\sin x)^{-1}\le \pi/2x$ on the interval $(0,\frac{\pi}{2}]$, we see that 
\begin{align*}
I_{1}&\le \frac{C_{1}}{n} \sum_{a=1}^{r} \sum_{n/2\ge m_{1}>\cdots >m_{r}>0} 
\frac{1}{m_{1}^{k_{1}} \cdots m_{a}^{k_{a}-1} \cdots m_{r}^{k_{r}}} \\ 
&\le \frac{C_{1}}{n} \sum_{a=1}^{r} \sum_{n/2\ge m_{1}>\cdots >m_{r}>0} 
\frac{1}{m_{1}^{k_{1}-1} m_{2}^{k_{2}} \cdots m_{r}^{k_{r}}} \\ 
&=\frac{C_{1} r}{n} \sum_{n/2\ge m>0}\frac{1}{m^{k_{1}-1}} 
\left(
\sum_{m>m_{2}>\cdots >m_{r}>0}\prod_{j=2}^{r}\frac{1}{m_{j}^{k_{j}}}  
\right) 
\end{align*}
for some positive constant $C_{1}$ which depends on $\kk$.
Using the estimation 
\begin{align*}
\sum_{m>m_{2}>\cdots >m_{r}>0}\prod_{j=2}^{r}\frac{1}{m_{j}^{k_{j}}} \le 
\left(\sum_{s=1}^{m-1}\frac{1}{s}\right)^{r-1} \le 
(2\log{m})^{r-1},  
\end{align*}
we get
\begin{align*}
I_{1}+I_{2} \le C_{2} \left(
\frac{1}{n}\sum_{n/2>m>0}\frac{(\log{m})^{r-1}}{m^{k_{1}-1}}+
\sum_{m>n/2}\frac{(\log{m})^{r-1}}{m^{k_{1}}}
\right)
\end{align*}
for some positive constant $C_{2}$ which depends on $\kk$. 
Since $k_{1} \ge 2$ it holds that 
\begin{align*}
\sum_{n/2>m>0}\frac{(\log{m})^{r-1}}{m^{k_{1}-1}}=O((\log{n})^{r}), \quad 
\sum_{m>n/2}\frac{(\log{m})^{r-1}}{m^{k_{1}}}=O\left(\frac{(\log{n})^{r-1}}{n}\right) 
\end{align*}
as $n \to +\infty$. 
This completes the proof. 
\end{proof}

\begin{lemma}[\cite{BTT1} Lemma 2.8]\label{lem:anplus1} We have
\begin{align}
A_{n}^{+}(1)=\log{\left(\frac{n}{\pi}\right)}+\gamma-\frac{\pi i}{2}+O\left(\frac{1}{n}\right) \qquad (n \to +\infty).    
\label{eq:A1-asymptotics}
\end{align}
\end{lemma}
\begin{proof}
From the definition of $A_{n}^{+}(1)$ we see that 
\begin{align*}
A_{n}^{+}(1)=\frac{\pi}{n}\sum_{n/2\ge m>0}\left(
\frac{\cos{\frac{m\pi}{n}}}{\sin{\frac{m\pi}{n}}}-i\right)=
\frac{\pi}{n}\sum_{n/2\ge m>0}
\frac{\cos{\frac{m\pi}{n}}}{\sin{\frac{m\pi}{n}}}-\frac{\pi i}{2}+O\left(\frac{1}{n}\right)
\end{align*} 
as $n \to +\infty$. 
Hence it suffices to show that 
\begin{align}
\frac{\pi}{n}\sum_{n/2\ge m>0}\frac{\cos{\frac{m\pi}{n}}}{\sin{\frac{m\pi}{n}}}=
\log{\frac{n}{\pi}}+\gamma+O\left(\frac{1}{n}\right) \qquad 
(n \to +\infty).  
\label{eq:T1-asym}
\end{align}

Since the function $f(x)=x^{-1}-(\tan{x})^{-1}$ is positive and increasing 
on the interval $(0,\pi)$, we see that 
\begin{align*}
\int_{0}^{\frac{n-1}{2}}f\Big(\frac{\pi x}{n}\Big)\,dx \le 
\sum_{n/2\ge m>0}\left(\frac{n}{\pi}\frac{1}{m}-\frac{\cos{\frac{m\pi}{n}}}{\sin{\frac{m\pi}{n}}}\right) \le 
\int_{1}^{\frac{n}{2}+1}f\Big(\frac{\pi x}{n}\Big)\,dx.
\end{align*}
Set $g(x)=\log{(1+x)}-\log{(\cos{\frac{\pi x}{2}})}$. 
By direct calculation we have 
\begin{align*}
& 
\int_{0}^{\frac{n-1}{2}}f\Big(\frac{\pi x}{n}\Big)\,dx=\frac{n}{\pi}\left( 
g\Big(-\frac{1}{n}\Big)+\log{\frac{\pi}{2}}
\right), \\ 
& 
\int_{1}^{\frac{n}{2}+1}f\Big(\frac{\pi x}{n}\Big)\,dx=\frac{n}{\pi}\left( 
g\Big(\frac{2}{n}\Big)+\log{\Big(\frac{n}{\pi}\sin{\frac{\pi}{n}}\Big)}+\log{\frac{\pi}{2}}
\right)\,.
\end{align*}
Since $g(x)=x+o(x) \,\, (x \to 0)$ and $\log{(x^{-1}\sin{x})}=o(x) \,\, (x \to +0)$, 
there exist positive constants $c_{1}$ and $c_{2}$ such that 
\begin{align*}
\int_{0}^{\frac{n-1}{2}}f\Big(\frac{\pi x}{n}\Big)\,dx \ge -c_{1}+\frac{n}{\pi}\log{\frac{\pi}{2}}, \quad 
\int_{1}^{\frac{n}{2}+1}f\Big(\frac{\pi x}{n}\Big)\,dx \le c_{2}+\frac{n}{\pi}\log{\frac{\pi}{2}}
\end{align*}
for $n \gg 0$. 
Therefore we find that 
\begin{align*}
\frac{\pi}{n} 
\sum_{n/2\ge m>0}\frac{\cos{\frac{m\pi}{n}}}{\sin{\frac{m\pi}{n}}}=
\sum_{n/2\ge m>0}\frac{1}{m}-\log{\frac{\pi}{2}}+O\Big(\frac{1}{n}\Big) \qquad (n \to +\infty). 
\end{align*}
Using the asymptotic expansion 
\begin{align*}
\sum_{n/2\ge m>0}\frac{1}{m}=\log{\frac{n}{2}}+\gamma+O\Big(\frac{1}{n}\Big) \qquad (n \to +\infty),  
\end{align*}
we get the formula \eqref{eq:T1-asym}. 
\end{proof}

\begin{proposition}\label{prop:asym-A+}
For any index $\kk$ it holds that 
\begin{align}
A_{n}^{\pm}(\kk)=\zeta^\ast\left(\kk ;\log{\left(\frac{n}{\pi}\right)}+\gamma\mp \frac{\pi i}{2}\right)+
O\left(\frac{(\log{n})^{J(\kk)}}{n}\right) \qquad (n \to +\infty),  
\label{eq:asymptotics-T}
\end{align} 
where $\gamma$ is Euler's constant and $J(\kk)$ is a positive integer which depends on $\kk$. 
\end{proposition} 
\begin{proof}
Lemma \ref{lem:T-asym-admissible} implies that 
the equality \eqref{eq:asymptotics-T} for $A_{n}^{+}(\kk)$ holds if $\kk$ is admissible. 
Let us prove that it holds also for the index $(\{1\}^{s}, \kk)$ 
with any $s \ge 0$ and any admissible index $\kk$. 

Using the equality 
\begin{align*}
\frac{e^{-\frac{\pi i}{n}m}}{\frac{n}{\pi}\sin{\frac{m\pi}{n}}} 
\frac{e^{\frac{\pi i}{n}(k-2)m}}{\left(\frac{n}{\pi}\sin{\frac{m\pi}{n}}\right)^{k}}=
\frac{e^{\frac{\pi i}{n}(k-1)m}}{\left(\frac{n}{\pi}\sin{\frac{m\pi}{n}}\right)^{k+1}}-\frac{2\pi i}{n}
\frac{e^{\frac{\pi i}{n}(k-2)m}}{\left(\frac{n}{\pi}\sin{\frac{m\pi}{n}}\right)^{k}},  
\end{align*} 
for $k \ge 1$ and $n/2\ge m>0$, we see that 
\begin{align*}
A_{n}^{+}(1)A_{n}^{+}(\{1\}^{s}, \kk)&=(s+1)A_{n}^{+}(\{1\}^{s+1}, \kk) \\ 
&+\sum_{b=1}^{s}\left( 
A_{n}^{+}(\{1\}^{b-1}, 2, \{1\}^{s-b}, \kk)-\frac{2\pi i}{n}A_{n}^{+}(\{1\}^{s}, \kk)
\right) \\ 
&+\sum_{a=1}^{r}\left( 
A_{n}^{+}(\{1\}^{s}, \kk'(a))+A_{n}^{+}(\{1\}^{s}, \kk''(a))-\frac{2\pi i}{n}A_{n}^{+}(\{1\}^{s}, \kk)
\right),  
\end{align*}
where $\kk'(a)$ and $\kk''(a)$ are the indices defined by \begin{align}
\mathbf{k}'(a)=(k_{1}, \ldots , k_{a}+1, \ldots , k_{r}), \quad 
\mathbf{k}''(a)=(k_{1}, \ldots , k_{a}, 1, k_{a+1}, \ldots , k_{r})\,. 
\label{eq:kk-primes}
\end{align}
We obtain the desired equality \eqref{eq:asymptotics-T} by induction on $s$ by the stuffle product and Lemma \ref{lem:anplus1}.
\end{proof}

\begin{proof}[Proof of Theorem \ref{thm:xilimit}]
The statement follows by combining Lemma \ref{lem:tindepend} and Proposition \ref{prop:asym-A+}.
\end{proof}

\subsection{Duality}
In this section we want to give the proof of Theorem \ref{thm:btthoffmandual}, i.e. we want to show that for all $n\geq 1$  and  all non-empty indices ${\bf k}$ we have for any primitive $n$-th root of unity $\zeta_n$
\[ H_{n-1}^\star({\bf k};\zeta_n) = (-1)^{\wt({\bf k})+1} H_{n-1}^\star(\overline{{\bf k}^\vee};\zeta_n) \,.\]

Here we define for an index ${\bf k}=(k_1,\dots,k_r)$ its reverse by $\overline{\kk} = (k_r,\dots,k_1)$ and ${\bf k}^\vee$ is the Hoffman dual, which we defined at the end of Section \ref{subsec:fmzv}.

We will use the following fact. 

\begin{lemma}\label{lem:q-binom-power}
Suppose that $n \ge 1$ and $\zeta_n$ is a primitive $n$-th root of unity. 
Then it holds that $(-1)^{n}\zeta_n^{n(n+1)/2}=-1$.  
\end{lemma}

\begin{proof}[Proof of Theorem \ref{thm:btthoffmandual}]
Notice that any index is uniquely written in the form 
\begin{align}
(\{1\}^{a_{1}-1}, b_{1}+1, \ldots , \{1\}^{a_{r-1}-1}, b_{r-1}+1, \{1\}^{a_{r}-1}, b_{r}),  
\label{eq:index-ab-rep}
\end{align}  
where $r$ and $a_{i}, b_{i} \, (1\le i \le r)$ are positive integers
\footnote{If $r=1$, \eqref{eq:index-ab-rep} should read as $(\{1\}^{a_{1}-1}, b_{1})$.}. 
Denote it by $[a_{1}, \ldots , a_{r}; b_{1}, \ldots , b_{r}]$. 
Then we see that 
\begin{align*}
\overline{[a_{1}, \ldots , a_{r}; b_{1}, \ldots , b_{r}]^{\vee}}=
[b_{r}, \ldots , b_{1}; a_{r}, \ldots , a_{1}]. 
\end{align*}

Now we fix a positive integer $r$ and introduce the generating function 
\begin{align*}
K(x_{1}, \ldots , x_{r}; y_{1}, \ldots , y_{r})=\sum
\frac{H_{n-1}^{\star}([a_{1}, \ldots , a_{r}; b_{1}, \ldots , b_{r}]; \zeta_n)}{(1-\zeta_n)^{a_1+\cdots+a_r+b_1+\cdots+b_r-1}}
\prod_{i=1}^{r}(x_{i}^{a_{i}-1}y_{i}^{b_{i}-1}),   
\end{align*}
where the sum is taken over all positive integers $a_{i}, b_{i} \, (1\le i \le r)$. 
Then Theorem \ref{thm:btthoffmandual} follows from the equality 
\begin{align}
K(x_{1}, \ldots , x_{r}; y_{1}, \ldots , y_{r})=K(-y_{r}, \ldots , -y_{1}; -x_{r}, \ldots , -x_{1}). 
\label{eq:zstar-duality-proof1}
\end{align}

Let us prove \eqref{eq:zstar-duality-proof1}. 
It holds that 
\begin{align*}
1+\sum_{a=2}^{\infty}\sum_{B\ge m_{1}\ge \cdots \ge m_{a-1} \ge A}
\frac{x^{a-1}}{\prod_{i=1}^{a-1}(1-\zeta_{n}^{m_{i}})}=
\prod_{i=A}^{B}\frac{1-\zeta_{n}^{i}}{1-x-\zeta_{n}^{i}} 
\end{align*}
for $n>B\ge A>0$, and that  
\begin{align*}
\sum_{b=1}^{\infty}\frac{\zeta_{n}^{bm}}{(1-\zeta_{n}^{m})^{b+1}}y^{b-1}=
\frac{1}{1-\zeta_{n}^{m}}\,\frac{\zeta_{n}^{m}}{1-\zeta_{n}^{m}(1+y)} 
\end{align*}
for $n>m>0$. 
Using the above formulas we have 
\begin{align*}
& 
K(x_{1}, \ldots , x_{r}; y_{1}, \ldots , y_{r})=\sum_{n>l_{1}\ge \cdots \ge l_{r}>0}
\prod_{i=l_{r}}^{n-1}(1-\zeta_n^{i}) \\ 
&\times 
\prod_{j=1}^{r-1}\left(\frac{\zeta_n^{l_{j}}}{1-\zeta_n^{l_{j}}(1+y_{j})}
\prod_{i=l_{j}}^{l_{j-1}}\frac{1}{1-x_{j}-\zeta_n^{i}}\right)
\frac{1}{1-\zeta_n^{l_{r}}(1+y_{r})}\prod_{i=l_{r}}^{l_{r-1}}\frac{1}{1-x_{r}-\zeta_n^{i}},  
\end{align*}
where $l_{0}=n-1$. 
Rewrite the right-hand side above by using the partial fraction expansion 
\begin{align*}
\prod_{i=A}^{B}\frac{1}{X-\zeta_{n}^{i}}&= 
\sum_{i=A}^{B} \frac{1}{X-\zeta_{n}^i} \prod_{j=A}^{i-1} \frac{1}{\zeta_{n}^i-\zeta_{n}^j} 
\prod_{j=i+1}^{B} \frac{1}{\zeta_{n}^i-\zeta_{n}^j}
\notag \\
&=\sum_{t=A}^{B}\frac{1}{X-\zeta_{n}^{t}}
\frac{(-1)^{B-t}\zeta_{n}^{-\binom{B+1}{2}+At-\binom{t}{2}}}
{\prod_{i=1}^{t-A}(1-\zeta_{n}^{-i})\prod_{i=1}^{B-t}(1-\zeta_{n}^{-i})} 
\end{align*}
for $n>B\ge A>0$. 
Then we find that 
\begin{align*}
&
K(x_{1}, \ldots , x_{r}; y_{1}, \ldots , y_{r}) \\ 
&=\sum_{n>t_{1}\ge l_{1}\ge \cdots \ge t_{r}\ge l_{r}>0}
\prod_{i=l_{r}}^{n-1}(1-\zeta_n^{i}) 
(-1)^{\sum_{j=1}^{r}(l_{j-1}-t_{j})}
\zeta_{n}^{\sum_{j=1}^{r}(-\binom{l_{j-1}+1}{2}+l_{j}t_{j}-\binom{t_{j}}{2})} 
\nonumber \\ 
&\qquad \qquad {}\times 
\prod_{j=1}^{r}\left(\prod_{i=1}^{t_{j}-l_{j}}\frac{1}{1-\zeta_n^{-i}}
\prod_{i=1}^{l_{j-1}-t_{j}}\frac{1}{1-\zeta_n^{-i}}\right) 
\nonumber \\ 
&\qquad \qquad {}\times 
\prod_{j=1}^{r-1}\left(\frac{\zeta_n^{l_{j}}}{1-\zeta_n^{l_{j}}(1+y_{j})}\frac{1}{1-x_{j}-\zeta_n^{t_{j}}}\right)
\frac{1}{1-\zeta_n^{l_{r}}(1+y_{r})}\frac{1}{1-x_{r}-\zeta_n^{t_{r}}}.   
\end{align*}
Now change the summation variables $t_{j}$ and $l_{j}$ to 
$n-l_{r+1-j}$ and $n-t_{r+1-j}$, respectively ($1\le j \le r$). 
As a result we get the desired equality \eqref{eq:zstar-duality-proof1} using Lemma \ref{lem:q-binom-power}.  
\end{proof}

\subsection{Reversal relation \& Linear shuffle relations}
We will now describe the linear shuffle relations for the multiple harmonic sums at roots of unity. Recall that we set $\overline{\N} = \{ \bar{1} \} \cup \Z_{\geq 1}$ and then define for $k\in \overline{\N} $ and $m\geq 1$ the following $q$-series
\begin{align}
    f_k(m) = \begin{cases}  \frac{q^m}{[m]_q},&\text{ if } k = \bar{1}\\
    \frac{q^{(k-1)m}}{[m]_q^k},&\text{ if } k \in \Z_{\geq 1}
\end{cases}\,.
\end{align}   
With this we can generalize the definition of the multiple harmonic $q$-series (and their star version) for $\kk = (k_1,\dots,k_r) \in \ON^r$ and $m\geq 1$ by
\begin{align}\label{eq:genmultipleharmq}
    H_m(\kk;q) &=  H_m(k_1,\dots,k_r;q) = \sum_{m\geq m_1 > \dots > m_r > 0} \prod_{j=1}^r f_{k_j}(m_j) \in \Q[[q]]\,,\\
    H^\star_m(\kk;q) &=  H^\star_m(k_1,\dots,k_r;q) = \sum_{m\geq m_1 \geq \dots \geq m_r > 0} \prod_{j=1}^r f_{k_j}(m_j) \in \Q[[q]]\,.
\end{align}
These can be viewed for any $m\geq 1$ as $\CC$-algebra homomorphisms defined on the generators by
\begin{align*}
    H_m(\cdot; q): \HHH^1_{\ast_q} &\longrightarrow \Q[[q]],\\
    e_{k_1}\dots e_{k_r} &\longmapsto H_m(k_1,\dots,k_r;q)\,,
\end{align*}
where we defined $\HHH^1 = \CC\langle e_{\bar{1}}, e_1, e_2, e_3, \dots \rangle$ and $\HHH^1_{\ast_q}$ is $\HHH^1 $ equipped with $q$-stuffle product. As in Section \ref{sec:qmzv} we view $\Q[[q]]$ here as a $\CC$-module, where the multiplication with $\hbar$ is given by the multiplication with $(1-q)$. Now recall that for $k\geq 2$ we defined 
\begin{align}
    e_{\overline{k}} &:= \sum_{j=2}^k \binom{k-2}{j-2} \hbar^{k-j} e_j\,,
\end{align}
which was inspired by the equation
\begin{align}
    \frac{q^{m}}{[m]_q^k} &= \sum_{j=2}^k \binom{k-2}{j-2} (1-q)^{k-j} \frac{q^{(j-1)m}}{[m]_q^j}\,,\qquad (k\geq 2).
\end{align}
The next lemma shows that the $e_{\overline{k}}$ will be used when considering the reverse of a word (doing the change of summation variables $m\rightarrow n-m$), when evaluating the $H_{n-1}$ at primitive $n$-th roots of unity. 

\begin{lemma}\label{lem:qreverse}Let $q=\zeta_n$ be a primitive $n$-th root of unity and $n>m$. Then we have for $k\geq 1$
\begin{align*}
\qquad \frac{q^{(k-1)(n-m)}}{[n-m]^k_q} = (-1)^k \frac{q^{m}}{[m]_q^k}\,.
\end{align*}
\end{lemma}
\begin{proof}
By direct calculation, we get $    [n-m]_q = - q^{-m} [m]_q$, from which the statement follows.
\end{proof}

Lemma \ref{lem:qreverse} can be seen as an analogue of the equation $\frac{1}{m^k} \equiv (-1)^k \frac{1}{(p-m)^k} \mod p$, which we used in the proof of the reversal relation and the linear shuffle relations of finite multiple zeta values. Inspired by this we define the following $q$-version of the reversal of a word.

\begin{definition}
We define the anti-automorphism\footnote{with respect to the concatenation product} $\psi: \HHH^1 \rightarrow \HHH^1$ for $k\geq 1$ by 
\begin{align*}
    \psi(e_k) = e_{\overline{k}},\quad \psi( e_{\overline{1}}) = e_{1}\,,
\end{align*}
i.e. for $k_1,\dots,k_r\in \ON$ we have $\psi(e_{k_1}\cdots e_{k_r}) = e_{\overline{k_r}}\cdots e_{\overline{k_1}}$ with the convention $e_{\overline{\overline{1}}}=e_1$.
\end{definition}
One can check that $\psi$ is an involution and as a consequence of Lemma \ref{lem:qreverse} we get the following analogue of the reversal formula:

\begin{proposition} For any primitive $n$-th root of unity $\zeta_n$ and $w \in \HHH^1$ we have 
\begin{align*}
    H_{n-1}(w;\zeta_n) = (-1)^{\wt(w)} H_{n-1}(\psi(w);\zeta_n)\,.
\end{align*}
\end{proposition}
\begin{proof}
This follows now by making the change of summation variables $m_j \rightarrow n - m_j$ (as in \eqref{eq:reversalfmzvcalc}) in \eqref{eq:genmultipleharmq}, together with Lemma \ref{lem:qreverse}.
\end{proof}

As for finite multiple zeta values, this reversal relation is just a special case of the linear shuffle relations.

\begin{theorem} \label{thm:linearshuffleq}For  $n\geq 2$ and any primitive $n$-th root of unity $\zeta_n$ and $w,v \in \HHH^1$ we have
 \begin{align*}
        H_{n-1}(w \shq v;\zeta_n) = (-1)^{\wt(w)}   H_{n-1}(\psi(w)v;\zeta_n)\,.
 \end{align*}
\end{theorem}
\begin{proof}
This can be shown in a similar way as the linear shuffle relations for finite multiple zeta values (Theorem \ref{thm:linearshufflefmzv}) together with Proposition \ref{prop:qpolylogalghom} (Exercise \ref{fex16}).
\end{proof}

\subsection{$\mathcal{Q}$-multiple zeta values}
For a prime $p$ we found a lot of relations among $H_{p-1}(\kk;\zeta_p)$ which are valid for any primitive $p$-th root of unity $\zeta_p$. Inspired by this, the authors of \cite{TT} introduced $\mq$-multiple zeta values, which capture these relations. Since 
\begin{align*}
    [p]_q = 1+q+\dots+q^{p-1}=\frac{1-q^p}{1-q} = \prod_{a=1}^{p-1}\left(q- e^{\frac{2\pi ia}{p}} \right)\,,
\end{align*}
a polynomial $f(q) \in \Z[q]$ satisfies $f(\zeta_p) =0 $ for all primitive $p$-th roots of unity if and only if $f(q) \equiv 0 \mod [p]_q$. We will consider the following ring, which can be seen as an analogue of the ring $\ma$.

\begin{definition} Define the ring $\mq$ by
\begin{align*}
    \mq = \quotient{\prod_{p \text{ prime}} \fpq \,}{\bigoplus_{p \text{ prime}} \fpq }\,,
\end{align*}
where $\Z_{(p)}=\{ \frac{a}{b} \mid a,b \in \Z, b\not \in p\Z\}$ denotes the localization of $\Z$ at the ideal $(p)$ and $([p]_q)$ is the ideal of $\Z_{(p)}[q]$ generated by the irreducible polynomial $[p]_q$.
\end{definition}
Notice that $\mq$ is a $\Q[q]$ algebra since we have the embedding
\begin{align*}
    \Q[q] &\longrightarrow \mq\\
    f &\longmapsto \left( f \mod [p]_q \right)_p\,,
\end{align*}
which is well-defined, since we can ignore those finitely many primes $p$, for which the denominators of the coefficients of $f$ are divisible by $p$.
The $\mq$-multiple zeta values, introduced by Takeyama and Tasaka in
\cite[Section~2.7, Definition~2.5]{TT} in the case $\kk \in \Z^r_{\geq 1}$, are
then defined as the following objects in $\mq$.

\begin{definition} For any index $\kk \in \ON^r$ we define the $\mq$-multiple zeta value $\zq(\kk) \in \mq$ by
\begin{align*}
    \zq(\kk) &= \left(  H_{p-1}(\kk; q) \mod [p]_q \right)_{p}\,,
\end{align*}
where the $H_{p-1}$ are given by \eqref{eq:genmultipleharmq}.
\end{definition}

\begin{proposition}[Algebraic realization]\label{prop:qtofinite}
There is a surjective $\Q$-algebra homomorphism
\begin{align*}
 \phi_{\ma}:\mq&\longrightarrow\ma,\\
 \left(f_p(q)\bmod[p]_q\right)_p&\longmapsto
 \left(f_p(1)\bmod p\right)_p.
\end{align*}
If $\kk^\circ$ denotes the index obtained from $\kk\in\ON^r$ by replacing every
$\overline1$ by $1$, then
\begin{align*}
 \phi_{\ma}\big(\zq(\kk)\big)=\za(\kk^\circ).
\end{align*}
In particular, $\phi_{\ma}(\zq(\kk))=\za(\kk)$ for every ordinary index
$\kk\in\Z_{\geq1}^r$.
\end{proposition}
\begin{proof}
Since $[p]_1=p$, evaluation at $q=1$ induces a map
\begin{align*}
 \Z_{(p)}[q]/([p]_q)\longrightarrow\Z/p\Z.
\end{align*}
Taking the product over all primes and passing to the quotients gives $\phi_{\ma}$.
Changing finitely many components does not change either side. The map is surjective
because it is the identity on constant representatives. Finally,
\begin{align*}
 \phi_{\ma}\big(\zq(\kk)\big)
 =\left(H_{p-1}(\kk;1)\bmod p\right)_p=\za(\kk^\circ),
\end{align*}
since both $f_1(m)$ and $f_{\overline1}(m)$ specialize to $1/m$ at $q=1$.
For ordinary indices this is \cite[Section~2.8 and Theorem~2.8]{TT}.
\end{proof}

The analytic realization, introduced in \cite[Section~2.8]{TT}, is not defined on
all of $\mq$. Let
\begin{align*}
 \mathcal O=\left\{\left(f_p(q)\bmod[p]_q\right)_p\in\mq\ \middle|\
 \lim_{p\rightarrow\infty}f_p\left(e^{2\pi i/p}\right)\text{ exists}\right\},
\end{align*}
where $p$ tends to infinity through the primes. Evaluation at the indicated root of
unity and passage to the limit define the $\Q$-algebra homomorphism
\begin{align*}
 \phi_{\ms}:\mathcal O&\longrightarrow\C,\\
 \left(f_p(q)\bmod[p]_q\right)_p&\longmapsto
 \lim_{p\rightarrow\infty}f_p\left(e^{2\pi i/p}\right).
\end{align*}
This is well-defined because $[p]_{e^{2\pi i/p}}=0$ and a change in finitely many
components does not affect the limit.

\begin{proposition}[Analytic realization]\label{prop:qtosymmetric}
For every ordinary index $\kk\in\Z_{\geq1}^r$, we have $\zq(\kk)\in\mathcal O$ and
\begin{align*}
 \phi_{\ms}\big(\zq(\kk)\big)
 &=\lim_{p\rightarrow\infty}H_{p-1}\left(\kk;e^{2\pi i/p}\right)=\xi(\kk),\\
 \operatorname{Re}\left(\phi_{\ms}\big(\zq(\kk)\big)\right)
 &\equiv\zs(\kk)\pmod{\pi^2\mz}.
\end{align*}
\end{proposition}
\begin{proof}
The existence of the limit and the first equality follow from Theorem
\ref{thm:xilimit}. See also \cite[Theorem~2.4~(ii) and Theorem~2.8]{TT}. The
congruence is the last statement of Theorem \ref{thm:xilimit}.
\end{proof}

The restriction to ordinary indices in Proposition \ref{prop:qtosymmetric} is necessary.
For example, \cite[Remark~2.11]{BTT1} gives
\begin{align*}
 H_{n-1}\left(\overline1,1;e^{2\pi i/n}\right)
 =2\zeta(2)+2\pi i\left(\log\frac{n}{2\pi}+\gamma\right)
 +\operatorname{O}\left(\frac{\log n}{n}\right),
\end{align*}
so $\zq(\overline1,1)$ does not belong to $\mathcal O$. Thus the algebraic
specialization works for the extended alphabet $\ON$, whereas the analytic specialization
does not. In particular, a $\Q$-linear relation among ordinary $\mq$-multiple zeta values specializes
under $\phi_{\ma}$ to the corresponding finite relation and, provided its terms belong to
$\mathcal O$, under $\phi_{\ms}$ to the corresponding symmetric relation after taking real
parts modulo $\pi^2\mz$. This explains the BTT philosophy above.

\vspace{1cm}
\begin{center}
\bf{ \Large {\color{qmzvlinecol}\ding{94}}~ Exercises ~{\color{qmzvlinecol}\ding{94}}}\end{center}

The following is a collection of exercises intended to deepen the reader's understanding of this chapter.

\begin{exer}\label{ffex14}
For $k\geq1 $ show that $\za(k)=0$ just using the equations \eqref{eq:conjast} and \eqref{eq:conjshuffle} in Conjecture \ref{conj:fmzvallrelations}.
\end{exer}

\begin{exer}\label{fex15}
Prove Lemma \ref{lem:tindepend}, i.e. show that for any index $\kk=(k_1,\dots,k_r)$ the polynomial 
\begin{align*}
R_\kk(X;T) = \sum_{j=0}^r (-1)^{k_1+\dots+k_j} \zeta^\ast(k_j,k_{j-1},\dots,k_1;T+X) \zeta^\ast(k_{j+1},\dots,k_r;T-X)
\end{align*}
does not depend on $T$.
\end{exer}

\begin{exer}\label{fex16} Prove the linear shuffle relations for multiple harmonic $q$-series at roots of unity (Theorem \ref{thm:linearshuffleq}), i.e. show that for any primitive $n$-th root of unity $\zeta_n$ and $w,v \in \HHH^1$ we have
 \begin{align*}
        H_{n-1}(w \shq v;\zeta_n) = (-1)^{\wt(w)}   H_{n-1}(\psi(w)v;\zeta_n)\,.
 \end{align*}
\end{exer}

%% file: chap_FamiliesOfRelations.tex
\chapter[Further relations \& variants of MZVs]{Further relations and variants \\of multiple zeta values} \label{sec:families}
	
In this chapter, we want to discuss several families of relations among multiple zeta values and also some of their $q$-analogues. Asking for linear relations among multiple zeta values is equivalent to asking for the kernel of the map $\zeta: \HH^0 \rightarrow \mz$. We will start with Seki-Yamamoto's connected sums, which give a nice proof of the duality relation and of Ohno's relation. After this, we discuss the Maesaka--Seki--Watanabe formula and the Drop1 operator, which will show up again when we study the derivatives of multiple Eisenstein series in Section \ref{sec:mesderiv}. We then have a look at parity results and give a small overview of the zoo of relations. At the end, we come back to the Schur multiple zeta values of Section \ref{subsec:schurmzvoverview} and prove a Jacobi--Trudi formula for them.

\section{Seki-Yamamoto's connected sums and Ohno's relation}

In this section, we want to present an alternative proof of the duality relation (Proposition \ref{prop:duality}), which was given by Seki and Yamamoto in \cite{SY}\footnote{In the work \cite{SY} the order of summation in the definition of multiple zeta values is reversed. Therefore one needs to be careful when comparing the results here and the ones in their work.}. In this nice work, they introduce the notion of connected sums, which recently also found their way into various other proofs of families of relations among multiple zeta values and some of their variants, such as finite multiple zeta values. In \cite{Se1} you can find an overview of different applications of connected sums. We will use this setup to present a proof of Ohno's relation for $q$-analogues of multiple zeta values as it was given in \cite{SY}.

\subsection{Connected sums and the duality relation for multiple zeta values}

We start by reformulating the duality relations on the level of indices. Recall that $\tau : \HH^0 \rightarrow \HH^0$ was defined as the anti-automorphism (with respect to the usual multiplication in $\Q\langle x,y \rangle$) satisfying $\tau(x)=y$ and $\tau(y)=x$. Now let $\kk = (k_1,\dots,k_r)$ be an admissible index, i.e. $k_1\geq 2$ and $z_\kk \in \HH^0$. Then there exist numbers $a_1,b_1,\dots,a_s,b_s \geq 1$ such that
\begin{align*}
\kk = (a_1 + 1, \{1\}^{b_1-1},a_2 +1 , \{1\}^{b_2-1},\dots, a_s+1, \{1\}^{b_s-1})\,.
\end{align*}
With these numbers we define the admissible index $\kk^\dagger$ by 
\begin{align*}
\kk^\dagger := (b_s + 1, \{1\}^{a_s-1},b_{s-1} +1 , \{1\}^{a_{s-1}-1},\dots, b_1+1, \{1\}^{a_1-1})\,.
\end{align*}
One can then see by induction on $s$ that we have (Exercise \ref{ex11})
\begin{align}
z_{\kk^\dagger} = \tau(z_\kk)\,,
\end{align} 
i.e. the duality relation of multiple zeta values can be stated as $\zeta(\kk) = \zeta(\kk^\dagger)$ for all admissible $\kk$. 

Now we will introduce the connected sum for multiple zeta values. The name comes from the fact that these sums look like the product of two multiple zeta values, which are connected at one point by a connector.

\begin{definition}  Let $\kk = (k_1,\dots,k_r)$, $\kl = (l_1,\dots,l_s)$ be two non-empty indices. Then we define the \emph{connected sum} $Z(\kk;\kl)$ by 
	\begin{align} \label{eq:defconsum}
	Z(\kk ;\kl) = Z(k_1,\dots ,k_r ; l_1, \dots, l_s) = \sum_{\substack{m_1 > m_2 > \cdots > m_r > 0 \\ n_1 > n_2 > \cdots > n_s >0}} \frac{m_1! \,n_1!}{(m_1+n_1)!}\, \frac{1}{m_1^{k_1} \cdots m_r^{k_r}} \,\frac{1}{n_1^{l_1} \cdots n_s^{l_s}}
	\end{align}
	and set $Z(\kk; \emptyset) = Z(\emptyset ;\kk ) = \zeta(\kk)$ for admissible $\kk$.
\end{definition}
One can check that the sum \eqref{eq:defconsum} converges for all non-empty indices $\kk$ and $\kl$. Also notice that  \eqref{eq:defconsum}  is essentially the product $\zeta(\kk) \zeta(\kl)$, which gets connected by the \emph{connector} $c(m_1,n_1)=\frac{m_1! \,n_1!}{(m_1+n_1)!}$ at the beginning. The relationship to the duality relations comes from the fact that one can show that $Z(\kk;\emptyset) =\dots = Z(\emptyset; \kk^\dagger)$, by using the following transport relations.

\begin{proposition}[Transport relations] \label{thm:transportrel} Let $(k_1,\dots,k_r)$ and $(l_1,\dots,l_s)$ be two indices.
	If $s>0$ then we have
	\begin{align*}
	Z(1,k_1,\dots,k_r ; l_1, \dots, l_s) = Z(k_1,\dots,k_r ; l_1+1, l_2, \dots, l_s)
	\end{align*}
	and if $r>0$, then we have 
	\begin{align*}
	Z(k_1+1,k_2,\dots,k_r ; l_1, \dots, l_s) = Z(k_1,\dots,k_r ; 1, l_1, l_2, \dots, l_s)\,.
	\end{align*}
\end{proposition}
\begin{proof} To prove the first equality we use for $m\geq 0$ the telescoping sum         
	\begin{align}\label{eq:transporttelescope}
	\sum_{a=m+1}^\infty \frac{1}{a} \frac{a! \,n!}{(a+n)!} &= \frac{1}{n} \sum_{a=m+1}^\infty \left(  \frac{(a-1)! n!}{(a-1+n)!} - \frac{a! n!}{(a+n)!}   \right) =\frac{1}{n} \frac{m! n!}{(m+n)!}\,,
	\end{align}
	from which we obtain 
	\begin{align*}
	Z(1,k_1,\dots,k_r ; l_1, \dots, l_s) &= \sum_{\substack{a>m_1 > m_2 > \cdots > m_r > 0 \\ n_1 > n_2 > \cdots > n_s >0}} \frac{1}{a} \frac{a! \,n_1!}{(a+n_1)!}\, \frac{1}{m_1^{k_1} \cdots m_r^{k_r}} \,\frac{1}{n_1^{l_1} \cdots n_s^{l_s}} \\
	&= \sum_{\substack{m_1 > m_2 > \cdots > m_r > 0 \\ n_1 > n_2 > \cdots > n_s >0}} \frac{1}{n_1}\frac{m_1! \,n_1!}{(m_1+n_1)!}\, \frac{1}{m_1^{k_1} \cdots m_r^{k_r}} \,\frac{1}{n_1^{l_1} \cdots n_s^{l_s}} \\
	&=  Z(k_1,\dots,k_r ; l_1+1, l_2, \dots, l_s)\,.
	\end{align*}
	Also notice that we obtain $    Z(1 ; l_1, \dots, l_s) = Z(\emptyset;l_1+1,l_2,\dots,l_s)=\zeta(l_1+1,l_2,\dots,l_s)$ here in the case $r=0$  by using the $m=0$ case of \eqref{eq:transporttelescope}.
	The second statement follows from the symmetry $Z(\kk;\kl) = Z(\kl ; \kk)$.
\end{proof}

Using these transport relations, we can therefore always transport one index from the left to the right and vice versa. For example, we have 
\begin{align*}
\zeta(3) = Z(3;\emptyset) = Z(2;1) = Z(1;1,1) = Z(\emptyset; 2,1 ) = \zeta(2,1)\,.
\end{align*}
In general the duality relation follows from this, since if we start with an admissible index $    \kk = (a_1 + 1, \{1\}^{b_1-1},a_2 +1 , \{1\}^{b_2-1},\dots, a_s+1, \{1\}^{b_s-1})$, we can apply the second transport relation $a_1$ times and the first one $b_1$ times, then the second one $a_2$ times and the first one $b_2$ times, etc., to get
\begin{align}\label{eq:dualityfromtransport}
\begin{split}
\zeta(\kk) &= Z(a_1 + 1, \{1\}^{b_1-1},a_2 +1 , \{1\}^{b_2-1},\dots, a_s+1, \{1\}^{b_s-1} ; \emptyset) = \dots \\
&= Z(\{1\}^{b_1-1},a_2 +1 , \{1\}^{b_2-1},\dots, a_s+1, \{1\}^{b_s-1}; 2, \{1\}^{a_1-1})= \dots \\
&= Z(a_2 +1 , \{1\}^{b_2-1},\dots, a_s+1, \{1\}^{b_s-1}; b_1+1, \{1\}^{a_1-1})= \dots \\
&= \dots \\
&= Z(\emptyset ; b_s + 1, \{1\}^{a_s-1},b_{s-1} +1 , \{1\}^{a_{s-1}-1},\dots, b_1+1, \{1\}^{a_1-1}) = Z(\emptyset; \kk^\dagger) = \zeta(\kk^\dagger)\,.
\end{split}
\end{align}

\subsection{The sum formula and Ohno's relation}
We will now present a family of linear relations, which are known as Ohno's relation. This type of relation generalizes the duality relation and the so-called sum formula, given by the following.

\begin{theorem}[Sum formula (Granville {\cite[equation~(1)]{G}}, Zagier independently)] \label{thm:sumformula} For all $k\geq 2$ and $1  \leq r < k$ we have
	\begin{align}\label{eq:sumformula}
	\sum_{\substack{k_1+\dots + k_r = k\\ k_1 \geq 2,\, k_2,\dots, k_r \geq 1}} \zeta(k_1,\dots,k_r) = \zeta(k)\,.    \end{align}
\end{theorem}
The sum formula therefore states that the sum over all multiple zeta values of weight $k$ in any fixed depth always gives the Riemann zeta value $\zeta(k)$. Notice again that our first relation $\zeta(2,1)=\zeta(3)$ is also the first (non-trivial) example of the sum formula.
There are various generalizations of the sum formula which also include weights and therefore are called weighted sum formulas.

To state Ohno's relation we first define for an admissible index $\kk=(k_1,\dots,k_r)$ the \emph{Ohno sum}    by
\begin{align*}
 \mathcal{O}^X(\kk) = \sum_{c\geq 0} O(\kk ; c)  X^c \in \mz[[X]]\,,
\end{align*}
where we write for $c\geq 0$
\begin{align*}
O(\kk ; c) = \sum_{\substack{c_1 + \dots + c_r = c\\c_1,\dots,c_r \geq 0}} \zeta(k_1+c_1,\dots,k_r+ c_r) \in \mz_{\wt(\kk)+c}\,.
\end{align*}
Notice that for $\kk=(2,\{1\}^{r-1})$ and $c=k-r-1$ we obtain the left-hand side of \eqref{eq:sumformula}
\begin{align*}
O(2,\{1\}^{r-1}  ; k-r-1) =     \sum_{\substack{k_1+\dots + k_r = k\\ k_1 \geq 2,\, k_2,\dots, k_r \geq 1}} \zeta(k_1,\dots,k_r)\,.
\end{align*}
Ohno's relation now states that the Ohno sum also satisfies the duality relation.

\begin{theorem}[Ohno's relation  (Ohno \cite{Oh})] \label{thm:ohno} For any admissible index $\kk$ we have
	\begin{align*}
	\mathcal{O}^X(\kk) =     \mathcal{O}^X(\kk^\dagger)\,.
	\end{align*}
\end{theorem}
Since     $\mathcal{O}^0(\kk) = O(\kk; 0) = \zeta(\kk)$ we obtain the duality as a special case by considering the constant term in Ohno's relation. Choosing for $r\geq 1$ the index $\kk=(2,\{1\}^{r-1})$ we have $\kk^\dagger = (r+1)$, and therefore the sum formula \eqref{eq:sumformula} follows by considering the coefficient of $X^{k-r-1}$ in $\mathcal{O}^X(2,1,\dots,1) =     \mathcal{O}^X(r+1)$. The formulation of the Ohno relation in the original work of Ohno is a bit different and we use here the formulation of \cite{HMOS}, where the authors also prove further relations of the Ohno sums besides the duality relation.

\begin{ex} For $\kk=(2,2,1)$ we have $\kk^\dagger = (3,2)$ and
	\begin{align*}
	\mathcal{O}^X(2,2,1)  = &\,\,\zeta(2,2,1) + \big( \zeta(3,2,1) + \zeta(2,3,1) + \zeta(2,2,2)  \big) X \\&+ \big( \zeta(4,2,1) + \zeta(2,4,1) + \zeta(2,2,3)+ \zeta(3,3,1)+\zeta(3,2,2)+\zeta(2,3,2)     \big) X^2 + \dots \,,\\
	\mathcal{O}^X(3,2)  = &\,\, \zeta(3,2) + \big( \zeta(4,2) + \zeta(3,3) \big) X+ \big( \zeta(5,2) + \zeta(3,4) + \zeta(4,3) \big) X^2 + \dots\,,
	\end{align*}
		which gives a linear relation among multiple zeta values of weight $5+c$ by comparing the coefficients of $X^c$ in $\mathcal{O}^X(2,2,1) = \mathcal{O}^X(3,2) $.
\end{ex}

The number of Ohno's relations is given by the following table (calculated by Tanaka  in \cite{Tan1}).
	\begin{center}
	\begin{tabular}{|c|c|c|c|c|c|c|c|c|c|c|c|c|c|}
		\hline
		weight $k$         & 3 & 4 & 5 & 6 & 7 & 8 & 9 & 10 & 11 & 12  \\ \hline
		\# all conjectured relations       &1& 3& 6& 14& 29& 60& 123& 249& 503& 1012\\ \hline
		\# Ohno's relations  & 1 & 2 & 5 & 10 & 23 & 46 & 98 & 199  & 411  & 830  \\ \hline
	\end{tabular}
\end{center}

We now want to state the $q$-analogue version of Ohno's relation. Recall that we defined for an admissible index $\kk = (k_1,\dots,k_r)$ the Bradley-Zhao $q$-analogues of multiple zeta values by 
\begin{align}
\zeta_q^{\rm{BZ}}(\kk)= \zeta_q^{\rm{BZ}}(k_1,\dots,k_r) = \sum_{m_1 > \cdots > m_r > 0} \frac{q^{(k_1-1)m_1} \cdots q^{(k_r-1)m_r}}{[m_1]_q^{k_1} \cdots [m_r]_q^{k_r} } \,,
\end{align} 
where $[m]_q = \frac{1-q^m}{1-q}$ denotes the $q$-integer. It was first shown by Bradley (\cite{Bra}) that these $q$-series also satisfy Ohno's relation. Define for an admissible index $\kk=(k_1,\dots,k_r)$ the \emph{$q$-Ohno sum}    by
\begin{align*}
\mathcal{O}^X_q(\kk) = \sum_{c\geq 0} O_q(\kk ; c)  X^c \in \Q[[q]][[X]]\,,
\end{align*}
where we write for $c\geq 0$
\begin{align*}
O_q(\kk ; c) = \sum_{\substack{c_1 + \dots + c_r = c\\c_1,\dots,c_r \geq 0}} \zeta_q^{\rm{BZ}}(k_1+c_1,\dots,k_r+ c_r)\,.
\end{align*}

The series $\mathcal{O}^X_q(\kk)$ is a formal power series in $X$ with coefficients given by formal power series in $q$. For real $0 < q,X < 1$ one can show that $\mathcal{O}^X_q(\kk)$ converges and gives a well-defined real number.
\begin{theorem}[Ohno's relation for $q$-MZV  (Bradley \cite{Bra})]\label{thm:qohno}  For any admissible index $\kk$ we have
	\begin{align*}
	\mathcal{O}^X_q(\kk) =     \mathcal{O}^X_q(\kk^\dagger)\,.
	\end{align*}
\end{theorem}
Notice that Theorem \ref{thm:qohno} implies that the $q$-analogues $\zeta_q^{\rm{BZ}}$ also satisfy the sum formula and the duality relation. Also, sending $q \rightarrow 1$ we obtain Ohno's relation for multiple zeta values.

The goal is now to give a proof of Theorem \ref{thm:qohno} by using Seki-Yamamoto's concept of connected sums. For this we first rewrite the $q$-Ohno sum as
\begin{align*}
\mathcal{O}^X_q(\kk) &= \sum_{c\geq 0} O_q(\kk ; c)  X^c  = \sum_{\substack{c_1,\dots,c_r \geq 0}} \zeta_q^{\rm{BZ}}(k_1+c_1,\dots,k_r+ c_r) X^{c_1+\dots+c_r}\\
&= \sum_{\substack{m_1 > \cdots > m_r > 0\\c_1,\dots,c_r \geq 0}} \frac{q^{(k_1+c_1-1)m_1} \cdots q^{(k_r+c_r-1)m_r}}{[m_1]_q^{k_1+c_1} \cdots [m_r]_q^{k_r+c_r} }   X^{c_1+\dots+c_r}\\
&= \sum_{m_1 > \cdots > m_r > 0} \frac{q^{(k_1-1)m_1}}{([m_1]_q-q^{m_1}X)[m_1]_q^{k_1-1}}\cdots  \frac{q^{(k_r-1)m_r}}{([m_r]_q-q^{m_r}X)[m_r]_q^{k_r-1}} \\
&= \sum_{m_1 > \cdots > m_r > 0} s^X_q(k_1,m_1) \cdots s^X_q(k_r,m_r)\,,
\end{align*}
where we set 
\begin{align}\label{eq:sxq}
s^X_q(k,m) = \frac{q^{(k-1)m}}{([m]_q-q^{m}X)[m]_q^{k-1}}\,. 
\end{align}
Notice that $s^0_1(k,m)=\frac{1}{m^k}$, i.e. the $q$-Ohno sum reduces to the multiple zeta value $\mathcal{O}^0_1(\kk)=\zeta(\kk)$ in this case. To define the connected sum we need to find the correct connector $c^X_q(m,n)$, which generalizes the connector $c^0_1(m,n) = c(m,n)=\frac{m! \,n!}{(m+n)!}$ we used for multiple zeta values. This was done by Seki and Yamamoto in \cite{SY}, by choosing the connector
\begin{align}\label{eq:cxq}
c^X_q(m,n) =\frac{q^{m n} f^X_q(m) f^X_q(n)}{f^X_q(m+n)}\,,
\end{align}
where $f^X_q(m) = \prod_{j=1}^m( [j]_q - q^j X)$, which can be seen as a variant of the factorial, since $f^0_1(m)=m!$.
\begin{definition} Let $\kk = (k_1,\dots,k_r)$, $\kl = (l_1,\dots,l_s)$ be two non-empty indices. Then we define the \emph{connected $q$-Ohno sum} $    \mathcal{O}^X_q(\kk ;\kl)$ by 
	\begin{align} 
	\mathcal{O}^X_q(\kk ;\kl) =  \sum_{\substack{m_1 > m_2 > \cdots > m_r > 0 \\ n_1 > n_2 > \cdots > n_s >0}} c^X_q(m_1,n_1) \prod_{i=1}^r s^X_q(k_i, m_i) \prod_{j=1}^s s^X_q(l_j, n_j)
	\end{align}
	and set $\mathcal{O}^X_q(\kk; \emptyset) = \mathcal{O}^X_q(\emptyset ;\kk ) = \mathcal{O}^X_q(\kk)$ for admissible $\kk$. The $s^X_q$ and $c^X_q$ are given by \eqref{eq:sxq} and \eqref{eq:cxq}.
\end{definition}

\begin{proposition}[$q$-Transport relations]\label{prop:qtransportrel}
	Let $(k_1,\dots,k_r)$ and $(l_1,\dots,l_s)$ be two indices. \\
	If $s>0$ then we have
	\begin{align*}
	\mathcal{O}^X_q(1,k_1,\dots,k_r ; l_1, \dots, l_s) = \mathcal{O}^X_q(k_1,\dots,k_r ; l_1+1, l_2, \dots, l_s)
	\end{align*}
	and if $r>0$, then we have 
	\begin{align*}
	\mathcal{O}^X_q(k_1+1,k_2,\dots,k_r ; l_1, \dots, l_s) = \mathcal{O}^X_q(k_1,\dots,k_r ; 1, l_1, l_2, \dots, l_s)\,.
	\end{align*}
\end{proposition}
\begin{proof} This proof is Exercise \ref{ex11} and it is similar to the proof of Proposition \ref{thm:transportrel}.
\end{proof}
With the same argument as in \eqref{eq:dualityfromtransport} we see that these transport relations imply Theorem \ref{thm:qohno}, i.e. $\mathcal{O}^X_q(\kk) =     \mathcal{O}^X_q(\kk^\dagger)$.

\begin{remark} We will see in Chapter \ref{sec:mes} that Theorem \ref{thm:qohno} also implies a version of Ohno's relation for our $q$-series $\g(\kk)$. More precisely, we will introduce double-indexed series $\g\bi{\kk}{\kl}$ whose values span the space $\mz_q$, introduced in Section \ref{subsubsec:generalmopdifiedqana}. We then obtain relations in this space since the coefficients of the modified version $(1-q)^{-\wt(\kk)}\mathcal{O}^{(1-q)^{-1}X}_q(\kk)$ are also elements in $\mz_q$.
\end{remark}

\section{The MSW formula and the Drop1 operator}

\subsection{The MSW formula}

The iterated integral expression for multiple zeta values has a finite
counterpart which, surprisingly, is equal to the multiple harmonic sum for
every truncation parameter. This was discovered by Maesaka, Seki and
Watanabe in \cite{MSW}. We state their formula in our convention for the
order of the summation variables.

\begin{definition}
For an index $\kk=(k_1,\dots,k_r)$ and $N\geq 1$, define
\begin{align}\label{eq:defmswflat}
\zeta^\flat_{<N}(\kk)
=\sum_{\substack{
N>n_{1,1}\geq\cdots\geq n_{1,k_1}>
n_{2,1}\geq\cdots\geq n_{2,k_2}>\cdots\\
{}\cdots>n_{r,1}\geq\cdots\geq n_{r,k_r}>0}}
\prod_{i=1}^r
\frac{1}{n_{i,1}\cdots n_{i,k_i-1}(N-n_{i,k_i})}\,.
\end{align}
If $k_i=1$, the product $n_{i,1}\cdots n_{i,k_i-1}$ is understood to be
empty.
\end{definition}

\begin{theorem}[MSW formula, {\cite[Theorem~1.3]{MSW}}]\label{thm:mswformula}
For every index $\kk=(k_1,\dots,k_r)$ and every $N\geq1$, we have
\begin{align}\label{eq:mswformula}
H_{N-1}(\kk)
=\sum_{N>m_1>\cdots>m_r>0}
\frac{1}{m_1^{k_1}\cdots m_r^{k_r}}
=\zeta^\flat_{<N}(\kk)\,.
\end{align}
\end{theorem}

\begin{proof}
We give the generating-series version of the connected-sum proof from
\cite[Lemmas~4.1--4.2]{MSW}, written throughout in our convention. For
$r\geq1$, define the two formal power series
\begin{align*}
F_N(X_1,\dots,X_r)
&=\sum_{N>n_1>\cdots>n_r>0}
\prod_{i=1}^r\frac{1}{n_i-X_i},\\
G_N(X_1,\dots,X_r)
&=\sum_{N=n_0>n_1>\cdots>n_r>0}
\prod_{i=1}^r\left(
\frac{1}{N-n_i}
\frac{(1-X_i)_{n_i-1}}{(n_i-1)!}
\frac{(n_{i-1}-1)!}{(1-X_i)_{n_{i-1}-1}}
\right),
\end{align*}
where $(a)_m=a(a+1)\cdots(a+m-1)$. The coefficient of
$X_1^{k_1-1}\cdots X_r^{k_r-1}$ in $F_N$ is $H_{N-1}(\kk)$. Moreover,
\begin{align}\label{eq:mswpochhammer}
\frac{(1-X)_{n-1}}{(n-1)!}\frac{(m-1)!}{(1-X)_{m-1}}
&=\prod_{a=n}^{m-1}\left(1-\frac{X}{a}\right)^{-1}\\
&=\sum_{d\geq0}X^d
\sum_{m>a_1\geq\cdots\geq a_d\geq n}
\frac{1}{a_1\cdots a_d}.\notag
\end{align}
Taking $m=n_{i-1}$, $n=n_i$ and $d=k_i-1$ in the $i$-th factor of
$G_N$, the variables $a_1,\dots,a_d,n_i$ form precisely the $i$-th block in
\eqref{eq:defmswflat}. Therefore, the same coefficient in $G_N$ is
$\zeta^\flat_{<N}(\kk)$.

It remains to prove $F_N=G_N$. For $0\leq n<N$ and $0<m\leq N$, set
\begin{align*}
C_N(n,m)=\frac{\binom{m-1}{n}}{\binom{N-1}{n}},
\end{align*}
where $C_N(n,m)=0$ if $m\leq n$. For $r,s\geq1$, define
\begin{align*}
&Z_N(X_1,\dots,X_r\mid Y_1,\dots,Y_s)\\
&=\sum_{N=m_0>m_1>\cdots>m_s>n_1>\cdots>n_r>0}
\prod_{i=1}^r\frac{1}{n_i-X_i}\,C_N(n_1,m_s)\\
&\hspace{2.3cm}\cdot
\prod_{j=1}^s\left(
\frac{1}{N-m_j}
\frac{(1-Y_j)_{m_j-1}}{(m_j-1)!}
\frac{(m_{j-1}-1)!}{(1-Y_j)_{m_{j-1}-1}}
\right).
\end{align*}
If one side of the vertical line is empty, we use the boundary values
\begin{align*}
Z_N(X_1,\dots,X_r\mid\emptyset)=F_N(X_1,\dots,X_r),\qquad
Z_N(\emptyset\mid Y_1,\dots,Y_s)=G_N(Y_1,\dots,Y_s).
\end{align*}

The required transport relation follows from the finite identity
\begin{align}\label{eq:mswgentransport}
\sum_{n<a<m}\frac{C_N(a,m)}{a-X}
=\sum_{n<b<m}\frac{C_N(n,b)}{N-b}
\frac{(1-X)_{b-1}}{(b-1)!}
\frac{(m-1)!}{(1-X)_{m-1}}.
\end{align}
To check it, put
\begin{align*}
D_N(a,b;X)
=\frac{(1-X)_{b-1}}{(b-a-1)!}
\frac{(N-a-1)!}{(1-X)_{N-1}}
\end{align*}
and set it equal to zero if $b\leq a$. A direct calculation gives
\begin{align*}
\frac{D_N(a,b+1;X)-D_N(a,b;X)}{a-X}
=\frac{D_N(a-1,b;X)-D_N(a,b;X)}{N-b}.
\end{align*}
Summing this equality over the triangular region $n<a\leq b<m$, cancelling the
common factor
\begin{align*}
\frac{(N-1)!}{(1-X)_{N-1}},
\end{align*}
and dividing by $(1-X)_{m-1}/(m-1)!$ gives
\eqref{eq:mswgentransport}. Inserting this identity into the
connected sum gives
\begin{align*}
Z_N(X,X_1,\dots,X_r\mid Y_1,\dots,Y_s)
=Z_N(X_1,\dots,X_r\mid Y_1,\dots,Y_s,X).
\end{align*}
This also holds in the two boundary cases, since
$C_N(a,N)=C_N(0,b)=1$. Transporting the variables one at a time now yields
\begin{align*}
F_N(X_1,\dots,X_r)
&=Z_N(X_1,\dots,X_r\mid\emptyset)\\
&=Z_N(X_2,\dots,X_r\mid X_1)
=\cdots\\
&=Z_N(\emptyset\mid X_1,\dots,X_r)
=G_N(X_1,\dots,X_r).
\end{align*}
Comparing coefficients proves \eqref{eq:mswformula}.
\end{proof}

The strength of Theorem \ref{thm:mswformula} is that it contains information
about real and finite multiple zeta values at the same time. For an
admissible index of weight $k$, apply the formula to $\kk$ and $\kk^\dagger$.
After the relabelling $n_i\mapsto N-n_{k+1-i}$, the strictly ordered parts of
the two flat sums agree. The remaining boundary sums tend to zero as
$N\to\infty$. This gives the usual duality relation. A precise estimate for
these boundary terms is given in \cite[Lemma~2.1 and the proof of
Theorem~1.1]{MSW}.

On the other hand, set $N=p$ and reduce \eqref{eq:mswformula} modulo $p$.
Using
\begin{align*}
\frac{1}{p-n}\equiv-\frac{1}{n}\pmod p
\end{align*}
for every distinguished factor in \eqref{eq:defmswflat} gives
\begin{align*}
 \za(\kk)=(-1)^{\dep(\kk)}\sum_{\kk\preceq\kl}\za(\kl),
\end{align*}
where $\kk\preceq\kl$ means that $\kk$ is obtained from $\kl$ by replacing
commas by plus signs. This is the non-starred version of Hoffman's finite
duality, which is equivalent to Theorem \ref{thm:fmzvhoffmandual} (see \cite[Theorem~3.1 and its proof]{MSW}). Thus the same finite
identity gives the two different dualities by two different limiting
procedures. Seki also used the MSW formula to prove the extended double
shuffle relations without using iterated integrals (see \cite[Theorem~3.1 and the proof of the extended double shuffle relation]{Se3}).

\subsection{The Drop1 operator}\label{sec:dropone}

The MSW formula is also one of the ingredients in the recent work of Hirose,
Maesaka, Seki and Watanabe \cite{HMSW}. They give an explicit algorithm which
removes all indices equal to $1$ from a multiple zeta value. The corresponding
map will also be important for the derivatives of multiple Eisenstein series.
Put
\begin{align*}
\HH^{\geq2}=\Q\langle z_k\mid k\geq2\rangle\subset\HH^0.
\end{align*}

The convention for the order of an index in \cite{HMSW} is opposite to ours.
We therefore use the equivalent form of the operator given in
\cite[Definition~4.1]{BKM}, which is already written in our convention. Every
non-empty word $w\in\HH^0$ has a unique block decomposition
\begin{align}\label{eq:droponeblocks}
w=x^{c_1}y^{c_2}\cdots x^{c_{2s-1}}y^{c_{2s}},
\qquad c_1,\dots,c_{2s}\geq1.
\end{align}
In terms of indices, this says
\begin{align*}
w=z_{c_1+1}z_1^{c_2-1}\cdots z_{c_{2s-1}+1}z_1^{c_{2s}-1}.
\end{align*}
For example, $z_3z_1z_2=x^2y^2xy$ corresponds to
$\mathbf{c}=(2,2,1,1)$.
For $n\geq0$, write $[n]=\{1,\dots,n\}$, where $[0]=\emptyset$, and set
\begin{align*}
[2s]^1_{\mathbf{c}}&=\{i\in[2s]\mid c_i=1\},&
[2s]^{>1}_{\mathbf{c}}&=\{i\in[2s]\mid c_i>1\}.
\end{align*}
A subset of $[2s]$ is called \emph{even--odd} if it is a union of pairs
$\{2r,2r+1\}$ with $1\leq r<s$, and \emph{odd--even} if it is a union of
pairs $\{2r-1,2r\}$ with $1\leq r\leq s$. If
$A\subset[2s]^1_{\mathbf{c}}$ and
$[2s]\setminus A=\{i_1<\dots<i_t\}$, define
\begin{align*}
\mathbf{c}_{(-A)}&=(c_{i_1},\dots,c_{i_t}),\qquad
\mathbf{c}_{(-A)}-\delta_B
=(c_{i_1}-\delta_{i_1\in B},\dots,c_{i_t}-\delta_{i_t\in B}),
\end{align*}
where $\delta_{i\in B}=1$ if $i\in B$ and $0$ otherwise.

To state the recursion more clearly, put $m(A,B)=\#A+\#B$ and, for
$\nu=1,2$, define
\begin{align*}
\mathcal E_\nu(\mathbf{c})
=\bigl\{(A,B)\ \big|\ &A\subset[2s]^1_{\mathbf{c}}\text{ is even--odd},\quad
B\subset[2s]^{>1}_{\mathbf{c}},\quad m(A,B)\geq\nu,\notag\\[-0.1cm]
&\{2r,2r+1\}\not\subset B\text{ for }1\leq r<s\bigr\}
\end{align*}
and
\begin{align*}
\mathcal O_2(\mathbf{c})
=\bigl\{(A,B)\ \big|\ &A\subset[2s]^1_{\mathbf{c}}\text{ is odd--even},\quad
B\subset[2s]^{>1}_{\mathbf{c}},\quad m(A,B)\geq2,\notag\\[-0.1cm]
&\{2r-1,2r\}\not\subset B\text{ for }1\leq r\leq s\bigr\}.
\end{align*}

\begin{definition}[{\cite[Section~1.2]{HMSW}, in our convention
\cite[Definition~4.1]{BKM}}]\label{def:dropone}
The \emph{Drop1 operator} is the $\Q$-linear map
\begin{align*}
\dropone:\HH^0\longrightarrow\HH^{\geq2}
\end{align*}
defined by $\dropone(1)=\mathfrak D(\emptyset)=1$ and
$\dropone(w)=\mathfrak D(\mathbf{c})$ for the block decomposition
\eqref{eq:droponeblocks}, where
\begin{align}\label{eq:droponeRecursion}
\mathfrak D(\mathbf{c})
={}&\sum_{(A,B)\in\mathcal E_1(\mathbf{c})}
(-1)^{\#B-1}x^{m(A,B)}
\mathfrak D\bigl(\mathbf{c}_{(-A)}-\delta_B\bigr)\notag\\
&+\sum_{(A,B)\in\mathcal E_2(\mathbf{c})}
(-1)^{\#B-1}x^{m(A,B)-1}y
\mathfrak D\bigl(\mathbf{c}_{(-A)}-\delta_B\bigr)\notag\\
&+\sum_{(A,B)\in\mathcal O_2(\mathbf{c})}
(-1)^{\#B}x^{m(A,B)-1}y
\mathfrak D\bigl(\mathbf{c}_{(-A)}-\delta_B\bigr).
\end{align}
\end{definition}

Notice that all factors in front of $\mathfrak D$ in
\eqref{eq:droponeRecursion} are multiplied from the left. This is where the
reversal of the convention in \cite{HMSW} becomes visible. Since $A$ is a
union of pairs and $B\subset[2s]^{>1}_{\mathbf{c}}$, every new tuple again
has an even number of positive entries. Its sum is smaller by
$m(A,B)\geq1$, so the recursion terminates.

For $s=1$, the recursion takes the simpler form
\begin{align*}
\mathfrak D(a,b)={}&x\mathfrak D(a-1,b)+x\mathfrak D(a,b-1)
-x^2\mathfrak D(a-1,b-1)-xy\mathfrak D(a-1,b-1)
\end{align*}
if $a,b>1$, together with
\begin{align*}
\mathfrak D(1,b)=x\mathfrak D(1,b-1),\qquad
\mathfrak D(a,1)=x\mathfrak D(a-1,1),\qquad
\mathfrak D(1,1)=xy.
\end{align*}
For example,
\begin{align*}
\dropone(z_2z_1)&=\mathfrak D(1,2)=x\mathfrak D(1,1)=z_3,\\
\dropone(z_3z_1)&=\mathfrak D(2,2)\notag\\
&=x\mathfrak D(1,2)+x\mathfrak D(2,1)
-x^2\mathfrak D(1,1)-xy\mathfrak D(1,1)
=z_4-z_2z_2.
\end{align*}
More generally, \cite[Theorem~4.7]{BKM} gives, for $p\geq2$,
\begin{align}\label{eq:droponepone}
\dropone(z_pz_1)=z_{p+1}-\sum_{a=2}^{p-1}z_a z_{p+1-a}.
\end{align}
Some examples with more than two indices are
\begin{align*}
\dropone(z_2z_1z_1)&=z_4,&
\dropone(z_2z_1z_2)&=z_2z_3,&
\dropone(z_2z_2z_1)&=z_3z_2,\\
\dropone(z_2z_1z_3)&=z_2z_2z_2+2z_3z_3.
\end{align*}

We briefly explain the main property of this operator. For an admissible index
$\kk=(k_1,\dots,k_r)$ and $N\geq1$, let
$[r]^1_{\kk}=\{i\in[r]\mid k_i=1\}$ and define the \emph{multiple zeta
diamond value}
\begin{align}\label{eq:defzetadiamond}
\zeta_N^\diamondsuit(\kk)
=\sum_{A\subset[r]^1_{\kk}}
\sum_{(n_1,\dots,n_r)\in S_{r,N}(A)}
\left(\prod_{i\in A}\frac{1}{N-n_i}\right)
\left(\prod_{i\in[r]\setminus A}\frac{1}{n_i^{k_i}}\right),
\end{align}
where $S_{r,N}(A)$ consists of all
$(n_1,\dots,n_r)\in[N-1]^r$ satisfying
\begin{align*}
n_{i-1}\geq n_i\quad(i\in A),\qquad
n_{i-1}>n_i\quad(i\in\{2,\dots,r\}\setminus A).
\end{align*}
There is no condition involving $n_0$, since admissibility gives $k_1\geq2$
and therefore $1\notin A$. If all $k_i\geq2$, then
$\zeta_N^\diamondsuit(\kk)$ is just the usual multiple harmonic sum. In
general, $\lim_{N\rightarrow\infty}\zeta_N^\diamondsuit(\kk)=\zeta(\kk)$.
The finite sums already satisfy some multiple zeta value relations. For
example,
\begin{align}\label{eq:zetadiamondeuler}
\zeta_N^\diamondsuit(2,1)
={}&\sum_{N>n_1>n_2>0}\frac{1}{n_1^2n_2}
+\sum_{N>n_1\geq n_2>0}\frac{1}{n_1^2(N-n_2)}
=\sum_{n=1}^{N-1}\frac{1}{n^3}
=\zeta_N^\diamondsuit(3).
\end{align}
This is the finite identity behind
$\dropone(z_2z_1)=z_3$ above.

Define the $\Q$-linear map
\begin{align*}
\zeta^\diamondsuit:\HH^0&\longrightarrow\Q^\N,\\
1&\longmapsto(1,1,\dots),\\
z_{k_1}\cdots z_{k_r}&\longmapsto
\bigl(\zeta_N^\diamondsuit(k_1,\dots,k_r)\bigr)_{N\geq1}.
\end{align*}
\begin{theorem}[{\cite[Theorems~1.8 and~2.5]{HMSW},
see also \cite[Theorem~4.2 and Remark~4.4]{BKM}}]\label{thm:dropone}
For every $w\in\HH^0$, we have
\begin{align}\label{eq:droponezetainvariance}
\zeta^\diamondsuit\bigl(\dropone(w)\bigr)
=\zeta^\diamondsuit(w),\qquad
\zeta\bigl(\dropone(w)\bigr)=\zeta(w).
\end{align}
Moreover, $\dropone(w)=w$ for $w\in\HH^{\geq2}$.
\end{theorem}

\begin{proof}
We just sketch the argument. The diamond values have an MSW-type discrete
iterated integral expression by \cite[Corollary~3.10]{HMSW}. Taking its
difference with respect to $N$ gives the three groups of terms in
\eqref{eq:droponeRecursion} (see \cite[Theorem~4.5]{HMSW}). Induction on the
weight and summation over $N$ then prove the first identity in
\eqref{eq:droponezetainvariance}, as in the proof of
\cite[Theorem~2.5]{HMSW}. For words in $\HH^{\geq2}$, the diamond values are
the ordinary multiple harmonic sums, whose sequences are linearly independent
by \cite[Lemma~2.1]{HMSW}. This gives $\dropone(w)=w$ on
$\HH^{\geq2}$. Finally, taking the limit $N\rightarrow\infty$ gives the
second identity in \eqref{eq:droponezetainvariance}.
\end{proof}

Since the recursion has integer coefficients, Theorem~\ref{thm:dropone} shows
that every multiple zeta value is a $\Z$-linear combination of multiple zeta
values $\zeta(k_1,\dots,k_r)$ with all $k_i\geq2$. The number of such indices
of weight $k$ is the Fibonacci number $F_{k-1}$, where $F_1=F_2=1$. In
particular,
\begin{align*}
\dim_\Q\mz_k\leq F_{k-1}.
\end{align*}
Another useful consequence is
\begin{align}\label{eq:droponeharmonic}
\dropone(u\ast v)=\dropone(u)\ast v
\qquad(u\in\HH^0,\ v\in\HH^{\geq2}),
\end{align}
see \cite[Lemma~4.3]{BKM}, which follows from
\cite[Proposition~3.1]{HMSW}.

The relation space
\begin{align*}
\mathsf{Drop1}
=\big\langle w-\dropone(w)\mid w\in\HH^0\big\rangle_\Q
\end{align*}
is exactly the kernel of the diamond evaluation by
\cite[Theorem~2.6]{HMSW}. It contains the stuffle extension of the linear part
of Kawashima's relations (see \cite[Theorem~1.17]{HMSW}). It is conjectured
that this inclusion is an equality (see \cite[Conjecture~1.16]{HMSW}). We will
use the operator itself in the derivative formula for multiple Eisenstein series in
Section~\ref{sec:mesderiv} and for formal multiple Eisenstein series in
Chapter~\ref{sec:formalspaces}.

\section{Parity}\label{sec:parity}

We already stated the parity theorem for multiple zeta values in the
overview. There are also parity results for finite and symmetric multiple
zeta values, but it is important that the parity condition is different in
the finite case. For the precise ordinary statement, put
\begin{align*}
 \widetilde{\operatorname{Fil}}^{\operatorname D}_r(\mz_k)
 =\sum_{\substack{n\geq0\\2n\leq k}}\zeta(2)^n\fild_r(\mz_{k-2n}).
\end{align*}
Thus powers of $\zeta(2)$ have modified depth zero.

\begin{theorem}[Parity for multiple zeta values]\label{thm:mzvparity}
Let $\kk=(k_1,\dots,k_r)$ be admissible and put $k=\wt(\kk)$. If
$k\not\equiv r\pmod 2$, then
\begin{align}\label{eq:mzvparitymodified}
 \zeta(\kk)\in
 \widetilde{\operatorname{Fil}}^{\operatorname D}_{r-1}(\mz_k).
\end{align}
In particular,
\begin{align}\label{eq:mzvparitymodpi}
\zeta(\kk)\in\fild_{r-1}(\mz_k)+\pi^2\mz_{k-2}.
\end{align}
\end{theorem}

\begin{proof}
Panzer's parity theorem states that
\begin{align*}
 \widetilde{\operatorname{Fil}}^{\operatorname D}_r(\mz_k)
 =\widetilde{\operatorname{Fil}}^{\operatorname D}_{r-1}(\mz_k)
 \qquad\text{if }k\not\equiv r\pmod2;
\end{align*}
see \cite[Theorem~1.1]{Pan}. Since $\zeta(\kk)$ belongs to the left-hand
side, this proves \eqref{eq:mzvparitymodified}. In its decomposition on the
right-hand side, the term with $n=0$ has depth at most $r-1$, whereas every
term with $n\geq1$ is divisible by $\zeta(2)=\pi^2/6$. This proves
\eqref{eq:mzvparitymodpi}, which is also \cite[Proposition~7]{KZ}. Apart from
the exceptional value $\zeta(2)$, the weaker parity theorem modulo
products follows from the linearized regularized double shuffle relations (see \cite[Corollary~8]{IKZ} and \cite[\S4]{Ts}).
\end{proof}

For the finite statement, define
\begin{align*}
\fild_r(\mza_k)
=\left\langle\za(\kl)\mid\wt(\kl)=k,\ \dep(\kl)\leq r\right\rangle_\Q.
\end{align*}

\begin{theorem}[Parity for finite multiple zeta values]\label{thm:fmzvparity}
Let $\kk=(k_1,\dots,k_r)$ be an index of weight $k$. If
$k\equiv r\pmod2$, then
\begin{align}\label{eq:fmzvparity}
\za(\kk)\in\fild_{r-1}(\mza_k)
+\sum_{\substack{k'+k''=k\\k',k''>0\\
r'+r''\leq r,\ 1\leq r',r''<r}}
\fild_{r'}(\mza_{k'})\fild_{r''}(\mza_{k''}).
\end{align}
\end{theorem}

\begin{proof}
Apply the antipode relation \eqref{eq:antipode-relation} for the stuffle
Hopf algebra and then the algebra homomorphism $\za$. This gives
\begin{align}\label{eq:fmzvantipode}
\sum_{i=0}^r
\za(k_1,\dots,k_i)(-1)^{r-i}
\za^\star(k_r,\dots,k_{i+1})=0.
\end{align}
Modulo finite multiple zeta values of length smaller than $r$, a star
value can be replaced by its strict term. The terms with $0<i<r$ in
\eqref{eq:fmzvantipode} are products of lower-weight values. We therefore
obtain
\begin{align*}
\za(k_1,\dots,k_r)+(-1)^r\za(k_r,\dots,k_1)\equiv0
\end{align*}
modulo the right-hand side of \eqref{eq:fmzvparity}. The reversal relation
\begin{align*}
\za(k_r,\dots,k_1)=(-1)^k\za(k_1,\dots,k_r)
\end{align*}
now gives
\begin{align*}
\bigl(1+(-1)^{k+r}\bigr)\za(\kk)\equiv0.
\end{align*}
The coefficient is $2$ precisely when $k\equiv r\pmod2$, which proves the
claim. This is \cite[Proposition~6]{KZ}.
\end{proof}

\begin{ex}
If $a+b+c$ is odd, the depth three case of
\eqref{eq:fmzvantipode} gives the explicit reduction
\begin{align}\label{eq:fmzvparitydepththree}
2\za(a,b,c)+\za(a+b,c)+\za(a,b+c)=0.
\end{align}
Indeed, expand $\za^\star(c,b,a)$ in \eqref{eq:fmzvantipode}, use the
vanishing of depth-one finite multiple zeta values, and apply the reversal
relation to the remaining terms.
\end{ex}

It is conjectured that the product terms in \eqref{eq:fmzvparity} are not
needed, i.e. that $\za(\kk)$ is a linear combination of finite multiple
zeta values of length smaller than $r$ whenever $k\equiv r\pmod2$. This
stronger statement is explicitly formulated in \cite[Remark~1 after Proposition~6]{KZ}.
Thus, unlike \eqref{eq:fmzvparity}, it is not known in general.
In the universal formal finite algebra introduced in
Section~\ref{sec:formalfinitemzv}, the stronger statement can be proved in
depth at most four. See Theorem~\ref{thm:formalfiniteparityfour}.

There is also a particularly clean statement for symmetric multiple zeta
values. The depth of a class in $\mz/\pi^2\mz$ is the smallest depth of a
representative of this class.

\begin{theorem}[Parity for symmetric multiple zeta values]\label{thm:smzvparity}
Let $\kk=(k_1,\dots,k_r)$ be a non-empty index of weight $k$. Then
$\zs(\kk)$ has depth at most $r-1$. If $k\equiv r\pmod2$, then it has
depth at most $r-2$ and it is a linear combination of symmetric multiple
zeta values of length at most $r-1$.
\end{theorem}

\begin{proof}
Choose the stuffle-regularized representative $\zs^\ast(\kk)$. The parity
argument above also holds for stuffle-regularized values (see \cite[proof of Theorem~4]{KZ}). Modulo lower depth and $\pi^2\mz$, a factor whose
weight and depth have opposite parity vanishes. In all remaining depth-$r$
terms, the weight signs in the definition of $\zs^\ast(\kk)$ can therefore
be replaced by the corresponding depth signs. Since a regularized star value
and its strict term differ only in lower depth, the depth-$r$ part of this
strict-term sum is the antipode cancellation in \eqref{eq:fmzvantipode}.
Hence it vanishes, which proves that the depth is at most $r-1$.

If $k\equiv r\pmod2$, then $k$ has opposite parity to $r-1$.
Applying Theorem \ref{thm:mzvparity} once more lowers the remaining
depth-$r-1$ part, so the depth is at most $r-2$. The last assertion is
Yasuda's refinement quoted in the corollary to Theorem~4 of \cite{KZ}.
This proves the theorem (see \cite[Theorem~4 and its corollary]{KZ}).
\end{proof}

For depth one, the last depth bound is interpreted as the vanishing of the
corresponding symmetric multiple zeta value. Notice that the parity conditions
are different for the three spaces. Ordinary multiple zeta values reduce in the
modified depth filtration when weight and depth have opposite parity. Finite
multiple zeta values reduce, modulo products, when they have the same parity.
Symmetric multiple zeta values always lose one depth, and they lose a second
depth in the same-parity case.

\section{The zoo of relations}
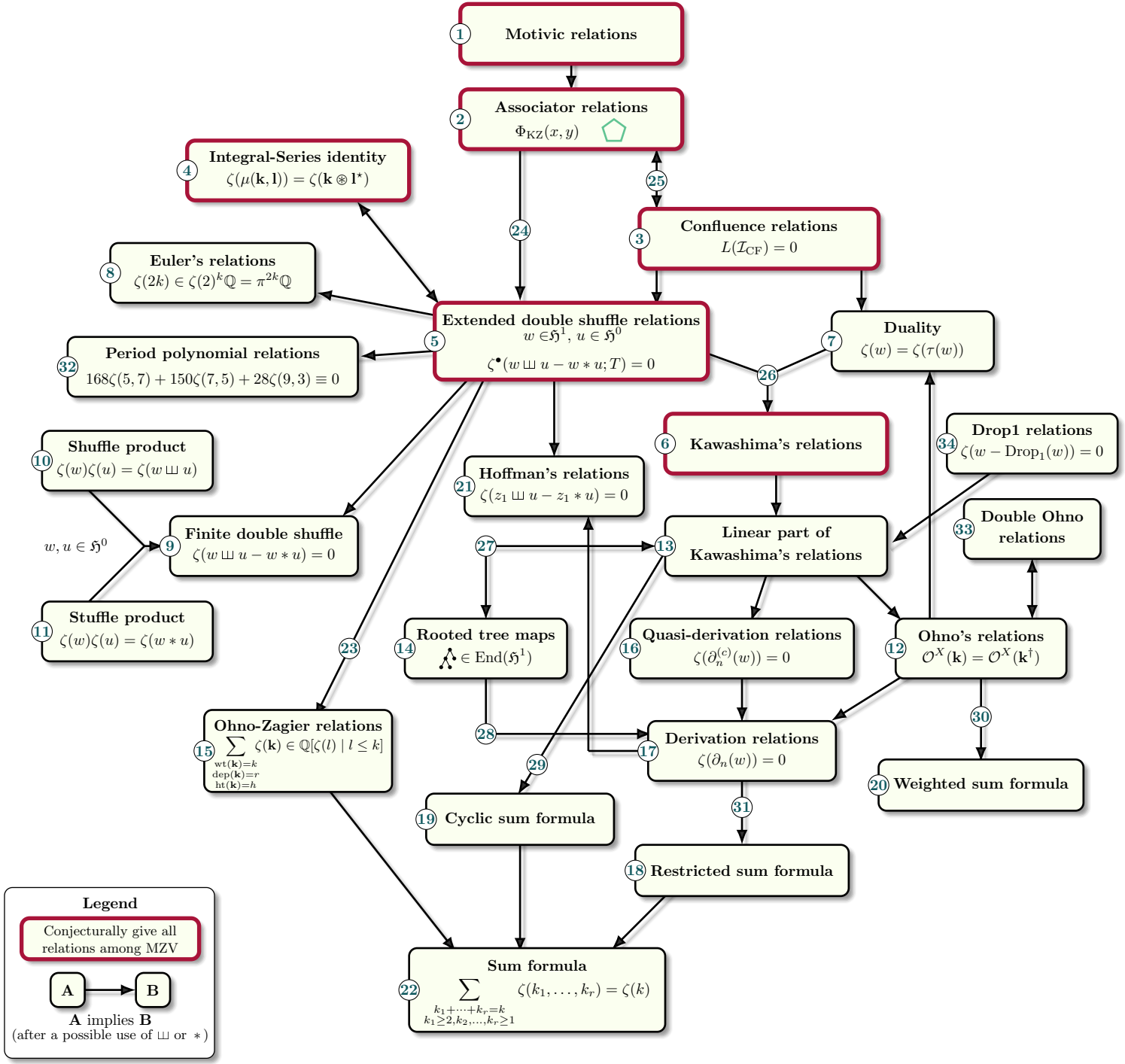
\begin{figure}[ht!]
	\centering
	\begin{adjustwidth}{-1.5cm}{} 
		\input{TikZ/mzv_map}
	\end{adjustwidth}
	\caption{Overview of some relations among multiple zeta values.}
	\label{fig:relationmap}
\end{figure}

We already studied various relations of multiple zeta values. But there are many more relations, which we will not be able to cover in these notes. Instead, we will provide a small overview of some of them and their relationships. In Figure \ref{fig:relationmap}, we give an overview of families of linear and algebraic relations among multiple zeta values and some of their implications. This overview was done together with the help of T. Tanaka, and for a more detailed list of results on multiple zeta values, one should have a look at the list of research papers collected by Hoffman in \cite{H5}.
Some of the presented relations already appeared before, and for some of them, we will provide a few explanations. For the remaining ones, we refer to the literature.

\newpage
In the following, we provide references and some explanations of Figure \ref{fig:relationmap}.
\begin{rel}{\bf Motivic relations:}
	See \cite{Br2,Br3}.
\end{rel}

\begin{rel}{\bf Associator relations:}
	See \cite{D} and \cite{F1}.
\end{rel}

\begin{rel}{\bf Confluence relations:}
	See \cite{HS}.
\end{rel}

\begin{rel}{\bf Integral-Series identity:}
	See \cite{KY}.
\end{rel}

\begin{rel}{\bf Extended double shuffle relations:}
	This is Theorem \ref{thm:eds}.
\end{rel}

\begin{rel}{\bf Kawashima's relations:}
	Define the automorphism $\varphi \in \Aut(\h)$ (with respect to the concatenation) on the generators by $\varphi(x) =  x+y$ and $\varphi(y) = -y$ and define for words $v,w \in \h y$ the operator $z_p v \circledast z_q w = z_{p+q} (v \ast w)$. With this the Kawashima relations can be stated as follows:
	\begin{theorem}[{\cite[Corollary 5.4]{Kaw}}]\label{thm:kawashimarel} For all $v,w \in \h y$ and $m\geq 1$ we have
		\begin{equation}\label{eq:kawashima}
		\sum_{\substack{i+j=m\\i,j \geq 1}} \zeta(\varphi(v) \circledast y^i) \zeta(\varphi(w) \circledast y^j)  = \zeta( \varphi(v \ast w) \circledast y^m) \,.
		\end{equation}
	\end{theorem}
	It is expected that Theorem \ref{thm:kawashimarel} gives all $\Q$-linear relations between multiple zeta values after evaluating the product on the left-hand side by the shuffle product formula. Moreover, numerical experiment suggests that the two cases $m=1,2$ are enough to obtain all linear relations. 
\end{rel}

\begin{rel}{\bf Duality:}
	This is Proposition \ref{prop:duality}.
\end{rel}

\begin{rel}{\bf Euler's relations:}
	This is Proposition \ref{prop:euler}.
\end{rel}

\begin{rel}{\bf Finite double shuffle relations:}
	This is Proposition \ref{prop:fds}.
\end{rel}

\begin{rel}{\bf Shuffle product:}
	This is Corollary \ref{cor:zetashufflemap}.
\end{rel}

\begin{rel}{\bf Stuffle product:}
	This is Corollary \ref{cor:zetastufflemap}.
\end{rel}

\begin{rel}{\bf Ohno's relations:}
	This is Theorem \ref{thm:ohno}.
\end{rel}

\begin{rel}{\bf Linear part of Kawashima's relations:}
	This is the $m=1$ case of Theorem \ref{thm:kawashimarel}. Notice that in this case the sum on the left-hand side  of \eqref{eq:kawashima} is zero and therefore we obtain the linear relation 
	\begin{align*} 
	\zeta( \varphi(v \ast w) \circledast y)=0\,.
	\end{align*}
	The number of relations obtained from the linear part of Kawashima's relations is given by the following table (calculated by Tanaka  in \cite{Tan1}).
	\begin{center}
		\begin{tabular}{|c|c|c|c|c|c|c|c|c|c|c|c|c|c|}
			\hline
			weight $k$         & 3 & 4 & 5 & 6 & 7 & 8 & 9 & 10 & 11 & 12  \\ \hline
			\# all conjectured relations       &1& 3& 6& 14& 29& 60& 123& 249& 503& 1012\\ \hline
			\# linear part of Kawashima's relations & 1 & 2 & 5 & 10 & 23 & 46 & 98 & 200  & 413  & 838  \\ \hline
		\end{tabular}
	\end{center}
\end{rel}

\begin{rel}{\bf Tanaka's rooted tree maps:} Rooted tree maps were introduced by Tanaka in \cite{Tan4}. To a rooted tree he assigns a map $f \in \End(\HH)$, which gives an element $f(w) \in \ker \zeta$ when evaluated at an admissible word $w\in \HH^0$. Before we can give the definition of the rooted tree maps, we need to recall some basics on rooted trees and the Connes-Kreimer coproduct.
	
	A rooted tree is a finite graph which is connected, has no
	cycles, and has a distinguished vertex called the root. We draw rooted trees with the root on top, and we just consider rooted trees with no plane structure, which means that we, for example, do not distinguish between  $\input{TikZ/tree_1_2_2}$ and $\input{TikZ/tree2_1_2_2}$. A product (given by the disjoint union) of rooted trees will be called a (rooted) forest, and by $\HHCK$ we denote the $\Q$-algebra of forests generated by all trees. The unit of $\HHCK$, given by the empty forest, will be denoted by $\II$. Since we just consider trees without plane structure, the algebra $\HHCK$ is commutative.  Due to the work of Connes and Kreimer (\cite{CK}), the space $\HHCK$ has the structure of a Hopf algebra. To define the coproduct on $\HHCK$, we first define the linear map $B_+$ on $\HHCK$, which connects all roots of the trees in a forest to a new root. For example we have $B_+\left(\input{TikZ/tree_1_2} \,\,\input{TikZ/tree_1}\right) = \input{TikZ/tree_1_2_2}$. Clearly for every non-empty tree $t \in \HHCK$ there exists a unique forest $f_t \in \HHCK$ with $t = B_+(f_t)$, which is just given by removing the root of $t$.
	The coproduct on $\HHCK$ can then be defined recursively for a tree $t \in \HHCK$ by
	\[ \Delta(t) = t \otimes \II + (\id \otimes B_{+}) \circ \Delta(f_t)\]
	and for a forest $f=g h$ with $g,h \in \HHCK$ multiplicatively by $\Delta(f) = \Delta(g)\Delta(h)$ and $\Delta(\II)= \II \otimes \II$. For example we have
	\[\Delta(\input{TikZ/tree_1_2}) = \input{TikZ/tree_1_2} \otimes \II + \input{TikZ/tree_1}\,\input{TikZ/tree_1}  \otimes \input{TikZ/tree_1} + 2 \,\input{TikZ/tree_1}  \otimes \input{TikZ/tree_2}+ \II\otimes \input{TikZ/tree_1_2}\,.\]
	In \cite{Tan4}  Tanaka uses the coproduct $\Delta$ to assign to a forest $f\in \HHCK$ a $\Q$-linear map on the space $\HH$, called a rooted tree map, by the following:
	
	\begin{definition}\label{def:rtm} The rooted tree map of the empty forest $\II$ is given by the identity map on $\h$. For any non-empty forest $f \in \HHCK$, we define a $\Q$-linear map on $\HH$, also denoted by $f$, recursively, where we set $f(1)=0$ on the empty word: For a word $w \in \HH$ and a letter $u \in \{x,y\}$ we set
		\begin{equation}\label{eq:defrtm}
		f(w u):= M(\Delta(f)(w \otimes u))\,,
		\end{equation}  
		where $M(w_1 \otimes w_2) = w_1 w_2$ denotes the multiplication on $\h$. This reduces the calculation to $f(u)$ for a letter $u \in \{x,y\}$, which is defined by the following:
		\begin{enumerate}[\textup{(}i\textup{)}] 
			\item If $f=\input{TikZ/tree_1}\,$, then $f(x) := xy$ and $f(y) := -xy$.
			\item For a tree $t = B_+(f)$ we set $t(u) := R_y R_{x+2y}R_y^{-1} f(u)$,
			where $R_v$ is the linear map given by $R_v(w)=wv$ ($v,w \in \h$).
			\item If $f = gh$ is a forest with $g,h \neq \II$, then $f(u):=g(h(u))$.
		\end{enumerate}
	\end{definition}
	
	\begin{theorem}[{\cite[Theorem 1.3]{Tan4}}]\label{thm:rtmrelation}  For any non-empty forest $f \in \HHCK$ we have $$f(\h^0) \subset \ker \zeta.$$
	\end{theorem}

	\begin{ex}
		For the tree $f=\input{TikZ/tree_2}$ and the word $w=xy$ we obtain for $f(w)$
		\begin{align*}
		\input{TikZ/tree_2}(xy) = M(\Delta(\input{TikZ/tree_2})(x \otimes y))= M(\input{TikZ/tree_2}(x) \otimes y + \input{TikZ/tree_1}(x) \otimes \input{TikZ/tree_1}(y) +x \otimes \input{TikZ/tree_2}(y)).
		\end{align*}
		Together with  $\input{TikZ/tree_1}(x)=xy$ and $\input{TikZ/tree_2}(x)=R_y R_{x+2y} R_y^{-1} \input{TikZ/tree_1}(x)=x(x+2y)y$ we get
		\[\input{TikZ/tree_2}(xy) = 2xyyy - xyxy-xxxy-xxyy = 2z_2 z_1 z_1 - z_2 z_2 - z_4 - z_3 z_1 \,,\]
		which by Theorem \ref{thm:rtmrelation} gives the linear relation $2\zeta(2,1,1)=\zeta(4)+\zeta(2,2)+\zeta(3,1)$.
	\end{ex}
	
\end{rel}

\begin{rel}{\bf Ohno-Zagier relations:} 
	We define the \emph{height}  of an index $\kk=(k_1,\dots,k_r)$ by $$\htt(\kk) = \#\{ i \mid k_i > 1\}\,,$$ i.e. the number of $k_i$ not equal to $1$.  The Ohno-Zagier relations give an explicit formula for the sum of all multiple zeta values of a fixed weight, depth and height as a polynomial in single zeta values:
	\begin{theorem}[\cite{OZ}] We have
		\begin{align*}
		\sum_{\substack{k\geq r+h\\r\geq h \geq 1}} \Bigg(\sum_{\substack{\kk \text{ adm.}\\ \wt(\kk)=k\\\dep(\kk)=r \\\htt(\kk)=h}} \zeta(\kk) \Bigg) X^{k-r-h} Y^{r-h} Z^{h-1} = \frac{1}{XY-Z}\left(1 - \exp \left( \sum_{n=2}^\infty \frac{\zeta(n)}{n} S_n(X,Y,Z)\right) \right)\,,
		\end{align*}
		where the $S_n(X,Y,Z) \in \Z[X,Y,Z]$ are defined by 
		\begin{align*}
		S_n(X,Y,Z) = X^n + Y^n - \alpha^n - \beta^n,\qquad \alpha, \beta = \frac{X+Y \pm \sqrt{(X+Y)^2-4Z}}{2}\,.
		\end{align*}
	\end{theorem}
For $n=2,3,4$ the $S_n(X,Y,Z)$ are given by 
\begin{align*}
S_2(X,Y,Z) &= -2XY+2Z,\quad S_3(X,Y,Z) =-3 X^2 Y-3 X Y^2+3 X Z+3 Y Z,\\  S_4(X,Y,Z)&= -4 X^3 Y-6 X^2 Y^2-4 X Y^3+4 X^2 Z+8 X Y Z+4 Y^2 Z-2 Z^2\,.
\end{align*}

	There are also Ohno-Zagier type relations for the $q$-MZV $\zeta_q^{\rm{BZ}}$ proven by Okuda and Takeyama in \cite{OT}.
\end{rel}

\begin{rel}{\bf Quasi-derivation relations:} 
	Quasi-derivation relations    were first proposed in \cite{K1}, and then it was shown in \cite{Tan1} that they give linear relations among multiple zeta values.
	\begin{definition} Let $c\in \Q$ and $H$ the derivation on $\h$ defined by $H(w) = \deg(w) w$ for any $w\in \h$. 
		For an integer $n\geq 1$, the $\Q$-linear map $\partial_n^{(c)} : \h \rightarrow \h$, called \emph{quasi-derivation}, is defined by\footnote{Here $\ad(\theta)(\partial) := \theta \partial - \partial \theta.$}
		\begin{align*}
		\partial_n^{(c)} = \frac{1}{(n-1)!} \ad\left( \theta^{(c)} \right)^{n-1}(\partial_1)\,,
		\end{align*}
		where $ \theta^{(c)} : \h \rightarrow \h$ is the $\Q$-linear map defined by $\theta^{(c)}(x) = \frac{1}{2}(xz + zx)$, $\theta^{(c)}(y) = \frac{1}{2}(yz + zy)$ with $z=x+y$ and the rule
		\begin{align*}
		\theta^{(c)}(w w') = \theta^{(c)}(w) w' + w \theta^{(c)}(w') + c \partial_1(w) H(w')
		\end{align*}
		for any $w,w' \in \h$.
	\end{definition}
	
	In \cite{Tan1} it was shown that $\partial_n^{(c)}$ evaluated at admissible words gives linear relations among multiple zeta values. Further, it was shown that these relations are consequences of the linear part of Kawashima's relations.
	
	\begin{theorem}[{Quasi-derivation relations, \cite[Theorem 3]{Tan1}}]\label{thm:quasideriv} For all $n\geq 1$ and $c\in \Q$ we have $$\partial^{(c)}_n(\h^0) \subset  \ker \zeta\,.$$
	\end{theorem}
	
	\begin{problem} Show that the quasi-derivations $\partial^{(c)}_n$ can be written in terms of rooted tree maps. In \cite{BTan2} it was shown that the linear part of Kawashima's relation ($m=1$ in Theorem \ref{thm:kawashimarel}) is equivalent to the rooted tree map relations (Theorem \ref{thm:rtmrelation}). Since the proof of the quasi-derivation relations uses only the linear part of Kawashima's relations, one might expect that there is an explicit relationship. For $n=1,2,3,4$ one can actually show that we have
		\begin{align*}
		\partial^{(c)}_1 &= \,\input{TikZ/tree_1}\,,\\
		3\partial^{(c)}_2 &=  \left( 2 \,\, \input{TikZ/tree_2} -  \input{TikZ/tree_1}\,\input{TikZ/tree_1}\right) + \left(  \,\, \input{TikZ/tree_2}  + \,\, \input{TikZ/tree_1}\,\input{TikZ/tree_1} \right) c \,,\\
		7 \partial^{(c)}_3 &= \left(3 \,\, \input{TikZ/tree_3}  - 3 \,\,\input{TikZ/tree_2} \, \input{TikZ/tree_1} + \,\input{TikZ/tree_1} \, \input{TikZ/tree_1} \, \input{TikZ/tree_1}\right) + \left( \frac{7}{3} \,\,\input{TikZ/tree_1_2}  + \frac{17}{6}\,\, \input{TikZ/tree_3} - \frac{1}{2}\,\,\input{TikZ/tree_2} \, \input{TikZ/tree_1}  - \,\input{TikZ/tree_1} \, \input{TikZ/tree_1} \, \input{TikZ/tree_1}   \right) c + \left( \frac{7}{6} \,\,\input{TikZ/tree_1_2}  + \frac{2}{3}\,\, \input{TikZ/tree_3} + \frac{1}{2}\,\,\input{TikZ/tree_2} \, \input{TikZ/tree_1}  + \,\input{TikZ/tree_1} \, \input{TikZ/tree_1} \, \input{TikZ/tree_1}   \right) c^2\,,\\
		15 \partial^{(c)}_4 &= \left( 2\, \input{TikZ/tree_1_3} - 4 \, \input{TikZ/tree_1_2_1} + 2 \, \input{TikZ/tree_1_1_2} + 4 \, \input{TikZ/tree_4} - 2 \, \input{TikZ/tree_1_2} \, \input{TikZ/tree_1} - 4 \, \input{TikZ/tree_3} \, \input{TikZ/tree_1} + 4 \, \input{TikZ/tree_2} \, \input{TikZ/tree_1} \, \input{TikZ/tree_1} - \, \input{TikZ/tree_1}\, \input{TikZ/tree_1} \, \input{TikZ/tree_1} \, \input{TikZ/tree_1} \right) \\
		&+\left(\frac{139}{63}\, \input{TikZ/tree_1_3} - \frac{53}{63} \, \input{TikZ/tree_1_2_1} + \frac{52}{9} \, \input{TikZ/tree_1_1_2} + \frac{37}{7} \, \input{TikZ/tree_4} - \frac{52}{9}\, \input{TikZ/tree_1_2}\, \input{TikZ/tree_1}-\frac{22}{7} \, \input{TikZ/tree_3}\, \input{TikZ/tree_1}+\frac{58}{63}\, \input{TikZ/tree_2}\, \input{TikZ/tree_1}\, \input{TikZ/tree_1}+\frac{53}{63}\, \input{TikZ/tree_1}\, \input{TikZ/tree_1}\, \input{TikZ/tree_1}\, \input{TikZ/tree_1}\right)c\\
		&+ \left( \frac{211}{63}\, \input{TikZ/tree_1_3}+ \frac{281}{126} \, \input{TikZ/tree_1_2_1} + \frac{71}{18} \, \input{TikZ/tree_1_1_2} + \frac{31}{14} \, \input{TikZ/tree_4} - \frac{41}{18} \, \input{TikZ/tree_1_2} \, \input{TikZ/tree_1} - \frac{4}{21} \, \input{TikZ/tree_3}\, \input{TikZ/tree_1} + \frac{82}{63}\, \input{TikZ/tree_2}\, \input{TikZ/tree_1}\, \input{TikZ/tree_1}-\frac{88}{63}\, \input{TikZ/tree_1}\, \input{TikZ/tree_1}\, \input{TikZ/tree_1}\, \input{TikZ/tree_1} \right) c^2\\
		&+ \left( \frac{173}{126} \, \input{TikZ/tree_1_3} + \frac{52}{63} \, \input{TikZ/tree_1_2_1} + \frac{7}{9} \, \input{TikZ/tree_1_1_2} + \frac{2}{7} \, \input{TikZ/tree_4} + \frac{1}{18} \, \input{TikZ/tree_1_2}\, \input{TikZ/tree_1}+\frac{4}{21} \, \input{TikZ/tree_3} \, \input{TikZ/tree_1}+ \frac{58}{63} \, \input{TikZ/tree_2} \, \input{TikZ/tree_1} \, \input{TikZ/tree_1} + \frac{53}{63} \, \input{TikZ/tree_1}\, \input{TikZ/tree_1}\, \input{TikZ/tree_1}\, \input{TikZ/tree_1} \right) c^3\,.
		\end{align*}
		But for $n\geq 5$ it is not clear how to express $\partial^{(c)}_n$ in terms of rooted tree maps and for $n\geq 4$ this representation is also not unique, since there are relations among rooted tree maps. For example we have 
		\[ \input{TikZ/tree_2}\,\input{TikZ/tree_2} = 2 \,\,\input{TikZ/tree_1_2_1} + \input{TikZ/tree_1_2}\, \input{TikZ/tree_1}  - \input{TikZ/tree_1_3} - \input{TikZ/tree_1_1_2} \,. \]
	\end{problem}
\end{rel}

\begin{rel}{\bf Derivation relations:} Define for $n\geq1$ the derivation $\partial_n$ on $\h$ by $\partial_n(x) = x(x+y)^{n-1}y$ and $\partial_n(y)=-x(x+y)^{n-1}y$. This is a derivation on $\HH$ with respect to the usual non-commutative multiplication and therefore it suffices to just define it on the generators $x$ and $y$. For any $w=uv$ with $u,v \in \HH$ it is then defined by using Leibniz's rule $\partial_n (uv) = \partial_n(u)v + u \partial_n(v)$.
	
	\begin{theorem}[{Derivation relation, \cite[Corollary~6]{IKZ}}] For all $n\geq 1$ we have $$\partial_n(\h^0) \subset  \ker \zeta\,.$$
	\end{theorem}
	Also notice that this is the $c=0$ case of the quasi-derivation relations (Theorem \ref{thm:quasideriv}), since $\partial^{(0)}_n = \partial_n$.

	\begin{ex} If $n=2$ we have
		\begin{align*}
		\partial_2(x) = -\partial_2(y) = xxy + xyy\,,
		\end{align*}
		i.e. we get for the admissible word $xy$:
		\begin{align*}
		\partial_2(xy) &= \partial_2(x)y + x \partial_2(y)\\
		&= xxyy + xyyy - xxxy - xxyy = xyyy - xxxy = z_{2,1,1} - z_4\,.
		\end{align*}
		The derivation relation then gives the relation $\zeta(2,1,1)=\zeta(4)$.
	\end{ex}
\end{rel}

\begin{rel}{\bf Restricted sum formula:} As a generalization of the sum formula Eie, Liaw, and Ong proved the following formula.
	\begin{theorem}[Restricted sum formula, \cite{ELO}] For integers $p\geq 0$ and $k> r \geq 1$ we have 
		\begin{align*}
		\sum_{\substack{\kk \text{ adm.}\\\wt(\kk)=k\\ \dep(\kk)= r}}  \zeta(\kk, \{1\}^p) = \sum_{\substack{\kl \text{ adm.} \\ \wt(\kl) = k + p\\\dep(\kl)=p+1\\ l_1 > k - r} } \zeta(\kl)\,.
		\end{align*}
	\end{theorem}
	Notice that this reduces to the sum formula (Theorem \ref{thm:sumformula}) in the case $p=0$.
\end{rel}

\begin{rel}{\bf Cyclic sum formula:}
	In \cite{HO} Hoffman and Ohno  proved the following relation.
	\begin{theorem}[Cyclic sum formula] Let $k_1,\dots,k_r \geq 1$ be integers with at least one $k_j \geq 2$. Then
		\begin{align*}
		\sum_{j=1}^r \zeta(k_j+1, k_{j+1}, \dots, k_r, k_1, \dots, k_{j-1}) &= \sum_{\substack{1\leq j \leq r\\ k_j \geq 2}}\sum_{m=0}^{k_j-2} \zeta(k_j-m, k_{j+1}, \dots, k_r, k_1, \dots, k_{j-1}, m+1)\,.
		\end{align*}
	\end{theorem}
	\begin{ex} If we take $k_1=1, k_2=2$ and $k_3=3$ the cyclic sum formula gives
		\begin{align*}
		\zeta(2,2,3) + \zeta(3,3,1)+\zeta(4,1,2) = \zeta(2,3,1,1)+ \zeta(3,1,2,1)+ \zeta(2,1,2,2)\,.
		\end{align*}
	\end{ex}
\end{rel}

\begin{rel}{\bf Weighted sum formula:} There are different types of weighted sum formulas, and one variant can be found in \cite{Kad}.
\end{rel}

\begin{rel}{\bf Hoffman's relations:}
	This is Proposition \ref{prop:hoffmanrel}.
\end{rel}

\begin{rel}{\bf Sum formula:}
	This is Theorem \ref{thm:sumformula}.
\end{rel}

\begin{rel}{\bf The extended double shuffle relations imply the Ohno-Zagier relations:}
	See \cite{L}.
\end{rel}

\begin{rel}{\bf The associator relations imply the extended double shuffle relations:}
	See \cite{F2}. A mould-theoretic generalization of this implication is given in \cite{FHK}.
\end{rel}

\begin{rel}{\bf The associator relations and the confluence relations are equivalent:}
	See \cite{F3}.
\end{rel}

\begin{rel}{\bf The extended double shuffle relations together with the duality imply Kawashima's relations:}
	See \cite{Kaw}.
\end{rel}

\begin{rel}{\bf The rooted tree map relations are equivalent to the linear part of Kawashima's relations:}
	See \cite{BTan2}.
\end{rel} 

\begin{rel}{\bf The rooted tree map relations imply the derivation relation:}
	This was shown in \cite{BTan1} by writing the derivations $\partial_n$ explicitly as rooted tree maps. 
	For this, only trees without any branches are needed, i.e. for $m \geq 1$ consider the ladder trees
	\[ \lambda_m = \input{TikZ/tree_ladder} \]
	and set $\lambda_0 = \II$. With this, the main result of \cite{BTan1} states:
	\begin{theorem}[\cite{BTan1}] For all $n \geq 1$ the derivation $\partial_n$ is given by
		\begin{equation}\label{eq:partialn}
		\partial_n = \frac{n}{2^n-1} \sum_{d=1}^n \frac{(-1)^{d+1}}{d} \sum_{\substack{m_1+\dots+m_d=n\\m_1,\dots,m_d \geq 1}}  \lambda_{m_1} \dots \lambda_{m_d}  \,.
		\end{equation}
	\end{theorem}
	By Definition \ref{def:rtm}(iii) we have $\lambda_{m_1} \lambda_{m_2} = \lambda_{m_2} \lambda_{m_1}$, so we get for the first few values of $n$  ($w \in \h$)
	\begin{align*}\partial_1(w) = \input{TikZ/tree_1}(w)\,,\quad \partial_2(w) = \frac{2}{3}  \,\input{TikZ/tree_2}(w) - \frac{1}{3}\,\input{TikZ/tree_1}\,\input{TikZ/tree_1}(w) \,,\quad
	\partial_3(w) = \frac{3}{7} \,\input{TikZ/tree_3}(w) -  \frac{3}{7} \,\input{TikZ/tree_2}\,\input{TikZ/tree_1}(w)+\frac{1}{7}\input{TikZ/tree_1}\,\input{TikZ/tree_1}\,\input{TikZ/tree_1}(w)\,.
	\end{align*}
\end{rel}

\begin{rel}{\bf The linear part of Kawashima's relations implies the cyclic sum formula:}
	See \cite{TW}.
\end{rel}

\begin{rel}{\bf Ohno's relations imply the weighted sum formula:}
	See \cite{Kad}.
\end{rel}

\begin{rel}{\bf The derivation relations imply the restricted sum formula:}
	See \cite{Tan3}.
\end{rel}

\begin{rel}{\bf Period polynomial relations:} This type of relation will be explained in the next chapter (see Corollary \ref{cor:exoticrel} or \cite{GKZ}).
\end{rel}

\begin{rel}{\bf Double Ohno relations:} See \cite{HMOS} and \cite{HSS}.
\end{rel}

\begin{rel}{\bf Drop1 relations:}
	These are the relations $w-\dropone(w)$ introduced in
	Section~\ref{sec:dropone}. They contain the stuffle extension of the linear
	part of Kawashima's relations, while equality is conjectural. This gives the
	corresponding arrow in Figure~\ref{fig:relationmap}.
\end{rel}

\section{Schur multiple zeta values}\label{sec:schurmzv}

We finish this chapter by returning to the Schur multiple zeta values discussed in Section~\ref{subsec:schurmzvoverview}. For this, we introduce the partition notation for skew Young diagrams and prove the Jacobi--Trudi formula for their Schur multiple zeta values.
Recall that they were introduced by Nakasuji, Phuksuwan and
Yamasaki in \cite{NPY}. They are indexed by skew Young diagrams whose boxes
are filled with positive integers, and contain
multiple zeta values and multiple zeta-star values as special cases. An
interpolated version for straight Young diagrams and a lattice path proof of
its Jacobi--Trudi formula are given in \cite{Ba4}.

Let $\mu\subseteq\lambda$ be partitions. Append zero parts as necessary and
write $\lambda=(\lambda_1,\dots,\lambda_d)$ and
$\mu=(\mu_1,\dots,\mu_d)$ with $\mu_i\leq\lambda_i$. The
associated skew Young diagram is
\begin{align*}
D(\lambda/\mu)
=\{(i,j)\mid 1\leq i\leq d,\ \mu_i<j\leq\lambda_i\}.
\end{align*}
For example, for $\lambda=(3,2)$ and $\mu=(1,0)$ we get
\begin{align*}
D(\lambda/\mu)=
{\ytableausetup{centertableaux,boxsize=1.4em}
\begin{ytableau}
\none & {} & {}\\
{} & {}
\end{ytableau}}.
\end{align*}
The missing box in the upper left corner belongs to $\mu$.
A \emph{semistandard Young tableau} of shape $\lambda/\mu$ is a collection
of positive integers
\begin{align*}
\mathbf m=(m_{i,j})_{(i,j)\in D(\lambda/\mu)}
\end{align*}
which is weakly increasing along rows and strictly increasing along columns,
i.e.
\begin{align*}
m_{i,j}\leq m_{i,j+1},\qquad m_{i,j}<m_{i+1,j}
\end{align*}
whenever the two entries in question belong to the diagram.
For the above shape, for example,
\begin{align*}
{\ytableausetup{centertableaux,boxsize=1.4em}
\begin{ytableau}
\none & 2 & 3\\
1 & 4
\end{ytableau}}
\end{align*}
is semistandard.

\begin{definition}
Let $\mathbf h=(h_{i,j})_{(i,j)\in D(\lambda/\mu)}$ be a collection of
positive integers. The corresponding \emph{Schur multiple zeta value} is
\begin{align}\label{eq:defschurmzv}
\zeta(\mathbf h)
=\sum_{\mathbf m\in\operatorname{SSYT}(\lambda/\mu)}
\prod_{(i,j)\in D(\lambda/\mu)}\frac{1}{m_{i,j}^{h_{i,j}}}.
\end{align}
\end{definition}

We will also write the filling $\mathbf h$ directly into its Young diagram in
the argument of $\zeta$.
A box $(i,j)$ is called a corner if neither $(i+1,j)$ nor $(i,j+1)$
belongs to the diagram. The series \eqref{eq:defschurmzv} converges
absolutely if $h_{i,j}\geq2$ at every corner and $h_{i,j}\geq1$ at all
other boxes (see \cite[Lemma~2.1 and Section~4.1]{NPY}). Because we use decreasing
summation variables, the first non-trivial column and row examples are
\begin{align*}
\zeta\left(
{\ytableausetup{centertableaux,boxsize=1.4em}
\begin{ytableau}
h_1\\h_2\\h_3
\end{ytableau}}
\right)
&=\zeta(h_3,h_2,h_1),\\
\zeta\left(
{\ytableausetup{centertableaux,boxsize=1.4em}
\begin{ytableau}
h_1&h_2&h_3
\end{ytableau}}
\right)
&=\zeta^\star(h_3,h_2,h_1).
\end{align*}

The analogue of the Jacobi--Trudi formula holds when the entries are
constant along diagonals. The skew version below was proved in
\cite[Theorem~4.3~(2)]{NPY}. For straight diagrams it is also the $t=0$
specialization of the interpolated formula in \cite[Theorem~1.1]{Ba4}.
Let $\lambda'$ and $\mu'$ denote the conjugate partitions and set $s=\lambda_1$.

\begin{theorem}[Jacobi--Trudi formula]
\label{thm:schurjacobitrudi}
Assume that
\begin{align*}
h_{i,j}=a_{j-i}
\end{align*}
for a family of integers $a_n\geq2$ indexed by $n\in\Z$. Put
\begin{align*}
L_{i,j}=\lambda'_i-\mu'_j-i+j.
\end{align*}
Then
\begin{align}\label{eq:schurjacobitrudi}
\zeta(\mathbf h)
=\det\left[
\zeta\bigl(a_{i-\lambda'_i},a_{i-\lambda'_i+1},
\dots,a_{j-\mu'_j-1}\bigr)
\right]_{1\leq i,j\leq s}.
\end{align}
An entry in the determinant is defined to be $1$ if $L_{i,j}=0$ and
$0$ if $L_{i,j}<0$.
\end{theorem}

\begin{proof}
We follow the lattice path proof of \cite[Theorem~4.7]{Ba4}. In the case $t=0$
the graph becomes particularly simple, and the choice of the initial points
below also gives the skew version.
For $N\geq1$, first restrict all tableau entries to be at most $N$. Consider
the directed graph with vertices
\begin{align*}
(c,m)\in\Z\times\{0,1,\dots,N\}.
\end{align*}
For $1\leq m\leq N$, it has vertical edges from $(c,m-1)$ to $(c,m)$
of weight $1$ and diagonal edges
\begin{align*}
(c,m-1)\longrightarrow(c+1,m)
\end{align*}
of weight $m^{-a_{-c}}$. For $1\leq i,j\leq s$, set
\begin{align*}
A_j=(\mu'_j-j+1,0),\qquad
B_i=(\lambda'_i-i+1,N).
\end{align*}
Every edge increases the second coordinate, and all relevant paths lie in a
finite strip. Thus it suffices to consider a finite acyclic subgraph.
Put $\alpha_j=j-\mu'_j-1$ and $\beta_i=i-\lambda'_i$. A path from $A_j$
to $B_i$ has $L_{i,j}=\alpha_j-\beta_i+1$ diagonal edges. If $L_{i,j}=0$,
there is one vertical path of weight $1$, and if $L_{i,j}<0$, there is no
such path. If $L_{i,j}>0$, the heights of its diagonal edges are strictly
increasing. Reading them in reverse order shows that the generating series
of these paths is
\begin{align*}
\zeta_{\leq N}\bigl(a_{\beta_i},a_{\beta_i+1},\dots,a_{\alpha_j}\bigr),
\end{align*}
where the subscript means that all summation variables are at most $N$.

For example, for $\lambda=(2,1)$, $\mu=\emptyset$ and $N=4$, the tableau
\begin{align*}
{\ytableausetup{centertableaux,boxsize=1.4em}
\begin{ytableau}
1&2\\3
\end{ytableau}}
\end{align*}
gives the following two paths. The diagonal steps of the solid path occur at
heights $1,3$ and record the first column, while the diagonal step of the
dashed path at height $2$ records the second column.
\begin{center}
\begin{tikzpicture}[x=0.75cm,y=0.55cm,>=stealth]
\foreach \m in {0,...,4}{\draw[gray!25] (-1.25,\m)--(2.25,\m);}
\foreach \x in {-1,...,2}{\foreach \m in {0,...,4}{\fill[gray!55] (\x,\m) circle (1.1pt);}}
\draw[->,line width=1pt,color=reallinecol]
 (0,0)--node[midway,below right] {\scriptsize $1$}(1,1)--(1,2)--node[midway,below right] {\scriptsize $3$}(2,3)--(2,4);
\draw[->,line width=1pt,dashed,color=numcol2]
 (-1,0)--(-1,1)--node[midway,below right] {\scriptsize $2$}(0,2)--(0,3)--(0,4);
\node[below] at (0,0) {\scriptsize $A_1$};
\node[below] at (-1,0) {\scriptsize $A_2$};
\node[above] at (2,4) {\scriptsize $B_1$};
\node[above] at (0,4) {\scriptsize $B_2$};
\end{tikzpicture}
\end{center}

By the Lindstr\"om--Gessel--Viennot lemma \cite{GV}, the determinant of
these path generating series is the signed sum over vertex-disjoint path
families. Both endpoint sequences are strictly ordered by their first
coordinates, so only the identity permutation can occur. The heights of the
diagonal edges of the $j$-th path
give the entries in the $j$-th column of a tableau. The strict increase of
these heights gives the strict increase down columns. Moreover, two neighbouring
paths are disjoint if and only if the corresponding entries are weakly
increasing along rows. This gives a
weight-preserving bijection with the semistandard Young tableaux of shape
$\lambda/\mu$ whose entries are at most $N$.

Thus the truncated version of \eqref{eq:defschurmzv} is the determinant in
\eqref{eq:schurjacobitrudi} with all multiple zeta values truncated at $N$.
Since all $a_n\geq2$, we can let $N\to\infty$ and obtain the formula.
\end{proof}

\begin{ex}
Here are three examples of the determinant formula.
\begin{enumerate}[\textup{(}i\textup{)}]
\item For a row of length three we get
\begin{align*}
\zeta\left(
{\ytableausetup{centertableaux,boxsize=1.4em}
\begin{ytableau}
a_0&a_1&a_2
\end{ytableau}}
\right)
&=\det\begin{pmatrix}
\zeta(a_0)&\zeta(a_0,a_1)&\zeta(a_0,a_1,a_2)\\
1&\zeta(a_1)&\zeta(a_1,a_2)\\
0&1&\zeta(a_2)
\end{pmatrix}\\
&=\zeta^\star(a_2,a_1,a_0).
\end{align*}
This is the $t=0$ case of \cite[Example~1.3~(i)]{Ba4}, written in our
convention.

\item For $\lambda=(2,1)$ and $\mu=\emptyset$ the formula gives
\begin{align*}
\zeta\left(
{\ytableausetup{centertableaux,boxsize=1.8em}
\begin{ytableau}
a_0&a_1\\a_{-1}
\end{ytableau}}
\right)
&=\det\begin{pmatrix}
\zeta(a_{-1},a_0)&\zeta(a_{-1},a_0,a_1)\\
1&\zeta(a_1)
\end{pmatrix}\\
&=\zeta(a_{-1},a_0)\zeta(a_1)-\zeta(a_{-1},a_0,a_1).
\end{align*}
In particular, if all entries are equal to $2$, then
\begin{align*}
\zeta\left(
{\ytableausetup{centertableaux,boxsize=1.4em}
\begin{ytableau}
2&2\\2
\end{ytableau}}
\right)
=\zeta(2,2)\zeta(2)-\zeta(2,2,2)=\frac{\pi^6}{840}.
\end{align*}

\item For the skew shape $\lambda=(3,2)$, $\mu=(1,0)$ and constant entries
equal to $2$, we have $\lambda'=(2,2,1)$ and $\mu'=(1,0,0)$. Therefore
\begin{align*}
\zeta\left(
{\ytableausetup{centertableaux,boxsize=1.4em}
\begin{ytableau}
\none&2&2\\2&2
\end{ytableau}}
\right)
&=\det\begin{pmatrix}
\zeta(2)&\zeta(2,2,2)&\zeta(2,2,2,2)\\
1&\zeta(2,2)&\zeta(2,2,2)\\
0&1&\zeta(2)
\end{pmatrix}\\
&=\frac{61\pi^8}{362880}.
\end{align*}
\end{enumerate}
\end{ex}

\subsection{Checkerboard style Schur multiple zeta values}

Let $a,b\geq1$ with $b\geq2$. A skew Young diagram is called
\emph{checkerboardable} if its boxes can be filled with $a$ and $b$ such
that the entries alternate along rows and columns and every corner contains
$b$. This filling is unique. We write $\zeta(a,b;\lambda/\mu)$ for the
corresponding Schur multiple zeta value.

Put $\delta_n=(n,n-1,\ldots,2,1)$ and
$\delta_0=\delta_{-1}=\emptyset$. The following four families of ribbons
will be useful:
\begin{align*}
A_{a,b}(n)
&=\zeta\bigl(a,b;(n+1,n+1,n,n-1,\ldots,3,2)/\delta_n\bigr),\\
B_{a,b}(n)
&=\zeta\bigl(a,b;\delta_{n+1}/\delta_{n-1}\bigr),\\
S_{a,b}(n)
&=\zeta\bigl(a,b;(n,n,n-1,\ldots,2,1)/\delta_{n-1}\bigr),\\
S^\star_{a,b}(n)
&=\zeta\bigl(a,b;(n+1,n,\ldots,3,2)/\delta_{n-1}\bigr).
\end{align*}
Here $A_{a,b}(n)$ is defined for $n\geq1$, $B_{a,b}(n)$ for $n\geq0$,
and the two $S$-families for $n\geq1$.

\begin{theorem}[{\cite[Theorems~3.4 and~3.5]{BY1}}]
For $n\geq1$ we have
\begin{align*}
S_{1,3}(n)
&=\frac{1}{4^n}\zeta^\star(\{4\}^n),\\
S^\star_{1,3}(n)
&=\sum_{j=0}^n\frac{1}{4^j}
\zeta^\star(\{4\}^j)\zeta(\{4\}^{n-j}),\\
A_{1,3}(n)
&=\frac{2}{4^n}\zeta(4n+1).
\end{align*}
For $n\geq0$ we also have
\begin{align*}
B_{1,3}(n)=\frac{1}{4^n}\zeta(4n+3).
\end{align*}
In particular, the $A$- and $B$-families together give all odd single zeta
values.
\end{theorem}

\begin{proof}[Idea of proof]
The harmonic product gives recursions for the four families. The generating
series of the alternating multiple polylogarithms factors into two Gauss
hypergeometric functions. Gauss's evaluation at $1$ gives the formulas for
$S_{1,3}$ and $S^\star_{1,3}$. For $A_{1,3}$ one also uses the asymptotic
formula for a zero-balanced hypergeometric function and the expansion of the
digamma function. The formula for $B_{1,3}$ follows from the evaluation of
$\zeta(3,\{1,3\}^n)$. See \cite[Sections~3.2 and~3.3]{BY1} for the details.
The convention for the order of ordinary multiple zeta indices in that paper
is opposite to ours, but the formulas above are unchanged.
\end{proof}

For example, the first two $A$- and $B$-ribbons give the following nice
pictures for four odd zeta values:
\begin{align*}
\zeta(5)
&=2\zeta\left(
{\ytableausetup{centertableaux,boxsize=1.2em}
\begin{ytableau}
\none&1\\
1&3
\end{ytableau}}
\right),&
\zeta(7)
&=4\zeta\left(
{\ytableausetup{centertableaux,boxsize=1.2em}
\begin{ytableau}
1&3\\
3
\end{ytableau}}
\right),\\
\zeta(9)
&=8\zeta\left(
{\ytableausetup{centertableaux,boxsize=1.2em}
\begin{ytableau}
\none&\none&1\\
\none&1&3\\
1&3
\end{ytableau}}
\right),&
\zeta(11)
&=16\zeta\left(
{\ytableausetup{centertableaux,boxsize=1.2em}
\begin{ytableau}
\none&1&3\\
1&3\\
3
\end{ytableau}}
\right)\,.
\end{align*}

The checkerboard example in Section~\ref{subsec:schurmzvoverview} is the
first case of the following Hankel determinant.

\begin{corollary}[{\cite[Corollary~4.6]{BY1}}]
Let $N\geq1$ and $0\leq n\leq N-1$. If $N\equiv n\pmod 2$, then
\begin{align*}
\zeta(1,3;\delta_N/\delta_n)
=4^{-\frac{(N+n)(N-n)}{4}}
\det\left[\zeta\bigl(4(i+j+n)-1\bigr)\right]_{
1\leq i,j\leq\frac{N-n}{2}}.
\end{align*}
If $N\not\equiv n\pmod 2$, then
\begin{align*}
\zeta(1,3;\delta_N/\delta_n)
&=(-1)^{\frac{n(n+1)}{2}}
4^{-\frac{(N+n+1)(N-n-1)}{4}}\\
&\quad\cdot
\det\left[\zeta\bigl(4(i+j-n)-5\bigr)\right]_{
1\leq i,j\leq\frac{N+n+1}{2}},
\end{align*}
where $\zeta(r)=0$ for $r<0$.
In particular, these values belong to
$\Q[\zeta(4m+3)\mid m\geq0]$.
\end{corollary}

\begin{proof}
Apply the admissible version of the column Jacobi--Trudi formula and the
matrix reduction in \cite[Section~4]{BY1}, and then use
$B_{1,3}(n)=4^{-n}\zeta(4n+3)$.
\end{proof}

The author and Charlton showed that every admissible checkerboard Schur
multiple zeta value with entries $1$ and $3$ belongs to
\begin{align*}
\Q[\pi^4,\zeta(3),\zeta(5),\zeta(7),\ldots]
\end{align*}
in \cite[Theorem~4.3]{BC}. Their generalized regularized Jacobi--Trudi
formula uses suitable stair-shaped outside decompositions. It shows more
precisely that pure $S$- or $S^\star$-stair tessellations give elements of
$\Q[\pi^4]$, while pure $A$- and $B$-stair tessellations give elements of
$\Q[\zeta(4n+1)\mid n\geq1]$ and
$\Q[\zeta(4n+3)\mid n\geq0]$, respectively
\cite[Proposition~4.4 and Corollary~4.5]{BC}.

The formula in Theorem~\ref{thm:schurjacobitrudi} is also the zeta
realization of the Jacobi--Trudi identity in the Young tableaux Hopf algebra
of \cite[Theorem~2.24 and Corollary~4.2]{Yu}. The integer-entry case of
Yu's multiple Schur series is the following.

\begin{definition}[{\cite[Definition~3.1]{Yu}}]
Let $\mathcal X=(X,\prec)$ be a finite totally ordered set, let $R$ be a
commutative $\Q$-algebra and choose an element $f_x(h)\in R$ for every
$x\in X$ and $h\geq1$. For a filling
$\mathbf h=(h_{i,j})_{(i,j)\in D(\lambda/\mu)}$ by positive integers, set
\begin{align*}
F_{\mathcal X,f}(\mathbf h)
=\sum_{\mathbf x\in\operatorname{SSYT}(\lambda/\mu,\mathcal X)}
\prod_{(i,j)\in D(\lambda/\mu)}f_{x_{i,j}}(h_{i,j})\,.
\end{align*}
Here $\operatorname{SSYT}(\lambda/\mu,\mathcal X)$ consists of the fillings
by elements of $X$ which are weakly increasing along rows and strictly
increasing down columns with respect to $\prec$. We call
$F_{\mathcal X,f}(\mathbf h)$ a \emph{multiple Schur series} and set
$F_{\mathcal X,f}(\emptyset)=1$.
\end{definition}

For $X=\{1,\ldots,N\}$ with its usual order and $f_m(h)=m^{-h}$, letting
$N\to\infty$ gives the Schur multiple zeta value in
\eqref{eq:defschurmzv}. Other choices give Schur multiple $L$-values and
$q$-analogues.

For the Eisenstein realization, let $\tau\in\Ha$ and order the lattice by
\begin{align*}
m_1\tau+n_1\prec_\tau m_2\tau+n_2
\quad\Longleftrightarrow\quad
m_1<m_2\quad\text{or}\quad(m_1=m_2\text{ and }n_1<n_2)\,.
\end{align*}
For $M,N\geq1$ put
\begin{align*}
X_{M,N}^{\tau,+}
=\{m\tau+n\mid m,n\in\Z,\ |m|\leq M,\ |n|\leq N,
\ 0\prec_\tau m\tau+n\}\,.
\end{align*}
Thus $\lambda\prec_\tau\mu$ is the same as $\mu\succ\lambda$ for the order
used in Chapter~\ref{sec:mes}.

\begin{definition}[{\cite[Definition/Proposition~4.7]{Yu}}]
Write
$\mathcal T_{M,N}^\tau(\lambda/\mu)
=\operatorname{SSYT}(\lambda/\mu,X_{M,N}^{\tau,+})$.
Let all entries of $\mathbf h$ be at least $2$. Define the truncated series by
\begin{align*}
\aG^{\mathrm{Sch}}_{M,N}(\mathbf h;\tau)
&=\sum_{\boldsymbol\ell\in
\mathcal T_{M,N}^\tau(\lambda/\mu)}
\prod_{(i,j)\in D(\lambda/\mu)}
\frac{1}{\ell_{i,j}^{h_{i,j}}}\,.
\end{align*}
The iterated limit
\begin{align*}
\lim_{M\to\infty}\lim_{N\to\infty}
\aG^{\mathrm{Sch}}_{M,N}(\mathbf h;\tau)\,.
\end{align*}
exists. It is called the \emph{Schur multiple Eisenstein series} and is
denoted by $\aG^{\mathrm{Sch}}(\mathbf h;\tau)$.
\end{definition}

The inner $N$-limit is taken first, just as for the Eisenstein summation in
\eqref{def:meseisensteinsum}. A column with entries
$h_1,\ldots,h_r$ from top to bottom gives
$\aG_{h_r,\ldots,h_1}(\tau)$ in the notation of Chapter~\ref{sec:mes}.
The Fourier expansion of $\aG^{\mathrm{Sch}}(\mathbf h;\tau)$ has constant
term $\zeta(\mathbf h)$ (see \cite[(4.3)]{Yu}).

When all entries are $2$, the author and Yu obtained the following
particularly nice consequence.

\begin{theorem}[{\cite[Theorem~2.3 and Proposition~2.5]{BY2}}]
Let $d=|D(\lambda/\mu)|$ and let $\{2\}^{\lambda/\mu}$ denote the filling
of $\lambda/\mu$ whose entries are all equal to $2$. Then
\begin{align*}
G^{\mathrm{Sch}}_{\lambda/\mu}(\tau)
:=(-2\pi i)^{-2d}
\aG^{\mathrm{Sch}}\bigl(\{2\}^{\lambda/\mu};\tau\bigr)
\in\qmf_{2d}\,.
\end{align*}
In particular, this is a homogeneous quasimodular form of weight $2d$.
For every $n\geq1$, the series $G^{\mathrm{Sch}}_\lambda$ with
$\lambda\vdash n$ span $\qmf_{2n}$ over $\Q$.
\end{theorem}

\vspace{1cm}
\begin{center}
\bf{ \Large {\color{qmzvlinecol}\ding{94}}~ Exercises ~{\color{qmzvlinecol}\ding{94}}}\end{center}

The following is a collection of exercises intended to deepen the reader's understanding of this chapter.

\begin{exer}\label{ex11}
\begin{enumerate}[\textup{(}i\textup{)}] 
    \item Show that $z_{\kk^\dagger} = \tau(z_\kk)$ for any admissible index $\kk$. 
    \item Give the proof of Proposition \ref{prop:qtransportrel}.
\end{enumerate}
\end{exer}

%% file: TikZ/mzv_map.tex
\tikzset{%
  block/.style={draw, top color=white!90!#1, bottom color=white, shading angle=45, text width=18pt, text centered, rounded corners, minimum height=10pt},
}

\begin{tikzpicture}[scale=0.3, >=Latex,every node/.style={scale=0.77}]

\def\linewidth{1pt}
\def\linewidthbig{2pt}
\def\arrowwidth{1pt}

\coordinate (Z) at (0,0);

\def\circrad{0.6}

\coordinate (PR1) at (32,58);

\coordinate (PR2) at (32,53);

\coordinate (PR3) at (43,46);

\coordinate (PR4) at (16,50);

\coordinate (PR5) at (32,40);

\coordinate (PR6) at (44,34);

\coordinate (PR7) at (52,40);

\coordinate (PR8) at (11,44);

\coordinate (PR9) at (14,28);

\coordinate (PR10) at (6,33);

\coordinate (PR11) at (6,23);

\coordinate (PR12) at (56,22);

\coordinate (PR13) at (44,28);

\coordinate (PR14) at (27,22);

\coordinate (PR15) at (16,16);

\coordinate (PR16) at (42,22);

\coordinate (PR17) at (42,16);

\coordinate (PR18) at (42,9);

\coordinate (PR19) at (29,12);

\coordinate (PR20) at (56,14);

\coordinate (PR21) at (31,31.5);

\coordinate (PR22) at (30,2);

\coordinate (PR23) at (11,38.5);

\coordinate (PR24) at (59,29);

\coordinate (PR25) at (59,34);

\newdimen\mydim
 \newcommand\getx[1]{
      \pgfextractx\mydim{\pgfpointanchor{#1}{center}}
 }

\draw[line width=\arrowwidth,shorten >=0.5cm,blur shadow,->] (PR1) -- (PR2);

\draw[line width=\arrowwidth,shorten >=0.5cm,shorten <=0.5cm,blur shadow,<->] ($(PR2)+(5,0)$) -- ($(PR3)+(-6,0)$);

\coordinate (CI) at ($(PR3)+(-6,3.4)$);
\draw[black,fill=white] (CI) circle (\circrad);
\node[text=numcol] at (CI) {{\bf 25}};

\draw[line width=\arrowwidth,shorten >=0.7cm, blur shadow,->]  ($(PR2)+(-3,0)$) -- ($(PR5)+(-3,0)$);

\coordinate (CI) at ($(PR5)+(-3,6.5)$);
\draw[black,fill=white] (CI) circle (\circrad);
\node[text=numcol] at (CI) {{\bf 24}};

\draw[line width=\arrowwidth,shorten >=0.83cm,shorten <=0.7cm,blur shadow,<->]  ($(PR4)+(2,0)$) -- ($(PR5)+(-6,0)$);

\draw[line width=\arrowwidth,shorten >=0.6cm,blur shadow,->]  ($(PR3)+(-6,0)$) -- ($(PR5)+(5,0)$);

\draw[line width=\arrowwidth,shorten >=0.5cm,blur shadow,->]  ($(PR3)+(6,0)$) -- ($(PR7)+(-3,0)$);

\coordinate (EDK) at ($(PR6)+(-0.5,4)$);
\draw[line width=\arrowwidth,blur shadow]  ($(PR5)+(6,0)$) -- (EDK);
\draw[line width=\arrowwidth,blur shadow]  ($(PR7)+(-4,0)$) -- (EDK);
\draw[line width=\arrowwidth,shorten >=0.5cm,blur shadow,->]  (EDK) -- ($(PR6)+(-0.5,0)$);

\draw[black,fill=white] (EDK) circle (\circrad);
\node[text=numcol] at (EDK) {{\bf 26}};

\draw[line width=\arrowwidth,shorten >=1.86cm,blur shadow,->]  ($(PR5)$) -- ($(PR8)$);

\draw[line width=\arrowwidth,shorten >=2.6cm,blur shadow,->]  ($(PR5)$) -- ($(PR23)$);

\coordinate (SHF) at ($(PR9)-(7,0)$);
\draw[line width=\arrowwidth,blur shadow]  ($(PR10)+(-4,0)$) -- (SHF);
\draw[line width=\arrowwidth,blur shadow]  ($(PR11)+(-4,0)$) -- (SHF);
\draw[line width=\arrowwidth,shorten >=1.15cm,blur shadow,->]  (SHF) -- ($(PR9)+(-2,0)$);

\draw[line width=\arrowwidth,shorten >=0.7cm,blur shadow,->]  ($(PR5)+(-4,0)$) -- ($(PR9)+(3,0)$);

\draw[line width=\arrowwidth,shorten >=0.65cm,blur shadow,->]  ($(PR5)+(-4,0)$) -- ($(PR15)$);

\coordinate (CI) at (19,22);
\draw[black,fill=white] (CI) circle (\circrad);
\node[text=numcol] at (CI) {{\bf 23}};

\draw[line width=\arrowwidth,shorten >=0.5cm,blur shadow,->]  ($(PR6)$) -- ($(PR13)$);

\draw[line width=\arrowwidth,shorten >=0.5cm,blur shadow,->]  ($(PR13)$) -- ($(PR16)$);

\draw[line width=\arrowwidth,shorten >=0.5cm,blur shadow,->]  ($(PR16)$) -- ($(PR17)$);

\draw[line width=\arrowwidth,shorten >=0.45cm,blur shadow,->]  ($(PR17)$) -- ($(PR18)$);

\coordinate (CI) at ($(PR18)+(0,3.7)$);
\draw[black,fill=white] (CI) circle (\circrad);
\node[text=numcol] at (CI) {{\bf 31}};

\draw[line width=\arrowwidth,shorten >=1cm,blur shadow,->]  ($(PR18)+(-3,0)$) -- ($(PR22)+(2,0)$);

\draw[line width=\arrowwidth,shorten >=0.5cm,blur shadow,->]
($(PR12)+(-3,0)$) -- ($(PR7)+(1,0)$);

\draw[line width=\arrowwidth,shorten >=0.5cm,shorten <=0.5cm,blur shadow,<->]  ($(PR12)+(3,0)$) -- ($(PR24)$);

\coordinate (RLK) at ($(PR14 |- PR13)$);
\draw[line width=\arrowwidth,blur shadow,shorten <=0.5cm,<-]  ($(PR14)$) -- ($(RLK)$);
\draw[line width=\arrowwidth,blur shadow,shorten <=2.1cm,<-]  ($(PR13)$) -- ($(RLK)$);

\draw[black,fill=white] (RLK) circle (\circrad);
\node[text=numcol] at (RLK) {{\bf 27}};

\coordinate (RTMD) at ($(PR14)+(0,-5)$);
\draw[line width=\arrowwidth,blur shadow]  ($(PR14)$) -- ($(RTMD)$);
\draw[line width=\arrowwidth,shorten >=0.5cm,blur shadow,->]  ($(RTMD)$) -- ($(PR17)+(-3.6,1)$);

\draw[black,fill=white] (RTMD) circle (\circrad);
\node[text=numcol] at (RTMD) {{\bf 28}};

\draw[line width=\arrowwidth,shorten >=0.7cm,blur shadow,->]  ($(PR13)+(3,0)$) -- ($(PR12)+(-3,0)$);

\draw[line width=\arrowwidth,shorten >=0.5cm,blur shadow,->]  ($(PR13)+(-6.3,0)$) -- ($(PR19)-(1,0)$);

\draw[line width=\arrowwidth,shorten >=0.1cm,blur shadow,->]
(PR25) -- ($(PR13)+(6.5,0)$);

\coordinate (CI) at (30,15.2);
\draw[black,fill=white] (CI) circle (\circrad);
\node[text=numcol] at (CI) {{\bf 29}};

\draw[line width=\arrowwidth,shorten >=0.65cm,blur shadow,->]  ($(PR19)$) -- ($(PR22)+(-1,0)$);

\draw[line width=\arrowwidth,shorten >=0.85cm,blur shadow,->]  ($(PR15)$) -- ($(PR22)+(-3,0)$);

\draw[line width=\arrowwidth,shorten >=1.0cm,blur shadow,->]  ($(PR12)+(-2,0)$) -- ($(PR17)+(2.4,0)$);

\draw[line width=\arrowwidth,shorten >= 0.4cm,blur shadow,->]  ($(PR12)$) -- ($(PR20)$);

\coordinate (CI) at ($(PR20)+(0,4)$);
\draw[black,fill=white] (CI) circle (\circrad);
\node[text=numcol] at (CI) {{\bf 30}};

\coordinate (DEH) at ($(PR17)+(-9,0)$);
\draw[line width=\arrowwidth,blur shadow,-]  (PR17) -- (DEH) ;
\draw[line width=\arrowwidth,shorten >=0.5cm,blur shadow,->]  (DEH)  -- ($(PR21)+(2,0)$);

\draw[line width=\arrowwidth,shorten >=0.5cm,blur shadow,->]  ($(PR5)-(1,0)$) -- ($(PR21)$);

\def\titletxt{Motivic relations};
\def\n{1}
\coordinate (dim) at (13,3.5);

\coordinate (CC) at (PR\n);
\coordinate (TT) at ($(CC)+0.5*(Z |- dim)+(0,-1)$);
\coordinate   (MA) at ($(CC)+(0,-0.5)$);
\coordinate (LL) at ($(CC)-0.5*(dim)$);
\coordinate (UR) at ($(CC)+0.5*(dim)$);
\coordinate (CI) at ($(CC)-0.5*(Z -| dim)$);

\draw[rounded corners,line width=\linewidthbig,color=bigcolor,fill=algrelcol,blur shadow] (LL) rectangle (UR) {};
\draw[black,fill=white] (CI) circle (\circrad);
\node[text=numcol] at (CI) {{\bf \n}};
\node at (CC) {{\bf \titletxt}};

\def\titletxt{Associator relations};
\def\n{2}
\coordinate (dim) at (13,3.5);
\def\math{$\begin{aligned}
        \Phi_{\text{KZ}}(x,y)
      \end{aligned}$};

\coordinate (CC) at (PR\n);
\coordinate (TT) at ($(CC)+0.5*(Z |- dim)+(0,-1)$);
\coordinate   (MA) at ($(CC)+(0,-0.5)$);
\coordinate (LL) at ($(CC)-0.5*(dim)$);
\coordinate (UR) at ($(CC)+0.5*(dim)$);
\coordinate (CI) at ($(CC)-0.5*(Z -| dim)$);

\draw[rounded corners,line width=\linewidthbig,color=bigcolor,fill=algrelcol,blur shadow] (LL) rectangle (UR) {};
\draw[black,fill=white] (CI) circle (\circrad);
\node[text=numcol] at (CI) {{\bf \n}};
\node at (TT) {{\bf \titletxt}};
\node at ($(MA)-(1.4,0.2)$) {\math};
\def\R{0.7}
\coordinate (PENT) at ($(MA)+(2.5,-0.2)$);
\draw[shift=(PENT),color=pentacol,line width=1pt] (90:\R) 
    \foreach \x in {162,234,...,449} {
      -- (\x:\R)
    }-- cycle (0:\R);

\def\titletxt{Confluence relations};
\def\n{3}
\coordinate (dim) at (14,3.5);
\def\math{$\begin{aligned}
         L(\mathcal{I}_{\text{CF}}) = 0 
      \end{aligned}$};

\coordinate (CC) at (PR\n);
\coordinate (TT) at ($(CC)+0.5*(Z |- dim)+(0,-1)$);
\coordinate   (MA) at ($(CC)+(0,-0.5)$);
\coordinate (LL) at ($(CC)-0.5*(dim)$);
\coordinate (UR) at ($(CC)+0.5*(dim)$);
\coordinate (CI) at ($(CC)-0.5*(Z -| dim)$);

\draw[rounded corners,line width=\linewidthbig,color=bigcolor,fill=algrelcol,blur shadow] (LL) rectangle (UR) {};
\draw[black,fill=white] (CI) circle (\circrad);
\node[text=numcol] at (CI) {{\bf \n}};
\node at (TT) {{\bf \titletxt}};
\node at (MA) {\math};

\def\titletxt{Integral-Series identity};
\def\n{4}
\coordinate (dim) at (13,3.5);
\def\math{$\begin{aligned}
         \zeta(\mu( {\bf k} , {\bf l} ))  = \zeta( {\bf k} \varoast {\bf l}^\star) 
      \end{aligned}$};

\coordinate (CC) at (PR\n);
\coordinate (TT) at ($(CC)+0.5*(Z |- dim)+(0,-1)$);
\coordinate   (MA) at ($(CC)+(0,-0.5)$);
\coordinate (LL) at ($(CC)-0.5*(dim)$);
\coordinate (UR) at ($(CC)+0.5*(dim)$);
\coordinate (CI) at ($(CC)-0.5*(Z -| dim)$);

\draw[rounded corners,line width=\linewidthbig,color=bigcolor,fill=linrelcol,blur shadow] (LL) rectangle (UR) {};
\draw[black,fill=white] (CI) circle (\circrad);
\node[text=numcol] at (CI) {{\bf \n}};
\node at (TT) {{\bf \titletxt}};
\node at (MA) {\math};

\def\titletxt{Extended double shuffle relations};
\def\n{5}
\coordinate (dim) at (16,4.5);
\def\math{$\begin{aligned}
w \in &\HH^1,\, u\in \HH^0\\
	\zeta^\bullet(w \sh u &- w\ast u; T) = 0
      \end{aligned}$};

\coordinate (CC) at (PR\n);
\coordinate (TT) at ($(CC)+0.5*(Z |- dim)+(0,-1)$);
\coordinate   (MA) at ($(CC)+(0,-0.5)$);
\coordinate (LL) at ($(CC)-0.5*(dim)$);
\coordinate (UR) at ($(CC)+0.5*(dim)$);
\coordinate (CI) at ($(CC)-0.5*(Z -| dim)$);

\draw[rounded corners,line width=\linewidthbig,color=bigcolor,fill=algrelcol,blur shadow] (LL) rectangle (UR) {};
\draw[black,fill=white] (CI) circle (\circrad);
\node[text=numcol] at (CI) {{\bf \n}};
\node at (TT) {{\bf \titletxt}};
\node at (MA) {\math};

\def\titletxt{Kawashima's relations};
\def\n{6}
\coordinate (dim) at (13,3.5);
\def\math{$\begin{aligned}
      \sum_{\substack{i+j=m\\i,j \geq 1}} \zeta(\varphi(v) \circledast y^i) \zeta(\varphi(w) \circledast y^j)  = \zeta( \varphi(v \ast w) \circledast y^m)
      \end{aligned}$};

\coordinate (CC) at (PR\n);
\coordinate (TT) at ($(CC)+0.5*(Z |- dim)+(0,-1)$);
\coordinate   (MA) at ($(CC)+(0,-0.5)$);
\coordinate (LL) at ($(CC)-0.5*(dim)$);
\coordinate (UR) at ($(CC)+0.5*(dim)$);
\coordinate (CI) at ($(CC)-0.5*(Z -| dim)$);

\draw[rounded corners,line width=\linewidthbig,color=bigcolor,fill=algrelcol,blur shadow] (LL) rectangle (UR) {};
\draw[black,fill=white] (CI) circle (\circrad);
\node[text=numcol] at (CI) {{\bf \n}};
\node at (CC) {{\bf \titletxt}};

\def\titletxt{Duality};
\def\n{7}
\coordinate (dim) at (9.5,3.5);
\def\math{$\begin{aligned}
         \zeta(w) = \zeta(\tau(w))
      \end{aligned}$};

\coordinate (CC) at (PR\n);
\coordinate (TT) at ($(CC)+0.5*(Z |- dim)+(0,-1)$);
\coordinate   (MA) at ($(CC)+(0,-0.5)$);
\coordinate (LL) at ($(CC)-0.5*(dim)$);
\coordinate (UR) at ($(CC)+0.5*(dim)$);
\coordinate (CI) at ($(CC)-0.5*(Z -| dim)$);

\draw[rounded corners,line width=\linewidth,fill=linrelcol,blur shadow] (LL) rectangle (UR) {};
\draw[black,fill=white] (CI) circle (\circrad);
\node[text=numcol] at (CI) {{\bf \n}};
\node at (TT) {{\bf \titletxt}};
\node at (MA) {\math};

\def\titletxt{Euler's relations};
\def\n{8}
\coordinate (dim) at (12,3.5);
\def\math{$\begin{aligned}
         \zeta(2k) \in \zeta(2)^k \mathbb{Q} = \pi^{2k} \Q 
      \end{aligned}$};

\coordinate (CC) at (PR\n);
\coordinate (TT) at ($(CC)+0.5*(Z |- dim)+(0,-1)$);
\coordinate   (MA) at ($(CC)+(0,-0.5)$);
\coordinate (LL) at ($(CC)-0.5*(dim)$);
\coordinate (UR) at ($(CC)+0.5*(dim)$);
\coordinate (CI) at ($(CC)-0.5*(Z -| dim)$);

\draw[rounded corners,line width=\linewidth,fill=algrelcol,blur shadow] (LL) rectangle (UR) {};
\draw[black,fill=white] (CI) circle (\circrad);
\node[text=numcol] at (CI) {{\bf \n}};
\node at (TT) {{\bf \titletxt}};
\node at (MA) {\math};

\def\titletxt{Finite double shuffle};
\def\n{9}
\coordinate (dim) at (11,3.5);
\def\math{$\begin{aligned}
        	\zeta(w \sh u - w\ast u) = 0
      \end{aligned}$};

\coordinate (CC) at (PR\n);
\coordinate (TT) at ($(CC)+0.5*(Z |- dim)+(0,-1)$);
\coordinate   (MA) at ($(CC)+(0,-0.5)$);
\coordinate (LL) at ($(CC)-0.5*(dim)$);
\coordinate (UR) at ($(CC)+0.5*(dim)$);
\coordinate (CI) at ($(CC)-0.5*(Z -| dim)$);

\draw[rounded corners,line width=\linewidth,fill=linrelcol,blur shadow] (LL) rectangle (UR) {};
\draw[black,fill=white] (CI) circle (\circrad);
\node[text=numcol] at (CI) {{\bf \n}};
\node at (TT) {{\bf \titletxt}};
\node at (MA) {\math};

\node at ($(PR9)+(-11,0)$) {$w,u \in \mathfrak{H}^0$};

\def\titletxt{Shuffle product};
\def\n{10}
\coordinate (dim) at (10,3.5);
\def\math{$\begin{aligned}
         \zeta(w) \zeta(u)= \zeta(w \shuffle u)
      \end{aligned}$};

\coordinate (CC) at (PR\n);
\coordinate (TT) at ($(CC)+0.5*(Z |- dim)+(0,-1)$);
\coordinate   (MA) at ($(CC)+(0,-0.5)$);
\coordinate (LL) at ($(CC)-0.5*(dim)$);
\coordinate (UR) at ($(CC)+0.5*(dim)$);
\coordinate (CI) at ($(CC)-0.5*(Z -| dim)$);

\draw[rounded corners,line width=\linewidth,fill=algrelcol,blur shadow] (LL) rectangle (UR) {};
\draw[black,fill=white] (CI) circle (\circrad);
\node[text=numcol] at (CI) {{\bf \n}};
\node at (TT) {{\bf \titletxt}};
\node at (MA) {\math};

\def\titletxt{Stuffle product};
\def\n{11}
\coordinate (dim) at (10,3.5);
\def\math{$\begin{aligned}
            \zeta(w) \zeta(u)= \zeta(w \ast u)
      \end{aligned}$};

\coordinate (CC) at (PR\n);
\coordinate (TT) at ($(CC)+0.5*(Z |- dim)+(0,-1)$);
\coordinate   (MA) at ($(CC)+(0,-0.5)$);
\coordinate (LL) at ($(CC)-0.5*(dim)$);
\coordinate (UR) at ($(CC)+0.5*(dim)$);
\coordinate (CI) at ($(CC)-0.5*(Z -| dim)$);

\draw[rounded corners,line width=\linewidth,fill=algrelcol,blur shadow] (LL) rectangle (UR) {};
\draw[black,fill=white] (CI) circle (\circrad);
\node[text=numcol] at (CI) {{\bf \n}};
\node at (TT) {{\bf \titletxt}};
\node at (MA) {\math};

\def\titletxt{Ohno's relations};
\def\n{12}
\coordinate (dim) at (10,3.5);
\def\math{$\begin{aligned}
\mathcal{O}^X(\kk) =     \mathcal{O}^X(\kk^\dagger)
      \end{aligned}$};

\coordinate (CC) at (PR\n);
\coordinate (TT) at ($(CC)+0.5*(Z |- dim)+(0,-1)$);
\coordinate   (MA) at ($(CC)+(0,-0.5)$);

\coordinate (LL) at ($(CC)-0.5*(dim)$);
\coordinate (UR) at ($(CC)+0.5*(dim)$);
\coordinate (CI) at ($(CC)-0.5*(Z -| dim)$);

\draw[rounded corners,line width=\linewidth,fill=linrelcol,blur shadow] (LL) rectangle (UR) {};
\draw[black,fill=white] (CI) circle (\circrad);
\node[text=numcol] at (CI) {{\bf \n}};
\node at (TT) {{\bf \titletxt}};
\node at (MA) {\math};

\def\titletxt{Linear part of};
\def\n{13}
\coordinate (dim) at (13,3.5);
\def\math{$\begin{aligned}
         ... = ...
      \end{aligned}$};

\coordinate (CC) at (PR\n);
\coordinate (TT) at ($(CC)+0.5*(Z |- dim)+(0,-1)$);
\coordinate   (MA) at ($(CC)+(0,-0.5)$);
\coordinate (LL) at ($(CC)-0.5*(dim)$);
\coordinate (UR) at ($(CC)+0.5*(dim)$);
\coordinate (CI) at ($(CC)-0.5*(Z -| dim)$);

\draw[rounded corners,line width=\linewidth,fill=linrelcol,blur shadow] (LL) rectangle (UR) {};
\draw[black,fill=white] (CI) circle (\circrad);
\node[text=numcol] at (CI) {{\bf \n}};
\node at (TT) {{\bf \titletxt}};
\node at (MA) {{\bf Kawashima's relations}};

\def\titletxt{Rooted tree maps};
\def\n{14}
\coordinate (dim) at (9.5,3.5);
\def\math{$\begin{aligned}
         \input{TikZ/tree_1_2_2} \in \operatorname{End}(\mathfrak{H}^1)
      \end{aligned}$};

\coordinate (CC) at (PR\n);
\coordinate (TT) at ($(CC)+0.5*(Z |- dim)+(0,-1)$);
\coordinate   (MA) at ($(CC)+(0,-0.5)$);
\coordinate (LL) at ($(CC)-0.5*(dim)$);
\coordinate (UR) at ($(CC)+0.5*(dim)$);
\coordinate (CI) at ($(CC)-0.5*(Z -| dim)$);

\draw[rounded corners,line width=\linewidth,fill=linrelcol,blur shadow] (LL) rectangle (UR) {};
\draw[black,fill=white] (CI) circle (\circrad);
\node[text=numcol] at (CI) {{\bf \n}};
\node at (TT) {{\bf \titletxt}};
\node at (MA) {\math};

\def\titletxt{Ohno-Zagier relations};
\def\n{15}
\coordinate (dim) at (11,4.8);
\def\math{$\begin{aligned}
         \sum_{\substack{\operatorname{wt}({\bf k})=k\\ \operatorname{dep}({\bf k})=r\\\operatorname{ht}({\bf k})=h}} \!\!\!\zeta({\bf k}) \in \Q[\zeta(l) \mid l\leq k]
      \end{aligned}$};

\coordinate (CC) at (PR\n);
\coordinate (TT) at ($(CC)+0.5*(Z |- dim)+(0,-1)$);
\coordinate   (MA) at ($(CC)+(0,-0.7)$);
\coordinate (LL) at ($(CC)-0.5*(dim)$);
\coordinate (UR) at ($(CC)+0.5*(dim)$);
\coordinate (CI) at ($(CC)-0.5*(Z -| dim)$);

\draw[rounded corners,line width=\linewidth,fill=algrelcol,blur shadow] (LL) rectangle (UR) {};
\draw[black,fill=white] (CI) circle (\circrad);
\node[text=numcol] at (CI) {{\bf \n}};
\node at (TT) {{\bf \titletxt}};
\node at (MA) {\small \math};

\def\titletxt{Quasi-derivation relations};
\def\n{16}
\coordinate (dim) at (13.1,3.5);
\def\math{$\begin{aligned}
         \zeta(\partial^{(c)}_n(w)) = 0
      \end{aligned}$};

\coordinate (CC) at (PR\n);
\coordinate (TT) at ($(CC)+0.5*(Z |- dim)+(0,-1)$);
\coordinate   (MA) at ($(CC)+(0,-0.5)$);
\coordinate (LL) at ($(CC)-0.5*(dim)$);
\coordinate (UR) at ($(CC)+0.5*(dim)$);
\coordinate (CI) at ($(CC)-0.5*(Z -| dim)$);

\draw[rounded corners,line width=\linewidth,fill=linrelcol,blur shadow] (LL) rectangle (UR) {};
\draw[black,fill=white] (CI) circle (\circrad);
\node[text=numcol] at (CI) {{\bf \n}};
\node at (TT) {{\bf \titletxt}};
\node at (MA) {\math};

\def\titletxt{Derivation relations};
\def\n{17}
\coordinate (dim) at (11,3.5);
\def\math{$\begin{aligned}
         \zeta(\partial_n(w)) = 0
      \end{aligned}$};

\coordinate (CC) at (PR\n);
\coordinate (TT) at ($(CC)+0.5*(Z |- dim)+(0,-1)$);
\coordinate   (MA) at ($(CC)+(0,-0.5)$);
\coordinate (LL) at ($(CC)-0.5*(dim)$);
\coordinate (UR) at ($(CC)+0.5*(dim)$);
\coordinate (CI) at ($(CC)-0.5*(Z -| dim)$);

\draw[rounded corners,line width=\linewidth,fill=linrelcol,blur shadow] (LL) rectangle (UR) {};
\draw[black,fill=white] (CI) circle (\circrad);
\node[text=numcol] at (CI) {{\bf \n}};
\node at (TT) {{\bf \titletxt}};
\node at (MA) {\math};

\def\titletxt{Restricted sum formula};
\def\n{18}
\coordinate (dim) at (12.5,3);

\coordinate (CC) at (PR\n);
\coordinate (TT) at ($(CC)+0.5*(Z |- dim)+(0,-1)$);
\coordinate   (MA) at ($(CC)+(0,-0.5)$);
\coordinate (LL) at ($(CC)-0.5*(dim)$);
\coordinate (UR) at ($(CC)+0.5*(dim)$);
\coordinate (CI) at ($(CC)-0.5*(Z -| dim)$);

\draw[rounded corners,line width=\linewidth,fill=linrelcol,blur shadow] (LL) rectangle (UR) {};
\draw[black,fill=white] (CI) circle (\circrad);
\node[text=numcol] at (CI) {{\bf \n}};
\node at (CC) {{\bf \titletxt}};

\def\titletxt{Cyclic sum formula};
\def\n{19}
\coordinate (dim) at (11,3);

\coordinate (CC) at (PR\n);
\coordinate (TT) at ($(CC)+0.5*(Z |- dim)+(0,-1)$);
\coordinate   (MA) at ($(CC)+(0,-0.5)$);
\coordinate (LL) at ($(CC)-0.5*(dim)$);
\coordinate (UR) at ($(CC)+0.5*(dim)$);
\coordinate (CI) at ($(CC)-0.5*(Z -| dim)$);

\draw[rounded corners,line width=\linewidth,fill=linrelcol,blur shadow] (LL) rectangle (UR) {};
\draw[black,fill=white] (CI) circle (\circrad);
\node[text=numcol] at (CI) {{\bf \n}};
\node at (CC) {{\bf \titletxt}};

\def\titletxt{Weighted sum formula};
\def\n{20}
\coordinate (dim) at (12,3);

\coordinate (CC) at (PR\n);
\coordinate (TT) at ($(CC)+0.5*(Z |- dim)+(0,-1)$);
\coordinate   (MA) at ($(CC)+(0,-0.5)$);
\coordinate (LL) at ($(CC)-0.5*(dim)$);
\coordinate (UR) at ($(CC)+0.5*(dim)$);
\coordinate (CI) at ($(CC)-0.5*(Z -| dim)$);

\draw[rounded corners,line width=\linewidth,fill=linrelcol,blur shadow] (LL) rectangle (UR) {};
\draw[black,fill=white] (CI) circle (\circrad);
\node[text=numcol] at (CI) {{\bf \n}};
\node at (CC) {{\bf \titletxt}};

\def\titletxt{Hoffman's relations};
\def\n{21}
\coordinate (dim) at (10.5,3.5);
\def\math{$\begin{aligned}
         \zeta(z_1 \shuffle u - z_1 \ast u) = 0
      \end{aligned}$};

\coordinate (CC) at (PR\n);
\coordinate (TT) at ($(CC)+0.5*(Z |- dim)+(0,-1)$);
\coordinate   (MA) at ($(CC)+(0,-0.5)$);
\coordinate (LL) at ($(CC)-0.5*(dim)$);
\coordinate (UR) at ($(CC)+0.5*(dim)$);
\coordinate (CI) at ($(CC)-0.5*(Z -| dim)$);

\draw[rounded corners,line width=\linewidth,fill=linrelcol,blur shadow] (LL) rectangle (UR) {};
\draw[black,fill=white] (CI) circle (\circrad);
\node[text=numcol] at (CI) {{\bf \n}};
\node at (TT) {{\bf \titletxt}};
\node at (MA) {\math};

\def\titletxt{Sum formula};
\def\n{22}
\coordinate (dim) at (15,4.9);
\def\math{$
	  \begin{aligned}
         \sum_{\substack{k_1+\dots+k_r=k\\k_1\geq 2, k_2,\dots,k_{r}\geq 1}} \zeta(k_1,\dots,k_r) = \zeta(k)
      \end{aligned}
$};

\coordinate (CC) at (PR\n);
\coordinate (TT) at ($(CC)+0.5*(Z |- dim)+(0,-1)$);
\coordinate   (MA) at ($(CC)+(0,-0.7)$);
\coordinate (LL) at ($(CC)-0.5*(dim)$);
\coordinate (UR) at ($(CC)+0.5*(dim)$);
\coordinate (CI) at ($(CC)-0.5*(Z -| dim)$);

\draw[rounded corners,line width=\linewidth,fill=linrelcol,blur shadow] (LL) rectangle (UR) {};
\draw[black,fill=white] (CI) circle (\circrad);
\node[text=numcol] at (CI) {{\bf \n}};
\node at (TT) {{\bf \titletxt}};
\node at (MA) {\math};

\def\titletxt{Period polynomial relations};
\def\n{23}
\coordinate (dim) at (17,3.5);
\def\math{$
	\begin{aligned}
168 \zeta(5,7)+150 \zeta(7,5)+28 \zeta(9,3) \equiv 0
	\end{aligned}
	$};

\coordinate (CC) at (PR\n);
\coordinate (TT) at ($(CC)+0.5*(Z |- dim)+(0,-1)$);
\coordinate   (MA) at ($(CC)+(0,-0.7)$);
\coordinate (LL) at ($(CC)-0.5*(dim)$);
\coordinate (UR) at ($(CC)+0.5*(dim)$);
\coordinate (CI) at ($(CC)-0.5*(Z -| dim)$);

\draw[rounded corners,line width=\linewidth,fill=linrelcol,blur shadow] (LL) rectangle (UR) {};
\draw[black,fill=white] (CI) circle (\circrad);
\node[text=numcol] at (CI) {{\bf 32}};
\node at (TT) {{\bf \titletxt}};
\node at (MA) {\math};

\def\titletxt{Double Ohno};
\def\n{24}
\coordinate (dim) at (8,3.5);
\def\math{relations};

\coordinate (CC) at (PR\n);
\coordinate (TT) at ($(CC)+0.5*(Z |- dim)+(0,-1)$);
\coordinate   (MA) at ($(CC)+(0,-0.5)$);

\coordinate (LL) at ($(CC)-0.5*(dim)$);
\coordinate (UR) at ($(CC)+0.5*(dim)$);
\coordinate (CI) at ($(CC)-0.5*(Z -| dim)$);

\draw[rounded corners,line width=\linewidth,fill=linrelcol,blur shadow] (LL) rectangle (UR) {};
\draw[black,fill=white] (CI) circle (\circrad);
\node[text=numcol] at (CI) {{\bf 33}};
\node at (TT) {{\bf \titletxt}};
\node at (MA) {{\bf \math}};

\def\titletxt{Drop1 relations};
\def\n{34}
\coordinate (dim) at (10,3.5);
\def\math{$\zeta(w-\dropone(w))=0$};

\coordinate (CC) at (PR25);
\coordinate (TT) at ($(CC)+0.5*(Z |- dim)+(0,-1)$);
\coordinate (MA) at ($(CC)+(0,-0.5)$);
\coordinate (LL) at ($(CC)-0.5*(dim)$);
\coordinate (UR) at ($(CC)+0.5*(dim)$);
\coordinate (CI) at ($(CC)-0.5*(Z -| dim)$);

\draw[rounded corners,line width=\linewidth,fill=linrelcol,blur shadow] (LL) rectangle (UR) {};
\draw[black,fill=white] (CI) circle (\circrad);
\node[text=numcol] at (CI) {{\bf \n}};
\node at (TT) {{\bf \titletxt}};
\node at (MA) {\math};

\draw[rounded corners,line width=0.5pt,fill=white,blur shadow] (-1.2,-2) rectangle (11.1,8) {};

\draw[rounded corners,line width=\linewidthbig,color=bigcolor,blur shadow,shading = axis,rectangle, left color=linrelcol, right color=algrelcol,shading angle=135, anchor=north] (-0.2,3.8) rectangle (10.2,6.2) {};
\node at (5,5.4) {\small Conjecturally give all};
\node at (5,4.4) {\small relations among MZV};

\draw[line width=\arrowwidth,shorten >=0.3cm,blur shadow,->] (2.5,2)  -- (7.5,2);
\draw[rounded corners,line width=\linewidth,blur shadow,shading = axis,rectangle, left color=linrelcol, right color=algrelcol,shading angle=135, anchor=north] (1.5,1) rectangle (3.5,3) {};
\draw[rounded corners,line width=\linewidth,blur shadow,shading = axis,rectangle, left color=linrelcol, right color=algrelcol,shading angle=135, anchor=north] (6.5,1) rectangle (8.5,3) {};
\node at (2.5,2) {{\bf A}};
\node at (7.5,2) {{\bf B}};

\node at (5,7) {{\bf Legend}};
\node at (5,0.2) [align=center]{{\bf A} implies {\bf B}};
\node at (5,-0.7) [align=center]{\small (after a possible use of $\shuffle$ or \,$\ast$\,)};

\end{tikzpicture}

%% file: TikZ/tree_1_2_2.tex
\begin{tikzpicture}[scale=0.2,baseline={([yshift=-.5ex]current bounding box.center)}]
\def\cz{5}
\def\wi{0.5}

\newcommand{\ci}[1]{	
	\fill[black] (#1) circle (\cz pt);
	\draw (#1) circle (\cz pt);
}

\coordinate (R) at (0,0);

\coordinate (r1) at (\wi,-1);
\coordinate (l1) at (-\wi,-1);

\coordinate (l11) at (-\wi-\wi,-2);
\coordinate (l12) at (-\wi+\wi,-2);

\draw (R) to (r1);
\draw (R) to (l1);

\draw (l1) to (l11);
\draw (l1) to (l12);

\ci{R}
\ci{r1}
\ci{l1}
\ci{l11}
\ci{l12}
\end{tikzpicture}

%% file: TikZ/tree2_1_2_2.tex
\begin{tikzpicture}[scale=0.25,baseline={([yshift=-.5ex]current bounding box.center)}]
\def\cz{5}
\def\wi{0.5}

\newcommand{\ci}[1]{	
	\fill[black] (#1) circle (\cz pt);
	\draw (#1) circle (\cz pt);
}

\coordinate (R) at (0,0);

\coordinate (r1) at (\wi,-1);
\coordinate (l1) at (-\wi,-1);

\coordinate (r11) at (\wi-\wi,-2);
\coordinate (r12) at (\wi+\wi,-2);

\draw (R) to (r1);
\draw (R) to (l1);

\draw (r1) to (r11);
\draw (r1) to (r12);

\ci{R}
\ci{r1}
\ci{l1}
\ci{r11}
\ci{r12}
\end{tikzpicture}

%% file: TikZ/tree_1_2.tex
\begin{tikzpicture}[scale=0.3,baseline={([yshift=-.5ex]current bounding box.center)}]
\def\cz{5}
\def\wi{0.5}

\newcommand{\ci}[1]{	
	\fill[black] (#1) circle (\cz pt);
	\draw (#1) circle (\cz pt);
}

\coordinate (R) at (0,0);
\coordinate (r1) at (\wi,-1);
\coordinate (l1) at (-\wi,-1);

\draw (R) to (l1);
\draw (R) to (r1);

\ci{R};
\ci{l1};
\ci{r1};
\end{tikzpicture}

%% file: TikZ/tree_1.tex
\begin{tikzpicture}[scale=0.3,baseline={([yshift=-.5ex]current bounding box.center)}]
\def\cz{5}
\def\wi{0.5}

\newcommand{\ci}[1]{	
	\fill[black] (#1) circle (\cz pt);
	\draw (#1) circle (\cz pt);
}

\coordinate (R) at (0,0);
\ci{R}
\end{tikzpicture}

%% file: TikZ/tree_2.tex
\begin{tikzpicture}[scale=0.3,baseline={([yshift=-.5ex]current bounding box.center)}]
\def\cz{5}
\def\wi{0.5}

\newcommand{\ci}[1]{	
	\fill[black] (#1) circle (\cz pt);
	\draw (#1) circle (\cz pt);
}

\coordinate (R) at (0,0);
\coordinate (r1) at (0,-1);

\draw (R) to (r1);

\ci{R};
\ci{r1};

\end{tikzpicture}

%% file: TikZ/tree_3.tex
\begin{tikzpicture}[scale=0.3,baseline={([yshift=-.5ex]current bounding box.center)}]
\def\cz{5}
\def\wi{0.5}

\newcommand{\ci}[1]{	
	\fill[black] (#1) circle (\cz pt);
	\draw (#1) circle (\cz pt);
}

\coordinate (R) at (0,0);
\coordinate (r1) at (0,-1);
\coordinate (r2) at (0,-2);

\draw (R) to (r1);
\draw (r1) to (r2);
\ci{R};
\ci{r1};
\ci{r2};
\end{tikzpicture}

%% file: TikZ/tree_1_3.tex
\begin{tikzpicture}[scale=0.3,baseline={([yshift=-.5ex]current bounding box.center)}]
\def\cz{5}
\def\wi{0.5}

\newcommand{\ci}[1]{	
	\fill[black] (#1) circle (\cz pt);
	\draw (#1) circle (\cz pt);
}

\coordinate (R) at (0,0);
\coordinate (a1) at (-\wi,-1);
\coordinate (a2) at (0,-1);
\coordinate (a3) at (\wi,-1);

\draw (R) to (a1);
\draw (R) to (a2);
\draw (R) to (a3);

\ci{R}
\ci{a1}
\ci{a2}
\ci{a3}

\end{tikzpicture}

%% file: TikZ/tree_1_2_1.tex
\begin{tikzpicture}[scale=0.3,baseline={([yshift=-.5ex]current bounding box.center)}]
\def\cz{5}
\def\wi{0.5}

\newcommand{\ci}[1]{	
	\fill[black] (#1) circle (\cz pt);
	\draw (#1) circle (\cz pt);
}

\coordinate (R) at (0,0);
\coordinate (l1) at (-\wi,-1);
\coordinate (r1) at (\wi,-1);
\coordinate (l11) at (-\wi,-2);

\draw (R) to (l1);
\draw (R) to (r1);
\draw (l1) to (l11);

\ci{R}
\ci{l1}
\ci{r1}
\ci{l11}

\end{tikzpicture}

%% file: TikZ/tree_1_1_2.tex
\begin{tikzpicture}[scale=0.3,baseline={([yshift=-.5ex]current bounding box.center)}]
\def\cz{5}
\def\wi{0.5}

\newcommand{\ci}[1]{	
	\fill[black] (#1) circle (\cz pt);
	\draw (#1) circle (\cz pt);
}

\coordinate (R) at (0,0);

\coordinate (l1) at (0,-1);
\coordinate (l11) at (\wi,-2);
\coordinate (l12) at (-\wi,-2);

\draw (R) to (l1);

\draw (l1) to (l11);
\draw (l1) to (l12);

\ci{R}
\ci{l1}
\ci{l11}
\ci{l12}
\end{tikzpicture}

%% file: TikZ/tree_4.tex
\begin{tikzpicture}[scale=0.3,baseline={([yshift=-.5ex]current bounding box.center)}]
\def\cz{5}
\def\wi{0.5}

\newcommand{\ci}[1]{	
	\fill[black] (#1) circle (\cz pt);
	\draw (#1) circle (\cz pt);
}

\coordinate (R) at (0,0);
\coordinate (r1) at (0,-1);
\coordinate (r2) at (0,-2);
\coordinate (r3) at (0,-3);

\draw (R) to (r1);
\draw (r1) to (r2);
\draw (r2) to (r3);
\ci{R};
\ci{r1};
\ci{r2};
\ci{r3};
\end{tikzpicture}

%% file: TikZ/tree_ladder.tex
\begin{tikzpicture}[scale=0.3,baseline={([yshift=-.5ex]current bounding box.center)}]
\def\cz{5}
\def\wi{0.5}

\newcommand{\ci}[1]{	
	\fill[black] (#1) circle (\cz pt);
	\draw (#1) circle (\cz pt);
}

\coordinate (R) at (0,0);
\coordinate (r1) at (0,-1);
\coordinate (r11) at (0,-1.4);

\coordinate (r21) at (0,-2.6);
\coordinate (r2) at (0,-3);
\coordinate (r3) at (0,-4);

\draw (R) to (r1);

\draw (r1) to (r11);
\draw[densely dotted] (r11) to (r21);
\draw (r21) to (r2);
\draw (r2) to (r3);

\ci{R};
\ci{r1};
\ci{r2};
\ci{r3};

\draw[decoration={brace,mirror,raise=5pt},decorate]
  (r3) -- node[right=6pt] {$\scriptstyle m$} (R);

\end{tikzpicture}

%% file: chap_FormalDoubleZetaSpace.tex
\chapter[Double zeta values and modular forms]{Double zeta values \\and modular forms} \label{sec:mdandmzv}

We now return to the connection between multiple zeta values and modular forms, which we already mentioned in Chapter \ref{sec:overview}. In this chapter, we want to make this connection precise in depth two, i.e. for double zeta values. For this, we introduce the formal double zeta space, which is given by formal symbols satisfying the same double shuffle relations as double zeta values. We will then study its realizations and see how period polynomials and cusp forms show up. The depth-two results in this chapter are contained in or inspired by the beautiful work \cite{GKZ}. In the last section, we will see that, with a shift in depth, the same formal double zeta space and the same period polynomials also appear for finite multiple zeta values of depth three and four.

\section{The formal double zeta space}
The idea of the formal double zeta space is to consider formal symbols satisfying the same depth-two double shuffle relations as double zeta values. 
We therefore start by recalling these relations before defining the formal double zeta space.  

\subsection{Double zeta values}
Recall that for $k_1,k_2\geq 2$ we have the finite double shuffle relations
\begin{align*}
\zeta(k_1) \zeta(k_2) &= \zeta(k_1,k_2) + \zeta(k_2,k_1) + \zeta(k_1+k_2)\\
&=\sum_{j=2}^{k_1+k_2-1} \left( \binom{j-1}{k_1-1} + \binom{j-1}{k_2-1} \right) \zeta(j, k_1+k_2-j)\,.
\end{align*}
Using the stuffle and shuffle regularized multiple zeta values, we have for all {\color{red}$k_1,k_2 \geq 1$}
\begin{align}\begin{split}
\label{eq:shufflereg1}
\zeta^\ast(k_1; T) \zeta^\ast(k_2; T) &= \zeta^\ast(k_1,k_2; T) + \zeta^\ast(k_2,k_1;T) + \zeta^\ast(k_1+k_2;T)\\
&=\sum_{{\color{red}j=1}}^{k_1+k_2-1} \left( \binom{j-1}{k_1-1} + \binom{j-1}{k_2-1} \right) \zeta^\sh(j, k_1+k_2-j;T)\,.
\end{split}
\end{align}
The comparison map $\rho$ (Theorem \ref{thm:comparezastzsh}), which gives $\zeta^\sh(\kk;T) = \rho(\zeta^\ast(\kk;T))$, satisfies $\rho(1)=1, \rho(T)=T$ and $\rho(T^2) = T^2 + \zeta(2)$. Therefore the  $\zeta^\sh(k_1,k_2;T)$ and  $\zeta^\ast(k_1,k_2;T)$ just differ in the case $k_1=k_2=1$ and we have 
\begin{align}\label{eq:diffsh11andst11}
\zeta^\sh(1,1;T) = \zeta^\ast(1,1;T) + \frac{1}{2} \zeta(2) \,,
\end{align}
and $\zeta^\sh(1,1;T) =  \frac{1}{2} T^2$, $\zeta^\ast(1,1;T)  = \frac{1}{2} T^2 - \frac{1}{2}\zeta(2)$.
Now define for $\bullet \in \{\sh, \ast \}$ their generating series
\begin{align*}
\gt^\bullet(X) = \sum_{k\geq 1} \zeta^\bullet(k;T) X^{k-1}\,, \qquad \gt^\bullet(X,Y) = \sum_{k_1, k_2\geq 1} \zeta^\bullet(k_1, k_2;T) X^{k_1-1} Y^{k_2-1}\,.
\end{align*}

Using $\frac{X^{k-1}- Y^{k-1}}{X-Y} = \sum_{k_1+k_2=k} X^{k_1-1} Y^{k_2-1}$ we see that \eqref{eq:shufflereg1} together with \eqref{eq:diffsh11andst11} can therefore be written as 
\begin{align}\begin{split}\label{eq:gensermzvdep2}
\gt^\bullet(X)\gt^\bullet(Y) &=\gt^\bullet(X,Y)+\gt^\bullet(Y,X)+\frac{\gt^\bullet(X)-\gt^\bullet(Y)}{X-Y}   - \delta_{\bullet, \sh} \zeta(2)\\
&=  \gt^\bullet(X+Y,Y) + \gt^\bullet(X+Y,X) + \delta_{\bullet, \ast} \zeta(2)\,,
\end{split}
\end{align}
where $\delta$ denotes the Kronecker delta.

\subsection{The formal double zeta space}
We will now define the formal double zeta space. It is spanned by formal symbols $ Z_k , Z_{k_1,k_2}$ and $P_{k_1,k_2} $ for $k,k_1,k_2\geq 1$, which satisfy similar relations to the regularized versions of $\zeta(k), \zeta(k_1,k_2)$ and $\zeta(k_1)\zeta(k_2)$. The only difference is that we ignore the correction term in the $k=2$ case. This causes no problems and makes things a little bit cleaner. 

\begin{definition}    We define for $k\geq 1$ the \emph{formal double zeta space} in weight $k$ as 
    \begin{align*}
    \dz_k = \quotient{\big \langle Z_k , Z_{k_1,k_2}, P_{k_1,k_2} \mid k_1 + k_2 = k, k_1,k_2 \geq 1 \big \rangle_\Q}{\eqref{eq:dzrel} }        
    \end{align*}
    where we divide out the following relations for $k_1,k_2 \geq 1$
    \begin{align}\begin{split}
    \label{eq:dzrel}
    P_{k_1,k_2} &= Z_{k_1, k_2} +  Z_{k_2, k_1} + Z_{k_1+k_2}
    = \sum_{j=1}^{k_1+k_2-1} \left( \binom{j-1}{k_1-1} + \binom{j-1}{k_2-1} \right) Z_{j, k_1+k_2-j}\,.
    \end{split}
    \end{align}
\end{definition}
\begin{remark}
    Notice that by definition the $P_{k_1,k_2}$ can always be expressed in terms of the $Z$'s and it would therefore be equivalent to define the space $\dz_k$ by the span of elements $Z_{k_1,k_2}$ and $Z_k$ modulo the relations
    \begin{align}\begin{split}
    \label{eq:dzrel2}
    Z_{k_1, k_2} +  Z_{k_2, k_1} + Z_{k_1+k_2} &= \sum_{j=1}^{k_1+k_2-1} \left( \binom{j-1}{k_1-1} + \binom{j-1}{k_2-1} \right) Z_{j, k_1+k_2-j}\,.
    \end{split}
    \end{align}
    But it is convenient to also work with the $P_{k_1,k_2}$, since they correspond to something like the product in most of the realizations (Definition \ref{def:dzrealization}) later. 
\end{remark}

For small weights the formal double zeta space is given by the following relations and basis elements. Since $P_{k_1,k_2} = P_{k_2,k_1}$ is symmetric, we will always just consider the case $k_1\leq k_2$.
\begin{figure}[ht!]
    \centering
    \renewcommand{\arraystretch}{1.8}
    \begin{tabular}{c|c|c|c}
        $k$ & Relations in $\dz_k $ &  Basis of $\dz_k$  &  $\dim_\Q \dz_k$  \\[3pt]  \hline 
        $1$&  - & $Z_1$ & $1$ \\[3pt] \hline
        $2$&  {$\!\begin{aligned}    Z_2 &= 0,\quad  P_{1,1} = 2 Z_{1,1} \end{aligned}$}  & $Z_{1,1}$  &  $1$ \\[3pt] \hline
        $3$&  {$\!\begin{aligned}    Z_{2,1} &= Z_3,\quad  Z_{1,2} = P_{1,2}-2Z_{3} \end{aligned}$}  & $Z_{3}, P_{1,2}$  &  $2$ \\[3pt] \hline
        $4$&  {$\!\begin{aligned}    Z_{4} &= 4Z_{3,1},\quad  Z_{2,2} = 3Z_{3,1},\\ P_{1,3} &= Z_{1,3}+ 5Z_{3,1},\quad  P_{2,2}=10 Z_{3,1} \end{aligned}$}  & $Z_{1,3}, Z_{3,1}$  &  $2$ \\[3pt] \hline
        $5$&  {$\!\begin{aligned}    Z_{4,1} &= 2 Z_5-P_{2,3},\quad Z_{3,2} =  -\frac{11}{2}Z_5 + 3 P_{2,3},\\ Z_{2,3} &=\frac{9}{2}Z_5-2 P_{2,3} ,\quad
            Z_{1,4} =  -3 Z_5 + P_{1,4} + P_{2,3}\end{aligned}$}  & $Z_{5}, P_{1,4}, P_{2,3}$  &  $3$ \\ \hline         
        $6$& $Z_6 = 4 Z_{3,3} + 4 Z_{5,1}\,,\quad \dots$ & $Z_{1,5}, Z_{3,3}, Z_{5,1}$  &  $3$ \\ \hline         
    \end{tabular}
    \renewcommand{\arraystretch}{1}
    \caption{Relations and bases for the formal double zeta space in small weights.}
    \label{fig:dktable}
\end{figure}

Observe in Figure \ref{fig:dktable} that for even weight $k$ we seem to have $\dim_{\Q} \dz_k = \frac{k}{2}$ with a basis given by $Z_{\text{odd},\text{odd}}$ and for odd weight $k$ it seems that $\dim_{\Q} \dz_k = \frac{k+1}{2}$ with a basis given by $Z_k$ and $P_{k_1,k_2}$. These observations are indeed correct and we will prove them in this section (Theorem \ref{thm:dzoddbasis} \& \ref{thm:dzevenbasis}). Notice that lower bounds of $\dim_{\Q} \dz_k$, which coincide with the observed dimension, already follow from the definition. Since \eqref{eq:dzrel2} is symmetric in $k_1$ and $k_2$ we obtain for $k$ even $\frac{k}{2}$ relations among the $k$ generators $Z_k, Z_{1,k-1},\dots,Z_{k-1,1}$ and therefore we have for even $k$
\begin{align}
\dim_{\Q} \dz_k \geq  \frac{k}{2}\,.\qquad (k \text{ even})
\end{align}
For $k$ odd we have $\frac{k-1}{2}$ relations and therefore 
\begin{align}
\dim_{\Q} \dz_k \geq  \frac{k+1}{2}\,.\qquad (k \text{ odd})
\end{align}

It is convenient to consider  generating series when working
with the formal double zeta space and therefore we define the following elements in $\dz_k[X,Y]$:
\begin{align*}
\gz_k(X,Y) &= \sum_{\substack{k_1+k_2=k\\k_1, k_2 \geq 1}} Z_{k_1,k_2} X^{k_1-1} Y^{k_2-1}\,,\\
\gp_k(X,Y) &=  \sum_{\substack{k_1+k_2=k\\k_1, k_2 \geq 1}}P_{k_1,k_2} X^{k_1-1} Y^{k_2-1}\,,\\
\gr_k(X,Y) &= Z_k \frac{X^{k-1} - Y^{k-1}}{X-Y}\,.
\end{align*}
With this the double shuffle relations \eqref{eq:dzrel} can be written as 
\begin{align}\begin{split}\label{eq:dzgenrel} 
\gp_k(X,Y)  &=  \gz_k(X,Y) + \gz_k(Y,X)  + \gr_k(X,Y)
= \gz_k(X+Y,Y) + \gz_k(X+Y,X) \,.  
\end{split}
\end{align}

We will not only be interested in relations in the space $\dz_k$, but also in ``real'' mathematical objects which satisfy these relations. Therefore we first introduce the following notation.

\begin{definition} \label{def:dzrealization} Let $A$ be a $\Q$-vector space. We define the \emph{$A$-valued points} $\dz_k(A)$ for $\dz_k$ by
    \begin{align*}
    \dz_k(A) = \Hom_\Q(\dz_k,A) = \left\{  (Z_k , Z_{k_1,k_2}) \in  A^{k} \mid \text{ satisfying } \eqref{eq:dzrel2}  \right\}\,.
    \end{align*}
    An element in $\dz_k(A)$, i.e. one particular choice of $Z_k, Z_{k_1,k_2} \in A$ for all $k_1+k_2=k$ which satisfy \eqref{eq:dzrel2}, will be called a \emph{realization} of the double zeta space in $A$.
\end{definition}
By comparing  \eqref{eq:gensermzvdep2} and \eqref{eq:dzgenrel} we see that one realization is given by the shuffle regularized multiple zeta values: For $A=\R[T]$ a realization of $\dz_k$ with $k\geq 1$ is given for $k_1,k_2 \geq 1 $ and $k_1+k_2=k$ by 
\begin{align*}
Z_{k} &\longmapsto \begin{cases}
\zeta^\sh(k; T) & k\neq 2\\
0 & k=2
\end{cases}\,,\\
Z_{k_1,k_2} &\longmapsto \zeta^\sh(k_1,k_2; T)\,,\\
P_{k_1,k_2} &\longmapsto \zeta^\sh(k_1; T) \zeta^\sh(k_2; T)\,.
\end{align*}
We will refer to this realization as the \emph{multiple zeta realization}.
Later we will also introduce realizations in the cases $A=\Q$, $A=\Q[[q]]$ and $A=\mathcal{O}(\Ha)$. But before doing so we will prove some results in $\dz_k$. Using the description in terms of generating series we obtain the following theorem. 
\begin{theorem}\label{thm:dzsumformula} \begin{enumerate}[\textup{(}i\textup{)}] 
        \item For all $k\geq 2$ we have 
        \begin{align*}
        \sum_{j=2}^{k-1} Z_{j,k-j} = Z_k\,.
        \end{align*}
        \item For $k \geq 2$  even, we have 
        \begin{align*}
        \sum_{\substack{j=2\\ j \text{ even}}}^{k-2} Z_{j,k-j} = \frac{3}{4} Z_k\,, \qquad     \qquad    \sum_{\substack{j=2\\ j \text{ odd}}}^{k-1} Z_{j,k-j} =  \frac{1}{4} Z_k\,.
        \end{align*}
    \end{enumerate}
\end{theorem}
\begin{proof} By \eqref{eq:dzgenrel} we have 
    \begin{align*}
    D(X,Y) := \gz_k(X+Y,Y) + \gz_k(X+Y,X) -  \gz_k(X,Y) - \gz_k(Y,X)  - \gr_k(X,Y) = 0\,.
    \end{align*}
    The first statement now follows by taking the case $(X,Y)=(1,0)$, since
    \begin{align*}
    0 = D(1,0) =  \gz_k(1,0) + \gz_k(1,1) - \gz_k(1,0) - \gz_k(0,1) - Z_k =     \sum_{j=1}^{k-1} Z_{j,k-j} - Z_{1,k-1} - Z_k\,.
    \end{align*}
    For the second statement first consider for even $k$
    \begin{align*}
    0 = D(1,-1) = \gz_k(0,-1) + \gz_k(0,1) -  \gz_k(1,-1) - \gz_k(-1,1)  -  Z_k = 2 \sum_{j=2}^{k-1} (-1)^j Z_{j,k-j} - Z_k\,.
    \end{align*}
    Taking $D(1,0) \pm \frac{1}{2} D(1,-1)$ we therefore obtain
    \begin{align*}
    0 &= D(1,0) + \frac{1}{2}D(1,-1) = 2    \sum_{\substack{j=2\\ j \text{ even}}}^{k-2} Z_{j,k-j} - \frac{3}{2} Z_k\\
    0 &= D(1,0) - \frac{1}{2}D(1,-1) = 2    \sum_{\substack{j=2\\ j \text{ odd}}}^{k-1} Z_{j,k-j} - \frac{1}{2} Z_k\,,    
    \end{align*} 
    from which the second statement follows after dividing by $2$.
\end{proof}

The polynomials $\gz_k$, $\gp_k$ and $\gr_k$ are all elements in $\dz_k \otimes_\Q V_k$, where $V_k\subset \Q[X,Y]$ denotes the space of all homogeneous polynomials of degree $k-2$. On $V_k$ we define a right-action of $\GL_2(\Z)$ for a $\gamma = \sabcd \in \GL_2(\Z)$ and $F\in V_k$ by 
\begin{align*}
(F{\mid\! \gamma})(X,Y) = F(aX + bY, cX + dY)\,.
\end{align*}
This action can then be extended linearly to an action of the group ring $\Z[\GL_2(\Z)]$ on $\dz_k \otimes_\Q V_k$.
The following elements in $\GL_2(\Z)$ will be of importance when working with the above group action. 
\begin{align*}
&\sigma= \pmatr{-1 & 0\\0 & -1}\,,\qquad
\epsilon = \pmatr{0 & 1 \\ 1 & 0}\,,\qquad \delta=\pmatr{-1 & 0\\0 & 1} \,,\\
&T= \pmatr{1 & 1\\ 0 & 1}\,,\qquad S=\pmatr{0 & -1\\1 & 0}\,,\qquad U= \pmatr{1 & -1\\1 & 0}\,.
\end{align*}

With this we can rewrite \eqref{eq:dzgenrel} even simpler as 
\begin{align}\begin{split}
\label{eq:dzactrel}
\gp_k   &=  \gz_k \mid\! (1 + \epsilon) + \gr_k
= \gz_k \mid\! T(1 +\epsilon)\,.
\end{split}
\end{align}

\begin{lemma} \label{lem:dzparitylemma} For $k\geq 1$ and $A= \epsilon U\epsilon $ we have 
    \begin{align*}
    \gz_k \mid (1 - \sigma)  &= \gp_k \mid\! (1-\delta)(1 + A - SA^2) - \gr_k \mid\! (1+A+A^2)\,.
    \end{align*}
\end{lemma}
\begin{proof}
    First notice that $A  = \epsilon U\epsilon = T \epsilon T^{-1} \epsilon =  \left( \begin{smallmatrix}
0 & 1 \\
-1 & 1
\end{smallmatrix} \right) $ and that we have $A^3 = \sigma$. By \eqref{eq:dzactrel} we get 
    \begin{align*}
    \gz_k \mid\! \epsilon &= -    \gz_k  +  \gp_k -  \gr_k \\
    \gz_k \mid\! T \epsilon T^{-1} &=  -    \gz_k  + \gp_k \mid\! T^{-1}  \,.
    \end{align*}
    and therefore 
    \begin{align*}
    \gz_k \mid\! A &= \gz_k \mid\!  (T \epsilon T^{-1}) \epsilon = \left( - \gz_k  +  \gp_k \mid\! T^{-1}  \right) \mid\! \epsilon = \gz_k + \underbrace{\gp_k \mid\! (T^{-1}\epsilon - 1) + \gr_k }_{=: \mathfrak{K}}\,.    
    \end{align*}
    Iterating this identity two more times we get
    \begin{align*}
    \gz_k \mid\! A^3 &= \gz_k + \mathfrak{K} \mid\! (1+A+A^2)\,.
    \end{align*}
    By direct calculation one can check that the action of  $(T^{-1}\epsilon - 1)(1+A+A^2)$ and $-(1-\delta)(1 + A - SA^2)$ is the same on the symmetric polynomial $\gp_k$.
\end{proof}
Since $(\gz_k \mid \sigma)(X,Y) =  \gz_k(-X,-Y) = (-1)^k\gz_k(X,Y)$ we have 
\begin{align*}
\gz_k \mid (1 - \sigma) = \begin{cases}
2 \, \gz_k &, k \text{ odd} \\
0 &, k \text{ even} \,.
\end{cases}
\end{align*}
Also notice that $\gp_k \mid\! (1-\delta)$ is the generating series of the $P_{j,k-j}$ with $j$ even, for which we write
\begin{align*}
\gp^{\text{ev}}_k(X,Y) := \frac{1}{2}\left( \gp_k \mid\! (1-\delta)\right)(X,Y)  =  \sum_{\substack{j=2\\ j \text{ even}} }^{k-1 }P_{j,k-j} X^{j-1} Y^{k-j-1} \,.
\end{align*}

\begin{theorem}[Parity]\label{thm:dzparity}
    For odd $k\geq 3$, every $Z_{k_1,k_2}$ with $k_1,k_2\geq 1$ and $k_1+k_2=k$ can be written as a linear combination of $P_{\text{ev},\text{od}}$ and $Z_k$. 
    More precisely we have 
    \begin{align*}
    Z_{k_1,k_2}
    &=(-1)^{k_2} \sum_{\substack{j=2\\j \text{ even}}}^{k-1} \left( \binom{k-j-1}{k_1-1} + \binom{k-j-1}{k_2-1} + (-1)^{k_2}\delta_{j,k_1} \right)  P_{j,k-j} \\
    &\quad + \frac{1}{2}\left(  (-1)^{k_1} \binom{k_1+k_2}{k_2}  -1\right) Z_k\,.
    \end{align*}
\end{theorem}
\begin{proof}
    This follows directly from Lemma \ref{lem:dzparitylemma} by considering the coefficient of $X^{k_1-1} Y^{k_2-1}$ in
    \begin{align}\label{eq:parityformula}
    \gz_k(X,Y) &=  \left(\gp_k^{\text{ev}}\mid\! (1 + A - SA^2)\right)(X,Y)  - \frac{1}{2} \left( \gr_k \mid\! (1+A+A^2) \right)(X,Y)\,,
    \end{align}
    and checking that 
    \begin{align*}
    \sum_{\substack{j=1}}^{k-1}\left( \binom{k-j-1}{k_1-1} + \binom{k-j-1}{k_2-1} \right) = \binom{k_1+k_2}{k_1}\,.
    \end{align*}
\end{proof}

\begin{ex} As a consequence of Theorem \ref{thm:dzparity} we get the following relations
    \begin{align*}
    Z_{1,2} &= P_{2, 1} - 2 Z_3\,,    \qquad &&Z_{2,3} = -2 P_{2, 3} + \frac{9}{2} Z_5 \,,\\
    Z_{2,1} &= Z_3\,,         \qquad &&Z_{3,2} = 3 P_{2, 3} - \frac{11}{2} Z_5 \,,\\
    Z_{1,4} &= P_{2, 3} + P_{4, 1} - 3 Z_5 \,,\qquad     &&Z_{4,1} = -P_{2, 3} + 2 Z_5 \,.\\
    \end{align*}
\end{ex}
Using Theorem \ref{thm:dzparity} we can now also prove the dimension formula for $\dz_k$ in the odd weight $k$ case.
\begin{theorem} \label{thm:dzoddbasis} For odd $k\geq 1$ we have $\dim_{\Q} \dz_k =  \frac{k+1}{2}$. For $k=1$ we have $\dz_1 = \Q Z_1$, and for odd $k\geq 3$ the sets
    \begin{align*}
    B_1=\left\{ Z_k, P_{2,k-2}, P_{4,k-4},\dots,P_{k-1,1} \right\},\qquad B_2 =\left\{ Z_k, Z_{1,k-1}, Z_{3,k-3},\dots,Z_{k-2,2} \right\},
    \end{align*}
    are both bases of $\dz_k$.
\end{theorem}
\begin{proof} We already saw that $\dim_{\Q} \dz_k \geq   \frac{k+1}{2}$ since there are just $\frac{k-1}{2}$ different relations among the $k$ generators of $\dz_k$. To prove equality it suffices therefore to show that any element can be expressed in terms of elements in $B_1$ or $B_2$. By Theorem \ref{thm:dzparity} we get that any element in $\dz_k$ can be expressed as linear combinations of elements in $B_1$.  
    For the second basis we just need to show, by symmetry, that any $P_{j,k-j}$ with $j$ even can be expressed by $Z_k$ and $Z_{k_1,k_2}$ with $k_1$ odd. The equation in Theorem \ref{thm:dzparity} modulo $Z_k$ for odd $k_1=1,3,\dots,k-2$ reads
    \begin{align*}
    Z_{k_1,k-k_1} \equiv \sum_{\substack{j=2\\j \text{ even}}}^{k-1} \left( \binom{k-j-1}{k_1-1} + \binom{k-j-1}{k-k_1-1}\right)  P_{j,k-j}   \mod \Q Z_k\,.
    \end{align*}
    We therefore need to show that the matrix $\left( \binom{k-2j-1}{2i-2} + \binom{k-2j-1}{k-2i} \right)_{1\leq i,j \leq \frac{k-1}{2}}$ is invertible. The binomial coefficient $\binom{m}{n}= \frac{m}{n} \binom{m-1}{n-1}$ is even, when $m$ is even and $n$ is odd. Therefore modulo $2$ the factor $\binom{k-2j-1}{k-2i} $ vanishes. Reversing the order of the columns, i.e. replacing the column index $j$ by $\frac{k+1}{2}-j$, the remaining entries become $\binom{2j-2}{2i-2}$ and we get a triangular matrix with $1$ on the diagonal, from which we deduce that the matrix is invertible and therefore we can express $P_{j,k-j}$ in terms of  $Z_{\text{odd},\text{even}}$ and $Z_k$.
\end{proof}

Period polynomials also give rise to relations among formal double zeta values of odd weight. Explicit formulas and the linear independence of the resulting relations are given in \cite{Ma}.

\subsection{The space $\dz_k$ in even weight}

We will now present consequences of Lemma \ref{lem:dzparitylemma} for the even weight case. The rest of this section will be devoted to the even weight case and the connection to modular forms. For even $k$, Lemma \ref{lem:dzparitylemma} implies relations among $P_{\text{ev},\text{ev}}$ and $Z_k$. 
\begin{theorem}[Relations among $P_{\text{ev},\text{ev}}$ and $Z_\text{ev}$]\label{thm:relpevevk} For all $k_1,k_2\geq 1$ with $k=k_1+k_2$ even we have 
    \begin{align}\label{eq:relpevevz}
    \frac{1}{2}\left(  \binom{k_1+k_2}{k_2} - (-1)^{k_1}\right) Z_k =  \sum_{\substack{j=2\\j \text{ even}}}^{k-2} \left( \binom{k-j-1}{k_1-1} + \binom{k-j-1}{k_2-1} - \delta_{j,k_1} \right)  P_{j,k-j}    \,.
    \end{align}
\end{theorem}
\begin{proof} This also follows from Lemma \ref{lem:dzparitylemma}, considering the coefficient of $X^{k_1-1} Y^{k_2-1}$ in
    \begin{align}\label{eq:pevevzk}
    \left(\gp_k^{\text{ev}}\mid\! (1 + A - SA^2)\right)(X,Y)  &= \frac{1}{2} \left( \gr_k \mid\! (1+A+A^2) \right)(X,Y)\,.
    \end{align}
\end{proof}
\begin{ex}  As a consequence of Theorem \ref{thm:relpevevk} we get the following relations by considering the coefficients of $X^5 Y$ and $X^4 Y^2$ in \eqref{eq:pevevzk}:
    \begin{align*}
    6 P_{2, 6} + 3 P_{4, 4} = \frac{27}{2} Z_8\,,\qquad 
    15 P_{2, 6} + 3 P_{4, 4} = \frac{57}{2} Z_8\,.
    \end{align*}
    Combining these two relations we obtain
    \begin{align*}
    P_{4,4} = \frac{7}{6} Z_8\,.
    \end{align*}
    Using the multiple zeta realization, this gives another proof of $\zeta(4)^2 = \frac{7}{6}\zeta(8)$.
\end{ex}

\begin{corollary} \label{cor:prodzk} For even $k$ we have
    \begin{align*}
    \sum_{\substack{j=2\\j\text{ even}}}^{k-2} P_{j, k-j} = \frac{k+1}{2} Z_k\,.
    \end{align*}
\end{corollary}
\begin{proof} This is the $(k_1,k_2) = (1,k-1)$ case in Theorem \ref{thm:relpevevk} but can also be obtained from the even sum formula in Theorem \ref{thm:dzsumformula} (ii) together with the relation $P_{j,k-j} = Z_{j,k-j} + Z_{k-j,j} + Z_k$.
\end{proof}

The even weight analogue of Theorem \ref{thm:dzoddbasis}  is given by the following 
\begin{theorem}\label{thm:dzevenbasis}  For even $k\geq 2$ we have $\dim_{\Q} \dz_k =  \frac{k}{2}$ and the set of $Z_{\text{od},\text{od}}$, i.e. 
    \begin{align*}
    \left\{ Z_{1,k-1}, Z_{3,k-3},\dots ,Z_{k-1,1} \right\},
    \end{align*}
    is a basis of $\dz_k$.
\end{theorem}
\begin{proof} We will give a proof of this below after introducing some further notation. Also, similar to the odd weight case, we will give explicit formulas to express the $P_{k_1,k_2}$ and $Z_{\text{ev},\text{ev}}$ in terms of $Z_{\text{od},\text{od}}$, which can be found in \eqref{eq:pevevzevevexplicitformula}. By Theorem \ref{thm:dzsumformula} (ii) this is already known for $Z_k$.
\end{proof}

Recall that by $V_k \subset \Q[X,Y]$ we denote the space of homogeneous polynomials of degree $k-2$. In this section we just consider the case when $k$ is even. In this case the action of $\GL_2(\Z)$ from the previous section induces an action of $\Gamma := \PGL_2(\Z)$ on $V_k$. On $V_k$ we define the following  pairing for $r,s,m,n \geq 1$, $r+s=m+n=k$
\begin{align}\label{eq:pairingdef}
\langle X^{r-1} Y^{s-1} \,,\, X^{m-1} Y^{n-1} \rangle = \frac{(-1)^{r-1}}{\binom{k-2}{m-1}} \delta_{(r,s), (n,m)}\,.
\end{align}
\begin{proposition} \label[proposition]{prop:pairing}\begin{enumerate}[\textup{(}i\textup{)}] 
        \item For even $k\geq 2$ the pairing \eqref{eq:pairingdef} on $V_k$  is bilinear, symmetric and non-degenerate.
        \item         For even $k\geq 2$ and $F,G \in V_k$ and $\gamma \in \Gamma$ we have 
        \begin{align*}
        \langle F\!\!\mid\!\gamma \,,\, G\!\!\mid\!\gamma  \rangle =     \langle F,  G\rangle \,.
        \end{align*}
    \end{enumerate}
    
\end{proposition}
\begin{proof}
    This is Exercise \ref{ex13}. 
\end{proof}
Recall that we have $\gz_k \mid\! (1 + \epsilon) + \gr_k = \gz_k \mid\! T(1 +\epsilon)$. Defining the element \begin{align*}
\Delta =(T-1)(\epsilon +1)
\end{align*}
we see that this relation then is equivalent to $\gz_k \!\mid\!\! \Delta = \gr_k$.
Now denote for $\xi \in \Z[\GL_2(\Z)]$ by $\xi^\ast$ its adjoint action given by the anti-automorphism induced by $\gamma \mapsto \gamma^{-1}$. By the $\Gamma$-invariance of the pairing we then obtain
\begin{align*}
\langle F\!\!\mid\!\xi \,,\, G  \rangle =     \langle F\,,\, G \!\!\mid\!\xi^\ast \rangle\,.
\end{align*}
This means that for any $F\in V_k$ we get a relation in $\dz_k$ with $\Delta^\ast = (1+\epsilon)(T^{-1}-1)$  by 
\begin{align}\label{eq:fdeltaast}
\langle F \!\!\mid\!\Delta^\ast \,,\, \gz_{k}  \rangle = \langle F,  \gr_k\rangle \,.
\end{align}
In the following we want to prove Theorem \ref{thm:dzevenbasis}, i.e. show that any $Z_{\text{ev}, \text{ev}}$ can be expressed in terms of $Z_{\text{odd}, \text{odd}}$. To prove this we want to find for $m,n\geq 2$ polynomials $F_{m,n}(X,Y) \in V_{m+n}$ such that
\begin{align}\label{eq:fmncondition}
F_{m,n}\mid\!(\Delta^\ast \piod) \in  X^{n-1} Y^{m-1} \Q^\times \,, 
\end{align}
where $ \piod = \frac{1}{2} (1-\delta)$ is the projection to the odd part. This would imply that for some $\alpha_{k_1,k_2}  \in \Q$ and $\lambda \in \Q^\times$ we have with $m+n=k$
\begin{align}     \label{eq:relamongzzev}
\langle F_{m,n}\!\!\mid\!\Delta^\ast \,,\, \gz_{k}  \rangle  = \lambda  Z_{m,n} + \sum_{\substack{k_1+k_2=k\\k_1, k_2 \text{ odd}}} \alpha_{k_1,k_2}  Z_{k_1,k_2}  = \langle F_{m,n} \,,\, \gr_{k}  \rangle \,,
\end{align}
which then gives an expression for $Z_{m,n}$ in terms of $Z_k$ and $Z_{\text{odd}, \text{odd}}$. We will deal with all $m,n\geq 2$ at the same time and therefore consider elements in 
\begin{align*}
\mathcal{V}^{\text{ev}} = \sum_{\substack{m\geq 2, n \geq 0\\ m,n \text{ even}}} V_{m+n}\, M^{m-2} N^{n} \subset \Q[[X,Y,M,N]]\,.
\end{align*}
On this space, we extend the action of $\Z[\Gamma]$ and define the following element in $\mathcal{V}^{\text{ev}}$ 
\begin{align}
f_{M,N}(X,Y) := -\frac{1}{4} YN\cosh (M Y) \left(\cosh (X N) \coth \left(\frac{Y N}{2}\right)+\sinh (X N)\right)\,,
\end{align}
where 
\begin{align*}
\cosh(X) = \sum_{\substack{n\geq 0\\n \text{ even}}} \frac{X^{n}}{n!}\,,\quad \sinh(X) =\cosh'(X) = \sum_{\substack{n\geq 1\\n \text{ odd}}}\frac{X^{n}}{n!},\quad \frac{X}{2}\coth\left( \frac{X}{2}\right)  = 1+\sum_{n\geq 2} \frac{B_n}{n!} X^n\,.
\end{align*}
With this we define the polynomials $F_{m,n}(X,Y) \in V_{m+n}$ for $m\geq 2, n\geq 0$ by 
\begin{align*}
F_{M,N}(X,Y) &:= \sum_{\substack{m \geq 2, n \geq 0\\m,n \text{ even}}} F_{m,n}(X,Y) \frac{M^{m-2}}{(m-2)!} \frac{N^n}{n!}  \\
&= f_{M,N}(X,Y) + f_{M,N}(Y,X) - f_{M,N}(-X,X+Y)\,,
\end{align*}
where the right-hand side is exactly $f_{M,N}\mid (1+\epsilon - ST)$. These $F_{m,n}$ give a solution to \eqref{eq:fmncondition}:
\begin{lemma}\label{lem:fmn}
    We have
    \begin{align*}
    (F_{M,N}\mid\!(\Delta^\ast \piod))(X,Y)&=N Y \cosh (M Y) \sinh (N X) =  \sum_{\substack{m,n \geq 2\\m,n \text{ even}}}  n X^{n-1} Y^{m-1} \frac{M^{m-2}}{(m-2)!}\frac{N^n}{n!}  \,,
    \end{align*}
    i.e. we have for even $m,n \geq 2$
    \begin{align*}
    F_{m,n}\mid\!(\Delta^\ast \piod) =  n X^{n-1} Y^{m-1}\,.
    \end{align*}
\end{lemma}
\begin{proof} This follows by a direct (but tedious) calculation using trigonometric identities. 
\end{proof}
Using \eqref{eq:relamongzzev} this gives explicit expressions for $Z_{\text{ev}, \text{ev}}$, which is given by the following.

\begin{proposition} \label{prop:explicitevevinodod} For $m, n \geq 2$ even and $k=m+n$ we have
    \begin{align*}
    Z_{m,n} &= \frac{2}{1-m} \sum_{\substack{k_1+k_2=k\\k_1,k_2 \geq 1 \text{ odd}}}\left( \sum_{j=0}^{\min\{k_1-2,n\}} \binom{k-2-j}{m-2}\binom{k_1-1}{j} B_{n-j}  \right)\left(Z_{k_1,k_2}  + \frac{1}{2} Z_k\right) - \frac{1}{2} Z_k\,.
    \end{align*}
\end{proposition}
\begin{proof}
    By Lemma \ref{lem:fmn} the polynomials $F_{m,n}$ give the relation 
    \begin{align} 
    \langle F_{m,n}\!\!\mid\!\Delta^\ast \,,\, \gz_{k}  \rangle  = \frac{n (-1)^{n-1}}{\binom{k-2}{m-1}}  Z_{m,n} + \sum_{\substack{k_1+k_2=k\\k_1, k_2 \text{ odd}}} \alpha_{k_1,k_2}  Z_{k_1,k_2}  = \lambda Z_k\,,
    \end{align}
    where the coefficients $\alpha_{k_1,k_2} \in \Q$ and $\lambda = \langle F_{m,n} \,,\, \gr_{k}  \rangle $ can be obtained by considering the coefficients of the generating series $F_{M,N}$. Notice that these $\alpha_{k_1,k_2}$ and $\lambda$ are not the $\alpha^{m,n}_{k_1,k_2}$ and $\lambda_{m,n}$ defined in \eqref{eq:defalphamn} below, but differ from them by the above normalization factor and a sign. The same formula can also be found in \cite[(7)]{GKZ}.
\end{proof}
Clearly Proposition \ref{prop:explicitevevinodod} also gives formulas for $P_{k_1,k_2} = Z_{k_1,k_2} + Z_{k_2,k_1} + Z_k$ in terms of $Z_{\text{od}, \text{od}}$. For this we can define for even $m,n\geq 2$ and odd $k_1,k_2\geq 1$ with $k_1+k_2=m+n$ the coefficients 
\begin{align}\begin{split}\label{eq:defalphamn}
\alpha^{m,n}_{k_1,k_2} &=  \frac{2}{1-m}  \sum_{j=0}^{\min\{k_1-2,n\}} \binom{k-2-j}{m-2}\binom{k_1-1}{j} B_{n-j} \,,\\
\lambda_{m,n} &=-\frac{1}{2} + \frac{1}{2} \sum_{\substack{k_1+k_2=k\\k_1,k_2 \geq 1 \text{ odd}}} \alpha^{m,n}_{k_1,k_2}   \,.
\end{split}
\end{align}
From the odd/odd sum formula in Theorem \ref{thm:dzsumformula} (ii) we therefore get for even $m,n\geq 2$ the following two explicit formulas
\begin{align}\begin{split}\label{eq:pevevzevevexplicitformula}
Z_{m,n} &=\sum_{\substack{k_1+k_2=k\\k_1\geq 3,k_2 \geq 1 \text{ odd}}}  \left( \alpha^{m,n}_{k_1,k_2} + 4 \lambda_{m,n}\right) Z_{k_1,k_2}  \,,  \\
P_{m,n} &= \sum_{\substack{k_1+k_2=k\\k_1\geq 3,k_2 \geq 1 \text{ odd}}}  \left( \alpha^{m,n}_{k_1,k_2} +\alpha^{n,m}_{k_1,k_2}+ 4\lambda_{m,n}+4\lambda_{n,m} + 4\right) Z_{k_1,k_2}\,.
\end{split}
\end{align}
This proves Theorem \ref{thm:dzevenbasis}, i.e., that the $Z_{\text{od}, \text{od}}$ form a basis for the space $\dz_k$ when $k$ is even.

\section{Realizations and combinatorial double Eisenstein series}\label{sec:realizations}
From now on we will be interested in realizations in a $\Q$-algebra $A$ and we will consider  all weights $k$ at the same time. Assume we have a family of homomorphisms $\varphi= ( \varphi_k )_{k\geq 1}$ where $\varphi_k$ is a realization of $\dz_k$ in $A$, i.e. $\varphi_k \in \Hom_\Q(\dz_k, A)$. Such a family $\varphi$ will be called a \emph{realization of $\dz$ in $A$}. For a realization $\varphi$ we will write $\varphi(P_{k_1,k_2}),\varphi(Z_{k_1,k_2})$ and $\varphi(Z_{k})$ instead of  $\varphi_{k_1+k_2}(P_{k_1,k_2}), \varphi_{k_1+k_2}(Z_{k_1,k_2}), \varphi_{k}(Z_{k})$, since the weight is always clear from the indices.

\begin{lemma} \label[lemma]{lem:realina} Assume we have power series $P(X,Y), Z(X,Y),Z(X)\in A[[X,Y]]$ such that 
    \begin{align*}
    P(X,Y) &= Z(X,Y) + Z(Y,X) + \frac{Z(X)- Z(Y)}{X-Y} - z(2)
    = Z(X+Y,Y) + Z(X+Y,X)\,,
    \end{align*}
    where $Z(X)=\sum_{k\geq 1} z(k) X^{k-1}$.
    Then $\varphi$ defined by 
    \begin{align*}
    \varphi(Z_k) &= z(k) - \delta_{k,2} z(2)\,,\\
    \varphi(Z_{k_1,k_2}) &= \text{coefficient of }    X^{k_1-1} Y^{k_2-1} \text{ in } Z(X,Y)\,,\\
    \varphi(P_{k_1,k_2}) &= \text{coefficient of }    X^{k_1-1} Y^{k_2-1} \text{ in } P(X,Y)
    \end{align*}
    gives a realization of $\dz$ in $A$.
\end{lemma}
\begin{proof} This is the same argument as for the realization $\varphi_\zeta$ given by shuffle regularized multiple zeta values.
\end{proof}

\begin{theorem}\label{thm:dzeuler} Let $\varphi$ be a realization of $\dz$ in a $\Q$-algebra $A$. For  $k\geq 1$ we write
    \[ \mathsf{Z}(k)=\varphi(Z_{k})  + \delta_{k,2} \mathsf{Z}(2)\]
    for some fixed element $\mathsf{Z}(2) \in A$.
    \begin{enumerate}[\textup{(}i\textup{)}] 
        \item Assume that for even $k_1,k_2 \geq 2$ we have 
        \begin{align*}
        \varphi(P_{k_1,k_2})= \mathsf{Z}(k_1) \mathsf{Z}(k_2) \,.
        \end{align*}
        Then for even $k\geq 2$ we have  $\mathsf{Z}(k) \in \Q[\mathsf{Z}(2)]$ and more precisely  we obtain for $m\geq 1$
        \begin{align}\label{eq:phieuler}
        \mathsf{Z}(2m) = -\frac{B_{2m} }{2 (2m)!}  \left( -24 \mathsf{Z}(2) \right)^m\,. 
        \end{align}
        \item Assume that there exists a derivation $\partial \in \Der(A)$ such that for all even $k_1,k_2 \geq 2$  
        \begin{align*}
        \varphi(P_{k_1,k_2})= \mathsf{Z}(k_1) \mathsf{Z}(k_2) + \frac{\delta_{k_1,2}}{2 k_2} \,\partial \mathsf{Z}(k_2) + \frac{\delta_{k_2,2}}{2 k_1} \ \,\partial \mathsf{Z}(k_1)\,.
        \end{align*}
        Then for even $k\geq 2$ we have $\mathsf{Z}(k) \in \Q[\mathsf{Z}(2), \mathsf{Z}(4),\mathsf{Z}(6)]$. Moreover, we get
        \begin{align}\begin{split}
        \label{eq:dzramaujandiff}
        \partial \mathsf{Z}(2) &= 5 \,\mathsf{Z}(4) - 2\, \mathsf{Z}(2)^2\,,\\ \partial \mathsf{Z}(4) &= 14\,\mathsf{Z}(6) -8\, \mathsf{Z}(2)\mathsf{Z}(4)\,,\\
        \partial\mathsf{Z}(6) &= \frac{120}{7} \mathsf{Z}(4)^2 -12\mathsf{Z}(2)\mathsf{Z}(6)\,,\\
        \partial^3 \mathsf{Z}(2) &= 36 (\partial \mathsf{Z}(2) )^2 -24 \mathsf{Z}(2) \partial^2 \mathsf{Z}(2) \,,
        \end{split}
        \end{align}
        and therefore the space $\Q[\mathsf{Z}(2), \mathsf{Z}(4),\mathsf{Z}(6)]=\Q[\mathsf{Z}(2), \partial \mathsf{Z}(2), \partial^2 \mathsf{Z}(2)]$ is closed under $\partial$.
        \item Assume that for even $k_1,k_2 \geq 4$ 
        \begin{align*}
        \varphi(P_{k_1,k_2})= \mathsf{Z}(k_1) \mathsf{Z}(k_2) \,.
        \end{align*}
        Then for even $k\geq 4$ we have  $ \mathsf{Z}(k) \in \Q[\mathsf{Z}(4),\mathsf{Z}(6)]$. 
    \end{enumerate}
\end{theorem}
\begin{proof} We first show (ii). The explicit formulas in \eqref{eq:dzramaujandiff} follow from the relations among the $P_{\text{ev},\text{ev}}$ and $Z_{\text{ev}}$ in Theorem \ref{thm:relpevevk}. In the smallest weights these relations read
    \begin{align*}
    P_{2,2} = \frac{5}{2} Z_4\,,\qquad P_{2,4} = \frac{7}{4} Z_6\,,\qquad P_{2,6} = \frac{5}{3} Z_8\,,\qquad P_{4,4} = \frac{7}{6} Z_8\,.
    \end{align*}
    Applying $\varphi$ to the first one gives $\mathsf{Z}(2)^2 + \frac{1}{2} \partial \mathsf{Z}(2) = \frac{5}{2} \mathsf{Z}(4)$, which is the first equation. The second one gives $\mathsf{Z}(2)\mathsf{Z}(4) + \frac{1}{8} \partial \mathsf{Z}(4) = \frac{7}{4} \mathsf{Z}(6)$, which is the second equation. For the third one we need both relations in weight eight: from $P_{4,4}$ we get $\mathsf{Z}(4)^2 = \frac{7}{6} \mathsf{Z}(8)$ and from $P_{2,6}$ we get $\mathsf{Z}(2)\mathsf{Z}(6) + \frac{1}{12} \partial \mathsf{Z}(6) = \frac{5}{3} \mathsf{Z}(8)$. Eliminating $\mathsf{Z}(8)$ gives the third equation. The Chazy equation then follows from these three by a direct calculation. One of the relations in Theorem \ref{thm:relpevevk} (Corollary \ref{cor:prodzk}) gives for $k\geq 4$
    \begin{align*}
    \mathsf{Z}(k) = \frac{2}{k+1} \sum_{\substack{j=2\\j\text{ even}}}^{k-2} \varphi( P_{j, k-j} ) =\frac{2}{k+1} \sum_{\substack{j=2\\j\text{ even}}}^{k-2} \mathsf{Z}(j) \mathsf{Z}(k-j) + \frac{2}{(k+1)(k-2)} \partial \mathsf{Z}(k-2) \,.
    \end{align*}
    Inductively we see, together with \eqref{eq:dzramaujandiff},  that $\mathsf{Z}(k) \in \Q[\mathsf{Z}(2), \mathsf{Z}(4),\mathsf{Z}(6)] $. The first statement in (i) is just a special case of (ii) for $\partial \equiv 0$.  The explicit formula \eqref{eq:phieuler} can be obtained by considering the generating series of the right hand side and showing that it satisfies the same recursion as the generating series of the left hand side (Exercise \ref{ex12}).
    For (iii) consider the equation \eqref{eq:relpevevz} in the cases $(k_1,k_2) = (m,2)$ and $(k_1,k_2)=(m-1,3)$. Taking $m-1$ times the first case and subtracting twice the second gives for even $k = m+2 \geq 8$ a relation of $Z_k$ as a linear combination of $P_{j,k-j}$ with $j,k-j \geq 4$ even. Recursively it follows that $ \mathsf{Z}(k) \in \Q[\mathsf{Z}(4),\mathsf{Z}(6)]$.
\end{proof}

Since the multiple zeta realization gives a realization of $\dz$ in $\R[T]$ which satisfies condition (i) with $\mathsf{Z}(2)=\zeta(2) = -\frac{(2\pi i)^2}{24}$, we obtain from \eqref{eq:phieuler} the Euler relation  $\zeta(2m) = -\frac{B_{2m}}{2 (2m)!} (2\pi i)^{2m}$. Notice that in depth one the regularized values do not involve $T$, so both sides of this relation are real numbers.

The first three equations in \eqref{eq:dzramaujandiff} are known as \emph{Ramanujan's differential equations} and the last one is called the \emph{Chazy equation} (e.g. see \cite[Section 5.1]{Za3}). One set of solutions for \eqref{eq:dzramaujandiff} (and actually a description of almost all solutions for the Chazy equation) is given by the Eisenstein series $\mathsf{Z}(k) = G_k$, which we defined for $k\geq 1$ by
\begin{align*}
G_k = \beta(k) + \g(k)\,.
\end{align*}
Recall that, for $k\geq 0$, we have $\beta(k)=-\frac{B_k}{2k!}$ when $k$ is even and $\beta(k)=0$ when $k$ is odd. In particular, $\beta(0)=-\frac{1}{2}$ and $\beta(1)=0$.

In the following, we want to introduce a realization which satisfies the condition in (ii) and is exactly given by the Eisenstein series in depth one. For this, we first introduce a realization with $\mathsf{Z}(k) =\g(k)$ and then a realization with $\mathsf{Z}(k) = \beta(k)$. Since $\Hom_\Q(\dz_k, A)$ are groups, we can add realizations and obtain another realization. Their sum gives the realization of the Eisenstein series. 
This realization satisfies the conditions of (ii) and (iii) in Theorem \ref{thm:dzeuler}. It therefore gives a proof that every Eisenstein series $G_k$ with $k\geq 4$ even can be written as a polynomial in $G_4$ and $G_6$.

\subsection{Combinatorial double Eisenstein series}
In this section, we want to introduce a realization $\varphi_G$ of $\dz$ in the $\Q$-algebra $A=\Q[[q]]$, which is given by the Eisenstein series $G(k) = \beta(k) + \g(k)$ in depth one. As mentioned before, we introduce two realizations, one for the constant term, denoted by $\varphi_\beta$, and one for the ``$g$-part'', denoted by $\varphi_{\g}$. Related modified double $q$-zeta values and their connection to modular forms are studied in \cite{Ba9}. The combinatorial double Eisenstein realization is the sum
\begin{align*}
\varphi_G = \varphi_\beta + \varphi_{\g}\,.
\end{align*}
We will start with the realization $\varphi_{\g}$, which in depth one will be given by $\varphi_{\g}(Z_k) = \g(k)$ for $k \neq 2$. To use Lemma \ref{lem:realina} we need to consider the generating series, which were defined in Section \ref{subsec:qmzv} by
\begin{align*}
\ggen(X_1,\dots,X_r) = \sum_{k_1,\dots,k_r\geq 1} \g(k_1,\dots,k_r) X_1^{k_1-1} \cdots X_r^{k_r-1}\,.
\end{align*}    
In Lemma \ref{lem:genexpr} we saw that we have the following two explicit expressions
\begin{align}\begin{split}
\label{eq:gtwoexpression2}
\ggen(X_1,\dots,X_r) 
&= \sum_{m_1> \dots > m_r > 0}  \frac{e^{X_1} q^{m_1}}{1-e^{X_1}q^{m_1}} \cdots  \frac{e^{X_r} q^{m_r}}{1-e^{X_r}q^{m_r}} \\
&= \sum_{m_1> \dots > m_r > 0}  \frac{e^{m_1 X_r} q^{m_1}}{1-q^{m_1}} \frac{e^{m_2 (X_{r-1}-X_r)} q^{m_2}}{1-q^{m_2}} \cdots \frac{e^{m_r (X_{1}-X_2)} q^{m_r}}{1-q^{m_r}}\,.
\end{split}
\end{align}
To express their product in a suitable way we introduce the following series
\begin{align*}
\gb(X)&=\frac{1}{2} \left(\frac{1}{X} - \frac{1}{e^X-1}-\frac{1}{2}\right) =\sum_{k\geq 1} \beta(k) X^{k-1} =-\sum_{k\geq 2} \frac{B_k}{2k!}X^{k-1} =  \sum_{m\geq 1}  \frac{\zeta(2m)}{(2\pi i)^{2m}} X^{2m-1}\,,
\end{align*}
which will also give the depth one part of the realization $\varphi_\beta$.
\begin{lemma} \label{lem:gprodprep} We have
    \begin{align*}
    \begin{split}
    \ggen(X) \ggen(Y)
    &=\ggen(X,Y) + \ggen(Y,X) +\frac{\ggen(X) - \ggen(Y)}{X-Y}\\
    &\quad+ \Big(\gb(Y-X) - \gb(X-Y)\Big)\left(\ggen(X)-\ggen(Y)\right)
    -\frac{1}{2}\left(\ggen(X)+\ggen(Y)\right) \\
    &= \ggen(X+Y,X) + \ggen(X+Y,Y) -\ggen(X+Y)+ (X+Y)\ggen'(X+Y)+\g(2)\,,
    \end{split}
    \end{align*}
    where we write $\ggen'(X)=  q \frac{d}{dq}\sum_{k\geq 1} \g(k) \frac{X^{k-1}}{k}$.
\end{lemma}    
\begin{proof}
    In Proposition \ref{prop:ggenproddep1} we showed by using \eqref{eq:gtwoexpression2} that
    \begin{align*}
    \ggen(X) \ggen(Y)  &= \ggen(X,Y)  + \ggen(Y,X) + \frac{1}{e^{X-Y}-1} \ggen(X) + \frac{1}{e^{Y-X}-1}\ggen(Y) \\
    &= \ggen(X+Y,X) + \ggen(X+Y,Y) -\ggen(X+Y)  +  q \frac{d}{dq}\sum_{k\geq 1} \g(k) \frac{(X+Y)^k}{k} + \g(2)\,,
    \end{align*}
    from which the statement follows by using the definition of $\gb(X)$.
\end{proof}

\begin{theorem}\label{thm:grealization} Define the following generating series
    \begin{align*}
    \gh(X,Y) &= \ggen(X,Y)-\left(\gb(X-Y)+\frac{1}{2}\right)\ggen(X) +\gb(Y) \ggen(X) +\gb(X-Y)\ggen(Y) \\
    &+ \frac{1}{2}(X - Y) \ggen'(Y)+  \frac{1}{2}X \ggen'(X)+ \frac{1}{2} \g(2)\,,
    \end{align*}
    where as before $\ggen'(X)=  q \frac{d}{dq}\sum_{k\geq 1} \g(k) \frac{X^{k-1}}{k}$. Then we have
    \begin{align}\begin{split}\label{eq:hrealization}
    P(X,Y) &= \gh(X,Y) + \gh(Y,X) + \frac{\ggen(X)- \ggen(Y)}{X-Y} - \g(2)
    = \gh(X+Y,Y) + \gh(X+Y,X)\,,    
    \end{split}
    \end{align}
    where 
    \begin{align*}
    P(X,Y)=\ggen(X)\ggen(Y) + \gb(X)\ggen(Y) + \gb(Y)\ggen(X) + \frac{1}{2}\left(\ggen'(X)Y + \ggen'(Y)X \right)\,.
    \end{align*}
    In particular  \eqref{eq:hrealization} gives a realization $\varphi_{\g}$ of $\dz$ in $q\Q[[q]]$  by  Lemma \ref{lem:realina}.
\end{theorem}
\begin{proof}
    This follows by a direct calculation from Lemma \ref{lem:gprodprep} but is also given in \cite[Theorem~7]{GKZ}. 
\end{proof}

We will now give the realization $\varphi_\beta$, which will give the constant term for the combinatorial double Eisenstein series.
\begin{theorem}\label{thm:betarealization}  With $\gb(X)=\sum_{k\geq 1} \beta(k) X^{k-1}= \frac{1}{2} \left(\frac{1}{X} - \frac{1}{e^X-1}-\frac{1}{2}\right)$ and 
    \begin{align*}
    &\gb(X,Y) = \sum_{k_1,k_2\geq 1} \beta(k_1,k_2) X^{k_1-1} Y^{k_2-1}\\
    &:= \frac{1}{3} \big(\gb(X)+\gb(X-Y) \big) \gb(Y) - \frac{5}{12} \frac{\gb(X)-\gb(Y)}{X-Y}+\frac{\gb(X)-\gb(X-Y)}{4Y} -\frac{\gb(Y)-\gb(Y-X)}{12X}- \frac{1}{96}
    \end{align*}
    we have 
    \begin{align*}
    \gb(X)\gb(Y) &= \gb(X,Y) +\gb(Y,X) + \frac{\gb(X)-\gb(Y)}{X-Y}   - \beta(2)
    = \gb(X+Y,Y) + \gb(X+Y,X)\,.
    \end{align*}
    In particular this gives a realization $\varphi_\beta$ of $\dz$ in $\Q$ by  Lemma \ref{lem:realina}.
\end{theorem}
\begin{proof}
    Again this can be checked explicitly by using the definition 
    of $\gb(X)$.
    A more systematic point of view, which we will make more explicit later, is that the hyperbolic cotangent  
    \begin{align*}
    F(X) = -\frac{1}{2}\frac{1}{X} + \gb(X) = -\frac{1}{4}\coth\left( \frac{X}{2} \right)
    \end{align*}
    satisfies the Fay identity 
    \begin{align}\label{eq:cotfay}
    F(X)F(Y) + F(X-Y)F(X)+F(-Y)F(X-Y) = \frac{1}{16}\,.
    \end{align}
    Writing $G(X,Y)=F(X)F(Y)$ the equation \eqref{eq:cotfay} can be written as $G\mid\!(1+U+U^2)=\frac{1}{16}$, which leads to a connection of period polynomials for modular forms. We will see later that basically any such $G$ gives rise to a realization of $\dz$ (Theorem \ref{thm:prop5supp}).
\end{proof}

\begin{definition} \begin{enumerate}[\textup{(}i\textup{)}]
        \item We define the \emph{combinatorial Eisenstein realization} of $\dz$ in $\Q[[q]]$ by
        \begin{align*}
        \varphi_G = \varphi_\beta + \varphi_{\g}\,,
        \end{align*}
        where the realizations $\varphi_\beta$ and $\varphi_{\g}$ are given by Theorems~\ref{thm:betarealization} and~\ref{thm:grealization}.
        \item     For $k_1,k_2\geq 1$ the \emph{combinatorial double Eisenstein series} $G(k_1,k_2)\in \Q[[q]]$ are defined by $G(k_1,k_2) = \varphi_G(Z_{k_1,k_2})$, i.e. they are explicitly given by
        \begin{align*}
        \sum_{k_1,k_2\geq 1}G(k_1,k_2) X^{k_1-1} Y^{k_2-1} := \gb(X,Y) + \gh(X,Y)\,,
        \end{align*}
        where $\gb(X,Y)$ and  $\gh(X,Y)$ are given in Theorem \ref{thm:betarealization}  and \ref{thm:grealization}. 
    \end{enumerate}
    
\end{definition}

Notice that we have 
\begin{align*}
\varphi_G(P_{k_1,k_2}) &=  \beta(k_1)\beta(k_2) + \g(k_1) \g(k_2)+\beta(k_1) \g(k_2) + \beta(k_2) \g(k_1) + \frac{\delta_{k_1,2}}{2k_2} q\frac{d}{dq}\g(k_2)+ \frac{\delta_{k_2,2}}{2k_1} q\frac{d}{dq}\g(k_1)\\
&= G_{k_1} G_{k_2} + \frac{\delta_{k_1,2}}{2k_2} q\frac{d}{dq}G_{k_2}+ \frac{\delta_{k_2,2}}{2k_1} q\frac{d}{dq} G_{k_1}\,,
\end{align*}
and therefore the combinatorial Eisenstein realization satisfies the conditions in (ii), with $\partial = q\frac{d}{dq}$, and (iii) of Theorem \ref{thm:dzeuler}. This gives a proof of Ramanujan's differential equations and the fact that for any even $k\geq 4$ we have $G_k \in \Q[G_4,G_6]$.
\begin{proposition}\label{prop:cdesareqana} The combinatorial double Eisenstein series are modified $q$-analogues of the double zeta values, i.e. for $k_1\geq 2, k_2\geq 1$ we have 
    \begin{align*}
    \lim_{q\rightarrow 1} (1-q)^{k_1+k_2} G(k_1,k_2) = \zeta(k_1,k_2)\,.
    \end{align*}
\end{proposition}
\begin{proof}
Put $K=k_1+k_2$. Since $\gb(X)\in X\Q[[X]]$, comparison of the
coefficients in Theorems \ref{thm:grealization} and
\ref{thm:betarealization} gives
\begin{align*}
G(k_1,k_2)-\g(k_1,k_2)
=\beta(k_1,k_2)+\sum_{j=1}^{K-1}c_j\g(j)
+c\,q\frac{d}{dq}\g(K-2)
\end{align*}
for some $c,c_1,\ldots,c_{K-1}\in\Q$. Indeed, the terms involving
$\gb$ give the $\g(j)$ with $j\leq K-1$, while multiplication by
$X-Y$ or $X$ in the two derivative terms leaves only the derivative
of $\g(K-2)$. The last term $\frac12\g(2)$ in
Theorem \ref{thm:grealization} does not contribute since $K\geq3$.

We show that all correction terms have zero modified limit of weight
$K$. Write $q=e^{-t}$ with $t>0$. Since a positive integer $n$ has
at most $2\sqrt n$ divisors, we have
\begin{align*}
\sigma_{j-1}(n)\leq 2n^{j-\frac12}.
\end{align*}
It follows from \eqref{eq:smallg}, by comparison with an integral, that
\begin{align*}
\g(j)&=O\bigl(t^{-j-\frac12}\bigr),&
q\frac{d}{dq}\g(j)&=O\bigl(t^{-j-\frac32}\bigr)
\qquad(t\longrightarrow0^+).
\end{align*}
Since $1-q\sim t$, we therefore obtain, for $j\leq K-1$,
\begin{align*}
(1-q)^K\g(j)&\longrightarrow0,&
(1-q)^Kq\frac{d}{dq}\g(K-2)&\longrightarrow0.
\end{align*}
The constant term $\beta(k_1,k_2)$ also vanishes after multiplication
by $(1-q)^K$. Finally, $(k_1,k_2)$ is admissible, and hence
Proposition \ref{prop:gisqmzv} gives
\begin{align*}
\lim_{q\rightarrow1}(1-q)^K\g(k_1,k_2)
=\zeta(k_1,k_2).
\end{align*}
This proves the claim.
\end{proof}

\begin{proposition}\label{prop:mkisgodod} Any modular form with rational coefficients can be written as a linear combination of $G(\text{odd},\text{odd})$, i.e.  for even $k\geq 4$ we have 
    \begin{align*}
    \mf^\Q_k   \subset \langle G(j,k-j) \mid j=3,5,\dots,k-1\rangle_{\Q}\,.
    \end{align*}
\end{proposition}
\begin{proof}
    This follows from the fact that any modular form can be written as a linear combination of products of Eisenstein series (Theorem \ref{thm:mf}) together with the explicit formulas for $P_{k_1,k_2}$ in  \eqref{eq:pevevzevevexplicitformula}.
\end{proof}

\begin{lemma} For even $k\geq 4$ the set 
    \begin{align*}
    \left\{ G_k \right\} \cup \left\{ G_{2j} G_{k-2j} \mid  \left\lfloor \frac{k-2}{6} \right\rfloor + 2 \leq j \leq \left \lfloor \frac{k}{4} \right\rfloor \right\}
    \end{align*}
    forms a basis of $\mf_k^\Q$.
\end{lemma}
\begin{proof}
    This is \cite[Corollary~2]{HST}.
\end{proof}
By the lemma we get a basis of $S^\Q_k$ by
\begin{align*}
\left\{ G_{2j} G_{k-2j} - \frac{\beta(2j)\beta(k-2j)}{\beta(k)}G_k \mid  \left\lfloor \frac{k-2}{6} \right\rfloor + 2 \leq j \leq \left \lfloor \frac{k}{4} \right\rfloor \right\}\,.
\end{align*}
Since $\varphi_G(P_{m,n})=G_mG_n$ for $m,n\geq4$, the formulas in
\eqref{eq:pevevzevevexplicitformula} make this basis completely explicit.
For even $m,n\geq4$, put $k=m+n$. For odd $3\leq r\leq k-1$ define
\begin{align*}
A_r^{m,n}
&=\frac{2}{1-m}\sum_{\ell=0}^{\min\{r-2,n\}}
\binom{k-2-\ell}{m-2}\binom{r-1}{\ell}B_{n-\ell}\\
&\quad+\frac{2}{1-n}\sum_{\ell=0}^{\min\{r-2,m\}}
\binom{k-2-\ell}{n-2}\binom{r-1}{\ell}B_{m-\ell}
\end{align*}
and
\begin{align}\label{eq:explicitgododcoeff}
C_r^{m,n}
=A_r^{m,n}
+2\sum_{\substack{s=3\\s\text{ odd}}}^{k-1}A_s^{m,n}
+2\binom{k}{m}\frac{B_mB_n}{B_k}\,.
\end{align}

\begin{theorem}\label{thm:explicitgododbasis}
For even $k\geq4$ the set
\begin{align*}
\left\{
\sum_{\substack{r=3\\r\text{ odd}}}^{k-1}
C_r^{2j,k-2j}G(r,k-r)
\ \middle|\ 
\left\lfloor\frac{k-2}{6}\right\rfloor+2
\leq j\leq\left\lfloor\frac{k}{4}\right\rfloor
\right\}
\end{align*}
forms a basis of $S_k^\Q$. More precisely, for even $m,n\geq4$ we have
\begin{align}\label{eq:explicitgododbasis}
G_mG_n-\frac{\beta(m)\beta(n)}{\beta(m+n)}G_{m+n}
=\sum_{\substack{r=3\\r\text{ odd}}}^{m+n-1}
C_r^{m,n}G(r,m+n-r)\,.
\end{align}
\end{theorem}

\begin{proof}
By \eqref{eq:pevevzevevexplicitformula}, the coefficient of
$G(r,k-r)$ in $G_mG_n=\varphi_G(P_{m,n})$ is
\begin{align*}
A_r^{m,n}+4\lambda_{m,n}+4\lambda_{n,m}+4
=A_r^{m,n}+2\sum_{\substack{s=3\\s\text{ odd}}}^{k-1}A_s^{m,n}\,.
\end{align*}
Here we used the definition of $\lambda_{m,n}$ in
\eqref{eq:defalphamn} and the fact that the coefficient with first
index $1$ vanishes. Moreover,
\begin{align*}
-\frac{\beta(m)\beta(n)}{\beta(k)}G_k
=2\binom{k}{m}\frac{B_mB_n}{B_k}
\sum_{\substack{r=3\\r\text{ odd}}}^{k-1}G(r,k-r)
\end{align*}
by $\beta(a)=-B_a/(2a!)$ and the odd sum formula in
Theorem~\ref{thm:dzsumformula}~(ii). Adding the two formulas gives
\eqref{eq:explicitgododbasis}. The left hand sides with
$m=2j$ and $n=k-2j$ form the basis displayed above, which proves the
statement.
\end{proof}

Using Proposition \ref{prop:cdesareqana} and the fact that modular forms are also modified $q$-analogues of their constant terms gives an explanation for the factor $\odd(X) \odd(X) - \Sx(X)$ in the Broadhurst-Kreimer conjecture (Conjecture \ref{conj:brokrei}).

\begin{ex}
After rescaling the basis elements, the first two cases are as follows.
    \begin{enumerate}[\textup{(}i\textup{)}] 
        \item 
        A basis for $S^\Q_{12}$ is given by
        \begin{align*}
        G(3,9)-\frac{23825}{5197} G(5,7)-\frac{41431 }{10394}G(7,5)+\frac{360
        }{5197}G(9,3)+G(11,1)\,.
        \end{align*}
        \item 
        A basis for $S^\Q_{16}$ is given by
        \begin{align*}
        G(3,13)&-\frac{279116 G(5,11)}{78967}-\frac{2125607 G(7,9)}{315868}
        -\frac{154671 G(9,7)}{22562}\\
        &\quad-\frac{1040507 G(11,5)}{315868}
        +\frac{38573 G(13,3)}{157934}+G(15,1)\,.
        \end{align*}
    \end{enumerate}
\end{ex}

\section{Period polynomials}\label{sec:periodpolynomials}
In this section, we will introduce period polynomials and their connection to relations in the formal double zeta space. Period polynomials are polynomials with complex coefficients associated to modular forms. By work of Eichler, Shimura, Manin and Zagier, the spaces of period polynomials are isomorphic to the space of modular forms. For more details on this topic, see \cite{La} and \cite{Za1}. Throughout this section we assume that $k$ is even.

For a cusp form $f \in S_k$ we define its \emph{period polynomial} as the following polynomial in $\C \otimes V_k$
\begin{align}\label{eq:defperiodpol}
P_f(X,Y) = \int_0^{i \infty} (X-Y\tau)^{k-2} f(\tau) \,d\tau \,.
\end{align}

Recall that we defined on $V_k$ a right-action of $\GL_2(\Z)$ for a $\gamma = \left( \begin{smallmatrix}
a & b \\
c & d
\end{smallmatrix} \right) \in \GL_2(\Z)$ and $F\in V_k$ by 
\begin{align*}
(F{\mid\! \gamma})(X,Y) = F(aX + bY, cX + dY)\,.
\end{align*}
This action can then be extended linearly to an action of the group ring $\Z[\GL_2(\Z)]$ on $\C \otimes_\Q V_k$.

\begin{lemma} For a cusp form $f\in S_k$ and $\gamma \in \SLZ$ we have
    \begin{align}
    (P_f \mid\! \gamma)(X,Y) = \int_{\gamma^{-1}(0)}^{\gamma^{-1}(i \infty)} (X-Y\tau)^{k-2} f(\tau) \,d\tau\,,
    \end{align}
    where $\gamma(z) = \frac{a z + b}{c z + d}$ for $\gamma = \begin{psmallmatrix}a&b\\c&d \end{psmallmatrix}$ as in Section \ref{subsec:mfandbk}.
\end{lemma}
\begin{proof}
    This is Exercise \ref{ex14} (i).
\end{proof}
In particular, we see that (Exercise \ref{ex14} (ii))
\begin{align}\begin{split}
\label{eq:1pS1UU}
(P_f \mid\! (1+S))(X,Y) &= \left(\int_{0}^{\infty i} + \int_{\infty i}^{0}   \right) (X-Y\tau)^{k-2} f(\tau) \,d\tau = 0\,,\\
(P_f \mid\! (1+U+U^2))(X,Y) &= 0     \,.
\end{split}
\end{align}
Motivated by \eqref{eq:1pS1UU} we define the following subspace of $V_k$ for even $k\geq 2$  
\begin{align}\label{eq:defwk}
W_k = \ker(1+S) \cap \ker(1+U+U^2) \subset V_k\,.
\end{align}
Let $V^\pm_k$ denote the $\pm 1$-eigenspaces under the action of $\epsilon$, i.e. $V^+_k$ are the symmetric and $V^-_k$ are the antisymmetric polynomials. By $V^\text{ev}$ and $V^\text{od}$ we denote the $\pm1$-eigenspaces of $\delta$, i.e. the even and odd polynomials. 
With this we define the symmetric ($+$), antisymmetric ($-$), even (ev) and odd (od) parts of $W_k$ for $\bullet \in \{+,-,\text{ev}, \text{od}\}$ by 
\begin{align*}
W^\bullet_k = W_k \cap V^\bullet_k\,.
\end{align*}

\begin{lemma}\label[lemma]{lem:lewis}
    \begin{enumerate}[\textup{(}i\textup{)}] 
        \item We have 
        \begin{align*}
        W^+_k = W^\text{od}_k\,, \qquad W^-_k = W^\text{ev}_k\,. 
        \end{align*}
        \item The spaces $W_k$ and $W^\pm_k$ can also be written as
        \begin{align}\label{eq:lewis}
        W_k = \ker(1- T - T')\,,
        \end{align}
        where $T' = -U^2 S = \begin{psmallmatrix}1&0\\1&1 \end{psmallmatrix}$ and
        \begin{align*}
        W^\pm_k = \ker(1- T \mp T \epsilon)\,.
        \end{align*}
    \end{enumerate}
\end{lemma}
\begin{proof}
    This is Exercise \ref{ex14} (iii). 
\end{proof}
The equation \eqref{eq:lewis} is also called the \emph{Lewis equation}. The period polynomials of cusp forms satisfy the Lewis equation, and the following theorem states that there are no more relations. More precisely, the Eichler-Shimura theorem gives isomorphisms of $S_k$ to $W^\pm_k$ (modulo a subspace of dimension one in the even case) by sending a cusp form $f$ to the odd part $P^+_f$ and even part $P^-_f$ of its period polynomial $P_f$.

\begin{theorem}[Eichler-Shimura isomorphism]\label{thm:eichlershimura} The map $p^\pm_f: f \mapsto P^\pm_f$ induces isomorphisms
    \begin{align*}
    p^+_f : S_k \xrightarrow{\raisebox{-0.25 em}{\smash{\ensuremath{\sim}}}} \C \otimes W^+_k\,,\qquad  \qquad
    p^-_f: S_k \xrightarrow{\raisebox{-0.25 em}{\smash{\ensuremath{\sim}}}} \C \otimes \quotient{W^-_k}{\Q(X^{k-2}-Y^{k-2})}\,.
    \end{align*}
\end{theorem}

Moreover, one can show that everything is ``defined over $\Q$'', i.e. we can always find a basis of $S_k$ whose images under $p^\pm_f$ are elements in $W_k^\pm$. Later we will see that we can also define period polynomials for Eisenstein series $G_k$, which are given exactly by a multiple of $X^{k-2}-Y^{k-2}$ in the antisymmetric (even) case. Therefore we will see that $W_k^-$ is isomorphic to the space of all modular forms $\mf_k$.

\begin{ex} The odd and even parts of the period polynomial of $f=c \Delta \in S_{12}$ for some explicit $c\in \C$ are given by 
    \begin{align*}
    P^+_f(X,Y) &= XY (X^2-Y^2)^2(X^2-4Y^2)(4X^2-Y^2)\,,\\
    P^-_f(X,Y) &= \frac{36}{691} (X^{10} - Y^{10}) - X^2 Y^2 (X^2-Y^2)^3\,.
    \end{align*}
\end{ex}

\subsection{The space $\mathcal{P}^\text{ev}_k$ and its connection to $W_k^\pm$}
Let $\mathcal{P}^\text{ev}_k \subset \dz_k$ denote the space spanned by all $P_{\text{ev},\text{ev}}$.
\begin{align*}
\mathcal{P}^\text{ev}_k  &= \langle P_{m,n} \mid m,n \geq 2 \text{ even}, m+n = k \rangle_{\Q}
= \Q Z_k + \langle P_{m,n} \mid m,n \geq 4 \text{ even}, m+n = k \rangle_{\Q}\,,
\end{align*}
where the second equality follows from Corollary \ref{cor:prodzk}.
By \eqref{eq:pevevzevevexplicitformula} we can write any $P_{m,n}$ explicitly in terms of $Z_{\text{od},\text{od}}$. In \cite{GKZ} it was observed that $W_k^-$ is canonically isomorphic to $\mathcal{P}^\text{ev}_k$, where the isomorphism can be written down explicitly by using the coefficients of period polynomials. 
For a $p \in W_k^-$ we define the coefficients $\beta^p_{k_1,k_2} \in \Q$ for $k_1,k_2 \geq 1$, $k_1+k_2=k$ by
\begin{align*}
\sum_{\substack{k_1+k_2=k\\ k_1,k_2 \geq 1}} \binom{k-2}{k_1-1} \beta^p_{k_1,k_2} X^{k_1-1} Y^{k_2-1} := p(X+Y,Y), 
\end{align*}
i.e. these are the coefficients of $p\mid T$ divided by $\binom{k-2}{k_1-1}$.
\begin{theorem} \label{thm:wkminuspkev} For even $k\geq 4$ the following map is an isomorphism of $\Q$-vector spaces
    \begin{align*}
    W_k^- &\longrightarrow \mathcal{P}^\text{ev}_k\\
    p &\longmapsto \sum_{\substack{k_1+k_2=k\\k_1, k_2 \text{ odd}}}\beta^p_{k_1,k_2}  Z_{k_1,k_2}\,.
    \end{align*}
    Moreover the image can be written in terms of the generators of $\mathcal{P}^\text{ev}_k$ explicitly as 
    \begin{align*}
    \sum_{\substack{k_1+k_2=k\\k_1, k_2 \text{ odd}}}\beta^p_{k_1,k_2}  Z_{k_1,k_2} \equiv     \frac{1}{6} \sum_{\substack{k_1+k_2=k\\k_1, k_2 \text{ even}}}\beta^p_{k_1,k_2}  P_{k_1,k_2}  \quad  \mod  \Q Z_k\,.
    \end{align*}		
\end{theorem}

Before we give a proof of this theorem, we give some applications of it for our multiple zeta realization $\varphi_\zeta$ and our combinatorial Eisenstein realization $\varphi_G$.

\begin{theorem} \label{thm:pkevisotomk} For $k\geq 4$ the combinatorial Eisenstein realization gives an isomorphism 
    \begin{align*}
    \varphi_G : \mathcal{P}^\text{ev}_k \xrightarrow{\raisebox{-0.25 em}{\smash{\ensuremath{\sim}}}}  \mf^\Q_k \,.
    \end{align*}
\end{theorem}

\begin{proof} By the Eichler-Shimura isomorphism (Theorem \ref{thm:eichlershimura}) we know that $\dim W_k^- = \dim S_k + 1 = \dim \mf^\Q_k$, so by Theorem \ref{thm:wkminuspkev} it suffices to show surjectivity. By Theorem \ref{thm:mf} we know that the set 
    \begin{align*}
    \Q G_k + \langle G_m G_n \mid m,n \geq 4 \text{ even}, m+n = k \rangle_{\Q}
    \end{align*}
    generates $\mf^\Q_k$. But this is exactly the image of the generating set of $\mathcal{P}^\text{ev}_k$ under $\varphi_G$.
\end{proof}

Another consequence of Theorem \ref{thm:wkminuspkev} is the following.
\begin{corollary} \label{cor:exoticrel} Let $k\geq 4$ be even and $p\in W_k^-$.
    \begin{enumerate}[\textup{(}i\textup{)}] 
        \item We have 
        \begin{align*}
        \sum_{\substack{k_1+k_2=k\\k_1\geq 3\text{ odd}\\k_2\geq 1\text{ odd}}} \beta^p_{k_1,k_2}  \zeta(k_1,k_2) \in \Q \pi^k\,.
        \end{align*}
        \item We have 
        \begin{align}\label{eq:goddquasimodular}
        \sum_{\substack{k_1+k_2=k\\k_1, k_2 \text{ odd}}} \beta^p_{k_1,k_2}  G(k_1,k_2) \in \mf_k^\Q\,.
        \end{align}
    \end{enumerate}
\end{corollary}
\begin{proof}
    Since $p$ is antisymmetric, we have $\beta^p_{1,k-1}=p(Y,Y)=0$. The two statements follow from $\varphi_\zeta(\mathcal{P}^\text{ev}_k)= \Q \pi^k$ and $\varphi_G(\mathcal{P}^\text{ev}_k) =  \mf_k^\Q$ (Theorem \ref{thm:pkevisotomk} above).
\end{proof}

\begin{remark} After extending the maps above to $\C$, starting with a cusp form $f \in S_k$ and taking $p=P_f^-\in \C\otimes W_k^-$ in \eqref{eq:goddquasimodular} gives a modular form. By Theorem \ref{thm:pkevisotomk}, every modular form is obtained from some $p\in \C\otimes W_k^-$ in this way. For a long time the relationship between $f$ and the modular form obtained from $P_f^-$ was unknown. Tasaka proved in \cite[Theorem~1]{Tas} that one can get back exactly $f$ in \eqref{eq:goddquasimodular} by, roughly speaking, considering a slight modification of the isomorphism in Theorem \ref{thm:wkminuspkev}.
\end{remark}

Now we will give an analogous statement for the space $W_k^+$ for odd period polynomials and see that it is canonically dual to $\mathcal{P}^\text{ev}_k $: 
\begin{theorem} \label{thm:wkpluspkev}  For even $k\geq 4$ the following map is an isomorphism of $\Q$-vector spaces
    \begin{align*}
    W_k^+ &\longrightarrow \left\{ \varphi \in \Hom( \mathcal{P}^\text{ev}_k , \Q) \mid  \varphi(Z_k) = 0 \right\}\\
    \sum_{\substack{k_1+k_2=k\\ k_1,k_2 \geq 2\text{ even}}} p_{k_1,k_2} X^{k_1-1} Y^{k_2-1}  &\longmapsto ( \varphi: P_{k_1,k_2} \mapsto p_{k_1,k_2})\,.
    \end{align*}
\end{theorem}

In other words any polynomial in $W_k^+$ gives a realization $\varphi$ of $\mathcal{P}^\text{ev}_k$ in $\Q$ with $\varphi(Z_k) = 0$. Later we will see that we can extend the space $W_k^+$ to a space from which we will actually obtain any realization of $\mathcal{P}^\text{ev}_k$ in $\Q$.

To prove Theorems \ref{thm:wkminuspkev}  and \ref{thm:wkpluspkev} we will use the pairing on $V_k$ defined in \eqref{eq:pairingdef} by 
\begin{align*}
\langle X^{r-1} Y^{s-1} \,,\, X^{m-1} Y^{n-1} \rangle = \frac{(-1)^{r-1}}{\binom{k-2}{m-1}} \delta_{(r,s), (n,m)}\,.
\end{align*}
Assume we have the following element in $V_k$
\begin{align} \label{eq:polA}
A(X,Y) = \sum_{\substack{k_1+k_2=k\\ k_1,k_2 \geq 1}} \binom{k-2}{k_1-1} a_{k_1,k_2} X^{k_1-1} Y^{k_2-1} \,.
\end{align}
Pairing $A\! \mid\! S$ with some polynomial therefore gives
\begin{align}\label{eq:pairingaandc}
\langle A\!\!\mid \! S ,  		\sum_{\substack{k_1+k_2=k\\ k_1,k_2 \geq 1}} c_{k_1,k_2} X^{k_1-1} Y^{k_2-1}  \rangle = \sum_{\substack{k_1+k_2=k\\ k_1,k_2\geq 1}} a_{k_1,k_2}  c_{k_1,k_2} \,.
\end{align}
This will be used to prove the following lemma. 
\begin{lemma}  \label{lem:Arealindk} Let $A$ and $a_{k_1,k_2}$ be as in \eqref{eq:polA} and let $k\geq 2$ be even.
    \begin{enumerate}[\textup{(}i\textup{)}] 
        \item In $\dz_k$ the relation
        \begin{align*}
        \sum_{\substack{k_1+k_2=k\\ k_1,k_2\geq 1}} a_{k_1,k_2}  Z_{k_1,k_2} = \lambda Z_k\,,
        \end{align*}
        holds for some $\lambda\in\Q$ if and only if $A = H \mid (T'-1)$ for some polynomial $H\in V_k^+$. In this case we have
        \begin{align*}
        \lambda = \frac{1}{2} \langle H \!\mid\! \delta , \frac{X^{k-1}-Y^{k-1}}{X-Y}\rangle = \frac{k-1}{2} \int_0^1 H(t,1-t) dt\,.
        \end{align*}
        
        \item Assume in addition that $A\in V_k^+$. Then in $\dz_k$ the relation
        \begin{align*}
        \sum_{\substack{k_1+k_2=k\\ k_1,k_2\geq 1}} a_{k_1,k_2}  P_{k_1,k_2} = \mu Z_k\,,
        \end{align*}
        holds for some $\mu\in\Q$ if and only if $A = H \mid (1-S)$ for some  polynomial $H\in V^U_k \cap V_k^+$. In this case we have
        \begin{align*}
        \mu = \langle H , \frac{X^{k-1}-Y^{k-1}}{X-Y}\rangle \,.
        \end{align*}		
        Here $ V^U_k $ denotes the space of $U$ invariant polynomials, i.e. $H\!\mid\! U = H$. 
    \end{enumerate}
    
\end{lemma}
\begin{proof} For (i) first notice that by \eqref{eq:pairingaandc} we have
    \begin{align*}
    \langle A\!\!\mid \! S , \gz_k \rangle = \sum_{\substack{k_1+k_2=k\\ k_1,k_2\geq 1}} a_{k_1,k_2}  Z_{k_1,k_2} \,.
    \end{align*}
    Now assume that $A \in V_k^+ \mid (T'-1)$. Since $V^+_k \mid S = V^+_k$ and $T^{-1}S = S T'$ we get 
    \begin{align*}
    V_k^+ \mid (T'-1) = V_k^+ \mid S(T'-1) = V_k^+ \mid (T^{-1} - 1)S
    \end{align*}
    Hence $A \in  V^+_k \mid (T^{-1} - 1)S$. Since $V^+_k  = V_k \mid (1+\epsilon)$ and $\Delta^\ast = (1+\epsilon)(T^{-1} - 1)$ we obtain $A\mid S \in V_k \mid \Delta^\ast$ if and only if $A \in V_k^+ \mid (T'-1)$. More precisely if $A= H \mid (T'-1)$ for some $H \in V_k^+$ we therefore have  by using $H = \frac{1}{2} H \! \mid (1+\epsilon)$ and $(1+\epsilon) S = \delta (1+\epsilon)$ 
    \begin{align*}
    A\! \mid \! S
    &=H\! \mid \! (T'-1) S
    =H\! \mid \! S(T^{-1}-1)\\
    &=\frac{1}{2} H\! \mid \! (1+\epsilon) S(T^{-1}-1)
    =\frac{1}{2} H\! \mid \! \delta(1+\epsilon)(T^{-1}-1)\\
    &=\frac{1}{2} H \!\mid \!\delta \Delta^* \,,
    \end{align*}
    which gives by the $\Gamma$-invariance of the pairing (Proposition \ref{prop:pairing})  and $\gz_k \!\mid\!\! \Delta = \gr_k = Z_k \frac{X^{k-1}-Y^{k-1}}{X-Y}$ 
    \begin{align*}
    \sum_{\substack{k_1+k_2=k\\ k_1,k_2\geq 1}} a_{k_1,k_2}  Z_{k_1,k_2} = \langle A\!\!\mid \! S , \gz_k \rangle  =\frac{1}{2} \langle   H \! \mid \!\delta \Delta^* , \gz_k \rangle = \frac{1}{2} \langle   H \! \mid \!\delta , \gz_k\!\mid\! \Delta \rangle  = \frac{1}{2} \langle   H \! \mid \!\delta , \gr_k \rangle  = \lambda Z_k\,.
    \end{align*}
    The second expression for $\lambda$ follows by using the beta integral 
    \begin{align*}
    \int_0^1 t^{k_1-1} (1-t)^{k_2-1} dt = \left( \binom{k-2}{k_1-1}   (k-1) \right)^{-1}\,.
    \end{align*}
    The proof for (ii) follows by considering a symmetric version of (i) together with $P_{k_1,k_2} = Z_{k_1,k_2} +Z_{k_2,k_1} + Z_k$ (Exercise \ref{ex15}). For this one shows that $A = H\mid (T'-1)$ for $H \in V_k^+$ is symmetric if and only if $H=H\mid U$.
\end{proof}

\begin{remark}
    Notice that there is a one-to-one correspondence between elements $H \in V_k^+$ and the relations in Lemma \ref{lem:Arealindk} (i), since $\dim V_k^+$ gives exactly the proven number of relations in $\dz_k$, which follow from the dimension formulas in Theorem \ref{thm:dzoddbasis} and \ref{thm:dzevenbasis}. In other words, for a given relation in $\dz_k$, there exists a unique symmetric polynomial $H \in V_k^+$ that generates this relation by the construction in Lemma \ref{lem:Arealindk}. 
    This can also be proven by showing that the pairing is non-degenerate, which provides another proof of the dimension formulas (without giving explicit bases).
\end{remark}

We now have all the necessary ingredients to prove Theorem \ref{thm:wkminuspkev}.

\begin{proof}[Proof of Theorem \ref{thm:wkminuspkev}]
    First we will show that the map 
    \begin{align*}
    \Theta : W_k^- &\longrightarrow \mathcal{P}^\text{ev}_k\\
    p &\longmapsto \sum_{\substack{k_1+k_2=k\\k_1, k_2 \text{ odd}}}\beta^p_{k_1,k_2}  Z_{k_1,k_2}\,
    \end{align*}
    is well defined, i.e. that the image is actually an element in $\mathcal{P}^\text{ev}_k$. For this we write $q= p\mid T$ and then see by direct calculation that with $A=q^{\text{od}} - 3 q^{\text{ev}}$ and $H = 2 (q^{\text{ev},+} - q^{\text{od},+}) \in V_k^+$ (Exercise \ref{ex15} (i)) we have 
    \begin{align*}
    A = H \mid (T'-1)\,.
    \end{align*}
    Here one uses that $p$ satisfies the Lewis equation and is antisymmetric to show that $q^{\text{ev},-} = \frac{1}{2}p$ and $q^{\text{od},-}=0$. Now we can apply Lemma \ref{lem:Arealindk} (i) and, since  the coefficients $\beta^p_{k_1,k_2}$ are (up to a binomial coefficient) defined by the coefficients of $q$, we  obtain the relation
    \begin{align*}
    \sum_{\substack{k_1+k_2=k\\k_1, k_2 \text{ odd}}}\beta^p_{k_1,k_2}  Z_{k_1,k_2} \equiv     \frac{1}{3} \sum_{\substack{k_1+k_2=k\\k_1, k_2 \text{ even}}}\beta^p_{k_1,k_2}  Z_{k_1,k_2}  \quad  \mod  \Q Z_k\,.
    \end{align*} 
    Further since $q^{\text{od},-}=0$ we have $\beta^p_{k_1,k_2} = \beta^p_{k_2,k_1}$ for even $k_1,k_2 \geq 2$. From this we obtain   
    \begin{align*}
    \sum_{\substack{k_1+k_2=k\\k_1, k_2 \text{ odd}}}\beta^p_{k_1,k_2}  Z_{k_1,k_2} \equiv     \frac{1}{6} \sum_{\substack{k_1+k_2=k\\k_1, k_2 \text{ even}}}\beta^p_{k_1,k_2}  P_{k_1,k_2}  \quad  \mod  \Q Z_k\,,
    \end{align*}
    i.e. the map $\Theta$ is well-defined and the image is given as stated in the theorem. The coefficient of $Z_k$ can be obtained explicitly by calculating the corresponding $\lambda$ in Lemma \ref{lem:Arealindk}. It remains to show that $\Theta$ is bijective. We will show in the proof of Theorem \ref{thm:wkpluspkev} that the dimension of $W_k^-$ and $\mathcal{P}^\text{ev}_k$ coincide and therefore we just need to show injectivity. This follows directly since if $\Theta(p) =0$, then $\beta^p_{k_1,k_2} = 0$ for all $k_1,k_2 \geq 1$ odd, since the $Z_{\text{od}, \text{od}}$ are linearly independent. But then we have $p=0$ since $0 = q^{\text{ev},-} = \frac{1}{2}p$. 
\end{proof}

Using the second part of Lemma \ref{lem:Arealindk}, we will now give the proof of Theorem \ref{thm:wkpluspkev}. For this, we will first show the following. 

\begin{lemma} \label{lem:Aorthlemma} We have $A \in  V_k^{U,+}  \mid (1-S)$ if and only if $A \in V_k^{\text{od}}\cap V_k^+$ and $A \perp W_k^+$.
\end{lemma}
\begin{proof}
    Let $v \in  (V_k^{U,+} \mid (1-S))^{\perp}$. Since  $V_k^{U,+} \mid (1-S) = V_k \! \mid (1+U+U^2)(1+\epsilon)(1-S)$ we obtain, by using again the $\Gamma$-invariance and the fact that the pairing is non-degenerate, that 
    \begin{align*}
    v \! \mid(1-S)(1+\epsilon)(1+U+U^2) = 0\,.
    \end{align*}
    But this means that 
    \begin{align*}
    v \! \mid(1-S)(1+\epsilon) \in \ker(1+U+U^2) \cap \ker(1+S) \cap \ker(1-\epsilon) = W_k^+\,,
    \end{align*}
    from which we get $v \in W_k^+ + V_k^- + V_k^{\text{ev}}$. The argument is reversible. Since $\left(V_k^-\right)^\perp=V_k^+$ and $\left(V_k^{\text{ev}}\right)^\perp=V_k^\text{od}$, we therefore have
    \begin{align*}
    V_k^{U,+} \mid (1-S) = \left(W_k^+ +V_k^-+V_k^{\text{ev}}\right)^\perp
    =\left(W_k^+\right)^\perp\cap V_k^+\cap V_k^\text{od}\,,
    \end{align*}
    from which the statement follows.
\end{proof}

\begin{proof}[Proof of Theorem \ref{thm:wkpluspkev}]
    Since $P_{k_1,k_2}=P_{k_2,k_1}$, an assignment $P_{k_1,k_2}\mapsto p_{k_1,k_2}$ must satisfy $p_{k_1,k_2}=p_{k_2,k_1}$. Given such numbers for $k_1,k_2\geq 2$ even with $k_1+k_2=k$, write
    \begin{align*}
    p=\sum_{\substack{k_1+k_2=k\\ k_1,k_2 \geq 2\text{ even}}} p_{k_1,k_2} X^{k_1-1} Y^{k_2-1}\in V_k^{\text{od}}\cap V_k^+\,.
    \end{align*}
    Every relation among the $P_{k_1,k_2}$ can also be written with symmetric coefficients. By Lemma \ref{lem:Arealindk} (ii), these numbers define a linear map $\varphi$ on $\mathcal{P}^{\mathrm{ev}}_k$ with $\varphi(P_{k_1,k_2})=p_{k_1,k_2}$ and $\varphi(Z_k)=0$ if and only if
    \begin{align*}
    \sum_{\substack{k_1+k_2=k\\ k_1,k_2\geq 2\text{ even}}} a_{k_1,k_2}  p_{k_1,k_2} = 0\,,
    \end{align*}
    for every $A\in V_k^{U,+}\mid(1-S)$, where the $a_{k_1,k_2}$ are the coefficients in \eqref{eq:polA}. Every such $A$ satisfies $A\mid S=-A$, since $S^2=\sigma$ acts trivially on $V_k$ for even $k$. By \eqref{eq:pairingaandc}, the last condition is therefore equivalent to $\langle A,p\rangle=0$ for every $A\in V_k^{U,+}\mid(1-S)$. Lemma \ref{lem:Aorthlemma} shows that this is equivalent to $p\in W_k^+$, which proves the theorem.
\end{proof}

\subsection{Extended period polynomials}
In this section we want to introduce period polynomials for all modular forms. For a modular form $f = \sum_{n\geq 0} a_n q^n \in M_k$, which is not a cusp form, i.e. $a_0 \neq 0$, the integral \eqref{eq:defperiodpol} does not converge. In \cite{Za1} Zagier introduces the \emph{extended period polynomial} for any $f=\sum_{n\geq 0} a_n q^n \in M_k$ by
\begin{align}\begin{split}
\label{eq:defextendedperiodpol}
\widehat{P}_f(X,Y) = &\int_{\tau_0}^{i \infty } (X- Y \tau)^{k-2} (f(\tau) - a_0) d\tau + \int_0^{\tau_0} (X-Y \tau)^{k-2} \left( f(\tau)  - \frac{a_0}{\tau^k}\right) d\tau \\
&+ \frac{a_0}{(k-1)}\left( \frac{1}{Y} - \frac{\tau_0^{1-k}}{X}\right) (X - Y \tau_0)^{k-1}\,.
\end{split}
\end{align}
Here $\tau_0 \in \Ha$ is arbitrary, and one can check that the definition of $\widehat{P}_f(X,Y)$ is independent of $\tau_0$ since the derivative of the right-hand side with respect to $\tau_0$ vanishes.

The Eisenstein series denoted by $G_k$ in \cite[Section~6]{GKZ} is our $\aG_k$. Taking into account the normalization $G_k=(-2\pi i)^{-k}\aG_k$ and our convention in \eqref{eq:defperiodpol}, a direct calculation from \eqref{eq:defextendedperiodpol} using \cite[Equations~(2) and~(9)]{Za1} gives

\begin{align*}
\widehat{P}_{G_k}(X,Y) = -\frac{\zeta(k-1)}{2(k-1)(2\pi i)^{k-1}} (X^{k-2} - Y^{k-2}) - \frac{1}{2(k-1)} \sum_{\substack{k_1+k_2 = k \\k_1,k_2 \geq 0}} \frac{B_{k_1}}{k_1!} \frac{B_{k_2}}{k_2!} X^{k_1-1}  Y^{k_2-1}\,.
\end{align*}
Notice that the first sign is opposite to the one printed in \cite[Section~6]{GKZ}. Substituting the Fourier expansion of $G_k$ into \eqref{eq:defextendedperiodpol} gives the sign above.

In general $\widehat{P}_f(X,Y)$ is not an element in $\C \otimes V_k$ anymore, since we can get poles in $X$ and $Y$. We therefore define the space
\begin{align*}
\widehat{V}_k = \bigoplus_{\substack{k_1+k_2=k \\ k_1,k_2\geq 0}} \Q X^{k_1-1} Y^{k_2-1}\,.
\end{align*}
With this, $\widehat{P}_f(X,Y) \in \C \otimes \widehat{V}_k$ for any $f \in \mf_k$. The group $\Gamma$ does not act on $\widehat{V}_k$ anymore, but still it makes sense to define the following subspace of $\widehat{V}_k$
\begin{align*}
\widehat{W}_k := \ker(1+U+U^2) \cap \ker(1+S) = \ker(1-T-T')\,,
\end{align*}
which contains all elements in $\widehat{V}_k$ which vanish under the actions $1+U+U^2$ and $1+S$ (defined in the same way as before). One can then also check that $\widehat{P}_f(X,Y) \in \C \otimes \widehat{W}_k$ and we have the following extended version of the Eichler-Shimura theorem, where we define $\widehat{W}_k^\pm$ again by the symmetric and antisymmetric parts of $\widehat{W}_k$.

\begin{theorem}[Eichler-Shimura, Zagier \cite{Za1}]
The map $\widehat{p}^\pm_f: f \mapsto \widehat{P}^\pm_f$ induces isomorphisms
\begin{align*}
\widehat{p}^+_f : \mf_k \xrightarrow{\raisebox{-0.25 em}{\smash{\ensuremath{\sim}}}} \C \otimes \widehat{W}^+_k\,,\qquad  \qquad
\widehat{p}^-_f: \mf_k \xrightarrow{\raisebox{-0.25 em}{\smash{\ensuremath{\sim}}}} \C \otimes \widehat{W}^-_k = \C \otimes W^-_k \,.
\end{align*}
\end{theorem}

We obtain an extended version of Theorem \ref{thm:wkpluspkev} given by the following.
\begin{theorem}[{\cite[Theorem~4]{GKZ}}] \label{thm:wkpluspkevextended}  For even $k\geq 4$ the following map is an isomorphism of $\Q$-vector spaces
\begin{align*}
\widehat{W}^+_k &\longrightarrow \left\{ \varphi \in \Hom( \mathcal{P}^\text{ev}_k , \Q) \right\}\\
\lambda \left( X^{k-1} Y^{-1} + X^{-1} Y^{k-1} \right) + \sum_{\substack{k_1+k_2=k\\ k_1,k_2 \geq 2\text{ even}}} p_{k_1,k_2} X^{k_1-1} Y^{k_2-1}  &\longmapsto ( \varphi: P_{k_1,k_2} \mapsto p_{k_1,k_2}, \, Z_k \mapsto -2\lambda)\,.
\end{align*}
\end{theorem}
\begin{proof}
Corresponding to the splitting $\mf_k = S_k \oplus \C G_k$ we also have $\widehat{W}^+_k = W_k^+ \oplus \Q \mathcal{E}_k$, where 
\begin{align*}
\mathcal{E}_k = \frac{1}{4\,k!}\sum_{\substack{k_1+k_2 = k \\k_1,k_2 \geq 0}} \binom{k}{k_1} B_{k_1}B_{k_2} X^{k_1-1}  Y^{k_2-1} = \sum_{\substack{k_1+k_2 = k \\k_1,k_2 \geq 0}} \beta(k_1)\beta(k_2) X^{k_1-1}  Y^{k_2-1}\,.
\end{align*}
The theorem follows then by Theorem \ref{thm:wkpluspkev} together with the fact that the image of $\mathcal{E}_k $ gives exactly the realization $\varphi_\beta$ given in Theorem \ref{thm:betarealization}, which satisfies $\varphi_\beta(Z_k) \neq 0$. 
\end{proof}

There is also a (non-unique) way to lift the realizations in Theorem~\ref{thm:wkpluspkevextended} of $\mathcal{P}^\text{ev}_k$ to realizations of $\dz_k$ in $\Q$.

\begin{theorem}[{\cite[Supplement to Proposition~5]{GKZ}}] \label{thm:prop5supp} Let $\widetilde{P} = P + \lambda (X^{k-1} Y^{-1} + X^{-1} Y^{k-1}) \in \widehat{W}_k^+$ where $P\in V_k$ and $\lambda \in \Q$, and define
\begin{align*}
Z =\frac{1}{3} P \! \mid (T^{-1}+1) + \frac{\lambda}{6} \frac{X^{k-1}-Y^{k-1}}{X-Y} \mid (5-3U + U\epsilon).
\end{align*}
Then we have
\begin{align*}
P &= Z \mid (1+\epsilon) - 2\lambda  \frac{X^{k-1}-Y^{k-1}}{X-Y} = Z \mid T(1+\epsilon)\,.
\end{align*}
In particular we get  a realization $\varphi$ of $\dz_k$ in $\Q$ by $\varphi(Z_k) = -2 \lambda$ and 
\begin{align*}
\varphi(Z_{k_1,k_2}) &= \text{coefficient of }    X^{k_1-1} Y^{k_2-1} \text{ in } Z(X,Y)\,,\\
\varphi(P_{k_1,k_2}) &= \text{coefficient of }    X^{k_1-1} Y^{k_2-1} \text{ in } P(X,Y)
\end{align*}
\end{theorem}
\begin{proof}
This can be checked by direct calculation.
\end{proof}
The realization $\varphi_\beta$ (Theorem \ref{thm:betarealization}) can be seen as a special case of Theorem \ref{thm:prop5supp} by choosing $\widetilde{P}  = \mathcal{E}_k$.

Further modular phenomena for regularized double zeta values, including non-admissible indices and period polynomials for a subgroup of index two, are studied by Hirose in \cite{Hi2}.

\input{chap_FiniteMZVModularForms}

In the next chapter, we return to multiple Eisenstein series in arbitrary
depth. Chapter~\ref{sec:formalspaces} will then lift the present depth-two
space by adjoining the lower indices which encode derivatives.

\vspace{1cm}
\begin{center}
\bf{ \Large {\color{qmzvlinecol}\ding{94}}~ Exercises ~{\color{qmzvlinecol}\ding{94}}}\end{center}

The following is a collection of exercises intended to deepen the reader's understanding of this chapter.

\begin{exer}\label{ex12}\begin{enumerate}[\textup{(}i\textup{)}] 
    \item Prove for $m\geq 1$ Euler's formula \eqref{eq:phieuler}
\begin{align*}
\mathsf{Z}(2m) = -\frac{B_{2m} }{2 (2m)!}  \left( -24 \mathsf{Z}(2) \right)^m\,,
\end{align*}
by assuming the condition in (i) of Theorem \ref{thm:dzeuler}.
\item[(ii*)] Assume the conditions in (ii) of Theorem \ref{thm:dzeuler} hold. Try to find an explicit formula for $\mathsf{Z}(2m)$ as a polynomial in $\partial^j \,\mathsf{Z}(2)$ for $j\geq 0$, which generalizes Euler's formula.

(Part (ii*) is a (hard) bonus exercise.)
\end{enumerate}
\end{exer}

\begin{exer}\label{ex13}
Show that the pairing $\langle \cdot , \cdot \rangle$ defined in \eqref{eq:pairingdef} is $\Gamma$-invariant, bilinear, symmetric and non-degenerate, i.e. prove Proposition \ref{prop:pairing}.
\end{exer}

\begin{exer}\label{ex14}\begin{enumerate}[\textup{(}i\textup{)}] 
    \item Show that for a cusp form $f\in S_k$ and $\gamma \in \SLZ$ we have
    \begin{align*}
    (P_f \mid\! \gamma)(X,Y) = \int_{\gamma^{-1}(0)}^{\gamma^{-1}(i \infty)} (X-Y\tau)^{k-2} f(\tau) \,d\tau\,,
    \end{align*}
    where $\gamma(z) = \frac{a z + b}{c z + d}$ for $\gamma = \sabcd$.
    \item Use (i) to show that for all $f \in S_k$ we have 
    \begin{align*}
    P_f \mid\! (1+S)  = P_f \mid\! (1+U+U^2) = 0 \,.	
    \end{align*}
    \item Prove Lemma \ref{lem:lewis}. 
\end{enumerate}
\end{exer}

\begin{exer}\label{ex15}
\begin{enumerate}[\textup{(}i\textup{)}] 
    \item Let $p\in W_k^-$ and define $q = p \!\mid\! T$. Show that $A = H \!\mid\! (T'-1)$, where  $A=q^{\text{od}} - 3 q^{\text{ev}}$ and $H = 2 (q^{\text{ev},+} - q^{\text{od},+} ) \in V_k^+$.  
            \item 
    Assume that $A\in V_k^+$ and give the proof of Lemma \ref{lem:Arealindk} (ii) by using (i), i.e. show that in $\dz_k$ the relation
    \begin{align*}
    \sum_{\substack{k_1+k_2=k\\ k_1,k_2\geq 1}} a_{k_1,k_2}  P_{k_1,k_2} = \mu Z_k\,,
    \end{align*}
    holds for some $\mu\in\Q$ if and only if $A = H \mid (1-S)$ for some  polynomial $H\in V^U_k \cap V_k^+$. In this case we have
    \begin{align*}
    \mu = \langle H , \frac{X^{k-1}-Y^{k-1}}{X-Y}\rangle \,.
    \end{align*}		
    Here $ V^U_k $ denotes the space of $U$ invariant polynomials, i.e. $H\!\mid\! U = H$. 
\end{enumerate}
\end{exer}

\begin{exer}\label{ex:finiteoperatorslecture}
Use \eqref{eq:finiteGQlecture} to prove the following properties of the maps $\mathcal K$ and $Q$ from \Cref{sec:fmzvmodularforms} on even polynomials.
\begin{enumerate}[(i)]
\item The map $\mathcal N^{-2}Q$ is a projection onto the even symmetric polynomials, its kernel is the graph $H^+=-\mathcal N^{-2}Q(H^-)$, and we have $Q^2=\mathcal N^2Q$.
\item If $Q(H)=0$, then
\begin{align*}
 H^+(0,Y)=\frac{Y(\partial_XH^-)(Y,Y)-H^-(0,Y)}{k-1}\,.
\end{align*}
\item The conditions $\mathcal K(H)|(1+U+U^2)=0$ and $Q(H)=0$ separately define spaces of dimensions $k/2-\lfloor(k-2)/6\rfloor$ and $\lfloor k/4\rfloor$, and their intersection has dimension $\dim_\C\ms_k+1$.
\end{enumerate}
Notice that it is still open whether $H^-\mapsto-\mathcal N^{-2}Q(H^-)$ is injective on the even antisymmetric polynomials.
\end{exer}

\begin{exer}\label{ex:finitedeformationlecture}
Prove the statement about the deformation $\mathcal K_t$ after \cref{ex:finiteH12lecture}. For this, notice first that on even antisymmetric polynomials the even part of $\mathcal K_t$ is $(k-2)\mathcal N$, so its restriction to this space is injective for every $t$. The explicit formula defines $\mathcal K_t$ on all of $V_k$. Using the polynomials $\Omega_j$ in \eqref{eq:finiteEblecture}, show that the even parts of the $\mathcal K_t(\Omega_j|\delta)$ are independent, where their lowest terms have the coefficients
\begin{align*}
 -2t(2j+1)(2j+2)(k-1-2j)&\quad\text{at }X^{2j}Y^{k-2-2j}\quad(t\neq0),\\
 (k-2)(2j+3)(k-3-2j)^2&\quad\text{at }X^{2j+2}Y^{k-4-2j}\quad(t=0)\,.
\end{align*}
Deduce that the cyclic condition has rank $\lfloor(k-2)/6\rfloor$, and use \eqref{eq:finiteperiodspacelecture} to identify the image for $t=0$ with $W_k^-$.
\end{exer}

\begin{exer}\label{ex:finiteperiodtriplelecture}
Let $f$ be the finite polynomial in \eqref{eq:finitepolynomialslecture}. Prove that for every $P\in W_k^-\oplus W_k^+$ we have
\begin{align*}
 P(\partial_X+\partial_Y,\partial_Y)f(X,Y)=0\,.
\end{align*}
For odd $P$ this already follows from the double shuffle equations, as in \cite[Proposition~2]{GKZ}. Use $P=X^2Y^2(X^2-Y^2)^3$ in weight $12$ to obtain
\begin{align*}
 14\za(3,1,8)+14\za(4,1,7)+5\za(5,1,6)+2\za(4,1)\za(6,1)=0\,.
\end{align*}
Notice that these are identities among the values $\za(n,1,k-1-n)$ and products, i.e. they do not say that these values span the depth-four space.
\end{exer}

\begin{exer}\label{ex:finiterisanoperatorlecture}
Let $\mathcal V_k$ be the space from \Cref{sec:finitemiddleone}. For $f\in\mathcal V_k$, define the polynomial $u_f$ by the right-hand side of \eqref{eq:finiteulecture}, with $B=f+f|\epsilon$ and $A=\partial_X(XB)$. Show that its coefficients with even $2\leq j\leq k-4$ detect every nonzero $f$ with $B=0$. (Hint: If $B=0$, the double shuffle equations simplify the right-hand side of \eqref{eq:finiteulecture} to $u_f=-(\partial_X+\partial_Y+Y^{-1})f$. For $f=\Omega_j$ from \eqref{eq:finiteEblecture}, the part of $u_f$ which is odd in $X$ has its first nonzero term at $X^{2j+1}$, with coefficient $2-2j(k-3-2j)$.) Deduce that any failure of this coefficient map to be injective must involve the product part $B\neq0$. Notice that this is a statement about the polynomial model, and it does not prove Risan's spanning conjecture for finite values.
\end{exer}

\begin{exer}\label{ex:finiterisancriterionlecture}
Let $q\in\Q[X,Y]_{k-3}$ be odd in $X$, even in $Y$, and satisfy $q(X,0)=0$. Its monomials correspond to Risan's family $\za(1,2a,2b+1)$ with $a,b\geq1$. Put
\begin{align*}
 J_q&=(X+Y)q+\int_0^Yq(X,t)\,dt,\\
 I_q&=\int_0^Xq(t,Y)\,dt-\int_0^Yq(X,t)\,dt,\\
 C_q&=Yq+I_q(X-Y,Y)-I_q(X,Y)\\
 &\quad-\frac12\bigl(J_q(X,Y)+J_q(X,Y-X)-J_q(X+Y,-Y)\bigr)\\
 &\quad-(k-1)X^{k-2}\int_0^1\int_0^1q(x,y)\,dx\,dy\,.
\end{align*}
With $[q,u]=q(\partial_X,\partial_Y)u$, prove that $[q,u_f]=0$ for all $f\in\mathcal V_k$ if and only if
\begin{align*}
 J_q^-|\delta(1+U+U^2)&=0,\\
 C_q|(1+\epsilon)(1+\delta)&=(C_q(1,0)+C_q(0,1))\bigl((X+Y)^{k-2}+(X-Y)^{k-2}\bigr)\,.
\end{align*}
(Hint: Decompose $f=f_{\mathfrak B}+F$ as in \eqref{eq:finitecanonicalliftlecture}. Moving the derivatives and substitutions in \eqref{eq:finiteulecture} to the left gives $[q,u_f]=-[J_q,F]+[C_q,B]$, and the two conditions annihilate the part $F$ and the even symmetric $B$ with $B(1,1)=0$.) In particular, these conditions give finite relations, and they describe all coefficient relations in this rational polynomial model. The necessity for formal finite relations also follows from \cref{thm:formalfinitebetaisomorphism}.
\end{exer}

%% file: chap_FiniteMZVModularForms.tex
\section{Finite MZVs and modular forms}\label{sec:fmzvmodularforms}
In this section, we want to come back to the relations among finite multiple zeta values which we saw in \Cref{subsec:modularfmzv}. The first one was a relation in weight $12$, i.e. it corresponds to the cusp form $\Delta$ just as in the case of double zeta values, but the finite multiple zeta values appearing there have depth four. We will give a direct construction of these relations, and we will see that the formal double zeta space and the period polynomials from the previous sections show up again, but now with a shift in depth. The formal explanation of this phenomenon will be given later in \Cref{sec:formalfinitemzv}.

Throughout this section, $k\geq4$ is even. Recall that $\za(a)=0$ for all $a\geq1$ (\cref{prop:fmzvdepth1}) and that by \cref{prop:fmzvdep2} we have
\begin{align*}
 \za(a,b)=(-1)^a\binom{a+b}{a}Z(a+b),\qquad Z(j)=\frac{\za(j-1,1)}j\quad(j\geq2)\,,
\end{align*}
where $Z(j)$ is the element defined in \eqref{eq:defZ}. We also set $Z(1)=0$. Since $Z(2j)=0$, all double values of even weight vanish, and we will use this several times below. Notice that there are several families of finite multiple zeta values with an entry $1$ in this section, and one has to be a bit careful to keep them apart. The values $\za(2a,1,2b,1)$ give the depth-four space in \cref{thm:fmzvmodularforms}. On the other hand, the values $\za(1,j,k-1-j)$ studied by Kaneko and Zagier and Risan's families in \cref{conj:fmzvmodularforms}\,(iii) have different index ranges, and the values $\za(n,1,k-1-n)$ with a $1$ in the middle behave differently again, as we will see in \Cref{sec:finitemiddleone}.

\subsection{Some identities}
We first need some identities among finite multiple zeta values. The first one is a finite version of Hoffman's relation, which follows directly from the linear shuffle relations.
\begin{lemma}\label[lemma]{lem:finitehoffmanlecture}
For positive $a_1,\ldots,a_r$, we have
\begin{align*}
 \sum_{j=1}^r\za(a_1,\ldots,a_j+1,\ldots,a_r)
 ={}&\za(1,a_1,\ldots,a_r)\\
 &+\sum_{j=1}^r\sum_{i=2}^{a_j}\za(a_1,\ldots,a_{j-1},i,a_j+1-i,a_{j+1},\ldots,a_r)\,.
\end{align*}
\end{lemma}
\begin{proof}
Put $w=z_{a_1}\cdots z_{a_r}$. By \cref{thm:linearshufflefmzv} we have $\za(z_1\shuffle w)=-\za(z_1w)$, while $\za(z_1*w)=0$. If we subtract the two expansions, the ordinary insertions of $z_1$ cancel. The remaining terms of the shuffle product split an $x$-block into $(i,a_j+1-i)$ with $2\leq i\leq a_j$, and the remaining terms of the stuffle product increase $a_j$ by one, which gives the formula.
\end{proof}
\noindent We also need the following two identities, which follow from the antipode of the stuffle algebra together with the reversal formula:
\begin{align}\label{eq:finiteantipodeslecture}
 2\za(a,b,c)&=-\za(a+b,c)-\za(a,b+c)\quad(a+b+c\text{ odd}),\notag\\
 2\za(a,b,c,d)&=\za(a,b)\za(c,d)-\za(a+b,c,d)\notag\\
 &\qquad-\za(a,b+c,d)-\za(a,b,c+d)\quad(a+b+c+d\text{ even})\,.
\end{align}
To see this, one can expand $\sum_{i=0}^r(-1)^i\za(a_i,\ldots,a_1)\za^\star(a_{i+1},\ldots,a_r)=0$ for $r=3,4$. Here all terms with a factor of depth one vanish, the double values of even weight vanish, and the reversal formula identifies the two remaining terms of depth $r$.

With these identities we can reduce every triple value of even weight to the values with a $1$ in the middle and products of the $Z(j)$.
\begin{lemma}\label[lemma]{lem:finitetriplereductionlecture}
For positive $a,b,c$ with $a+b+c=k$, we have
\begin{align}\label{eq:finitetriplereductionlecture}
 \za(a,b,c)={}&(-1)^{b-1}\sum_{n+m=k-1}\binom na\binom mc\za(n,1,m)\notag\\
 &+(-1)^{a+c}\sum_{i=1}^{b-1}\binom{a+i}a\binom{b+c-i}c Z(a+i)Z(b+c-i)\,,
\end{align}
where $n,m$ in the first sum are positive.
\end{lemma}
\begin{proof}
We write $\Sigma(t;c)=\sum_{u+v=t}\za(u,v,c)$. If we apply \cref{lem:finitehoffmanlecture} to a double value of weight $k-1$ and use the reversal formula and the vanishing of the double values of weight $k$, we get $\Sigma(n+1;m)+\Sigma(m+1;n)=\za(n,1,m)$. Now let $a+b+c=k-1$ and apply the same lemma to $(a,b,c)$ together with \eqref{eq:finiteantipodeslecture}. The terms containing $\Sigma$ cancel by the identity we just obtained, and the complete double sums vanish, so we are left with
\begin{align}\label{eq:finitetriplerecurrencelecture}
 &(a+1)\za(a+1,b,c)+b\za(a,b+1,c)+(c+1)\za(a,b,c+1)\notag\\
 &\qquad=\sum_{i=2}^{b}\za(a,i)\za(b+1-i,c)-\za(a,b)\za(1,c)\,.
\end{align}
By the depth-two formula and $(a+1)\binom{u}{a+1}=(u-a)\binom ua$, the product sum in \eqref{eq:finitetriplereductionlecture} satisfies this recurrence, and after subtracting it the recurrence becomes homogeneous. The first sum satisfies the homogeneous recurrence by the same binomial identity in $a$ and $c$, and for $b=1$ it agrees with $\za(a,1,c)$. Therefore the statement follows by induction on the middle entry $b$.
\end{proof}

We now want to put all these identities into one polynomial, which we will use for the relations in depth three and in depth four. For this, we define
\begin{align}\label{eq:finitepolynomialslecture}
 M(X,Y)&=\sum_{n+m=k-1}\za(n,1,m)X^{n-1}Y^{m-1},\notag\\
 \mathfrak B(X,Y)&=\sum_{r+s=k}Z(r)Z(s)X^{r-1}Y^{s-1},\notag\\
 B(X,Y)&=\mathfrak B(X,Y)-\frac{\mathfrak B(1,1)}2(X^{k-2}+Y^{k-2}),\notag\\
 f(X,Y)&=\frac{Y-X}{X}\bigl(YM(Y-X,Y)-B(Y-X,Y)\bigr)\,,
\end{align}
where all sums run over positive indices. Although the last expression has a denominator, it is a polynomial, as we will see in the following lemma.
\begin{lemma}\label[lemma]{lem:finitedspolynomiallecture}
The polynomial $f$ is homogeneous of degree $k-2$, satisfies $f(0,Y)=0$, and
\begin{align}\label{eq:finitedspolynomiallecture}
 f(X,Y)+f(Y,X)=f(X+Y,Y)+f(X+Y,X)=B(X,Y)\,.
\end{align}
\end{lemma}
\begin{proof}
By the reversal formula, $M$ is symmetric. The stuffle product with a value of depth one gives
\begin{align*}
 \sum_{a+b+c=k}\za(a,b,c)=0\,,
\end{align*}
and summing \cref{lem:finitehoffmanlecture} over all double values then gives $M(1,1)=0$. By construction, $B$ is symmetric, even in both variables and vanishes at $(1,1)$. Notice that a symmetric polynomial which vanishes on the diagonal has both first partial derivatives equal to zero there. Therefore $YM(Y-X,Y)-B(Y-X,Y)$ is divisible by $X^2$, which proves the first assertion, and the second equality in \eqref{eq:finitedspolynomiallecture} follows by substitution.

For the first equality, we multiply \eqref{eq:finitetriplereductionlecture} by $X^aY^bW^c$ and sum over all $a+b+c=k$. This gives
\begin{align*}
 XYW\sum_{a+b+c=k}\za(a,b,c)X^{a-1}Y^{b-1}W^{c-1}
 ={}&Y\sum_{n+m=k-1}\za(n,1,m)\bigl((X-Y)^n-(-Y)^n\bigr)\\
 &\qquad\cdot\bigl((W-Y)^m-(-Y)^m\bigr)\\
 &+\sum_{r+s=k}Z(r)Z(s)\bigl((Y-X)^r-Y^r+X^r\bigr)\\
 &\qquad\cdot\bigl((Y-W)^s-Y^s+W^s\bigr)\,.
\end{align*}
If we write $\mathcal Z(X,Y,W)=\sum_{a+b+c=k}\za(a,b,c)X^{a-1}Y^{b-1}W^{c-1}$ for the generating polynomial of the triple values, the stuffle product and the reversal formula give $\mathcal Z(X,Y,X)+2\mathcal Z(Y,X,X)=0$. We insert the displayed expression into this identity multiplied by $X^2Y$, and expanding the products gives $2XY(f(X,Y)+f(Y,X)-B(X,Y))=0$, as required.
\end{proof}
In the language of \Cref{sec:realizations}, \cref{lem:finitedspolynomiallecture} says that the coefficients of $f$ and $B$ give a realization of $\dz_k$ in $\ma$, which sends $Z_{r,s}$ to the coefficient of $X^{r-1}Y^{s-1}$ in $f$, $P_{r,s}$ to the coefficient of $X^{r-1}Y^{s-1}$ in $B$, and $Z_k$ to zero (see \cref{lem:realina}). So the formal double zeta space really shows up for finite triple values. Up to a boundary correction, the coefficients of $f$ are the finite double-zeta coordinates constructed by Kina and Zhao (see \cite[Appendix~A]{Ki} and \cite[Proposition~5.6(i)]{Zh4}). The argument above gives the double shuffle equations for them directly from identities among finite multiple zeta values.

\subsection{Depth-four relations}
The following direct construction of the relations in depth four is new and appears here for the first time. For this, we use the polynomial space $V_k$ with the action of $\GL_2(\Z)$ and the pairing \eqref{eq:pairingdef}, which is non-degenerate and invariant under this action (\cref{prop:pairing}). Notice that an even polynomial pairs to zero with every polynomial which is odd in $X$. For $q\in V_k$ we write $q^\pm=\frac12(q\pm q|\epsilon)$ for its symmetric and antisymmetric part, and we denote by $\mathcal N$ the operator $q\mapsto\partial_X\partial_Y(XYq)$, which multiplies $X^{r-1}Y^{s-1}$ by $rs$. Now let $f,\mathfrak B\in V_k$ be polynomials satisfying \eqref{eq:finitedspolynomiallecture}, where $B$ is obtained from $\mathfrak B$ as in \eqref{eq:finitepolynomialslecture} and $\mathfrak B$ is even and symmetric. Then we define
\begin{align}\label{eq:finitedepthfourpolynomiallecture}
 \Phi(f,\mathfrak B)=\frac12\mathcal N^2(\mathfrak B)
 -\Bigl(\partial_X\partial_Y\bigl((XY\mathfrak B)|T+Y^2f+XY(X+Y)\partial_Xf\bigr)\Bigr)^-\,.
\end{align}
We denote by $\mathcal H_k$ the space of even $H\in V_k$ with $\langle H,\Phi(f,\mathfrak B)\rangle=0$ for every such rational pair $(f,\mathfrak B)$.

\begin{theorem}\label{thm:finitedirectperiodrelation}
For every $H\in\mathcal H_k$, we have $\langle H,\mathcal F_k\rangle=0$, where
\begin{align}\label{eq:finiteFlecture}
 \mathcal F_k(X,Y)={}&\sum_{\substack{r+s=k\\r,s\geq2}}rs\za(r-1,1,s-1,1)X^{r-1}Y^{s-1}\notag\\
 &+\frac{k-1}{2}\za(1,k-2,1)(X^{k-2}-Y^{k-2})\,.
\end{align}
\end{theorem}
\begin{proof}
We take the finite polynomials from \eqref{eq:finitepolynomialslecture} and claim that the even part of $\Phi(f,\mathfrak B)$ equals the even part of $\mathcal F_k$. For odd $r,s\geq3$, the antipode identity and the stuffle product give
\begin{align}\label{eq:finitefourbridgelecture}
 2\za(r-1,1,s-1,1)=rsZ(r)Z(s)+\za(1,r-1,s)-\za(1,s-1,r)\,,
\end{align}
and therefore the symmetric part is $\frac12\mathcal N^2(\mathfrak B)$. To compute the antisymmetric part, we write $f_{r,s}$ for the coefficient of $X^{r-1}Y^{s-1}$ in $f$ and set $f_{k,0}=0$. If $c_{k-r}=-c_r$ are antisymmetric coefficients, then \eqref{eq:finitetriplereductionlecture} gives
\begin{align*}
 \sum_r c_r\za(r-1,1,k-r-1,1)
 ={}&-\sum_r c_r\bigl((r-1)f_{r,k-r}+(r+1)f_{r+1,k-r-1}\bigr)\\
 &-\sum_r\left(\sum_a\binom{k-r}{a}c_a\right)Z(r)Z(k-r)\,,
\end{align*}
where $r$ and $a$ first run over the odd integers from $3$ to $k-3$. This identity also holds for $r=1,k-1$, if we read the two depth-four terms there as $-\za(1,k-2,1)/2$ and $\za(1,k-2,1)/2$. Namely, \eqref{eq:finitetriplereductionlecture} gives $\za(1,k-2,1)/2=f_{2,k-2}+(k-1)\mathfrak B(1,1)/2$, while $f_{1,k-1}=0$ and $f_{k-1,1}=-\mathfrak B(1,1)/2$. The full mixed derivative in \eqref{eq:finitedepthfourpolynomiallecture} has the coefficient
\begin{align*}
 rs\left((r-1)f_{r,s}+(r+1)f_{r+1,s-1}+\sum_a\binom{k-a}{r}Z(a)Z(k-a)\right)
\end{align*}
at $X^{r-1}Y^{s-1}$, so its antisymmetric part has the coefficient $-rs(v_r-v_s)/2$, where $v_r=\za(r-1,1,s-1,1)$ in the interior and the two endpoint values are the ones above. This proves the claim, including the boundary terms. Since the equations defining the rational pairs have rational coefficients, the finite pair is an $\ma$-linear combination of rational pairs, and the statement follows from the definition of $\mathcal H_k$ and the evenness of $H$.
\end{proof}
\noindent If we write $H=\sum_{r+s=k}\binom{k-2}{r-1}h_rX^{r-1}Y^{s-1}$, then the relation in \cref{thm:finitedirectperiodrelation} reads
\begin{align}\label{eq:finiteexplicitHrelationlecture}
 \sum_{\substack{r+s=k\\r,s\geq3\text{ odd}}}rsh_r\za(s-1,1,r-1,1)
 +\frac{k-1}{2}(h_1-h_{k-1})\za(1,k-2,1)=0\,.
\end{align}
Only the values $\za(2a,1,2b,1)$ with even entries occur, since $H$ is even. Notice that $H^-(0,Y)=0$ removes the boundary triple $\za(1,k-2,1)$. Also \eqref{eq:finitefourbridgelecture} gives
\begin{align*}
 \za(r-1,1,s-1,1)+\za(s-1,1,r-1,1)=\za(r-1,1)\za(s-1,1)\,,
\end{align*}
i.e. the symmetric part of $H$ only contributes products.

\subsection{The space \texorpdfstring{$\mathcal H_k$}{Hk} and period polynomials}
We now want to understand the space $\mathcal H_k$ better. Every even and symmetric $\mathfrak B$ has the lift
\begin{align}\label{eq:finitecanonicalliftlecture}
 f_{\mathfrak B}(X,Y)=\frac12\bigl(B(X,Y)+B(X-Y,Y)-B(X,X-Y)\bigr)\,,
\end{align}
where $B$ is the boundary-corrected polynomial from \eqref{eq:finitepolynomialslecture}. One checks by substitution that $f_{\mathfrak B}$ satisfies \eqref{eq:finitedspolynomiallecture} and $f_{\mathfrak B}(0,Y)=0$. Any other lift differs from $f_{\mathfrak B}$ by an antisymmetric polynomial $F$ with $F|U=F$, and for such $F$ we have $\Phi(F,0)=-\mathcal L(F)/2$, where
\begin{align*}
 \mathcal L(q)=\partial_X\partial_Y\bigl((X+Y)(\partial_X+\partial_Y)-2\bigr)(XYq)\,.
\end{align*}
This gives two maps on even polynomials. We define $\mathcal K(H)$ to be antisymmetric and $Q(H)$ to be even and symmetric with
\begin{align}\label{eq:finiteadjointslecture}
 \langle\mathcal K(H),F\rangle&=\langle H,\mathcal L(F)\rangle\quad(F\text{ antisymmetric}),\notag\\
 \langle Q(H),\mathfrak B\rangle&=2\langle H,\Phi(f_{\mathfrak B},\mathfrak B)\rangle\quad(\mathfrak B\text{ even and symmetric})\,.
\end{align}
Since the restricted pairings are non-degenerate, these definitions determine the two maps uniquely.

\begin{proposition}\label[proposition]{prop:finiteGQlecture}
An even polynomial $H$ belongs to $\mathcal H_k$ if and only if
\begin{align*}
 \mathcal K(H)|(1+U+U^2)=0,\qquad Q(H)=0\,.
\end{align*}
The two maps are given by
\begin{align}\label{eq:finiteGQlecture}
 \mathcal K(H)={}&\partial_X\partial_Y\bigl((X-Y)(\partial_X-\partial_Y)-2\bigr)(XYH^-),\notag\\
 Q(H)={}&\mathcal N^2(H^+)+\frac12(X\partial_X-Y\partial_Y)\mathcal N(H^-)\notag\\
 &+\frac14\Bigl(\mathcal K(H)|U^2-2\partial_X\partial_Y\bigl((XYH^-)|T\bigr)\Bigr)|(1+\epsilon)(1+\delta)\,.
\end{align}
\end{proposition}
\begin{proof}
By the decomposition of the lifts above, $H$ has to annihilate the lift part and all pairs $(F,0)$ with $F|\epsilon=-F$ and $F|U=F$. The first condition is $Q(H)=0$. On antisymmetric polynomials, $(1+U+U^2)/3$ is a self-adjoint projection onto the $U$-invariant polynomials, and this gives the other condition. To compute the adjoints, we use $[p,q]=p(\partial_X,\partial_Y)q$ and
\begin{align}\label{eq:finiteFischerlecture}
 (k-2)!\langle p,q\rangle=[p|\delta\epsilon,q]\,.
\end{align}
With respect to $[\cdot,\cdot]$, the multiplication by $X$ and $Y$ is adjoint to differentiation, and a substitution is adjoint to the transposed substitution. Moving these operators to the other side in \eqref{eq:finiteadjointslecture} gives the displayed formulas, where the projection $(1+\epsilon)(1+\delta)/4$ takes the even symmetric part.
\end{proof}
In particular, $\mathcal K(H)$ only depends on $H^-$, and since $\mathcal N^2$ is invertible, the condition $Q(H)=0$ determines $H^+$ uniquely. So $\mathcal H_k$ is the graph of a map from a subspace of the even antisymmetric polynomials to the even symmetric polynomials.

Recall from \Cref{sec:periodpolynomials} and \cref{lem:lewis} that the space $W_k^-$ of even period polynomials is given by the polynomials $P\in V_k$ with
\begin{align}\label{eq:finiteperiodspacelecture}
 P(-X,Y)=P(X,Y),\qquad P(X,Y)=P(X+Y,Y)-P(X+Y,X)\,.
\end{align}
By the Eichler-Shimura isomorphism (\cref{thm:eichlershimura}) it has dimension $\dim_\C\ms_k+1$, where the extra direction $X^{k-2}-Y^{k-2}$ corresponds to the Eisenstein series. The following theorem makes the connection with $\mathcal H_k$ precise.
\begin{theorem}\label{thm:finiteHperiodlecture}
There is an isomorphism $W_k^-\longrightarrow\mathcal H_k$ taking the subspace $P(0,Y)=0$ onto the subspace $H^-(0,Y)=0$.
\end{theorem}
\begin{proof}
Put $d=\lfloor(k-2)/6\rfloor$ and let $r^\ast$ be the largest odd integer smaller than $k/2$. We write
\begin{align*}
 P&=\sum_{r+s=k}\binom{k-2}{r-1}p_rX^{r-1}Y^{s-1},\\
 H^-&=\sum_{r+s=k}\binom{k-2}{r-1}c_rX^{r-1}Y^{s-1}\,.
\end{align*}
The isomorphism is specified by
\begin{align}\label{eq:finiteHnormalizationlecture}
 c_1&=p_1,\notag\\
 r(k-r)c_r-\frac{2r-k}{2r^\ast-k}r^\ast(k-r^\ast)c_{r^\ast}
 &=k(k-1)p_{r+2}\quad(2d+1\leq r<r^\ast,\ r\text{ odd})\,.
\end{align}
For $k\geq8$, the antisymmetric $U$-invariant polynomials have the basis
\begin{align}\label{eq:finiteEblecture}
 \Omega_j={}&XY(X-Y)(X+Y)(2X-Y)(X-2Y)\notag\\
 &\cdot(X^2-XY+Y^2)^{(k-8-6j)/2}\bigl(X^2Y^2(X-Y)^2\bigr)^j\quad(0\leq j<d)\,.
\end{align}
To see this, divide an alternating polynomial in $X$, $Y-X$, $-Y$ by the Vandermonde determinant and use symmetric polynomials. If we put $\omega_j(t)=\Omega_j(t,1)$ and $\widetilde\omega_j(t)=((k-1)t^2+(k-2)t+1)\omega_j(t)+(t-t^3)\omega_j'(t)$, their leading terms are
\begin{align*}
 \omega_j(t)=-2t^{2j+1}+(k-3-2j)t^{2j+2}+\cdots,\qquad \widetilde\omega_j(t)=-4(j+1)t^{2j+1}+\cdots\,.
\end{align*}
The equations for $H^-$ are $\sum_{r\leq r^\ast}r(k-r)c_r[t^r]\widetilde\omega_j=0$, where $r$ is odd. We start with $c_1=p_1$, set $c_r=k(k-1)p_{r+2}/(r(k-r))$ in the range of \eqref{eq:finiteHnormalizationlecture} and set the other interior coefficients to zero. The leading terms allow us to adjust $c_{2j+1}$ successively for $j=d-1,\ldots,1$. For the last correction we use
\begin{align}\label{eq:finiteomomentlecture}
 \sum_{\substack{3\leq r\leq r^\ast\\r\text{ odd}}}(2r-k)[t^r]\widetilde\omega_j=-6(k-2)\delta_{j,0}\,,
\end{align}
which follows by pairing the coefficients at $r$ and $k-r$ and evaluating $\widetilde\omega_j'$ at $1$ and $-1$. Adding a multiple of $((2r-k)/(r(k-r)))_r$ for odd $3\leq r\leq r^\ast$, we can satisfy the last equation. By the moment identity this does not change the other equations, and all these correction directions form exactly the kernel of the interior normalization equations. Therefore the solution is unique.

It remains to see that $p_1$ and the $p_{r+2}$ in \eqref{eq:finiteHnormalizationlecture} determine $P$. The period equation implies $P|\epsilon=-P$, and transposing it gives $[P,\Omega_j]=2[P,\Omega_j|T']=2[P,\Omega_j|T'\delta]=-2[P,\Omega_j]$, i.e. $[P,\Omega_j]=0$. Now the leading terms of $\omega_j$ show that $p_1$ and the coefficients $p_{r+2}$ in \eqref{eq:finiteHnormalizationlecture} determine $P$ uniquely. There are $\lfloor k/4\rfloor-\lfloor(k-2)/6\rfloor=\dim_\C\ms_k+1$ such coordinates, and therefore the map is an isomorphism. In the cases $k=4,6$ there is just the coordinate $c_1=p_1$. Finally, $Q(H)=0$ gives the unique symmetric part.
\end{proof}
So we get $\dim_\C\ms_k$ independent coefficient relations without the boundary triple. Notice that we do not lose independence when we remove the boundary: if all interior coefficients of such an $H$ vanish, then $H$ is symmetric, and $Q(H)=\mathcal N^2(H)=0$ forces $H=0$. The remaining Eisenstein relation expresses $\za(1,k-2,1)$ in terms of the depth-four values.

\begin{ex}\label[example]{ex:finiteH12lecture}
In weight $12$, the polynomial
\begin{align*}
 H=X^2Y^2(40X^6+81X^4Y^2+162X^2Y^4-5Y^6)+\frac{54}{11}(X^{10}+Y^{10})
\end{align*}
belongs to $\mathcal H_{12}$ and satisfies $H^-(0,Y)=0$. For this polynomial, \eqref{eq:finiteexplicitHrelationlecture} gives
\begin{align*}
 16\za(2,1,8,1)+9\za(4,1,6,1)+18\za(6,1,4,1)-2\za(8,1,2,1)=0\,,
\end{align*}
which is, up to sign and the reversal formula, exactly the relation in weight $12$ from \Cref{subsec:modularfmzv}. The Eisenstein direction gives instead
\begin{align*}
 45\za(1,10,1)=112\za(2,1,8,1)+48\za(4,1,6,1)+36\za(6,1,4,1)+76\za(8,1,2,1)\,.
\end{align*}
\end{ex}
Here $\dim_\Q\mathcal H_{12}=2$. In general, \cref{thm:fmzvmodularforms} gives
\begin{align*}
 \fild_4(\mza_k)&=\left\langle\za(2a,1,2b,1):2(a+b+1)=k,\ a,b\geq1\right\rangle_\Q,\\
 \dim_\Q\fild_4(\mza_k)&\leq\frac{k}{2}-2-\dim_\C\ms_k\,.
\end{align*}
Equality and the other depth-four placements remain the claims of \cref{conj:fmzvmodularforms}. Of course, the independent coefficient relations above do not show that there are no further relations among the actual finite values. This is the injectivity question for the finite map in \Cref{sec:formalfinitemzv}.

There is also a nice deformation of this picture. For $t\in\Q$ and $H\in V_k$ we put
\begin{align*}
 \mathcal K_t(H)=\partial_X\partial_Y\bigl((k-2)-t(X\partial_Y+Y\partial_X)\bigr)(XYH^-)\,.
\end{align*}
The image under $\mathcal K_t$ of the even $H$ with $\mathcal K_t(H)|(1+U+U^2)=0$ has dimension $\dim_\C\ms_k+1$ for every $t$. For $t=0$ it is $W_k^-$, and for $t=1$ it is $\mathcal K(\mathcal H_k)$, since $\mathcal K_1=\mathcal K$. So $\mathcal K(\mathcal H_k)$ is a deformation of the space of even period polynomials. We leave the proof as \cref{ex:finitedeformationlecture}.

\subsection{The Kaneko--Zagier triple relations}\label{sec:finiteKZtriples}
We now come back to the triple values $\za(1,j,k-1-j)$. Kaneko and Zagier conjectured the following generating-series formula in \cite[\S3(g), (57) and Conjecture~1(2)]{KZ}. The proof below is new and uses the same finite polynomial $f$.
\begin{theorem}\label{thm:finiteKZtripleslecture}
For an integer $\ell\geq k/6$, define $a_n$ by
\begin{align*}
 \frac{(1-t+t^2)^{3\ell-k/2}}{(1-t)^{2\ell}}=\sum_{n\geq0}a_nt^n,\qquad a_n=0\quad(n<0)\,.
\end{align*}
Then $\sum_{j=1}^{k-2}a_{2\ell-j}\za(1,j,k-1-j)=0$.
\end{theorem}
For example, for $k=12$ and $\ell=2$ we obtain
\begin{align*}
 20\za(1,1,10)+10\za(1,2,9)+4\za(1,3,8)+\za(1,4,7)=0\,.
\end{align*}
If $k/6\leq\ell\leq(k-4)/4$, only the indices $1\leq j<k/2-1$ occur. There are exactly $\dim_\C\ms_k$ such integers $\ell$, and their coefficient vectors are independent, since the last nonzero coefficient is $a_0=1$ at $j=2\ell$. So the theorem proves the relation part of Kaneko--Zagier's conjecture, while the following spanning statement is still open.
\begin{conjecture}[Kaneko--Zagier {\cite[Conjecture~1(1)]{KZ}}]\label[conjecture]{conj:finiteKZtriplespanlecture}
For every even $k\geq4$, we have
\begin{align*}
 \left\langle\za(a,b,c):a+b+c=k,\ a,b,c\geq1\right\rangle_\Q
 =\left\langle\za(1,j,k-1-j):1\leq j<\frac{k}{2}-1\right\rangle_\Q\,.
\end{align*}
\end{conjecture}
This is different from Risan's spanning statement in \cref{conj:fmzvmodularforms}\,(iii). For example, his relation in weight $12$ is
\begin{align*}
 1224\za(1,2,9)-171\za(1,4,7)+158\za(1,6,5)-426\za(1,8,3)=0\,.
\end{align*}
In \cref{conj:finiteKZtriplespanlecture}, we have $j=1,2,3,4$ here, while Risan's family uses $j=2,4,6,8$. Notice that neither spanning statement follows just from the number of coefficient relations.

For the proof of \cref{thm:finiteKZtripleslecture}, we put $u(X,Y)=\sum_{j=1}^{k-2}\za(1,j,k-1-j)X^{j-1}Y^{k-j-2}$ and $A=\partial_X(XB)$. The reduction formula and \eqref{eq:finitepolynomialslecture} give
\begin{align}\label{eq:finiteulecture}
 u(X,Y)={}&\frac{A(X,Y-X)-A(X,Y)}X-\frac{f(Y,Y-X)}Y-(\partial_Xf)(Y,Y-X)\notag\\
 &-f(1,0)\frac{(Y-X)^{k-1}-Y^{k-1}+X^{k-1}}{XY}\,.
\end{align}
To see this, set the first variable equal to zero in the triple generating identity in the proof of \cref{lem:finitedspolynomiallecture}, and then differentiate $XYM(X,Y)=(Y-X)f(Y-X,Y)+XB(X,Y)$.
\begin{proof}[Proof of \cref{thm:finiteKZtripleslecture}]
All residues are taken coefficientwise in $\ma$. For a polynomial $v$ we write $\operatorname{Res}_K(v)=\operatorname{Res}_{t=0}K(t)v(t)\,dt$ with $K(t)=(1-t+t^2)^{3\ell-k/2}/(t^{2\ell}(1-t)^{2\ell})$. The only poles of $K$ are $0,1,\infty$, and we have
\begin{align*}
 K(1-t)=K(t),\qquad K(1/t)=t^kK(t),\qquad K(t/(t-1))=(1-t)^kK(t)\,.
\end{align*}
On polynomials of degree at most $k-2$, we write $\mathsf s v(t)=v(1-t)$ and $\mathsf r v(t)=t^{k-2}v(1/t)$, where products denote composition. Changes of variables and the residue theorem give
\begin{align*}
 \operatorname{Res}_K(\mathsf s v)&=\operatorname{Res}_K(v)-\operatorname{Res}_K(\mathsf r v),\qquad \operatorname{Res}_K(\mathsf s\mathsf r\mathsf s v)=-\operatorname{Res}_K(v),\\
 \operatorname{Res}_K((1-\mathsf s\mathsf r\mathsf s)(2+\mathsf s+\mathsf r)v)&=6\operatorname{Res}_K(v)\,.
\end{align*}
Now put $g(t)=f(1,t)$. The double shuffle equations and the evenness of $B$ say that $\widetilde g=(1+\mathsf s)(1-\mathsf r\mathsf s)g=0$, and the chain rule gives
\begin{align*}
 (1-\mathsf s\mathsf r\mathsf s)(2+\mathsf s+\mathsf r)((tg+\mathsf r g)')={}&((k-2)t+2)\widetilde g+(k-1)\mathsf r \widetilde g\\
 &+t(2-t)\widetilde g'+(2-t)(\mathsf r \widetilde g)'\,.
\end{align*}
Therefore $\operatorname{Res}_K((tg+\mathsf r g)')=0$. Next we apply $\operatorname{Res}_K$ to \eqref{eq:finiteulecture} at $(X,Y)=(t,1)$. The boundary polynomial is fixed by $\mathsf s\mathsf r\mathsf s$, so its residue is zero. Under the substitution $t\mapsto t/(t-1)$, the residue $\operatorname{Res}_K(A(t,1)/t)$ becomes $\operatorname{Res}_K((1-t)A(t,1-t)/t)$, and after differentiating \eqref{eq:finitedspolynomiallecture} the remaining residue is
\begin{align*}
 \operatorname{Res}_K\bigl((1+(k-2)t)g+t(1-t)g'-(k-2)(1-t)g(1-t)+(1-t)^2g'(1-t)\bigr)\,.
\end{align*}
Take $v=-(k-2)tg+t^2g'$. Its degree is at most $k-2$ and $\mathsf r v=-(\mathsf r g)'$, so the residue identities reduce the last display to $\operatorname{Res}_K((tg+\mathsf r g)')=0$. Finally, $\operatorname{Res}_K(t^{j-1})=a_{2\ell-j}$, which proves the theorem.
\end{proof}

\subsection{The values \texorpdfstring{$\za(n,1,k-1-n)$}{zA(n,1,k-1-n)}}\label{sec:finitemiddleone}
Let $\mathcal V_k$ be the space of homogeneous $f\in V_k$ with $f(0,Y)=0$ such that $B=f+f|\epsilon=f|T(1+\epsilon)$ is even. Its dimension is
\begin{align}\label{eq:finitepolynomialdimensionlecture}
 \dim_\Q\mathcal V_k=\frac{k}{2}-2-\dim_\C\ms_k\,,
\end{align}
since these are exactly the equations dual to the reduced formal double zeta space $\widetilde{\dz}_k$ in \eqref{eq:reducedformaldoublezeta}, whose dimension is computed in \cref{prop:reducedformaldoublezetadim}. Notice that this is a statement about rational polynomials and not about the dimension of their finite realizations. The map to the generating polynomial of the values $\za(n,1,k-1-n)$ is
\begin{align}\label{eq:finitemiddlemaplecture}
 \mu_k(f)=\frac{(Y-X)f(Y-X,Y)+X(f(X,Y)+f(Y,X))}{XY}\,.
\end{align}
Since the numerator vanishes on both axes, this is a polynomial. By the double shuffle equations it is symmetric with $\mu_k(f)(1,1)=0$, and by \eqref{eq:finitepolynomialslecture} the finite polynomial satisfies $\mu_k(f)=M$. Recall that $W_{k+2}^+$ is the space of odd period polynomials of weight $k+2$, i.e. of the homogeneous polynomials $P$ of degree $k$ which are odd in $X$ and satisfy $P(X,Y)=P(Y,X)=P(X+Y,Y)+P(X+Y,X)$ (see \cref{lem:lewis}). Surprisingly, these describe the kernel of $\mu_k$.
\begin{theorem}\label{thm:finitemiddlekernellecture}
The following maps are inverse isomorphisms:
\begin{align*}
 W_{k+2}^+&\longrightarrow\ker\mu_k,&\ker\mu_k&\longrightarrow W_{k+2}^+,\\
 P&\longmapsto\frac{P(X-Y,Y)}{XY},&f&\longmapsto XY(f(X,Y)+f(Y,X))\,.
\end{align*}
\end{theorem}
\begin{proof}
Let $\mu_k(f)=0$ and put $B=f+f|\epsilon$. Then $f(X,Y)=(X-Y)B(X-Y,Y)/X$, hence $P=XYB$ is symmetric and odd in both variables, and $f+f|\epsilon=B$ gives $P(X,Y)=P(X-Y,Y)-P(X-Y,X)$. Replacing $X$ by $X+Y$ gives the equation for $W_{k+2}^+$ in \cref{lem:lewis}.

Conversely, let $P\in W_{k+2}^+$ and put $B=P/(XY)$. The period equation gives $P(X,X)=0$. Since $B$ is symmetric and even, it is divisible by $(X-Y)^2$ and $(X+Y)^2$. Therefore $f=P(X-Y,Y)/(XY)$ is a polynomial with $f(0,Y)=0$. The period equation shows that $f$ satisfies both double shuffle equations with symmetric part $B$, and by substitution we get $\mu_k(f)=0$. Finally, the two constructions are inverse to each other.
\end{proof}
Since $\dim_\Q W_{k+2}^+=\dim_\C\ms_{k+2}$ by \cref{thm:eichlershimura}, this proves the following bound for the actual finite values.
\begin{corollary}\label[corollary]{cor:finitemiddleboundlecture}
For every even $k\geq4$, we have
\begin{align*}
 \dim_\Q\left\langle\za(n,1,k-1-n):1\leq n\leq k-2\right\rangle_\Q
 \leq\frac{k}{2}-2-\dim_\C\ms_k-\dim_\C\ms_{k+2}\,.
\end{align*}
\end{corollary}
\begin{proof}
The finite polynomial $f$ belongs to $\mathcal V_k\otimes_\Q\ma$. Thus the span of the coefficients of $M=\mu_k(f)$ has dimension at most $\dim_\Q\mu_k(\mathcal V_k)$, and the statement follows from \cref{thm:finitemiddlekernellecture} and \eqref{eq:finitepolynomialdimensionlecture}.
\end{proof}
For example, the first extra relation occurs in weight $10$, although $\ms_{10}=0$:
\begin{align*}
 \za(3,1,6)+2\za(5,1,4)=0\,.
\end{align*}
It was observed by Risan in \cite[p.~53]{R}, and the extra relation comes from $\ms_{12}$. This is why the values $\za(n,1,k-1-n)$ also involve cusp forms of weight $k+2$. Period polynomials also give identities involving these triple values and products (\cref{ex:finiteperiodtriplelecture}). The formal comparison in \Cref{sec:formalfinitemzv} explains the depth-four spanning theorem, but it does not settle either of the two triple-spanning conjectures or the equality of the finite dimensions. There we will also use a reduced quotient of the formal double zeta space to transfer the period-polynomial relations to finite multiple zeta values.

%% file: chap_MES.tex
\chapter{Multiple Eisenstein series}\label{sec:mes}
After the depth-two constructions in the previous chapter, we now want to study multiple Eisenstein series in arbitrary depth. We will first calculate their Fourier expansion, which will also give the connection to the double-indexed $q$-analogues of multiple zeta values. Since multiple Eisenstein series are just defined for $k_1\geq 3, k_2,\dots,k_r \geq 2$, we then want to regularize them, i.e. extend the definition to all admissible indices. At the end of this chapter, we will have a look at their derivatives and at the $\mathfrak{sl}_2$-algebra structure of multiple Eisenstein series. Multiple Eisenstein series were introduced in \cite{GKZ} in depth $2$, and the author calculated their Fourier expansion in arbitrary depths in \cite{Ba1} (see also \cite{Ba5}). Later, the authors of \cite{BT} extended the definition to all indices and gave a connection to the so-called Goncharov coproduct of formal iterated integrals. Some variants which will not be discussed in this chapter are mentioned at the end of Section~\ref{subsec:mesoverview}. Throughout this chapter we will write $q=e^{2\pi i \tau}$ for $\tau \in \Ha$.
Recall that we defined in Section \ref{subsubsec:basicofmf} the Eisenstein series $\aG_k$ for even $k\geq 4$ and $\tau \in \Ha$  by
\begin{equation}\label{eq:defgklattice}
\aG_k(\tau) = \frac{1}{2}\sum_{\substack{m,n \in \Z \\ (m,n)\neq (0,0)}} \frac{1}{(m\tau + n)^k}\,.
\end{equation}	
Splitting the summation into the parts $m=0$ and $m \in \Z \backslash \{0\}$ we obtain for even $k$
\[ \aG_k(\tau) = \frac{1}{2} \sum_{n \neq 0} \frac{1}{n^k} + \sum_{m = 1}^\infty \left( \sum_{n \in \Z} \frac{1}{(m\tau + n)^k} \right) \,. \]
To calculate the Fourier expansion of the sum on the right, one uses the following.
\begin{proposition}(Lipschitz summation formula) For $k\geq 2$ and $q=e^{2\pi i \tau}$ we have 
    \begin{equation}\label{lipschitz-sum}
    \sum_{n \in \Z} \frac{1}{(\tau + n)^k} = \frac{(-2 \pi i)^k}{(k-1)!} \sum_{d=1}^\infty d^{k-1} q^d \,.
    \end{equation}		
\end{proposition}
\begin{proof}
    Consider the partial fraction expansion and Fourier expansion of the cotangent function
    \begin{align*}
    \frac{1}{\tau} + \sum_{n=1}^\infty \left( \frac{1}{\tau - n} + \frac{1}{\tau + n} \right) = \pi \cot(\pi \tau) = -\pi i - 2 \pi i \sum_{d=1}^\infty q^d\,.
    \end{align*}
    Taking the $k\!-\!1$-th derivative on both sides yields \eqref{lipschitz-sum}. 
\end{proof}
With \eqref{lipschitz-sum} we obtain
\begin{align}\begin{split}\label{eq:gk}
\aG_k(\tau) &= \zeta(k) + \frac{(-2\pi i)^k}{(k-1)!} \sum_{m = 1}^\infty \sum_{d=1}^\infty d^{k-1} q^{md} \\&= \zeta(k) + \frac{(-2\pi i)^k}{(k-1)!} \sum_{n = 1}^\infty \sigma_{k-1}(n) q^{n} \\&= \zeta(k) + (-2\pi i)^k \g(k)\,,
\end{split}
\end{align}
Here and in the following we will always suppress the dependence of $\g(k)$ on $\tau$ and always view it as a function in $\tau$ with $q=e^{2\pi i \tau}$. Formula \eqref{eq:gk} also makes sense for odd $k$ but does not give a modular form, since there are no non-trivial modular forms of odd weight. 
In the following, we want to construct a ``multiple version'' of $\aG_k$, such that its constant term gives multiple zeta values.  

The sum in \eqref{eq:defgklattice} vanishes for odd $k$, therefore instead of summing over the whole lattice, we restrict the summation to the positive lattice points, with positivity coming from an order on the lattice $\Z \tau + \Z$. This, in turn, will also enable us to give a multiple version of the Eisenstein series in an obvious way. 

\begin{definition}Let $\tau \in \Ha$. We define the order $\succ$ on $\Z \tau + \Z$ 	for $\lambda_1,\lambda_2 \in \Z \tau + \Z$ by  
    \[ \lambda_1 \succ \lambda_2 :\Leftrightarrow \lambda_1 - \lambda_2 \in P\,,\]
where $P$, the \emph{set of positive lattice points}, is defined by 
    \[ P := \left\{ m \tau + n \in \Z \tau + \Z \mid {\color{blue}{m > 0}}  \vee  \left( {\color{red}{ m=0 \wedge n>0 }}\right) \right\}  = \UU \cup \RR \,.\]
    \begin{figure}[H]
        \begin{center}
            \begin{tikzpicture}[scale=0.4]
            \draw[dotted,step=1,color=gray,thin] (-4.9,-4.9) grid (4.9,4.9); %
            \draw [->,thick] (0,-5.5) -- (0,5.5) node (yaxis) [above] {$m$};
            \draw [->,thick] (-5.5,0) -- (5.5,0) node (xaxis) [right] {$n$};
            \fill[black] (-5,1) circle (4pt);
            \fill[black] (-4,1) circle (4pt);
            \fill[black] (-3,1) circle (4pt);
            \fill[black] (-2,1) circle (4pt);
            \fill[black] (-1,1) circle (4pt);
            \fill[black] (0,1) circle (4pt);
            \fill[black] (1,1) circle (4pt);
            \fill[black] (2,1) circle (4pt);
            \fill[black] (3,1) circle (4pt);
            \fill[black] (4,1) circle (4pt);
            \fill[black] (5,1) circle (4pt);
            \fill[black] (-5,2) circle (4pt);
            \fill[black] (-4,2) circle (4pt);
            \fill[black] (-3,2) circle (4pt);
            \fill[black] (-2,2) circle (4pt);
            \fill[black] (-1,2) circle (4pt);
            \fill[black] (0,2) circle (4pt);
            \fill[black] (1,2) circle (4pt);
            \fill[black] (2,2) circle (4pt);
            \fill[black] (3,2) circle (4pt);
            \fill[black] (4,2) circle (4pt);
            \fill[black] (5,2) circle (4pt);
            \fill[black] (-5,3) circle (4pt);
            \fill[black] (-4,3) circle (4pt);
            \fill[black] (-3,3) circle (4pt);
            \fill[black] (-2,3) circle (4pt);
            \fill[black] (-1,3) circle (4pt);
            \fill[black] (0,3) circle (4pt);
            \fill[black] (1,3) circle (4pt);
            \fill[black] (2,3) circle (4pt);
            \fill[black] (3,3) circle (4pt);
            \fill[black] (4,3) circle (4pt);
            \fill[black] (5,3) circle (4pt);
            \fill[black] (-5,4) circle (4pt);
            \fill[black] (-4,4) circle (4pt);
            \fill[black] (-3,4) circle (4pt);
            \fill[black] (-2,4) circle (4pt);
            \fill[black] (-1,4) circle (4pt);
            \fill[black] (0,4) circle (4pt);
            \fill[black] (1,4) circle (4pt);
            \fill[black] (2,4) circle (4pt);
            \fill[black] (3,4) circle (4pt);
            \fill[black] (4,4) circle (4pt);
            \fill[black] (5,4) circle (4pt);
            \fill[black] (-5,5) circle (4pt);
            \fill[black] (-4,5) circle (4pt);
            \fill[black] (-3,5) circle (4pt);
            \fill[black] (-2,5) circle (4pt);
            \fill[black] (-1,5) circle (4pt);
            \fill[black] (0,5) circle (4pt);
            \fill[black] (1,5) circle (4pt);
            \fill[black] (2,5) circle (4pt);
            \fill[black] (3,5) circle (4pt);
            \fill[black] (4,5) circle (4pt);
            \fill[black] (5,5) circle (4pt);
            \fill[black] (1,0) circle (4pt);
            \fill[black] (2,0) circle (4pt);
            \fill[black] (3,0) circle (4pt);
            \fill[black] (4,0) circle (4pt);
            \fill[black] (5,0) circle (4pt);
            \draw[red,very thin] (1,0.5) arc (90:270:0.5);
            \draw[red,very thin] (1,0.5) -- (7,0.5);
            \draw[red,very thin] (1,-0.5) -- (7,-0.5);
            \coordinate [label=right:\textcolor{red}{$R$}] (R) at (6.5,0);
            \draw[blue,very thin] (-6,0.7) -- (7,0.7);
            \coordinate [label=right:\textcolor{blue}{$U$}] (R) at (6.5,3);
            \end{tikzpicture}\\
            The set $P$ for the case $\tau =i$.
        \end{center}
    \end{figure}
    In other words we have $m_1 \tau + n_1 \succ m_2 \tau + n_2$ if $m_1>m_2$ or if $m_1=m_2$ and $n_1 > n_2$.		\end{definition}

Since $P \cup (-P) = \Lambda_\tau \backslash \{0\}$ the Eisenstein series can be written for even $k \geq 4$ as
\begin{align*}
\aG_{k}(\tau) = \frac{1}{2}\sum_{\substack{m,n \in \Z \\ (m,n)\neq (0,0)}} \frac{1}{(m\tau + n)^k} = \sum_{\substack{\lambda\succ 0\\ \lambda \in \Lambda_\tau}} \frac{1}{\lambda^{k}}   \,. \qquad (k \geq 4 \text{ even})
\end{align*}
But now the right-hand side does not vanish anymore for odd $k$ and we can use it to define $\aG_{k}$ for all $k\geq 3$ (since the sum converges absolutely for $k>2$).  With the same argument as before we obtain
\begin{align*}
\aG_{k}(\tau) := \sum_{\substack{\lambda\succ 0\\ \lambda \in \Lambda_\tau}} \frac{1}{\lambda^{k}} 	= 	\zeta(k) + (-2\pi i)^k \g(k)\,. \qquad (k \geq 3)
\end{align*}
In general we will define the multiple version of these objects as follows. 
\begin{definition}
    For $k_1 \geq 3, k_2,\dots,k_r \geq 2$ the \emph{multiple Eisenstein series} are defined by
    \begin{align}\label{eq:defmes}
    \aG_{k_1,\dots,k_r}(\tau) := \sum_{\substack{\lambda_1 \succ \dots \succ \lambda_r \succ 0\\ \lambda_i \in \Z \tau + \Z}} \frac{1}{\lambda_1^{k_1} \dots \lambda_r^{k_r}}   \,.
    \end{align}
    By $k= k_1 + \dots + k_r$ we denote its weight and by $r$ its depth. 
\end{definition}
That the sum \eqref{eq:defmes} converges absolutely for $k_1 \geq 3, k_2,\dots,k_r \geq 2$ can be checked by generalizing the proofs of Theorems 4.3 and B.1 in \cite{C}.
The multiple Eisenstein series are holomorphic functions in the upper half-plane and they satisfy the stuffle product formula, i.e. for example we have 
\[ \aG_3(\tau) \cdot \aG_4(\tau) = \aG_{4,3}(\tau) + \aG_{3,4}(\tau) + \aG_{7}(\tau)\,. \] 
But, as we already see in depth one, they are in general not modular forms. The Eisenstein series $\aG_2$ does not converge absolutely but is conditionally convergent. We define it by 
\begin{align*}
\aG_2(\tau) := \lim_{M \rightarrow \infty}   \lim_{N \rightarrow \infty} \left( \sum_{0<n<N} \frac{1}{n^2} +  \sum_{M>m>0} \sum_{-N < n < N } \frac{1}{(m\tau + n)^2} \right) = \zeta(2) + \ag(2)
\end{align*}
i.e., we sum first ``horizontally'' the $n$ and then ``vertically''	the $m$. 
This is usually called \emph{Eisenstein summation}. In general set $\Z_M = \{m \in \Z \mid |m|<M\}$  for an integer $M>0$. With this we define the multiple Eisenstein series for all $k_1,\dots,k_r \geq 2$ by 
\begin{align}\label{def:meseisensteinsum}
\aG_{k_1,\dots,k_r}(\tau) :=\lim_{M \rightarrow \infty}   \lim_{N \rightarrow \infty}  \sum_{\substack{\lambda_1 \succ \dots \succ \lambda_r \succ 0\\ \lambda_i \in \Z_M \tau + \Z_N}} \frac{1}{\lambda_1^{k_1} \dots \lambda_r^{k_r}}   \,.
\end{align}
These also satisfy for all $k_1,\dots,k_r \geq 2$  the stuffle product formula, since one can check easily that this is already true for fixed $M$ and $N$. 

\section{The Fourier expansion of multiple Eisenstein series} \label{subsec:fourier}
By definition it is $\aG_{k_1,\dots,k_r}(\tau + 1) = \aG_{k_1,\dots,k_r}(\tau)$, i.e. there exists a Fourier expansion of $\aG_{k_1,\dots,k_r}$. In depth one we already saw that this is the case:
\begin{align*}
\aG_{k}(\tau) = 	\zeta(k) + (-2\pi i)^k \g(k)\,.
\end{align*}
We will see that it will be possible in general to write $\aG_{k_1,\dots,k_r}$ in terms of multiple zeta values and $(-2\pi i)^{k_1+\dots  + k_r} \g(k_1,\dots,k_r)$. Therefore we define for all $k_1,\dots,k_r\geq 1$
\begin{align}
\ag(k_1,\dots,k_r):=	(-2\pi i)^{k_1+\dots + k_r} \g(k_1,\dots,k_r)
\end{align}
which then gives $	\aG_{k}(\tau) = 	\zeta(k)  + \ag(k)$. Again in everything that follows $\ag$ can be viewed as a holomorphic function in $\tau$ but we will avoid writing $ \ag(k_1,\dots,k_r; \tau)$ in order to make notations shorter.

\begin{theorem}\label{thm:Fourier}
    The $\aG_{k_1,\dots,k_r}(\tau)$ can be written as a $\mz$-linear combination of $\ag$. 
    \\
    More precisely, for  $k_1,\dots,k_r \geq 2$ there exist rational numbers $\alpha^{(l_1,\dots,l_r)}_j \in \Q$, for $l_1,\dots,l_r \geq 2$ and $1\leq j \leq r-1$ with  $k=k_1+\dots+k_r = l_1+ \dots + l_r$, such that
    \begin{align*} \aG_{k_1,\dots,k_r}(\tau) \,= \,\zeta(k_1,\dots,k_r) \;+\; \sum_{\substack{1\leq j \leq r-1 \\ l_1,\dots,l_r \geq 2\\ l_1 +\dots+l_r = k}} \alpha^{(l_1,\dots,l_r)}_j 
    \cdot  \zeta(l_1,\dots,l_j) \cdot \ag(l_{j+1},\dots,l_r) \;+\; \ag(k_1,\dots,k_r) \,. 
    \end{align*}
\end{theorem}
The rest of this section will be devoted to proving this Theorem. 	
For example, the triple Eisenstein series $\aG_{3,2,2}$ can be written as
\begin{align*}
\aG_{3,2,2}(\tau) =& \zeta(3,2,2) + \left( \frac{54}{5} \zeta(2,3) + \frac{51}{5} \zeta(3,2) \right) \ag(2) + \frac{16}{3} \zeta  (2,2)\ag(3)  \\
&+3 \zeta(3) \ag(2,2) + 4 \zeta(2) \ag(3,2)  + \ag(3,2,2)\,.
\end{align*}

To derive the Fourier expansion, we introduce the following functions, which can be seen as a multiple version of the term $\sum_{n \in \Z} \frac{1}{(x+n)^k}$ appearing in the calculation of the Fourier expansion of classical Eisenstein series. 

\begin{definition}\label{def:multitangent}
    For $k_1,\dots,k_r \geq 2$ and $x \in \C \backslash \Z$ we define the multitangent function of depth $r$ by
    \[  \Psi_{k_1,\dots,k_r}(x) = \sum_{\substack{ n_1 > \cdots > n_r \\ n_i \in \Z}} \frac{1}{(x+n_1)^{k_1} \dots (x+n_r)^{k_r}} \,.  \]
    In the depth $r=1$ case we also refer to these as monotangent functions. 
\end{definition}

These functions were introduced and studied in detail by Bouillot in  \cite{Bo}. One of the main results in \cite{Bo}, which is crucial for the calculation of the Fourier expansion presented here, is the following theorem, which reduces the multitangent functions into monotangent functions.

\begin{theorem}(\cite[Theorem 3]{Bo}, Reduction of multitangent into monotangent)\label{thm:reduction}
For $k_1,\dots,k_r \geq 2$ with $k=k_1+\dots+k_r$ the multitangent function can be written as \\
\scalebox{0.96}{\parbox{\linewidth}{%
\begin{align*}
    \Psi_{k_1,\dots,k_r}(\tau) = \sum_{\substack{1\leq j \leq r\\ l_1+\dots+l_r = k\\l_1,\dots,l_r\geq1}} (-1)^{l_1+\dots+l_{j-1}+k_j+k} \prod_{\substack{1\leq i \leq r\\i\neq j}}\binom{l_i-1}{k_i-1} \zeta(l_1,\dots,l_{j-1}) \,\Psi_{l_j}(\tau)\, \zeta(l_r,l_{r-1},\dots,l_{j+1}).
\end{align*} }} \\
Moreover, the terms containing $\Psi_1(\tau)$ vanish. 
\end{theorem}
\begin{proof}
 The proof uses partial fraction decomposition and the antipode relation for (shuffle regularized) multiple zeta values (see \cite[Lemma 2.4]{BT}) to argue that the coefficient in front of $\Psi_1(\tau)$ vanishes. 
    For example, in depth two, we have
    \begin{align}\label{eq:psi32}
    \begin{split}
    &\Psi_{3,2}(x) = \sum_{n_1 > n_2} \frac{1}{(x+n_1)^3 (x+n_2)^2} \\
    &= \sum_{n_1 > n_2}  \left( \frac{1}{(n_1-n_2)^2 (x+n_1)^3} +\frac{2}{(n_1-n_2)^3 (x+n_1)^2} + \frac{3}{(n_1-n_2)^4 (x+n_1)} \right) \\
    &+\sum_{n_1 > n_2} \left( \frac{1}{(n_1-n_2)^3 (x+n_2)^2}  -  \frac{3}{(n_1-n_2)^4 (x+n_2)} \right) \\
    &=  3 \zeta(3) \Psi_2(x) +  \zeta(2) \Psi_3(x)  \,. 
    \end{split}
    \end{align}
\end{proof}

The connection between the functions $\ag$ and the monotangent functions is given by the following 
\begin{proposition}\label{prop:ghatasmonotangent}
    For $k_1,\dots,k_r \geq 2$ the functions $\ag$ can be written as 
    \begin{align*}
    \ag(k_1,\ldots,k_r) = \sum_{m_1 > \dots > m_r > 0} \Psi_{k_1}(m_1\tau) \dots \Psi_{k_r}(m_r\tau) \,.
    \end{align*}
\end{proposition}
\begin{proof}
    This follows directly from the Lipschitz formula \eqref{lipschitz-sum} and the definition of the functions $\ag$, since 
    \begin{align*} 
    \Psi_{k}(m\tau)  = \sum_{n \in \Z} \frac{1}{(m\tau + n)^k} = \frac{(-2 \pi i)^k}{(k-1)!} \sum_{d=1}^\infty d^{k-1} q^{md}
    \end{align*}
    and therefore 
    \begin{align*}
    \sum_{m_1 > \dots > m_r > 0} \Psi_{k_1}(m_1\tau) \dots \Psi_{k_r}(m_r\tau) &= (-2 \pi i)^{k_1+\dots+k_r} \sum_{\substack{m_1 > \dots > m_r > 0\\d_1,\dots,d_r>0}} \\
    &\qquad{}\cdot \frac{d_1^{k_1-1}}{(k_1-1)!} \cdots \frac{d_r^{k_r-1}}{(k_r-1)!} q^{m_1 d_1 + \dots + m_r d_r} \\
    &= (-2 \pi i)^{k_1+\dots+k_r}  \g(k_1,\dots,k_r) = \ag(k_1,\dots,k_r) \,.
    \end{align*}		
\end{proof}

{\bf Preparation for the Proof of Theorem \ref{thm:Fourier}:} We will now recall the construction of the Fourier expansion of multiple Eisenstein series introduced in \cite{Ba1}, in order to prove Theorem \ref{thm:Fourier}.
To calculate the Fourier expansion we rewrite the multiple Eisenstein series as 
\begin{align*}
\aG_{k_1,\dots,k_r}(\tau) &= \sum_{\lambda_1 \succ \dots \succ \lambda_r \succ 0} \frac{1}{\lambda_1^{k_1} \dots \lambda_r^{k_r}} \\
&= \sum_{(\lambda_1,\dots,\lambda_r) \in P^r} \frac{1}{(\lambda_1+\dots+\lambda_r)^{k_1} (\lambda_2 + \dots + \lambda_r )^{k_2} \dots \lambda_r^{k_r} } \,.
\end{align*}
We decompose the set of tuples of positive lattice points $P^r$  into the $2^r$ distinct subsets $A_1 \times \dots \times A_r \subset P^r$ with $A_i \in \{R,U\}$ and write 
{\small
    \[ \aG^{A_1 \dots A_r}_{k_1,\dots,k_r}(\tau) := \sum_{(\lambda_1,\dots,\lambda_r) \in A_1 \times \dots \times A_r} \frac{1}{(\lambda_1+\dots+\lambda_r)^{k_1} (\lambda_2 + \dots + \lambda_r )^{k_2} \dots \lambda_r^{k_r} } \]
}
this gives the decomposition
\[ \aG_{k_1,\dots,k_r} = \sum_{A_1,\dots,A_r \in \{R,U\}} \aG^{A_1\dots A_r}_{k_1,\dots,k_r} \,. \]
In the following, we identify the $A_1 \dots A_r$ with words in the alphabet $\{R,U\}$. We first illustrate the general algorithm in depth one and two.

\begin{ex} \begin{enumerate}[\textup{(}i\textup{)}] 
        \item 	In depth $r=1$ we have $\aG_{k}(\tau) =\aG_k^R(\tau) + \aG_k^U(\tau)$ and
        \begin{align*}
        \aG_k^R(\tau) &= \sum_{\substack{m_1 = 0 \\ n_1 > 0}} \frac{1}{(0\tau + n_1)^k} = \zeta(k) \,,\\
        \aG_k^U(\tau) &= 	\sum_{\substack{m_1 > 0 \\ n_1 \in \Z}} \frac{1}{(m_1\tau + n_1)^k} = \sum_{m_1>0} \Psi_k(m_1\tau) = \ag(k)\,, 
        \end{align*}
        which gives $\aG_{k}(\tau) = \zeta(k) + \ag(k)$.
        \item 
        In depth $2$ we have $ \aG_{k_1,k_2} = \aG^{RR}_{k_1,k_2}+\aG^{UR}_{k_1,k_2} +\aG^{RU}_{k_1,k_2}+\aG^{UU}_{k_1,k_2}$. The $RR$ and $UU$ part is similar to the depth one case and we get
        \begin{align*}
        \aG^{RR}_{k_1,k_2} &= \sum_{\substack{m_1 = m_2=0 \\ n_1 > n_2 > 0}} \frac{1}{(0\tau + n_1)^{k_1} (0\tau + n_2)^{k_2}} = \zeta(k_1,k_2)\,,\\ 
        \aG^{UU}_{k_1,k_2} &= \sum_{\substack{m_1 > m_2 >0 \\ n_1,n_2 \in \Z}} \frac{1}{(m_1\tau + n_1)^{k_1} (m_2\tau + n_2)^{k_2}} =\sum_{m_1>m_2>0} \Psi_{k_1}(m_1\tau) \Psi_{k_2}(m_2\tau)  = \ag(k_1,k_2)\,.
        \end{align*}
        For the word $UR$ we get
        \begin{align*}
        \aG^{UR}_{k_1,k_2} &= \sum_{\substack{m_1>0, m_2=0 \\ n_1 \in \Z, n_2 > 0}} \frac{1}{(m_1 \tau + n_1)^{k_1} (0\tau + n_2)^{k_2}} 
        =\sum_{m_1>0 }\Psi_{k_1}(m_1 \tau) \sum_{n_2>0} \frac{1}{n_2^{k_2}} = \ag(k_1) \zeta(k_2)  \,.
        \end{align*}
        Finally, to evaluate the $RU$ part, we need to use Theorem \ref{thm:reduction} which gives
        \begin{align*}
        \Psi_{k_1,k_2}(x) = \sum_{j=2}^{k_1+k_2} c^{k_1,k_2}_j  \Psi_j (x)\,.
        \end{align*}
        Using this we can write 
        \begin{align*}
        \aG^{RU}_{k_1,k_2}(\tau) &= \sum_{\substack{m_1=0, m_2>0 \\ n_1 > n_2 \\ n_i \in \Z}} \frac{1}{(m_2 \tau + n_1)^{k_1} (m_2 \tau + n_2)^{k_2}} = \sum_{m>0} \Psi_{k_1,k_2}(m\tau) = \sum_{j=2}^{k} c^{k_1,k_2}_j  \sum_{m>0} \Psi_{j}(m\tau) \\
        & = \sum_{j=2}^{k_1+k_2} c^{k_1,k_2}_j  \,\ag(j)\,.
        \end{align*}
        In total we therefore obtain 
        \begin{align*}
        \aG_{k_1,k_2}(\tau) = \zeta(k_1,k_2) +\zeta(k_2)  \ag(k_1)  + \sum_{j=2}^{k_1+k_2} c^{k_1,k_2}_j \, \ag(j) + \ag(k_1,k_2)
        \end{align*}
        with $c^{k_1,k_2}_j  \in \mz_{k_1+k_2-j}$.
    \end{enumerate}
\end{ex}

In general, by using Proposition \ref{prop:ghatasmonotangent} the $\aG^{U^r}_{k_1,\dots,k_r}$ can be written as
\begin{align*}
\aG^{U^r}_{k_1,\dots,k_r}(\tau) &= \sum_{\substack{m_1>\dots>m_r>0 \\ n_1,\dots,n_r \in \Z}} \frac{1}{(m_1\tau+n_1)^{k_1} \dots (m_r\tau + n_r)^{k_r}}\\
&= \sum_{m_1 > \dots > m_r>0} \Psi_{k_1}(m_1\tau) \dots \Psi_{k_r}(m_r \tau) = \ag(k_1,\dots,k_r)\,.
\end{align*}

The other special case $\aG^{R^r}_{k_1,\dots,k_r}$ can also be written down explicitly:
\begin{align*}
\aG^{R^r}_{k_1,\dots,k_r}(\tau) =  \sum_{\substack{m_1 =\dots=m_r=0 \\ n_1 > \dots > n_r > 0}} \frac{1}{(0\tau + n_1)^{k_1} \dots (0\tau + n_r)^{k_r}} = \zeta(k_1,\dots,k_r) \,.
\end{align*}
In the case $\aG^{UR}$ we saw that we could write it as $\aG^U$ multiplied with a zeta value. 
In general, for a word $w\neq R^r$ of depth $r$ ending in the letter $R$, there is a word $w'$ ending in $U$ with $w = w' R^l$ and $1 \leq l < r$, and we can write
\begin{align*}
\aG^{w}_{k_1,\dots,k_r}(\tau) =  \aG^{w'}_{k_1,\dots,k_{r-l}}(\tau) \cdot \zeta(k_{r-l+1},\dots,k_r) \,.
\end{align*}
For example we have $\aG^{RUURR}_{3,4,5,6,7} = \aG^{RUU}_{3,4,5} \cdot \zeta( 6,7)$. Hence one can concentrate on the words ending in $U$ when calculating the Fourier expansion of a multiple Eisenstein series.
Let $w$ be a word ending in $U$ then there are integers $r_1,\dots,r_j \geq 1$ with $w = R^{r_1-1} U R^{r_2-1} U \dots R^{r_j-1} U$. With this one can write
\[ \aG^w_{k_1,\dots,k_r}(\tau) = \sum_{m_1 > \dots > m_j > 0} \Psi_{k_1,\dots,k_{r_1}}(m_1 \tau) \cdot \Psi_{k_{r_1+1},\dots,k_{r_1+r_2}}(m_2 \tau) \dots \Psi_{k_{r-r_j+1},\dots,k_r}(m_j \tau) \,,\] 
which gives a $\mz$-linear combination of $\ag$ by using Theorem \ref{thm:reduction}.

\begin{ex} For example for the word $w=RURRU$ we have
    \begin{align*}
    \aG_{k_1,\dots,k_5}^{RURRU}(\tau) &= \sum_{m_1 > m_2 >0} \Psi_{k_1,k_2}(m_1 \tau) \Psi_{k_3,k_4,k_5}(m_2 \tau) \\
    &= \sum_{m_1 > m_2 >0} \sum_{\substack{2 \leq j_1 \leq k_1+ k_2\\2 \leq j_2 \leq k_3+k_4+k_5}} c^{k_1,k_2}_{j_1} c^{k_3,k_4,k_5}_{j_2}  \Psi_{j_1}(m_1 \tau)  \Psi_{j_2}(m_2 \tau) \\
    &= \sum_{\substack{2 \leq j_1 \leq k_1+ k_2\\2 \leq j_2 \leq k_3+k_4+k_5}} c^{k_1,k_2}_{j_1} c^{k_3,k_4,k_5}_{j_2} \,   \ag(j_1, j_2)\,.
    \end{align*}
    See Figure \ref{fig:RURRU} for an example of a summand of  $\aG_{k_1,\dots,k_5}^{RURRU}(\tau)$.
    \begin{figure}[h!] 
        \begin{center}
            \begin{tikzpicture}[scale=0.6]
            \draw[densely dotted,step=1,color=gray,thin] (-4.9,-4.9) grid (4.9,4.9);
            \draw [->] (0,-5.5) -- (0,5.5) node (yaxis) [above] {$m$};
            \draw [->] (-5.5,0) -- (5.5,0) node (xaxis) [right] {$n$};
            \draw[blue] (0,0) -- node[anchor=west] {$\lambda_5$} (1,1);
            \draw[red] (1,1) -- node[anchor=north] {$\lambda_4$} (3,1); 
            \draw[red] (3,1)-- node[anchor=north] {$\lambda_3$} (4,1);
            \draw[blue] (4,1) --  node[anchor=west] {$\lambda_2$} (-1,4);
            \draw[red] (-1,4) --  node[anchor=south] {$\lambda_1$} (2,4);
            \fill[black] (1,1) circle (3pt);
            \fill[black] (3,1) circle (3pt);
            \fill[black] (4,1) circle (3pt);
            \fill[black] (-1,4) circle (3pt);
            \fill[black] (2,4) circle (3pt); 
            \end{tikzpicture} \linebreak
        \end{center}
        \caption{			A summand of $\aG^{\RR \UU \RR \RR \UU}_{k_1,\dots,k_5}(\tau)$ for $\tau=i$.}
        \label{fig:RURRU}
    \end{figure}
\end{ex}

\begin{proof}[Proof of Theorem \ref{thm:Fourier}:]
    For $k_1,\dots,k_r \geq 2$ the Fourier expansion of the multiple Eisenstein series $\aG_{k_1,\dots,k_r}$ can be computed in the following way
    \begin{enumerate}[\textup{(}i\textup{)}] 
        \item Split up the summation into $2^r$ distinct parts $\aG^w_{k_1,\dots,k_r}$ where $w$ are words in $\{R,U\}$. 
        \item For $w$ being a word  ending in  $R$ one can write $\aG^w_{k_1,\dots,k_r}$  as $\aG^{w'}_{k_1,\dots} \cdot \zeta(\dots,k_r)$ with  $w'$ ending in $U$. 
        \item For $w = R^{r_1-1} U R^{r_2-1} U \dots R^{r_j-1} U$ being a word ending in $U$ one can write $\aG^w_{k_1,\dots,k_r}$ as
        \[ \aG^w_{k_1,\dots,k_r}(\tau) = \sum_{m_1 > \dots > m_j > 0} \Psi_{k_1,\dots,k_{r_1}}(m_1 \tau) \cdot \Psi_{k_{r_1+1},\dots,k_{r_1+r_2}}(m_2 \tau) \dots \Psi_{k_{r-r_j+1},\dots,k_r}(m_j \tau) \,. \] 
        \item Using the Theorem \ref{thm:reduction} we can write the  multitangent functions in iii) as a $\mz$-linear combination of monotangents. We therefore just have $\mz$-linear combinations with sums of the form
        \[    \sum_{m_1 > \dots > m_j >0} \Psi_{k_1}(m_1 \tau) \dots \Psi_{k_j}(m_j \tau)=\ag(k_1,\dots,k_j)\,. \]
    \end{enumerate}
\end{proof}

An explicit formula for the Fourier expansion of the multiple Eisenstein series for arbitrary depth can be found in \cite[Proposition~2.4]{BT} (with the order of indices reversed). Here we just give the Fourier expansion in depths two and three. For this we define for $n_1,n_2,k > 0$ the numbers $C_{n_1,n_2}^k$ by
\[ C_{n_1,n_2}^k = (-1)^{n_2} \binom{k-1}{n_2-1} + (-1)^{k-n_1} \binom{k-1}{n_1-1} \,. \]
\begin{proposition}\label{prop:fourier23}
    \begin{enumerate}[\textup{(}i\textup{)}] 
        \item (\cite[Formula (52)]{GKZ}, \cite{Ba1}, \cite{BT}) For $k_1,k_2 \geq 2$ the Fourier expansion of the double Eisenstein series is given by
        \begin{align*}
        \aG_{k_1,k_2}(\tau) = \zeta(k_1,k_2) + \zeta(k_2) \ag(k_1) + \sum_{\substack{l_1+l_2 = k_1+k_2 \\ l_1,l_2 \geq 2}} C^{l_2}_{k_1,k_2} \zeta(l_2) \ag(l_1) + \ag(k_1,k_2) \,.
        \end{align*}
        \item (\cite{Ba1}, \cite{BT}) For $k_1,k_2,k_3 \geq 2$ and $k=k_1+k_2+k_3$ the Fourier expansion of the triple Eisenstein series can be written as
        \begin{align*}
        \aG_{k_1,k_2,k_3}(\tau) &=  \zeta(k_1,k_2,k_3) + \zeta(k_2,k_3) \ag(k_1) + \zeta(k_3) \ag(k_1,k_2) + \ag(k_1,k_2,k_3) \\
        &+\zeta(k_3) \sum_{l_1+l_2 = k_1 + k_2} C_{k_1,k_2}^{l_1} \zeta(l_1) \ag(l_2) \\
        &+\sum_{ l_1+l_2 = k_1 + k_2} C_{k_1,k_2}^{l_2}  \zeta(l_2) \ag(l_1,k_3)+\sum_{ l_1+l_2 = k_2 + k_3} C_{k_2,k_3}^{l_2} \zeta(l_2) \ag(k_1,l_1) \\
        &+\sum_{l_1+l_2+l_3 = k} (-1)^{k_2+k_3} \binom{l_2-1}{k_2-1}\binom{l_3-1}{k_3-1} \zeta(l_3,l_2) \ag(l_1)   \\
        &+\sum_{l_1+l_2+l_3 = k}(-1)^{k_1+k_2+l_2+l_3} \binom{l_2-1}{k_2-1}\binom{l_3-1}{k_1-1} \zeta(l_3,l_2) \ag(l_1) \\
        &+(-1)^{k_1+k_3}\sum_{l_1+l_2+l_3 = k} (-1)^{l_2} \binom{l_2-1}{k_1-1}\binom{l_3-1}{k_3-1} \zeta(l_3) \zeta(l_2)  \ag(l_1) \,,
        \end{align*}
        where in the sums we sum over all $l_i \geq 2$. 
    \end{enumerate}
\end{proposition} 
We finish this section with a closer look at the stuffle product of two Eisenstein series. 
Since the product of multiple Eisenstein series can be written using the stuffle product, we have $\aG_2 \cdot \aG_3 = \aG_{2,3} + \aG_{3,2} + \aG_5$. On the other hand, we have
\begin{align*}
\aG_2(\tau) \cdot \aG_3(\tau) &= \left(\zeta(2) + \ag(2) \right) \left( \zeta(3) + \ag(3) \right) = \zeta(2) \zeta(3) + \zeta(3) \ag(2) + \zeta(2) \ag(3) + \ag(2) \cdot \ag(3) \,. 
\end{align*}
By Proposition \ref{prop:fourier23} we obtain
\begin{align*}
\aG_{2,3}(\tau)  &= \zeta(2,3) - 2 \zeta(3) \ag(2) + \zeta(2) \ag(3) + \ag(2,3) \,, \\
\aG_{3,2}(\tau) &= \zeta(3,2) + 3 \zeta(3) \ag(2) + 2 \zeta(2) \ag(3) + \ag(3,2) \,.
\end{align*}
In conclusion, we obtain $\ag(2) \cdot \ag(3) = \ag(3,2) + \ag(2,3) + \ag(5) + 2 \zeta(2) \ag(3)$. Dividing out $(-2\pi i)^5$ gives
\[ \g(2) \cdot \g(3) = \g(3,2) + \g(2,3) + \g(5) - \frac{1}{12} \g(3) \,,\] 
which we already saw in \eqref{eq:g2g3prod}
as a special case of Proposition \ref{prop:gstuffle1}.

\section{Regularized multiple Eisenstein series}

The multiple Eisenstein series $\aG_{k_1,\dots,k_r}$ were just defined for $k_1,\dots,k_r \geq 2$ in the previous section and we saw that they have a Fourier expansion of the form
\begin{align}\label{eq:messimplefourier}
\aG_{k_1,\dots,k_r}(\tau) = \zeta(k_1,\dots,k_r) + \sum_{n=1}^\infty a_n q^n\,.\qquad  \qquad (a_n \in \mz[\pi i] = \mz + \pi i \mz)
\end{align}

A natural question therefore is whether there exists a ``good'' extension of these objects for all admissible indices $k_1\geq 2, k_2, \dots, k_r\geq 1$. By ``good'' one could have different properties in mind that should be satisfied by these extended objects. One certainly is that they also should have a Fourier expansion of the form \eqref{eq:messimplefourier}. Another property could be that the extended version also satisfies the shuffle or stuffle product formula. We want to present two types of regularization: the shuffle regularized multiple Eisenstein series (\cite{BT}) and the stuffle regularized multiple Eisenstein series (\cite{Ba7,Ba10}). The definition of shuffle regularized multiple Eisenstein series uses a beautiful connection between the Fourier expansion of multiple Eisenstein series and the coproduct of formal iterated integrals. The other regularization, the stuffle regularized multiple Eisenstein series, uses the construction of the Fourier expansion of multiple Eisenstein series and a coproduct on $\h^1_\ast$ together with a result on regularization of multitangent functions by O. Bouillot (\cite{Bo}).
First recall from the last section that for $k_1,\dots,k_r \geq 2$ the $\aG_{k_1,\dots,k_r}$ satisfy the stuffle product formula. We define the $\Q$-vector space $\h^2 = \Q\langle z_2,z_3, \dots \rangle \subset \h^1$. Equipped with the stuffle product $\ast$, one can see easily that we obtain a subalgebra $\hz_\ast \subset  \h^1_\ast$.  Now we can view the multiple Eisenstein series as an algebra homomorphism defined on the generators by 
\begin{align*}
\aG : \hz_\ast  &\longrightarrow \MZB \\
z_{k_1} \cdots z_{k_r} &\longmapsto \aG_{k_1,\dots,k_r}\,.
\end{align*}
By abuse of notation we will use $\aG$ for both the map and the multiple Eisenstein series. The rough idea to extend the map $\aG$ to $\h^1$ will be as follows. On the algebra $\h^1_\sh$ one can define the Goncharov coproduct $\Delta_G$ and on $\h^1_\ast$ one can define the deconcatenation coproduct $\Delta_H$.  In addition to this we will construct algebra homomorphisms  ($\bullet \in \{ \ast , \shuffle \}$)
\begin{align*}
\ag^\bullet : \h^1_\bullet \rightarrow \Q[2\pi i][[q]]\,.
\end{align*}	 
With this we then can define two algebra homomorphisms $\aG^\bullet  :\h^1_\bullet \rightarrow \MZB$ as follows\footnote{The $\aG^{\ast}(w)$ will depend on an integer $M$, for which we can take the limit $M\rightarrow \infty$ in the case when $w \in \h^0$.}. 	 \begin{align*}
&
\xymatrix{  
    \h^1_\sh  \ar@{->}[r]^-{\Delta_G} \ar@{->}[d]_{\aG^\sh} & \h^1_\sh  \otimes \h^1_\sh \ar@{->}[d]^{\ag^{\sh} \otimes  \zeta^{\sh} }  \\ 
    \mz[\pi i][\![q]\!]  &   \Q[2\pi i][\![q]\!] \otimes \,\mz  \ar@{->}[l]^-{m}\\
} 
&
\xymatrix{  
    \h^1_\ast  \ar@{->}[r]^-{\Delta_H} \ar@{->}[d]_{\aG^\ast} & \h^1_\ast  \otimes \h^1_\ast  \ar@{->}[d]^{\ag^{\ast} \otimes  \zeta^{\ast} }  \\ 
    \mz[\pi i][\![q]\!]  &   \Q[2\pi i][\![q]\!] \otimes \,\mz  \ar@{->}[l]^-{m}\,,\\
} 
\end{align*}
where in both cases $m$ denotes the usual multiplication.  These are both extensions of the original multiple Eisenstein series, in the sense that we have the following
\begin{align} \label{eq:regandnormalmes}
\aG = \aG^\sh{\mid_{\hz}} = \aG^\ast{\mid_{\hz}} \,.
\end{align}

We start by reviewing the definition of formal iterated integrals and the coproduct defined by Goncharov.  In an explicit example in depth two, we will see the connection of this coproduct and the calculation of the Fourier expansion of multiple Eisenstein series from the previous section. This will give an indication of why  \eqref{eq:regandnormalmes} holds. 
After this, we give the definition of shuffle and stuffle regularized multiple Eisenstein series as presented in \cite{BT}, \cite{Ba7} and \cite{Ba10}. At the end of this section, we compare these two regularizations with the help of a few examples.

\subsection{Formal iterated integrals} \label{sec:formaliteratedintegrals}
Following Goncharov (Section 2 in \cite{Gon}) we consider the $\Q$-algebra ${\mathcal I}$ generated by the elements
\[ \mathbb{I}(a_0;a_1,\ldots,a_N;a_{N+1}), \quad a_i\in\{0,1\}, N\ge0. \]
together with the following relations
\begin{enumerate}[\textup{(}i\textup{)}] 
    \item For any $a,b\in\{0,1\}$ the unit is given by $\mathbb{I}(a;b):=\mathbb{I}(a;\emptyset;b)=1$.
    \item The product is given by the shuffle product $\sh$
    \begin{align*}
    & \mathbb{I}(a_0;a_1,\ldots,a_M;a_{M+N+1}) \mathbb{I}(a_0;a_{M+1},\ldots,a_{M+N};a_{M+N+1})\\
    &=\sum_{\sigma\in sh_{M,N}} \mathbb{I}(a_0;a_{\sigma^{-1}(1)},\ldots,a_{\sigma^{-1}(M+N)};a_{M+N+1}),
    \end{align*}
    where $sh_{M,N}$ is the set of $\sigma \in \mathfrak{S}_{M+N}$ such that $\sigma(1)<\cdots<\sigma(M)$ and $\sigma(M+1)<\cdots<\sigma(M+N)$.
    \item The path composition formula holds: for any $N\ge0$ and $a_i,x\in\{0,1\}$, one has
    \[ \mathbb{I}(a_0;a_1,\ldots,a_N;a_{N+1}) = \sum_{k=0}^N \mathbb{I}(a_0;a_1,\ldots,a_k;x)  \mathbb{I}(x;a_{k+1},\ldots,a_N;a_{N+1}).\]
    \item For $N\ge1$ and $a_i,a\in\{0,1\}$ it is $\mathbb{I} (a;a_1,\ldots,a_N;a)=0$.
    \item The path inversion formula holds:
    \[ \mathbb{I}(a_0;a_1,\ldots,a_N;a_{N+1}) = (-1)^N \mathbb{I}(a_{N+1};a_N,\ldots,a_1;a_0)\,. \]
\end{enumerate}

\begin{definition}(Goncharov coproduct) 
    Define the coproduct $\Delta_G$ on ${\mathcal I}$ by
    \begin{align*} 
    &\Delta_G \left( \mathbb{I}(a_0;a_1,\ldots,a_N;a_{N+1}) \right) :=\\
    &\sum \left( \mathbb{I}(a_0;a_{i_1},\ldots,a_{i_k};a_{N+1}) \otimes \prod_{p=0}^k \mathbb{I}(a_{i_p};a_{i_p+1},\ldots,a_{i_{p+1}-1};a_{i_{p+1}}) \right)  ,
    \end{align*}
    where the sum on the right runs over all $i_0=0<i_1<\cdots<i_k<i_{k+1}=N+1$ with $0\le k \le N$.
\end{definition}
\begin{proposition}(\cite[Proposition~2.2]{Gon})
    The triple $({\mathcal I},\sh, \Delta_G)$ is a commutative graded Hopf algebra. 
\end{proposition}

To calculate $\Delta_G \left( \mathbb{I}(a_0;a_1,\dots,a_8;a_{9}) \right)$ one sums over all possible diagrams of the following form. 
\begin{figure}[H]
    \begin{center}
        \begin{tikzpicture}[scale=1.5]
        \def \ra {2}
        \def \noteradius {2.2}
        \def \circsize {0.03}	
        \def \braceradius {2.4}
        \draw (-\ra,0) -- (\ra,0)
        (\ra,0) arc (0:180:\ra);

        \fill (20:\ra) circle (\circsize); 
        \node at (20:\noteradius) {\small $a_8$}; 
        
        \fill (40:\ra) circle (\circsize); 
        \node at (40:\noteradius) {\small $a_7$}; 
        
        \fill(60:\ra) circle (\circsize); 
        \node at (60:\noteradius) {\small $a_6$}; 
        
        \fill(80:\ra) circle (\circsize); 
        \node at (80:\noteradius) {\small $a_5$}; 
        
        \fill (100:\ra) circle (\circsize); 
        \node at (100:\noteradius) {\small $a_4$}; 								
        
        \fill (120:\ra) circle (\circsize); 
        \node at (120:\noteradius) {\small $a_3$}; 
        
        \fill (140:\ra) circle (\circsize); 
        \node at (140:\noteradius) {\small $a_2$}; 	
        
        \fill (160:\ra) circle (\circsize); 
        \node at (160:\noteradius) {\small $a_1$}; 	
        
        \fill(-\ra,0) circle (\circsize ); 				
        \draw (-\ra,0) circle (\circsize ); 
        
        \fill (\ra,0) circle (\circsize ); 
        \node at (180:\noteradius) {\small $a_0$}; 	
        \node at (0:\noteradius) {\small $a_9$};

        \draw  (0:\ra) to [bend left] (40:\ra) ;
        \draw  (40:\ra) to [bend left] (100:\ra) ;	
        \draw  (100:\ra) to [bend left] (160:\ra) ;		
        \draw  (160:\ra) to [bend left] (180:\ra) ;		
        
        \braceme[thick]{\braceradius}{0}{40}{br1}{$I(a_7;a_8;a_9)$}
        \braceme[thick]{\braceradius}{40}{100}{br1}{$I(a_4;a_5,a_6;a_7)$}
        \braceme[thick]{\braceradius}{100}{160}{br1}{$I(a_1;a_2,a_3;a_4)$}
        \braceme[thick]{\braceradius}{160}{180}{br1}{$I(a_0;a_1)$}
        \end{tikzpicture}
        \caption{One diagram for the calculation of  $\Delta_G \left( \mathbb{I}(a_0;a_1,\dots,a_8;a_{9}) \right)$. It gives the term 
            $I(a_0;a_1,a_4,a_7;a_9)  \otimes I(a_0;a_1)I(a_1;a_2,a_3;a_4)I(a_4;a_5,a_6;a_7)I(a_7;a_8;a_9) \,.$}
    \end{center}
\end{figure}
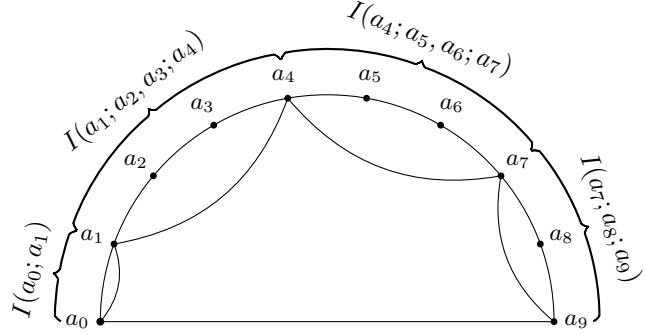
For our purpose it will be important to consider the quotient space\footnote{If one likes to interpret the integrals as real integrals, then the passage from $\mathcal{I}$ to $\mathcal{I}^1$  regularizes these integrals such that ``$-\log(0) = \int_{1 > t > 0} \frac{dt}{t} := 0$''. }
\[ \mathcal{I}^1=\mathcal{I}/\mathbb{I}(1;0;0)\mathcal{I} \,.\]
Let us denote by 
\[ I(a_0;a_1,\ldots,a_N;a_{N+1})\]
the image of $\mathbb{I}(a_0;a_1,\ldots,a_N;a_{N+1})$ in $\mathcal{I}^1$. The quotient map $\mathcal{I} \rightarrow \mathcal{I}^1$ induces a Hopf algebra structure on $\mathcal{I}^1$, but for our application we just need that for any $w_1,w_2\in\mathcal{I}^1$, one has $\Delta_G(w_1\ \sh \ w_2)=\Delta_G(w_1)\ \sh \ \Delta_G(w_2)$. The coproduct on $\mathcal{I}^1$ is given by the same formula as before by replacing $\mathbb{I}$ with $I$.
For integers $n\ge0,k_1,\ldots,k_r\ge1$, we set 
\[I_{n}(k_1,\ldots,k_r):=I(1;\underbrace{0,0,\ldots,1}_{k_1},\ldots,\underbrace{0,0,\ldots,1}_{k_r},\underbrace{0,\ldots,0}_{n};0).\]
In particular, we write\footnote{This notion fits well with the iterated integral expression of multiple zeta values. Recall that 
    \[ \zeta(2,3) = \int_{\scriptstyle 1 > t_1 > \dots > t_5 > 0} \underbrace{\frac{dt_1}{t_1} \cdot \frac{dt_2}{1-t_2}}_{2} \cdot \underbrace{ \frac{dt_3}{t_3}\cdot  \frac{dt_4}{t_4} \cdot \frac{dt_5}{1-t_5}}_{3} \,. \]
    This corresponds to $I(2,3)$ (but is of course not the same since the $I$ are formal symbols).} $I(k_1,\ldots,k_r)$ to denote $I_0(k_1,\ldots,k_r)$.
\begin{proposition}(\cite[Eq. (3.5),(3.6) and Prop. 3.5]{BT})\label{prop:basisfori}
    \begin{enumerate}[\textup{(}i\textup{)}] 
        \item We have $I_n(\emptyset)=0$ if $n\ge1$ or $1$ if $n=0$.
        \item For integers $n\ge0,k_1,\ldots,k_r\ge1$, 
        \begin{equation*}
        I_{n}(k_1,\ldots,k_r) = (-1)^n \sum^*\bigg(\prod_{j=1}^{r}\binom{l_j-1}{k_j-1} \bigg) I (l_1,\ldots,l_r) \,,
        \end{equation*}
        where the sum runs over all $l_1+\cdots+l_r=k_1+\cdots +k_r+n$ with $l_1,\ldots,l_r\ge1$.
        \item The set $\{I(k_1,\ldots,k_r)\mid r\ge0,k_i\ge1\}$ forms a basis of the space $\mathcal{I}^1$.
    \end{enumerate}
\end{proposition}

We give an example for (ii). In $\mathcal{I}^1$ it is $I(1;0;0)=0$ and therefore
\begin{align*}
0 &= I(1;0;0)  I(1;0,1;0) \\
&=I(1;0,0,1;0) + I(1;0,0,1;0) + I(1;0,1,0;0) \\
&= 2 I(3) + I_1(2) 
\end{align*}
which gives $I_1(2) = -2 I(3) = (-1)^1 \binom{2}{1} I(3)$. 

\begin{remark} \label{rmk:ih}
    Statement iii) in Proposition \ref{prop:basisfori} basically states that we can identify $\mathcal{I}^1$ with $\h^1$ by sending $I(k_1,\dots,k_r)$ to $z_{k_1}\dots z_{k_r}$, which is an algebra isomorphism with respect to the shuffle product. In other words, we can equip $\h^1$ with the coproduct $\Delta_G$. Instead of working with $I$ we will use this identification in the next section when defining the shuffle regularized multiple Eisenstein series. 
\end{remark}

\begin{ex}\label{ex:coproductfourier32}
    In the following we are going to calculate $\Delta_G(I(3,2)) = \Delta_G(I(1;0,0,1,0,1;0))$. Therefore we have to determine all possible markings of the diagram
    \begin{figure}[H]
        \centering
        \begin{tikzpicture}[scale=0.60]
        \def \ra {2}
        \def \circsize {0.1}		
        \draw (-\ra,0) -- (\ra,0)
        (\ra,0) arc (0:180:\ra);
        \fill[black] (180:\ra) circle (\circsize); 
        
        \fill[black] (30:\ra) circle (\circsize); 
        \draw (30:\ra) circle (\circsize); 	
        
        \fill[white] (60:\ra) circle (\circsize); 
        \draw (60:\ra) circle (\circsize); 	
        
        \fill[black] (90:\ra) circle (\circsize); 
        
        \fill[white] (120:\ra) circle (\circsize); 
        \draw (120:\ra) circle (\circsize); 
        
        \fill[white] (150:\ra) circle (\circsize);	
        \draw (150:\ra) circle (\circsize); 
        
        \fill[white] (0:\ra) circle (\circsize); 	
        \draw (0:\ra) circle (\circsize); 
        \end{tikzpicture}
    \end{figure}where the corresponding summand in the coproduct does not vanish. For simplicity we draw  $\circ$ to denote a $0$ and $\bullet$ to denote a $1$. We will consider the $4 = 2^2$ ways of marking the two $\bullet$ in the top part of the circle separately. As explained in Section \ref{subsec:fourier}, we want to compare the coproduct to the Fourier expansion of multiple Eisenstein series. Therefore, in this case we also calculate the expansion of $\aG_{3,2}(\tau)$ using the construction described there. Recall that we also had the $4$ different parts $\aG^{RR}_{3,2}$, $\aG^{UR}_{3,2}$, $\aG^{RU}_{3,2}$ and $\aG^{UU}_{3,2}$. We will see that the number and positions of the marked $\bullet$ correspond to the number and positions of the letter $U$ in the word $w$ of $\aG^w$. 
    \begin{enumerate}[\textup{(}i\textup{)}] 
        \item Diagrams with no marked $\bullet$: 
        \begin{figure}[H]
            \centering
            \begin{tikzpicture}[scale=0.6]
            \def \ra {2}
            \def \circsize {0.1}	
            \draw  (0:\ra) to [bend right] (180:\ra) ;			
            
            \draw (-\ra,0) -- (\ra,0)
            (\ra,0) arc (0:180:\ra);
            \draw (-\ra,0) -- (\ra,0)
            (\ra,0) arc (0:180:\ra);
            \fill[black] (180:\ra) circle (\circsize); 
            
            \fill[black] (30:\ra) circle (\circsize); 
            \draw (30:\ra) circle (\circsize); 	
            
            \fill[white] (60:\ra) circle (\circsize); 
            \draw (60:\ra) circle (\circsize); 	
            
            \fill[black] (90:\ra) circle (\circsize); 
            
            \fill[white] (120:\ra) circle (\circsize); 
            \draw (120:\ra) circle (\circsize); 
            
            \fill[white] (150:\ra) circle (\circsize);	
            \draw (150:\ra) circle (\circsize); 
            
            \fill[white] (0:\ra) circle (\circsize); 	
            \draw (0:\ra) circle (\circsize); 
            \end{tikzpicture}
        \end{figure}
        Corresponding sum in the coproduct: 
        \[ I(0;\emptyset;1) \otimes I(1;0,0,1,0,1;0) = 1 \otimes I(3,2) \,. \] 
        The part of the Fourier expansion of $\aG_{3,2}$ which is associated to this is the one with no $U$ ``occurring'', i.e.  $\aG^{RR}_{3,2}(\tau) = \zeta(3,2)$.  
        \item  Diagrams with the first $\bullet$ marked: 
        \begin{figure}[H]
            \centering
            \begin{tikzpicture}[scale=0.6]
            \def \ra {2}
            \def \circsize {0.1}

            \draw  (0:\ra) to [bend left] (90:\ra) ;	
            \draw  (90:\ra) to [bend left] (120:\ra) ;	
            \draw  (120:\ra) to [bend left] (150:\ra) ;	
            \draw  (150:\ra) to [bend left] (180:\ra) ;				
            
            \draw (-\ra,0) -- (\ra,0)
            (\ra,0) arc (0:180:\ra);
            \fill[black] (180:\ra) circle (\circsize); 
            
            \fill[black] (30:\ra) circle (\circsize); 
            \draw (30:\ra) circle (\circsize); 	
            
            \fill[white] (60:\ra) circle (\circsize); 
            \draw (60:\ra) circle (\circsize); 	
            
            \fill[black] (90:\ra) circle (\circsize); 
            
            \fill[white] (120:\ra) circle (\circsize); 
            \draw (120:\ra) circle (\circsize); 
            
            \fill[white] (150:\ra) circle (\circsize);	
            \draw (150:\ra) circle (\circsize); 
            
            \fill[white] (0:\ra) circle (\circsize); 	
            \draw (0:\ra) circle (\circsize); 
            
            \end{tikzpicture}
        \end{figure}
        Corresponding sum in the coproduct: 
        \begin{align*}
        I(1;0,0,1;0) \otimes \big( I(1;0) \cdot I(0;0) \cdot I(0;1) \cdot I(1;0,1;0) \big) = I(3) \otimes I(2) \,.
        \end{align*}
        The associated part of the Fourier expansion of $\aG_{3,2}$ is $\aG^{UR}_{3,2}(\tau) = \ag(3) \cdot \zeta(2)$.  
        \item Diagrams with the second $\bullet$ marked: 
        \begin{figure}[H]
            \centering
            \begin{tikzpicture}[scale=0.6]
            \def \ra {2}
            \def \circsize {0.1}

            \draw  (0:\ra) to [bend left] (30:\ra) ;			
            \draw  (30:\ra) to [bend left] (60:\ra) ;	
            \draw  (60:\ra) to [bend left] (180:\ra) ;				
            
            \draw (-\ra,0) -- (\ra,0)
            (\ra,0) arc (0:180:\ra);
            \fill[black] (180:\ra) circle (\circsize); 
            
            \fill[black] (30:\ra) circle (\circsize); 
            \draw (30:\ra) circle (\circsize); 	
            
            \fill[white] (60:\ra) circle (\circsize); 
            \draw (60:\ra) circle (\circsize); 	
            
            \fill[black] (90:\ra) circle (\circsize); 
            
            \fill[white] (120:\ra) circle (\circsize); 
            \draw (120:\ra) circle (\circsize); 
            
            \fill[white] (150:\ra) circle (\circsize);	
            \draw (150:\ra) circle (\circsize); 
            
            \fill[white] (0:\ra) circle (\circsize); 	
            \draw (0:\ra) circle (\circsize); 
            
            \end{tikzpicture}
            \,\,\,\,\,\,\,\,
            \begin{tikzpicture}[scale=0.6]
            \def \ra {2}
            \def \circsize {0.1}

            \draw  (0:\ra) to [bend left] (30:\ra) ;			
            \draw  (30:\ra) to [bend left] (150:\ra) ;	
            \draw  (150:\ra) to [bend left] (180:\ra) ;				
            
            \draw (-\ra,0) -- (\ra,0)
            (\ra,0) arc (0:180:\ra);
            \fill[black] (180:\ra) circle (\circsize); 
            
            \fill[black] (30:\ra) circle (\circsize); 
            \draw (30:\ra) circle (\circsize); 	
            
            \fill[white] (60:\ra) circle (\circsize); 
            \draw (60:\ra) circle (\circsize); 	
            
            \fill[black] (90:\ra) circle (\circsize); 
            
            \fill[white] (120:\ra) circle (\circsize); 
            \draw (120:\ra) circle (\circsize); 
            
            \fill[white] (150:\ra) circle (\circsize);	
            \draw (150:\ra) circle (\circsize); 
            
            \fill[white] (0:\ra) circle (\circsize); 	
            \draw (0:\ra) circle (\circsize); 
            
            \end{tikzpicture}
            \,\,\,\,\,\,\,\,
            \begin{tikzpicture}[scale=0.6]
            \def \ra {2}
            \def \circsize {0.1}

            \draw  (0:\ra) to [bend left] (30:\ra) ;	
            \draw  (30:\ra) to [bend left] (120:\ra) ;	
            \draw  (120:\ra) to [bend left] (150:\ra) ;	
            \draw  (150:\ra) to [bend left] (180:\ra) ;				
            
            \draw (-\ra,0) -- (\ra,0)
            (\ra,0) arc (0:180:\ra);
            \fill[black] (180:\ra) circle (\circsize); 
            
            \fill[black] (30:\ra) circle (\circsize); 
            \draw (30:\ra) circle (\circsize); 	
            
            \fill[white] (60:\ra) circle (\circsize); 
            \draw (60:\ra) circle (\circsize); 	
            
            \fill[black] (90:\ra) circle (\circsize); 
            
            \fill[white] (120:\ra) circle (\circsize); 
            \draw (120:\ra) circle (\circsize); 
            
            \fill[white] (150:\ra) circle (\circsize);	
            \draw (150:\ra) circle (\circsize); 
            
            \fill[white] (0:\ra) circle (\circsize); 	
            \draw (0:\ra) circle (\circsize); 
            
            \end{tikzpicture}
        \end{figure}
        Corresponding sum in the coproduct: 
        \begin{align*}
        & I(1;0,1;0) \otimes \big(  I(1;0,0,1;0) \cdot I(0;1) \cdot I(1;0)\big) \\
        + & I(1;0,1;0) \otimes \big(  I(1;0) \cdot I(0;0,1,0;1) \cdot I(1;0)\big) \\
        + & I(1;0,0,1;0) \otimes \big(  I(1;0) \cdot I(0;0) \cdot I(0;1,0;1) \cdot I(1;0)\big) \\
        & = I(2) \otimes I(3) - I(2) \otimes I_1(2) + I(3) \otimes I(2) \,,
        \end{align*}
        where we used $I(0;0,1,0;1)=-I_1(2)$ and $I(0;1,0;1) = (-1)^2 I(1;0,1;0) = I(2)$. Together with $I_1(2) = -2 I(3)$ this gives
        \[ 3 I(2) \otimes I(3) + I(3) \otimes I(2) \,.\]
        The associated part of the Fourier expansion is also the most complicated one. We have
        \[ \aG^{RU}_{3,2}(\tau) = \sum_{m>0} \Psi_{3,2}(m\tau) \]
        and with \eqref{eq:psi32} we derived $\Psi_{3,2}(x) = 3 \Psi_2(x) \cdot \zeta(3) + \Psi_3(x) \cdot \zeta(2)$, i.e. 
        \[ \aG^{RU}_{3,2}(\tau) = 3 \ag(2) \cdot \zeta(3) + \ag(3) \cdot \zeta(2) \,.\] 
        \item  Diagrams with both $\bullet$ marked: 
        \begin{figure}[H]
            \centering
            \begin{tikzpicture}[scale=0.6]
            \def \ra {2}
            \def \circsize {0.1}	
            
            \draw  (0:\ra) to [bend left] (30:\ra) ;	
            \draw  (30:\ra) to [bend left] (60:\ra) ;	
            \draw  (60:\ra) to [bend left] (90:\ra) ;	
            \draw  (90:\ra) to [bend left] (120:\ra) ;				
            \draw  (120:\ra) to [bend left] (150:\ra) ;	
            \draw  (150:\ra) to [bend left] (180:\ra) ;					
            
            \draw (-\ra,0) -- (\ra,0)
            (\ra,0) arc (0:180:\ra);
            \fill[black] (180:\ra) circle (\circsize); 
            
            \fill[black] (30:\ra) circle (\circsize); 
            \draw (30:\ra) circle (\circsize); 	
            
            \fill[white] (60:\ra) circle (\circsize); 
            \draw (60:\ra) circle (\circsize); 	
            
            \fill[black] (90:\ra) circle (\circsize); 
            
            \fill[white] (120:\ra) circle (\circsize); 
            \draw (120:\ra) circle (\circsize); 
            
            \fill[white] (150:\ra) circle (\circsize);	
            \draw (150:\ra) circle (\circsize); 
            
            \fill[white] (0:\ra) circle (\circsize); 	
            \draw (0:\ra) circle (\circsize); 
            
            \end{tikzpicture}
        \end{figure}
        Corresponding sum in the coproduct: $ I(3,2) \otimes 1$. The associated part of the Fourier expansion of $\aG_{3,2}$ is $\aG^{UU}_{3,2}(\tau) = \ag(3,2)$.  
    \end{enumerate}
    Summing all $4$ parts together we obtain for the coproduct
    \[ \Delta_G( I(3,2) ) =  1 \otimes I(3,2)   +  3 I(2) \otimes I(3) + 2 I(3) \otimes I(2)+  I(3,2)  \otimes 1  \]
    and for the Fourier expansion of $\aG_{3,2}(\tau)$: 
    \begin{align*}
    \aG_{3,2}(\tau) =\zeta(3,2) + 3  \ag(2) \zeta(3) + 2  \ag(3) \zeta(2) + \ag(3,2) \,. 
    \end{align*}
\end{ex}

This shows that the left factors of the terms in the coproduct correspond to the functions $\g$ and the right factors to the multiple zeta values. We will use this in the next section to define shuffle regularized multiple Eisenstein series. 

\subsection{The $q$-series $\g^\sh$}

In this section we want to construct for $k_1,\dots,k_r\geq 1$ $q$-series $\g^\sh(k_1,\dots,k_r)$ which satisfy the shuffle product formula. By this we mean that the following map is a $\Q$-algebra homomorphism
\begin{align*}
\g^\sh : \h^1_\sh &\longrightarrow \Q[[q]]\\
z_{k_1}\dots z_{k_r} &\longmapsto \g^\sh(k_1,\dots,k_r)\,.
\end{align*}	 
In smallest depths this means that we have for $k_1,k_2\geq 1$ and $k=k_1+k_2$
\begin{align}\label{eq:gsh_profdep1}
\g^\sh(k_1) \g^\sh(k_2) = &\sum_{j=1}^{k-1} \left( \binom{j-1}{k_1-1} + \binom{j-1}{k_2-1} \right) \g^\sh(j, k-j) \,.
\end{align}

In Section \ref{subsec:qmzv} we saw in Proposition \ref{prop:gshuffle1} that for $k_1,k_2\geq1$ and $k=k_1+k_2>2$ the $q$-series $\g$ satisfy a similar formula
\begin{align}\label{eq:gshuffle1_again}\begin{split}
\g(k_1) \g(k_2) = &\sum_{j=1}^{k-1} \left( \binom{j-1}{k_1-1} + \binom{j-1}{k_2-1} \right) \g(j, k-j) \\
&+ \binom{k-2}{k_1-1} \left(q\frac{d}{dq} \frac{\g(k-2)}{k-2} - \g(k-1) \right)\,.
\end{split}
\end{align}
For the remaining case $k_1=k_2=1$ we have
\begin{align}\label{eq:gshuffle11_again}
\g(1)^2=2\g(1,1)+\g(2)-\g(1)\,.
\end{align}
Indeed, we will see that the $\g^\sh$ are given by the $\g$ if all indices are greater than $1$. Comparing \eqref{eq:gsh_profdep1} with \eqref{eq:gshuffle1_again} and \eqref{eq:gshuffle11_again} we see that we could define the $\g^\sh$ in depths one and two by
\begin{align*}
\g^\sh(k) &= \g(k)\,,\\
\g^\sh(k_1,k_2) &=
\begin{cases}
\displaystyle \g(k_1,k_2),&k_2>1,\\[1mm]
\displaystyle \g(k_1,1) + \frac{1}{2} \left( q\frac{d}{dq} \frac{\g(k_1-1)}{k_1-1}  - \g(k_1) \right),&k_2=1,\ k_1\geq2,\\[2mm]
\displaystyle \g(1,1)+\frac{1}{2}\bigl(\g(2)-\g(1)\bigr),&(k_1,k_2)=(1,1).
\end{cases}
\end{align*}
These series then satisfy the equation \eqref{eq:gsh_profdep1}.
For higher depths the derivatives $ q\frac{d}{dq} $ will not be sufficient to correct the $q$-series $\g$ in order to get $\g^\sh$. We will need to introduce the following double-indexed version of $\g$, on which we will focus in more detail in Section \ref{subsec:doubleindexedqana}. 
\begin{definition}\label{def:doubleg}For $k_1,\dots,k_r \geq 1$, $d_1,\dots,d_r\geq 0$ we define the $q$-series
    \begin{align*}
    \mb{k_1,\dots,k_r}{d_1,\dots,d_r} &:= \sum_{m_1 > \cdots > m_r > 0} m_1^{d_1}\frac{P_{k_1}(q^{m_1})}{(1-q^{m_1})^{k_1}} \cdots m_r^{d_r} \frac{P_{k_r}(q^{m_r})}{(1-q^{m_r})^{k_r}} \\
    &= \sum_{\substack{m_1>\dots>m_r > 0\\n_1,\dots,n_r > 0}} m_1^{d_1} \frac{n_1^{k_1-1}}{(k_1-1)!} \cdots m_r^{d_r} \frac{n_r^{k_r-1}}{(k_r-1)!}  q^{m_1 n_1 + \dots + m_r n_r}\,.
    \end{align*}
\end{definition}
Notice that the original definition in \cite{BK1} included the factors $\frac{1}{d_1!\cdots d_r!}$.

These generalize the $q$-series $\g$, since by definition we have 
\begin{align*}
\mb{k_1,\dots,k_r}{0,\dots,0} = \g(k_1,\dots,k_r)\,.
\end{align*}
In depth one we can also see that these give the derivatives with respect to $q\frac{d}{dq}$, since
\begin{align*}
q \frac{d}{dq} \g(k) = q \frac{d}{dq}  \sum_{\substack{m>0\\n>0}} \frac{n^{k-1}}{(k-1)!} q^{m n} = \sum_{\substack{m>0\\n>0}} \frac{m n^{k}}{(k-1)!} q^{m n} =  k \sum_{\substack{m>0\\n>0}} \frac{m n^{k}}{k!} q^{m n}=  k \mb{k+1}{1}\,.
\end{align*}

We will deal with the operator $q \frac{d}{dq}$ in more generality in Section \ref{subsec:doubleindexedqana}.
There we will also study the algebraic structure of these double-indexed $\g$ and prove that they span the space of modified $q$-analogues, i.e. the space $\mz_q$ of modified $q$-analogues (Eq. \eqref{eq:defzq}) is spanned by the $q$-series $\mb{k_1,\dots,k_r}{d_1,\dots,d_r}$ (see Theorem \ref{thm:zqspannedbydoubleg}).	

We now want to construct the $\g^\sh$ in general as elements in $\mz_q$ by defining them in terms of the double-indexed $\g$. This will be done by using their generating series
\begin{align*}
\ggen\bi{X_1,\dots,X_r}{Y_1,\dots,Y_r} := \sum_{\substack {k_1,\dots,k_r\geq  1\\d_1,\dots,d_r \geq 0}} \mb{k_1,\dots,k_r}{d_1,\dots,d_r} X_1^{k_1-1} \frac{Y_1^{d_1}}{d_1!}\cdots X_r^{k_r-1} \frac{Y_r^{d_r}}{d_r!}\,.
\end{align*}	
Similar to Lemma \ref{lem:genexpr} we obtain the following explicit expression.
\begin{lemma} \label{lem:bigenexpr}We have 
    \begin{align}\label{eq:bigen}
    \ggen\bi{X_1,\dots,X_r}{Y_1,\dots,Y_r}
    &= \sum_{m_1> \dots > m_r > 0} e^{Y_1 m_1} \frac{e^{X_1} q^{m_1}}{1-e^{X_1}q^{m_1}} \cdots  e^{Y_r m_r }\frac{e^{X_r} q^{m_r}}{1-e^{X_r}q^{m_r}} \,.
    \end{align}
\end{lemma}

For $e_1,\dots,e_r \geq 1$ we generalize these generating series to the following 
\begin{align}\label{def:tribracket} 
\mtt{X_1, & ... & ,\, X_r}{Y_1, & ... & ,\,  Y_r}{e_1, & ... & ,\, e_r} = \sum_{m_1>\dots>m_r>0} \prod_{j=1}^r e^{m_j Y_j} \left(\frac{e^{X_j} q^{m_j}}{1-e^{X_j} q^{m_j}} \right)^{e_j}  \,. 
\end{align}
In particular for $e_1 = \dots = e_r=1$ these are the generating series of the $\g$ in \eqref{eq:bigen}. To show that the coefficients of these series are in $\mz_q$ for arbitrary $e_j$ we need to define the differential operator $\mathcal{D}^Y_{e_1,\dots,e_r} := D_{Y_1,e_1}  D_{Y_2,e_2} \dots  D_{Y_r,e_r}$ with
\begin{align*}
D_{Y_j,e} &=  \prod_{k=1}^{e-1} \left(\frac{1}{k}\left(\frac{\partial}{\partial Y_{r-j+1}} - \frac{\partial}{\partial Y_{r-j+2}} \right) - 1\right)\,.
\end{align*}
where we set $\frac{\partial}{\partial Y_{r+1}}=0$.

\begin{proposition}\label{prop:generalpartition}
    The coefficients of \eqref{def:tribracket} are in $\mz_q$ and we have
    \[\mathcal{D}^Y_{e_1,\dots,e_r}\mt{X_1, \dots, X_r}{Y_1,\dots,Y_r}  = \mtt{Y_1+\dots+Y_r,\, & ... & ,\, Y_1}{X_r,\, X_{r-1}-X_r,\,& ... & ,\,  X_1-X_2}{e_1, & ... & ,\, e_r} \,. \]
\end{proposition}
\begin{proof}
    By $\frac{\partial}{\partial X} L_n(X) = L_n(X)^2 + L_n(X)$ one inductively obtains 
    \[  L_n(Y)^{e+1}  =  \left( \frac{1}{e} \frac{\partial}{\partial Y}  -   1\right) L_n(Y)^e \,, \qquad\text{and therefore}\qquad  L_n(Y)^{e} = \prod_{k=1}^{e-1} \left( \frac{1}{k} \frac{\partial}{\partial Y}  -   1 \right) L_n(Y) \,, \]
    from which the statement follows after a suitable change of variables. 
\end{proof}

\begin{lemma} \label{lem:shuffleprodforgenseries}
    Let $A$ be an algebra spanned by elements $a_{k_1,\dots,k_r}$ with $k_1,\dots, k_r \geq 1$, and let
    \[ H(X_1,\dots,X_r) = \sum_{k_j} a_{k_1,\dots,k_r} X_1^{k_1-1} \dots X_r^{k_r-1} \]
    be the generating function of these elements. For $f \in \Q[[X_1,\dots,X_r]]$ define
    \[ f^{\sharp}(X_1,\dots,X_r) = f(X_1+\dots+X_r, X_2+\dots+X_r, \dots, X_r) \,. \] 
    Then the following two statements are equivalent.
    \begin{enumerate}[\textup{(}i\textup{)}] 
        \item The map $\h^1_{\sh} \rightarrow A$ given by $z_{k_1} \dots z_{k_r} \mapsto a_{k_1,\dots,k_r}$ is an algebra homomorphism. 
        \item For all $r,s \geq 1$  we have
        \begin{multline*}
        H^{\sharp}(X_1,\dots,X_r) \cdot H^{\sharp}(X_{r+1}, \dots, X_{r+s}) \\
        = H^{\sharp}(X_1,\dots,X_{r+s})_{\vert sh_r^{(r+s)}}\,,
        \end{multline*}
        where $sh_r^{(r+s)} = \sum_{\sigma\in\Sigma(r,s)}\sigma$ in the group ring $\Z[\mathfrak{S}_{r+s}]$ and the symmetric group $\mathfrak{S}_r$ acts on $\Q[[X_1,\ldots,X_r]]$ by $(f\big|\sigma)(X_1,\ldots,X_r)= f(X_{\sigma^{-1}(1)},\ldots,X_{\sigma^{-1}(r)})$\,.
    \end{enumerate}
\end{lemma}
\begin{proof} This can be proven by induction over $r$ together with Proposition 8 in \cite{I}.\end{proof}

\begin{definition}\label{def:gshuffle}
    For $k_1,\dots,k_r\geq 1$ define $\g^\sh(k_1,\dots,k_r) \in \mz_q$ as the coefficients of the following generating function
    \begin{align*}
    &H_{\sh}(X_1,\dots,X_r)  = \sum_{k_1,\dots,k_r \geq 1} \g^\sh(k_1,\dots,k_r) X_1^{k_1-1} \dots X_r^{k_r-1} \\
    &:=\sum_{\substack{1\leq m \leq r\\i_1 + \dots + i_m = r}}\frac{1}{i_1! \dots i_m!} \mathcal{D}^Y_{i_1,\dots,i_m} \mt{X_1,X_{i_m+1},X_{i_{m-1}+i_m+1}, \dots, X_{i_2+\dots+i_m+1}}{Y_1,\dots,Y_r} _{\big\vert Y=0} \,.
    \end{align*}
\end{definition}

\begin{theorem}(\cite[Thm. 5.7]{Ba7}) \label{thm:shufflebracket} We have
    \begin{enumerate}[\textup{(}i\textup{)}] 
        \item The $\g^\sh(k_1,\dots,k_r)$ satisfy the shuffle product formula, i.e. 
        \[ H^{\sharp}_{\sh}(X_1,\dots,X_r) \cdot H^{\sharp}_{\sh}(X_{r+1}, \dots, X_{r+s}) = H^{\sharp}_{\sh}(X_1,\dots,X_{r+s})_{\vert sh_r^{(r+s)}} \,. \]
        \item For $k_1\geq 1,\, k_2,\dots,k_r \geq 2$ we have $\g^\sh(k_1,\dots,k_r) = \g(k_1,\dots,k_r)$.
    \end{enumerate}
\end{theorem} 
\begin{proof}
    The first part of the proof is basically the same as in the discussion in Section 4.1 in \cite{BT} but with a reverse order and some changes in notation.
    Consider the alphabet $A=\left\{ \bi{y}{n} \mid n \geq 1 \,, y \in Y_\Z  \right\}$, where $Y_{\Z}$ is the set of finite sums of the elements in $Y=\{Y_1, Y_2,\dots\}$. We denote a word in these letters by $\bi{y_1,\dots,y_r}{n_1,\dots,n_r}$. For two letters $a, b \in A$ define $a \diamond b \in A$ as the component-wise sum. With this we can equip $\Q\langle A \rangle$ with the quasi-shuffle product $\qsh$ (see Section \ref{subsec:quasishuffle}) and therefore obtain a quasi-shuffle algebra $( \Q\langle A \rangle , \qsh )$. It is easy to see that the map $( \Q\langle A \rangle , \qsh ) \rightarrow \varinjlim_j \Q[[q]][[X_1,\dots,X_j,Y_1,\dots,Y_j]]$ given by
    \[ \bi{y_1,\dots,y_r}{n_1,\dots,n_r} \longmapsto  \mtt{0, \dots  ,\, 0}{y_1, \dots ,\,  y_r}{n_1, \dots,\, n_r} \] 
    is an algebra homomorphism. Using now Theorem \ref{thm:explog} the series $h$ defined by the exponential map
    \[  h(X_1,\dots,X_r) =\sum_{\substack{1\leq m \leq r\\ i_1+\dots+i_m = r}} \frac{1}{i_1! \dots i_m!} \mtt{0, & ... & ,\, 0}{Y_1, & ... & ,\,  Y_m}{i_1, & ... & ,\, i_m} \,, \] 
    where $s_0=0$, $s_j=i_1+\dots+i_j$ and $Y_j=X_{s_{j-1}+1}+\dots+X_{s_j}$, satisfies the (index-)shuffle product i.e.
    \[ h(X_1,\dots,X_r) \cdot h(X_{r+1}, \dots, X_{r+s}) = h(X_1,\dots,X_{r+s})_{\vert sh_r^{(r+s)}} \,.\]
    We now set $H_{\sh}(X_1,\dots,X_r) := h(X_r,X_{r-1}-X_r,\dots, X_1-X_2)$ and by the same argument as in Theorem 4.3 in \cite{BT} it is 
    \[ H^{\sharp}_{\sh}(X_1,\dots,X_r) \cdot H^{\sharp}_{\sh}(X_{r+1}, \dots, X_{r+s}) = H^{\sharp}_{\sh}(X_1,\dots,X_{r+s})_{\vert sh_r^{(r+s)}} \,. \]
    Combining the definition of $h$ and $H_{\sh}$ we observe that $H_{\sh}(X_1,\dots,X_r)$ equals
    \[ \sum_{\substack{1\leq m \leq r\\ i_1+\dots+i_m = r}} \frac{1}{i_1! \dots i_m!} \mtt{0, & ... & ,\, 0}{X_{r-i_1+1},\, X_{r-i_1-i_2+1} - X_{r - i_1 + 1} , &...& ,\,  X_1 - X_{r-i_1-\dots-i_{m-1}+1}}{i_1, & ... & ,\, i_m} \,.\] 
    We now apply Proposition \ref{prop:generalpartition} to this and obtain (i) of the Theorem. To prove (ii) one checks that the only summand on the right hand side, where {\bf all} variables $X_2,\dots,X_r$ appear, is the one with $i_1= \dots = i_m=1$ which is exactly $\g(k_1,\dots,k_r) X_1^{k_1-1} \dots X_r^{k_r-1}$. Therefore the $q$-series $\g^\sh(k_1,\dots,k_r)$ where $k_2,\dots,k_r\geq 2$ are given by $\g(k_1,\dots,k_r)$. 
\end{proof}

For small depths we obtain the following explicit expressions for $\g^\sh$.
\begin{corollary}%
    We have $\g^\sh(k_1) = \g(k_1)$ and for $r=2,3,4$  the $\g^\sh(k_1,\dots,k_r)$ are given by\footnote{Here $\delta_{a,b}$ again denotes the Kronecker delta, i.e. $\delta_{a,b}$ is $1$ for $a=b$ and $0$ otherwise.}
    \begin{enumerate}[\textup{(}i\textup{)}] 
        \item $\begin{aligned}[t]
        \g^\sh(k_1,k_2) &= \g(k_1,k_2) + \delta_{k_2,1}\cdot \frac{1}{2} \left( \mb{k_1}{1} - \g(k_1) \right) \,,\\
        \end{aligned}$
        \item $\begin{aligned}[t]
        \g^\sh(k_1,k_2,k_3) &= \g(k_1,k_2,k_3) + \delta_{k_3,1}\cdot \frac{1}{2} \left( \mb{k_1,k_2}{0,1} - \g(k_1,k_2) \right)\\
        &+\delta_{k_2,1}\cdot \frac{1}{2}\left(  \mb{k_1,k_3}{1,0} - \mb{k_1,k_3}{0,1}- \g(k_1,k_3) \right)\\
        &+\delta_{k_2 \cdot k_3,1}\cdot \frac{1}{6}\left(  \frac{1}{2}\mb{k_1}{2} - \frac{3}{2}\mb{k_1}{1}+ \g(k_1) \right) \,,\\
        \end{aligned}$
        \item \resizebox{\linewidth}{!}{$\begin{aligned}[t]
        \g^\sh(k_1,k_2,k_3&,k_4) = \g(k_1,k_2,k_3,k_4) + \delta_{k_4,1}\cdot \frac{1}{2} \left( \mb{k_1,k_2,k_3}{0,0,1} - \g(k_1,k_2,k_3) \right) \\
        + \delta_{k_3,1} &\cdot\frac{1}{2} \left( \mb{k_1,k_2,k_4}{0,1,0}-\mb{k_1,k_2,k_4}{0,0,1} - \g(k_1,k_2,k_4) \right)\\
        + \delta_{k_2,1} &\cdot\frac{1}{2} \left( \mb{k_1,k_3,k_4}{1,0,0}-\mb{k_1,k_3,k_4}{0,1,0} - \g(k_1,k_3,k_4) \right)\\
        + \delta_{k_2 \cdot k_4,1} &\cdot\frac{1}{4} \left( \mb{k_1,k_3}{1,1}-\mb{k_1,k_3}{0,2}-\mb{k_1,k_3}{1,0} + \g(k_1,k_3) \right)\\
        + \delta_{k_3 \cdot k_4,1} &\cdot\frac{1}{6} \left( \frac{1}{2}\mb{k_1,k_2}{0,2}-\frac{3}{2}\mb{k_1,k_2}{0,1} + \g(k_1,k_2) \right)\\
        + \delta_{k_2 \cdot k_3,1} &\cdot\frac{1}{6} \left( \frac{1}{2}\mb{k_1,k_4}{0,2}-\mb{k_1,k_4}{1,1}+\frac{3}{2}\mb{k_1,k_4}{0,1}+\frac{1}{2}\mb{k_1,k_4}{2,0}-\frac{3}{2}\mb{k_1,k_4}{1,0} + \g(k_1,k_4) \right)\\
        + \delta_{k_2 \cdot k_3 \cdot k_4,1} &\cdot\frac{1}{24} \left( \frac{1}{6}\mb{k_1}{3}-\mb{k_1}{2}+\frac{11}{6} \mb{k_1}{1} - \g(k_1) \right) \,.
        \end{aligned}$}
    \end{enumerate}
\end{corollary}
\begin{proof}
    This follows by calculating the coefficients of the series $H_\sh$.
\end{proof}

\subsection{Shuffle regularized multiple Eisenstein series}%
In this section, we present the definition of shuffle regularized multiple Eisenstein series as was done in \cite{BT}. We use the observation of the previous section and use the coproduct $\Delta_G$ of formal iterated integrals to define these series. As mentioned in Remark \ref{rmk:ih} we can equip the space $\h^1$ with the coproduct $\Delta_G$ instead of working with the space  $\mathcal{I}^1$. In analogy to the map $\zeta^\sh: \h^1_\sh\rightarrow \mz$ of shuffle regularized multiple zeta values, the map  $$\ag^{\sh} : \h^1_\sh \rightarrow \Q[2\pi i][\![q]\!]$$ defined on the generators $z_{k_1} \dots z_{k_r}$ by 
\[ \ag^{\sh}(z_{k_1} \dots z_{k_r}) := (-2\pi i)^{k_1+\dots+k_r}\g^\sh( k_1,\dots,k_r)\,, \] 
is also an algebra homomorphism by Theorem \ref{thm:shufflebracket}. With this we can give the definition of $\aG^\sh$ (\cite{BT}).
\begin{definition}
    For $k_1,\ldots,k_r\ge1$ define the  \emph{shuffle regularized multiple Eisenstein series} by 
    \[ \aG^\sh_{k_1,\ldots,k_r}(\tau) := m\left( (  \ag^{\sh} \otimes  \zeta^{\sh})\circ \Delta_G \big(z_{k_1} \dots z_{k_r} \big)\right)\,,\]
    where $m$ denotes the multiplication given by $m: a \otimes b \mapsto a\cdot b$ and $\zeta^\sh$ denotes shuffle regularized multiple zeta values (Definition \ref{def:regmzv} with $T=0$).
\end{definition}
We can view $\aG^\sh$ as an algebra homomorphism $\aG^\sh: \h^1_\sh\rightarrow \MZB$ such that the following diagram commutes
\[
\xymatrix{  
    \h^1_\sh  \ar@{->}[r]^-{\Delta_G} \ar@{->}[d]_{\aG^\sh} & \h^1_\sh \otimes \h^1_\sh \ar@{->}[d]^{\ag^{\sh} \otimes  \zeta^{\sh} }  \\ 
    \MZB &   \Q[2\pi i][\![q]\!] \otimes \,\mz  \ar@{->}[l]^-{m}\\
} 
\]
\begin{theorem}(\cite[Thm. 1.1, 1.2]{BT}) \label{thm:messh}
    For all $k_1,\ldots,k_r\ge1$ the shuffle regularized multiple Eisenstein series $\aG^{\sh}_{k_1,\ldots,k_r}$ have the following properties:
    \begin{enumerate}[\textup{(}i\textup{)}] 
        \item They are holomorphic functions on the upper half-plane having a Fourier expansion with the shuffle regularized multiple zeta values as the constant term, i.e. they can be written as
        \begin{align}\label{eq:messhsimplefourier}
        \aG^\sh_{k_1,\dots,k_r}(\tau) = \zeta^\sh(k_1,\dots,k_r) + \sum_{n=1}^\infty a_n q^n\,.\qquad  \qquad (a_n \in \mz[\pi i] = \mz + \pi i \mz)
        \end{align}
        \item They satisfy the shuffle product formula.
        \item For integers $k_1,\ldots,k_r\ge2$ they equal the multiple Eisenstein series
        \[ \aG^{\sh}_{k_1,\ldots,k_r}(\tau)=\aG_{k_1,\ldots,k_r}(\tau) \]
        and therefore satisfy the stuffle product formula in these cases. 
    \end{enumerate} 
\end{theorem}
Parts i) and ii) in this theorem follow directly by definition. The important part here is iii), which states that the connection of the Fourier expansion and the coproduct, as illustrated in Example \ref{ex:coproductfourier32}, holds in general. It also shows that the shuffle regularized multiple Eisenstein series satisfy the stuffle product formula in many cases, though the exact failure of the stuffle product of these series is still unknown. 

\subsection{Stuffle regularized multiple Eisenstein series}\label{sec:stufflereg}
Now we want to construct the stuffle regularized multiple Eisenstein series. 
Motivated by the calculation of the Fourier expansion of multiple Eisenstein series, we consider the following construction. 

\begin{constr} \label{const}
    Given a $\Q$-algebra $(A,\cdot)$ and a family of homomorphisms
    \[ \left\{ w \mapsto f_w(m) \right\}_{m \geq 1} \]
    from  $\h^1_\ast$ to $(A,\cdot)$, we set $F_{\emptyset}(M):=1$. For a non-empty word $w \in \h^1$ and $M \geq 1$ we define
    \[ F_w(M) := \sum_{\substack{1 \leq k \leq \len(w)\\w_1 \dots w_k = w\\ M> m_1 > \dots >m_k >0}}  f_{w_1}(m_1) \dots f_{w_k}(m_k) \in A \,, \]
    where $\len(w)$ denotes the length of the word $w$ and $w_1 \dots w_k = w$ is a decomposition of $w$ into $k$ non-empty words in $\h^1$.
\end{constr}
\begin{proposition}(\cite[Prop. 6.8]{Ba7})\label{prop:construction}
    For all $M \geq 1$ the assignment $w \mapsto F_w(M)$, described above, determines an algebra homomorphism from $\h^1_\ast$ to $(A,\cdot)$. In particular $\left\{ w \mapsto F_w(m) \right\}_{m \geq 1}$ is again a family  of homomorphisms as used in Construction \ref{const}. 
\end{proposition}

For a word $w=z_{k_1} \dots z_{k_r} \in \h^1$ we also write in the following $f_{k_1,\dots,k_r}(m):=f_w(m)$ and similarly $ F_{k_1,\dots,k_r}(M):=F_w(M)$. 

\begin{ex} Let $f_w(m)$ be as in Construction \ref{const}. In small depths the $F_w$ are given by
    \[ F_{k_1}(M) = \sum_{M> m_1 >0} f_{k_1}(m_1) \,, \quad F_{k_1,k_2}(M) = \sum_{M> m_1 >0} f_{k_1,k_2}(m_1) + \sum_{M> m_1 > m_2>0} f_{k_1}(m_1)  f_{k_2}(m_2) \,\]
    and one can check directly by the use of the stuffle product for the $f_w$ that 
    \begin{align*}
    &F_{k_1}(M) \cdot F_{k_2}(M) =  \sum_{M> m_1>0} f_{k_1}(m_1) \cdot  \sum_{M> m_2 >0} f_{k_2}(m_2) \\
    &=   \sum_{M> m_1 > m_2>0} f_{k_1}(m_1)  f_{k_2}(m_2) +  \sum_{M> m_2 > m_1>0}   f_{k_2}(m_2) f_{k_1}(m_1) +  \sum_{M>  m_1> 0}   f_{k_1}(m_1) f_{k_2}(m_1)  \\
    &=  \sum_{M> m_1 > m_2>0} f_{k_1}(m_1)  f_{k_2}(m_2) +  \sum_{M> m_2 > m_1>0}   f_{k_2}(m_2) f_{k_1}(m_1) \\
    &+ \sum_{M> m_1>0} \left(f_{k_1,k_2}(m_1) +  f_{k_2,k_1}(m_1) + f_{k_1+k_2}(m_1) \right)  \\
    &=  F_{k_1,k_2}(M)+ F_{k_2,k_1}(M) + F_{k_1+k_2}(M)  \,.
    \end{align*}
\end{ex}

Let us now give an explicit example for maps $f_w$ in which we are interested. Recall (Definition \ref{def:multitangent}) that for integers $k_1,\dots,k_r \ge 2$ we defined the multitangent function by
\[ \Psi_{k_1,\ldots,k_r}(x) = \sum_{\substack{n_1>\cdots>n_r\\ n_j \in \Z}} \frac{1}{(x+n_1)^{k_1}\cdots (x+n_r)^{k_r}}.\]
In \cite{Bo}, where these functions were introduced, the author uses the notation $\mathcal{T}e^{k_1,\ldots,k_r}(x)$ which corresponds to our notation $\Psi_{k_1,\ldots,k_r}(x)$. It was shown there that the series $\Psi_{k_1,\ldots,k_r}(x)$ converges absolutely when $k_1,\ldots,k_r\ge2$. These functions fulfill (for the cases they are defined) the stuffle product. As explained in Section \ref{subsec:fourier} the multitangent functions appear in the calculation of the Fourier expansion of the multiple Eisenstein series $\aG_{k_1,\dots,k_r}$, for example in depth two we have
\begin{align*}
\aG_{k_1,k_2}(\tau) ={}&\zeta(k_1,k_2) + \zeta(k_2) \sum_{m_1 > 0} \Psi_{k_1}(m_1 \tau) + \sum_{m_1>0} \Psi_{k_1,k_2}(m_1 \tau)
+\sum_{m_1 > m_2 > 0} \Psi_{k_1}(m_1\tau) \Psi_{k_2}(m_2\tau) \,.
\end{align*}
One nice result of \cite{Bo} is a regularization of the multitangent functions which defines $\Psi_{k_1,\ldots,k_r}(x)$ for all $k_1,\dots,k_r \geq 1$. We will use this result together with the above construction to recover the Fourier expansion of the multiple Eisenstein series. 
\begin{theorem}(\cite{Bo})\label{thm:olivier}
    For all $k_1,\dots,k_r \geq 1$  there exist holomorphic functions $\Psi_{k_1,\dots,k_r}$ on $\Ha$ with the following properties
    \begin{enumerate}[\textup{(}i\textup{)}] 
        \item Setting $q=e^{2\pi i \tau}$ for $\tau \in \Ha$ the map  $w \mapsto \Psi_w(\tau)$ defines an algebra homomorphism from $(\h^1, \ast)$ to $(\C[\![q]\!],\cdot)$.
        \item In the case $k_1,\dots,k_r \geq 2$ the  $\Psi_{k_1,\dots,k_r}$ are given by the multitangent functions in Definition \ref{def:multitangent}.
        \item The monotangent functions have the $q$-expansion given  by
        \begin{align*}
        \Psi_1(\tau) &= \frac{\pi}{\tan(\pi \tau)} =  (-2 \pi i)\left(\frac{1}{2} + \sum_{n>0} q^n \right)
        \end{align*}
        and for $k\geq 2$ by
        \begin{align*}
        \Psi_k(\tau) &=  \frac{(-2\pi i)^k}{(k-1)!} \sum_{n>0} n^{k-1} q^n\,.
        \end{align*}
        \item (Reduction into monotangent function) Every $\Psi_{k_1,\dots,k_r}(\tau)$ can be written as a $\mz$-linear combination of monotangent functions. There are explicit $\epsilon^{k_1,\dots,k_r}_{i,k} \in \mz$ such that
        \[\Psi_{k_1,\dots,k_r}(\tau) = \delta^{k_1,\dots,k_r} + \sum_{i=1}^r \sum_{k=1}^{k_i} \epsilon^{k_1,\dots,k_r}_{i,k} \Psi_k(\tau) \,,\]
        where $ \delta^{k_1,\dots,k_r} = \frac{(\pi i)^r}{r!}$ if $k_1=\dots=k_r=1$ and $r$ even and $\delta^{k_1,\dots,k_r} = 0$ otherwise. 
        For $k_1>1$ and $k_r>1$ the sum on the right starts at $k=2$, i.e. there are no $\Psi_1(\tau)$ appearing and therefore there is no constant term in the $q$-expansion.
    \end{enumerate}
\end{theorem}
\begin{proof}
    This is just a summary of the results in Sections~6 and~7 of \cite{Bo}. The last statement, (iv), is given by \cite[Theorem~6]{Bo}. 
\end{proof} 
Due to (iv) in the Theorem, the calculation of the Fourier expansion of multiple Eisenstein series, where ordered sums of multitangent functions appear, reduces to ordered sums of monotangent functions. We have seen in Proposition \ref{prop:ghatasmonotangent} that for $k_1,\dots,k_r \geq 2$ we have
\begin{align*}
&\ag(k_1,\ldots,k_r) = \sum_{m_1 >\dots >m_r>0} \Psi_{k_1}(m_1\tau) \dots \Psi_{k_r}(m_r\tau) \,.
\end{align*}
For $w \in \h^1$ we now use Construction \ref{const} with $A=\C[\![q]\!]$ and the family of homomorphisms $\{ w \mapsto \Psi_w(n \tau) \}_{n\geq 1}$ (see Theorem \ref{thm:olivier} (i)) to define
\[ \ag^{\ast,M}(w) := \sum_{\substack{1 \leq k \leq l(w)\\w_1 \dots w_k = w}} \sum_{M>m_1 > \dots > m_k >0} \Psi_{w_1}(m_1 \tau) \dots \Psi_{w_k}(m_k \tau)  \,. \] 
From Proposition \ref{prop:construction} it follows that for all $M\geq 1$ the map $\ag^{\ast,M}$ is an algebra homomorphism from $(\h^1,\ast)$ to $\C[\![q]\!]$.

To define stuffle regularized  multiple Eisenstein series we need the following: For an arbitrary quasi-shuffle algebra $\Q\langle A \rangle$ define the following coproduct for a word $w$ 
\[ \Delta_H(w) = \sum_{u v = w} u \otimes v\,. \] 
Then it is known due to Hoffman (See \cite{HI}) that the space $\left( \Q\langle A \rangle, \qsh , \Delta_H \right)$ has the structure of a Hopf algebra. With this, we try to mimic the definition of the $\aG^\sh$ and use the coproduct structure on the space $(\h^1, \ast, \Delta_H)$ to define for $M\geq 0$ the function $\aG^{\ast,M}$ and then take the limit $M\rightarrow \infty$ to obtain the stuffle regularized multiple Eisenstein series. For this, we consider the following diagram
\[
\xymatrix{  
    (\h^1,\ast)  \ar@{->}^-{\Delta_H} [r]\ar@{->}[d]_{\aG^{\ast,M}}  &(\h^1,\ast) \otimes (\h^1,\ast) \ar@{->}[d]^{\ag^{\ast,M} \otimes\, \zeta^{\ast} }   \\ 
    \C[\![q]\!] & \C[\![q]\!] \otimes  \mz  \ar@{->}[l]^-{m}   &\\
} \]
with the above algebra homomorphism $\ag^{\ast,M} : (\h^1, \ast) \rightarrow \C[\![q]\!]$ and the map  $\zeta^\ast$ for stuffle regularized multiple zeta values given in Definition \ref{def:regmzv} (with $T=0$).

\begin{definition}\label{def:gast}
    For integers $k_1,\ldots,k_r\ge1$ and $M \geq 1$, we define the $q$-series $\aG^{\ast,M}_{k_1,\ldots,k_r}\in\C[\![q]\!]$ as the image of the word $w= z_{k_1}\dots z_{k_r} \in \h^1$ under the algebra homomorphism $(\ag^{\ast,M} \otimes \zeta^{\ast})\circ \Delta_H$:
    \[ \aG^{\ast,M}_{k_1,\ldots,k_r}(\tau) :=m\left( (\ag^{\ast,M} \otimes \zeta^{\ast})\circ \Delta_H \big( w \big) \right) \in \C[\![q]\!]\,.\]
\end{definition}

For $k_1,\dots,k_r \geq 2$ the limit
\begin{equation}\label{eq:stufflimit}
\aG^\ast_{k_1,\ldots,k_r}(\tau) :=  \lim_{M\to\infty}   \aG^{\ast,M}_{k_1,\ldots,k_r}(\tau) \, 
\end{equation} 
exists and we have $\aG_{k_1,\ldots,k_r} = \aG^\ast_{k_1,\ldots,k_r}=\aG^\sh_{k_1,\ldots,k_r}$ (\cite[Prop. 6.13]{Ba7}).

\begin{remark} The limit in \eqref{eq:stufflimit} exists for $k_1\geq2$ and $k_2,\dots,k_r\geq 1$ by \cite[Proposition~4.8]{Ba10}. For $k_1=1$ this limit need not exist, since the monotangent function $\Psi_{1}$ has a constant term in its Fourier expansion. Thus the direct limit above cannot be used for all indices. A harmonic regularization for arbitrary indices, based on regularized multitangent functions, is constructed in \cite{Ba10}.
\end{remark}

\begin{theorem}(\cite{Ba7}) For all $k_1,\dots,k_r \geq 1$ and $M \geq 1$ the $\aG_{k_1,\dots,k_r}^{\ast,M} \in \C[\![q]\!]$ have the following properties:
    \begin{enumerate}[\textup{(}i\textup{)}] 
        \item Their product can be expressed in terms of the stuffle product.
        \item In the case where the limit $\aG^\ast_{k_1,\dots,k_r} :=  \lim_{M\to\infty} \aG_{k_1,\dots,k_r}^{\ast,M}$ exists, the functions $\aG^\ast_{k_1,\dots,k_r}$ are elements in $\MZB$. 
        \item For $k_1,\dots,k_r \geq 2$ the $\aG^\ast_{k_1,\dots,k_r}$ exist and equal the classical multiple Eisenstein series
        \[ \aG_{k_1,\dots,k_r}(\tau) = \aG^\ast_{k_1,\dots,k_r}(\tau) \,.\] 
    \end{enumerate}
\end{theorem} 

The admissible limits define an algebra homomorphism $\aG^\ast:\h^0_\ast\rightarrow\MZB$. Since $\h^1_\ast\simeq\h^0_\ast[z_1]$, we can extend this map by harmonic regularization. Following \cite[Proposition~4.8 and Definition~4.9]{Ba10}, we choose
\begin{align*}
\aG^\ast(z_1)=\ag(1)=(-2\pi i)\g(1)
\end{align*}
and obtain an algebra homomorphism $\aG^\ast:\h^1_\ast\rightarrow\MZB$. For arbitrary positive indices, $\aG^\ast_{k_1,\dots,k_r}$ will from now on denote the values of this harmonic-regularized extension.

\subsection{Double shuffle relations for regularized multiple Eisenstein series}

By Theorem \ref{thm:messh} we know that the product of two shuffle regularized multiple Eisenstein series $\aG^\sh_{k_1,\dots,k_r}$ with $k_1,\dots,k_r \geq 1$ can be expressed by using the shuffle product formula. This means we can for example replace every $\zeta$ by $\aG^\sh$ in the shuffle product (Example \ref{ex:z2z3}) of multiple zeta values and obtain
\begin{equation}\label{eq:shmessh}
\aG^\sh_2 \cdot \aG^\sh_3 =  \aG^\sh_{2,3} + 3 \aG^\sh_{3,2} + 6 \aG^\sh_{4,1} \,.
\end{equation}
Due to Theorem \ref{thm:messh} iii) we know that $\aG^\sh_{k_1,\dots,k_r}=\aG_{k_1,\dots,k_r}$ whenever $k_1,\dots,k_r \geq 2$. Since the product of two multiple Eisenstein series $\aG_{k_1,\dots,k_r}$ can be expressed using the stuffle product formula we also have
\begin{align}\label{eq:stmessh}
\begin{split}
\aG^\sh_2 \cdot \aG^\sh_3 &= \aG_2 \cdot \aG_3 = \aG_{2,3} + \aG_{3,2} + \aG_{5} \\
&= \aG^\sh_{2,3} + \aG^\sh_{3,2} + \aG^\sh_{5}   \,.
\end{split}
\end{align}
Combining \eqref{eq:shmessh} and \eqref{eq:stmessh} we obtain the relation $\aG^\sh_5 = 2 \aG^\sh_{3,2} + 6 \aG^\sh_{4,1}$. In the following we will call these relations, i.e. the relations obtained by writing the product of two $\aG^\sh_{k_1,\dots,k_r}$ with $k_1,\dots,k_r \geq 2$ as the stuffle and shuffle product, \emph{restricted double shuffle relations}.

We know that multiple zeta values fulfill even more linear relations, in particular we can express the product of two multiple zeta values $\zeta(k_1,\dots,k_r)$ in two different ways whenever $k_1\geq 2$ and $k_2,\dots,k_r\geq 1$. A natural question therefore is in which cases the $\aG^\sh$ also fulfill these additional relations. The answer to this question is that some are satisfied and some are not, as the following will show. 

The four examples in \cite[Example 6.15]{Ba7} are part of the following general comparison theorem.
\begin{theorem}[{\cite[Main Theorem~A]{BKM}}]\label{thm:mesregcomparison}
Suppose that $k_1,\dots,k_r\geq1$, with all entries at least $2$ except for at most one entry $k_i=1$ with $2\leq i\leq r$. Then
\begin{align*}
\aG^\sh_{k_1,\dots,k_r}=\aG^\ast_{k_1,\dots,k_r}.
\end{align*}
\end{theorem}
In particular, $\aG^\sh_{2,1,2}=\aG^\ast_{2,1,2}$, $\aG^\sh_{2,1}=\aG^\ast_{2,1}$, $\aG^\sh_{2,2,1}=\aG^\ast_{2,2,1}$ and $\aG^\sh_{4,1}=\aG^\ast_{4,1}$. Since the product of two $\aG^\ast$ can be expressed using the stuffle product we obtain
\begin{align}\label{eq:stufflew5}
\begin{split}
\aG^\sh_2 \cdot \aG^\sh_{2,1} &= \aG^\ast_2 \cdot \aG^\ast_{2,1} \\
&= \aG^\ast_{2,1,2} + 2 \aG^\ast_{2,2,1} + \aG^\ast_{4,1}+ \aG^\ast_{2,3} \\
&=  \aG^\sh_{2,1,2} + 2 \aG^\sh_{2,2,1} + \aG^\sh_{4,1}+ \aG^\sh_{2,3} \,.
\end{split}
\end{align}
Using also the shuffle product to express $\aG^\sh_2 \cdot \aG^\sh_{2,1}$ we obtain a linear relation in weight $5$ which is not covered by the restricted double shuffle relations. This linear relation was numerically observed in \cite{BT} but could not be proved there. Theorem~\ref{thm:mesregcomparison} proves the comparison whenever there is at most one internal index $1$. Outside this range, the exact relationship between $\aG^\ast$ and $\aG^\sh$ is still open. Similarly to multiple zeta values we expect, but cannot prove yet, that the space spanned by all $\aG^\shuffle$ equals the one spanned by $\aG^\ast$.

Notice that \cref{thm:mesregcomparison} also implies $\ds(u,v)\in\ker\aG^\sh$ for all $u\in\HH^{\geq2}$ and all $v$ which have at most one entry $1$, not at the beginning (see \cite[Main Theorem~A(ii)]{BKM}), and the relation in weight $5$ obtained from \eqref{eq:stufflew5} is an example of this. As a nice consequence, one gets a family of relations among the multiple Eisenstein series $\aG_{k_1,\dots,k_r}$ with $k_1,\dots,k_r\geq2$, which was proved by the author, Kanno and Maesaka. For this, we put for $u,v\in\HH^{\geq2}$
\begin{align}\label{eq:mesRuv}
R(u,v)=\ds(u,z_2)\ast v+u\ast\ds(v,z_2)-\ds(u\ast v,z_2)\,,
\end{align}
i.e. $R(u,v)$ measures how far the map $w\mapsto\ds(w,z_2)$ is from being a derivation with respect to the stuffle product.
\begin{theorem}[{\cite[Main Theorem~B and Lemma~A.1]{BKM}}]\label{thm:mesRrelations}
For all $u,v\in\HH^{\geq2}$ we have $R(u,v)\in\HH^{\geq2}$ and $\aG\bigl(R(u,v)\bigr)=0$.
\end{theorem}
\begin{proof}
For $w\in\HH^{\geq2}$, the element $\ds(w,z_2)$ is a linear combination of words with at most one letter $z_1$, which is not at the beginning. By the restricted double shuffle relations we have $\aG^\sh(\ds(w,z_2))=0$, and by \cref{thm:mesregcomparison} therefore also $\aG^\ast(\ds(w,z_2))=0$. Since $\aG^\ast$ is an algebra homomorphism with respect to the stuffle product, this gives $\aG^\ast(R(u,v))=0$. The fact that $R(u,v)\in\HH^{\geq2}$, i.e. that all words containing $z_1$ cancel, is \cite[Lemma~A.1]{BKM}.
\end{proof}
For example, we have
\begin{align*}
R(z_2,z_2)&=6z_3z_3-3z_4z_2-z_6\,,\\
R(z_3,z_2)&=4z_3z_4+3z_4z_3-2z_5z_2-z_7\,,
\end{align*}
which give exactly the two relations in weight $6$ and $7$ from \Cref{subsec:mesoverview}. If \cref{conj:mesdim} is true, then the ideal generated by the relations in \cref{thm:mesRrelations} contains all relations among multiple Eisenstein series up to weight $16$, but not in weight $17$ (see \cite{BKM}). We will come back to this at the end of \Cref{sec:messl2}.

Define the space of regularized multiple Eisenstein series of weight $k\geq 0$ by
\begin{align*}
\rmes_k = \langle \aG^\ast_{\kk} \mid \kk \in \Z_{\geq 1}^r, r\geq 0, \wt(\kk)=k \rangle_\Q \,
\end{align*}
and set $\rmes = \sum_{k\geq 0} \rmes_k$. We end this section by giving a dimension conjecture for these spaces. Recall that the Hilbert--Poincar\'e series of modular forms, cusp forms and quasimodular forms are given by
\begin{align*}
\mathsf{M}(X) &= \frac{1}{(1-X^4)(1-X^6)},\\
\mathsf{S}(X) &= X^{12}\mathsf{M}(X) = \frac{X^{12}}{(1-X^4)(1-X^6)},\\
\widetilde{\mathsf{M}}(X) &= \mathsf{D}(X) \mathsf{M}(X) = \frac{1}{(1-X^2)(1-X^4)(1-X^6)}.
\end{align*}
In addition, define the Hilbert--Poincar\'e series of the space of period polynomials $W_k$ (defined in \eqref{eq:defwk}) with  even $k\geq 2$ by
\begin{align*}
  \mathsf{W}(X) = \sum_{\substack{k\geq 2\\ k \text{ even}}} \dim_\Q W_k X^k = \mathsf{M}(X)+\mathsf{S}(X)-1 = \frac{X^4}{1-X^2} + 2\, \mathsf{S}(X)
\end{align*}
and recall that $\mathsf{D}(X)= \frac{1}{1-X^2}$ and $\mathsf{O}(X) = \frac{X^3}{1-X^2}$.
Then we have the following conjecture.

\begin{conjecture}\label{conj:rmesdim}
We have
\begin{align*}
    \sum_{k\geq 0} \dim_\Q \rmes_k X^k &= \widetilde{\mathsf{M}}(X) \cdot  \frac{1}{1 - \mathsf{D}(X) (X+\mathsf{O}(X)) + \mathsf{D}(X) \mathsf{W}(X)}\\
&= \mathsf{M}(X) \cdot  \frac{1}{1 - X - X^2 - \mathsf{O}(X) + \mathsf{W}(X)}\\
    &= \frac{1}{1 - X- X^2 - X^3 + X^6 + X^7 + X^8 + X^9}.
\end{align*}
\end{conjecture}

This conjecture was proposed in \cite{BK2} for the associated weight graded spaces of $q$-analogues. As it is expected that regularized multiple Eisenstein series satisfy the same relation as weight graded $q$-analogues, we obtain the above conjecture.

\section{Double-indexed \texorpdfstring{$q$}{q}-analogues of MZVs} \label{subsec:doubleindexedqana}

In this section, we want to study the double-indexed version of $\g$, which we defined in Definition \ref{def:doubleg} for $k_1,\dots,k_r \geq 1$, $d_1,\dots,d_r\geq 0$ by
\begin{align*}
\mb{k_1,\dots,k_r}{d_1,\dots,d_r} &:= \sum_{m_1 > \cdots > m_r > 0} m_1^{d_1}\frac{P_{k_1}(q^{m_1})}{(1-q^{m_1})^{k_1}} \cdots m_r^{d_r} \frac{P_{k_r}(q^{m_r})}{(1-q^{m_r})^{k_r}} \\
&= \sum_{\substack{m_1>\dots>m_r > 0\\n_1,\dots,n_r > 0}} m_1^{d_1} \frac{n_1^{k_1-1}}{(k_1-1)!} \cdots m_r^{d_r} \frac{n_r^{k_r-1}}{(k_r-1)!}  q^{m_1 n_1 + \dots + m_r n_r}\,.
\end{align*}
By $k_1+\dots + k_r + d_1 + \dots + d_r$ we denote its weight and by $r$ its depth. These $q$-series are also modified $q$-analogues of multiple zeta values.

\begin{proposition} \label{prop:doublegasqana}For $k_1 \geq d_1+2$ and $k_j \geq d_j + 1$ for $j=2,\dots,r$ we have 
    \begin{align*}
    \lim_{q\rightarrow 1} (1-q)^{k_1+\dots+k_r} \mb{k_1,\dots,k_r}{d_1,\dots,d_r}  = \zeta(k_1-d_1, k_2-d_2,\dots,k_r-d_r)\,.
    \end{align*} 
\end{proposition}
\begin{proof}
    This proof is similar to the proof of Proposition \ref{prop:gisqmzv}, since we have 
    \begin{align*}
    \lim\limits_{q\rightarrow 1} (1-q)^k\frac{ m^d P_k(q^m)}{(1-q^m)^k} =     \lim\limits_{q\rightarrow 1} \frac{m^d P_k(q^m)}{[m]_q^k} = \frac{m^d P_k(1)}{m^k} = \frac{1}{m^{k-d}}\,.
    \end{align*}
\end{proof}
In fact these $q$-series span the space of all modified $q$-analogues, which we defined in Section \ref{subsec:qmzv} by 
\begin{align*}
\mz_q := \Q + \Big \langle   \zeta_q(k_1,\dots,k_r ; Q_1 , \dots , Q_r) \,\,\big|\,\, r\ge 1 ,\,k_1,\dots,k_r\geq1 ,\,\deg(Q_j) \leq k_j 
\Big\rangle_\Q\,.
\end{align*}
Here we defined the $\zeta_q(k_1,\dots,k_r ; Q_1 , \dots , Q_r)$ for $k_1,\dots,k_r \geq 1$ and polynomials $Q_1(X) \in X \Q[X]$ and $Q_2(X),\dots,Q_r(X) \in \Q[X]$ by
\begin{align*} 
\zeta_q(k_1,\dots,k_r ; Q_1 , \dots , Q_r) = \sum_{m_1>\cdots> m_r >0} \frac{Q_1(q^{m_1}) \dots Q_r(q^{ m_r})}{(1-q^{m_1})^{k_1}\cdots (1-q^{m_r})^{k_r}} \,.
\end{align*}
These  series are modified $q$-analogues of $\zeta(k_1,\dots,k_r)$, since we have for $k_1 \geq 2$  \[ \lim\limits_{q\rightarrow 1} (1-q)^{k_1+\dots+k_r} \zeta_q(k_1,\ldots,k_r; Q_1 , \dots , Q_r) = Q_1(1) \dots Q_r(1)\cdot \zeta(k_1,\ldots,k_r) \,. \]

\begin{theorem} \label{thm:zqspannedbydoubleg}The space $\mz_q$ of modified $q$-analogues is spanned by the $q$-series $\mb{k_1,\dots,k_r}{d_1,\dots,d_r}$, i.e. 
    \begin{align*}
    \mz_q = \Q + \Big \langle   \mb{k_1,\dots,k_r}{d_1,\dots,d_r} \,\,\big|\,\, k_1,\dots,k_r \geq 1, d_1,\dots,d_r\geq 0
    \Big\rangle_\Q\,.
    \end{align*}
\end{theorem}
\begin{proof} We give a variation of the proof as given in \cite{BK2}.
    First we show the inclusion '$\subseteq$', i.e. that every $\zeta_q(k_1,\dots,k_r;Q_1,\dots,Q_r)$ can be written in terms of $\g$. 
    For all $k\geq 1$ we have $P_k(1) =1$ and $P_k(0)=0$ and therefore the polynomials $P_j(X) (1-X)^{k-j}$ with $j=1,\dots,k$ form a basis of the space $\{ Q \in X \Q[X] \mid \deg Q \leq k \} $. In particular for every polynomial $Q$ in this space there exist coefficients $\alpha_j \in \Q$ with 
    \begin{equation} \label{eq:QasP}
    \frac{Q(X)}{(1-X)^k} = \sum_{j=1}^k \alpha_j \frac{ P_j(X)}{(1-X)^j}\,.
    \end{equation}
    Therefore we just need to see what happens if one of the $Q_2,\dots,Q_r$ has a constant term. Without loss of generality we can focus on the cases $Q_i(X) = 1$ for a $2 \leq i \leq r$.  Since for all $k\geq 1$ we have
    \[ \frac{1}{(1-X)^k} = 1 + \sum_{m=1}^k \frac{X}{(1-X)^m}\,,\]
    we can write 
    \begin{align*} \sum_{n_1>\cdots>n_r>0} \prod_{j=1}^r \frac{Q_j(q^{n_j})}{(1-q^{n_j})^{s_j}} &=  \sum_{n_1>\cdots>n_r>0} \prod_{\substack{j=1\\ j \neq i}}^{r} \frac{Q_j(q^{n_j})}{(1-q^{n_j})^{s_j}} +  \sum_{\substack{n_1>\cdots>n_r>0\\ 1 \leq m \leq s_i}} \frac{q^{n_i}}{(1-q^{n_i})^m} \prod_{\substack{j=1\\ j \neq i}}^{r} \frac{Q_j(q^{n_j})}{(1-q^{n_j})^{s_j}} \,.
    \end{align*}
    For the second sum on the right-hand side we can again use \eqref{eq:QasP}. For the first sum we obtain (by setting $n_{r+1} =0$)
    \begin{align*}
    \sum_{n_1>\cdots>n_r>0} \prod_{\substack{j=1\\ j \neq i}}^{r} \frac{Q_j(q^{n_j})}{(1-q^{n_j})^{s_j}} = \sum_{n_1>\cdots>n_{i-1} > n_{i+1} > \dots >n_r>0} (n_{i-1} - n_{i+1}-1) \prod_{\substack{j=1\\ j \neq i}}^{r} \frac{Q_j(q^{n_j})}{(1-q^{n_j})^{s_j}}\,.
    \end{align*}
    Repeating this for all  $2 \leq i \leq r$ with $Q_i(X) = 1$ we obtain linear combinations of $\g$, from which we deduce '$\subseteq$'.
    To prove '$\supseteq$' we first define for $m\geq 0$ the polynomials $p_m(n)$ by $p_0(n)=1$ and 
    \begin{equation}\label{eq:defpm}
    p_m(n) = \sum_{n > N_1 > \dots > N_m > 0} 1 = \binom{n-1}{m} \,. 
    \end{equation}
    The polynomial $p_m(n)$ has degree $m$ and therefore we can always find $c_m(r) \in \Q$ with $n^r = \sum_{m=0}^r c_m(r)\, p_m(n)$. The idea is now to replace $n_j^{r_j}$ in the definition of $\g$ by $\sum_{m_j=0}^{r_j} c_{m_j}(r_j)\, p_{m_j}(n_j)$ and then use \eqref{eq:defpm} to get sums which can be written in terms of the $\zeta_q$. We illustrate this in the depth two case from which the general case becomes clear. 
    In depth two, we have
    \begin{align*}
    \mb{k_1,k_2}{d_1,d_2} &= \sum_{n_1 > n_2 > 0} \frac{n_1^{d_1}  P_{k_1}(q^{n_1})}{(1-q^{n_1})^{k_1}} \frac{n_2^{d_2} P_{k_2}(q^{n_2})}{(1-q^{n_2})^{k_2}} \\
    &= \sum_{0 \leq m_2 \leq d_2 } c_{m_2}(d_2) \sum_{n_1 > n_2 > N_1 > \dots > N_{m_2}> 0} \frac{n_1^{d_1} P_{k_1}(q^{n_1})}{(1-q^{n_1})^{k_1}} \frac{ P_{k_2}(q^{n_2})}{(1-q^{n_2})^{k_2}} \\
    &= \sum_{\substack{0 \leq m_1 \leq d_1 \\0 \leq m_2 \leq d_2 }} c_{m_1}(d_1) c_{m_2}(d_2) \sum_{\substack{n_1 > n_2 > N_1 > \dots > N_{m_2}> 0 \\ n_1 >  N'_1 > \dots > N'_{m_1}> 0}} \frac{P_{k_1}(q^{n_1})}{(1-q^{n_1})^{k_1}} \frac{ P_{k_2}(q^{n_2})}{(1-q^{n_2})^{k_2}}\,.
    \end{align*}
    Now considering all the possible shuffles, and possible equalities of the $N$ and the $N'$ it is clear that this sum can be written as a linear combination of $\zeta_q$ by interpreting appearing $1$ as $(1-q^N) (1-q^N)^{-1}$. For general depth $r$ the idea is the same and therefore we obtain '$\supseteq$'.
\end{proof}

Since the double-indexed $\g$ generalize the single-indexed $\g$ we clearly have $\gs \subset \mz_q$. In Proposition \ref{prop:gevensubspace} we saw that the following spaces
\begin{align*}
\gstwo &= \Big \langle   \g(k_1,\dots,k_r) \,\,\big|\,\, r\ge 0 ,\,k_1,\dots,k_r\geq 2 
\Big\rangle_\Q \,,\\
\gsev &= \Big \langle   \g(k_1,\dots,k_r) \,\,\big|\,\, r\ge 0 ,\,k_1,\dots,k_r\geq 2  \text{ even} 
\Big\rangle_\Q\,,
\end{align*} 
are also subalgebras of $\gs$. Further we saw in Proposition \ref{prop:qmprop} that the algebra of quasi-modular forms (with rational coefficients) $\qmf = \Q[\g(2), \g(4),\g(6)]$ is a subalgebra of $\gsev$.
Combining all this we obtain the following inclusions of $\Q$-algebras 
\begin{align}\label{eq:algebrainclusion}
\mf^\Q \subset \qmf \subset \gsev  \subset \gstwo  \subset \gs \subset \mz_q \,.
\end{align}

We will see that $\g$ form a nice generating set of modified $q$-analogues since they behave well under the operator $q \frac{d}{dq}$, and they satisfy the partition relations, which we can use to describe an analogue of double shuffle relations for them. Conjecturally these relations give all relations among elements in $\mz_q$ and by \eqref{eq:algebrainclusion} they are sufficient to prove any relation among (quasi-)modular forms.

\subsection{The operator $q \frac{d}{dq}$}

In this section we will study the operator $q \frac{d}{dq}$. We saw already that the space of quasimodular forms $\qmf$ is closed under this operator. In this section we will see that this is also the case for $\mz_q$, $\gs$ and $\gstwo$. First notice that this operator on a $q$-series is given as follows
\begin{align*}
q \frac{d}{dq} \sum_{n=1}^\infty a_n q^n = \sum_{n=1}^\infty n \,a_n q^n\,.
\end{align*}
Applying this to $\g$ therefore gives 
\begin{align}\begin{split}\label{eq:derivedoubleg}
q \frac{d}{dq} \mb{k_1,\dots,k_r}{d_1,\dots,d_r} &= \sum_{\substack{m_1>\dots>m_r > 0\\n_1,\dots,n_r > 0}} m_1^{d_1} \frac{n_1^{k_1-1}}{(k_1-1)!} \cdots m_r^{d_r} \frac{n_r^{k_r-1}}{(k_r-1)!} (m_1 n_1 + \dots + m_r n_r) q^{m_1 n_1 + \dots + m_r n_r}\\
&=\sum_{j=1}^r  k_j \mb{k_1,\dots,k_j+1,\dots,k_r}{d_1,\dots,d_j+1,\dots,d_r}\,. 
\end{split}
\end{align}
In particular we obtain the following:
\begin{proposition}The space $\mz_q$ is closed under $q \frac{d}{dq}$.
\end{proposition}
\begin{proof}
    This follows by Theorem \ref{thm:zqspannedbydoubleg} together with \eqref{eq:derivedoubleg}.
\end{proof}

In \cite{BK1} it was shown that also the subspace $\gs$ generated by $\g(k_1,\dots,k_r)$ is closed under $q\frac{d}{dq}$. This is not obvious, since by \eqref{eq:derivedoubleg} we have 
\begin{align*}
q \frac{d}{dq} \g(k_1,\dots,k_r) = \sum_{j=1}^r k_j \mb{k_1,\dots,k_j+1,\dots,k_r}{0,\dots,1,\dots,0}\,. 
\end{align*}
A priori it is not clear why the right-hand side can also be written in terms of single indexed $\g$. In the following we will give a proof of this fact by using generating series. 
For this we consider the following series (which is a special case of the $\mtt{.}{.}{.}$ in \eqref{def:tribracket})
\begin{align*}
H\bi{n_1,\dots,n_r}{X_1,\dots,X_r} &= \sum_{m_1 > \dots > m_r > 0} e^{m_1 X_1} \left( \frac{q^{m_1}}{1-q^{m_1}}\right)^{n_1} \cdots  e^{m_r X_r} \left( \frac{q^{m_r}}{1-q^{m_r}}\right)^{n_r} \,.
\end{align*} 
By Lemma \ref{lem:genexpr} we have the following relationship between $H$ and the generating series of $\g(k_1,\dots,k_r)$
\begin{align*}
H\bi{1,\dots,1}{X_1,\dots,X_r} &= \ggen(X_1+\dots+X_r, X_1+\dots+X_{r-1},\dots,X_1)\,.
\end{align*}
Also notice that the $H$ satisfy the stuffle product formula, e.g. 
\begin{align*}
H\bi{n_1}{X_1} H\bi{n_2}{X_2} = H\bi{n_1,n_2}{X_1,X_2} + H\bi{n_2,n_1}{X_2,X_1} + H\bi{n_1+n_2}{X_1+X_2}\,.
\end{align*} 
Therefore we get 
\begin{align}\begin{split}\label{eq:hderiv1}
\ggen(X_1) \ggen(X_2) &= H\bi{1}{X_1} H\bi{1}{X_2} = H\bi{1,1}{X_1,X_2} + H\bi{1,1}{X_2,X_1} + H\bi{2}{X_1+X_2} \\
&= \ggen(X_1+X_2,X_1) + \ggen(X_1+X_2,X_2) + H\bi{2}{X_1+X_2} \,.
\end{split}
\end{align}
Now we use $q \frac{d}{dq} \frac{q^m}{1-q^m} = m \left(  \frac{q^m}{1-q^m} \right)^2 + m  \frac{q^m}{1-q^m}$ to obtain 
\begin{align}\begin{split}\label{eq:hderiv2}
q \frac{d}{dq} \ggen(X)  &= q \frac{d}{dq} H\bi{1}{X} = q \frac{d}{dq}  \sum_{m>0} e^{mX} \frac{q^m}{1-q^m} = \sum_{m>0} m e^{mX}  \left(  \frac{q^m}{1-q^m} \right)^2 + \sum_{m>0} m e^{mX}  \frac{q^m}{1-q^m} \\
&= \frac{d}{dY}\left( H\bi{2}{X+Y} + H\bi{1}{X+Y} \right)_{Y=0}\,.
\end{split}
\end{align}
Combining \eqref{eq:hderiv1} and \eqref{eq:hderiv2} we obtain
\begin{align*}
q \frac{d}{dq} \ggen(X) = \frac{d}{dY}\big( \ggen(X) \ggen(Y) - \ggen(X+Y,X)- \ggen(X+Y,Y) +  \ggen(X+Y)\big)_{Y=0}\,.
\end{align*}
Since $\frac{d}{dY} \ggen(Y)_{Y=0} = \g(2)$ the above formula somehow states that $q \frac{d}{dq} \g(k)$ measures, up to lower weight terms, the failure of the shuffle product formula for $\g(k) \g(2)$. Also notice that for $k_1,k_2\geq1$ and $k=k_1+k_2$ this is just a special case of the shuffle product formula
\begin{align*}
\g(k_1) \g(k_2) = &\sum_{j=1}^{k-1} \left( \binom{j-1}{k_1-1} + \binom{j-1}{k_2-1} \right) \g(j, k-j) \\
&+ \binom{k-2}{k_1-1} \left(\mb{k-1}{1} - \g(k-1) \right)\,.
\end{align*}
The same strategy works in arbitrary depths, and we obtain the following explicit formula for the derivative of $\g$ in arbitrary depths.

\begin{theorem} \label{thm:gclosedunderqdq} For $k_1,\dots,k_r \geq 1$ we have 
\begin{align*}
    \qdq \g(k_1,\dots,k_r) =&\,  \g(2,k_1,\dots,k_r) + \sum_{j=2}^r \g(k_1,\dots,k_{j-1},2,k_{j},\dots,k_r) + \sum_{j=1}^r \g(k_1,\dots,k_j+2,\dots,k_r)\\
    &- \sum_{\substack{1 \leq j \leq r\\a+b = k_j+2}} (a-1)  \g(k_1,\dots,k_{j-1},a,b,k_{j+1},\dots,k_r) \\
 & - \sum_{\substack{1 \leq i < j \leq r\\a+b = k_j+1}}  k_i\, \g(k_1,\dots,k_i+1,\dots,k_{j-1},a,b,k_{j+1},\dots,k_r)\\
 &- \sum_{j=1}^r  k_j\, \g(k_1,\dots,k_j+1,\dots,k_r,1) \\
 &+ \sum_{j=1}^r (r-j+1) k_j \g(k_1,\dots,k_j+1,\dots,k_r)\\
 &+ \sum_{j=1}^r \sum_{a=1}^{k_j+1} \beta_{a,k_j}\, \g(k_1,\dots,k_{j-1},a,k_{j+1},\dots,k_r),
\end{align*}
where the coefficients $\beta_{a,k} \in \Q$ are given for $1 \leq a \leq k+1$ by 
\begin{align*}
    \beta_{a,k} = \left( (-1)^{k-1} \binom{k+1-a}{2-a} - \binom{k+1-a}{k-a} \right) \frac{B_{k+2-a}}{(k+2-a)!}.
\end{align*}
In particular, the space $\gs$ is closed under $q \frac{d}{dq}$.
\end{theorem}
\begin{proof}
With the same idea as above (see \cite[(3.2)]{BK1}) one shows that the action of the operator $\qdq$ on $\ggen$ can be written as\footnote{Here we set $X_{r+1}=0$, i.e. the case $j=r+1$ gives the term $\ggen(X_1+Y,\dots,X_r+Y,Y)$ for first sum on the right. }
\begin{align*}
  \qdq \ggen(X_1,\dots,X_r) &= \g(2) \ggen(X_1,\dots,X_r) - \sum_{j=1}^{r+1} \frac{d}{dY}\ggen(X_1+Y,\dots,X_j+Y,X_j,\dots,X_r)_{\big|Y=0}\\
  &+ \sum_{j=1}^r \frac{d}{dY}\ggen(X_1+Y,\dots,X_j+Y,X_{j+1},\dots,X_r)_{\big|Y=0}.
\end{align*}
Notice that this is equivalent to 
\begin{align*}
  \qdq \ggen(X_1,\dots,X_r) &= \g(2) \ggen(X_1,\dots,X_r) - \sum_{j=1}^{r+1} \frac{d}{dY}\ggen(X_1+Y,\dots,X_j+Y,X_j,\dots,X_r)_{\big|Y=0}\\
  &+ \sum_{j=1}^r (r-j+1) \frac{d}{dX_j}\ggen(X_1,\dots,X_r).
\end{align*}
Therefore the action of $\qdq$ on the $q$-series $\g$ is given by
\begin{align}\begin{split}\label{eq:qdgexplicitformula}
    \qdq \g(k_1,\dots,k_r) =&\, \g(2) \g(k_1,\dots,k_r)- \g(k_1,\dots,k_r,2)\\
    &- \sum_{\substack{1 \leq j \leq r\\a+b = k_j+2}} (a-1)  \g(k_1,\dots,k_{j-1},a,b,k_{j+1},\dots,k_r) \\
 & - \sum_{\substack{1 \leq i < j \leq r\\a+b = k_j+1}}  k_i \g(k_1,\dots,k_i+1,\dots,k_{j-1},a,b,k_{j+1},\dots,k_r)\\
 &- \sum_{j=1}^r  k_j \g(k_1,\dots,k_j+1,\dots,k_r,1) \\
 &+ \sum_{j=1}^r (r-j+1) k_j \g(k_1,\dots,k_j+1,\dots,k_r).
 \end{split}
\end{align}
Moreover, the product can be evaluated by Lemma \ref{lem:rproduct} as 
\[
\resizebox{\linewidth}{!}{$\begin{aligned}
    \g(2) &\g(k_1,\dots,k_r) = \sum_{j=0}^r \g(k_1,\dots,k_j,2,k_{j+1},\dots,k_r) + \sum_{j=1}^r \g(k_1,\dots,k_j+2,\dots,k_r) \\
    & + \sum_{j=1}^r \sum_{a=1}^{k_j+1} \left( (-1)^{k_j-1} \binom{k_j+1-a}{2-a} - \binom{k_j+1-a}{k_j-a} \right) \frac{B_{k_j+2-a}}{(k_j+2-a)!} \g(k_1,\dots,k_{j-1},a,k_{j+1},\dots,k_r).
\end{aligned}$}
\]
Combining this with the above gives the statement in the theorem.
\end{proof}

By Lemma \ref{lem:shuffleprodforgenseries} we can also interpret the right-hand side of \eqref{eq:qdgexplicitformula} as the measure of the failure of the shuffle product formula for $\g(k_1,\dots,k_r) \g(2)$ (up to lower weight terms).

\begin{theorem}[{\cite[Corollary~1.7]{BKM}}]\label{thm:g2closedunderdif}
    The space $\gstwo$ is closed under $\qdq$.
\end{theorem}
We will give the proof at the end of Section~\ref{sec:mesderiv}.

\begin{proposition} \label{prop:derivgkev1}For $k\geq 1$ we have 
    \begin{align*}
    q\frac{d}{dq} \g(k) = &(2k-1) \g(k+2) - \sum_{j=2}^k (k+j-1) \g(k+2-j,j)  - \g(2,k)\\
    &+ \sum_{j=2}^k \frac{B_{k+2-j}}{(k+2-j)!} (3k-j+1) \g(j) + (-1)^k \frac{B_k}{k!} \g(2)\,.
    \end{align*}
    In particular $q\frac{d}{dq} \g(k)  \in \gstwo$ for $k\geq 2$.
\end{proposition}
\begin{proof} This follows from Proposition \ref{prop:gshuffle1} and can also be found in \cite{Ba2}.
\end{proof}

The theorem is closely connected to multiple Eisenstein series. In Theorem \ref{thm:Fourier} we saw that the multiple Eisenstein series $\aG_{k_1,\dots,k_r}$ can be written as a $\C$-linear combination of $\g$. In particular we have
\begin{align*}
\mes^\C = \C + \langle \aG_{k_1,\dots,k_r} \mid r\geq 1, k_1, \dots,k_r \geq 2 \rangle_\C= \C \otimes_\Q \gstwo\,.
\end{align*}
Thus Theorem~\ref{thm:g2closedunderdif} gives the corresponding statement after extending the coefficients to $\C$. The stronger rational statement $\GD\mes\subset\mes$, together with an explicit formula for $\GD=(2\pi i)\frac{d}{d\tau}$, will be proved in Section~\ref{sec:mesderiv}. The depth-one formula is the following.

\begin{theorem}\label{thm:turanthm}
    For $k \geq 2$  we have 
    \begin{align*}
    (2\pi i) \frac{d}{d\tau}\aG_{k}(\tau) = (2k-1) \aG_{k+2}(\tau) - \sum_{j=2}^k (k+j-1) \aG_{k+2-j,j}(\tau)  - \aG_{2,k}(\tau)\,. 
    \end{align*}
\end{theorem}
\begin{proof}
    This is a consequence of Proposition \ref{prop:derivgkev1} together with the formula for the Fourier expansion of the double Eisenstein series given in Proposition \ref{prop:fourier23}. The details of this proof are worked out in the Bachelor thesis of Turan in \cite{Tu}.
\end{proof}

One does not expect that the space $\gsev$ is closed under $\qdq$, since it seems already to be the case that $\qdq \g(4,2) \notin \gsev$. Summarizing everything we get the following overview. 
\begin{figure}[h!]
    \begin{center}
        \begin{tikzcd}[column sep=8ex]
        \mf^{\Q} \arrow[r,hook] \arrow[r,bend right=50,"\qdq"]& 
        \qmf \arrow[loop,"\qdq",out=125,in=45,looseness = 5] \arrow[r,hook]&
        \gsev  \arrow[r,bend right=50,"\qdq ?",dashrightarrow]  \arrow[r,hook] &
        \gstwo  \arrow[loop,"\qdq",out=125,in=45,looseness = 5] \arrow[r,hook]&  
        \gs \arrow[loop,"\qdq",out=125,in=45,looseness = 5]  \arrow[r,hook]&
        \mz_q  \arrow[loop,"\qdq",out=125,in=45,looseness = 5] &
        \end{tikzcd}
    \end{center}
    \caption{Overview of the behavior of the operator $\qdq$. The dashed arrow is conjectural.}
\end{figure}
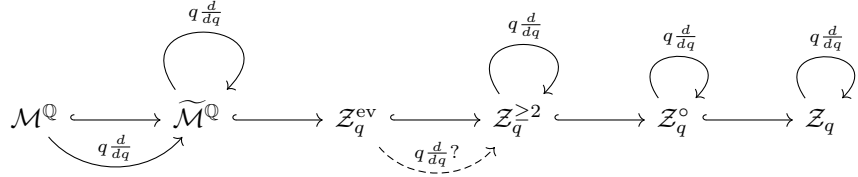

\subsection{Partition relation}
In this section we will introduce a family of linear relations among the $\g$, which follows immediately when viewing the coefficient of $q^N$ as a sum over all partitions of $N$.
By a partition of a natural number $N$ with $r$ distinct part sizes we denote a representation of $N$ which involves exactly $r$ distinct natural numbers. For example, $15 = 4 + 4 + 3 + 2 + 1 + 1$ is a partition of $15$ with the four distinct part sizes $4,3,2,1$. We identify such a partition with a tuple $(m,n) \in \N^r \times \N^r$ where the $m_j$'s are the distinct part sizes and the $n_j$'s count their appearance in the sum. The above partition of $15$ is therefore given by the tuple $(m,n) = ((4,3,2,1),(2,1,1,2))$. By $P_r(N)$ we denote all \emph{partitions of $N$ with $r$ distinct part sizes} and write
\[ P_r(N) := \left\{ (m,n) \in \N^r \times \N^r \, \mid \, N = m_1 n_1 + \dots + m_r n_r \, \text{ and } \, m_1 > \dots > m_r > 0  \right\}   \] 
On the set $P_r(N)$ we have an involution given by the conjugation $\rho$ of partitions which can be obtained by reflecting the corresponding Young diagram across the main diagonal. 
\begin{figure}[H]
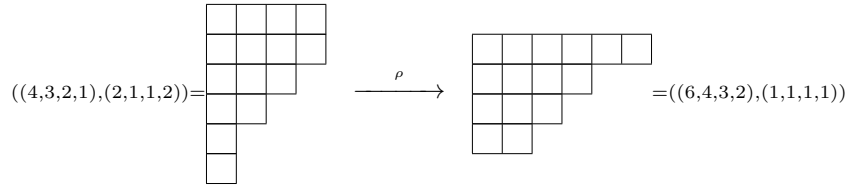

    \[ {\scriptstyle ((4,3,2,1),(2,1,1,2)) = \yng(4,4,3,2,1,1)  \quad  \overset{\rho}{\xrightarrow{\hspace*{1cm}}} \quad \yng(6,4,3,2) = ((6,4,3,2),(1,1,1,1))} \]
    \caption{The conjugation of the partition $15 = 4 + 4 + 3 + 2 + 1 + 1$ is given by $\rho( ((4,3,2,1),(2,1,1,2)) ) =  ((6,4,3,2),(1,1,1,1))$ which can be seen by reflecting the corresponding Young diagram across the main diagonal.  }
\end{figure}
On the set $P_r(N)$ the conjugation $\rho$ is explicitly given by $\rho( (m,n) ) = (m',n')$ where $m'_j = n_1 + \dots + n_{r-j+1}$ and $n'_j = m_{r-j+1}-m_{r-j+2}$ with $m_{r+1} := 0$, i.e.  
\begin{equation} \label{eq:conj}
\rho : \binom{m_1, \dots , m_r}{n_1, \dots ,n_r} {\longmapsto} \binom{n_1+\dots+n_r, \dots,n_1+n_2,n_1}{m_r,m_{r-1}-m_r,\dots,m_1-m_2} \,.
\end{equation}

We can interpret the coefficient of $q^N$ in $\g$ as a sum over all partitions of $N$, since we have 
\begin{align*}
\mb{k_1,\dots,k_r}{d_1,\dots,d_r} &= \sum_{\substack{m_1>\dots>m_r > 0\\n_1,\dots,n_r > 0}} m_1^{d_1} \frac{n_1^{k_1-1}}{(k_1-1)!} \cdots m_r^{d_r} \frac{n_r^{k_r-1}}{(k_r-1)!}  q^{m_1 n_1 + \dots + m_r n_r}\\
&= \sum_{N>0} \left(  \sum_{(m,n) \in P_r(N)} m_1^{d_1} \frac{n_1^{k_1-1}}{(k_1-1)!} \cdots m_r^{d_r} \frac{n_r^{k_r-1}}{(k_r-1)!} \right) q^N 
\end{align*}
The coefficients are given by sums over all elements in $P_r(N)$ and are therefore invariant under the action of $\rho$, which implies relations among $\g$. These relations can be written down nicely in terms of generating series. Recall that we denote these by
\begin{align*}
\ggen\bi{X_1,\dots,X_r}{Y_1,\dots,Y_r} := \sum_{\substack {k_1,\dots,k_r\geq  1\\d_1,\dots,d_r \geq 0}} \mb{k_1,\dots,k_r}{d_1,\dots,d_r} X_1^{k_1-1} \frac{Y_1^{d_1}}{d_1!}\cdots X_r^{k_r-1} \frac{Y_r^{d_r}}{d_r!}\,.
\end{align*}
\begin{lemma}\label{lem:doubleggenfct}
    For $m\geq 1$ set
    \[ E_m(X) := e^{mX}\,\quad \text{ and } \quad L_m(X) := \frac{e^X q^m}{1-e^X q^m} \in \Q[[q,X]]\,. \]
    Then for all $r\geq 1$ we have the following two different expressions for the generating functions:
    \begin{align*}
    \mt{X_1,\dots,X_r}{Y_1,\dots,Y_r} &= \sum_{m_1>\dots>m_r>0} \prod_{j=1}^r E_{m_j}(Y_j) L_{m_j}(X_j) \\
    &=  \sum_{m_1>\dots>m_r>0} \prod_{j=1}^r E_{m_j}(X_{r+1-j}-X_{r+2-j}) L_{m_j}(Y_1+\dots+Y_{r-j+1})
    \end{align*}
    (with $X_{r+1} := 0$).
\end{lemma} 
\begin{proof}
    This follows from a change of variables similar to the one used in Lemma \ref{lem:genexpr} and given exactly by \eqref{eq:conj}.
\end{proof}

As a direct consequence of Lemma \ref{lem:doubleggenfct} we obtain the following. 
\begin{theorem}(Partition relation) We have\footnote{The  generating series are  examples of a \emph{bimould}. In the language of moulds  the partition relation \eqref{eq:partition} states that the bimould of the generating series $\mathfrak{g}$ is swap invariant.} 
    \begin{equation} \label{eq:partition}
    \mt{X_1,\dots,X_r}{Y_1,\dots,Y_r}  = \mt{ Y_1 + \dots + Y_r ,\dots ,Y_1 + Y_2, Y_1}{X_r , X_{r-1} - X_r, \dots , X_{1}-X_{2}}\,. 
    \end{equation} 
\end{theorem} 
Comparing the coefficients of both sides yields linear relations among $\g$, which are exactly those coming from the above interpretation as a sum over all partitions of $N$ for the coefficient of $q^N$. 
\begin{corollary}(Partition relation in depth one and two) \label{cor:partitionlowdepth}
    For $k,k_1,k_2 \geq 1$ and $d,d_1,d_2 \geq 0$ we have the  following relations in length one and two
    \begin{align*}
    (k-1)!\mb{k}{d} &= d!\mb{d+1}{k-1} \,,\\
    \mb{k_1,k_2}{d_1,d_2} &= \sum_{\substack{0 \leq j \leq d_1\\0 \leq i \leq k_2-1}} (-1)^i \frac{d_1!(d_2+j)!}{(k_1-1)!\,i!\,j!\,(k_2-1-i)!} \mb{d_2+1+j \,,d_1+1-j}{k_2-1-i \,, k_1-1+i} \,.
    \end{align*} 
\end{corollary}

\subsection{Double shuffle relations for $\g$}
In the following we will describe the algebraic structure of the $\g$. For this we will again use the notion of quasi-shuffle algebras. Recall that we introduced the Quasi-shuffle algebra $(\Q \langle A_z \rangle, \gqsh)$, where the product $\gqsh$ is the quasi-shuffle product induced by the following product on $\Q A_z$
\begin{align}\label{eq:defgagaindiamond}
z_{k_1 }\gdiamond z_{k_2}&=     z_{k_1+k_2}+ \sum_{j=1}^{k_1+k_2-1}\left( \lambda^j_{k_1,k_2}+ \lambda^j_{k_2,k_1} \right) z_j\,,
\end{align}
where the rational numbers $\lambda^j_{k_1,k_2}$ are given by
\begin{align*}
\lambda^j_{k_1,k_2}  = (-1)^{k_2-1} \binom{k_1+k_2-1-j}{k_1-j}  \frac{B_{k_1+k_2-j}}{(k_1+k_2-j)!} \,.
\end{align*}
We then viewed $\g$ as an algebra homomorphism from $(\Q \langle A_z \rangle, \gqsh)$ to $\gs$, by sending $z_{k_1}\cdots z_{k_r}$ to $\g(k_1,\dots,k_r)$.  This setup can be generalized in an obvious way for the double-indexed version of $\g$. For this we first define the double-indexed version of $A_z = \left\{z_1,z_2,\dots \right\}$ by 
\begin{align*}
\azbi := \left\{ \ai{k}{d} \mid k\geq 1, d\geq 0 \right\}\,.
\end{align*}
We can view $A_z$ as a subset of $\azbi$ by identifying $z_k$ with $\ai{k}{0}$.
On $\Q \azbi$ we define the following product	\begin{align}\label{eq:defgbidiamond}
\ai{k_1}{d_1} \gdiamond \ai{k_2}{d_2} &:=     \ai{k_1+k_2}{d_1+d_2}+ \sum_{j=1}^{k_1+k_2-1}\left( \lambda^j_{k_1,k_2}+ \lambda^j_{k_2,k_1} \right) \ai{j}{d_1+d_2}\,.
\end{align}
Clearly \eqref{eq:defgbidiamond} reduces to \eqref{eq:defgagaindiamond} in the case $d_1=d_2=0$. One can check that this is again a commutative and associative product and we therefore obtain a quasi-shuffle algebra  $(\Q \langle \azbi \rangle, \gqsh)$ (Theorem \ref{thm:quasishufflealgebra}). For convenience we write for a word in the alphabet $\azbi$ \[\ai{k_1,\dots,k_r}{d_1,\dots,d_r} := \ai{k_1}{d_1}\cdots \ai{k_r}{d_r}\,.\]
Again we can view $\g$ as a $\Q$-linear map defined as follows
\begin{align*}
\g: \Q \langle \azbi \rangle &\longrightarrow \mz_q\\
\ai{k_1,\dots,k_r}{d_1,\dots,d_r}  &\longmapsto \mb{k_1,\dots,k_r}{d_1,\dots,d_r}\,.
\end{align*}

\begin{proposition}\label{prop:gqshhom} The map $\g$ is an algebra homomorphism from $(\Q \langle \azbi \rangle, \gqsh)$ to $\mz_q$.
\end{proposition}
\begin{proof}
This is a direct consequence of Lemma \ref{lem:fmtoFM} and 
\begin{align*}
m^{d_1} \frac{P_{k_1}(X)}{(1-X)^{k_1}} \cdot m^{d_2} \frac{P_{k_2}(X)}{(1-X)^{k_2}} =& m^{d_1+d_2} \frac{P_{k_1+k_2}(X)}{(1-X)^{k_1+k_2}}\\
&+ \sum_{j=1}^{k_1+k_2-1}\left( \lambda^j_{k_1,k_2}+ \lambda^j_{k_2,k_1} \right) m^{d_1+d_2} \frac{P_{j}(X)}{(1-X)^{j}}\,,
\end{align*}
which follows from Lemma \ref{lem:rproduct}.
\end{proof}

The product $\gqsh$ on $\Q \langle \azbi \rangle$ can be viewed as the analogue of the stuffle product. In the following we want to describe a product $\gshu$, which gives an analogue of the shuffle product. To define this product we will use the partition relations described in the previous subsection. Consider the following generating series of words in $\azbi$
\begin{align*}
\mathfrak{w} \bi{X_1,\dots,X_r}{Y_1,\dots,Y_r} := \sum_{\substack {k_1,\dots,k_r\geq  1\\d_1,\dots,d_r \geq 0}} \ai{k_1,\dots,k_r}{d_1,\dots,d_r} X_1^{k_1-1} \frac{Y_1^{d_1}}{d_1!}\cdots X_r^{k_r-1} \frac{Y_r^{d_r}}{d_r!} \in \Q \langle \azbi \rangle[[X_1,Y_1,\dots,X_r,Y_r]]\,.
\end{align*}
\begin{definition} We define the $\Q$-linear map $P: \Q \langle \azbi \rangle \rightarrow \Q \langle \azbi \rangle$ on the generators by 
\begin{align*}
\sum_{\substack {k_1,\dots,k_r\geq  1\\d_1,\dots,d_r \geq 0}} P\left(\ai{k_1,\dots,k_r}{d_1,\dots,d_r}\right) X_1^{k_1-1} \frac{Y_1^{d_1}}{d_1!}\cdots X_r^{k_r-1} \frac{Y_r^{d_r}}{d_r!}  :=  \mathfrak{w} \bi{ Y_1 + \dots + Y_r ,\dots ,Y_1 + Y_2, Y_1}{X_r , X_{r-1} - X_r, \dots , X_{1}-X_{2}} \,.
\end{align*}
\end{definition}
Notice that $P$ is an involution and that it corresponds exactly to the partition relation, i.e. we have for example $P\left( \ai{k}{d} \right)  = \frac{d!}{(k-1)!}\ai{d+1}{k-1}$. In particular we get the following. 
\begin{proposition} The map $\g$ is $P$-invariant, i.e. we have $\g(P(w)) = \g(w)$ for all $w \in \Q \langle \azbi \rangle$.
\end{proposition}
Using the map $P$ we can define a new product on the space $\Q \langle \azbi \rangle$ as follows. 
\begin{definition} We define the product $\gshu$ for $w,v \in \Q \langle \azbi \rangle$ by
\begin{align*}
w \,\gshu\, v := P\left( P(w)\, \gqsh \,P(v) \right)\,. 
\end{align*}
\end{definition}

\begin{proposition} \label{prop:gshuhom}Equipped with the product $\gshu$ the space $\Q \langle \azbi \rangle$ becomes a commutative $\Q$-algebra and $\g$ gives a $\Q$-algebra homomorphism $\g: (\Q \langle \azbi \rangle, \gshu) \rightarrow \mz_q$. 
\end{proposition}
\begin{proof} The first statement follows directly from the commutativity and associativity of $\gqsh$. The second statement follows from the $P$-invariance of $\g$.
\end{proof}

\begin{theorem}(Double shuffle relations for $\g$). For all $w,v \in \Q \langle \azbi \rangle$ we have 
\begin{align*}
\g(w \,\gshu \,v - w \,\gqsh \,v) = 0\,.
\end{align*}
\end{theorem}
\begin{proof}
This is a consequence of Proposition \ref{prop:gqshhom} and \ref{prop:gshuhom}.
\end{proof}

We will now give an example, which indicates why $\gshu$ could be seen as an analogue of the shuffle product.

\begin{ex} First we give the following example for the map $P$  (Using  Corollary \ref{cor:partitionlowdepth})
\begin{align*}
P\left( \ai{1,1}{1,2} \right) &= 2\ai{3,2}{0,0} + 6\ai{4,1}{0,0}\,,\\
P\left( \ai{1,1}{2,1} \right) &= 2\ai{2,3}{0,0} + 4\ai{3,2}{0,0} + 6\ai{4,1}{0,0} \,.
\end{align*}
Using this we obtain 
\begin{align*}
\ai{2}{0} \,\gshu\,  \ai{3}{0}  &= P\left( P\left(\ai{2}{0}  \right) \, \gqsh \, P\left( \ai{3}{0} \right)  \right) = P\left( \ai{1}{1}  \, \gqsh \,  \frac{1}{2}\ai{1}{2}  \right) \\
&= \frac{1}{2}P\left( \ai{1,1}{1,2} +  \ai{1,1}{2,1}  + \ai{2}{3} - \ai{1}{3} \right) \\
&=   \ai{2,3}{0,0}  + 3 \ai{3,2}{0,0} + 6 \ai{4,1}{0,0} + 3 \ai{4}{1} - 3 \ai{4}{0}\,.
\end{align*}
Together with 
\begin{align*}
\ai{2}{0} \,\gqsh\,  \ai{3}{0} &=  \ai{2,3}{0,0} +\ai{3,2}{0,0}+\ai{5}{0} -\frac{1}{12} \ai{3}{0}\,
\end{align*}
this gives the following linear relation among the $q$-series $\g$
\begin{align}\label{eq:gdoubledshrelex}
\mb{5}{0}= 2\mb{3,2}{0,0}+6\mb{4,1}{0,0} + \frac{1}{12}\mb{3}{0} - 3 \mb{4}{0}+3 \mb{4}{1}\,.
\end{align}
Notice that this gives exactly the same formula as in Example \ref{ex:gdshdep1}, but the new setup also allows us to consider these products in higher depths. 
\end{ex}
In general it was shown in \cite{Ba8} and \cite{Zu5} that the product $\gshu$ corresponds to the classical stuffle and shuffle products of multiple zeta values after considering the limit $q\rightarrow 1$.

In the case of multiple zeta values, we saw that conjecturally the extended double shuffle relations give all relations. A similar conjecture exists for the modified $q$-analogues $\g$ given as follows.

\begin{conjecture}\label{conj:allrelg} All linear relations among $\g$ are a consequence of the double shuffle relations and the partition relation, i.e. we have
\begin{align*}
\ker(\g) = \langle w \, \gshu \,v - w \,\gqsh\, v \mid w,v\in \Q \langle \azbi \rangle \rangle_\Q  + \langle P(w) - w \mid w \in \Q \langle \azbi \rangle \rangle_\Q  \,.
\end{align*}
\end{conjecture}
This conjecture was first proposed by the author in a talk in Bristol in 2014.\footnote{\url{https://www.math.uni-hamburg.de/home/bachmann/talks/bristol14_multdiv.pdf}}

In \eqref{eq:gdoubledshrelex}, we see that the double-indexed $\mb{4}{1}$ can be written in terms of single-indexed $\g$. As mentioned in the beginning, this is the case in general, i.e., any double-indexed $\g$ can be expressed in the original (single-indexed) $\g$ and therefore $\mz_q = \gs$.

\begin{theorem}[{Hirose--Maesaka--Watanabe \cite[Theorem~1.1]{HMW}}] \label{thm:hmwqspaces}We have $\mz_q = \gs$ and $\mz^\circ_{q,1}=\mz_{q,1}$. More precisely any $\mb{k_1,\dots,k_r}{d_1,\dots,d_r}$ can be expressed as a linear combination of $\g(s_1,\dots,s_l)$ with $1 \leq l \leq r + d_1+\dots + d_r$ and $s_1+\dots+s_l \leq k_1+\dots+k_r+d_1+\dots+d_r$.
\end{theorem}
The first equality was a conjecture of the author, first stated in \cite[Conjecture 4.3]{Ba7}, for a long time. It had been checked up to weight $9$ (using the double shuffle relations) and was proven for all weights in depth one (as a consequence of Theorem \ref{thm:gclosedunderqdq}), together with a few special cases in depth two (\cite{Ba7}). The second equality was another conjecture of the author based on numerical experiments. Recently Hirose, Maesaka and Watanabe proved both equalities in general, and in a much stronger form: they construct for every generator in the larger spaces an explicit expression with integer coefficients in terms of the smaller families. Their proof first establishes the corresponding identities for suitable finite $q$-analogues, which can be described recursively by generating series, and then passes to the limit.

\begin{ex} Using the double shuffle relations for $\g$ one can show that 
\begin{align}\label{eq:bdmdappl1}\begin{split}
\mb{3, 2}{1, 0} = &-\frac{1}{8}\g(6)-\g(5, 1)+\frac{5}{2}\g(4, 2)-\frac{3}{4}\g(4, 1, 1)-\frac{5}{4}\g(3, 2, 1)\\&-\frac{1}{1440} \g(2)+\frac{7}{96}\g(3)-\frac{1}{3}\g(4)-\frac{7}{48}\g(3, 1)+\frac{1}{2}\g(5)+\frac{3}{4}\g(4, 1)+\frac{5}{8}\g(3, 2)\,.
\end{split}
\end{align}
\end{ex}

We saw in Proposition \ref{prop:doublegasqana} that $	\mb{3, 2}{1, 0} $ can be viewed as a $q$-analogue of $\zeta(2,2)$, since 
\begin{align*}
\lim_{q\rightarrow 1} (1-q)^5\mb{3, 2}{1, 0} = \zeta(2,2)\,.
\end{align*}
But according to Theorem \ref{thm:hmwqspaces} the $q$-series is actually something of weight $6$ and not $5$. But clearly $\lim_{q\rightarrow 1} (1-q)^6\mb{3, 2}{1, 0} = 0$, and therefore \eqref{eq:bdmdappl1} implies a relation among multiple zeta values given by
\begin{align*}
\zeta(6) = 20 \zeta(4,2) - 8 \zeta(5,1) - 6 \zeta(4,1,1) - 10 \zeta(3,2,1)\,.
\end{align*} 
In general, almost all double-indexed $\mb{k_1,\dots,k_r}{d_1,\dots,d_r}$ vanish when considering the limit
\begin{align*}
\lim_{q\rightarrow 1} (1-q)^{k_1+\dots+k_r+d_1+\dots+d_r}.
\end{align*}
Theorem \ref{thm:hmwqspaces} therefore gives relations among multiple zeta values, which conjecturally give all relations among multiple zeta values.

\section{Derivatives of multiple Eisenstein series}\label{sec:mesderiv}

We now give the derivative formula from \cite{BKM}. The formula in arbitrary
depth uses the Drop1 operator from Section~\ref{sec:dropone}. We will also give
explicit formulas in depths one and two and prove that $\gstwo$ is closed under
$q\frac{d}{dq}$. Recall the algebra $\HH^{\geq2}$ from that section and extend
the map $\aG$ linearly to it by
$\aG(z_{k_1}\cdots z_{k_r})=\aG_{k_1,\dots,k_r}$.

The proof uses the formal multiple Eisenstein algebra which will be developed
in Chapter~\ref{sec:formalspaces}. For now, we write $\fmes$ for this algebra
and $\Gf(w)$ for the element associated with a word $w$. Only its formal
derivative and its analytic and $q$-series realizations are needed below.

\subsection{The derivative formula}

Although the recursive definition of Drop1 is needed in general, the derivative
formula only produces words containing at most one $z_1$. For these words there
is a more direct formula in \cite[Theorem~4.7]{BKM}.

We will use the notation $z_{k_1,\dots,k_r}=z_{k_1}\cdots z_{k_r}$. For the
depth-two formula we need the following two special cases.
\begin{lemma}\label{lem:droponeDepthTwo}
For $p,q\geq2$, put $M=p+q+1$. Then
\begin{align}
\dropone(z_pz_qz_1)={}&z_pz_{q+1}-z_qz_{p+1}
+\sum_{a=q+1}^{p+1}z_a z_{M-a}-\sum_{a=p+2}^{q}z_a z_{M-a}\notag\\
&-\sum_{\substack{a+b=q+1\\a,b\geq2}}
\bigl(z_a z_pz_b+z_pz_a z_b\bigr)
+\sum_{\substack{a+b=p+1\\a,b\geq2}}z_qz_a z_b,
\label{eq:droponepqone}\\
\dropone(z_pz_1z_q+z_pz_qz_1)={}&z_{p+1}z_q+z_pz_{q+1}\notag\\
&-\sum_{\substack{a+b=p+1\\a,b\geq2}}z_a z_bz_q
-\sum_{\substack{a+b=q+1\\a,b\geq2}}z_pz_a z_b.
\label{eq:droponeponeq}
\end{align}
\end{lemma}
\begin{proof}
Specializing the formula for a word with exactly one $z_1$ in
\cite[Theorem~4.7]{BKM} gives
\begin{align*}
\dropone(z_pz_qz_1)
={}&z_pz_qz_1+(z_qz_p+z_{p+q})\sh z_1
-z_p\ast(z_q\sh z_1)+z_pz_q\ast z_1,\\
\dropone(z_pz_1z_q)
={}&z_pz_1z_q-(z_p\ast z_q)\sh z_1+z_p\ast(z_q\sh z_1).
\end{align*}
To expand these expressions, one uses, for $m_1,\dots,m_r\geq2$,
\begin{align*}
z_{m_1}\cdots z_{m_r}\sh z_1={}&z_1z_{m_1}\cdots z_{m_r}
+2\sum_{i=1}^r z_{m_1}\cdots z_{m_i}z_1z_{m_{i+1}}\cdots z_{m_r}\\
&+\sum_{i=1}^r\sum_{a=2}^{m_i-1}
z_{m_1}\cdots z_{m_{i-1}}z_a z_{m_i+1-a}z_{m_{i+1}}\cdots z_{m_r}.
\end{align*}
All words containing $z_1$ cancel. For the first expression, the depth-three terms
are the last two sums in \eqref{eq:droponepqone}. For $2\leq a\leq M-2$, the
coefficient of $z_a z_{M-a}$ in the remaining depth-two sums is $-1$ for
$p+1\leq a\leq q$, is $1$ for $q+1\leq a\leq p$, and is zero otherwise.
Combining this with the isolated terms
$z_pz_{q+1}-z_qz_{p+1}+z_{p+1}z_q$ gives \eqref{eq:droponepqone}.
Adding the two expanded expressions gives
\eqref{eq:droponeponeq}.
\end{proof}

Recall that our sign convention for double shuffle differences is
\begin{align*}
\ds(u,v)=u\sh v-u\ast v.
\end{align*}
For $w\in\HH^{\geq2}$ define
\begin{align}\label{eq:defmesdertheta}
\theta(w)=\dropone\bigl(\ds(w,z_2)\bigr)\in\HH^{\geq2}.
\end{align}
Notice that \cite{BKM} uses the opposite sign for the double shuffle difference. This is why there is no minus sign in \eqref{eq:defmesdertheta}.

\begin{theorem}[{\cite[Main Theorem~D(ii)]{BKM}}]\label{thm:mesderformula}
For every $w\in\HH^{\geq2}$ we have
\begin{align}\label{eq:mesderformula}
(2\pi i)\frac{d}{d\tau}\aG(w)=\aG\bigl(\theta(w)\bigr).
\end{align}
In particular, $\mes$ is closed under the derivation $\GD=(2\pi i)\frac{d}{d\tau}$.
\end{theorem}
\begin{proof}
The formal multiple Eisenstein series developed in Chapter~\ref{sec:formalspaces} carry the derivation $\GD$. The main formal identity is \cite[Main Theorem~D(i)]{BKM}:
\begin{align}\label{eq:fmesdroponeidentity}
\GD\Gf(w)=\Gf\bigl(\theta(w)\bigr)
\qquad\bigl(w\in\HH^{\geq2}\bigr).
\end{align}
We prove this identity. If
$\mathbf{a}=(a_1,\dots,a_p)$ and
$\mathbf{b}=(b_1,\dots,b_s)$ have all entries at least $2$, then the formula
for $\dropone$ on a word containing one $z_1$, together with the swap
relation, gives
\begin{align}\label{eq:droponeformalswap}
&\Gf\bigl((\dropone-\operatorname{id})
(z_{a_1}\cdots z_{a_p}z_1z_{b_1}\cdots z_{b_s})\bigr)\notag\\
&\quad=
\Gf\bi{\mathbf{a},\mathbf{b}}
{\{0\}^{p-1},1,\{0\}^{s}}
-\Gf\bi{\mathbf{a},\mathbf{b}}
{\{0\}^{p},1,\{0\}^{s-1}},
\end{align}
where the last term is zero if $s=0$ (see \cite[Lemma~4.8]{BKM}). Now write $w=z_{k_1}\cdots z_{k_r}$. Modulo the
span of words whose indices are all at least $2$, the shuffle product is
\begin{align*}
w\sh z_2\equiv
2\sum_{1\leq j<i\leq r+1}k_j
z_{k_1}\cdots z_{k_j+1}\cdots z_{k_{i-1}}z_1
z_{k_i}\cdots z_{k_r}.
\end{align*}
Since $\dropone$ is the identity on $\HH^{\geq2}$, we may apply
\eqref{eq:droponeformalswap}. For each fixed $j$, the two terms telescope as
$i$ runs from $j+1$ to $r+1$. Hence
\begin{align*}
&\Gf\bigl(\dropone(w\sh z_2)-w\sh z_2\bigr)\\
&\quad=2\sum_{j=1}^rk_j
\Gf\bi{k_1,\dots,k_j+1,\dots,k_r}
{0,\dots,0,1,0,\dots,0}
=2\GD\Gf(w),
\end{align*}
where the $1$ in the lower row occurs in the $j$-th position. On the other
hand, $\dropone(w\ast z_2)=w\ast z_2$, and
\eqref{eq:fmesDonfilzero} gives
\begin{align*}
\Gf\bigl(\theta(w)\bigr)
&=\Gf\bigl(\dropone(w\sh z_2)-w\ast z_2\bigr)\\
&=2\GD\Gf(w)+\Gf(w\sh z_2-w\ast z_2)
=\GD\Gf(w).
\end{align*}
This proves~\eqref{eq:fmesdroponeidentity}.

The bi-multiple Eisenstein series constructed in \cite{BKM} give a realization of the formal algebra as holomorphic functions on $\Ha$. This realization sends $\Gf(k_1,\dots,k_r)$ to $\aG_{k_1,\dots,k_r}$ when all $k_j\geq2$ and intertwines the formal derivation with $(2\pi i)\frac{d}{d\tau}$. Applying it to \eqref{eq:fmesdroponeidentity} gives \eqref{eq:mesderformula}. Since $\theta(w)\in\HH^{\geq2}$, the right-hand side belongs to $\mes$.
\end{proof}

For $w=z_k$, formula \eqref{eq:defmesdertheta} gives
\begin{align*}
\theta(z_k)=(2k-1)z_{k+2}-\sum_{j=2}^k(k+j-1)z_{k+2-j}z_j-z_2z_k,
\end{align*}
which gives Theorem~\ref{thm:turanthm}. We also obtain the following short
formula for repeated twos.
\begin{proposition}[{\cite[Example~4.11]{BKM}}]\label{prop:mesderrepeatedtwo}
For $r\geq1$ we have
\begin{align*}
\GD\aG_{\{2\}^r}
=3\aG_{\{2\}^r}\aG_2-(r+1)(2r+3)\aG_{\{2\}^{r+1}}.
\end{align*}
\end{proposition}
\begin{proof}
In this case the Drop1 calculation in \cite[Lemma~4.6]{BKM} gives
\begin{align*}
\theta(z_2^r)=3z_2^r\ast z_2-(r+1)(2r+3)z_2^{r+1}.
\end{align*}
Now apply Theorem~\ref{thm:mesderformula} and the stuffle product formula.
\end{proof}

We finish the depth-two calculation announced in Section~\ref{subsec:mesoverview}.
\begin{proof}[Proof of Theorem~\ref{thm:mesderdepthtwo}]
Put $N=k+l+2$. A direct shuffle calculation gives
\begin{align}\label{eq:mesderdepthtwods}
\ds(z_kz_l,z_2)={}&-z_{2,k,l}-z_{k,2,l}-z_{k,l+2}-z_{k+2,l}\notag\\
&+\sum_{a=2}^{k}(a-1)z_{a,k+2-a,l}
+\sum_{a=2}^{l}(a-1)z_{k,a,l+2-a}\notag\\
&+k\sum_{a=2}^{l-1}z_{k+1,a,l+1-a}
+2k\bigl(z_{k+1}z_1z_l+z_{k+1}z_lz_1\bigr)
+2l z_kz_{l+1}z_1.
\end{align}
The terms without $z_1$ are fixed by $\dropone$. By
\eqref{eq:droponeponeq} and \eqref{eq:droponepqone}, respectively, the last two
parts become
\begin{align*}
&2k\left(z_{k+2}z_l+z_{k+1}z_{l+1}
-\sum_{\substack{a+b=k+2\\a,b\geq2}}z_a z_bz_l
-\sum_{\substack{a+b=l+1\\a,b\geq2}}z_{k+1}z_a z_b\right),\\
&2l\left(z_kz_{l+2}-z_{l+1}z_{k+1}
+\sum_{a=l+2}^{k+1}z_a z_{N-a}-\sum_{a=k+2}^{l+1}z_a z_{N-a}\right.\\
&\hspace{1.5cm}\left.
-\sum_{\substack{a+b=l+2\\a,b\geq2}}(z_a z_kz_b+z_kz_a z_b)
+\sum_{\substack{a+b=k+1\\a,b\geq2}}z_{l+1}z_a z_b\right).
\end{align*}
Substitution into \eqref{eq:mesderdepthtwods} gives \eqref{eq:mesderdepthtwo}.
For example, if $a+b=k+2$, then the coefficient of $z_a z_bz_l$ is
$(a-1)-2k=-(k+b-1)$. The case $a+b=l+2$ is analogous. Applying $\aG$ and
Theorem~\ref{thm:mesderformula} gives the result.
\end{proof}

\subsection{Derivatives of Okounkov's $q$-analogues}

We finally prove Theorem~\ref{thm:g2closedunderdif}. Recall that $\gstwo$ is spanned by $1$ and the $q$-series $\g(k_1,\dots,k_r)$ with all $k_j\geq2$. This is the space denoted by $\mathsf{qMZV}$ in \cite{BKM} and introduced by Okounkov in \cite{Ok}.

\begin{proof}[Proof of Theorem~\ref{thm:g2closedunderdif}]
Write $\operatorname{ev}_G:\fmes\rightarrow\Q\llbracket q\rrbracket$ for the algebraic realization given by the combinatorial bi-multiple Eisenstein series of \cite{BaBu}, and set
\begin{align*}
G(k_1,\dots,k_r)=\operatorname{ev}_G\bigl(\Gf(k_1,\dots,k_r)\bigr).
\end{align*}
This realization satisfies
\begin{align}\label{eq:cmesintertwineder}
\qdq\operatorname{ev}_G(f)=\operatorname{ev}_G(\GD f)\qquad(f\in\fmes).
\end{align}
Moreover, the proof of \cite[Lemma~4.12]{BKM}, using \cite[Proposition~6.15]{BaBu}, gives the triangular formula
\begin{align}\label{eq:cmestriangular}
G(k_1,\dots,k_r)=\g(k_1,\dots,k_r)
+\sum_{\substack{0\leq s<r\\l_1,\dots,l_s\geq2\\l_1+\dots+l_s\leq k_1+\dots+k_r}}
\lambda^{k_1,\dots,k_r}_{l_1,\dots,l_s}\g(l_1,\dots,l_s)
\end{align}
for some rational numbers $\lambda^{k_1,\dots,k_r}_{l_1,\dots,l_s}$, where $\g(\emptyset)=1$. The correction terms have smaller depth and no larger weight, so the sum is finite. It follows by induction on the depth that $\gstwo$ is also spanned by the $G(k_1,\dots,k_r)$ with all indices at least $2$.

For $w\in\HH^{\geq2}$, equations \eqref{eq:cmesintertwineder} and \eqref{eq:fmesdroponeidentity} now give
\begin{align*}
\qdq\operatorname{ev}_G\bigl(\Gf(w)\bigr)
=\operatorname{ev}_G\bigl(\GD\Gf(w)\bigr)
=\operatorname{ev}_G\bigl(\Gf(\theta(w))\bigr)\in\gstwo.
\end{align*}
This proves the claim.
\end{proof}

\section{The \texorpdfstring{$\sltwo$}{sl2}-algebra structure}\label{sec:messl2}
In this section, we want to explain the proof of \cref{thm:mesgraded,thm:messl2algebra}, which was given by the author and Kanno in \cite{BKa}. Since the derivative $\GD$ was already treated in the previous section, it remains to show that $\GW$ and $\Gdelta$ are well-defined derivations and that the commutator relations hold. The main idea is to compare the usual summation of multiple Eisenstein series with sums over a family of orders on the lattice, which depend on an integer $N$ and which come from the action of $\SLZ$. This comparison needs an analytic estimate, whose proof is quite technical and which we will not give here. Everything else is elementary and only uses results from the previous sections and the formal multiple Eisenstein series of Chapter~\ref{sec:formalspaces}.

\subsection{Polynomials in \texorpdfstring{$\aG_2$}{G2}}
We denote by $\HH^{\mathrm{ac}}\subset\HH^{\geq2}$ (``ac'' for absolutely convergent) the subspace spanned by $1$ and the words $z_{k_1}\cdots z_{k_r}$ with $k_1\geq3$ and $k_2,\dots,k_r\geq2$. By the definition of the stuffle product, this is a subalgebra of $\HH^{\geq2}_\ast$, and we set $\mes^{\mathrm{ac}}=\aG(\HH^{\mathrm{ac}})$, i.e. $\mes^{\mathrm{ac}}$ is the subalgebra of $\mes$ spanned by $1$ and the absolutely convergent multiple Eisenstein series.
\begin{proposition}[{\cite[Proposition~2.1]{BKa}}]\label[proposition]{prop:mesHadmpolynomial}
Every $w\in\HH^{\geq2}$ can be written uniquely as
\begin{align}\label{eq:mesHadmdecomposition}
w=\sum_{j=0}^d a_j\ast z_2^{\ast j}\qquad(a_j\in\HH^{\mathrm{ac}})\,,
\end{align}
where $z_2^{\ast j}$ denotes the $j$-th power with respect to $\ast$ and $z_2^{\ast0}=1$, i.e. we have $\HH^{\geq2}_\ast\cong\HH^{\mathrm{ac}}_\ast[z_2]$. If $w$ is homogeneous of weight $k$, then $a_j$ is homogeneous of weight $k-2j$.
\end{proposition}
\begin{proof}
Every word can be written uniquely as $z_2^mv$, where $v=1$ or $v$ starts with a letter $z_k$ with $k\geq3$. By induction on $m$, the definition of the stuffle product gives
\begin{align*}
z_2^{\ast m}\ast v=m!\,z_2^mv+\bigl(\text{words with fewer than $m$ letters $z_2$ at the beginning}\bigr)\,.
\end{align*}
Since there are only finitely many words of a fixed weight, the two bases $\{z_2^mv\}$ and $\{v\ast z_2^{\ast m}\}$ of the weight $k$ part of $\HH^{\geq2}$ are related by an invertible triangular matrix, which gives the statement.
\end{proof}
\noindent Now let $\bar\delta:\HH^{\geq2}\rightarrow\HH^{\geq2}$ be the linear map given by $\bar\delta(1)=0$ and
\begin{align*}
\bar\delta(z_{k_1}\cdots z_{k_r})=\begin{cases}-\frac12z_{k_2}\cdots z_{k_r},&k_1=2,\\0,&k_1>2\,,\end{cases}
\end{align*}
i.e. $\bar\delta$ is the map which should induce $\Gdelta$ on $\mes$. It is the restriction of the formal operator $\Gdelta$ in \eqref{eq:fmesdeltafilzero}.
\begin{lemma}[{\cite[Lemma~2.2]{BKa}}]\label[lemma]{lem:mesdeltabarderivation}
The map $\bar\delta$ is a derivation of $\HH^{\geq2}_\ast$, and with \eqref{eq:mesHadmdecomposition} we have
\begin{align*}
\bar\delta(w)=-\frac12\sum_{j=1}^d j\,a_j\ast z_2^{\ast(j-1)}\,,
\end{align*}
i.e. $\bar\delta=-\frac12\frac{\partial}{\partial z_2}$.
\end{lemma}
\begin{proof}
A term in the stuffle product of two words which starts with $z_2$ is obtained by taking the first letter $z_2$ of one of the two words without merging it. Removing this letter gives exactly $\bar\delta(u\ast v)=\bar\delta(u)\ast v+u\ast\bar\delta(v)$. Since $\bar\delta$ vanishes on $\HH^{\mathrm{ac}}$ and $\bar\delta(z_2)=-\frac12$, the second statement follows from the uniqueness in \cref{prop:mesHadmpolynomial}.
\end{proof}
Applying $\aG$ to \eqref{eq:mesHadmdecomposition}, we see that $\mes=\mes^{\mathrm{ac}}[\aG_2]$. But a priori it is not clear that the representation of an element of $\mes$ as a polynomial in $\aG_2$ is unique, and this is exactly the question whether $\Gdelta$ is well-defined.
\begin{proposition}[{\cite[Proposition~2.3]{BKa}}]\label[proposition]{prop:mesdeltaequivalent}
The map $\bar\delta$ preserves the kernel of $\aG$, i.e. it induces a well-defined derivation $\Gdelta$ on $\mes$, if and only if $\aG_2$ is transcendental over $\mes^{\mathrm{ac}}$. In this case, $\mes$ is a polynomial ring in $\aG_2$ over $\mes^{\mathrm{ac}}$ and $\Gdelta=-\frac12\frac{\partial}{\partial\aG_2}$.
\end{proposition}
\begin{proof}
If $\Gdelta$ is well-defined, then it vanishes on $\mes^{\mathrm{ac}}$ and satisfies $\Gdelta(\aG_2)=-\frac12$. Now suppose that $b_0+b_1\aG_2+\dots+b_d\aG_2^d=0$ with $b_j\in\mes^{\mathrm{ac}}$. Then applying $\Gdelta^d$ gives $d!\,(-\frac12)^db_d=0$, i.e. $b_d=0$, and repeating this argument gives $b_j=0$ for all $j$. Conversely, if $\aG_2$ is transcendental over $\mes^{\mathrm{ac}}$, then $\mes\cong\mes^{\mathrm{ac}}[T]$ and we can transport the derivation $-\frac12\partial/\partial T$ to $\mes$. By \cref{lem:mesdeltabarderivation}, this derivation sends $\aG(w)$ to $\aG(\bar\delta w)$, so $\bar\delta$ preserves the kernel of $\aG$.
\end{proof}
Therefore the following theorem, which we mentioned already in \Cref{subsec:mesoverview}, is equivalent to the well-definedness of $\Gdelta$, and we will prove it in \Cref{sec:messl2proofs}.
\begin{theorem}[{\cite[Main Theorem~A(iii)]{BKa}}]\label{thm:mespolynomialG2}
The algebra $\mes$ is a polynomial ring in $\aG_2$ over $\mes^{\mathrm{ac}}$, and $\Gdelta=-\frac12\frac{\partial}{\partial\aG_2}$.
\end{theorem}
Notice that the analogous statement for multiple zeta values is false, since $\zeta(2)^2=\frac52\zeta(4)$. So surprisingly, the Eisenstein series $\aG_2$ behaves like a free variable, while $\zeta(2)$ does not.

\subsection{Moving lattice orders}
We fix $\tau\in\Ha$ and write $\Lambda_\tau=\Z\tau+\Z$ as before. For an integer $N\geq1$ and $\lambda=m\tau+n\in\Lambda_\tau$ we put $H_N(\lambda)=Nm+n$ and define
\begin{align*}
P_N=\left\{\lambda\in\Lambda_\tau\mid H_N(\lambda)>0\right\}\cup\left\{m\tau+n\in\Lambda_\tau\mid H_N(m\tau+n)=0,\ m<0\right\}\,.
\end{align*}
As for $\succ$, we write $\lambda_1\succ_N\lambda_2$ if $\lambda_1-\lambda_2\in P_N$. For $\tau=i$, the sets $P_1$ and $P_3$ look as follows.
\begin{center}
\begin{tikzpicture}[scale=0.32]
\foreach \NN/\ss in {1/0,3/14}{
\begin{scope}[shift={(\ss,0)}]
\draw[dotted,step=1,color=gray,thin] (-4.9,-4.9) grid (4.9,4.9);
\draw[->,thick] (0,-5.5) -- (0,5.5) node[above] {$m$};
\draw[->,thick] (-5.5,0) -- (5.5,0) node[right] {$n$};
\foreach \mm in {-4,...,4}{\foreach \nn in {-4,...,4}{
\pgfmathtruncatemacro{\hval}{\NN*\mm+\nn}
\ifnum\hval>0 \fill[black] (\nn,\mm) circle (4pt); \fi
\ifnum\hval=0 \ifnum\mm<0 \fill[black] (\nn,\mm) circle (4pt); \fi\fi
}}
\node at (0,-6.8) {$N=\NN$};
\end{scope}}
\end{tikzpicture}
\end{center}
For two fixed lattice points, the orders $\succ_N$ and $\succ$ agree for $N$ large enough, since the sign of $H_N(\lambda_1-\lambda_2)=N(m_1-m_2)+(n_1-n_2)$ is then determined by $m_1-m_2$, and by $n_1-n_2$ if $m_1=m_2$. So the cones $P_N$ approach the set $P$ of positive lattice points as $N\rightarrow\infty$. For $k_1\geq3$ and $k_2,\dots,k_r\geq2$ we define
\begin{align*}
\aG^{(N)}_{k_1,\dots,k_r}(\tau)=\sum_{\substack{\lambda_1\succ_N\dots\succ_N\lambda_r\succ_N0\\\lambda_i\in\Lambda_\tau}}\frac{1}{\lambda_1^{k_1}\cdots\lambda_r^{k_r}}\,.
\end{align*}
These sums are nothing else than multiple Eisenstein series evaluated at another point.
\begin{lemma}[{\cite[Lemma~3.1]{BKa}}]\label[lemma]{lem:mesmovingorder}
Let $k_1\geq3$, $k_2,\dots,k_r\geq2$ and $k=k_1+\dots+k_r$. With $\tau_N=-1/(\tau-N)$ we have
\begin{align*}
\tau_N^k\aG_{k_1,\dots,k_r}(\tau_N)=\aG^{(N)}_{k_1,\dots,k_r}(\tau)\,.
\end{align*}
\end{lemma}
\begin{proof}
For a lattice point $a\tau_N+b\in\Lambda_{\tau_N}$ we put $\lambda=a-b(\tau-N)=m\tau+n$ with $(m,n)=(-b,a+Nb)$. This change of coordinates is given by a matrix in $\SLZ$, so it is a bijection, and we have $H_N(\lambda)=a$ and $a\tau_N+b=-\lambda/(\tau-N)$. The condition that $a>0$, or $a=0$ and $b>0$, becomes $H_N(\lambda)>0$, or $H_N(\lambda)=0$ and $m<0$, i.e. the order $\succ$ on $\Lambda_{\tau_N}$ becomes the order $\succ_N$ on $\Lambda_\tau$. Since $k_1\geq3$, both sums converge absolutely and we can reorder them, and the change of variables gives the factor $(-1)^k(\tau-N)^k=\tau_N^{-k}$.
\end{proof}
In other words, the moving orders come from the action of $T^NS$ on the lattice, and $\tau_N=(T^NS)^{-1}\tau$ (see \cite[Remark~3.2]{BKa}). Of course, we expect that $\aG^{(N)}_{k_1,\dots,k_r}\rightarrow\aG_{k_1,\dots,k_r}$ as $N\rightarrow\infty$. If all entries are at least $3$, this follows from the dominated convergence theorem, since then $\prod_{i=1}^r\sum_{\lambda\neq0}|\lambda|^{-k_i}$ gives a majorant independent of $N$, and for two fixed lattice points the two orders agree for large $N$. The difficulty comes from the entries equal to $2$, and for these one needs the following estimate.
\begin{theorem}[{\cite[Theorem~3.3]{BKa}}]\label{thm:mesmovingestimate}
For $k_1\geq3$ and $k_2,\dots,k_r\geq2$ there is an $M\geq0$ with
\begin{align*}
\aG^{(N)}_{k_1,\dots,k_r}(\tau)-\aG_{k_1,\dots,k_r}(\tau)=O\left(\frac{(1+\log(2+N))^M}{N^2}\right)\qquad(N\rightarrow\infty)
\end{align*}
locally uniformly in $\tau\in\Ha$.
\end{theorem}
The proof of this estimate is the analytic heart of \cite{BKa}. It needs more than ten pages of careful estimates of lattice sums of Eisenstein type, and we will not give it here.

\subsection{Proof of the main theorems}\label{sec:messl2proofs}
For a homogeneous $w\in\HH^{\geq2}$ of weight $k$ we put $\aG^{(N)}(w)(\tau)=\tau_N^k\aG(w)(\tau_N)$. By \cref{lem:mesmovingorder}, this agrees with the sums over the moving orders for $w\in\HH^{\mathrm{ac}}$, and since $\aG$ is an algebra homomorphism, we have $\aG^{(N)}(u\ast v)=\aG^{(N)}(u)\aG^{(N)}(v)$ for homogeneous $u,v$. For $w=z_2$, we use the transformation formula of $\aG_2$ from \Cref{sec:qmfsl2algebras}, which gives $\aG_2(-1/z)=z^2\aG_2(z)-\pi iz$. Together with $\aG_2(\tau-N)=\aG_2(\tau)$ we obtain
\begin{align}\label{eq:mesmovingG2}
\aG^{(N)}(z_2)(\tau)=\aG_2(\tau)-\frac{\pi i}{\tau-N}\,.
\end{align}
\begin{proposition}[{\cite[Proposition~4.1]{BKa}}]\label[proposition]{prop:mestwoterm}
For homogeneous $w\in\HH^{\geq2}$ there is an $M\geq0$ with
\begin{align*}
\aG^{(N)}(w)(\tau)=\aG(w)(\tau)+\frac{2\pi i}{\tau-N}\aG(\bar\delta w)(\tau)+O\left(\frac{(1+\log(2+N))^M}{N^2}\right)
\end{align*}
locally uniformly in $\tau\in\Ha$.
\end{proposition}
\begin{proof}
We write $w=\sum_ja_j\ast z_2^{\ast j}$ with homogeneous $a_j\in\HH^{\mathrm{ac}}$ as in \cref{prop:mesHadmpolynomial}. By the multiplicativity of $\aG^{(N)}$ and \eqref{eq:mesmovingG2} we have
\begin{align*}
\aG^{(N)}(w)=\sum_{j=0}^d\aG^{(N)}(a_j)\left(\aG_2-\frac{\pi i}{\tau-N}\right)^j\,.
\end{align*}
By \cref{lem:mesmovingorder} and \cref{thm:mesmovingestimate} we can replace $\aG^{(N)}(a_j)$ by $\aG(a_j)$ up to an error of the claimed size. Since $|\tau-N|$ grows like $N$ locally uniformly, the same is true for all terms containing $(\tau-N)^{-2}$ or a higher power, and we obtain
\begin{align*}
\aG^{(N)}(w)=\aG(w)-\frac{\pi i}{\tau-N}\sum_{j=1}^dj\,\aG(a_j)\aG_2^{j-1}+O\left(\frac{(1+\log(2+N))^M}{N^2}\right)\,.
\end{align*}
By \cref{lem:mesdeltabarderivation} we have $\aG(\bar\delta w)=-\frac12\sum_{j\geq1}j\,\aG(a_j)\aG_2^{j-1}$, which gives the statement.
\end{proof}
\noindent In particular, $\aG^{(N)}(w)\rightarrow\aG(w)$ locally uniformly for every homogeneous $w\in\HH^{\geq2}$. It is quite nice that the term $-\pi i/(\tau-N)$, which comes from the fact that $\aG_2$ is not modular, is exactly what detects the operator $\Gdelta$.

\begin{proof}[Proof of \cref{thm:mesgraded}]
We show the statement directly for the complex span of all multiple Eisenstein series. Suppose that $F_0+\dots+F_{k_0}=0$, where each $F_k$ is a complex linear combination of multiple Eisenstein series of weight $k$, and $F_{k_0}\neq0$. We choose homogeneous $w_k\in\C\otimes\HH^{\geq2}$ of weight $k$ with $\aG(w_k)=F_k$, where we extend $\aG$ and $\aG^{(N)}$ linearly. Evaluating the relation at $\tau_N$ and multiplying by $\tau_N^{k_0}$, we get
\begin{align*}
0=\aG^{(N)}(w_{k_0})(\tau)+\sum_{k=0}^{k_0-1}\tau_N^{k_0-k}\aG^{(N)}(w_k)(\tau)\,.
\end{align*}
By \cref{prop:mestwoterm}, all $\aG^{(N)}(w_k)$ are locally bounded and the first term converges to $F_{k_0}(\tau)$, while all other terms tend to zero, since $\tau_N\rightarrow0$. Therefore $F_{k_0}=0$, which is a contradiction. So the complex span of multiple Eisenstein series is graded by weight, and in particular $\mes=\bigoplus_{k\geq0}\mes_k$.
\end{proof}
Notice how the weights are separated in this argument: a multiple Eisenstein series of weight $k$ just gets the factor $\tau_N^k$, and $\tau_N\rightarrow0$.

\begin{proof}[Proof of \cref{thm:messl2algebra}]
By \cref{thm:mesgraded} the kernel of $\aG$ is spanned by homogeneous elements, i.e. $\GW$ is well-defined, and it is a derivation since $\mes_k\mes_l\subset\mes_{k+l}$ (\cref{prop:mesissubalgebra}). The operator $\GD$ is well-defined by \cref{thm:mesderformula}, and it is a derivation since it is a differential operator. Now let $w\in\ker\aG$ be homogeneous. Then $\aG^{(N)}(w)=0$ for all $N$, and if we multiply the expansion in \cref{prop:mestwoterm} by $\tau-N$ and let $N\rightarrow\infty$, we get $\aG(\bar\delta w)=0$. So $\bar\delta$ preserves $\ker\aG$, and by \cref{lem:mesdeltabarderivation} it induces a derivation $\Gdelta$ of $\mes$. By \cref{prop:mesdeltaequivalent}, this also proves \cref{thm:mespolynomialG2}.

For the commutator relations, we use the formal multiple Eisenstein series from Chapter~\ref{sec:formalspaces}. By \cref{thm:famessl2}, the algebra $\fames$ is an $\sltwo$-algebra and $\Gf(w)\mapsto\aG(w)$ gives a surjective algebra homomorphism $\fames\rightarrow\mes$. For $w\in\HH^{\geq2}$ homogeneous of weight $k$, the formal operators are given by $\GW\Gf(w)=k\Gf(w)$, $\Gdelta\Gf(w)=\Gf(\bar\delta w)$ by \eqref{eq:fmesdeltafilzero} and $\GD\Gf(w)=\Gf(\theta(w))$ by \eqref{eq:fmesdroponeidentity}. The relation $[\Gdelta,\GD]=\GW$ in $\fames$ therefore gives $\Gf([\bar\delta,\theta]w)=k\Gf(w)$, and applying the realization we get $\aG([\bar\delta,\theta]w)=k\aG(w)$. Together with \cref{thm:mesderformula} this gives
\begin{align*}
[\Gdelta,\GD]\aG(w)=\aG\bigl([\bar\delta,\theta]w\bigr)=\GW\aG(w)\,,
\end{align*}
and the other two commutator relations follow, since $\GD$ raises and $\Gdelta$ lowers the weight by two.
\end{proof}
In particular, the realization map $\fames\rightarrow\mes$ in \cref{thm:famessl2} commutes with all three operators, i.e. it is a surjective morphism of $\sltwo$-algebras. Its injectivity (\cref{conj:famesrealization}) is still open.

\subsection{Multiple Eisenstein series and \texorpdfstring{$q$}{q}-analogues}
We now want to explain the connection to the $q$-series $\g$ from \Cref{subsec:doubleindexedqana}. For $k\geq0$, we denote by $\fil{\mathrm W}{k}(\gstwo)$ the $\Q$-span of $1$ and all $\g(k_1,\dots,k_r)$ with $k_1,\dots,k_r\geq2$ and $k_1+\dots+k_r\leq k$, and we set $\fil{\mathrm W}{-1}(\gstwo)=0$. Since the product of two such $q$-series is a linear combination of $q$-series of this form of weight at most the sum of the two weights (see \cite{BK1}), this gives a filtration of the algebra $\gstwo$, and we write
\begin{align*}
\operatorname{gr}^{\mathrm W}\gstwo=\bigoplus_{k\geq0}\fil{\mathrm W}{k}(\gstwo)\big/\fil{\mathrm W}{k-1}(\gstwo)
\end{align*}
for the associated graded algebra.
\begin{theorem}[{\cite[Main Theorem~B(ii)]{BKa}}]\label{thm:mesgrqmzv}
The following map is a well-defined isomorphism of graded $\Q$-algebras:
\begin{align*}
\mes&\longrightarrow\operatorname{gr}^{\mathrm W}\gstwo\,,\\
\aG_{k_1,\dots,k_r}&\longmapsto\g(k_1,\dots,k_r)\bmod\fil{\mathrm W}{k_1+\dots+k_r-1}(\gstwo)\,.
\end{align*}
\end{theorem}
\noindent By \cite[Main Theorem~B(ii)]{BKa}, also the natural map from $\C\otimes_\Q\mes$ to the complex span of all multiple Eisenstein series is an isomorphism.
\begin{proof}[Sketch of proof]
Let $w\in\HH^{\geq2}$ be homogeneous of weight $k$. By \cref{thm:Fourier}, the difference
\begin{align*}
(-2\pi i)^{-k}\aG(w)-\g(w)
\end{align*}
is a complex linear combination of $q$-series $\g(l_1,\dots,l_s)$ with $l_1,\dots,l_s\geq2$ of weight smaller than $k$. If $\aG(w)=0$, this shows that $\g(w)$ lies in the complex span of $\fil{\mathrm W}{k-1}(\gstwo)$, and since $\g(w)$ has rational coefficients, a short linear algebra argument with finitely many Fourier coefficients shows that $\g(w)\in\fil{\mathrm W}{k-1}(\gstwo)$. So the map is well-defined, and it is clearly surjective. Conversely, if $\g(w)\in\fil{\mathrm W}{k-1}(\gstwo)$, then by the same triangularity $\aG(w)$ is a complex linear combination of multiple Eisenstein series of smaller weight, and the grading of the complex span shown in the proof of \cref{thm:mesgraded} gives $\aG(w)=0$. The compatibility with the products and the statement about the complex span follow in a similar way (see \cite[Section~4]{BKa}).
\end{proof}
This explains the remark after \cref{conj:mesdim}: the dimension conjecture for multiple Eisenstein series is the same as the one for the associated weight-graded algebra of $\gstwo$. Notice also that the statement about the complex span shows that the two definitions of $\mes^\C$ in \Cref{subsec:mesoverview} and \Cref{subsec:doubleindexedqana} agree.

\subsection{A conjecture on all relations}
Let us finish this chapter with a conjecture which describes all relations among multiple Eisenstein series. By \cref{thm:mesRrelations} and \cref{thm:mesderformula}, the ideal $\mathsf{DR}_\ast$ of $\HH^{\geq2}_\ast$ generated by all $\theta^n(R(u,v))$ with $n\geq0$ and $u,v\in\HH^{\geq2}$, with $R(u,v)$ as in \eqref{eq:mesRuv}, is contained in the kernel of $\aG$.
\begin{conjecture}[{\cite[Conjecture~1.8]{BKM}}]\label[conjecture]{conj:mesallrelations}
The kernel of $\aG:\HH^{\geq2}_\ast\rightarrow\mes$ is equal to $\mathsf{DR}_\ast$.
\end{conjecture}
This can be seen as an analogue of \cref{conj:extdshprecise} for multiple zeta values. Numerically, the ideal $\mathsf R_\ast$ generated by the $R(u,v)$ alone gives the number of relations predicted by \cref{conj:mesdim} up to weight $16$, but in weight $17$ one relation is missing. The ideal $\mathsf{DR}_\ast$ gives the predicted number of relations up to weight $20$ (see \cite{BKM}):
\begin{center}
\small
\setlength{\tabcolsep}{3pt}
\begin{tabular}{|c|c|c|c|c|c|c|c|c|c|c|c|c|c|c|c|}
    \hline
    weight $k$ & 6 & 7 & 8 & 9 & 10 & 11 & 12 & 13 & 14 & 15 & 16 & 17 & 18 & 19 & 20 \\ \hline
    \# of relations $\overset{?}{=}$ & 1 & 1 & 4 & 6 & 13 & 23 & 42 & 74 & 129 & 224 & 382 & 651 & 1098 & 1847 & 3089 \\ \hline
    relations in $\mathsf R_\ast$ & 1 & 1 & 4 & 6 & 13 & 23 & 42 & 74 & 129 & 224 & 382 & 650 & 1098 & 1845 & 3086 \\ \hline
    relations in $\mathsf{DR}_\ast$ & 1 & 1 & 4 & 6 & 13 & 23 & 42 & 74 & 129 & 224 & 382 & 651 & 1098 & 1847 & 3089 \\ \hline
\end{tabular}
\end{center}
In \cite[Main Theorem~E]{BKM} it is also shown that the quotient $\HH^{\geq2}_\ast/\mathsf{DR}_\ast$ is an $\sltwo$-algebra, which fits nicely with \cref{thm:messl2algebra}.

The proofs in the last two sections used the formal multiple Eisenstein algebra only through its derivations and realizations. In the next chapter, after introducing formal multiple zeta values, we return to the depth-two precursor of this algebra and then develop the formal framework in arbitrary depth.

%% file: chap_FormalSpaces.tex
\chapter{Formal and combinatorial multiple Eisenstein series}\label{sec:formalspaces}
\chaptermark{Formal and combinatorial MES}

In this chapter, we want to introduce formal versions of multiple zeta values
and multiple Eisenstein series, and then construct a $q$-series realization of
the latter. The idea behind all of these constructions is always the same:
instead of studying the objects themselves, we study formal symbols which
satisfy exactly the relations we expect, and then ask which realizations such
symbols have.

We start with the formal multiple zeta values, which are given by formal
symbols satisfying the extended double shuffle relations. After this, we come
back to depth two once more and introduce the formal double Eisenstein space,
which lifts the formal double zeta space of Chapter~\ref{sec:mdandmzv} by also
allowing lower indices. Using the stuffle product and the partition relation for
the double-indexed $q$-series of the previous chapter, we then define formal
multiple Eisenstein series in arbitrary depth. This algebra will have a lot of
structure: a formal constant term, an $\sltwo$-action and, with this, also
formal quasimodular, modular and cusp forms. We will give three realizations of
it, i.e. in $q$-series, in holomorphic functions and in real numbers. The first
of these is constructed in Section~\ref{sec:cmes} by the combinatorial multiple
Eisenstein series, where we also study their products, limits and derivatives.
At the end, in Section~\ref{sec:formalfinitemzv}, we come back to finite and
symmetric multiple zeta values, introduce their universal formal counterpart
and connect the depth-four part of finite multiple zeta values with the period
polynomials of Chapter~\ref{sec:mdandmzv}.

\section{Formal multiple zeta values}\label{sec:fmz}

Recall the extended double shuffle relations from
Theorem~\ref{thm:eds}. They are conjectured to give all relations among
multiple zeta values. This suggests defining formal multiple zeta values by
imposing precisely these relations.

Notice that Chapter~\ref{sec:mdandmzv} considered a formal object in depth two.
There the symbols $Z_k$, $Z_{k_1,k_2}$ and $P_{k_1,k_2}$ encode the
depth-two double shuffle relations, including separate symbols for products.
The algebra below is a quotient of the full stuffle algebra by all extended
double shuffle relations and contains the regularization parameter
$\fzeta(1)$. Thus the two constructions serve different purposes and will not
be identified.

Let $\mathfrak E$ be the ideal of the stuffle algebra $(\HH^1,\ast)$ generated
by
\begin{align}\label{eq:fmzedsideal}
w\ast v-w\sh v\qquad(w\in\HH^1,\ v\in\HH^0)\,.
\end{align}

\begin{definition}\label{def:fmz}
The algebra of \emph{formal multiple zeta values} is
\begin{align*}
\fmz=\quotient{(\HH^1,\ast)}{\mathfrak E}\,.
\end{align*}
For $k_1,\dots,k_r\geq1$, we denote the class of
$z_{k_1}\cdots z_{k_r}$ by $\fzeta(k_1,\dots,k_r)$.
\end{definition}

By definition, the elements $\fzeta(k_1,\dots,k_r)$ satisfy exactly the
extended double shuffle relations. The regularized multiple zeta values give a
surjective algebra homomorphism
\begin{align}\label{eq:fmzrealization}
\fmz&\longrightarrow\mz[T],\\
\fzeta(k_1,\dots,k_r)&\longmapsto
\zeta^\ast(k_1,\dots,k_r;T)\,.
\end{align}
The extended double shuffle conjecture predicts that this map is an
isomorphism. Notice that $\fzeta(1)$ is non-zero and maps to $T$. Thus it plays
the role of the regularization parameter.

\begin{problem}\label{prob:fmzduality}
Show that formal multiple zeta values satisfy the duality relation, i.e., for
every admissible index $\kk$,
\begin{align*}
\fzeta(\kk^\dagger)=\fzeta(\kk)\,.
\end{align*}
Equivalently, show that
\begin{align*}
\tau(v)-v\in\mathfrak E\qquad(v\in\HH^0)\,.
\end{align*}
This is related to Problem~\ref{prob:dualityeds}, where it is asked whether
$\tau(v)-v$ is already a linear combination of regularized extended double
shuffle relations.
\end{problem}

For $k_1,k_2\geq1$, the depth-two relations are
\begin{align}\label{eq:fmzdepthtwo}
\fzeta(k_1)\fzeta(k_2)
={}&\fzeta(k_1,k_2)+\fzeta(k_2,k_1)+\fzeta(k_1+k_2)\\
={}&\sum_{l_1+l_2=k_1+k_2}
\left(\binom{l_1-1}{k_1-1}+\binom{l_1-1}{k_2-1}\right)
\fzeta(l_1,l_2)
+\begin{cases}\fzeta(2),&k_1=k_2=1,\\0,&\text{otherwise.}\end{cases}
\end{align}
For example,
\begin{align*}
\fzeta(3)=\fzeta(2,1),\qquad
\fzeta(4)=\frac{2}{5}\fzeta(2)^2,
\qquad
\fzeta(6)=\frac{8}{35}\fzeta(2)^3\,.
\end{align*}
More generally, Euler's relation holds formally:
\begin{align}\label{eq:fmzeuler}
\fzeta(2m)=-\frac{B_{2m}}{2(2m)!}\bigl(-24\fzeta(2)\bigr)^m
\qquad(m\geq1)\,.
\end{align}
These are all extended double shuffle relations, so the proofs of the
corresponding identities for multiple zeta values carry over without any change.
See also \cite[Proposition~5.20]{BCK}.

The quotient of $\fmz$ obtained by setting $\fzeta(1)=0$ is the version of
formal multiple zeta values studied in detail by Burmester, Confurius and
K\"uhn in \cite[Definition~5.17]{BCK}. They construct the Goncharov--Brown
coaction on this quotient (see \cite[Theorem~5.40]{BCK}). Assuming the free odd
generation conjecture for the double shuffle Lie algebra, they also prove that
the formal multiple zeta values $\fzeta(k_1,\dots,k_r)$ with
$k_j\in\{2,3\}$ form a basis of this quotient (see \cite[Theorems~6.18 and~6.19]{BCK}).
Two explicit polar solutions of the double shuffle equations modulo products, due to \'Ecalle and Brown, are compared in \cite{MatT}.

\section{The formal double Eisenstein space}\label{sec:fdes}

The formal double zeta space of Chapter~\ref{sec:mdandmzv} only involves the
upper indices $k_1,k_2$. The double-indexed $q$-series of the previous chapter
suggest that one should also allow lower indices. This leads to the formal
double Eisenstein space introduced by the author, K\"uhn and Matthes in
\cite{BKuM}. It is a vector space in every weight, and not an algebra.\footnote{To avoid
a conflict with our notation $\mes$ for multiple Eisenstein series, we write
$\fdes_K$ for the space denoted by $\mathcal E_K$ in \cite{BKuM}. So one needs to be
careful when comparing the notation here with the one in their work.}

\begin{definition}[{\cite[Definition~2.1]{BKuM}}]\label{def:fdes}
For $K\geq1$, the \emph{formal double Eisenstein space} $\fdes_K$ is the
$\Q$-vector space spanned by the symbols
\begin{align*}
\Gde{k}{d},\qquad
\Gde{k_1,k_2}{d_1,d_2},\qquad
\Pde{k_1,k_2}{d_1,d_2}\,,
\end{align*}
where
\begin{align*}
k+d=k_1+k_2+d_1+d_2=K,\qquad
k,k_1,k_2\geq1,\quad d,d_1,d_2\geq0\,,
\end{align*}
subject to the relations
{\small
\begin{align}\begin{split}\label{eq:fdesrelation}
\Pde{k_1,k_2}{d_1,d_2}
={}&\Gde{k_1,k_2}{d_1,d_2}
+\Gde{k_2,k_1}{d_2,d_1}
+\Gde{k_1+k_2}{d_1+d_2}\\
={}&\sum_{\substack{l_1+l_2=k_1+k_2\\e_1+e_2=d_1+d_2}}
\left(
\binom{l_1-1}{k_1-1}\binom{d_1}{e_1}(-1)^{d_1-e_1}
+\binom{l_1-1}{k_2-1}\binom{d_2}{e_1}(-1)^{d_2-e_1}
\right)
\Gde{l_1,l_2}{e_1,e_2}\\
&+\frac{d_1!d_2!}{(d_1+d_2+1)!}
\binom{k_1+k_2-2}{k_1-1}
\Gde{k_1+k_2-1}{d_1+d_2+1},
\end{split}\end{align}}
where $l_1,l_2\geq1$ and $e_1,e_2\geq0$ in the sum.
\end{definition}

The normalization is the same as for the double-indexed $q$-series: there is
no factorial in the symbols, but we use divided powers in the generating
series. Put
\begin{align*}
\mathfrak G^{\mathrm{DE}}_1\bi{X}{Y}
&=\sum_{\substack{k\geq1\\d\geq0}}
\Gde{k}{d}X^{k-1}\frac{Y^d}{d!},\\
\mathfrak G^{\mathrm{DE}}_2\bi{X_1,X_2}{Y_1,Y_2}
&=\sum_{\substack{k_1,k_2\geq1\\d_1,d_2\geq0}}
\Gde{k_1,k_2}{d_1,d_2}
X_1^{k_1-1}X_2^{k_2-1}
\frac{Y_1^{d_1}}{d_1!}\frac{Y_2^{d_2}}{d_2!}\,,
\end{align*}
and define $\mathfrak P^{\mathrm{DE}}$ in the same way using the symbols
$\Pde{k_1,k_2}{d_1,d_2}$. Then~\eqref{eq:fdesrelation} is equivalent to
{\small
\begin{align}\begin{split}\label{eq:fdesgenrelation}
\mathfrak P^{\mathrm{DE}}\bi{X_1,X_2}{Y_1,Y_2}
={}&\mathfrak G^{\mathrm{DE}}_2\bi{X_1,X_2}{Y_1,Y_2}
+\mathfrak G^{\mathrm{DE}}_2\bi{X_2,X_1}{Y_2,Y_1}
+\frac{\mathfrak G^{\mathrm{DE}}_1\bi{X_1}{Y_1+Y_2}
-\mathfrak G^{\mathrm{DE}}_1\bi{X_2}{Y_1+Y_2}}{X_1-X_2}\\
={}&\mathfrak G^{\mathrm{DE}}_2\bi{X_1+X_2,X_2}{Y_1,Y_2-Y_1}
+\mathfrak G^{\mathrm{DE}}_2\bi{X_1+X_2,X_1}{Y_2,Y_1-Y_2}
+\frac{\mathfrak G^{\mathrm{DE}}_1\bi{X_1+X_2}{Y_1}
-\mathfrak G^{\mathrm{DE}}_1\bi{X_1+X_2}{Y_2}}{Y_1-Y_2}.
\end{split}\end{align}}

The first difference from the formal double zeta space occurs in weight two, where we have
\begin{align}\label{eq:fdesweighttwo}
\Pde{1,1}{0,0}
=2\Gde{1,1}{0,0}+\Gde{2}{0}
=2\Gde{1,1}{0,0}+\Gde{1}{1}\,.
\end{align}
Thus $\Gde{2}{0}=\Gde{1}{1}$ and $\dim_\Q\fdes_2=2$, whereas
$Z_2=0$ in $\dz_2$.

\subsection{The formal double zeta space and the derivative}

There are canonical maps in both directions between $\fdes_K$ and the formal
double zeta space. We use $s_K^{\mathrm{DE}}$ for the map into $\fdes_K$,
since the letter $\sigma$ will later denote the swap on words.

\begin{theorem}[{\cite[Propositions~2.5 and~2.6]{BKuM}}]
\label{thm:fdesdzmaps}
For every $K\geq1$, there is a well-defined linear map
\begin{align*}
\pi_K^{\mathrm{DE}}:\fdes_K\longrightarrow\dz_K
\end{align*}
given on the $G$-symbols by
\begin{align*}
\pi_K^{\mathrm{DE}}\left(\Gde{k}{d}\right)
={}&\delta_{d,0}Z_k+\delta_{k,1}d!Z_{d+1},\\
\pi_K^{\mathrm{DE}}\left(\Gde{k_1,k_2}{d_1,d_2}\right)
={}&\delta_{(d_1,d_2),(0,0)}Z_{k_1,k_2}
+\delta_{(k_1,k_2),(1,1)}d_1!d_2!
\sum_{\substack{a+b=d_1+d_2+2\\a,b\geq1}}
\binom{a-1}{d_2}Z_{a,b}\,.
\end{align*}
In the other direction, there is a well-defined linear map
$s_K^{\mathrm{DE}}:\dz_K\longrightarrow\fdes_K$ given by
\begin{align*}
s_K^{\mathrm{DE}}(Z_K)
={}&\Gde{K}{0}-\delta_{K,2}\Gde{2}{0},\\
s_K^{\mathrm{DE}}(Z_{k_1,k_2})
={}&\Gde{k_1,k_2}{0,0}
+\frac12\left(
\delta_{k_2,1}\Gde{k_1}{1}
-\delta_{k_1,1}\Gde{k_2}{1}
+\delta_{k_1,2}\Gde{k_2+1}{1}
\right),\\
s_K^{\mathrm{DE}}(P_{k_1,k_2})
={}&\Pde{k_1,k_2}{0,0}
+\frac12\left(
\delta_{k_1,2}\Gde{k_2+1}{1}
+\delta_{k_2,2}\Gde{k_1+1}{1}
\right)-\delta_{k_1,1}\delta_{k_2,1}\Gde{2}{0}\,.
\end{align*}
For $K\geq3$, we have
$\pi_K^{\mathrm{DE}}\circ s_K^{\mathrm{DE}}=\operatorname{id}_{\dz_K}$.
In particular, there is a split exact sequence
\begin{equation}\label{eq:fdessplit}
\begin{tikzcd}
0\arrow[r]&\ker(\pi_K^{\mathrm{DE}})\arrow[r]
&\fdes_K\arrow[r,two heads,"{\pi_K^{\mathrm{DE}}}"]
&\dz_K\arrow[r]\arrow[bend left=32,l,"{s_K^{\mathrm{DE}}}"]&0.
\end{tikzcd}
\end{equation}
\end{theorem}

\begin{proof}
Write
\begin{align*}
\mathfrak Z_1(X)&=\sum_{k\geq1}Z_kX^{k-1},\\
\mathfrak Z_2(X_1,X_2)&=\sum_{k_1,k_2\geq1}
Z_{k_1,k_2}X_1^{k_1-1}X_2^{k_2-1}\,,
\end{align*}
and define $\mathfrak P$ similarly. The two formulas defining
$\pi^{\mathrm{DE}}$ are equivalent to
\begin{align*}
\pi^{\mathrm{DE}}\left(\mathfrak G^{\mathrm{DE}}_1\bi{X}{Y}\right)
&=\mathfrak Z_1(X)+\mathfrak Z_1(Y),\\
\pi^{\mathrm{DE}}\left(\mathfrak G^{\mathrm{DE}}_2
\bi{X_1,X_2}{Y_1,Y_2}\right)
&=\mathfrak Z_2(X_1,X_2)+\mathfrak Z_2(Y_1+Y_2,Y_1)\,.
\end{align*}
The image of $\mathfrak P^{\mathrm{DE}}$ is
$\mathfrak P(X_1,X_2)+\mathfrak P(Y_1,Y_2)$. Substituting these three
identities into~\eqref{eq:fdesgenrelation} gives two copies of the formal
double zeta relation~\eqref{eq:dzgenrel}. This proves that
$\pi_K^{\mathrm{DE}}$ is well-defined.

For the other direction, substitute the displayed formulas for
$s_K^{\mathrm{DE}}$ into~\eqref{eq:dzrel}. We give the boundary calculation
explicitly. Put $K=k_1+k_2$ and
\begin{align*}
c_j=\binom{j-1}{k_1-1}+\binom{j-1}{k_2-1},\qquad
c=\binom{K-2}{k_1-1}\,.
\end{align*}
In the stuffle expression
$Z_{k_1,k_2}+Z_{k_2,k_1}+Z_K$, the correction terms in the images of the
two double symbols cancel except for
\begin{align*}
\frac12\left(
\delta_{k_1,2}\Gde{k_2+1}{1}
+\delta_{k_2,2}\Gde{k_1+1}{1}
\right)\,.
\end{align*}
This is the correction in the formula for
$s_K^{\mathrm{DE}}(P_{k_1,k_2})$. When $K=2$, the two additional terms
$-\Gde{2}{0}$ agree as well. For $K\geq3$, the correction in the shuffle
expression is
\begin{align*}
&\sum_{j=1}^{K-1}c_j\left(
s_K^{\mathrm{DE}}(Z_{j,K-j})-\Gde{j,K-j}{0,0}
\right)\notag\\
&\qquad=\left(c+\frac12\delta_{k_1,2}
+\frac12\delta_{k_2,2}\right)\Gde{K-1}{1}\,.
\end{align*}
Here only $j=1,2,K-1$ contribute: the terms containing
$\delta_{k_1,1}$ or $\delta_{k_2,1}$ cancel, while the term with $j=K-1$
gives $c\Gde{K-1}{1}$. The second equality in
\eqref{eq:fdesrelation} with $d_1=d_2=0$ contains precisely the first part
$c\Gde{K-1}{1}$. The case $K=2$ follows separately from
$\Gde{2}{0}=\Gde{1}{1}$. Thus $s_K^{\mathrm{DE}}$ is well-defined. Applying
$\pi_K^{\mathrm{DE}}$ to its three formulas gives the identity on the
generators of $\dz_K$ for $K\geq3$, which proves the last statement.
\end{proof}

The lower indices also give a formal derivative.

\begin{proposition}[{\cite[Proposition~2.7]{BKuM}}]
\label{prop:fdesderivative}
For every $K\geq1$, there is a well-defined linear map
\begin{align*}
D_K^{\mathrm{DE}}:\fdes_K\longrightarrow\fdes_{K+2}
\end{align*}
given by
\begin{align*}
D_K^{\mathrm{DE}}\Gde{k}{d}
={}&k\Gde{k+1}{d+1},\\
D_K^{\mathrm{DE}}\Gde{k_1,k_2}{d_1,d_2}
={}&k_1\Gde{k_1+1,k_2}{d_1+1,d_2}
+k_2\Gde{k_1,k_2+1}{d_1,d_2+1}\,.
\end{align*}
\end{proposition}

\begin{proof}
On the generating series, the two formulas are given by
\begin{align*}
\frac{\partial^2}{\partial X\partial Y},\qquad
\frac{\partial^2}{\partial X_1\partial Y_1}
+\frac{\partial^2}{\partial X_2\partial Y_2}\,,
\end{align*}
respectively. The second operator commutes with the two linear changes of
variables in~\eqref{eq:fdesgenrelation}. When it is applied to either
difference quotient, the terms in which a derivative hits the denominator
cancel. Therefore applying the operator to~\eqref{eq:fdesgenrelation} gives
the same relation with every coefficient replaced by its stated image. Hence
the map is well-defined.
\end{proof}

\subsection{Realizations from the Fay identity}

Let $A$ be a commutative $\Q$-algebra. A realization of $\fdes_K$ in $A$ is
just a linear map $\fdes_K\rightarrow A$. The generating-series relation
allows us to construct such maps from solutions of the Fay identity. We first
need a right action of $\GL_2(\Z)$. For
$\gamma=\left(\begin{smallmatrix}a&b\\c&d\end{smallmatrix}\right)$, put
\begin{align}\label{eq:fdesgroupaction}
(R\mid\gamma)\bi{X_1,X_2}{Y_1,Y_2}
={}&R\bi{aX_1+bX_2,cX_1+dX_2}
{\det(\gamma)(dY_1-cY_2),\det(\gamma)(-bY_1+aY_2)}\,.
\end{align}
We extend this action $\Z$-linearly to $\Z[\GL_2(\Z)]$ and use the matrices
\begin{align*}
\epsilon=\pmatr{0&1\\1&0},\qquad
T=\pmatr{1&1\\0&1},\qquad
S=\pmatr{0&-1\\1&0},\qquad
U=\pmatr{1&-1\\1&0}\,.
\end{align*}
On the $X$-variables, this is the four-variable lift of the action used in
Chapter~\ref{sec:mdandmzv}. The action on the $Y$-variables is the
contragredient one. Identities containing polar terms are interpreted in the
localization in which the linear forms occurring in them are invertible.
For $B_1(X,Y)\in A[[X,Y]]$, define
\begin{align*}
R_{B_1}^{\ast}\bi{X_1,X_2}{Y_1,Y_2}
&=\frac{B_1(X_1,Y_1+Y_2)-B_1(X_2,Y_1+Y_2)}{X_1-X_2},\\
R_{B_1}^{\sh}\bi{X_1,X_2}{Y_1,Y_2}
&=\frac{B_1(X_1+X_2,Y_1)-B_1(X_1+X_2,Y_2)}{Y_1-Y_2}\,.
\end{align*}

\begin{theorem}[{\cite[Theorem~3.3 and Corollary~3.4]{BKuM}}]
\label{thm:fdesfayrealization}
Suppose that $B_1(-X,-Y)=-B_1(X,Y)$ and that
\begin{align*}
F(X,Y)=-\frac12\left(\frac1X+\frac1Y\right)+B_1(X,Y)
\end{align*}
satisfies the Fay identity
\begin{align}\label{eq:fdesfay}
0={}&F(X_1,Y_1)F(X_2,Y_2)
+F(X_1-X_2,-Y_2)F(X_1,Y_1+Y_2)\\
&+F(-X_2,-Y_1-Y_2)F(X_1-X_2,Y_1).\notag
\end{align}
Set
\begin{align*}
P_B\bi{X_1,X_2}{Y_1,Y_2}=B_1(X_1,Y_1)B_1(X_2,Y_2)
\end{align*}
and
{\small
\begin{align}\label{eq:fdesBtwo}
B_2={}&\frac13P_B\mid(1+T^{-1})
-\frac1{12}R_{B_1}^{\ast}\mid(5-3U+U\epsilon)
-\frac1{12}R_{B_1}^{\sh}\mid\bigl(T^{-1}(5-3\epsilon+U)\bigr).
\end{align}}
Then
\begin{align}\label{eq:fdesBrelation}
P_B
=B_2\mid(1+\epsilon)+R_{B_1}^{\ast}
=B_2\mid T(1+\epsilon)+R_{B_1}^{\sh}\,.
\end{align}
Therefore, the coefficients of $B_1,B_2$ and $P_B$ give a realization of
$\fdes_K$ in $A$ for every $K\geq1$, where the lower variables are expanded
in divided powers as in~\eqref{eq:fdesgenrelation}.
\end{theorem}

\begin{proof}
Put
\begin{align*}
P_B^{\mathrm{pol}}
={}&-\frac12\left(\frac1{X_2}+\frac1{Y_2}\right)B_1(X_1,Y_1)
-\frac12\left(\frac1{X_1}+\frac1{Y_1}\right)B_1(X_2,Y_2),\\
\widetilde P_B={}&P_B^{\mathrm{pol}}+P_B\,.
\end{align*}
Thus
\begin{align*}
\widetilde P_B
=F(X_1,Y_1)F(X_2,Y_2)
-\frac14\left(\frac1{X_1}+\frac1{Y_1}\right)
\left(\frac1{X_2}+\frac1{Y_2}\right)\,.
\end{align*}
Oddness gives the equations with $1+S$ and $1-\epsilon$, while the Fay
identity gives the equation with $1+U+U^2$. The pure polar product satisfies
the same three equations by direct calculation. This gives
\begin{align}\label{eq:fdesperiodrelations}
\widetilde P_B\mid(1+U+U^2)
=\widetilde P_B\mid(1+S)
=\widetilde P_B\mid(1-\epsilon)=0\,.
\end{align}

We first record the group-ring calculation which will be used below. Let $R$
satisfy~\eqref{eq:fdesperiodrelations}, and put
\begin{align*}
\delta=\pmatr{-1&0\\0&1},\qquad
\sigma=\pmatr{-1&0\\0&-1},\qquad
T'=U^2S^{-1}\,.
\end{align*}
Since
\begin{align*}
1-T-T'=(1+S)S^{-1}-(1+U+U^2)S^{-1}\,,
\end{align*}
we have $R\mid(1-T-T')=0$. We also have
$R\mid\delta=-R$, $R\mid\sigma=R$, $T'=\epsilon T\epsilon$ and
$\sigma T'T^{-1}=\delta T^{-1}\epsilon$. It follows that
\begin{align}\label{eq:fdesgroupringcalculation}
R\mid(1+T^{-1})(1+\epsilon)
&=2R+R\mid(T^{-1}-\delta T^{-1}\epsilon)\\
&=3R+R\mid(1-T-T')T^{-1}=3R,\notag\\
R\mid(1+T^{-1})T(1+\epsilon)
&=2R+R\mid(T+T')=3R.\notag
\end{align}

Finally, set
\begin{align*}
\beta_B={}&\frac14R_{B_1}^{\ast}\mid(5-3U+U\epsilon)
+\frac14R_{B_1}^{\sh}\mid
\bigl(T^{-1}(5-3\epsilon+U)\bigr)\,.
\end{align*}
A direct substitution in the action~\eqref{eq:fdesgroupaction}, using the
oddness of $B_1$, gives
\begin{align}\label{eq:fdesbetacalculation}
\beta_B\mid(1+\epsilon)
&=3R_{B_1}^{\ast}
+P_B^{\mathrm{pol}}\mid(1-T^{-1}-T^{-1}\epsilon),\\
\beta_B\mid T(1+\epsilon)
&=3R_{B_1}^{\sh}
+P_B^{\mathrm{pol}}\mid(1-T-T\epsilon).\notag
\end{align}
Since $B_2=\frac13P_B\mid(1+T^{-1})-\frac13\beta_B$, the first identities
in~\eqref{eq:fdesgroupringcalculation} and
\eqref{eq:fdesbetacalculation}, together with
$P_B^{\mathrm{pol}}\mid\epsilon=P_B^{\mathrm{pol}}$, give
\begin{align*}
B_2\mid(1+\epsilon)
={}&\frac13(\widetilde P_B-P_B^{\mathrm{pol}})
\mid(1+T^{-1})(1+\epsilon)-\frac13\beta_B\mid(1+\epsilon)\\
={}&\widetilde P_B-P_B^{\mathrm{pol}}-R_{B_1}^{\ast}
=P_B-R_{B_1}^{\ast}\,.
\end{align*}
The second identities give in the same way
$B_2\mid T(1+\epsilon)=P_B-R_{B_1}^{\sh}$. This proves
\eqref{eq:fdesBrelation}. Comparing coefficients with
\eqref{eq:fdesgenrelation} proves the last statement.
\end{proof}

The main example is obtained from the Kronecker function. With the Eisenstein
series $G_k$ of~\eqref{eq:defgk}, put
\begin{align}\label{eq:fdeskronecker}
\mathfrak K_q(X,Y)
={}&-\frac12\left(\frac1X+\frac1Y\right)
+\sum_{\substack{r,s\geq0\\r+s\text{ odd}}}
\frac{|r-s|!}{r!}\left(\qdq\right)^{\min(r,s)}
G_{|r-s|+1}X^r\frac{Y^s}{s!}\,.
\end{align}

\begin{proposition}[{\cite[Lemma~3.5 and Proposition~3.7]{BKuM}}]
\label{prop:fdeskronecker}
The series $\mathfrak K_q$ satisfies the Fay identity. It therefore gives a
surjective Kronecker realization
\begin{align*}
\rho_K^{\mathfrak K}:\fdes_K\twoheadrightarrow\qmf_K\,.
\end{align*}
On the depth-one symbols it is given explicitly by
\begin{align}\label{eq:fdeskroneckerdepthone}
\rho_K^{\mathfrak K}\left(\Gde{k}{d}\right)
=\begin{cases}
\displaystyle
\frac{|k-d-1|!}{(k-1)!}
\left(\qdq\right)^{\min(k-1,d)}G_{|k-d-1|+1},
&k+d\text{ even},\\[0.3cm]
0,&k+d\text{ odd}\,.
\end{cases}
\end{align}
We also have
\begin{align*}
\rho_K^{\mathfrak K}\left(\Pde{k_1,k_2}{d_1,d_2}\right)
=\rho_{k_1+d_1}^{\mathfrak K}\left(\Gde{k_1}{d_1}\right)
\rho_{k_2+d_2}^{\mathfrak K}\left(\Gde{k_2}{d_2}\right)\,,
\end{align*}
and the images of the depth-two $G$-symbols are the coefficients of the series
$B_2$ in~\eqref{eq:fdesBtwo} for the non-polar part of $\mathfrak K_q$.
\end{proposition}

\begin{proof}
The Fay identity for $\mathfrak K_q$ is a form of Riemann's theta identity (see \cite[Lemma~3.5]{BKuM} and \cite[Section~5, equation~(21)]{Za1}). Theorem
\ref{thm:fdesfayrealization} now gives the realization and the product formula.
Formula~\eqref{eq:fdeskroneckerdepthone} follows by putting $r=k-1$ and $s=d$
in~\eqref{eq:fdeskronecker}. Finally, the products of derivatives of
Eisenstein series occurring in these coefficients span $\qmf_K$, which gives
the surjectivity (see \cite[Proposition~3.7]{BKuM}).
\end{proof}

The formal derivative becomes the usual $q$-derivative under this
realization.

\begin{proposition}[{\cite[Proposition~3.8]{BKuM}}]
\label{prop:fdeskroneckerderivative}
For every $K\geq1$, the following diagram commutes:
\begin{center}
\begin{tikzcd}[column sep=large,row sep=large]
\fdes_K\arrow[r,"{\rho_K^{\mathfrak K}}"]
\arrow[d,"{D_K^{\mathrm{DE}}}"']
&\qmf_K\arrow[d,"{\qdq}"]\\
\fdes_{K+2}\arrow[r,"{\rho_{K+2}^{\mathfrak K}}"']
&\qmf_{K+2}.
\end{tikzcd}
\end{center}
\end{proposition}

\begin{proof}
If $B_1$ is the non-polar part of $\mathfrak K_q$, then
\begin{align*}
\qdq B_1(X,Y)=\frac{\partial^2}{\partial X\partial Y}B_1(X,Y)\,.
\end{align*}
The four-variable differential operator in the proof of
Proposition~\ref{prop:fdesderivative} commutes with the action
in~\eqref{eq:fdesgroupaction} and with the two difference quotients. Applying
it to~\eqref{eq:fdesBtwo} therefore gives $\qdq B_2$. Comparing coefficients
proves the commutativity.
\end{proof}

For example, the weight-two relation gives
\begin{align*}
\rho_2^{\mathfrak K}\left(\Gde{2}{0}\right)
=\rho_2^{\mathfrak K}\left(\Gde{1}{1}\right)=G_2,\qquad
\rho_2^{\mathfrak K}\left(\Gde{1,1}{0,0}\right)=-\frac12G_2\,.
\end{align*}
and
\begin{align*}
\rho_4^{\mathfrak K}\left(D_2^{\mathrm{DE}}\Gde{2}{0}\right)
=2\rho_4^{\mathfrak K}\left(\Gde{3}{1}\right)=\qdq G_2\,.
\end{align*}

Taking the constant term of the Kronecker realization recovers the Bernoulli
realization of Theorem~\ref{thm:betarealization}. More precisely, for
$K\geq3$, comparison of the coefficients gives
\begin{align}\label{eq:fdesconstanttermidentity}
\operatorname{ct}\circ\rho_K^{\mathfrak K}
=\varphi_\beta\circ\pi_K^{\mathrm{DE}}
\end{align}
and hence the commutative diagram
\begin{center}
\begin{tikzcd}[column sep=large,row sep=large]
\fdes_K\arrow[r,"{\rho_K^{\mathfrak K}}"]
\arrow[d,"{\pi_K^{\mathrm{DE}}}"']
&\qmf_K\arrow[d,"{q\mapsto0}"]\\
\dz_K\arrow[r,"{\varphi_\beta}"']&\Q.
\end{tikzcd}
\end{center}
The same diagram is trivially commutative for $K=1$, where both compositions
vanish.
The weight-two case is exceptional:
\begin{align*}
\pi_2^{\mathrm{DE}}\left(\Gde{2}{0}\right)=Z_2=0,\qquad
\left.\rho_2^{\mathfrak K}\left(\Gde{2}{0}\right)\right|_{q=0}
=\beta(2)=-\frac1{24}\,.
\end{align*}

The formal double Eisenstein space also satisfies lower-index extensions of
the sum formula and the parity theorem (see \cite[Proposition~4.1 and Theorem~4.3]{BKuM}). These are stronger than the
formal double zeta statements of Chapter~\ref{sec:mdandmzv}, but they will not
be needed below. The even-weight consequences which lead to the Ramanujan
equations will reappear in the formal multiple Eisenstein algebra.

\section{Formal multiple Eisenstein series}\label{sec:fmes}

We now pass from depth two to arbitrary depth. Recall that in the previous
chapter we introduced the double-indexed $q$-series
$\mb{k_1,\dots,k_r}{d_1,\dots,d_r}$. They satisfy the partition relation and a
quasi-shuffle product, but this product contains terms of lower weight. We now
want to define formal multiple Eisenstein series. For this, we use the top-weight
part of the stuffle product and impose the partition relation, which is called
the \emph{swap} in this context. The definition and the general formal setup are
from \cite{BIM}. This is also the formal algebra used in the derivative proof
in Section~\ref{sec:mesderiv}.

Let us explain where this definition comes from, since this was also the original
motivation for it. Multiple Eisenstein series with all entries at least $2$ satisfy the same relations as the corresponding $\g$ modulo lower weight in $\gstwo$, which was expected in \cite{Ba3} and is now proved in \cite[Main Theorem~B(ii)]{BKa} (see \cref{thm:mesgrqmzv}). By Theorem~\ref{thm:hmwqspaces} the
double-indexed $\g$ span the same space as the single-indexed ones, and by
Conjecture~\ref{conj:allrelg} all relations among them should come from the swap
and the product. Putting these together, we expect that all relations among
multiple Eisenstein series come from the swap and the product modulo lower
weight, i.e. from the usual stuffle product. The algebra $\fmes$ is built exactly
from these two relations.

Notice that the formal algebra is different from the space $\rmes$ of
regularized multiple Eisenstein series introduced in the previous chapter. We
use $\fmes$ for the formal algebra and $\Gf$ for its elements.

\subsection{The algebra of formal multiple Eisenstein series}

Recall the alphabet
\begin{align*}
\azbi=\left\{\ai{k}{d}\ \middle|\ k\geq 1,\ d\geq 0\right\}
\end{align*}
from Section~\ref{subsec:doubleindexedqana}. For words we use the abbreviation
\begin{align*}
\ai{k_1,\dots,k_r}{d_1,\dots,d_r}
=\ai{k_1}{d_1}\cdots\ai{k_r}{d_r}\,.
\end{align*}
The product needed here is the top-weight part of the product $\gqsh$ used for the
double-indexed $q$-series.

\begin{definition}\label{def:fmesstuffle}
The \emph{bi-stuffle product} $\bist$ on $\Q\langle\azbi\rangle$ is the
$\Q$-bilinear product determined by $1\bist w=w\bist1=w$ and
\begin{align*}
\ai{k_1}{d_1}w\bist\ai{k_2}{d_2}v
={}&\ai{k_1}{d_1}\bigl(w\bist\ai{k_2}{d_2}v\bigr)
+\ai{k_2}{d_2}\bigl(\ai{k_1}{d_1}w\bist v\bigr)
+\ai{k_1+k_2}{d_1+d_2}(w\bist v)
\end{align*}
for letters $\ai{k_1}{d_1},\ai{k_2}{d_2}$ and words $w,v$.
\end{definition}

Thus $\bist$ is the quasi-shuffle product associated with the componentwise
product
\begin{align*}
\ai{k_1}{d_1}\diamond\ai{k_2}{d_2}
=\ai{k_1+k_2}{d_1+d_2}\,.
\end{align*}
In particular, $(\Q\langle\azbi\rangle,\bist)$ is a commutative algebra. Notice
that $\bist$ is different from $\gqsh$: the latter contains the additional
lower-weight terms in~\eqref{eq:defgbidiamond}.

For $r\geq1$, consider the generating series of words
\begin{align*}
\mathfrak a_r\bi{X_1,\dots,X_r}{Y_1,\dots,Y_r}
=\sum_{\substack{k_1,\dots,k_r\geq1\\d_1,\dots,d_r\geq0}}
\ai{k_1,\dots,k_r}{d_1,\dots,d_r}
\prod_{j=1}^r X_j^{k_j-1}\frac{Y_j^{d_j}}{d_j!}\,.
\end{align*}

\begin{definition}\label{def:fmesswap}
The \emph{swap} $\sigma$ is the linear involution on
$\Q\langle\azbi\rangle$ defined by comparing coefficients in
\begin{align}\label{eq:fmesswap}
\sigma\mathfrak a_r\bi{X_1,\dots,X_r}{Y_1,\dots,Y_r}
={}&\mathfrak a_r
\bi{Y_1+\dots+Y_r,\dots,Y_1+Y_2,Y_1}
{X_r,X_{r-1}-X_r,\dots,X_1-X_2}\,.
\end{align}
\end{definition}

This is exactly the map $P$ used for the partition relation in
Section~\ref{subsec:doubleindexedqana}. The notation $\sigma$ emphasizes that we
now regard it as the swap of two sets of variables.

\begin{definition}\label{def:fmes}
Let $\mathfrak I_\sigma$ be the ideal of
$(\Q\langle\azbi\rangle,\bist)$ generated by $\sigma(w)-w$ for all words $w$.
The algebra of \emph{formal multiple Eisenstein series} is
\begin{align*}
\fmes=\quotient{(\Q\langle\azbi\rangle,\bist)}{\mathfrak I_\sigma}\,.
\end{align*}
We write
\begin{align*}
\Gf\bi{k_1,\dots,k_r}{d_1,\dots,d_r}
\end{align*}
for the class of $\ai{k_1,\dots,k_r}{d_1,\dots,d_r}$, and set
\begin{align*}
\Gf(k_1,\dots,k_r)
=\Gf\bi{k_1,\dots,k_r}{0,\dots,0}\,.
\end{align*}
By $\Q$-linearity, $\Gf(w)$ then denotes the image of any
$w\in\Q\langle\azbi\rangle$ under the quotient map.
For a word $w=\ai{k_1,\dots,k_r}{d_1,\dots,d_r}$, define
\begin{align*}
\wt(w)&=k_1+\dots+k_r+d_1+\dots+d_r,\\
\lwt(w)&=d_1+\dots+d_r,\qquad \dep(w)=r\,.
\end{align*}
These give a weight grading and lower-weight and depth filtrations on $\fmes$.
We denote the lower-weight-zero part by
\begin{align*}
\fmesz=\big\langle\Gf(k_1,\dots,k_r)\mid r\geq0,\ k_1,\dots,k_r\geq1\big\rangle_\Q\,.
\end{align*}
\end{definition}

Already in depth one, the swap gives
\begin{align}\label{eq:fmesswapdepthone}
\Gf\bi{k}{d}=\frac{d!}{(k-1)!}\Gf\bi{d+1}{k-1}\,.
\end{align}
For example,
\begin{align*}
\Gf\bi{1}{1}=\Gf(2),\qquad
\Gf\bi{3}{1}=\frac{1}{2}\Gf\bi{2}{2}\,.
\end{align*}
This is the formal counterpart of Corollary~\ref{cor:partitionlowdepth}.

The double indices are essential in Definition~\ref{def:fmes}. Nevertheless, it
is expected that they are not needed to span the resulting algebra.

\begin{conjecture}\label{conj:fmeslowerweightzero}
The inclusion $\fmesz\subset\fmes$ is an equality. More precisely,
\begin{align*}
\fil{\lwt,\dep}{d,r}\fmes
\subset \fil{\lwt,\dep}{0,d+r}\fmes
\qquad(d,r\geq0)\,.
\end{align*}
\end{conjecture}

This is the formal version of Theorem~\ref{thm:hmwqspaces}. The expected Hilbert
series is the same one as in Conjecture~\ref{conj:rmesdim}:
\begin{align*}
\sum_{k\geq0}\dim_\Q\fmes_kX^k
\overset{?}{=}
\frac{1}{1-X-X^2-X^3+X^6+X^7+X^8+X^9}\,.
\end{align*}

\subsection{Swap and stuffle relations}

The swap gives a second product on words by
\begin{align*}
w\bist^\sigma v=\sigma\bigl(\sigma(w)\bist\sigma(v)\bigr)\,.
\end{align*}
Since every class in $\fmes$ is swap invariant, the products $\bist$ and
$\bist^\sigma$ have the same image in $\fmes$. This gives
\begin{align}\label{eq:fmesdoublestuffle}
\Gf(w\bist v-w\bist^\sigma v)=0
\qquad(w,v\in\Q\langle\azbi\rangle)\,.
\end{align}
These relations are the formal analogue of the double shuffle relations.

\begin{proposition}[{\cite[Proposition~2.9]{BIM}}]\label{prop:fmesdshdepthone}
For $k_1,k_2\geq1$ and $d_1,d_2\geq0$, we have
{\small
\begin{align*}
\Gf\bi{k_1}{d_1}\Gf\bi{k_2}{d_2}
={}&\Gf\bi{k_1,k_2}{d_1,d_2}
+\Gf\bi{k_2,k_1}{d_2,d_1}
+\Gf\bi{k_1+k_2}{d_1+d_2}\\
={}&\sum_{\substack{l_1+l_2=k_1+k_2\\e_1+e_2=d_1+d_2}}
\left(
\binom{l_1-1}{k_1-1}\binom{d_1}{e_1}(-1)^{d_1-e_1}
+\binom{l_1-1}{k_2-1}\binom{d_2}{e_1}(-1)^{d_2-e_1}
\right)
\Gf\bi{l_1,l_2}{e_1,e_2}\\
&+\frac{d_1!d_2!}{(d_1+d_2+1)!}
\binom{k_1+k_2-2}{k_1-1}
\Gf\bi{k_1+k_2-1}{d_1+d_2+1},
\end{align*}}
where the second sum is over $l_1,l_2\geq1$ and $e_1,e_2\geq0$.
\end{proposition}

\begin{proof}
The first equality is the definition of $\bist$. For the second equality, apply
the swap to both depth-one factors, take their bi-stuffle product and apply the
swap again. Comparing coefficients in~\eqref{eq:fmesswap} gives the stated
formula.
\end{proof}

Comparing this proposition with Definition~\ref{def:fdes} gives a canonical
link between the depth-two space above and the arbitrary-depth algebra.

\begin{proposition}\label{prop:fdesfmesmap}
For every $K\geq1$, there is a canonical linear map
\begin{align*}
\iota_K^{\mathrm{DE}}:\fdes_K\longrightarrow\fmes_K
\end{align*}
given by
\begin{align*}
\Gde{k}{d}&\longmapsto\Gf\bi{k}{d},\\
\Gde{k_1,k_2}{d_1,d_2}&\longmapsto
\Gf\bi{k_1,k_2}{d_1,d_2},\\
\Pde{k_1,k_2}{d_1,d_2}&\longmapsto
\Gf\bi{k_1}{d_1}\Gf\bi{k_2}{d_2}\,.
\end{align*}
\end{proposition}

\begin{proof}
The two equal expressions in~\eqref{eq:fdesrelation} are sent precisely to the
two expressions in Proposition~\ref{prop:fmesdshdepthone}. Thus the map is
well-defined.
\end{proof}

No injectivity statement is made here. The formal multiple Eisenstein algebra
also imposes the individual swap relations in every depth, whereas
Definition~\ref{def:fdes} only imposes their combined depth-two product
relation.

The difference between the formal and the actual product is already visible in
a small example. In $\fmes$ we have
\begin{align}\label{eq:fmesplainproductexample}
\Gf\bi{2}{1}\Gf\bi{3}{2}
={}&\Gf\bi{2,3}{1,2}+\Gf\bi{3,2}{2,1}+\Gf\bi{5}{3}\,,
\end{align}
whereas the $q$-series satisfy
\begin{align}\label{eq:fmesactualproductexample}
\mb{2}{1}\mb{3}{2}
={}&\mb{2,3}{1,2}+\mb{3,2}{2,1}+\mb{5}{3}
-\frac{1}{12}\mb{3}{3}\,.
\end{align}
The last term has weight $6$, while all other terms have weight $8$. Thus
$\bist$ is the top-weight part of $\gqsh$, but the map
$\Gf\binom{\mathbf{k}}{\mathbf{d}}\mapsto
\mb{\mathbf{k}}{\mathbf{d}}$ is not an algebra homomorphism.

The depth-one relations imply the following useful formulas.

\begin{corollary}\label{cor:fmesevenrelations}
For even $k\geq4$, we have
\begin{align}\label{eq:fmesevenrelationone}
\frac{k+1}{2}\Gf(k)
={}&\Gf\bi{k-1}{1}
+\sum_{\substack{k_1+k_2=k\\k_1,k_2\geq2\text{ even}}}
\Gf(k_1)\Gf(k_2)\,,
\end{align}
and, for even $k\geq6$,
\begin{align}\label{eq:fmesevenrelationtwo}
\frac{(k+1)(k-1)(k-6)}{12}\Gf(k)
={}&\sum_{\substack{k_1+k_2=k\\k_1,k_2\geq4\text{ even}}}
(k_1-1)(k_2-1)\Gf(k_1)\Gf(k_2)\,.
\end{align}
\end{corollary}

\begin{proof}
These are obtained by taking the same linear combinations of the relations in
Proposition~\ref{prop:fmesdshdepthone} as for the formal double Eisenstein
space (see \cite[Theorem~2.10 and Corollary~2.11]{BIM}).
\end{proof}

For example, the first relation in lower-weight zero is
\begin{align}\label{eq:fmesfirstrelation}
\Gf(4)=2\Gf(2,2)-2\Gf(3,1)\,.
\end{align}
Applying the realizations below to \eqref{eq:fmesfirstrelation} gives the
corresponding relation for several constructions of double Eisenstein series.
Notice that their individual depth-two values need not agree.

\subsection{The balanced alphabet}

There is another presentation in which the swap becomes a simple reversal. Let
$\mathcal B=\{b_0,b_1,b_2,\dots\}$ and let
$\Q\langle\mathcal B\rangle^0$ be the span of words which do not start with
$b_0$. Define a quasi-shuffle product $\ast_b$ by
\begin{align*}
b_i\diamond_b b_j=
\begin{cases}
b_{i+j},&i,j>0,\\
0,&ij=0\,,
\end{cases}
\end{align*}
and the usual quasi-shuffle recursion. On words written in blocks, define
\begin{align*}
\tau_b\bigl(b_{k_1}b_0^{m_1}\cdots b_{k_s}b_0^{m_s}\bigr)
=b_{m_s+1}b_0^{k_s-1}\cdots b_{m_1+1}b_0^{k_1-1}\,.
\end{align*}
Then $\tau_b$ is an involution. Let $\mathfrak T_b$ be the ideal generated by
the differences $\tau_b(w)-w$, i.e.,
\begin{align*}
\mathfrak T_b
=\bigl(\tau_b(w)-w\mid w\in\Q\langle\mathcal B\rangle^0\bigr)
\end{align*}
in $(\Q\langle\mathcal B\rangle^0,\ast_b)$. The balanced presentation of
\cite{Bu2,Bu3,BIM} gives an isomorphism
\begin{align}\label{eq:fmesbalanced}
\fmes\cong
\quotient{(\Q\langle\mathcal B\rangle^0,\ast_b)}
{\mathfrak T_b}\,.
\end{align}
Notice that in this presentation the swap reverses the order of the blocks and
sends each block $b_kb_0^m$ to $b_{m+1}b_0^{k-1}$. This presentation is used
for the balanced multiple $q$-zeta values in \cite{Bu2,BIM}.

\subsection{The formal constant term}

The constant term of a multiple Eisenstein series is a multiple zeta value. A
formal version of this operation is obtained by taking a quotient of $\fmes$.
Set
\begin{align*}
\mathcal A_0&=\left\{\ai{k}{0}\mid k\geq1\right\},&
\mathcal A^1&=\left\{\ai{1}{d}\mid d\geq0\right\}\,.
\end{align*}
Let $\mathfrak N$ be the ideal in $\fmes$ generated by the classes of all words
which are not contained in $(\mathcal A^1)^*(\mathcal A_0)^*$. Thus the words
which can survive modulo $\mathfrak N$ are of the form
\begin{align*}
\ai{1,\dots,1,k_1,\dots,k_r}{d_1,\dots,d_s,0,\dots,0}\,.
\end{align*}

\begin{theorem}[{\cite[Definition~2.13, Proposition~2.14 and
Theorem~2.17]{BIM}}]\label{thm:fmesfmzquotient}
There is a canonical isomorphism
\begin{align}\label{eq:fmesfmzquotient}
\quotient{\fmes}{\mathfrak N}&\longrightarrow\fmz,\\
\Gf(k_1,\dots,k_r)&\longmapsto\fzeta(k_1,\dots,k_r)\,.
\end{align}
In particular, the canonical formal constant-term map
\begin{align*}
\pi_f:\fmes\longrightarrow\fmz
\end{align*}
is surjective already on $\fmesz$.
\end{theorem}

The proof uses generating series of bimoulds and is given in the appendix of
\cite{BIM}. Notice that after taking the quotient by $\mathfrak N$, the stuffle
relations remain unchanged and the swap gives the regularized shuffle
relations. Therefore the quotient is isomorphic to $\fmz$.
This also explains why we introduced $\fmz$ first: it is the formal target of
the constant-term map.

\subsection{Derivations of quasi-shuffle algebras}

Before defining the derivations on $\fmes$, we recall a general construction for
quasi-shuffle algebras. Let $A$ be an alphabet
and let $\diamond$ be a commutative and associative product on $\Q A$. A
derivation $\varphi$ of $(\Q A,\diamond)$ acts on words by
\begin{align}\label{eq:fmesletterderivation}
\Theta^\varphi(a_1\cdots a_r)
=\sum_{j=1}^r a_1\cdots a_{j-1}\varphi(a_j)a_{j+1}\cdots a_r\,.
\end{align}

\begin{proposition}\label{prop:fmesletterderivation}
The map $\Theta^\varphi$ is a derivation of the quasi-shuffle algebra
$(\Q\langle A\rangle,\ast_\diamond)$. We also have
\begin{align*}
[\Theta^\varphi,\Theta^\psi]=\Theta^{[\varphi,\psi]}\,.
\end{align*}
\end{proposition}

\begin{proof}
The first statement follows by induction on the sum of the lengths of two words.
In the merged term of the quasi-shuffle recursion, one uses
$\varphi(a\diamond b)=\varphi(a)\diamond b+a\diamond\varphi(b)$. The second
statement follows directly from~\eqref{eq:fmesletterderivation}.
\end{proof}

For the lowering operator we will also need derivations which remove a specified
letter or act on two consecutive letters. They can be obtained from the
corresponding shuffle derivations by Hoffman's exponential isomorphism (see \cite[Section~3]{BIM}). Proposition~\ref{prop:fmesletterderivation} already gives
the leading part of the lowering operator.

A derivation of $(\Q\langle\azbi\rangle,\bist)$ descends to $\fmes$ if it
preserves the swap ideal $\mathfrak I_\sigma$. In particular, every derivation
which commutes with $\sigma$ descends.

\subsection{The \texorpdfstring{$\sltwo$}{sl2}-structure}

Recall from Section~\ref{sec:qmfsl2algebras} that an $\sltwo$-algebra is an
algebra equipped with three derivations $\GD,\GW,\Gdelta$ satisfying
\begin{align}\label{eq:fmessl2relations}
[\GW,\GD]=2\GD,\qquad
[\GW,\Gdelta]=-2\Gdelta,\qquad
[\Gdelta,\GD]=\GW\,.
\end{align}
It is quite surprising that the purely formal algebra $\fmes$, which was defined
just by the bi-stuffle product and the swap, carries such a structure at all.
We first define these operators on the free word algebra and then show that
they descend to $\fmes$. Put $\mathfrak a_0=1$ and let
$\mathfrak a=(\mathfrak a_r)_{r\geq0}$ be the word generating series above.
The operators are obtained by applying the following formulas to
$\mathfrak a$ and comparing coefficients. More generally, let
$F=(F_r)_{r\geq0}$ be a family of such power series and write
$\mathbf X=(X_1,\dots,X_r)$ and $\mathbf Y=(Y_1,\dots,Y_r)$. Set
\begin{align*}
p_j(\mathbf X)=(X_1,\dots,X_{j-1},X_{j+1},\dots,X_r)\,.
\end{align*}
For $j<r$, let $c_{j,j+1}(\mathbf Y)$ be obtained by replacing
$Y_j,Y_{j+1}$ by $Y_j+Y_{j+1}$, and set
$c_{r,r+1}(\mathbf Y)=(Y_1,\dots,Y_{r-1})$. For $2\leq j\leq r$, let
$c_{j-1,j}(\mathbf Y)$ be obtained by replacing $Y_{j-1},Y_j$ by
$Y_{j-1}+Y_j$. Define
\begin{align*}
(\varphi_j^+F)_r(\mathbf X,\mathbf Y)
&=\begin{cases}
F_{r-1}\bigl(p_j(\mathbf X),c_{j,j+1}(\mathbf Y)\bigr),&1\leq j\leq r,\\
0,&\text{otherwise,}
\end{cases}\\
(\varphi_j^-F)_r(\mathbf X,\mathbf Y)
&=\begin{cases}
F_{r-1}\bigl(p_j(\mathbf X),c_{j-1,j}(\mathbf Y)\bigr),&2\leq j\leq r,\\
0,&\text{otherwise.}
\end{cases}
\end{align*}

The raising and weight operators are
\begin{align}\label{eq:fmesDWseries}
(\GD F)_r&=\sum_{j=1}^r\frac{\partial^2F_r}{\partial X_j\partial Y_j},\\
(\GW F)_r&=\sum_{j=1}^r
\left(X_j\frac{\partial}{\partial X_j}
+Y_j\frac{\partial}{\partial Y_j}+1\right)F_r\,.
\end{align}
The lowering operator is given by
\begin{align}\label{eq:fmesdeltaseries}
(\Gdelta F)_r={}&\left(\sum_{j=1}^rX_jY_j\right)F_r\\
&+\sum_{j\geq1}\left(
-\frac{X_j-X_{j+1}+Y_j}{2}(\varphi_j^+F)_r
-\frac{X_{j-1}-X_j+Y_j}{2}(\varphi_j^-F)_r\right.\\
&\hspace{42mm}\left.
+\frac{1}{4}\bigl((\varphi_j^+)^2F\bigr)_r
-\frac{1}{4}\bigl((\varphi_j^-)^2F\bigr)_r
\right)\,,
\end{align}
where terms outside the ranges above are zero and $X_0=X_{r+1}=0$. Applied to
$\mathfrak a$, these formulas define $\Q$-linear coefficient operators on
$\Q\langle\azbi\rangle$.

\begin{theorem}[{\cite[Propositions~4.9 and~4.11 and Theorem~4.12]{BIM}}]\label{thm:fmessl2}
The coefficient operators $\GD,\GW$ and $\Gdelta$ are $\sigma$-equivariant
derivations of $(\Q\langle\azbi\rangle,\bist)$. Hence they descend to
derivations of $\fmes$, where they satisfy~\eqref{eq:fmessl2relations}. In
particular, $\fmes$ is an $\sltwo$-algebra.
\end{theorem}

\begin{proof}
On letters, $\GD$ and $\GW$ are induced by the derivations
\begin{align*}
\ai{k}{d}&\longmapsto k\ai{k+1}{d+1},&
\ai{k}{d}&\longmapsto(k+d)\ai{k}{d}
\end{align*}
of the componentwise product $\diamond$. Hence
Proposition~\ref{prop:fmesletterderivation} shows that they are derivations for
$\bist$. Their generating-series expressions commute with the change of
variables~\eqref{eq:fmesswap}.

On coefficients, we have
\begin{align*}
\Gdelta=\Gdelta^1-\frac{1}{2}
(\Gdelta^2+\Gdelta^3+\Gdelta^4+\Gdelta^5)\,,
\end{align*}
where each $\Gdelta^i$ is a derivation obtained from the constructions mentioned
above. The generating-series formula for $\Gdelta$ is $\sigma$-equivariant.
We also have $[\Gdelta^1,\GD]=\GW$, while
$[\Gdelta^i,\GD]=0$ for $i=2,3,4,5$. The remaining relations follow from the
weights of the operators. See \cite[Section~4]{BIM} for the coefficient
formulas and the full verification.
\end{proof}

On coefficients, the first two operators have the particularly simple form
\begin{align}\label{eq:fmesDoncoefficients}
\GD\Gf\bi{k_1,\dots,k_r}{d_1,\dots,d_r}
={}&\sum_{j=1}^rk_j
\Gf\bi{k_1,\dots,k_j+1,\dots,k_r}
{d_1,\dots,d_j+1,\dots,d_r},\\
\GW\Gf\bi{k_1,\dots,k_r}{d_1,\dots,d_r}
={}&\bigl(k_1+\dots+k_r+d_1+\dots+d_r\bigr)
\Gf\bi{k_1,\dots,k_r}{d_1,\dots,d_r}\,.
\end{align}
In particular, the map from Proposition~\ref{prop:fdesfmesmap} is compatible
with the derivatives:
\begin{align}\label{eq:fdesfmesderivative}
\iota_{K+2}^{\mathrm{DE}}\circ D_K^{\mathrm{DE}}
=\GD\circ\iota_K^{\mathrm{DE}}\,.
\end{align}
This follows directly on the single and double $G$-symbols from
Proposition~\ref{prop:fdesderivative}. These symbols span after eliminating
the $P$-symbols.
In particular,
\begin{align}\label{eq:fmesDdepthonebasic}
\GD\Gf(k)=k\Gf\bi{k+1}{1}\,.
\end{align}

The operator $\Gdelta$ also has a simple expression on the lower-weight-zero
part. With the convention $\Gf(\emptyset)=1$, and with terms involving
non-existent indices understood to be zero, we have
\begin{align}\label{eq:fmesdeltafilzero}
\Gdelta\Gf(k_1,\dots,k_r)
={}&-\frac{1}{2}\delta_{k_1,2}\Gf(k_2,\dots,k_r)
+\frac{1}{4}\delta_{k_1,1}\delta_{k_2,1}\Gf(k_3,\dots,k_r)\\
&+\frac{1}{2}
\sum_{\substack{1\leq j<r\\k_j=1,\ k_{j+1}>1}}
\Gf(k_1,\dots,k_{j-1},k_{j+1}-1,k_{j+2},\dots,k_r)\\
&-\frac{1}{2}
\sum_{\substack{2\leq j\leq r\\k_j=1,\ k_{j-1}>1}}
\Gf(k_1,\dots,k_{j-2},k_{j-1}-1,k_{j+1},\dots,k_r)\,.
\end{align}
Thus $\fmesz$ is stable under $\Gdelta$ and $\GW$. It is also stable under
$\GD$, since for $w=z_{k_1}\cdots z_{k_r}\in\HH^1$ we have
\begin{align}\label{eq:fmesDonfilzero}
\GD\Gf(k_1,\dots,k_r)
=\Gf\bigl(z_2\ast w-z_2\sh w\bigr)\,,
\end{align}
where the right-hand side is interpreted using the map from $\HH^1$ to
$\fmesz$. We have therefore proved the following.

\begin{proposition}\label{prop:fmeszerosl2}
The lower-weight-zero part $\fmesz$ is an $\sltwo$-subalgebra of $\fmes$.
\end{proposition}

For example, formulas~\eqref{eq:fmesDdepthonebasic} and
\eqref{eq:fmesDonfilzero} give
\begin{align}\label{eq:fmesDexamples}
\GD\Gf(1)&=\Gf(3)-\Gf(2,1),\\
\GD\Gf(3)&=\Gf(5)-2\Gf(3,2)-6\Gf(4,1)
=5\Gf(5)-4\Gf(3,2)-6\Gf(2,3)\,.
\end{align}
More generally, for $k\geq2$,
\begin{align}\label{eq:fmesDgeneraldepthone}
\GD\Gf(k)
={}&(2k-1)\Gf(k+2)-\Gf(2,k)
-\sum_{j=2}^k(k+j-1)\Gf(k+2-j,j)\,.
\end{align}
This expression only involves indices at least $2$.

We now consider the following formal version of the algebra $\mes$:
\begin{align}\label{eq:deffames}
\fames=\big\langle\Gf(k_1,\dots,k_r)
\mid r\geq0,\ k_1,\dots,k_r\geq2\big\rangle_\Q\,.
\end{align}
It is stable under $\GW$ and, by~\eqref{eq:fmesdeltafilzero}, under
$\Gdelta$. The formal identity~\eqref{eq:fmesdroponeidentity}, which uses the
Drop1 operator from Section~\ref{sec:dropone}, shows that it is also stable
under $\GD$ in arbitrary depth. We therefore obtain the following.

\begin{theorem}[{\cite[Main Theorems~C and D(i)]{BKM}}]\label{thm:famessl2}
The algebra $\fames$ is an $\sltwo$-subalgebra of $\fmes$. Moreover, the
analytic realization of \cite{BKM} restricts to a surjective algebra
homomorphism
\begin{align}\label{eq:famesanalyticrealization}
\fames&\longrightarrow\mes,\\
\Gf(k_1,\dots,k_r)&\longmapsto\aG_{k_1,\dots,k_r}
\end{align}
which intertwines $\GD$ with $(2\pi i)\frac{d}{d\tau}$.
\end{theorem}
\begin{proof}
The stability under $\GW$ and $\Gdelta$ follows from the formulas above, and
the stability under $\GD$ follows from~\eqref{eq:fmesdroponeidentity}. Hence the
$\sltwo$-relations on $\fmes$ restrict to $\fames$. The second statement is the
restriction of the bi-multiple Eisenstein realization in
Theorem~\ref{thm:fmesanalyticrealization}. It is surjective by the definition of
$\mes$.
\end{proof}

\begin{conjecture}\label[conjecture]{conj:famesrealization}
The map~\eqref{eq:famesanalyticrealization} is injective and hence an algebra
isomorphism.
\end{conjecture}

Notice that by \cref{thm:messl2algebra} the operators $\GW$ and $\Gdelta$ are also well-defined on $\mes$, and the map \eqref{eq:famesanalyticrealization} commutes with all three operators, i.e. it is a surjective morphism of $\sltwo$-algebras (see \Cref{sec:messl2}).

Finally, the raising operator satisfies
\begin{align}\label{eq:fmesDkerpi}
\GD\fmes\subset\ker(\pi_f)\,,
\end{align}
because every term in~\eqref{eq:fmesDoncoefficients} contains a letter
$\ai{k}{d}$ with $k\geq2$ and $d\geq1$. Hence $\GD$ induces the zero
derivation on $\fmz$. The lowering operator does not descend to $\fmz$. For
example,
\begin{align*}
\Gdelta\Gf\bi{2}{1}=\Gf(1)\,,
\end{align*}
where $\pi_f\Gf\binom{2}{1}=0$ but
$\pi_f\Gdelta\Gf\binom{2}{1}=\fzeta(1)\neq0$.

\subsection{Formal quasimodular forms}

The $\sltwo$-structure allows us to define formal versions of quasimodular
forms.

\begin{definition}\label{def:fqmf}
The algebra of \emph{formal quasimodular forms} $\fqmf$ is the smallest
$\sltwo$-subalgebra of $\fmes$ containing $\Gf(2)$.
\end{definition}

Combining~\eqref{eq:fmesevenrelationone} with
$\GD\Gf(k)=k\Gf\binom{k+1}{1}$ gives the formal Ramanujan equations.

\begin{corollary}\label{cor:fmesramanujan}
In $\fmes$, we have
\begin{align}\label{eq:fmesramanujan}
\GD\Gf(2)&=5\Gf(4)-2\Gf(2)^2,\\
\GD\Gf(4)&=14\Gf(6)-8\Gf(2)\Gf(4),\\
\GD\Gf(6)&=\frac{120}{7}\Gf(4)^2-12\Gf(2)\Gf(6)\,.
\end{align}
This gives
\begin{align}\label{eq:fmeschazy}
\GD^3\Gf(2)+24\Gf(2)\GD^2\Gf(2)
-36\bigl(\GD\Gf(2)\bigr)^2=0\,.
\end{align}
\end{corollary}

\begin{proof}
The first two identities follow by taking $k=4$ and $k=6$ in
\eqref{eq:fmesevenrelationone}. For the third, first use
\eqref{eq:fmesevenrelationtwo} with $k=8$ to obtain
$\Gf(8)=\frac{6}{7}\Gf(4)^2$, and then take $k=8$ in
\eqref{eq:fmesevenrelationone}. Differentiating the first identity twice and
using all three Ramanujan equations gives~\eqref{eq:fmeschazy}.
\end{proof}

\begin{theorem}[{\cite[Theorem~5.4 and Corollary~5.5~(ii)]{BIM}}]\label{thm:fqmf}
There is an isomorphism of $\sltwo$-algebras
\begin{align*}
\fqmf=\Q[\Gf(2),\Gf(4),\Gf(6)]\cong\qmf\,.
\end{align*}
We also have
\begin{align*}
\fqmf=\Q[\Gf(2),\GD\Gf(2),\GD^2\Gf(2)]\,.
\end{align*}
\end{theorem}

\begin{proof}
The Ramanujan equations show that the algebra generated by
$\Gf(2),\Gf(4),\Gf(6)$ is stable under $\GD$. It is stable under $\GW$ by
weight homogeneity and under $\Gdelta$ by~\eqref{eq:fmesdeltafilzero}. Thus it
is an $\sltwo$-subalgebra generated as an $\sltwo$-algebra by $\Gf(2)$. The
$q$-series realization in
Theorem~\ref{thm:fmesqrealization} sends these three elements to
$G(2),G(4),G(6)$. Since the latter are algebraically independent and generate
the algebra of quasimodular forms, there can be no further relation among the
formal generators.
\end{proof}

The following is the formal version of Euler's formula.

\begin{corollary}\label{cor:fmeseuler}
For every $m\geq1$, there is an element
$Q_{2m}\in\GD\fqmf$ such that
\begin{align}\label{eq:fmeseuler}
\Gf(2m)=-\frac{B_{2m}}{2(2m)!}
\bigl(-24\Gf(2)\bigr)^m+Q_{2m}\,.
\end{align}
\end{corollary}

After applying the formal constant-term map, the term $Q_{2m}$ vanishes by
\eqref{eq:fmesDkerpi}, and~\eqref{eq:fmeseuler} becomes
Euler's relation~\eqref{eq:fmzeuler}. One can prove the corollary directly from
\eqref{eq:fmesevenrelationone}: the generating series of its non-derivative
part satisfies the same differential equation as
\begin{align*}
\frac{1}{2}-\frac{\sqrt{-24\Gf(2)}X}{4}
\coth\left(\frac{\sqrt{-24\Gf(2)}X}{2}\right)\,.
\end{align*}

\begin{definition}\label{def:fmf}
The algebra of \emph{formal modular forms} is
\begin{align*}
\fmf=\ker\bigl(\Gdelta|_{\fqmf}\bigr)\,.
\end{align*}
We denote its weight-$k$ part by $\fmf_k$.
\end{definition}

Since~\eqref{eq:fmesdeltafilzero} gives
\begin{align*}
\Gdelta\Gf(2)=-\frac{1}{2},\qquad
\Gdelta\Gf(4)=\Gdelta\Gf(6)=0\,,
\end{align*}
Theorem~\ref{thm:fqmf} immediately gives
\begin{align}\label{eq:fmfpolynomial}
\fmf=\Q[\Gf(4),\Gf(6)]\cong\mf^\Q\,.
\end{align}

For even $k\geq2$, define the formal normalized Eisenstein series by
\begin{align*}
E_k^f=-\frac{2k!}{B_k}\Gf(k)\,.
\end{align*}
Its image under the $q$-series realization is the normalized Eisenstein series
$E_k$ with constant term $1$.

\begin{definition}\label{def:fcusp}
The algebra of \emph{formal cusp forms} is
\begin{align*}
\fcusp=\ker\bigl(\pi_f|_{\fmf}\bigr)\,.
\end{align*}
\end{definition}

The first non-zero example occurs in weight $12$.

\begin{proposition}\label{prop:fDelta}
The element
\begin{align}\label{eq:fDelta}
\fDelta
&=\frac{(E_4^f)^3-(E_6^f)^2}{1728}
=2400\cdot6!\,\Gf(4)^3-420\cdot7!\,\Gf(6)^2
\end{align}
belongs to $\fcusp$ and satisfies
\begin{align}\label{eq:fDeltaderivative}
\GD\fDelta=E_2^f\fDelta\,.
\end{align}
\end{proposition}

\begin{proof}
Equation~\eqref{eq:fmfpolynomial} shows that $\fDelta$ is a formal modular
form. After applying $\pi_f$, Euler's relation implies
$(E_4^f)^3=(E_6^f)^2$, so $\pi_f(\fDelta)=0$. Finally,
\eqref{eq:fDeltaderivative} follows by differentiating~\eqref{eq:fDelta} and
using the Ramanujan equations.
\end{proof}

Under the isomorphism~\eqref{eq:fmfpolynomial}, $\fDelta$ corresponds to the
classical discriminant
\begin{align*}
\Delta(q)=q\prod_{n\geq1}(1-q^n)^{24}\,.
\end{align*}
Notice that this construction starts just from the bi-stuffle product and the
swap.

\subsection{Realizations}

Let $A$ be a $\Q$-algebra. A \emph{realization} of $\fmes$ in $A$ is an
algebra homomorphism $\fmes\rightarrow A$. Equivalently, it is a
swap-invariant algebra homomorphism
\begin{align*}
(\Q\langle\azbi\rangle,\bist)\longrightarrow A\,.
\end{align*}
We will now give three realizations of $\fmes$.

\begin{theorem}[{\cite[Definition~6.4, Theorem~6.5, Example~6.6(i), and Proposition~6.29]{BaBu}}]\label{thm:fmesqrealization}
There exists an algebra homomorphism
\begin{align*}
\operatorname{ev}_G:\fmes\longrightarrow\Q[[q]]
\end{align*}
such that, for $k\geq2$,
\begin{align}\label{eq:fmesqrealizationdepthone}
\operatorname{ev}_G(\Gf(k))
=G(k)=-\frac{B_k}{2k!}
+\frac{1}{(k-1)!}\sum_{n>0}\sigma_{k-1}(n)q^n\,.
\end{align}
Moreover,
\begin{align*}
\operatorname{ev}_G(\GD f)=\qdq\operatorname{ev}_G(f)
\qquad(f\in\fmes)\,.
\end{align*}
\end{theorem}

The values in higher depth are the combinatorial bi-multiple Eisenstein series
of \cite[Definition~6.4]{BaBu}. The raw series
$\mb{\mathbf{k}}{\mathbf{d}}$ are swap invariant but, as
\eqref{eq:fmesactualproductexample} shows, they are not a realization of
$\fmes$. The construction corrects their lower-weight terms while preserving
the swap and depends on a rational solution of the extended double shuffle
equations. We will construct it in Section~\ref{sec:cmes}. In depth two, it is
closely related to, but not identical with, the realization of
Chapter~\ref{sec:mdandmzv}. The complete comparison, including all boundary
corrections and the exceptional weight-two term, is given in
Theorem~\ref{thm:cmesfdzbbcomparison}.

Composing this realization with $\iota_K^{\mathrm{DE}}$ gives a realization
of the formal double Eisenstein space, but it is not the Kronecker
realization. Equations~\eqref{eq:fmesswapdepthone} and
\eqref{eq:fmesDoncoefficients} show that
\begin{align*}
(\operatorname{ev}_G\circ\iota_K^{\mathrm{DE}})
\left(\Gde{k}{d}\right)
=\rho_K^{\mathfrak K}\left(\Gde{k}{d}\right)
\qquad(k+d=K\text{ even})\,.
\end{align*}
Their depth-two values are different in general, however, and the image of
the Kronecker realization is always contained in $\qmf$, while no such
restriction is imposed on $\operatorname{ev}_G$.

\begin{theorem}[{\cite[Main Theorem~C]{BKM}}]\label{thm:fmesanalyticrealization}
There exists an algebra homomorphism
\begin{align*}
\operatorname{ev}_{\mes}:\fmes\longrightarrow\mathcal O(\Ha)
\end{align*}
such that, for $k_1,\dots,k_r\geq2$,
\begin{align*}
\operatorname{ev}_{\mes}\bigl(\Gf(k_1,\dots,k_r)\bigr)
=\aG_{k_1,\dots,k_r}\,.
\end{align*}
Moreover,
\begin{align*}
\operatorname{ev}_{\mes}(\GD f)
=(2\pi i)\frac{d}{d\tau}\operatorname{ev}_{\mes}(f)
\qquad(f\in\fmes)\,.
\end{align*}
\end{theorem}

The values of this realization on general double indices are the
\emph{bi-multiple Eisenstein series} of \cite{BKM}. Notice that under this
realization the formal derivative gives Theorem~\ref{thm:mesderformula}. Under
the combinatorial realization it gives
Theorem~\ref{thm:cmesbiderivative}.

\begin{theorem}[{\cite[Theorem~4.23]{BI}}]\label{thm:fmesmzvrealization}
There exists an algebra homomorphism
\begin{align*}
\operatorname{ev}_\zeta:\fmes\longrightarrow\mz
\end{align*}
such that
\begin{align*}
\operatorname{ev}_\zeta\bigl(\Gf(k_1,\dots,k_r)\bigr)
=\zeta^\ast(k_1,\dots,k_r;0)\,.
\end{align*}
\end{theorem}

This realization can also be obtained by composing $\pi_f$ with the map which
sends the regularization parameter $\fzeta(1)$ to zero. The following partial
diagram records the $q$-series realization and the two constant-term maps.
\begin{center}
\begin{tikzcd}[column sep=large,row sep=large]
\fmes \arrow[r,"{\operatorname{ev}_G}"] \arrow[d,"{\pi_f}"']
& \Q[[q]]\\
\fmz \arrow[r,two heads] & \mz[T].
\end{tikzcd}
\end{center}

Thus these realizations return from the formal algebra to the real-number,
$q$-series and holomorphic-function worlds of the picture in the Introduction.
The next section constructs the $q$-series realization explicitly and studies
its limits, derivatives and relation with the earlier depth-two construction.

\input{chap_CombinatorialMES}

\input{chap_FormalFiniteMZV}

\vspace{1cm}
\begin{center}
\bf{\Large{\color{qmzvlinecol}\ding{94}}~Exercises~{\color{qmzvlinecol}\ding{94}}}
\end{center}

The following is a collection of exercises intended to deepen the reader's
understanding of this chapter.

\input{chap_FormalFiniteExercises}

\begin{exer}\label{exer:fmzdepthtwo}
Use~\eqref{eq:fmzdepthtwo} with $(k_1,k_2)=(1,2)$ to prove
$\fzeta(3)=\fzeta(2,1)$. Derive the displayed formulas for $\fzeta(4)$ and
$\fzeta(6)$ from the extended double shuffle relations.
\end{exer}

\begin{exer}\label{exer:fdesweighttwo}
Use~\eqref{eq:fdesrelation} in weight two to prove
$\Gde{2}{0}=\Gde{1}{1}$ and $\dim_\Q\fdes_2=2$. Then use the product part of
the Kronecker realization to show that
\begin{align*}
\rho_2^{\mathfrak K}\left(\Gde{1,1}{0,0}\right)=-\frac12G_2\,.
\end{align*}
\end{exer}

\begin{exer}\label{exer:fdeskroneckerderivative}
Extract the coefficient of $X^{k-1}Y^d/d!$ in
\eqref{eq:fdeskronecker} and derive~\eqref{eq:fdeskroneckerdepthone}. Verify
the commutative diagram in Proposition~\ref{prop:fdeskroneckerderivative}
directly on every depth-one symbol.
\end{exer}

\begin{exer}\label{exer:fmesswap}
Use~\eqref{eq:fmesswap} to prove~\eqref{eq:fmesswapdepthone}. Then verify the
corresponding identity for the $q$-series $\mb{k}{d}$ directly from their
definition.
\end{exer}

\begin{exer}\label{exer:fmesproducts}
Prove~\eqref{eq:fmesplainproductexample} from the bi-stuffle product, and
derive the correction term in~\eqref{eq:fmesactualproductexample} from
\eqref{eq:defgbidiamond}.
\end{exer}

\begin{exer}\label{exer:fmeschazy}
Derive the three Ramanujan equations~\eqref{eq:fmesramanujan} from
Corollary~\ref{cor:fmesevenrelations}, and use them to prove the Chazy
equation~\eqref{eq:fmeschazy}.
\end{exer}

\begin{exer}\label{exer:fDelta}
Show directly from the Ramanujan equations that the element $\fDelta$ in
\eqref{eq:fDelta} satisfies $\GD\fDelta=E_2^f\fDelta$.
\end{exer}

\begin{exer}\label{exer:cmesbeta11}
Use the stuffle and shuffle equations in depth two to prove
\eqref{eq:cmesbeta11}. Check directly from~\eqref{eq:cmesgammaZ} that
$\gamma_2=\beta(2)/2=-1/48$ and that the depth-two symmetral equation contains
the term $2\gamma_2$.
\end{exer}

\begin{exer}\label{exer:cmesdepthoneswap}
Use~\eqref{eq:cmesdepthone} to prove the depth-one swap relation
\begin{align*}
\Gb{k}{d}=\frac{d!}{(k-1)!}\Gb{d+1}{k-1}\,.
\end{align*}
Pay attention to the divided powers $Y^d/d!$ in the generating series.
\end{exer}

\begin{exer}\label{exer:cmesG2G3}
Derive both product expressions for $G(2)G(3)$ from symmetrility and swap
invariance, and prove~\eqref{eq:cmesG5relation}.
\end{exer}

\begin{exer}\label{exer:cmescomparison}
Extract coefficients from~\eqref{eq:cmesfdzbbgen} and derive all six cases
in~\eqref{eq:cmesfdzbbcoeff}. Then rewrite the result in the form
\eqref{eq:cmesbkumcorrection} when $k_1+k_2\geq3$.
\end{exer}

\begin{exer}\label{exer:cmesderivatives}
Prove~\eqref{eq:cmesG1derivative} directly from the stuffle and shuffle
products of $z_2$ and $z_1$. Then use~\eqref{eq:cmesderivativekernel} with
$u=v=z_1$ to derive~\eqref{eq:cmesweightfourrelation}.
\end{exer}

\begin{exer}\label{exer:formalfinitedepthtwo}
Use~\eqref{eq:formalfinitedepthtwo} to prove
\eqref{eq:formalfiniteweightfiverelation}. Apply the two maps in
\eqref{eq:formalfiniteclassicalrealizations} and compare the resulting
identities with the examples in Chapter~\ref{sec:overview}.
\end{exer}

\begin{exer}\label{exer:formalfiniteweakparity}
Starting from the stuffle antipode relation, derive
\eqref{eq:formalfiniteweakdepththree} and
\eqref{eq:formalfiniteweakdepthfour}. Explain why the product of two
depth-two values is the only obstruction to strong parity in depth four.
\end{exer}

%% file: chap_CombinatorialMES.tex
\begingroup
\renewcommand{\sectionmark}[1]{\markright{Combinatorial MES}}
\section{Combinatorial multiple Eisenstein series}\label{sec:cmes}
\endgroup

In Theorem~\ref{thm:fmesqrealization}, we stated that the algebra of formal
multiple Eisenstein series has a realization in $\Q[[q]]$. We now want to
construct this realization. Recall that the Fourier expansion of a multiple
Eisenstein series is built out of multiple zeta values and the $q$-series
$\ag$. The idea is to replace the multiple zeta values by a rational solution
$\beta$ of the extended double shuffle relations and $\ag$ by $\g$. This gives
$q$-series with rational coefficients, which we call \emph{combinatorial
multiple Eisenstein series}. The construction is due to the author and
Burmester \cite{BaBu}.

The resulting $q$-series should not be confused with either the combinatorial
double Eisenstein series of Chapter~\ref{sec:mdandmzv} or the multiple
Eisenstein series of Chapter~\ref{sec:mes}. The two depth-two combinatorial
constructions are closely related, but they are not identical. We will compare
them in Section~\ref{sec:cmesdepthtwocomparison}. Unlike the analytic
realization $\operatorname{ev}_{\mes}$ and the Kronecker realization of the
formal double Eisenstein space, the realization constructed here takes values
in $\Q[[q]]$ without being restricted to quasimodular forms. At the end of the
section, we will see that it interpolates between a rational solution of the
extended double shuffle relations at $q=0$ and multiple zeta values as
$q\rightarrow1$. This is one of the nicest features of these $q$-series.

\subsection{Moulds and rational extended double shuffle solutions}
\label{sec:cmesmoulds}

We first recall the mould notation needed for the construction. Mould calculus
was introduced by \'Ecalle \cite{Ec}. For an introduction to its application to
multiple zeta values and double shuffle relations, see \cite{Sc}. Let $A$ be a
commutative $\Q$-algebra. A \emph{mould} $Z$ with values in $A$ is a family of power
series
\begin{align*}
Z(X_1,\dots,X_r)\in A[[X_1,\dots,X_r]]\qquad(r\geq0)\,,
\end{align*}
where we always assume that the depth-zero part is $Z(\emptyset)=1$. We write
its coefficients as
\begin{align}\label{eq:cmesmouldcoeff}
Z(X_1,\dots,X_r)
=\sum_{k_1,\dots,k_r\geq1}z(k_1,\dots,k_r)
X_1^{k_1-1}\cdots X_r^{k_r-1}\,.
\end{align}
Similarly, a \emph{bimould} $B$ with values in $A$ is a family
\begin{align*}
B\bi{X_1,\dots,X_r}{Y_1,\dots,Y_r}
\in A[[X_1,\dots,X_r,Y_1,\dots,Y_r]]
\end{align*}
with $B\bi{\emptyset}{\emptyset}=1$. Its coefficients are defined by
\begin{align}\label{eq:cmesbimouldcoeff}
B\bi{X_1,\dots,X_r}{Y_1,\dots,Y_r}
={}&\sum_{\substack{k_1,\dots,k_r\geq1\\d_1,\dots,d_r\geq0}}
b\bi{k_1,\dots,k_r}{d_1,\dots,d_r}
\prod_{j=1}^rX_j^{k_j-1}\frac{Y_j^{d_j}}{d_j!}\,.
\end{align}
Thus our convention is the same as in Chapter~\ref{sec:mes}: there is no
factor $1/d_j!$ in the coefficient
$b\bi{k_1,\dots,k_r}{d_1,\dots,d_r}$, but divided powers are used in the
generating series.

The coefficient map of $B$ is the linear map
\begin{align*}
\varphi_B:\Q\langle\azbi\rangle&\longrightarrow A,\\
\ai{k_1,\dots,k_r}{d_1,\dots,d_r}
&\longmapsto b\bi{k_1,\dots,k_r}{d_1,\dots,d_r}\,.
\end{align*}
We call $B$ \emph{symmetril} if $\varphi_B$ is an algebra homomorphism for
the bi-stuffle product $\bist$ from Definition~\ref{def:fmesstuffle}. It is
\emph{swap invariant} if it satisfies the functional equation
\begin{align}\label{eq:cmesswap}
B\bi{X_1,\dots,X_r}{Y_1,\dots,Y_r}
={}&B\bi{Y_1+\dots+Y_r,\dots,Y_1+Y_2,Y_1}
{X_r,X_{r-1}-X_r,\dots,X_1-X_2}\,.
\end{align}
This is exactly the swap $\sigma$ from Definition~\ref{def:fmesswap}, or,
equivalently, the partition relation $P$ from
Section~\ref{subsec:doubleindexedqana}. Therefore, a symmetril and swap
invariant bimould with values in $A$ gives a realization $\fmes\rightarrow A$.

For two moulds, or two bimoulds, we define the product $\times$ by
deconcatenating the variables. For bimoulds it is given by
\begin{align}\label{eq:cmesmouldproduct}
(B\times C)\bi{X_1,\dots,X_r}{Y_1,\dots,Y_r}
={}&\sum_{j=0}^r
B\bi{X_1,\dots,X_j}{Y_1,\dots,Y_j}
C\bi{X_{j+1},\dots,X_r}{Y_{j+1},\dots,Y_r}\,.
\end{align}

\begin{proposition}[{\cite[Proposition~3.9]{BaBu}}]\label{prop:cmesproductsymmetril}
If $B$ and $C$ are symmetril, then $B\times C$ is symmetril.
\end{proposition}

\begin{proof}
Let $\Delta$ denote the deconcatenation coproduct on
$\Q\langle\azbi\rangle$. By~\eqref{eq:cmesmouldproduct}, the coefficient map
of $B\times C$ is the convolution
\begin{align*}
\varphi_{B\times C}
=m_A\circ(\varphi_B\otimes\varphi_C)\circ\Delta\,,
\end{align*}
where $m_A$ is multiplication in $A$. The deconcatenation coproduct is an
algebra homomorphism for the bi-stuffle product, and the convolution of two
algebra homomorphisms with values in a commutative algebra is again an algebra
homomorphism. This proves the claim.
\end{proof}

There is an analogous terminology for moulds. We call $Z$ symmetril if its
coefficient map is an algebra homomorphism for the stuffle product, and
symmetral if it is an algebra homomorphism for the shuffle product. To include
the regularization at $z_1$, define
\begin{align}\label{eq:cmesgammaZ}
\Gamma^Z(T)=\sum_{j\geq0}\gamma_j^ZT^j
:=\exp\left(\sum_{n\geq2}\frac{(-1)^n}{n}z(n)T^n\right)
\end{align}
and
\begin{align}\label{eq:cmesZgamma}
Z_\gamma(X_1,\dots,X_r)
:=\sum_{j=0}^r\gamma_j^Z
Z(X_1+\dots+X_{r-j},\dots,X_1+X_2,X_1)\,.
\end{align}
For $j=r$, the last factor in~\eqref{eq:cmesZgamma} is the depth-zero part
$Z(\emptyset)=1$. We say that $Z$ satisfies the \emph{extended double shuffle
equations} if $Z$ is symmetril and $Z_\gamma$ is symmetral. This formulation is
equivalent to the extended double shuffle relations of
Theorem~\ref{thm:eds} (see \cite[Definition~3.5]{BaBu} and \cite{IKZ,Rac}).

We now fix a $\Q$-valued mould $\gb$ satisfying the extended double shuffle
equations and whose depth-one part is
\begin{align}\label{eq:cmesbetadepthone}
\gb(X)
&=-\sum_{k\geq2}\frac{B_k}{2k!}X^{k-1}
=\frac12\left(\frac1X-\frac1{e^X-1}-\frac12\right)\,.
\end{align}
Write $\beta(k_1,\dots,k_r)$ for its coefficients. Together with
$\beta(0)=-\frac12$, its depth-one coefficients agree with the convention
\begin{align*}
\beta(k)=
\begin{cases}
-\dfrac{B_k}{2k!}\,,&k\text{ even},\\[0.2em]
0\,,&k\text{ odd}
\end{cases}
\qquad(k\geq0)\,.
\end{align*}
Such rational solutions exist by the work of Racinet, or by combining the
existence of rational associators with Furusho's theorem (see \cite[Remark~4.1]{BaBu}, \cite{Rac,D,F2}). The choice is not canonical, and
different choices occur starting in weight $8$. All combinatorial multiple
Eisenstein series below depend on this fixed choice.
In particular, $\beta$ defines a realization
\begin{align}\label{eq:cmesevbeta}
\operatorname{ev}_\beta:\fmz&\longrightarrow\Q,\\
\fzeta(k_1,\dots,k_r)&\longmapsto\beta(k_1,\dots,k_r)\,.
\end{align}

Already in depth two the stuffle relation determines
\begin{align}\label{eq:cmesbeta11}
0=\beta(1)^2=2\beta(1,1)+\beta(2),
\qquad\text{and therefore}\qquad \beta(1,1)=\frac1{48}\,.
\end{align}
This normalization will be important when we compare the construction with
the formal double zeta space in Section~\ref{sec:cmesdepthtwocomparison}.
Equivalently, the depth-two symmetral equation contains $2\gamma_2$.
The factor $2$ is missing in the corresponding displayed equation
\cite[(3.7)]{BaBu}, but follows directly from Definitions~3.4 and~3.5 there.

Next, we associate a bimould to any mould $Z$. Define
\begin{align}\label{eq:cmesBZ}
B^Z\bi{X_1,\dots,X_r}{Y_1,\dots,Y_r}
={}&\sum_{j=0}^r
Z_\gamma(Y_1,\dots,Y_j)Z(X_{j+1},\dots,X_r)\\
={}&\sum_{0\leq i\leq j\leq r}\gamma_i^Z
Z(Y_1+\dots+Y_{j-i},\dots,Y_1)Z(X_{j+1},\dots,X_r)\,.
\end{align}

\begin{proposition}[{\cite[Propositions~4.3 and~4.4]{BaBu}}]\label{prop:cmesBZ}
For every mould $Z$, the bimould $B^Z$ is swap invariant. If $Z$ satisfies the
extended double shuffle equations, then $B^Z$ is also symmetril.
\end{proposition}

\begin{proof}
After applying the swap in~\eqref{eq:cmesswap}, the second expression
in~\eqref{eq:cmesBZ} becomes
\begin{align*}
\sum_{0\leq i\leq j\leq r}\gamma_i^Z
Z(X_{r-j+i+1},\dots,X_r)
Z(Y_1+\dots+Y_{r-j},\dots,Y_1)\,.
\end{align*}
The change of variables $j'=r-j+i$ turns this into
\begin{align*}
\sum_{0\leq i\leq j'\leq r}\gamma_i^Z
Z(Y_1+\dots+Y_{j'-i},\dots,Y_1)Z(X_{j'+1},\dots,X_r)\,,
\end{align*}
which is~\eqref{eq:cmesBZ}.

Suppose now that $Z$ satisfies the extended double shuffle equations. The
bimould which is given by $Z(X_1,\dots,X_r)$ is symmetril. The bimould which
is given by $Z_\gamma(Y_1,\dots,Y_r)$ is symmetral, and since it is independent
of the $X$-variables, it is also symmetril. Formula~\eqref{eq:cmesBZ} is the
product of these two bimoulds. The claim therefore follows from
Proposition~\ref{prop:cmesproductsymmetril}.
\end{proof}

For our fixed rational solution, we write again $\gb$ for the bimould
$B^\gb$. Thus
\begin{align}\label{eq:cmesbetabimould}
\gb\bi{X_1,\dots,X_r}{Y_1,\dots,Y_r}
=\sum_{0\leq i\leq j\leq r}\gamma_i
\gb(Y_1+\dots+Y_{j-i},\dots,Y_1)
\gb(X_{j+1},\dots,X_r)\,,
\end{align}
where
\begin{align*}
\sum_{j\geq0}\gamma_jT^j
=\exp\left(\sum_{n\geq2}\frac{(-1)^n}{n}\beta(n)T^n\right)\,.
\end{align*}
By Proposition~\ref{prop:cmesBZ}, this bimould is symmetril and swap invariant.

\subsection{The construction}\label{sec:cmesconstruction}

We now combine the rational bimould $\gb$ with the double-indexed $q$-series
of Chapter~\ref{sec:mes}. For $m\geq1$, set
\begin{align}\label{eq:cmesLm}
L_m\bi{X}{Y}:=\frac{e^{X+mY}q^m}{1-e^Xq^m}\,.
\end{align}
Then the generating series of the $q$-series
$\mb{k_1,\dots,k_r}{d_1,\dots,d_r}$ from
Definition~\ref{def:doubleg} is
\begin{align}\label{eq:cmesrawbimould}
\ggen\bi{X_1,\dots,X_r}{Y_1,\dots,Y_r}
=\sum_{m_1>\dots>m_r>0}\prod_{j=1}^rL_{m_j}\bi{X_j}{Y_j}\,.
\end{align}
As we saw in Chapter~\ref{sec:mes}, this bimould is swap invariant, but its
coefficient map is an algebra homomorphism for $\gqsh$, not for the
bi-stuffle product $\bist$. We therefore need to correct its lower-weight
terms.

Following \cite[Definition~6.1]{BaBu}, first define the bimould $\gbr$ by
\begin{align}\label{eq:cmesbtilde}
\gbr\bi{X_1,\dots,X_r}{Y_1,\dots,Y_r}
:=\sum_{i=0}^r\frac{(-1)^i}{2^ii!}
\gb\bi{X_{i+1},\dots,X_r}{-Y_1,\dots,-Y_{r-i}}\,.
\end{align}
For each $m\geq1$, define the bimould $\gL_m$ by
\begin{align}\label{eq:cmesmathfrakLm}
&\gL_m\bi{X_1,\dots,X_r}{Y_1,\dots,Y_r}\\
&\quad:=\sum_{j=1}^r
\gb\bi{X_1-X_j,\dots,X_{j-1}-X_j}{Y_1,\dots,Y_{j-1}}
L_m\bi{X_j}{Y_1+\dots+Y_r}
\gbr\bi{X_r-X_j,\dots,X_{j+1}-X_j}{Y_r,\dots,Y_{j+1}}.
\nonumber
\end{align}
Its depth-one part is just $L_m\bi{X}{Y}$. The formula is the combinatorial
counterpart of the reduction of multitangent functions used in the Fourier
expansion of multiple Eisenstein series. Compare
Theorem~\ref{thm:reduction}.

The ordered sums of these bimoulds are collected in $\ggen^\ast$, defined as
in \cite[Definition~6.3]{BaBu} by
\begin{align}\label{eq:cmesgstar}
\ggen^\ast\bi{X_1,\dots,X_r}{Y_1,\dots,Y_r}:={}&
\sum_{\substack{1\leq j\leq r\\
0=r_0<r_1<\dots<r_j=r\\m_1>\dots>m_j>0}}
\prod_{i=1}^j
\gL_{m_i}\bi{X_{r_{i-1}+1},\dots,X_{r_i}}
{Y_{r_{i-1}+1},\dots,Y_{r_i}}\,.
\end{align}
This is the analogue of the ordered multitangent sums $\g^{\ast,M}$ from
Section~\ref{sec:stufflereg}.

\begin{definition}[{\cite[Definitions~6.3 and~6.4]{BaBu}}]\label{def:cmes}
Define the bimould
\begin{align}\label{eq:cmesGdefinition}
\gG:=\ggen^\ast\times\gb\,.
\end{align}
The \emph{combinatorial bi-multiple Eisenstein series} are its coefficients,
\begin{align}\label{eq:cmesGgenerating}
\gG\bi{X_1,\dots,X_r}{Y_1,\dots,Y_r}
={}&\sum_{\substack{k_1,\dots,k_r\geq1\\d_1,\dots,d_r\geq0}}
\Gb{k_1,\dots,k_r}{d_1,\dots,d_r}
\prod_{j=1}^rX_j^{k_j-1}\frac{Y_j^{d_j}}{d_j!}\,.
\end{align}
Their weight is $k_1+\dots+k_r+d_1+\dots+d_r$ and their depth is $r$. For
$d_1=\dots=d_r=0$, we write
\begin{align*}
G(k_1,\dots,k_r):=\Gb{k_1,\dots,k_r}{0,\dots,0}
\end{align*}
and call these the \emph{combinatorial multiple Eisenstein series}.
\end{definition}

In depth one, $\ggen^\ast=\ggen$, and hence
$\gG\bi{X}{Y}=\gb\bi{X}{Y}+\ggen\bi{X}{Y}$. Comparing coefficients gives
the following explicit formula.

\begin{proposition}[{\cite[Example~6.6(i)]{BaBu}}]\label{prop:cmesdepthone}
For $k\geq1$ and $d\geq0$, we have
\begin{align}\label{eq:cmesdepthone}
\Gb{k}{d}
={}&\delta_{d,0}\beta(k)+\delta_{k,1}d!\beta(d+1)
+\frac1{(k-1)!}\sum_{m,n\geq1}m^dn^{k-1}q^{mn}\,.
\end{align}
If $k>d$, then
\begin{align}\label{eq:cmesdepthonederivative}
\Gb{k}{d}
=\frac{(k-d-1)!}{(k-1)!}\left(\qdq\right)^dG(k-d)\,.
\end{align}
In particular, $G(k)=G_k$ for $k\geq2$, with the notation
of~\eqref{eq:defgk}.
\end{proposition}

\begin{proof}
Formula~\eqref{eq:cmesdepthone} follows from
\eqref{eq:cmesbetabimould}, \eqref{eq:cmesrawbimould} and the depth-one
identity $\gG=\gb+\ggen$. Applying $(\qdq)^d$ to the Fourier expansion of
$G(k-d)$ proves~\eqref{eq:cmesdepthonederivative}.
\end{proof}

Notice that~\eqref{eq:cmesdepthone} also gives
$\Gb{1}{0}=\g(1)$.\footnote{The Bernoulli-number formula printed in
\cite[Example~6.6(i)]{BaBu} does not cover this case, so one needs to be careful
here when using the convention $B_1=-1/2$.}

For example,
\begin{align*}
G(2)&=-\frac1{24}+q+3q^2+4q^3+7q^4+6q^5+O(q^6),\\
G(3)&=\frac12q+\frac52q^2+5q^3+\frac{21}{2}q^4+13q^5+O(q^6),\\
G(4)&=\frac1{1440}+\frac16q+\frac32q^2+\frac{14}{3}q^3
+\frac{73}{6}q^4+21q^5+O(q^6)\,.
\end{align*}

Before proving the main properties, let us also write the construction in
depth two. From~\eqref{eq:cmesGdefinition} and~\eqref{eq:cmesgstar}, we obtain
\cite[Example~6.6(ii)]{BaBu}
\begin{align}\label{eq:cmesdepthtwogen}
\gG\bi{X_1,X_2}{Y_1,Y_2}
={}&\ggen\bi{X_1,X_2}{Y_1,Y_2}
+\gb\bi{X_1-X_2}{Y_1}\ggen\bi{X_2}{Y_1+Y_2}\\
&+\ggen\bi{X_1}{Y_1+Y_2}\gbr\bi{X_2-X_1}{Y_2}
+\ggen\bi{X_1}{Y_1}\gb\bi{X_2}{Y_2}
+\gb\bi{X_1,X_2}{Y_1,Y_2}.\nonumber
\end{align}
At $Y_1=Y_2=0$, this becomes
\begin{align}\label{eq:cmesdepthtwomono}
\gG(X,Y)={}&\gb(X,Y)+\ggen(X,Y)+\gb(X-Y)\ggen(Y)
+\ggen(X)\left(\gb(Y-X)-\frac12\right)+\ggen(X)\gb(Y)\,.
\end{align}
Extracting the coefficient of $X^{k_1-1}Y^{k_2-1}$ gives the following
formula, which is not written out in \cite[Example~6.6(ii)]{BaBu}. If
$K=k_1+k_2$, then
\begin{align}\label{eq:cmesdepthtwoexplicit}
G(k_1,k_2)
={}&\beta(k_1,k_2)+\g(k_1,k_2)
-\frac12\delta_{k_2,1}\g(k_1)+\beta(k_2)\g(k_1)\nonumber\\
&+\sum_{j=k_1}^{K-1}(-1)^{j-k_1}
\binom{j-1}{k_1-1}\beta(j)\g(K-j)
+\sum_{j=k_2}^{K-1}(-1)^{j-k_2}
\binom{j-1}{k_2-1}\beta(j)\g(K-j)\,.
\end{align}
For example, the extended double shuffle relations for $\beta$ and
\eqref{eq:cmesdepthtwoexplicit} give
\begin{align}\label{eq:cmeslowdepthexamples}
G(1,1)&=\frac1{48}-\frac12\g(1)+\g(1,1),\nonumber\\
G(1,2)&=-\frac1{24}\g(1)+\g(1,2),\qquad
G(2,1)=-\frac12\g(2)+\g(2,1),\nonumber\\
G(2,2)&=\frac1{1920}-\frac18\g(2)+\g(2,2),\nonumber\\
G(2,3)&=-\frac1{24}\g(3)+\g(2,3),\nonumber\\
G(3,3)&=\frac1{120960}-\frac1{240}\g(2)+\g(3,3)\,.
\end{align}
Here $\beta(2,1)=\beta(3)=0$. The stuffle products
$\beta(1)\beta(2)$ and $\beta(2)\beta(3)$ give
$\beta(1,2)=0$ and $\beta(2,3)=-\beta(3,2)$, while the weight-five extended
double shuffle relations give $\beta(3,2)=0$. Finally,
\begin{align*}
2\beta(2,2)&=\beta(2)^2-\beta(4)=\frac1{960},\\
2\beta(3,3)&=\beta(3)^2-\beta(6)=\frac1{60480}\,.
\end{align*}
The first three examples show that the boundary terms in
\eqref{eq:cmesdepthtwoexplicit} are essential.

For the reversed pair in weight five, \cite[Example~6.6(iv)]{BaBu} gives
\begin{align*}
G(3,2)=\g(3,2)-\frac1{12}\g(3)\,,
\end{align*}
and therefore
\begin{align}\label{eq:cmesG32}
G(3,2)=-\frac1{24}q-\frac5{24}q^2+\frac1{12}q^3
+\frac58q^4+\frac{53}{12}q^5+O(q^6)\,.
\end{align}
In contrast, the analytic multiple Eisenstein series $\aG_{3,2}$ contains
zeta values and the rescaled series $\ag$. Compare
Example~\ref{ex:coproductfourier32}.

A first example in depth three is still short:
\begin{align}\label{eq:cmesG222}
G(2,2,2)
={}&\beta(2,2,2)
 +\bigl(3\beta(2,2)+3\beta(2)^2\bigr)\g(2)
 +5\beta(2)\g(2,2)+\g(2,2,2)\nonumber\\
={}&-\frac1{322560}+\frac{13}{1920}\g(2)
-\frac5{24}\g(2,2)+\g(2,2,2)\,.
\end{align}
For this, notice that the stuffle product gives
\begin{align*}
\beta(2,2,2)
=\frac{\beta(2)^3-3\beta(2)\beta(4)+2\beta(6)}6
=-\frac1{322560}\,.
\end{align*}
The remaining coefficients can be checked by extracting the coefficient of
$T^7$ in the repeated-two identity following
\cite[Example~6.14]{BaBu},
\begin{align*}
\sum_{r\geq0}G(\underbrace{2,\dots,2}_{r})T^{2r+1}
=\sum_{r\geq0}\g(\underbrace{2,\dots,2}_{r})
\left(2\sin\left(\frac T2\right)\right)^{2r+1}\,.
\end{align*}

\subsection{Symmetrility and swap invariance}\label{sec:cmesmaintheorem}

It remains to prove the two properties which make $\gG$ a realization of the formal
multiple Eisenstein series. The main ingredients are that every $\gL_m$ is
symmetril and that $\gG$ can be decomposed into swap invariant pieces.

\begin{proposition}[{\cite[Lemma~6.20]{BaBu}}]\label{prop:cmesLsymmetril}
For every $m\geq1$, the bimould $\gL_m$ is symmetril.
\end{proposition}

\begin{proof}
We give the main step and refer to \cite[Lemma~6.20]{BaBu} for the full
coefficient calculation. Write $q^m=e^{-T}$ and put
\begin{align*}
L_T(X)=\frac{e^{X-T}}{1-e^{X-T}}
=2\gb(X-T)-\frac1{X-T}-\frac12\,.
\end{align*}
Let $\gL_T$ be obtained from~\eqref{eq:cmesmathfrakLm} by replacing $L_m$ by
$L_T$ and omitting the factor $e^{m(Y_1+\dots+Y_r)}$. Introduce
\begin{align*}
\gb_T\bi{X_1,\dots,X_r}{Y_1,\dots,Y_r}
&=\gb\bi{X_1-T,\dots,X_r-T}{Y_1,\dots,Y_r},\\
\gbr_T\bi{X_1,\dots,X_r}{Y_1,\dots,Y_r}
&=\gbr\bi{X_r-T,\dots,X_1-T}{Y_r,\dots,Y_1}\,,
\end{align*}
and let $M_T$ have positive-depth part
\begin{align*}
M_T\bi{X_1,\dots,X_r}{Y_1,\dots,Y_r}
=\begin{cases}
\dfrac1{T-X_1},&r=1,\\
0,&r>1\,.
\end{cases}
\end{align*}
The extended double shuffle equations imply that the three bimoulds on the
right of
\begin{align}\label{eq:cmesLTfactorization}
\gL_T=\gb_T\times M_T\times\gbr_T\,.
\end{align}
are symmetril. To verify the factorization, one applies symmetrility to the
factors on the left of $L_T$ and the shuffle antipode to those on the right.
After the nonempty terms cancel, the remaining identity is
\begin{align*}
\left(\sum_{j\geq0}\beta(\underbrace{1,\dots,1}_{j})X^j\right)^2
=\frac{2}{X}\sinh\left(\frac X2\right)\,,
\end{align*}
which follows from the inverse series \cite[(3.6)]{BaBu}, Euler's formula for
$\beta(2n)$ and the reflection formula for the Gamma function. The complete
cancellation is given in \cite[Lemma~6.20]{BaBu}. In the first display of that proof, the
$Y$-arguments of the right-hand $\gbr$-factor are printed with minus signs.
\cite[Definition~6.1]{BaBu} and~\eqref{eq:cmesLTfactorization} require the
positive signs used in~\eqref{eq:cmesmathfrakLm}. Substituting $e^{-T}=q^m$ and restoring the
factor $e^{m(Y_1+\dots+Y_r)}$ now gives $\gL_m$ and proves the claim.
\end{proof}

The following elementary construction allows us to pass from the individual
$\gL_m$ to the ordered sum $\ggen^\ast$.

\begin{lemma}[{\cite[Lemma~6.21]{BaBu}}]\label{lem:cmesorderedsum}
Let $(B_m)_{m\geq1}$ be a family of symmetril bimoulds. For $M\geq1$, define
$C_M\bi{\emptyset}{\emptyset}=1$ and, in positive depth,
\begin{align*}
C_M\bi{X_1,\dots,X_r}{Y_1,\dots,Y_r}
:={}&\sum_{\substack{1\leq j\leq r\\0=r_0<\dots<r_j=r\\
M>m_1>\dots>m_j>0}}
\prod_{i=1}^jB_{m_i}\bi{X_{r_{i-1}+1},\dots,X_{r_i}}
{Y_{r_{i-1}+1},\dots,Y_{r_i}}\,.
\end{align*}
Then $C_M$ is symmetril.
\end{lemma}

\begin{proof}
For $M=1$, the positive-depth part is zero, so the assertion is clear. By
separating the terms with $m_1=M$ from those with $m_1<M$, one obtains
\begin{align*}
C_{M+1}=B_M\times C_M\,.
\end{align*}
The result now follows by induction from
Proposition~\ref{prop:cmesproductsymmetril}.
\end{proof}

Taking $B_m=\gL_m$ and passing coefficientwise to the limit $M\rightarrow
\infty$ gives the following.

\begin{corollary}[{\cite[Proposition~6.22]{BaBu}}]\label{cor:cmesgstarsymmetril}
The bimould $\ggen^\ast$ is symmetril.
\end{corollary}

It remains to understand the swap. We first record the transformation of the
right-hand correction in~\eqref{eq:cmesmathfrakLm}.

\begin{lemma}[{\cite[Lemma~6.24]{BaBu}}]\label{lem:cmesbtildeswap}
For $r\geq1$, we have
\begin{align*}
\gbr\bi{X_1,\dots,X_r}{Y_1,\dots,Y_r}=\gbr
\bi{-Y_1-\dots-Y_r,\dots,-Y_1-Y_2,-Y_1}
{-X_r,-X_{r-1}+X_r,\dots,-X_1+X_2}\,.
\end{align*}
\end{lemma}

\begin{proof}
Apply the swap invariance of $\gb$ to each summand in
\eqref{eq:cmesbtilde}. This gives
\begin{align*}
\gbr\bi{Y_1+\dots+Y_r,\dots,Y_1}
{X_r,X_{r-1}-X_r,\dots,X_1-X_2}
={}&\sum_{i=0}^r\frac{(-1)^i}{2^ii!}
\gb\bi{-X_{i+1},\dots,-X_r}{Y_1,\dots,Y_{r-i}}
=\gbr\bi{-X_1,\dots,-X_r}{-Y_1,\dots,-Y_r}\,.
\end{align*}
Replacing the variables by their inverse swap gives the stated formula.
\end{proof}

For $j\geq0$, define a family of power series $\gG_j$ as follows. Set
$\gG_0=\gb$ and set the depth-$r$ part of $\gG_j$ equal to zero for $j>r$.
For $1\leq j\leq r$, put
\begin{align}\label{eq:cmesGj}
\gG_j\bi{X_1,\dots,X_r}{Y_1,\dots,Y_r}:={}&
\sum_{\substack{0=r_0<r_1<\dots<r_j\leq r\\
m_1>\dots>m_j>0}}
\left(\prod_{i=1}^j
\gL_{m_i}\bi{X_{r_{i-1}+1},\dots,X_{r_i}}
{Y_{r_{i-1}+1},\dots,Y_{r_i}}\right)
\gb\bi{X_{r_j+1},\dots,X_r}{Y_{r_j+1},\dots,Y_r}\,.
\end{align}

\begin{proposition}[{\cite[Theorem~6.26 and Proposition~6.27]{BaBu}}]
\label{prop:cmesGjswap}
Each $\gG_j$ is swap invariant, and
\begin{align}\label{eq:cmesGsumGj}
\gG=\sum_{j=0}^r\gG_j
\end{align}
in depth $r$.
\end{proposition}

\begin{proof}
The second statement follows immediately by expanding
$\gG=\ggen^\ast\times\gb$ and collecting the terms according to the number
$j$ of positive integers $m_1>\dots>m_j$.

For the first statement, expand each $\gL_{m_i}$ in~\eqref{eq:cmesGj} by
choosing the position $n_i$, with $r_{i-1}<n_i\leq r_i$, of its factor
$L_{m_i}$. The sum over the resulting
$L$-factors is a component of the raw bimould $\ggen$ in
\eqref{eq:cmesrawbimould}. Under the swap, this component is fixed by the
partition relation for $\ggen$. The factors to its left are transformed by
the swap invariance of $\gb$, and the factors to its right by
Lemma~\ref{lem:cmesbtildeswap}. Reversing the order of the blocks and making
the change of summation variables
\begin{align*}
n'_i=r-r_{j-i+1}+1,\qquad r'_i=r-n_{j-i+1}+1
\end{align*}
returns exactly the summation range in~\eqref{eq:cmesGj}. Hence $\gG_j$ is
fixed by the swap. The full bookkeeping of the block endpoints is given in
\cite[Theorem~6.26]{BaBu}.
\end{proof}

The main result now follows.

\begin{theorem}[{\cite[Theorem~6.5]{BaBu}}]\label{thm:cmesmain}
The bimould $\gG$ is symmetril and swap invariant.
\end{theorem}

\begin{proof}
By Corollary~\ref{cor:cmesgstarsymmetril}, the bimould $\ggen^\ast$ is
symmetril, and $\gb$ is symmetril by Proposition~\ref{prop:cmesBZ}. Hence
$\gG=\ggen^\ast\times\gb$ is symmetril by
Proposition~\ref{prop:cmesproductsymmetril}. Its swap invariance follows from
the decomposition~\eqref{eq:cmesGsumGj} and
Proposition~\ref{prop:cmesGjswap}.
\end{proof}

\begin{corollary}\label{cor:cmesrealization}
There is an algebra homomorphism
\begin{align*}
\operatorname{ev}_G:\fmes&\longrightarrow\Q[[q]],\\
\Gf\bi{k_1,\dots,k_r}{d_1,\dots,d_r}
&\longmapsto\Gb{k_1,\dots,k_r}{d_1,\dots,d_r}\,.
\end{align*}
This is the $q$-series realization of Theorem~\ref{thm:fmesqrealization}.
It is also the realization used in the proof of
Theorem~\ref{thm:g2closedunderdif}.
\end{corollary}

\begin{proof}
Symmetrility says that the coefficient map is an algebra homomorphism for
$\bist$, and swap invariance says that it vanishes on the ideal
$\mathfrak I_\sigma$. Thus it factors through the quotient defining $\fmes$.
\end{proof}

\subsection{Double shuffle relations and the depth-two comparison}
\label{sec:cmesdepthtwocomparison}

The symmetrility of $\gG$ gives the stuffle product, whereas the combination
of symmetrility and swap invariance gives an analogue of the shuffle product.
In depth two both products can be written explicitly.

\begin{proposition}[{\cite[Proposition~6.7]{BaBu}}]\label{prop:cmesdepthtwoproduct}
For $k_1,k_2\geq1$, we have
\begin{align}\label{eq:cmesdepthtwoproduct}
G(k_1)G(k_2)
={}&G(k_1,k_2)+G(k_2,k_1)+G(k_1+k_2)\\
={}&\sum_{j=1}^{k_1+k_2-1}
\left(\binom{j-1}{k_1-1}+\binom{j-1}{k_2-1}\right)
G(j,k_1+k_2-j)+R_G(k_1,k_2),\nonumber
\end{align}
where
\begin{align}\label{eq:cmesRG}
R_G(k_1,k_2)=
\begin{cases}
\displaystyle
\frac{(k_1+k_2-3)!}{(k_1-1)!(k_2-1)!}\,
\qdq G(k_1+k_2-2),&k_1+k_2\geq3,\\[0.3cm]
G(2),&k_1=k_2=1\,.
\end{cases}
\end{align}
\end{proposition}

\begin{proof}
The first equality is the coefficient of
$\ai{k_1}{0}\bist\ai{k_2}{0}$ in the symmetrility of $\gG$. For the second
equality, apply the swap to the two depth-one factors, take their bi-stuffle
product, and apply the swap again. The two ways of interleaving the factors
give the binomial coefficients in~\eqref{eq:cmesdepthtwoproduct}. The merged
letter gives
\begin{align*}
\binom{k_1+k_2-2}{k_1-1}\Gb{k_1+k_2-1}{1}\,.
\end{align*}
If $k_1+k_2\geq3$, formula~\eqref{eq:cmesdepthonederivative} turns this into
the first case of~\eqref{eq:cmesRG}. If $k_1=k_2=1$, the depth-one swap gives
$\Gb{1}{1}=G(2)$.
\end{proof}

For $(k_1,k_2)=(2,3)$, the two lines of
\eqref{eq:cmesdepthtwoproduct} give
\begin{align*}
G(2)G(3)
&=G(2,3)+G(3,2)+G(5)\\
&=G(2,3)+3G(3,2)+6G(4,1)+\qdq G(3)\,.
\end{align*}
This gives
\begin{align}\label{eq:cmesG5relation}
G(5)=2G(3,2)+6G(4,1)+\qdq G(3)\,.
\end{align}
The derivative is precisely the obstruction to an ordinary double shuffle
relation.

We now compare the depth-two series with the combinatorial double Eisenstein
series constructed in Chapter~\ref{sec:mdandmzv}. These two constructions
agree away from the boundary, but they are not identical in general. To make
the distinction visible, write
\begin{align*}
G^{\mathrm{FDZ}}(k_1,k_2)
\end{align*}
for the series defined by Theorems~\ref{thm:grealization}
and~\ref{thm:betarealization}, and write $G^{\mathrm{BB}}(k_1,k_2)$ for the
depth-two specialization of Definition~\ref{def:cmes}. For
$\bullet\in\{\mathrm{FDZ},\mathrm{BB}\}$, set
\begin{align*}
G^\bullet(X,Y)
:=\sum_{k_1,k_2\geq1}G^\bullet(k_1,k_2)X^{k_1-1}Y^{k_2-1}\,.
\end{align*}

The rational depth-two solution in Chapter~\ref{sec:mdandmzv} satisfies
\begin{align*}
\beta^{\mathrm{FDZ}}(1,1)=0\,,
\end{align*}
whereas the stuffle normalization~\eqref{eq:cmesbeta11} gives
$\beta^{\mathrm{BB}}(1,1)=1/48$. We therefore match the two rational
depth-two solutions by setting
\begin{align}\label{eq:cmesbetamatching}
\gb^{\mathrm{BB}}(X,Y)
=\gb^{\mathrm{FDZ}}(X,Y)+\frac1{48}\,.
\end{align}
For the comparison, $G^{\mathrm{BB}}$ will from now on denote the depth-two
series obtained with this matched constant term. Only the depth-two solution
in~\eqref{eq:cmesbetamatching} is used here. No choice of a full rational
mould with this prescribed depth-two part is needed. For a different full
rational solution, the positive-$q$ part of the following comparison is
unchanged, and one merely adds the difference of the constant terms.

Put
\begin{align}\label{eq:cmesAderivative}
A(T):=\qdq\sum_{k\geq1}G(k)\frac{T^{k-1}}{k}\,.
\end{align}

\begin{theorem}\label{thm:cmesfdzbbcomparison}
With the matching~\eqref{eq:cmesbetamatching}, we have
\begin{align}\label{eq:cmesfdzbbgen}
G^{\mathrm{FDZ}}(X,Y)-G^{\mathrm{BB}}(X,Y)
=\frac12(X-Y)A(Y)+\frac12XA(X)+\frac12G(2)\,.
\end{align}
Equivalently,
\begin{align}\label{eq:cmesfdzbbcoeff}
G^{\mathrm{FDZ}}(k_1,k_2)-G^{\mathrm{BB}}(k_1,k_2)={}&
\begin{cases}
\frac12G(2),&(k_1,k_2)=(1,1),\\[0.1cm]
-\dfrac{1}{2(k_2-1)}\qdq\g(k_2-1),&k_1=1,\ k_2\geq2,\\[0.3cm]
\qdq\g(1),&(k_1,k_2)=(2,1),\\[0.1cm]
\dfrac{1}{2k_2}\qdq\g(k_2),&k_1=2,\ k_2\geq2,\\[0.3cm]
\dfrac{1}{2(k_1-1)}\qdq\g(k_1-1),&k_1\geq3,\ k_2=1,\\[0.3cm]
0,&k_1\geq3,\ k_2\geq2\,.
\end{cases}
\end{align}
\end{theorem}

\begin{proof}
The positive-$q$ part of the FDZ construction is the series $\gh(X,Y)$ in
Theorem~\ref{thm:grealization}. The positive-$q$ part of the BB construction
is obtained from~\eqref{eq:cmesdepthtwomono} by omitting the first term
$\gb(X,Y)$. Since $\gb$ is odd, all non-derivative terms cancel in their
difference. The remaining terms are
\begin{align*}
\frac12(X-Y)A(Y)+\frac12XA(X)+\frac12\g(2)\,.
\end{align*}
By~\eqref{eq:cmesbetamatching}, the difference of the constant terms is
$-1/48=\frac12\beta(2)$. This proves~\eqref{eq:cmesfdzbbgen}. Extracting the
coefficient of $X^{k_1-1}Y^{k_2-1}$ gives the six cases in
\eqref{eq:cmesfdzbbcoeff}.
\end{proof}

For weights at least $3$, the correction in~\eqref{eq:cmesfdzbbcoeff} can be
written more conceptually as
\begin{align}\label{eq:cmesbkumcorrection}
\frac12\left(
\delta_{k_2,1}\Gb{k_1}{1}
-\delta_{k_1,1}\Gb{k_2}{1}
+\delta_{k_1,2}\Gb{k_2+1}{1}
\right)\,.
\end{align}
This is exactly the splitting correction in
Theorem~\ref{thm:fdesdzmaps}. Compare \cite[Proposition~2.6]{BKuM}. To see
this, notice that for $K=k_1+k_2\geq3$, applying the
combinatorial realization after the map
$\iota_K^{\mathrm{DE}}$ to the formula for
$s_K^{\mathrm{DE}}(Z_{k_1,k_2})$ gives
\begin{align*}
(\operatorname{ev}_G\circ\iota_K^{\mathrm{DE}})
\left(s_K^{\mathrm{DE}}(Z_{k_1,k_2})\right)
={}&G^{\mathrm{BB}}(k_1,k_2)+\frac12\left(
\delta_{k_2,1}\Gb{k_1}{1}
-\delta_{k_1,1}\Gb{k_2}{1}
+\delta_{k_1,2}\Gb{k_2+1}{1}
\right)\\
={}&G^{\mathrm{FDZ}}(k_1,k_2)\,.
\end{align*}
Thus the abstract splitting~\eqref{eq:fdessplit} is realized by the boundary
correction~\eqref{eq:cmesbkumcorrection}. The case $(1,1)$ is exceptional
and gives $G(2)/2$.

There is also a direct comparison with the double Eisenstein series of
Gangl--Kaneko--Zagier. Let $Z^{\mathrm{GKZ}}_{k_1,k_2}$ denote their
zero-constant-term series and put
\begin{align*}
K(k_1,k_2):=(-1)^{k_1+k_2}Z^{\mathrm{GKZ}}_{k_1,k_2}\,.
\end{align*}
Then
\begin{align}\label{eq:cmesgkzcomparison}
G^{\mathrm{FDZ}}(k_1,k_2)
=K(k_1,k_2)+\beta^{\mathrm{FDZ}}(k_1,k_2)
+\frac14\delta_{k_1,1}\g(k_2)\,.
\end{align}
In the convergence range $k_1\geq3$, $k_2\geq2$, this follows from
\cite[(17)]{GKZ}. The boundary terms are given by
\cite[(54)--(57)]{GKZ}. Combining~\eqref{eq:cmesgkzcomparison} with
Theorem~\ref{thm:cmesfdzbbcomparison} gives the complete comparison of all
three depth-two constructions. Notice that this comparison concerns lower
indices $d_1=d_2=0$. There is no corresponding three-way identification for
arbitrary bi-indices.

For illustration, the first terms in four characteristic cases are
\begin{align*}
(3,2):\quad
K&=G^{\mathrm{FDZ}}=G^{\mathrm{BB}}
=-\frac1{24}q-\frac5{24}q^2+\frac1{12}q^3+O(q^4),\\[0.1cm]
(2,2):\quad
K&=\frac18q+\frac98q^2+\frac72q^3+O(q^4),\\
G^{\mathrm{FDZ}}&=\frac1{1920}+K,\qquad
G^{\mathrm{BB}}=\frac1{1920}-\frac18q-\frac38q^2+\frac12q^3+O(q^4),\\[0.1cm]
(2,1):\quad
K&=G^{\mathrm{FDZ}}=\frac12q+\frac52q^2+5q^3+O(q^4),\\
G^{\mathrm{BB}}&=-\frac12q-\frac32q^2-q^3+O(q^4),\\[0.1cm]
(1,1):\quad
K&=-\frac14q+\frac32q^3+O(q^4),\\
G^{\mathrm{FDZ}}&=\frac12q^2+2q^3+O(q^4),\qquad
G^{\mathrm{BB}}=\frac1{48}-\frac12q-q^2+O(q^4)\,.
\end{align*}
These expansions make both the boundary corrections and the exceptional
weight-two constant term visible.

\subsection{The \texorpdfstring{$q$}{q}-series spaces and their limits}
\label{sec:cmeslimits}

We now return to arbitrary depth. By construction, every combinatorial
bi-multiple Eisenstein series is a
rational linear combination of the double-indexed $q$-series $\mb{\kk}{\kl}$.
The converse also holds.

\begin{theorem}[{\cite[Proposition~6.15 and Remark~6.16(i)]{BaBu}}]
\label{thm:cmesspans}
We have
\begin{align}\label{eq:cmesbispan}
\mz_q
&=\Q+\left\langle
\Gb{k_1,\dots,k_r}{d_1,\dots,d_r}
\ \middle|\ r\geq1,\ k_j\geq1,\ d_j\geq0
\right\rangle_\Q,\\
\gs
&=\Q+\left\langle
G(k_1,\dots,k_r)
\ \middle|\ r\geq1,\ k_j\geq1
\right\rangle_\Q.
\label{eq:cmesmonospan}
\end{align}
\end{theorem}

\begin{proof}
The definition of $\gG$ gives a triangular expansion
\begin{align}\label{eq:cmesbitriangular}
\Gb{k_1,\dots,k_r}{d_1,\dots,d_r}
=\mb{k_1,\dots,k_r}{d_1,\dots,d_r}
+\text{terms involving $\mb{\kk'}{\kl'}$ of smaller depth and weight}\,.
\end{align}
Thus every combinatorial bi-multiple Eisenstein series belongs to $\mz_q$.
Conversely, formula~\eqref{eq:cmesdepthone} starts the induction, and
\eqref{eq:cmesbitriangular} expresses every $\mb{\kk}{\kl}$ as a linear
combination of combinatorial bi-multiple Eisenstein series. Equation
\eqref{eq:cmesbispan} now follows from
Theorem~\ref{thm:zqspannedbydoubleg}.

After setting all $Y$-variables equal to zero, the decomposition used in the
proof of Proposition~\ref{prop:cmesGjswap} contains only the single-indexed
series $\g(\kk)$. It gives the same triangular argument in both directions,
which proves~\eqref{eq:cmesmonospan}.
\end{proof}

In particular, asking whether the combinatorial multiple Eisenstein series
already span all combinatorial bi-multiple Eisenstein series is exactly the
question whether
\begin{align*}
\gs=\mz_q\,.
\end{align*}
By Theorem~\ref{thm:hmwqspaces} this is the case, i.e. the combinatorial multiple
Eisenstein series do span all combinatorial bi-multiple Eisenstein series. Notice
that this was still open when \cite{BaBu} was written, so the corresponding
statement there is formulated as a conjecture. It is
conjectured that all relations among the combinatorial bi-multiple Eisenstein
series follow from symmetrility and swap invariance. Since these defining
relations are homogeneous, one therefore expects all relations to be
homogeneous in weight (see \cite[Remark~6.11]{BaBu} and \cite[Conjecture~3]{BK2}).

It remains to make precise in which sense the combinatorial multiple Eisenstein
series interpolate between $\beta$ and multiple zeta values. For an admissible
index $k_1\geq2$, the ordinary limit is
\begin{align}\label{eq:cmesadmissiblelimit}
\lim_{q\rightarrow1}(1-q)^{k_1+\dots+k_r}G(k_1,\dots,k_r)
=\zeta(k_1,\dots,k_r)\,.
\end{align}
If $k_1=1$, a regularization is needed. Regard $G$ as the algebra
homomorphism
\begin{align*}
G:(\HH^1,\ast)&\longrightarrow\gs,\\
z_{k_1}\cdots z_{k_r}&\longmapsto G(k_1,\dots,k_r)\,.
\end{align*}
By the polynomial decomposition $\HH^1=\HH^0[z_1]$ for the stuffle product,
every $w\in\HH^1$ has a unique expression
\begin{align*}
w=\sum_{j\geq0}w_j\ast z_1^{\ast j},\qquad w_j\in\HH^0\,.
\end{align*}
For a word $w=z_{k_1}\cdots z_{k_r}$ of weight $k$, define
\begin{align}\label{eq:cmesregularizedlimit}
{\lim_{q\rightarrow1}}^{\!*}(1-q)^kG(w)
:=\lim_{q\rightarrow1}(1-q)^kG(w_0)=\zeta(w_0)\,.
\end{align}

\begin{proposition}[{\cite[Proposition~6.17]{BaBu}}]\label{prop:cmeslimits}
For all $k_1,\dots,k_r\geq1$, we have
\begin{align}\label{eq:cmesbothlimits}
\lim_{q\rightarrow0}G(k_1,\dots,k_r)
&=\beta(k_1,\dots,k_r),\\
{\lim_{q\rightarrow1}}^{\!*}(1-q)^{k_1+\dots+k_r}
G(k_1,\dots,k_r)
&=\zeta^\ast(k_1,\dots,k_r;0)\,.
\end{align}
For $k_1\geq2$, the second identity is the ordinary limit
\eqref{eq:cmesadmissiblelimit}.
\end{proposition}

\begin{proof}
At $q=0$, all $L_m$ vanish, and hence the constant term of $\gG$ is the
bimould $\gb$. Since $\gb_\gamma$ is symmetral and its depth-one value at zero
is $\beta(1)=0$, the shuffle relation gives
$\gb_\gamma(0,\dots,0)=0$ in every positive depth. Setting all $Y$-variables
equal to zero in~\eqref{eq:cmesbetabimould} therefore leaves precisely the
mould $\gb(X_1,\dots,X_r)$. This proves the first identity.

Put $K=k_1+\dots+k_r$. The mono-index expansion used in the proof of
\cite[Proposition~6.15]{BaBu} has leading term $\g(k_1,\dots,k_r)$, while
every other nonconstant term has weight less than $K$. After multiplication
by $(1-q)^K$, these lower-weight terms vanish as $q\rightarrow1$. This is the
highest-weight argument preceding \cite[Proposition~6.17]{BaBu}. The limit of
the leading term is Proposition~\ref{prop:gisqmzv}, which proves
\eqref{eq:cmesadmissiblelimit}.

Finally, $G$ and $\zeta^\ast(\,\cdot\,;0)$
are both algebra homomorphisms for the stuffle product. Applying them to the
polynomial decomposition in $z_1$ proves the regularized statement.
\end{proof}

The construction and its two limits can be summarized by
\begin{center}
\begin{tikzcd}[column sep=large,row sep=large]
\fmes \arrow[r,"{\operatorname{ev}_G}"] \arrow[d,"{\pi_f}"']
& \Q[[q]] \arrow[d,"{q\to0}"] \arrow[r,dashed,"{q\to1\ (d_j=0)}"]
& \mz\\
\fmz \arrow[r,"{\operatorname{ev}_\beta}"'] & \Q
\end{tikzcd}
\end{center}
The lower map is the realization~\eqref{eq:cmesevbeta}.
The dashed arrow is shorthand for the weight-scaled, regularized limit on the
lower-index-zero values $G(k_1,\dots,k_r)$. It is not a map defined on all of
$\Q[[q]]$.

The two limits also explain what happens when a linear combination of
combinatorial multiple Eisenstein series is a modular form.

\begin{proposition}[{\cite[Proposition~6.19]{BaBu}}]
Suppose that $k\geq4$ and that
\begin{align*}
F=\sum_{\substack{r\geq1\\k_1+\dots+k_r=k}}
c_{k_1,\dots,k_r}G(k_1,\dots,k_r)
\end{align*}
is a modular form of weight $k$. Then
\begin{align*}
\sum_{\substack{r\geq1\\k_1+\dots+k_r=k}}
c_{k_1,\dots,k_r}\zeta^\ast(k_1,\dots,k_r;0)
=(2\pi i)^k
\sum_{\substack{r\geq1\\k_1+\dots+k_r=k}}
c_{k_1,\dots,k_r}\beta(k_1,\dots,k_r)\,.
\end{align*}
\end{proposition}

If $c_{k_1,\dots,k_r}=0$ whenever $k_1=1$, the left-hand side is the
corresponding sum of ordinary multiple zeta values. In
\cite[Proposition~6.19]{BaBu}, it is written as $\zeta(k_1,\dots,k_r)$ also
when $k_1=1$. Its proof uses the regularized statement of
\cite[Proposition~6.17]{BaBu}, and therefore has to be read in the sense above.

\begin{proof}
Put
\begin{align*}
W=\sum_{\substack{r\geq1\\k_1+\dots+k_r=k}}
c_{k_1,\dots,k_r}z_{k_1}\cdots z_{k_r}
=\sum_{j\geq0}W_j\ast z_1^{\ast j},\qquad W_j\in\HH^0\,.
\end{align*}
Since $G$ is an algebra homomorphism,
$F=\sum_jG(W_j)G(1)^j$. We also have
\begin{align*}
G(1)=\sum_{m\geq1}\frac{q^m}{1-q^m},\qquad
(1-q)G(1)=\log\left(\frac1{1-q}\right)+O(1)\,.
\end{align*}
Together with the admissible limits, this shows by descending powers of the
logarithm that $\zeta(W_0)=\zeta^\ast(W;0)$ is the finite part of the
weight-scaled expansion of $F$ at $q=1$. A modular form has no logarithmic
terms at a cusp, so this finite part agrees with the ordinary limit.

The first limit in Proposition~\ref{prop:cmeslimits} gives
\begin{align*}
a_0=\sum_{\substack{r\geq1\\k_1+\dots+k_r=k}}
c_{k_1,\dots,k_r}\beta(k_1,\dots,k_r)\,.
\end{align*}
The transformation $F(-1/\tau)=\tau^kF(\tau)$ gives
\begin{align*}
\lim_{q\rightarrow1}(1-q)^kF(q)=(2\pi i)^ka_0\,.
\end{align*}
The second limit in Proposition~\ref{prop:cmeslimits} identifies the left-hand
side with the stated sum of regularized multiple zeta values.
\end{proof}

\subsection{Quasimodularity and derivatives}\label{sec:cmesderivatives}

Finally, we return to two structures from Chapter~\ref{sec:mes} and
Section~\ref{sec:fmes}: quasimodular forms and derivatives. Some combinatorial
bi-multiple Eisenstein series are quasimodular forms. The simplest systematic
family consists of repeated equal bi-indices.

\begin{proposition}[{\cite[Proposition~6.13]{BaBu}}]\label{prop:cmesrepeatedmodular}
Let $k\geq1$ and $d\geq0$ with $k+d$ even. For every $r\geq1$, the series
\begin{align*}
\Gb{\overbrace{k,\dots,k}^{r}}{\underbrace{d,\dots,d}_{r}}
\end{align*}
is a quasimodular form of weight $r(k+d)$ and depth at most
$r\min(d+1,k)$, where the depth of a quasimodular form is its degree in
$G(2)$.
\end{proposition}

\begin{proof}
The repeated-letter identity in a quasi-shuffle algebra gives
\begin{align}\label{eq:cmesrepeatedexp}
1+\sum_{r\geq1}\Gb{\overbrace{k,\dots,k}^{r}}
{\underbrace{d,\dots,d}_{r}}T^r
=\exp\left(\sum_{r\geq1}(-1)^{r-1}\Gb{rk}{rd}\frac{T^r}{r}\right),
\end{align}
compare \cite[(32)]{HI} and Theorem~\ref{thm:explog}. By
Proposition~\ref{prop:cmesdepthone}, together with the depth-one swap, each
$\Gb{nk}{nd}$ is a derivative of an Eisenstein series of even weight and is
therefore quasimodular. More precisely, its depth is at most
\begin{align*}
\min(nd+1,nk)\,.
\end{align*}
If $k>d$, it is an $nd$-th derivative of $G(n(k-d))$. If $k\leq d$,
apply the depth-one swap first and obtain a derivative of order $nk-1$ of
$G(n(d-k)+2)$. In both cases, one has to allow depth one when the Eisenstein
series at the beginning is $G(2)$. Since
\begin{align*}
\min(nd+1,nk)\leq n\min(d+1,k)\,,
\end{align*}
comparing coefficients of powers of $T$ in~\eqref{eq:cmesrepeatedexp} gives
the stated bound.
\end{proof}

The smaller bound $r\min(d,k-1)$ printed at the end of the proof of
\cite[Proposition~6.13]{BaBu} does not hold for the usual depth of
quasimodular forms. For example,
\begin{align*}
G(2,2)=\frac12\bigl(G(2)^2-G(4)\bigr)
\end{align*}
has depth two.

For example, the case $k=d=1$ and $r=2$ gives
\begin{align}\label{eq:cmesbiquasimodularexample}
\Gb{1,1}{1,1}
=\frac12G(2)^2-\frac12\qdq G(2)\,.
\end{align}
This is \cite[Example~6.14]{BaBu}. It also follows directly from
symmetrility, the depth-one swap and Proposition~\ref{prop:cmesdepthone}.

The derivative of an arbitrary combinatorial bi-multiple Eisenstein series
has a particularly simple form.

\begin{theorem}[{\cite[Proposition~6.29]{BaBu}}]\label{thm:cmesbiderivative}
For $k_1,\dots,k_r\geq1$ and $d_1,\dots,d_r\geq0$, we have
\begin{align}\label{eq:cmesbiderivative}
\qdq\Gb{k_1,\dots,k_r}{d_1,\dots,d_r}
=\sum_{i=1}^rk_i
\Gb{k_1,\dots,k_i+1,\dots,k_r}
{d_1,\dots,d_i+1,\dots,d_r}\,.
\end{align}
\end{theorem}

\begin{proof}
It is enough to prove the generating-series identity
\begin{align}\label{eq:cmesGderivativegen}
\qdq\gG\bi{X_1,\dots,X_r}{Y_1,\dots,Y_r}
=\sum_{i=1}^r\frac{\partial^2}{\partial X_i\partial Y_i}
\gG\bi{X_1,\dots,X_r}{Y_1,\dots,Y_r}\,.
\end{align}
The building block~\eqref{eq:cmesLm} satisfies
\begin{align*}
\qdq L_m\bi{X}{Y}
=\frac{\partial^2}{\partial X\partial Y}L_m\bi{X}{Y}\,.
\end{align*}
Expand $\gG$ into the pieces $\gG_j$ from~\eqref{eq:cmesGj} and mark in each
$\gL_m$ the position of its $L_m$-factor, as in the proof of
Proposition~\ref{prop:cmesGjswap}. The remaining factors contain only the sum
and difference variables occurring in $\gb$ and $\gbr$. When the mixed
derivatives are summed over all positions, their product-rule terms cancel
and their own mixed derivatives vanish. Thus only the marked $L_m$-factors
remain, and for them the preceding identity is exactly $\qdq$. This proves
\eqref{eq:cmesGderivativegen}. The complete bookkeeping is given in
\cite[Proposition~6.29]{BaBu}. Extracting coefficients gives
\eqref{eq:cmesbiderivative}, where the factor $k_i$ comes from differentiating
$X_i^{k_i}$.
\end{proof}

For the mono-indexed series, the derivative can be expressed entirely in the
language of the usual stuffle and shuffle products. Extend $G$ linearly to
$\HH^1$ as above.

\begin{corollary}[{\cite[Corollary~6.31]{BaBu}}]\label{cor:cmeswordderivative}
For every $w\in\HH^1$, we have
\begin{align}\label{eq:cmeswordderivative}
\qdq G(w)=G(z_2\ast w-z_2\sh w)\,.
\end{align}
If $h:\HH^1\rightarrow\HH^1$ is defined by
\begin{align*}
h(w)=z_2\ast w-z_2\sh w\,,
\end{align*}
then, for all $u,v\in\HH^1$,
\begin{align}\label{eq:cmesderivativekernel}
h(u\ast v)-h(u)\ast v-u\ast h(v)\in\ker G\,.
\end{align}
\end{corollary}

\begin{proof}
Set the lower indices equal to zero in
Theorem~\ref{thm:cmesbiderivative} and use the swap to return the resulting
bi-indices to lower weight zero. Comparing the resulting insertions and merged
letters with the shuffle of $z_2=xy$ and an arbitrary word gives
\eqref{eq:cmeswordderivative}. Equivalently, this is the coefficient form of
\cite[Proposition~6.30]{BaBu}. For~\eqref{eq:cmesderivativekernel}, use that
$G$ is an algebra homomorphism and $\qdq$ is a derivation:
\begin{align*}
G(h(u\ast v))
&=\qdq\bigl(G(u)G(v)\bigr)\\
&=G(h(u))G(v)+G(u)G(h(v))\,.
\end{align*}
\end{proof}

Taking $w=z_1$ in~\eqref{eq:cmeswordderivative} gives
\begin{align}\label{eq:cmesG1derivative}
\qdq G(1)=G(3)-G(2,1)\,.
\end{align}
Taking $u=v=z_1$ in~\eqref{eq:cmesderivativekernel} gives the weight-four
relation
\begin{align}\label{eq:cmesweightfourrelation}
G(4)=2G(2,2)-2G(3,1)\,.
\end{align}
On the other hand, combinatorial multiple Eisenstein series do not satisfy all
relations of multiple zeta values. For example,
\begin{align}\label{eq:cmesG211}
G(2,1,1)
=\beta(2,1,1)+\frac16\g(2)-\g(2,1)+\g(2,1,1)\,,
\end{align}
as in \cite[Example~6.6(v)]{BaBu},
and hence $G(2,1,1)\neq G(4)$, even though
$\beta(2,1,1)=\beta(4)$. More explicitly,
\begin{align*}
G(2,1,1)-G(4)=-q^2+O(q^3)\,.
\end{align*}
Thus duality is not built into the construction.

Theorem~\ref{thm:cmesbiderivative} is the concrete $q$-series version of the
formal derivative $\GD$ in Section~\ref{sec:fmes}. Together with the
analytic realization in Theorem~\ref{thm:fmesanalyticrealization}, it shows
that the same formal derivative controls both the derivatives of
combinatorial multiple Eisenstein series and those of the analytic multiple
Eisenstein series studied in Section~\ref{sec:mesderiv}.

This completes the construction of the $q$-series realization announced in
Theorem~\ref{thm:fmesqrealization}: its products, its two limits and its
derivative are now explicit.

Let us finally mention that the balanced multiple $q$-zeta values introduced
by Burmester give another presentation of these $q$-series. The Hopf algebra
isomorphism in \cite[Theorem~7.11]{Bu2} identifies the bi-stuffle product and
swap with $\ast_b$ and $\tau_b$. Under the corresponding change of variables,
$\gG$ becomes the generating series of the balanced multiple $q$-zeta values (see \cite[Definition~10.1]{Bu2}). In particular,
$\zeta_q(k_1,\dots,k_r)=G(k_1,\dots,k_r)$ for $k_1,\dots,k_r\geq1$ (see \cite[Example~10.6(i)]{Bu2}). In this sense, the balanced multiple $q$-zeta
values and the combinatorial bi-multiple Eisenstein series are two coordinate
systems for the same space $\mz_q$.

The last section returns to the finite side of the four-world picture and
introduces formal finite multiple zeta values.

%% file: chap_FormalFiniteMZV.tex
\section{Formal finite multiple zeta values}\label{sec:formalfinitemzv}

In Chapter~\ref{sec:fmzv}, we saw that finite and symmetric multiple zeta
values satisfy the stuffle product and the same linear shuffle relations, and
Conjecture~\ref{conj:fmzvallrelations} predicts that these are all
relations among finite multiple zeta values. This suggests introducing a
universal algebra in which precisely these relations hold. The material in
this section is based on joint work of the author and Risan~\cite{BaR}, which
grew out of Risan's master's thesis~\cite{R}.

\subsection{The formal finite algebra}

Recall that $\overline{w}$ denotes the reverse of a word $w\in\HH^1$. Let
$\mathfrak L$ be the ideal of the stuffle algebra $\HH^1_\ast$ generated by
\begin{align}\label{eq:formalfiniteideal}
w\sh v-(-1)^{\wt(w)}\overline{w}v\qquad(w,v\in\HH^1)\,,
\end{align}
where the product on the right-hand side is concatenation.

\begin{definition}[{\cite[Definition~4.1]{R}}]\label{def:formalfinitemzv}
The algebra of \emph{formal finite multiple zeta values} is
\begin{align*}
\fmza=\quotient{\HH^1_\ast}{\mathfrak L}\,.
\end{align*}
For $k_1,\dots,k_r\geq1$, the class of $z_{k_1}\cdots z_{k_r}$ is denoted by
$\fza(k_1,\dots,k_r)$, and we set $\fza(\emptyset)=1$.
\end{definition}

This algebra is graded by weight. It is also equipped with the depth
filtration induced by the length of a word. By construction, its elements
satisfy the stuffle product and the formal linear shuffle relations
\begin{align}\label{eq:formalfiniteproducts}
\fza(w\ast v)&=\fza(w)\fza(v),\notag\\
\fza(w\sh v)&=(-1)^{\wt(w)}\fza(\overline{w}v)\,.
\end{align}
Taking $v=1$ in the second relation gives
\begin{align}\label{eq:formalfinitereversal}
\fza(k_r,\dots,k_1)=(-1)^{k_1+\dots+k_r}
\fza(k_1,\dots,k_r)\,.
\end{align}

\begin{proposition}[{\cite[Proposition~4.2]{R}}]
There are surjective algebra homomorphisms
\begin{align}\label{eq:formalfiniteclassicalrealizations}
\pi_{\ma}:\fmza&\longrightarrow\mza,\notag\\
\fza(\kk)&\longmapsto\za(\kk),\notag\\
\pi_{\ms}:\fmza&\longrightarrow\quotient{\mz}{\pi^2\mz},\notag\\
\fza(\kk)&\longmapsto\zs(\kk)\,.
\end{align}
\end{proposition}

\begin{proof}
Theorem~\ref{thm:linearshufflefmzv} and
Theorem~\ref{thm:linearshufflesmzv} show that the generators of
$\mathfrak L$ are in the kernels of the two maps. Both maps also respect the
stuffle product. The first one is surjective by the definition of $\mza$,
whereas the surjectivity of the second one follows from
Theorem~\ref{thm:symspanmz}.
\end{proof}

Thus every relation proved in $\fmza$ immediately gives a relation for
both finite and symmetric multiple zeta values. In contrast to the classical
finite algebra, the formal algebra is defined entirely by generators and
relations, so non-vanishing questions do not enter its construction.

The first calculations already explain why the depth-one symbols are not the
right formal analogues of single zeta values.

\begin{proposition}[{\cite[Propositions~4.3, 4.5 and~4.6]{R}}]\label{prop:formalfinitedepthtwo}
For $a,b\geq1$, we have
\begin{enumerate}[(i)]
\item $\fza(a)=0$.
\item $\fza(a,b)=0$ if $a+b$ is even.
\end{enumerate}
Define
\begin{align}\label{eq:formaltruezeta}
\fZ(1)=0,\qquad
\fZ(k)=\frac{1}{k}\fza(k-1,1)\qquad(k\geq2)\,.
\end{align}
Then, for all $a,b\geq1$,
\begin{align}\label{eq:formalfinitedepthtwo}
\fza(a,b)=(-1)^a\binom{a+b}{a}\fZ(a+b)\,.
\end{align}
\end{proposition}

\begin{proof}
For $k=1$, the depth-one vanishing follows from reversal. In weight two,
the linear shuffle and stuffle products of $z_1$ with itself give
\begin{align*}
3\fza(1,1)=0,\qquad
0=\fza(1)^2=2\fza(1,1)+\fza(2)\,.
\end{align*}
Now suppose that the depth-one statement is known below weight $k$. The
linear shuffle relation for $z_{k-1}$ and $z_1$, together with reversal,
gives
\begin{align}\label{eq:formalfinitedepthoneinduction}
0=2\fza(k-1,1)+\sum_{i=2}^{k-2}\fza(k-i,i)
+2\fza(1,k-1)\,.
\end{align}
Pairing the terms in this sum and using the stuffle product of two
depth-one values shows that the right-hand side is
$-(k+1)\fza(k)/2$. This proves (i) by induction. Part~(ii) now follows
from the stuffle product and reversal:
\begin{align*}
0=\fza(a,b)+\fza(b,a)=2\fza(a,b)
\end{align*}
if $a+b$ is even.

For $k,d\geq1$, we will also use the formal sum formula
\begin{align}\label{eq:formalfinitecompositionsum}
\sum_{\substack{k_1+\dots+k_d=k\\k_1,\dots,k_d\geq1}}
\fza(k_1,\dots,k_d)=0\,.
\end{align}
To see this, notice that Hoffman's symmetric-sum identity writes the left-hand side as a
polynomial, for the stuffle product, in depth-one values (see \cite{H3}).
These values have just been shown to vanish.

For the last statement, the linear shuffle relation of $z_1$ with
$z_az_b$ gives
\begin{align*}
0={}&\fza(1,a,b)
+\sum_{u+v=a+1}\fza(u,v,b)
+\sum_{u+v=b+1}\fza(a,u,v)+\fza(a,b,1)\,,
\end{align*}
where $u,v\geq1$ in both sums. Expand in addition the three zero products
\begin{align*}
\fza(1)\fza(a,b),\qquad
\sum_{u+v=a+1}\fza(u)\fza(v,b),\qquad
\sum_{u+v=b+1}\fza(u)\fza(a,v)\,.
\end{align*}
After subtracting twice the preceding linear shuffle relation, the complete
composition sums vanish by~\eqref{eq:formalfinitecompositionsum}. The
remaining terms give
\begin{align}\label{eq:formalfinitedoubleremainder}
0={}&(a+1)\fza(a+1,b)+(b+1)\fza(a,b+1)\notag\\
&+\fza(a,1,b)-\fza(1,a,b)-\fza(a,b,1)\,.
\end{align}
If $a+b+1$ is even, the first line vanishes by part~(ii). If $a+b+1$ is
odd, put
\begin{align*}
T=\fza(b,a,1)+\fza(1,b,a)+\fza(a,1,b)\,.
\end{align*}
Expanding stuffle products with one depth-one factor and using reversal gives
\begin{align*}
0&=2T-\bigl(\fza(1,a,b)+\fza(a,b,1)+\fza(b,1,a)\bigr)\\
&=\bigl(2+(-1)^{a+b}\bigr)T\,.
\end{align*}
Thus $T=0$, and in odd total weight $T$ is precisely the second line of
\eqref{eq:formalfinitedoubleremainder}. So in both parities we
obtain the recurrence
\begin{align}\label{eq:formalfinitedoublerecurrence}
(a+1)\fza(a+1,b)+(b+1)\fza(a,b+1)=0\,.
\end{align}
Iterating~\eqref{eq:formalfinitedoublerecurrence} from
$\fza(k-1,1)=k\fZ(k)$ gives~\eqref{eq:formalfinitedepthtwo}.
\end{proof}

For example,
\begin{align}\label{eq:formalfinitedepthtwoexamples}
\fza(1,2)&=-3\fZ(3),&
\fza(2,1)&=3\fZ(3),\notag\\
\bigl(\fza(1,4),\fza(2,3),\fza(3,2),\fza(4,1)\bigr)
&=\bigl(-5,10,-10,5\bigr)\fZ(5)\,.
\end{align}
In particular, the relation
\begin{align}\label{eq:formalfiniteweightfiverelation}
2\fza(4,1)+\fza(3,2)=0
\end{align}
is the formal version of the relations in
\eqref{eq:za41z32} and the paragraph following
Theorem~\ref{thm:symspanmz}. Under $\pi_{\ma}$, the true single value
$\fZ(k)$ maps to the element $Z(k)$ in~\eqref{eq:defZ}.

\subsection{Formal symmetric values and the Kaneko--Zagier map}

The formal multiple zeta algebra in Section~\ref{sec:fmz} contains the
regularization parameter $\fzeta(1)$. For the construction of formal
symmetric multiple zeta values, it is important to set this parameter equal
to zero first. Put
\begin{align}\label{eq:formaladmissiblequotients}
\fmz_0&=\quotient{\fmz}{\fzeta(1)\fmz},&
\overline{\fmz}_0
&=\quotient{\fmz_0}{\fzeta(2)\fmz_0}
=\quotient{\fmz}{(\fzeta(1),\fzeta(2))}\,.
\end{align}
Since $\HH^1_\ast\simeq\HH^0_\ast[z_1]$, the first quotient is the
admissible, or $T=0$, formal multiple zeta algebra used in
\cite[Definition~1.19]{R}.
Namely, after setting $z_1$ equal to zero, the generators in
\eqref{eq:fmzedsideal} become the stuffle-regularized extended double
shuffle relations in $\HH^0$.

There are formal stuffle and shuffle regularizations
\begin{align*}
Z^{f,\bullet}:\HH^1_\bullet\longrightarrow\fmz_0[T]
\qquad(\bullet\in\{\ast,\sh\})
\end{align*}
which extend the quotient map on $\HH^0$ and send $z_1$ to $T$. We write
$\zeta^{f,\bullet}(k_1,\dots,k_r;T)$ for the image of
$z_{k_1}\cdots z_{k_r}$.

\begin{definition}[{\cite[Definition~4.13]{R}}]
For an index $\kk=(k_1,\dots,k_r)$, set
\begin{align}\label{eq:formalsymmetricdefinition}
\zeta_{\ms}^{f,\bullet}(\kk)
=\sum_{j=0}^r(-1)^{k_1+\dots+k_j}
\zeta^{f,\bullet}(k_j,\dots,k_1;T)
\zeta^{f,\bullet}(k_{j+1},\dots,k_r;T)\,.
\end{align}
The value attached to an empty index in this formula is $1$.
\end{definition}

\begin{proposition}[{\cite[Theorems~4.11--4.12 and Corollary~4.16]{R}}]
The expression~\eqref{eq:formalsymmetricdefinition} is independent of $T$,
and
\begin{align}\label{eq:formalsymmetricregularization}
\zeta_{\ms}^{f,\ast}(\kk)\equiv\zeta_{\ms}^{f,\sh}(\kk)
\pmod{\fzeta(2)\fmz_0}\,.
\end{align}
Their common class $\fzs(\kk)\in\overline{\fmz}_0$ satisfies the stuffle
product and the formal linear shuffle relations
\begin{align}\label{eq:formalsymmetricrelations}
\fzs(w\ast v)&=\fzs(w)\fzs(v),\notag\\
\fzs(w\sh v)&=(-1)^{\wt(w)}\fzs(\overline{w}v)\,.
\end{align}
\end{proposition}

\begin{proof}
The proof is the formal version of the argument in Chapter~\ref{sec:fmzv}.
For an admissible index $\kk$, define
\begin{align*}
G_f^\bullet(\kk;X,T)
=\sum_{j\geq0}\zeta^{f,\bullet}(\{1\}^j,\kk;T)X^j\,.
\end{align*}
Formal regularization gives
\begin{align*}
G_f^\bullet(\kk;X,T)&=e^{XT}G_f^\bullet(\kk;X,0),\notag\\
G_f^\sh(\kk;X,0)&=A^f(X)G_f^\ast(\kk;X,0)\,,
\end{align*}
where
\begin{align*}
A^f(X)=\exp\left(\sum_{n\geq2}\frac{(-1)^n}{n}\fzeta(n)X^n\right)\,.
\end{align*}
Group the cuts in~\eqref{eq:formalsymmetricdefinition} according to the block
of consecutive $1$'s through which the cut passes. Every such block is the
coefficient of a product
\begin{align*}
(-1)^{\wt(\kl)}G_f^\bullet(\kl;-X,T)
G_f^\bullet({\bf m};X,T)
\end{align*}
for two admissible indices $\kl$ and ${\bf m}$. The factors $e^{-XT}$ and
$e^{XT}$ cancel, which proves the independence of $T$. We also have
\begin{align*}
A^f(-X)A^f(X)
=\exp\left(\sum_{m\geq1}\frac{\fzeta(2m)}{m}X^{2m}\right)
\equiv1\pmod{\fzeta(2)\fmz_0[[X]]}
\end{align*}
by the formal Euler relation~\eqref{eq:fmzeuler}. This proves
\eqref{eq:formalsymmetricregularization}.

For the stuffle product, write the alternating cut in
\eqref{eq:formalsymmetricdefinition} as the convolution of the two stuffle
homomorphisms obtained before and after the cut. The deconcatenation
coproduct is compatible with the stuffle product, so their convolution is
again an algebra homomorphism. For the linear shuffle relation, use the
shuffle representative and the purely combinatorial cut--shuffle identity
\eqref{eq:symmetricshuffleoperator}. This gives~\eqref{eq:formalsymmetricrelations}.
\end{proof}

In depth two the definition can be evaluated explicitly.

\begin{proposition}[{\cite[Proposition~4.14]{R}}]\label{prop:formalsymmetricdepthtwo}
If $a,b\geq1$ and $k=a+b$, then
\begin{align}\label{eq:formalsymmetricdepthtwo}
\fzs(a,b)=
\begin{cases}
0,&k\text{ even},\\
(-1)^a\binom{k}{a}\fzeta(k),&k\text{ odd}
\end{cases}
\qquad\text{in }\overline{\fmz}_0\,.
\end{align}
\end{proposition}

\begin{proof}
Set $T=0$ in~\eqref{eq:formalsymmetricdefinition}. If $k$ is even, the
stuffle relation gives
\begin{align*}
\fzs(a,b)=\bigl(1+(-1)^a\bigr)\fzeta(a)\fzeta(b)-\fzeta(k)\,,
\end{align*}
which is divisible by $\fzeta(2)$ by~\eqref{eq:fmzeuler}. For odd $k$, the
depth-two extended double shuffle relation gives
\begin{align*}
\zeta^{f,\ast}(a,b;0)
\equiv\frac12\left((-1)^a\binom{k}{a}-1\right)\fzeta(k)
\pmod{\fzeta(2)\fmz_0}\,.
\end{align*}
Substituting this formula and its version with $a,b$ interchanged into
\eqref{eq:formalsymmetricdefinition} gives the result. The same calculation,
including the boundary case $a=1$, is given in
\cite[proof of Proposition~4.14]{R}.
\end{proof}

Now the formal analogue of the Kaneko--Zagier map can be stated.

\begin{theorem}[Formal Kaneko--Zagier map, {\cite[Theorem~4.18 and Corollary~4.19]{R}}]
\label{thm:formalkanekozagier}
The assignment
\begin{align}\label{eq:formalkanekozagier}
\varphi_{\mathrm{KZ}}^f:\fmza&\longrightarrow\overline{\fmz}_0,\notag\\
\fza(\kk)&\longmapsto\fzs(\kk)
\end{align}
defines a surjective homomorphism of graded $\Q$-algebras. Moreover,
for $k\geq2$,
\begin{align}\label{eq:formaltruesingleimage}
\varphi_{\mathrm{KZ}}^f\bigl(\fZ(k)\bigr)
=\fzeta(k)\qquad\text{in }\overline{\fmz}_0\,.
\end{align}
\end{theorem}

\begin{proof}
The relations~\eqref{eq:formalsymmetricrelations} show that the map from
$\HH^1_\ast$ to $\overline{\fmz}_0$ annihilates the ideal $\mathfrak L$,
so it factors through $\fmza$. The non-trivial point is surjectivity.
Yasuda's proof that symmetric multiple zeta values span multiple zeta values
is an induction in weight and depth whose polynomial step uses only the
linearized extended double shuffle relations. It can therefore be carried
out in the formal algebra. This shows that the formal symmetric values span
$\overline{\fmz}_0$. See \cite[Proposition~4.17]{R} for the
formal polynomial argument. Finally,
Propositions~\ref{prop:formalfinitedepthtwo} and
\ref{prop:formalsymmetricdepthtwo} give
\begin{align*}
\varphi_{\mathrm{KZ}}^f\bigl(\fZ(k)\bigr)
=\frac{1}{k}\fzs(k-1,1)=\fzeta(k)
\end{align*}
for odd $k$. For even $k$, both sides vanish in
$\overline{\fmz}_0$ by~\eqref{eq:fmzeuler}.
\end{proof}

\begin{conjecture}[{\cite[Chapter~5(i)]{R}}]\label{conj:formalkanekozagier}
The homomorphism $\varphi_{\mathrm{KZ}}^f$ is an isomorphism.
\end{conjecture}

This is a formal version of the Kaneko--Zagier conjecture, but it does not
yet give a map $\mza\rightarrow\mz/\pi^2\mz$. The problem is that the realization
$\fmza\rightarrow\mza$ might have additional relations, and it is not known
that all of them vanish under the symmetric realization. This is exactly the
difficulty in the classical conjecture.

\subsection{Parity in depth at most four}

For $d,k\geq0$, define
\begin{align}\label{eq:formalfinitedepthfiltration}
\fild_d(\fmza_k)
=\left\langle\fza(k_1,\dots,k_r)\mid
r\leq d,\ k_1+\dots+k_r=k\right\rangle_\Q\,.
\end{align}
The weak parity theorem from Section~\ref{sec:parity} holds already in the
formal algebra.

\begin{theorem}[{\cite[Theorem~4.22]{R}}]\label{thm:formalfiniteweakparity}
If $k,d\geq1$ and $k\equiv d\pmod2$, then
\begin{align}\label{eq:formalfiniteweakparity}
\fild_d(\fmza_k)\subseteq{}&\fild_{d-1}(\fmza_k)\notag\\
&+\sum_{\substack{k_1+k_2=k,\ d_1+d_2=d\\d_1,d_2\geq2}}
\fild_{d_1}(\fmza_{k_1})\fild_{d_2}(\fmza_{k_2})\,.
\end{align}
\end{theorem}

\begin{proof}
Apply the stuffle antipode relation~\eqref{eq:antipode-relation} to a word
$z_{a_1}\cdots z_{a_d}$ and then pass to $\fmza$. Modulo depth $d-1$, all
contractions disappear, and the two end terms add up to
\begin{align*}
\bigl(1+(-1)^{k+d}\bigr)\fza(a_1,\dots,a_d)
\end{align*}
by~\eqref{eq:formalfinitereversal}. The coefficient is $2$ if
$k\equiv d\pmod2$. Every interior term is a product of two smaller-depth
values, and the terms having a depth-one factor vanish by
Proposition~\ref{prop:formalfinitedepthtwo}. This gives exactly the sum in
\eqref{eq:formalfiniteweakparity}.
\end{proof}

The proof also gives useful explicit formulas. If $a+b+c$ is odd, then
\begin{align}\label{eq:formalfiniteweakdepththree}
\fza(a,b,c)
=-\frac12\bigl(\fza(a+b,c)+\fza(a,b+c)\bigr)\,.
\end{align}
If $a+b+c+d$ is even, then
\begin{align}\label{eq:formalfiniteweakdepthfour}
\fza(a,b,c,d)=-\frac12\bigl(&\fza(a+b,c,d)+\fza(a,b+c,d)
+\fza(a,b,c+d)\notag\\
&+(-1)^{a+b}\fza(a,b)\fza(c,d)\bigr)\,.
\end{align}
These are obtained by expanding the antipode relation in depths three and
four. Notice in particular the sign $(-1)^{a+b}$ in the product term. It
comes from applying reversal to $\fza(b,a)$.

It is natural to expect that the product term in
Theorem~\ref{thm:formalfiniteweakparity} can always be removed. This is the
formal version of the strong parity conjecture stated after
Theorem~\ref{thm:fmzvparity}. In the formal algebra this is known through
depth four.

\begin{theorem}[{\cite[Theorem~4.30]{R}}]\label{thm:formalfiniteparityfour}
Let $k,d\geq1$, $k\equiv d\pmod2$ and $d\leq4$. Then
\begin{align}\label{eq:formalfiniteparityfour}
\fild_d(\fmza_k)\subseteq\fild_{d-1}(\fmza_k)\,.
\end{align}
\end{theorem}

\begin{proof}
Only depth four needs to be considered. Fix an even weight $k$ and put
\begin{align*}
P_{r,s}=\fZ(r)\fZ(s)=P_{s,r}
\end{align*}
for odd $r,s\geq3$ with $r+s=k$. If $r=a+b$ and $s=c+d$, then
\eqref{eq:formalfiniteweakdepthfour} and
\eqref{eq:formalfinitedepthtwo} give
\begin{align}\label{eq:formalfiniteparitydictionary}
\fza(a,b,c,d)\equiv
\begin{cases}
\dfrac12(-1)^{a+c}\binom{r}{a}\binom{s}{c}P_{r,s},&r\text{ odd},\\
0,&r\text{ even}
\end{cases}
\pmod{\fild_3(\fmza_k)}\,.
\end{align}
It remains to lower the products $P_{r,s}$.

For odd $r\geq3$, set
\begin{align*}
f_r(X,Y)=\frac{(Y-X)^r-Y^r+X^r}{XY}\,.
\end{align*}
The generating series of the depth-four values satisfies the five-term
relation
\begin{align}\label{eq:formalfinitefiveterm}
0={}&\mathfrak F_k(X,Y,Z,W)+\mathfrak F_k(Y-X,Y,Z,W)\notag\\
&+\mathfrak F_k(Y-X,Z-X,Z,W)
+\mathfrak F_k(Y-X,Z-X,W-X,W)\notag\\
&+\mathfrak F_k(Y-X,Z-X,W-X,-X)\,,
\end{align}
where
\begin{align*}
\mathfrak F_k(X,Y,Z,W)
=\sum_{a+b+c+d=k}\fza(a,b,c,d)
X^{a-1}Y^{b-1}Z^{c-1}W^{d-1}\,.
\end{align*}
This is the generating-series form of the linear shuffle relation between a
word of depth one and a word of depth three. Insert
\eqref{eq:formalfiniteparitydictionary} into
\eqref{eq:formalfinitefiveterm} and set $X=W=0$. Comparing the coefficient
of $Y^{\alpha-1}Z^{k-\alpha-3}$ gives
\begin{align}\label{eq:formalfiniteparitytriangular}
\sum_{\substack{r+s=k,\ r,s\geq3\text{ odd}\\\alpha\leq r-2}}
s\binom{r}{\alpha}P_{r,s}\equiv0
\pmod{\fild_3(\fmza_k)}\,.
\end{align}
Using $P_{r,s}=P_{s,r}$, choose successively $\alpha=r_0-2$ for
$r_0=k-3,k-5,\dots$ down to the smallest odd $r_0\geq k/2$. Pairs with
$\max(r,s)<r_0$ do not occur in~\eqref{eq:formalfiniteparitytriangular},
whereas pairs with $\max(r,s)>r_0$ have already been treated. The remaining
coefficient of $P_{r_0,k-r_0}$ is non-zero. Hence every $P_{r,s}$ has depth
at most three. Formula~\eqref{eq:formalfiniteparitydictionary} now proves
the theorem. Depths one and two follow from
Proposition~\ref{prop:formalfinitedepthtwo}, and depth three follows from
\eqref{eq:formalfiniteweakdepththree}.
\end{proof}

Risan proved Theorem~\ref{thm:formalfiniteparityfour} in his master's thesis
by a generating-series calculation. The proof above is a shorter uniform
version of this argument. The thesis also gives a closed reduction
formula~\cite[Corollary~4.29]{R}. For example, it gives the identity
\begin{align}\label{eq:formalfinite2121}
\fza(2,1,2,1)=\frac32\fza(1,2,3)\,.
\end{align}
For the general formula we refer to the thesis.

\subsection{The formal double zeta space and modular forms}

The depth-four part of $\fmza$ is closely related to the formal double zeta space from \Cref{sec:mdandmzv}. So the cusp forms, which we
met first in the world of real multiple zeta values, also show up in the finite
world. This is quite mysterious, since the two worlds have nothing to do with
each other at first sight. For even $k\geq4$, define
its reduced quotient by
\begin{align}\label{eq:reducedformaldoublezeta}
\widetilde{\dz}_k
=\quotient{\dz_k}{
\Q Z_k+\Q Z_{1,k-1}
+\left\langle P_{r,s}\mid r,s\text{ even}\right\rangle_\Q}\,.
\end{align}
Here $Z_k$ and the even product symbols model the part divisible by
$\zeta(2)$, while $Z_{1,k-1}$ removes the regularized boundary case.

\begin{proposition}[{\cite[Proposition~4.33(i)]{R}}]\label[proposition]{prop:reducedformaldoublezetadim}
For even $k\geq4$, we have
\begin{align}\label{eq:reducedformaldoublezetadim}
\dim_\Q\widetilde{\dz}_k=\frac{k}{2}-2-\dim_\C\ms_k\,.
\end{align}
\end{proposition}

\begin{proof}
In even weight, the odd--odd symbols form a basis of $\dz_k$ by
Theorem~\ref{thm:dzevenbasis}. The subspace
\begin{align*}
\mathcal P_k^{\mathrm{ev}}
=\Q Z_k+\left\langle P_{r,s}\mid r,s\text{ even}\right\rangle_\Q
\end{align*}
has dimension $\dim_\C\ms_k+1$ by Theorems~\ref{thm:eichlershimura} and
\ref{thm:wkminuspkev}. Notice that $Z_{1,k-1}$ does not belong to this
subspace: in the odd--odd basis its coefficient is zero on
$\mathcal P_k^{\mathrm{ev}}$, as follows from the period-polynomial formula
in Theorem~\ref{thm:wkminuspkev}. Quotienting by its class therefore removes
one further dimension. This leaves
$k/2-(\dim_\C\ms_k+1)-1$ dimensions.
\end{proof}

We will now describe the proved connection with finite multiple zeta values.
For odd $r,s\geq3$ with $r+s=k$, define the following element of
$\widetilde{\dz}_k$:
\begin{align}\label{eq:formalfiniteHrs}
H_{r,s}=\frac12\Bigg(&
(rs-r+1)Z_{r,s}+(rs+s-1)Z_{s,r}
+(s+1)Z_{s+1,r-1}-(r+1)Z_{r+1,s-1}\notag\\
&+\sum_{\substack{3\leq\ell\leq k-3\\\ell\text{ odd}}}
\left(\binom{k-\ell}{s}-\binom{k-\ell}{r}\right)
P_{\ell,k-\ell}\Bigg)\,.
\end{align}

\begin{theorem}[Bachmann--Risan {\cite{BaR}}]
\label{thm:finitebetasurjection}
For every even $k\geq4$, the elements $H_{r,s}$ with $r,s\geq3$ odd and
$r+s=k$ span $\widetilde{\dz}_k$. The rule
\begin{align}\label{eq:finitebetasurjection}
\beta_k:\widetilde{\dz}_k&\longrightarrow\fild_4(\mza_k),\notag\\
H_{r,s}&\longmapsto\za(r-1,1,s-1,1)
\end{align}
defines a surjective $\Q$-linear map. We also have
\begin{align}\label{eq:finitebetaspan}
\fild_4(\mza_k)
&=\left\langle\za(2a,1,2b,1)\mid
a,b\geq1,\ 2(a+b+1)=k\right\rangle_\Q\notag\\
&=\left\langle\za(1,2a,1,2b)\mid
a,b\geq1,\ 2(a+b+1)=k\right\rangle_\Q\,.
\end{align}
The dimension satisfies
\begin{align}\label{eq:finitebetadimension}
\dim_\Q\fild_4(\mza_k)
\leq\frac{k}{2}-2-\dim_\C\ms_k\,.
\end{align}
In particular, the generators in either line of~\eqref{eq:finitebetaspan}
satisfy at least $\dim_\C\ms_k$ linearly independent relations.
\end{theorem}

\begin{proof}
We first replace $\za$ by $\fza$ and $Z(j)$ by $\fZ(j)$ in the construction of \Cref{sec:fmzvmodularforms}. Notice that the proofs of \cref{lem:finitehoffmanlecture,lem:finitetriplereductionlecture,lem:finitedspolynomiallecture} only use the stuffle product, the linear shuffle relations, the reversal formula and the formulas in low depth, and all of these hold in $\fmza$. Therefore the formal polynomial $f=\sum_{r+s=k}f_{r,s}X^{r-1}Y^{s-1}$ satisfies the same two double shuffle equations, i.e. sending $Z_{r,s}$ to $f_{r,s}$ and $Z_k$ to zero defines a realization of $\dz_k$ as in \Cref{sec:fmzvmodularforms}. This realization kills $Z_{1,k-1}$, since $f(0,Y)=0$, and it kills the even product symbols, since $B=f+f|\epsilon$ has zero coefficients for even $r,s$. So it factors through $\widetilde{\dz}_k$.

The binomial theorem in \eqref{eq:finitepolynomialslecture} together with the reduction formula gives
\begin{align*}
 f_{r,s}=(-1)^r\sum_{n=1}^{k-2}\left(\binom nr-\frac12\binom{k-1}{r}\right)
 \bigl(\fza(n,1,k-1-n)-\fZ(n)\fZ(k-n)\bigr)-\delta_{r,1}\mu,
\end{align*}
where $\mu=\frac12\sum_{n=1}^{k-1}\fZ(n)\fZ(k-n)$. Substituting this into \eqref{eq:formalfiniteHrs}, we get
\begin{align*}
 H_{r,s}\longmapsto\frac12\bigl(rs\fZ(r)\fZ(s)+\fza(1,r-1,s)-\fza(1,s-1,r)\bigr)
 =\fza(r-1,1,s-1,1),
\end{align*}
by the formal version of \eqref{eq:finitefourbridgelecture}, and this constructs $\beta_k^{\mathrm f}$. Conversely, the identity $XYM(X,Y)=(Y-X)f(Y-X,Y)+XB(X,Y)$ recovers every value $\fza(n,1,k-1-n)$ from its image, and since $B=f+f|\epsilon$, all product coefficients are in the image as well. So by the formal triple reduction the image is $\fild_3(\fmza_k)$, which equals $\fild_4(\fmza_k)$ by \cref{thm:formalfiniteparityfour}.

By the double shuffle calculation in \cite{BaR}, the elements $H_{r,s}$ span $\widetilde{\dz}_k$, where the antisymmetric part uses the triangular basis \eqref{eq:finiteEblecture} and the symmetric part uses $H_{r,s}+H_{s,r}=rsP_{r,s}$. After the finite realization, this gives $\beta_k$ and the first equality in \eqref{eq:finitebetaspan}, and the second one follows from the reversal formula. The dimension bound follows from \cref{prop:reducedformaldoublezetadim}.
\end{proof}

\input{chap_FormalFiniteComparison}

So only the injectivity after the realization in actual finite multiple zeta values remains conjectural.
\begin{conjecture}[Bachmann--Risan {\cite{BaR}}]\label{conj:formalfinitedoublezeta}
For every even $k\geq4$, the finite map
\begin{align}\label{eq:formalfinitedoublezeta}
 \beta_k:\widetilde{\dz}_k\longrightarrow\fild_4(\mza_k)
\end{align}
is an isomorphism of $\Q$-vector spaces.
\end{conjecture}
This is equivalent to equality in \eqref{eq:finitebetadimension}, i.e. \cref{conj:fmzvmodularforms}\,(i).

We next want to see explicitly how the period polynomials from
Chapter~\ref{sec:mdandmzv} give relations among finite multiple zeta values.
For $p\in W_k^-$, write
\begin{align}\label{eq:formalfiniteperiodcoefficients}
p(X+Y,Y)=\sum_{r+s=k}\binom{k-2}{r-1}q_{r,s}
X^{r-1}Y^{s-1}\,.
\end{align}

\begin{proposition}[{\cite{BaR}}]
\label[proposition]{prop:formalfiniteperiodrelation}
For every even $k\geq4$ and $p\in W_k^-$, we have
\begin{align}\label{eq:formalfiniteperiodrelation}
\sum_{\substack{r+s=k\\r,s\geq3\text{ odd}}}
\bigl(q_{r,s}-q_{k-1,1}\bigr)Z_{r,s}=0
\qquad\text{in }\widetilde{\dz}_k\,.
\end{align}
Every relation among the classes $Z_{r,s}$ with $r,s\geq3$ odd is obtained
in this way. The kernel of the assignment from $W_k^-$ to these relations is
$\Q(X^{k-2}-Y^{k-2})$.
\end{proposition}

\begin{proof}
The coefficients $q_{r,s}$ are the coefficients denoted by
$\beta^p_{r,s}$ in Theorem~\ref{thm:wkminuspkev}. That theorem shows that
the sum of $q_{r,s}Z_{r,s}$ over all odd $r,s$ belongs to
$\mathcal P_k^{\mathrm{ev}}$, so it vanishes in
$\widetilde{\dz}_k$. There we also have $Z_{1,k-1}=0$, and the odd sum
formula in Theorem~\ref{thm:dzsumformula} gives
\begin{align*}
Z_{k-1,1}=-\sum_{\substack{3\leq r\leq k-3\\r\text{ odd}}}Z_{r,k-r}\,.
\end{align*}
Substituting this gives~\eqref{eq:formalfiniteperiodrelation}. The final two
statements follow from the converse part of the period-polynomial theorem and
the dimension formulas in Theorem~\ref{thm:eichlershimura} and
Proposition~\ref{prop:reducedformaldoublezetadim}.
\end{proof}

\input{chap_FormalFinitePeriodTransfer}

To write the corresponding finite relation down, consider the slightly larger
quotient
\begin{align*}
\overline{\dz}_k
=\quotient{\dz_k}{\Q Z_k+\Q Z_{1,k-1}}\,.
\end{align*}
Theorem~\ref{thm:dzevenbasis} and the odd sum formula show that the classes
$Z_{r,k-r}$ with $3\leq r\leq k-3$ odd form a basis of
$\overline{\dz}_k$. Read the right-hand side of~\eqref{eq:formalfiniteHrs}
in this quotient and define the square matrix $A_k=(A_{r,r'})$, whose rows
and columns are indexed by the odd integers from $3$ to $k-3$, by
\begin{align}\label{eq:formalfiniteAmatrix}
2H_{r,k-r}
=\sum_{\substack{3\leq r'\leq k-3\\r'\text{ odd}}}
A_{r,r'}Z_{r',k-r'}\,.
\end{align}
Its entries are explicit, since the even--even symbols occurring in
$H_{r,k-r}$ are given in terms of the odd--odd basis by
\eqref{eq:pevevzevevexplicitformula}. For $p\in W_k^-$ and $q_{r,s}$ as
in~\eqref{eq:formalfiniteperiodcoefficients}, let $\widehat c_r(p)$ be the
determinant of the matrix obtained from $A_k$ by replacing its row indexed by
$r$ with
\begin{align*}
\bigl(2(q_{r',k-r'}-q_{k-1,1})\bigr)_{r'}\,.
\end{align*}

\begin{theorem}[Bachmann--Risan {\cite{BaR}}]
\label{thm:formalfiniteperiodtransfer}
For every even $k\geq6$ and every $p\in W_k^-$, we have
\begin{align}\label{eq:formalfiniteperiodtransfer}
\sum_{\substack{3\leq r\leq k-3\\r\text{ odd}}}
\widehat c_r(p)\za(r-1,1,k-r-1,1)=0\,.
\end{align}
If $\det A_k\neq0$, these relations, as $p$ varies, give all relations among
the elements $H_{r,k-r}$, and their coefficient vectors vanish exactly on
$\Q(X^{k-2}-Y^{k-2})$.
\end{theorem}

\begin{proof}
Put $b_{r'}=q_{r',k-r'}-q_{k-1,1}$. Expansion along the replaced row and
the adjugate identity give
\begin{align*}
\sum_r\widehat c_r(p)A_{r,r'}=2\det(A_k)b_{r'}\,.
\end{align*}
Combining this with~\eqref{eq:formalfiniteAmatrix} and
Proposition~\ref{prop:formalfiniteperiodrelation} gives
\begin{align*}
\sum_r\widehat c_r(p)H_{r,k-r}=0
\qquad\text{in }\widetilde{\dz}_k\,.
\end{align*}
Now apply the map $\beta_k$ from
Theorem~\ref{thm:finitebetasurjection}. If $A_k$ is invertible, Cramer's rule
and Proposition~\ref{prop:formalfiniteperiodrelation} give the last statement.
\end{proof}

It is conjectured that $\det A_k\neq0$ for every even $k\geq6$, or equivalently that the $H_{r,s}$ form a basis of $\overline{\dz}_k$. The public version of \cite{BaR} records checks through weight $100$, and later exact computations reach every even weight up to $500$, but of course these are just finite checks. On the other hand, the cofactor formula holds unconditionally, but it might produce zero coefficients if the determinant vanishes. In contrast, \cref{thm:completeformalfiniteperiodtransfer} already gives all formal relations without this conjecture. The remaining inverse problem and its full $Z_k$ correction are described in \cref{ex:formalfiniteinversionlecture}.

The first cusp-form relation occurs in weight twelve. In this case it even
holds in the universal formal algebra, before either classical realization.

\begin{samepage}
\begin{theorem}[{\cite[Proposition~4.35(i)]{R}}]\label{thm:formalfiniteweighttwelve}
In $\fmza_{12}$, we have
\begin{align}\label{eq:formalfiniteweighttwelve}
2\fza(1,2,1,8)-18\fza(1,4,1,6)
-9\fza(1,6,1,4)-16\fza(1,8,1,2)=0\,.
\end{align}
\end{theorem}

\begin{proof}
This follows by combining the stuffle product and the linear shuffle relations
in weight twelve. Only words of depth at most four are needed.
\end{proof}
\end{samepage}

The finite realization of this identity is the first relation displayed in
Section~\ref{subsec:modularfmzv}. After reversal and rescaling, it is also the
case of
Theorem~\ref{thm:formalfiniteperiodtransfer} obtained from the period
polynomial of $\Delta$ (see \cite{BaR}). On the formal double
zeta side, this polynomial gives
\begin{align*}
28Z_{9,3}+150Z_{7,5}+168Z_{5,7}=0
\qquad\text{in }\widetilde{\dz}_{12}\,.
\end{align*}

Therefore the formal depth-three reduction and the formal isomorphism are proved in every even weight. Notice that the determinant question concerns a larger quotient, before the even product symbols are killed, so it is a separate question. The injectivity for the actual finite values, the two triple-spanning conjectures and the injectivity of the formal Kaneko--Zagier map in general are still open.

%% file: chap_FormalFiniteComparison.tex
We now want to show that the formal map $\beta_k^{\mathrm f}:\widetilde{\dz}_k\twoheadrightarrow\fild_4(\fmza_k)$ from the preceding proof is also injective. For this, we use motivic multiple zeta values. We write $\overline{\mz}^{\mathfrak m}=\mz^{\mathfrak m}/(\zeta^{\mathfrak m}(2))$ and denote by $\fild_r\overline{\mz}^{\mathfrak m}$ the usual motivic depth filtration. If we compose the formal Kaneko--Zagier map $\varphi_{\mathrm{KZ}}^f$ from \cref{thm:formalkanekozagier} with the motivic realization of formal multiple zeta values, we obtain a map
\begin{align*}
 \varphi_{\mathrm{KZ}}^{\mathfrak m}:\fmza\longrightarrow\overline{\mz}^{\mathfrak m}\,,
\end{align*}
and in particular we have $\varphi_{\mathrm{KZ}}^{\mathfrak m}(\fZ(2j+1))=\zeta^{\mathfrak m}(2j+1)$.

\begin{lemma}\label[lemma]{lem:formalfinitemotivicimage}
For every even $k\geq4$, we have $\varphi_{\mathrm{KZ}}^{\mathfrak m}(\fild_4(\fmza_k))=\fild_2\overline{\mz}^{\mathfrak m}_k$.
\end{lemma}
\begin{proof}
Let $\Pi_k$ be the span of the products of two odd motivic single zeta values of weight at least three and total weight $k$. These products are in the image. By the depth-three calculation of Yasuda \cite[Proposition~3.1 and Remark~3.2]{Yas}, the image of $\fza(c,b,a)$ modulo $\Pi_k$ is given by
\begin{align*}
 &-(-1)^a\sum_{i+j=a}\binom{b+i-1}{i}\binom{c+j-1}{j}\zeta^{\mathfrak m}(c+j,b+i)\\
 &\qquad+(-1)^c\sum_{i+j=c}\binom{a+i-1}{i}\binom{b+j-1}{j}\zeta^{\mathfrak m}(b+j,a+i)\,,
\end{align*}
where $i,j\geq0$ and a first entry $1$ is interpreted by regularization at $T=0$. Since \cref{thm:formalfiniteparityfour} gives $\fild_4(\fmza_k)=\fild_3(\fmza_k)$, the image is contained in $\fild_2$.

For the reverse inclusion, let $\lambda$ be a rational functional on $\fild_2\overline{\mz}^{\mathfrak m}_k$ which annihilates these expressions and $\Pi_k$, and put $F(X,Y)=\sum_{r+s=k}\lambda(\zeta^{\mathfrak m}(s,r))X^{r-1}Y^{s-1}$. The double shuffle relations give $F(Y,X)=-F(X,Y)$ and $F(X,X+Y)+F(Y,X+Y)=0$, and by the binomial theorem the annihilation of the displayed expressions becomes
\begin{align*}
 \frac{F(Y-X,Z-X)-F(Y,Z)}X
 =\frac{F(X-Z,Y-Z)-F(X,Y)}Z\,,
\end{align*}
where both quotients are polynomials. For $Z=0$, the two double shuffle equations give $F+X(\partial_X+\partial_Y)F=F(Y,0)$, and evaluating at $X=0$ gives $F(Y,0)=0$. Interchanging $X$ and $Y$ then gives $(X-Y)(\partial_X+\partial_Y)F=0$, so the preceding identity forces $F=0$. Thus the functional is zero, which proves the surjectivity onto $\fild_2$.
\end{proof}

\begin{theorem}[Bachmann--Risan {\cite{BaR}}]\label{thm:formalfinitebetaisomorphism}
For every even $k\geq4$, the map $\beta_k^{\mathrm f}$ is an isomorphism.
\end{theorem}
\begin{proof}
By Li's depth-two dimension theorem \cite[Theorem~1.6(i)]{JLi} we have
\begin{align*}
 \dim_\Q\fild_2\overline{\mz}^{\mathfrak m}_k=\frac{k}{2}-2-\dim_\C\ms_k\,,
\end{align*}
which by \cref{prop:reducedformaldoublezetadim} is also the dimension of $\widetilde{\dz}_k$. Therefore the preceding lemma and the surjectivity of $\beta_k^{\mathrm f}$ force both $\beta_k^{\mathrm f}$ and the restriction of $\varphi_{\mathrm{KZ}}^{\mathfrak m}$ to be isomorphisms.
\end{proof}
Notice that this theorem uses the revised version of \cite{BaR}.\footnote{The public August 2026 version proves the map to actual finite values, but leaves the formal comparison conjectural. The formal reduction, isomorphism and complete period transfer used here are from the later revision.} In particular, the formal Kaneko--Zagier map is injective on the part of even weight and depth at most four. Unfortunately, this does not give a lower bound for the actual finite values, since the motivic lower bound only applies to the formal algebra and no factorization of $\varphi_{\mathrm{KZ}}^{\mathfrak m}$ through $\mza$ is known.

\begin{corollary}\label[corollary]{cor:formaltruesinglesindependent}
The elements $\fZ(3),\fZ(5),\fZ(7),\ldots$ are algebraically independent, and the formal Kaneko--Zagier map is injective on the algebra generated by them.
\end{corollary}
\begin{proof}
Their images under $\varphi_{\mathrm{KZ}}^{\mathfrak m}$ are the odd single zeta values, which can be chosen as the letters of odd weight in Brown's shuffle description of motivic multiple zeta values \cite{Br2}. Shuffle monomials in these letters are linearly independent, since different multiplicities have disjoint word support. Therefore every polynomial relation among the formal elements is mapped to a polynomial relation among algebraically independent motivic elements, i.e. it is trivial. Since $\varphi_{\mathrm{KZ}}^{\mathfrak m}$ factors through the formal Kaneko--Zagier map, the same argument shows that the latter is injective on this subalgebra.
\end{proof}
The injectivity of the formal Kaneko--Zagier map in general and the strong parity beyond depth four are still open.

%% file: chap_FormalFinitePeriodTransfer.tex
\subsubsection{A complete period transfer}\label{sec:completeformalfiniteperiodtransfer}
The construction in \Cref{sec:fmzvmodularforms} gives relations among the actual finite values. We now want to give a construction which produces all relations among the formal values $\fza(2a,1,2b,1)$ and which does not need any assumption on a determinant. It consists of a triangular correction, followed by the even--odd conversion from \Cref{sec:mdandmzv}.

For every even $k\geq4$ we put $I=\{r:3\leq r\leq k-3,\ r\text{ odd}\}$, and for $k\geq8$ we also put $m=k-2$, $d=\lfloor m/6\rfloor$ and $I^-=\{r\in I:2r<k\}$. Let $p\in W_k^-$ with $p=\sum_{r+s=k}\binom m{r-1}p_rX^{r-1}Y^{s-1}$. The coefficients of the relation coming from $p$ are defined in the following three steps.
\begin{enumerate}[(i)]
\item For $r\in I^-$, start with $c_r=p_{r+2}$ if $2d+1\leq r<k/2-2$ and $c_r=0$ otherwise. Let $\omega_j$ and $\widetilde\omega_j$ be the polynomials from the proof of \cref{thm:finiteHperiodlecture}, i.e. $\omega_j(t)=\Omega_j(t,1)$ with $\Omega_j$ as in \eqref{eq:finiteEblecture} and $\widetilde\omega_j(t)=((m+1)t^2+mt+1)\omega_j(t)+(t-t^3)\omega_j'(t)$. For $j=d-1,d-2,\ldots,1$, successively replace
\begin{align*}
 c_{2j+1}\longleftarrow c_{2j+1}+\frac1{4(j+1)}\sum_{r\in I^-}c_r[t^r]\widetilde\omega_j(t)\,.
\end{align*}
Then replace all coefficients simultaneously by
\begin{align*}
 c_r\longleftarrow c_r+\frac{2r-k}{6m}\sum_{u\in I^-}c_u[t^u]\widetilde\omega_0(t)\,.
\end{align*}
Finally, extend by $c_{k-r}=-c_r$ and set $c_{k/2}=0$ if $k/2\in I$.
\item Form
\begin{align*}
 h(X,Y)&=\sum_{r\in I}\binom kr c_rX^rY^{k-r},\\
 F(X,Y)&=m![t^m]\left(e^{Yt}\frac{Xt}{e^{Xt}-1}\sum_{r\in I}\frac{r+1}{r!}c_r(Xt)^{r-1}\right)\,.
\end{align*}
In the following expression, we keep only the monomials $X^aY^{m-a}$ with $a$ even and $2\leq a\leq m-2$, and call the resulting polynomial $R(X,Y)$:
\begin{align}\label{eq:formalrecursiveRlecture}
 &\frac{\partial_X\partial_Y}{k(k-1)}\bigl(h(X+Y,Y)+h(X+Y,X)-X\partial_Xh(X,Y)+h(X,Y)\bigr)\notag\\
 &\qquad+2\bigl(F(X+Y,X)-F(Y,X)-F(1,1)(X+Y)^m\bigr)\,.
\end{align}
\item Write
\begin{align*}
 R(X+Y,Y)-R(Y,X+Y)-R(X,Y)=\sum_{r+s=k}\binom m{r-1}b_{r,s}X^{r-1}Y^{s-1}
\end{align*}
and put $C_r(p)=c_r+b_{r,k-r}/(r(k-r))$ for $r\in I$.
\end{enumerate}
For $k=4,6$, we set all $C_r(p)=0$.

\begin{theorem}[Bachmann--Risan {\cite{BaR}}]\label{thm:completeformalfiniteperiodtransfer}
For every even $k\geq4$ and $p\in W_k^-$, we have
\begin{align}\label{eq:completeformalfiniteperiodtransfer}
 \sum_{r\in I}C_r(p)\fza(r-1,1,k-r-1,1)=0\,.
\end{align}
The map $p\mapsto(C_r(p))_r$ identifies $W_k^-/\Q(X^{k-2}-Y^{k-2})$ with the space of all relations among these formal values.
\end{theorem}
\begin{proof}
We take the quotient of $\widetilde{\dz}_k$ by the odd product symbols. Its dual has the basis $\Omega_j$ from \eqref{eq:finiteEblecture}, where the functional corresponding to $\Omega_j$ sends $Z_{r,k-r}$ to $[t^{r-1}]\omega_j$ and $H_{r,k-r}$ to $-[t^r]\widetilde\omega_j/2$. The leading terms and the moment identity in the proof of \cref{thm:finiteHperiodlecture} show that step (i) makes every sum $\sum_{r\in I^-}c_r[t^r]\widetilde\omega_j$ vanish, i.e. $\sum_{r\in I}c_rH_{r,k-r}$ vanishes modulo odd products. Applying the same triangular argument to the $\omega_j$, we see that a restricted period polynomial is determined by the free coefficients $p_{r+2}$ used in step (i). Since the low-coordinate corrections and the vector $(2r-k)_r$ span a complementary space, step (i) is an isomorphism onto these antisymmetric relation vectors modulo products.

Now we write $R=\sum_{r\in I}\binom m{r-1}\nu_rX^{r-1}Y^{k-r-1}$. The even--odd conversion \eqref{eq:pevevzevevexplicitformula} gives
\begin{align}\label{eq:formalrecursiveconversionlecture}
 \sum_{r\in I}c_rH_{r,k-r}=\sum_{r\in I}\nu_rZ_{r,k-r}\qquad\text{in }\widetilde{\dz}_k\,.
\end{align}
To check this coefficient calculation, notice that the binomial theorem gives the non-Bernoulli coefficient $-(r-1)c_r-\sum_u c_u\binom{k-r}{u}-\sum_u c_u\binom ru$. If we write $F=\sum_{n=0}^m\binom mn f_nX^nY^{m-n}$, the Bernoulli series gives
\begin{align*}
 f_n=\sum_{\substack{u\in I\\u-1\leq n}}\frac{u+1}{u}c_u\binom n{u-1}B_{n-u+1}\,.
\end{align*}
Substituting this in $F(X+Y,X)-F(Y,X)$ gives exactly the coefficients of the two even--even symbols in \eqref{eq:formalfiniteHrs} by the same conversion formula, and the term $F(1,1)$ removes $Z_{k-1,1}$ by the odd sum formula. This proves \eqref{eq:formalrecursiveconversionlecture}.

Since the element in \eqref{eq:formalrecursiveconversionlecture} vanishes modulo odd products, $R-R|\epsilon$ is a restricted period polynomial. Its relation in \cref{prop:formalfiniteperiodrelation} shows that the $b_{r,s}$ in step (iii) are symmetric for odd $r,s$ and that $\sum_r(\nu_r+b_{r,k-r})Z_{r,k-r}=0$. But in the reduced quotient we have $H_{r,s}+H_{s,r}=rsP_{r,s}$, and therefore
\begin{align*}
 \sum_{r\in I}\frac{b_{r,k-r}}{r(k-r)}H_{r,k-r}=\sum_{r\in I}b_{r,k-r}Z_{r,k-r}\,.
\end{align*}
Together with \eqref{eq:formalrecursiveconversionlecture}, this gives $\sum_rC_r(p)H_{r,k-r}=0$, and we can apply \cref{thm:formalfinitebetaisomorphism}. Finally, $C_r(p)-C_{k-r}(p)=2c_r$ proves the injectivity modulo the Eisenstein direction. Both spaces have dimension $\dim_\C\ms_k$, so the map gives all relations, and for $k=4,6$ both spaces are zero.
\end{proof}
After realization, the same coefficients give relations among the actual finite values, but we do not know whether these are all relations there. For example, $p=X^2Y^2(X^2-Y^2)^3$ in weight $12$ gives
\begin{align*}
 (c_3,c_5)=\left(\frac1{175},-\frac1{350}\right),\qquad
 (C_3,C_5,C_7,C_9)=\frac1{1575}(16,9,18,-2)\,.
\end{align*}
So the relation in \cref{ex:finiteH12lecture} also holds in the universal formal algebra. Notice that no matrix inversion is needed for this construction.

\subsubsection{The determinant question}

%% file: chap_FormalFiniteExercises.tex
\begin{exer}\label{ex:formalfiniteinversionlecture}
Let $k\geq6$ be even and $I=\{3,5,\ldots,k-3\}$, and read the right-hand side of \eqref{eq:formalfiniteHrs} in the full space $\dz_k$. Prove that
\begin{align*}
 H_{r,s}+H_{s,r}&=rs(P_{r,s}-Z_k),\\
 Z_{k-1,1}&=\frac14Z_k-\sum_{r\in I}\frac{H_{r,k-r}}{r(k-r)}\,.
\end{align*}
Let $A_k$ be the matrix in \eqref{eq:formalfiniteAmatrix}. Use \eqref{eq:pevevzevevexplicitformula} and the odd sum formula to show that
\begin{align*}
 \sum_{u\in I}A_{r,u}(Z_{u,k-u}+Z_k/2)=2H_{r,k-r}+\kappa_{r,k-r}Z_k,
\end{align*}
where
\begin{align*}
 \kappa_{r,s}=rs+s-r-\frac{k-1}{2}\binom{k-2}{r-1}
 \left(\frac{r+1}{r}B_{s-1}-\frac{s+1}{s}B_{r-1}\right)\,.
\end{align*}
If $\det A_k\neq0$, deduce the conditional inverse for $t\in I$:
\begin{align*}
 Z_{t,k-t}=2\sum_{r\in I}(A_k^{-1})_{t,r}H_{r,k-r}
 +\left(-\frac12+\sum_{r\in I}(A_k^{-1})_{t,r}\kappa_{r,k-r}\right)Z_k\,.
\end{align*}
Then the boundary identity and the even--odd conversion reconstruct every admissible $Z_{r,s}$. Show that $Z_{1,k-1}$ is independent of these elements. Finally, if $\det(1+xA_k)=\sum_{j=0}^Ne_jx^j$, use Cayley--Hamilton to write
\begin{align*}
 A_k^{-1}=\frac{(-1)^{N-1}}{e_N}\sum_{j=0}^{N-1}(-1)^je_jA_k^{N-1-j},\qquad N=|I|\,.
\end{align*}
Substitute this expression into the formula for $Z_{t,k-t}$ and multiply by $e_N=\det A_k$. Notice that the resulting identity involves only finite Bernoulli sums, and that it holds unconditionally, even if the determinant vanishes.
\end{exer}

\begin{exer}\label{ex:formalfinitedeterminantreductions}
Put $m=k-2$ and $I^-=\{r\in I:2r<k\}$. Define
\begin{align*}
 d_{r,u}&=\frac{u+1}{u}\binom{r-1}{u-1}B_{r-u}\quad(u\leq r),\qquad d_{r,u}=0\quad(u>r),\\
 (\mathsf B_k)_{r,u}&=-m\delta_{r,u}+2(d_{r,u}-d_{r,k-u}-d_{k-r,u}+d_{k-r,k-u})\quad(r,u\in I^-).
\end{align*}
\begin{enumerate}[(i)]
\item Prove
\begin{align*}
 \det A_k=\prod_{\substack{r\in I\\r\leq k/2}}2r(k-r)\cdot\det\mathsf B_k\,.
\end{align*}
(Hint: The symmetric part is given by the first identity in \cref{ex:formalfiniteinversionlecture}. For antisymmetric $c_{k-r}=-c_r$, the difference of the coefficients of $Z_{r,k-r}$ and $Z_{k-r,r}$ in $\sum_uc_uH_{u,k-u}$ is $(\mathsf B_kc)_r$.)
\item For an odd polynomial $o(t)$ of degree at most $m-1$, put
\begin{align*}
 d_o(t)&=o'(t)+o(t)/t-2o'(0),\\
 z_o(t)&=o(t)+(t^md_o(1/t)-d_o(t))/m,\\
 \mathcal J_m(o)&=z_o(t+1)+t^mz_o(1+1/t)-z_o(t)-t^mz_o(1/t)\,.
\end{align*}
Show that all these expressions are polynomials and that $\det A_k\neq0$ is equivalent to the injectivity of $\mathcal J_m$. If $\mathcal J_m(o)=0$, interpret the coefficients of $z_o$ as a functional on $\overline{\dz}_k$ which annihilates every $H_{r,s}$.
\item Suppose $k\geq8$. Let $p\in W_k^-$ be restricted, i.e. $p(X,0)=0$, and let $c_r(p)$ be its coefficients after step (i) of the complete period transfer. Show that there is a unique restricted period polynomial
\begin{align*}
 \Upsilon_k(p)=\sum_{r+s=k}\binom{k-2}{r-1}y_rX^{r-1}Y^{s-1},\qquad
 (y_r)_{r\in I^-}=\mathsf B_k(c_r(p))_{r\in I^-}\,.
\end{align*}
In $\overline{\dz}_k$, prove
\begin{align*}
 \sum_{r\in I}C_r(p)H_{r,k-r}
 =\sum_{\substack{r+s=k\\r,s\text{ odd}}}q_{r,s}(\Upsilon_k(p))Z_{r,s}\,.
\end{align*}
Here the $q_{r,s}$ are defined by \eqref{eq:formalfiniteperiodcoefficients}. Deduce that $\det A_k\neq0$ is equivalent to the invertibility of $\Upsilon_k$ on the restricted period space, i.e. the remaining question is reduced to a space of dimension $\dim_\C\ms_k$.
\end{enumerate}
Of course, none of these reductions proves the nonvanishing of the determinant in arbitrary weight.
\end{exer}

%% file: main.bbl
\begin{thebibliography}{99}
\setlength{\itemsep}{2pt plus 1pt minus 1pt}
\bibitem[Ap]{Ap} R.~Ap\'ery:
{\itshape Irrationalit\'e de $\zeta(2)$ et $\zeta(3)$}, in ``Journ\'ees Arithm\'etiques de Luminy'', Ast\'erisque \textbf{61} (1979), 11--13.
\url{https://numdam.org/item/AST_1979__61__11_0/}

\bibitem[AK]{AK} T. Arakawa, M. Kaneko, {\itshape \begin{CJK}{UTF8}{min}多重ゼータ値入門\end{CJK} (Introduction to multiple zeta values)}, Kyushu University MI Lecture Note Series \textbf{23} (2010, corrected and revised January 2016), Japanese. (available at \url{https://www2.math.kyushu-u.ac.jp/~mkaneko/papers/MZV_LectureNotes.pdf})

\bibitem[Ba1]{Ba1} H.~Bachmann:
{\itshape Multiple Zeta-Werte und die Verbindung zu Modulformen durch Multiple Eisensteinreihen}, Master thesis (German), Hamburg University (2012).\\(available at \url{http://www.henrikbachmann.com})

\bibitem[Ba2]{Ba2} H.~Bachmann:
{\itshape Derivatives of q-analogues of multiple zeta values}, Talk at Kindai University - February 2017. (available at \url{https://www.henrikbachmann.com/uploads/7/7/6/3/77634444/qderiv_bachmann_osaka.pdf}) 

\bibitem[Ba3]{Ba3} H.~Bachmann:
{\itshape Double shuffle relations for q-analogues of multiple zeta values, their derivatives and the connection to multiple Eisenstein series}, RIMS K\^oky\^uroku No.~{\bf 2015} (2017), 22--43.

\bibitem[Ba4]{Ba4} H.~Bachmann:
{\itshape Interpolated Schur multiple zeta values}, J. Aust. Math. Soc. {\bf 104} (2018), no.~3, 289--307,
\href{https://doi.org/10.1017/S1446788717000209}{doi:10.1017/S1446788717000209}.

\bibitem[Ba5]{Ba5} H. Bachmann, {\itshape Multiple zeta values and their relations }(poster),\\
\url{https://www.henrikbachmann.com/uploads/7/7/6/3/77634444/mzv_poster_3.pdf}

\bibitem[Ba6]{Ba6} H. Bachmann, {\itshape Introduction to modular forms}, Lecture notes, Nagoya University (Fall 2018). \\
(available at \url{https://www.henrikbachmann.com/mf2018.html})

\bibitem[Ba7]{Ba7} H.~Bachmann:
{\itshape The algebra of bi-brackets and regularized multiple Eisenstein series}, J. Number Theory {\bf 200}, 2019, 260--294.

\bibitem[Ba8]{Ba8} H. Bachmann, {\itshape Multiple Eisenstein series and q-analogues of multiple zeta values}, in ``Periods in Quantum Field Theory and Arithmetic'', Springer Proceedings in Mathematics \& Statistics {\bf 314} (2020), 173--235,
\href{https://doi.org/10.1007/978-3-030-37031-2_8}{doi:10.1007/978-3-030-37031-2\_8}.

\bibitem[Ba9]{Ba9} H.~Bachmann:
{\itshape Modular forms and $q$-analogues of modified double zeta values},
Abh. Math. Semin. Univ. Hambg. {\bf 90} (2020), no.~2, 201--213,
\href{https://doi.org/10.1007/s12188-020-00227-7}{doi:10.1007/s12188-020-00227-7}.

\bibitem[Ba10]{Ba10} H.~Bachmann:
{\itshape Stuffle regularized multiple Eisenstein series revisited},
RIMS K\^oky\^uroku {\bf 2238} (2023), 73--86,
\href{https://arxiv.org/abs/2212.10700}{arXiv:2212.10700}.

\bibitem[Ba11]{Ba11} H.~Bachmann:
{\itshape MacMahon's sums-of-divisors and their connection to multiple Eisenstein series},
Res. Number Theory {\bf 10} (2024), no.~2, Paper No.~50,
\href{https://doi.org/10.1007/s40993-024-00537-2}{doi:10.1007/s40993-024-00537-2}.

\bibitem[Ba12]{Ba12} H.~Bachmann:
{\itshape A weighted sum formula for double Eisenstein series},
preprint, \href{https://arxiv.org/abs/2607.21358}{arXiv:2607.21358} (2026).

\bibitem[BaBu]{BaBu}  H.~Bachmann, A.~Burmester:
\newblock {\itshape Combinatorial multiple Eisenstein series}, 
\newblock Res. Math. Sci. {\bf 10} (2023), Paper No.~35,
\href{https://doi.org/10.1007/s40687-023-00398-8}{doi:10.1007/s40687-023-00398-8}.

\bibitem[BC]{BC} H.~Bachmann, S.~Charlton:
{\itshape Generalized Jacobi--Trudi determinants and evaluations of Schur multiple zeta values}, European J. Combin. {\bf 87} (2020), Paper No.~103133,
\href{https://doi.org/10.1016/j.ejc.2020.103133}{doi:10.1016/j.ejc.2020.103133}.

\bibitem[BI]{BI}  H.~Bachmann, J.-W.~van Ittersum:
\newblock {\itshape Partitions, multiple zeta values and the $q$-bracket}, 
\newblock Selecta Math. (N.S.) {\bf 30} (2024), Paper No.~3, 46~pp.,
\href{https://doi.org/10.1007/s00029-023-00893-4}{doi:10.1007/s00029-023-00893-4}.

\bibitem[BIM]{BIM} H.~Bachmann, J.-W.~van Ittersum:
\newblock {\itshape Formal multiple Eisenstein series and their derivations}, with an appendix by N.~Matthes,
\newblock Adv. Math. {\bf 487} (2026), Paper No.~110739,
\href{https://doi.org/10.1016/j.aim.2025.110739}{doi:10.1016/j.aim.2025.110739}.

\bibitem[BKa]{BKa} H.~Bachmann, H.~Kanno:
{\itshape The $\mathfrak{sl}_2$-algebra structure of multiple Eisenstein series},
preprint, \href{https://arxiv.org/abs/2609.03777}{arXiv:2609.03777} (2026).

\bibitem[BKM]{BKM} H.~Bachmann, H.~Kanno, T.~Maesaka:
{\itshape Relations and derivatives of multiple Eisenstein series},
preprint, \href{https://arxiv.org/abs/2602.08176}{arXiv:2602.08176} (2026).

\bibitem[BK1]{BK1} H.~Bachmann, U.~K\"uhn:
{\itshape The algebra of generating functions for multiple divisor sums and applications to multiple zeta values}, Ramanujan J. {\bf 40} (2016), 605--648. 

\bibitem[BK2]{BK2} H.~Bachmann, U.~K\"uhn:
{\itshape A dimension conjecture for q-analogues of multiple zeta values}, in ``Periods in Quantum Field Theory and Arithmetic'', Springer Proceedings in Mathematics \& Statistics {\bf 314} (2020), 237--258,
\href{https://doi.org/10.1007/978-3-030-37031-2_9}{doi:10.1007/978-3-030-37031-2\_9}.

\bibitem[BKuM]{BKuM} H.~Bachmann, U.~K\"uhn, N.~Matthes:
{\itshape Realizations of the formal double Eisenstein space},
to appear in Tohoku Math. J.,
\href{https://arxiv.org/abs/2109.04267}{arXiv:2109.04267}.

\bibitem[BaR]{BaR} H.~Bachmann, Risan:
{\itshape Formal finite multiple zeta values},
preprint, \href{https://arxiv.org/abs/2608.15480}{arXiv:2608.15480} (2026).

\bibitem[BTT1]{BTT1} H.~Bachmann, Y.~Takeyama, K.~Tasaka:
{\itshape Cyclotomic analogues of finite multiple zeta values}, Compositio Math., {\bf 154} (12), 2018, 2701--2721,
\href{https://doi.org/10.1112/S0010437X18007583}{doi:10.1112/S0010437X18007583}.\vspace{0.1cm}

\bibitem[BTT2]{BTT2} H.~Bachmann, Y.~Takeyama, K.~Tasaka:
{\itshape Special values of finite multiple harmonic q-series at roots of unity}, Algebraic Combinatorics, Resurgence, Moulds and Applications (CARMA), Volume 2, IRMA Lectures in Mathematics and Theoretical Physics 32, 2020, 1--18,
\href{https://doi.org/10.4171/205-1/1}{doi:10.4171/205-1/1}.\vspace{0.1cm}

\bibitem[BTan1]{BTan1} H.~Bachmann, T.~Tanaka:
{\itshape Rooted tree maps and the derivation relation for multiple zeta values}, Int. J. Number Theory {\bf 14} (2018), no. 10, 2657--2662.

\bibitem[BTan2]{BTan2} H. Bachmann and T. Tanaka: {\itshape Rooted tree maps and the Kawashima relations for multiple zeta values}, Kyushu J. Math. {\bf 74} (2020), no.~1, 169--176,
\href{https://doi.org/10.2206/kyushujm.74.169}{doi:10.2206/kyushujm.74.169}.

\bibitem[BT]{BT} H. Bachmann and K. Tasaka, {\itshape The double shuffle relations for multiple Eisenstein series}, Nagoya Math. J. {\bf 230} (2018), 180--212,
\href{https://doi.org/10.1017/nmj.2017.9}{doi:10.1017/nmj.2017.9}.

\bibitem[BTs]{BTs} H.~Bachmann, H.~Tsumura:
{\itshape On multiple series of Eisenstein type}, Ramanujan J. {\bf 42} (2017), no.~2, 479--489,
\href{https://doi.org/10.1007/s11139-015-9738-0}{doi:10.1007/s11139-015-9738-0}.

\bibitem[BY1]{BY1} H.~Bachmann, Y.~Yamasaki:
{\itshape Checkerboard style Schur multiple zeta values and odd single zeta values},
Math. Z. {\bf 290} (2018), no.~3--4, 1173--1197,
\href{https://doi.org/10.1007/s00209-018-2058-5}{doi:10.1007/s00209-018-2058-5}.

\bibitem[BY2]{BY2} H.~Bachmann, J.~Yu:
{\itshape Schur Eisenstein series and Schur MacMahon series},
preprint, \href{https://arxiv.org/abs/2607.27702}{arXiv:2607.27702} (2026).

\bibitem[BR]{BR} K.~Ball, T.~Rivoal:
{\itshape Irrationalit\'e d'une infinit\'e de valeurs de la fonction z\^eta aux entiers impairs}, Invent. Math. {\bf 146} (2001), no.~1, 193--207,
\href{https://doi.org/10.1007/s002220100168}{doi:10.1007/s002220100168}.

\bibitem[BS]{BS} S.~Baumard, L.~Schneps:
{\itshape Period polynomial relations between double zeta values}, Ramanujan J. {\bf 32} (2013), no.~1, 83--100,
\href{https://doi.org/10.1007/s11139-013-9466-2}{doi:10.1007/s11139-013-9466-2}.

\bibitem[BoBr]{BoBr} J. M. Borwein and D. M. Bradley, {\itshape Thirty-two Goldbach variations}, Intl. J. Number Theory {\bf 2} (2006), 65--103.

\bibitem[Bo]{Bo} O.~Bouillot: 
{\itshape The algebra of multitangent functions}, J. Algebra {\bf 410}, 2014, 148--238.

\bibitem[Bra]{Bra} D. M.~Bradley:
{\it Multiple q-zeta values}, J. Algebra {\bf 283} (2005), 752--798.

\bibitem[Bri]{Bri} B. Brindle: {\itshape{A unified approach to qMZVs}}, INTEGERS {\bf 24} (2024), Paper No.~A6, 41~pp.

\bibitem[BroK]{BroK} D. Broadhurst, D. Kreimer:
{\it Association of multiple zeta values with positive knots via Feynman diagrams up to 9 loops}, Phys. Lett. B {\bf 393} (1997), 403--412.

\bibitem[Br1]{Br1} F.~Brown:
{\itshape Multiple zeta values and periods of moduli spaces $\overline{\mathfrak{M}}_{0,n}$}, Ann. Sci. \'Ec. Norm. Sup\'er. (4) {\bf 42} (2009), no.~3, 371--489,
\href{https://doi.org/10.24033/asens.2099}{doi:10.24033/asens.2099}.

\bibitem[Br2]{Br2} F.~Brown:
{\itshape Mixed Tate motives over $\Z$}, Ann. of Math. (2) {\bf 175} (2012), no.~2, 949--976.

\bibitem[Br3]{Br3} F.~Brown:
{\itshape Depth-graded motivic multiple zeta values}, Compos. Math. {\bf 157} (2021), no.~3, 529--572,
\href{https://doi.org/10.1112/S0010437X20007654}{doi:10.1112/S0010437X20007654}.

\bibitem[BrZ]{BrZ} F.~Brown, W.~Zudilin:
{\itshape On cellular rational approximations to $\zeta(5)$}, preprint,
\href{https://arxiv.org/abs/2210.03391}{arXiv:2210.03391v3} (2022, revised 2026), 32~pp.

\bibitem[BF]{BF} J.~I. Burgos Gil, J. Fres\'an: {\itshape Multiple zeta values: From numbers to motives}, with contributions by U.~K\"uhn, Clay Mathematics Proceedings, to appear.
(available at \url{http://javier.fresan.perso.math.cnrs.fr/mzv.pdf}) 

\bibitem[Bu1]{Bu1}  A.~Burmester:
\newblock {\itshape An algebraic approach to multiple q-zeta values}, 
\newblock \href{https://ediss.sub.uni-hamburg.de/handle/ediss/10158}{PhD Thesis, Universit\"at Hamburg} (2023), 232~pp.

\bibitem[Bu2]{Bu2} A.~Burmester:
{\itshape Balanced multiple q-zeta values}, Adv. Math. {\bf 439} (2024), Paper No.~109487,
\href{https://doi.org/10.1016/j.aim.2024.109487}{doi:10.1016/j.aim.2024.109487}.

\bibitem[Bu3]{Bu3} A.~Burmester:
{\itshape A generalization of formal multiple zeta values related to multiple Eisenstein series and multiple q-zeta values},
J. Number Theory {\bf 269} (2025), 106--137,
\href{https://doi.org/10.1016/j.jnt.2024.09.011}{doi:10.1016/j.jnt.2024.09.011}.

\bibitem[BCK]{BCK} A.~Burmester, N.~Confurius, U.~K\"uhn:
{\itshape AGZT-Lectures on formal multiple zeta values},
preprint, \href{https://arxiv.org/abs/2406.13630}{arXiv:2406.13630} (2024), 72~pp.


\bibitem[Car]{Car} L.~Carlitz:
{\itshape A degenerate Staudt--Clausen theorem}, Arch. Math. (Basel) {\bf 7} (1956), 28--33,
\href{https://doi.org/10.1007/BF01900520}{doi:10.1007/BF01900520}.

\bibitem[ChK]{ChK} S.~Charlton, A.~Keilthy:
{\itshape Evaluation of $\zeta(2,\ldots,2,4,2,\ldots,2)$ and period polynomial relations},
Forum Math. Sigma {\bf 12} (2024), Paper No.~e46, 50~pp.,
\href{https://doi.org/10.1017/fms.2024.16}{doi:10.1017/fms.2024.16}.

\bibitem[CK]{CK} A.~Connes, D.~Kreimer:
{\itshape Hopf algebras, renormalization and noncommutative geometry}, Commun. Math. Phys. {\bf 199} (1998), 203--242.

\bibitem[C]{C} K. Conrad, {\itshape Modular forms}, Lecture notes. 
(available at \url{https://ctnt-summer.math.uconn.edu/wp-content/uploads/sites/1632/2016/02/CTNTmodularforms.pdf}) 

\bibitem[DG]{DG} P.~Deligne, A.~B. Goncharov:
{\itshape Groupes fondamentaux motiviques de Tate mixte}, Ann. Sci. \'Ec. Norm. Sup\'er. (4) {\bf 38} (2005), no.~1, 1--56,
\href{https://doi.org/10.1016/j.ansens.2004.11.001}{doi:10.1016/j.ansens.2004.11.001}.

\bibitem[D]{D} V.~Drinfeld, {\itshape On quasitriangular quasi-Hopf algebras and a group closely connected with $\operatorname{Gal}(\overline{\Q}/\Q)$}, Leningrad Math. J. {\bf 2} (1991), no. 4, 829--860.

\bibitem[Ec]{Ec} J.~\'Ecalle:
{\itshape ARI/GARI, la dimorphie et l'arithm\'etique des multiz\^etas: un premier bilan}, J. Th\'eor. Nombres Bordeaux {\bf 15} (2003), no.~2, 411--478,
\href{https://doi.org/10.5802/jtnb.410}{doi:10.5802/jtnb.410}.

\bibitem[ELO]{ELO} M.~Eie, W.~Liaw, Y.~Ong, {\itshape A restricted sum formula among multiple zeta values}, J. Number Theory {\bf 129} (2009), 908--921.

\bibitem[Eu]{Eu} L.~Euler:
{\itshape De summis serierum reciprocarum}, Comment. Acad. Sci. Petropolitanae {\bf 7} (1740), 123--134. Written in 1735.
\url{https://scholarlycommons.pacific.edu/euler-works/41/}

\bibitem[Fi]{Fi} S.~Fischler:
{\itshape Linear independence of odd zeta values using Siegel's lemma}, J. London Math. Soc. (2) {\bf 113} (2026), no.~4, Paper No.~e70535,
\href{https://doi.org/10.1112/jlms.70535}{doi:10.1112/jlms.70535}.

\bibitem[FSZ]{FSZ} S.~Fischler, J.~Sprang, W.~Zudilin:
{\itshape Many odd zeta values are irrational}, Compos. Math. {\bf 155} (2019), no.~5, 938--952,
\href{https://doi.org/10.1112/S0010437X1900722X}{doi:10.1112/S0010437X1900722X}.

\bibitem[F1]{F1} H.~Furusho, {\itshape The multiple zeta value algebra and the stable derivation algebra}, Publ. Res. Inst. Math. Sci., {\bf 39}(4) (2003), 695--720.

\bibitem[F2]{F2} H.~Furusho, {\itshape Double shuffle relation for associators}, Ann. of Math. (2) {\bf 174} (2011), no. 1, 341--360.

\bibitem[F3]{F3} H.~Furusho, {\itshape The pentagon equation and the confluence relations}, Amer. J. Math. {\bf 144} (2022), no.~4, 873--894,
\href{https://doi.org/10.1353/ajm.2022.0018}{doi:10.1353/ajm.2022.0018}.

\bibitem[FHK]{FHK} H.~Furusho, M.~Hirose, N.~Komiyama:
{\itshape Associators in mould theory}, to appear in J. Lie Theory,
\href{https://arxiv.org/abs/2312.15423}{arXiv:2312.15423} (2023, revised 2026), 81~pp.

\bibitem[GKZ]{GKZ} H.~Gangl, M. Kaneko, D. Zagier:
{\it Double zeta values and modular forms}, in ``Automorphic forms and zeta functions'', World Sci. Publ., Hackensack, NJ (2006), 71--106,
\href{https://doi.org/10.1142/9789812774415_0004}{doi:10.1142/9789812774415\_0004}.

\bibitem[GV]{GV} I.~Gessel, G.~Viennot:
{\itshape Binomial determinants, paths, and hook length formulae}, Adv. Math. {\bf 58} (1985), no.~3, 300--321,
\href{https://doi.org/10.1016/0001-8708(85)90121-5}{doi:10.1016/0001-8708(85)90121-5}.

\bibitem[Gon]{Gon} A.~B. Goncharov:
{\itshape Galois symmetries of fundamental groupoids and noncommutative geometry}, Duke Math. J. {\bf 128} (2005), no.~2, 209--284,
\href{https://doi.org/10.1215/S0012-7094-04-12822-2}{doi:10.1215/S0012-7094-04-12822-2}.

\bibitem[G]{G} A. Granville, {\itshape A decomposition of Riemann's zeta function}, in ``Analytic Number Theory'', London Mathematical Society Lecture Note Series {\bf 247} (1997), 95--102,
\href{https://doi.org/10.1017/CBO9780511666179.009}{doi:10.1017/CBO9780511666179.009}.

\bibitem[HaST]{HaST} T.~Hara, K.~Sakugawa and K.~Tasaka:
{\itshape Symmetric multiple Eisenstein series}, preprint,
\href{https://arxiv.org/abs/2601.13626}{arXiv:2601.13626} (2026), 43~pp.

\bibitem[Hi1]{Hi1} M.~Hirose:
{\itshape Double shuffle relations for refined symmetric multiple zeta values}, Doc. Math. {\bf 25} (2020), 365--380,
\href{https://doi.org/10.4171/DM/750}{doi:10.4171/DM/750}.

\bibitem[Hi2]{Hi2} M.~Hirose:
{\itshape Modular phenomena for regularized double zeta values}, Israel J. Math. {\bf 261} (2024), no.~2, 501--547,
\href{https://doi.org/10.1007/s11856-023-2587-4}{doi:10.1007/s11856-023-2587-4}.

\bibitem[HMSW]{HMSW} M.~Hirose, T.~Maesaka, S.~Seki and T.~Watanabe:
{\itshape The $\Z$-module of multiple zeta values is generated by ones for indices without ones},
preprint, \href{https://arxiv.org/abs/2505.07221}{arXiv:2505.07221} (2025).

\bibitem[HMW]{HMW} M.~Hirose, T.~Maesaka and T.~Watanabe:
{\itshape A unified proof of conjectures on the spaces of multiple $q$-zeta values},
preprint, \href{https://arxiv.org/abs/2605.28584}{arXiv:2605.28584} (2026).

\bibitem[HMOS]{HMOS} M. Hirose, H. Murahara, T. Onozuka and N. Sato: 
{\itshape Linear relations of Ohno sums of multiple zeta values}, Indag. Math. (N.S.) {\bf 31} (2020), no.~4, 556--567,
\href{https://doi.org/10.1016/j.indag.2020.04.004}{doi:10.1016/j.indag.2020.04.004}.

\bibitem[HS]{HS} M.~Hirose, N.~Sato: {\itshape Iterated integrals on $\mathbb{P}^1 \backslash \{0, 1,\infty, z\}$ and a class of relations among multiple zeta values}, Adv. Math. {\bf 348} (2019), 163--182,
\href{https://doi.org/10.1016/j.aim.2019.03.005}{doi:10.1016/j.aim.2019.03.005}.

\bibitem[HSS]{HSS} M.~Hirose, N.~Sato and S.~Seki: {\itshape The connector for the double Ohno relation}, Acta Arith. {\bf 201} (2021), no.~2, 109--118,
\href{https://doi.org/10.4064/aa200621-18-5}{doi:10.4064/aa200621-18-5}.

\bibitem[HST]{HST} M.~Hirose, N.~Sato and K.~Tasaka: {\itshape Eisenstein series identities based on partial fraction decomposition}, Ramanujan J., {\bf 38}(3) (2015), 455--463.

\bibitem[H1]{H1} M.~E. Hoffman:
{\itshape Multiple harmonic series}, Pacific J. Math. {\bf 152} (1992), no.~2, 275--290,
\href{https://doi.org/10.2140/pjm.1992.152.275}{doi:10.2140/pjm.1992.152.275}.

\bibitem[H2]{H2} M. E. Hoffman: {\itshape The algebra of multiple harmonic series}, J. Algebra {\bf 194} (1997), 477--495.

\bibitem[H3]{H3} M. E. Hoffman: {\itshape Quasi-shuffle products}, J. Algebraic Combin. {\bf 11} (2000), 49--68.

\bibitem[H4]{H4} M. E. Hoffman: {\itshape Quasi-symmetric functions and mod p multiple harmonic sums}, Kyushu J. Math. {\bf 69} (2015), 345--366.

\bibitem[H5]{H5} M. E. Hoffman: {\itshape References on multiple zeta values and Euler sums}, (webpage)\\ \url{https://www.usna.edu/Users/math/meh/biblio.php}.

\bibitem[HI]{HI} M. E. Hoffman and K. Ihara: {\itshape Quasi-shuffle products revisited}, J. Algebra {\bf 481} (2017), 293--326. 

\bibitem[HO]{HO} M.~Hoffman and Y.~Ohno, {\itshape Relations of multiple zeta values and their algebraic expression}, J. Algebra {\bf 262} (2003), 332--347.

\bibitem[I]{I} K.~Ihara: 
{\itshape Derivation and double shuffle relations for multiple zeta values, joint work with M. Kaneko, D. Zagier}, RIMS K\^oky\^uroku {\bf 1549} (2007), 47--63.

\bibitem[IKZ]{IKZ} K.~Ihara, M.~Kaneko and D.~Zagier, 
{\itshape Derivation and double shuffle relations for multiple zeta values}, Compositio Math. {\bf 142} (2006), 307--338,
\href{https://doi.org/10.1112/S0010437X0500182X}{doi:10.1112/S0010437X0500182X}.

\bibitem[J]{J} D.~Jarossay, \textit{Double m\'elange des multiz\^etas finis et multiz\^etas sym\'etris\'es}, 
\textit{Comptes Rendus Math\'ematique} \textbf{352} (2014), no.~10, 767--771,
\href{https://doi.org/10.1016/j.crma.2014.08.005}{doi:10.1016/j.crma.2014.08.005}.


\bibitem[Kad]{Kad} S.~Kadota, {\itshape Certain weighted sum formulas for multiple zeta values with some parameters}, Comment. Math. Univ. St. Pauli {\bf 66} (2017), no. 1-2, 1--13.

\bibitem[K1]{K1} M. Kaneko:
{\itshape On an extension of the derivation relation for multiple zeta values}, in L.~Weng and M.~Kaneko (eds.), ``The Conference on L-Functions'' (Fukuoka, 18--23 February 2006), World Sci. Publ., Hackensack, NJ (2007), 89--94,
\href{https://doi.org/10.1142/9789812772398_0005}{doi:10.1142/9789812772398\_0005}.

\bibitem[K2]{K2} M. Kaneko: {\itshape \begin{CJK}{UTF8}{min}有限多重ゼータ値\end{CJK} (Finite multiple zeta values)}, RIMS K\^oky\^uroku Bessatsu {\bf B68} (2017), 175--190 (Japanese).

\bibitem[K3]{K3} M.~Kaneko:
{\itshape An introduction to classical and finite multiple zeta values}, Publications Math\'ematiques de Besan\c{c}on. Alg\`ebre et th\'eorie des nombres, 2019/1 (2019), 103--129,
\href{https://doi.org/10.5802/pmb.31}{doi:10.5802/pmb.31}.

\bibitem[K4]{K4} M.~Kaneko:
{\itshape \begin{CJK}{UTF8}{min}多重ゼータ値入門\end{CJK} (Introduction to Multiple Zeta Values)},
Iwanami Studies in Advanced Mathematics, Iwanami Shoten, Tokyo, 2026, ix+209~pp. (Japanese), ISBN 978-4-00-029940-4.

\bibitem[KNT]{KNT} M. Kaneko, M. Noro and K. Tsurumaki: {\itshape
On a conjecture for the dimension of the space of the multiple zeta values}, in M.~Stillman, J.~Verschelde and N.~Takayama (eds.), ``Software for Algebraic Geometry'', IMA Volumes in Mathematics and its Applications {\bf 148}, Springer, New York (2008), 47--58,
\href{https://doi.org/10.1007/978-0-387-78133-4_4}{doi:10.1007/978-0-387-78133-4\_4}.

\bibitem[KT]{KT} M.~Kaneko, K.~Tasaka:
{\itshape Double zeta values, double Eisenstein series, and modular forms of level 2},
Math. Ann. {\bf 357} (2013), no.~3, 1091--1118,
\href{https://doi.org/10.1007/s00208-013-0930-5}{doi:10.1007/s00208-013-0930-5}.

\bibitem[KY]{KY} M.~Kaneko, S.~Yamamoto, {\itshape A new integral-series identity of multiple zeta values and regularizations}, Selecta Math. (N.S.) {\bf 24} (2018), no.~3, 2499--2521,
\href{https://doi.org/10.1007/s00029-018-0400-8}{doi:10.1007/s00029-018-0400-8}.

\bibitem[KZ]{KZ} M. Kaneko and D. Zagier:
{\itshape Finite multiple zeta values}, preprint, 50~pp., to appear in the Proceedings of the 17th MSJ-SI conference ``Modular Forms and Multiple Zeta Values'',
available at \url{https://www2.math.kyushu-u.ac.jp/~mkaneko/papers/FMZV.pdf}.

\bibitem[Ka]{Ka} H.~Kanno:
{\itshape Shuffle regularization for multiple Eisenstein series of level $N$},
Ramanujan J. {\bf 67} (2025), no.~4, Paper No.~95,
\href{https://doi.org/10.1007/s11139-025-01136-0}{doi:10.1007/s11139-025-01136-0}.

\bibitem[Kaw]{Kaw} G.~Kawashima:
{\itshape A class of relations among multiple zeta values},
J. Number Theory {\bf 129} (2009), no.~4, 755--788,
\href{https://doi.org/10.1016/j.jnt.2008.11.002}{doi:10.1016/j.jnt.2008.11.002}.

\bibitem[Ki]{Ki} K.~Kina:
{\itshape An explicit expression for MZVs in terms of symmetric MZVs},
preprint, \href{https://arxiv.org/abs/2607.05795v2}{arXiv:2607.05795v2} (2026).

\bibitem[KoZ]{KoZ} W. Kohnen, D. Zagier: {\itshape Modular forms with rational periods}, In: Rankin, R. A. (ed.): Modular Forms, pp. 197--249, Horwood, Chichester, 1984.

\bibitem[Lai]{Lai} L.~Lai:
{\itshape A note on the number of irrational odd zeta values, II}, to appear in J. Th\'eor. Nombres Bordeaux,
\href{https://arxiv.org/abs/2501.05321}{arXiv:2501.05321} (2025).

\bibitem[LZ]{LZ} L.~Lai, L.~Zhou:
{\itshape At least two of $\zeta(5),\zeta(7),\ldots,\zeta(35)$ are irrational}, Publ. Math. Debrecen {\bf 101} (2022), no.~3--4, 353--372,
\href{https://doi.org/10.5486/PMD.2022.9252}{doi:10.5486/PMD.2022.9252}.

\bibitem[La]{La} S. Lang: {\itshape Introduction to modular forms}, Grundlehren der mathematischen Wissenschaften, No.~222, Springer-Verlag, Berlin--New York, 1976,
\href{https://doi.org/10.1007/978-3-642-51447-0}{doi:10.1007/978-3-642-51447-0}.

\bibitem[JLi]{JLi} J.~Li:
{\itshape The depth structure of motivic multiple zeta values},
Math. Ann. {\bf 374} (2019), 179--209,
\href{https://arxiv.org/abs/1710.06135v4}{arXiv:1710.06135v4}.

\bibitem[L]{L} Z.~Li, {\itshape Regularized double shuffle and Ohno-Zagier relations of multiple zeta values}, J. Number Theory {\bf 133} (2013), no.~2, 596--610,
\href{https://doi.org/10.1016/j.jnt.2012.08.005}{doi:10.1016/j.jnt.2012.08.005}.

\bibitem[Lin]{Lin} F.~Lindemann:
{\itshape Ueber die Zahl $\pi$}, Math. Ann. {\bf 20} (1882), 213--225,
\href{https://doi.org/10.1007/BF01446522}{doi:10.1007/BF01446522}.

\bibitem[Ma]{Ma} D.~Ma:
{\itshape Period polynomial relations between formal double zeta values of odd weight}, Math. Ann. {\bf 365} (2016), no.~1--2, 345--362,
\href{https://doi.org/10.1007/s00208-015-1308-7}{doi:10.1007/s00208-015-1308-7}.

\bibitem[MaT]{MaT} D.~Ma, K.~Tasaka:
{\itshape Relationships between multiple zeta values of depths 2 and 3 and period polynomials}, Israel J. Math. {\bf 242} (2021), no.~1, 359--400,
\href{https://doi.org/10.1007/s11856-021-2139-8}{doi:10.1007/s11856-021-2139-8}.

\bibitem[MSW]{MSW} T.~Maesaka, S.~Seki and T.~Watanabe:
{\itshape Deriving two dualities simultaneously from a family of identities for multiple harmonic sums}, preprint, \href{https://arxiv.org/abs/2402.05730}{arXiv:2402.05730} (2024).

\bibitem[MatT]{MatT} N.~Matthes, K.~Tasaka:
{\itshape On \'Ecalle's and Brown's polar solutions to the double shuffle equations modulo products}, Kyushu J. Math. {\bf 73} (2019), no.~2, 337--356,
\href{https://doi.org/10.2206/kyushujm.73.337}{doi:10.2206/kyushujm.73.337}.

\bibitem[MJOP]{MJOP} H. N. Minh, G. Jacob, N. E. Oussous and M. Petitot: {\itshape Aspects combinatoires des polylogarithmes et des sommes d'Euler-Zagier}, S\'em. Lothar. Combin. {\bf 43} (2000), Art.~B43e, 29~pp.

\bibitem[M]{M} H. Murahara: {\itshape A study on relations among finite multiple zeta values}, Doctoral thesis, Graduate School of Mathematics, Kyushu University, January 2016.

\bibitem[NPY]{NPY} M.~Nakasuji, O.~Phuksuwan and Y.~Yamasaki:
{\itshape On Schur multiple zeta functions: a combinatoric generalization of multiple zeta functions}, Adv. Math. {\bf 333} (2018), 570--619,
\href{https://doi.org/10.1016/j.aim.2018.05.014}{doi:10.1016/j.aim.2018.05.014}.

\bibitem[Oh]{Oh} Y. Ohno, {\itshape A generalization of the duality and sum formulas on the multiple zeta values}, J. Number
Theory {\bf 74} (1999), no.~1, 39--43,
\href{https://doi.org/10.1006/jnth.1998.2314}{doi:10.1006/jnth.1998.2314}.

\bibitem[OZ]{OZ} Y. Ohno and D. Zagier, {\itshape Multiple zeta values of fixed weight, depth, and height}, Indag. Math. (N.S.) {\bf 12} (2001), no.~4, 483--487,
\href{https://doi.org/10.1016/S0019-3577(01)80037-9}{doi:10.1016/S0019-3577(01)80037-9}.

\bibitem[Ok]{Ok} A.~Okounkov. {\itshape Hilbert schemes and multiple q-zeta values}, Funct. Anal. Appl. {\bf 48} (2014), no.~2, 138--144,
\href{https://doi.org/10.1007/s10688-014-0054-z}{doi:10.1007/s10688-014-0054-z}.

\bibitem[OT]{OT} J. Okuda and Y. Takeyama, {\itshape On relations for the multiple q-zeta values}, Ramanujan J. {\bf 14} (2007), no.~3, 379--387,
\href{https://doi.org/10.1007/s11139-007-9053-5}{doi:10.1007/s11139-007-9053-5}.

\bibitem[On]{On} T. Onozuka, {\itshape Analytic Properties of Multiple Zeta Functions}, PhD Thesis, Nagoya University 2014.
(available at \url{https://hdl.handle.net/2237/20845})

\bibitem[Pan]{Pan} E.~Panzer:
{\itshape The parity theorem for multiple polylogarithms}, J. Number Theory {\bf 172} (2017), 93--113,
\href{https://doi.org/10.1016/j.jnt.2016.08.004}{doi:10.1016/j.jnt.2016.08.004}.

\bibitem[P]{P} Yu. A. Pupyrev, {\itshape Linear and algebraic independence of q-zeta values}, Math. Notes {\bf 78} (2005), no.~3--4, 563--568,
\href{https://doi.org/10.1007/s11006-005-0155-3}{doi:10.1007/s11006-005-0155-3}. Translated from Mat. Zametki {\bf 78} (2005), no.~4, 608--613.

\bibitem[Q]{Q}  Zhenbo Qin, {\itshape Multiple q-zeta values and traces}, preprint, \href{https://arxiv.org/abs/2505.14614}{arXiv:2505.14614} (2025).


\bibitem[Rac]{Rac} G.~Racinet:
{\itshape Doubles m\'elanges des polylogarithmes multiples aux racines de l'unit\'e},
Publ. Math. Inst. Hautes \'Etudes Sci. {\bf 95} (2002), 185--231,
\href{https://doi.org/10.1007/s102400200004}{doi:10.1007/s102400200004}.


\bibitem[R]{R} Risan, {\itshape Formal finite multiple zeta values and modular forms}, Master's thesis, Graduate School of Mathematics, Nagoya University (2025), 61~pp.
(available at \url{https://www.risan.io/} and at
\href{https://drive.google.com/file/d/1p1eCMkUlHQtxHLg4nW_z8kQHnx4iYWbT/view}{drive.google.com})

\bibitem[Sa]{Sa} M. Satoh, {\itshape Generalized quasi-shuffle products}, preprint, \href{https://arxiv.org/abs/2303.12386}{arXiv:2303.12386} (2023).

\bibitem[Sc]{Sc} L.~Schneps:
{\itshape ARI, GARI, Zig and Zag: An introduction to Ecalle's theory of multiple zeta values}, preprint,
\href{https://arxiv.org/abs/1507.01534}{arXiv:1507.01534} (2015, revised 2025), 79~pp.

\bibitem[Se1]{Se1} S. Seki, {\itshape Connectors}, RIMS K\^oky\^uroku {\bf 2160} (2020), 15--27.

\bibitem[Se2]{Se2} S. Seki: {\itshape Regular primes, non-Wieferich primes, and finite multiple zeta values of level $N$}, Integers 24, Paper No. A22, 2024.

\bibitem[Se3]{Se3} S.~Seki:
{\itshape A proof of the extended double shuffle relation without using integrals},
Kyushu J. Math. {\bf 79} (2025), no.~1, 191--198,
\href{https://doi.org/10.2206/kyushujm.79.191}{doi:10.2206/kyushujm.79.191}.

\bibitem[SY]{SY} S. Seki and S. Yamamoto, {\itshape A new proof of the duality of multiple zeta values and its generalizations}, Intern. J. Number Theory {\bf 15} (2019), no.~6, 1261--1265,
\href{https://doi.org/10.1142/S1793042119500702}{doi:10.1142/S1793042119500702}.
 

\bibitem[Tak1]{Tak1} Y. Takeyama: {\itshape The algebra of a q-analogue of multiple harmonic series}, SIGMA {\bf 9} (2013), Paper No.~061, 15~pp.,
\href{https://doi.org/10.3842/SIGMA.2013.061}{doi:10.3842/SIGMA.2013.061}.

\bibitem[Tak2]{Tak2} Y. Takeyama: {\itshape Derivations on the algebra of multiple harmonic q-series and their applications}, Ramanujan J. {\bf 52} (2020), no.~1, 41--65,
\href{https://doi.org/10.1007/s11139-019-00139-y}{doi:10.1007/s11139-019-00139-y}.

\bibitem[TT]{TT} Y. Takeyama and K. Tasaka: {\itshape Supercongruences of multiple harmonic q-sums and generalized finite/symmetric multiple zeta values}, Kyushu J. Math. {\bf 77}(1) (2023), 75--120,
\href{https://doi.org/10.2206/kyushujm.77.75}{doi:10.2206/kyushujm.77.75}.

\bibitem[Tan1]{Tan1} T. Tanaka:
{\itshape On the quasi-derivation relation for multiple zeta values}, J. Number Theory {\bf 129} (2009), no.~9, 2021--2034,
\href{https://doi.org/10.1016/j.jnt.2009.03.003}{doi:10.1016/j.jnt.2009.03.003}.

\bibitem[Tan2]{Tan2} T. Tanaka:
{\itshape Duality for multiple zeta values and related topics} (survey in Japanese), Proceedings of the second MZV seminar, 33--43, 2010.2. (available at \url{https://www.cc.kyoto-su.ac.jp/~tatsushi/doc.pdf})

\bibitem[Tan3]{Tan3} T. Tanaka: {\itshape Restricted sum formula and derivation relation for multiple zeta values}, preprint, \href{https://arxiv.org/abs/1303.0398}{arXiv:1303.0398} (2013).

\bibitem[Tan4]{Tan4} T. Tanaka:
{\itshape Rooted tree maps}, Commun. Number Theory Phys. {\bf 13} (2019), no.~3, 647--666,
\href{https://doi.org/10.4310/CNTP.2019.v13.n3.a6}{doi:10.4310/CNTP.2019.v13.n3.a6}.

\bibitem[TW]{TW} T.~Tanaka, N.~Wakabayashi, {\itshape An algebraic proof of the cyclic sum formula for multiple zeta values}, J. Algebra {\bf 323} (2010), no.~3, 766--778,
\href{https://doi.org/10.1016/j.jalgebra.2009.11.016}{doi:10.1016/j.jalgebra.2009.11.016}.

\bibitem[Tas]{Tas} K. Tasaka: 
{\itshape Hecke eigenform and double Eisenstein series},
Proc. Amer. Math. Soc. {\bf 148} (2020), no.~1, 53--58,
\href{https://doi.org/10.1090/proc/14680}{doi:10.1090/proc/14680}.

\bibitem[Ter]{Ter} T.~Terasoma:
{\itshape Mixed Tate motives and multiple zeta values}, Invent. Math. {\bf 149} (2002), no.~2, 339--369,
\href{https://doi.org/10.1007/s002220200218}{doi:10.1007/s002220200218}.

\bibitem[Ts]{Ts} H.~Tsumura, {\itshape Combinatorial relations for Euler-Zagier sums}, Acta Arith. {\bf 111} (2004), no.~1, 27--42,
\href{https://doi.org/10.4064/aa111-1-3}{doi:10.4064/aa111-1-3}.

\bibitem[Tu]{Tu} C.~Turan, {\itshape Ableitungen von Eisensteinreihen und Erzeugendenreihen von Multiplen Teilersummen}, Bachelor thesis, Universit\"at Hamburg, 2021. (German)
(available at \url{https://www.henrikbachmann.com/uploads/7/7/6/3/77634444/bachelorarbeit_turan.pdf})

\bibitem[W]{W} M. Waldschmidt: {\itshape Lectures on multiple zeta values}, preliminary draft of eight lectures at IMSc, Chennai, April 2011 (updated 21 April 2011). (available at \url{https://webusers.imj-prg.fr/~michel.waldschmidt/articles/pdf/MZV2011IMSc.pdf})


\bibitem[Yam]{Yam} S. Yamamoto: 
{\itshape Multivariable Hoffman-Ihara operators and the operad of formal power series}, 
J. Algebra {\bf 556} (2020), 634--648,
\href{https://doi.org/10.1016/j.jalgebra.2020.04.006}{doi:10.1016/j.jalgebra.2020.04.006}.

\bibitem[Yas]{Yas} S. Yasuda: {\itshape Finite real multiple zeta values generate the whole space $Z$}, Int. J. Number Theory {\bf 12} (2016), no.~3, 787--812,
\href{https://doi.org/10.1142/S1793042116500512}{doi:10.1142/S1793042116500512}.

\bibitem[Yu]{Yu} J.~Yu:
{\itshape The Young Tableaux Hopf algebra and multiple Schur series},
preprint, \href{https://arxiv.org/abs/2607.09157}{arXiv:2607.09157} (2026), 50~pp.

\bibitem[YZ]{YZ} H.~Yuan, J.~Zhao:
{\itshape Double shuffle relations of double zeta values and the double Eisenstein series at level $N$},
J. Lond. Math. Soc. (2) {\bf 92} (2015), no.~3, 520--546,
\href{https://doi.org/10.1112/jlms/jdv042}{doi:10.1112/jlms/jdv042}.

\bibitem[Za1]{Za1} D.~Zagier:
{\itshape Periods of modular forms and Jacobi theta functions},
Invent. Math. {\bf 104} (1991), no.~3, 449--465,
\href{https://doi.org/10.1007/BF01245085}{doi:10.1007/BF01245085}.

\bibitem[Za2]{Za2} D.~Zagier:
{\itshape Values of zeta functions and their applications}, in ``First European Congress of Mathematics'', Volume II, Progress in Math. {\bf 120}, Birkh\"auser-Verlag, Basel (1994), 497--512,
\href{https://doi.org/10.1007/978-3-0348-9112-7_23}{doi:10.1007/978-3-0348-9112-7\_23}.

\bibitem[Za3]{Za3} D. Zagier: {\itshape Elliptic modular forms and their applications}, in ``The 1-2-3 of modular forms'', Universitext, Springer-Verlag, Berlin, 2008, 1--103,
\href{https://doi.org/10.1007/978-3-540-74119-0_1}{doi:10.1007/978-3-540-74119-0\_1}.

\bibitem[Za4]{Za4} D.~Zagier:
{\itshape Evaluation of the multiple zeta values $\zeta(2,\ldots,2,3,2,\ldots,2)$}, Ann. of Math. (2) {\bf 175} (2012), no.~2, 977--1000,
\href{https://doi.org/10.4007/annals.2012.175.2.11}{doi:10.4007/annals.2012.175.2.11}.

\bibitem[Zh1]{Zh1} J.~Zhao:
{\itshape Multiple q-zeta functions and multiple q-polylogarithms}, Ramanujan J. {\bf 14} (2007), no.~2, 189--221,
\href{https://doi.org/10.1007/s11139-007-9025-9}{doi:10.1007/s11139-007-9025-9}.

\bibitem[Zh2]{Zh2} J. Zhao: {\itshape Multiple zeta functions, multiple polylogarithms and their special values}, Series on Number Theory and Its Applications {\bf 12}, World Scientific, Hackensack, NJ, 2016, xxi+595~pp.

\bibitem[Zh3]{Zh3} J.~Zhao:
{\itshape Uniform approach to double shuffle and duality relations of various q-analogs of multiple zeta values via Rota-Baxter algebras}, in J.~I. Burgos Gil, K.~Ebrahimi-Fard and H.~Gangl (eds.), ``Periods in Quantum Field Theory and Arithmetic'', Springer Proceedings in Mathematics \& Statistics {\bf 314}, Springer, Cham (2020), 259--292,
\href{https://doi.org/10.1007/978-3-030-37031-2_10}{doi:10.1007/978-3-030-37031-2\_10}.

\bibitem[Zh4]{Zh4} J.~Zhao:
{\itshape Effective Euler sums}, with an appendix by R.~Shi and J.~Zhao,
preprint, \href{https://arxiv.org/abs/2609.09260}{arXiv:2609.09260} (2026).

\bibitem[Zu1]{Zu1} W.~Zudilin:
{\itshape One of the numbers $\zeta(5)$, $\zeta(7)$, $\zeta(9)$, $\zeta(11)$ is irrational}, Russian Math. Surveys {\bf 56} (2001), no.~4, 774--776,
\href{https://doi.org/10.1070/RM2001v056n04ABEH000427}{doi:10.1070/RM2001v056n04ABEH000427}.

\bibitem[Zu2]{Zu2} W.~Zudilin: {\it Diophantine problems for q-zeta values}, Mat. Zametki {\bf 72} (2002), no.~6, 936--940 (Russian). Translation in Math. Notes {\bf 72} (2002), no.~6, 858--862,
\href{https://doi.org/10.1023/A:1021450231834}{doi:10.1023/A:1021450231834}.

\bibitem[Zu3]{Zu3} W.~Zudilin:
{\it Algebraic relations for multiple zeta values}, Russian Math. Surveys {\bf 58} (2003), no.~1, 1--29,
\href{https://doi.org/10.1070/RM2003v058n01ABEH000592}{doi:10.1070/RM2003v058n01ABEH000592}.

\bibitem[Zu4]{Zu4} W.~Zudilin:
{\itshape Arithmetic of linear forms involving odd zeta values}, J. Th\'eor. Nombres Bordeaux {\bf 16} (2004), no.~1, 251--291,
\href{https://doi.org/10.5802/jtnb.447}{doi:10.5802/jtnb.447}.

\bibitem[Zu5]{Zu5} W.~Zudilin:
{\itshape Multiple $q$-zeta brackets}, Mathematics {\bf 3} (2015), no.~1, 119--130,
\href{https://doi.org/10.3390/math3010119}{doi:10.3390/math3010119}.

\bibitem[Zu6]{Zu6} W.~Zudilin:
{\itshape Multiple zeta values}, tasting notes (14 April 2025).
(available at \url{https://www.math.ru.nl/~wzudilin/PS/MZV.pdf})

\end{thebibliography}
